\documentclass[12pt]{amsbook}
\usepackage{amssymb}
\usepackage[all]{xy}

\usepackage[colorlinks=true,linkcolor=magenta,citecolor=magenta]{hyperref}
\makeindex

\newtheorem{theorem}{Theorem}[chapter]
\newtheorem{comment}[theorem]{Comment}
\newtheorem{conjecture}[theorem]{Conjecture}
\newtheorem{definition}[theorem]{Definition}
\newtheorem{exercise}[theorem]{Exercise}
\newtheorem{proposition}[theorem]{Proposition}
\newtheorem{question}[theorem]{Question}

\begin{document}

\title{Linear algebra and group theory}

\author{Teo Banica}
\address{Department of Mathematics, University of Cergy-Pontoise, F-95000 Cergy-Pontoise, France. {\tt teo.banica@gmail.com}}


\begin{abstract}
This is an introduction to linear algebra and group theory. We first review the linear algebra basics, namely the determinant, the diagonalization procedure and more, and with the determinant being constructed as it should, as a signed volume. We discuss then the basic applications of linear algebra to questions in analysis. Then we get into the study of the closed groups of unitary matrices $G\subset U_N$, with some basic algebraic theory, and with a number of probability computations, in the finite group case. In the general case, where $G\subset U_N$ is compact, we explain how the Weingarten integration formula works, and we present some basic $N\to\infty$ applications.
\end{abstract}

\maketitle

\chapter*{Preface}

Linear algebra is the source of many good things in this world. First of all, everything algebra, for sure. But also geometry and analysis, because any smooth function or manifold, taken locally, perturbes a certain linear transformation of $\mathbb R^N$. And finally probability too, remember indeed that Gauss integral needed for talking about normal laws, which can only be computed by using polar coordinates and their Jacobian.

\bigskip

The purpose of this book is to talk about linear algebra in a large sense, theory and applications, at a somewhat more advanced level than the beginner one, and by insisting on beautiful things. And with some graduate level mathematics, and quantum physics too, in mind. We will particularly insist on the groups of matrices, which are extremely useful for all sorts of mathematics and physics, and which are perhaps the most beautiful topic one could study, once the basics of linear algebra and matrices understood.

\bigskip

The first half of the book is concerned with linear algebra and its applications. Part I is a quick journey through basic linear algebra, from basic definitions and fun with $2\times2$ matrices, up to the Spectral Theorem in its most general form, for the normal matrices $A\in M_N(\mathbb C)$. Among the features of our presentation, the determinant will be introduced as it should, as a signed volume of a system of vectors. And also, we will discuss all sorts of useful matrix tricks, which are more advanced, and good to know. 

\bigskip

As a continuation of this, Part II deals with various applications of linear algebra, to questions in analysis. After a quick look at differentiation and integration, which in several variables are intimately related to matrix theory, via the Jacobian, Hessian and so on, we will develop some useful probability theory, in relation with the normal and hyperspherical laws, by using spherical coordinates and their Jacobian. We will also discuss some other analytic topics, such as special matrices and spectral theory.

\bigskip

The second half of the book is concerned with matrix groups. As already mentioned, this is perhaps the most beautiful topic one could study, once the basics of linear algebra understood. The subject is however huge, and Part III will be a modest introduction to it. Our philosophy will be that of talking about all sorts of interesting closed subgroups $G\subset U_N$, finite and continuous alike, and by using very basic methods, coming from standard calculus, combinatorics and probability, for their study.

\bigskip

As a conclusion to this, the finite group case will appear to be reasonably understood, while the continuous case, not. Part IV will be dedicated to the study of the closed subgroups $G\subset U_N$, and more specifically the continuous ones, by using heavy machinery, as heavy as it gets. We will discuss here the basics of representation theory, then the existence of the Haar measure, and the Peter-Weyl theory, and then more advanced topics, such as Tannakian duality, Brauer theorems, and Weingarten calculus.

\bigskip

In the hope that you will find this book useful. At the level of things which are not done here, notable topics include the Jordan decomposition, which is the nightmare of everyone involved, teacher or student, and this remains between us, as well as some basic Lie algebra theory, which would have perfectly make sense to include, but that we preferred to replace by representation theory, and its relation with combinatorics and probability, which are somewhat more elementary, and fitting better with the rest.

\bigskip

Let us also mention that this way of presenting things has its origins in some recent research work on the quantum groups, and more specifically on the so-called easy quantum groups. The idea there is that there is no much smoothness and geometry, with the main tools belonging to combinatorics and probability. Thus, as main philosophy, the present book, while dealing with classical topics, is written with a ``quantum'' touch.

\bigskip

This book remains an introductory text, and for more, we will recommend some reading at the end. Among others, for some help with the preliminaries, you have my general mathematics book \cite{ba1}, for more linear algebra, you have my advanced linear algebra book \cite{ba2}, and for more about groups, you have my group theory book \cite{ba3}.

\bigskip

Most of this book is based on lecture notes from various classes at Cergy, and I would like to thank my students. The final part goes into research topics, and I am grateful to Beno\^it Collins, Steve Curran and Jean-Marc Schlenker, for our joint work on the subject. Many thanks go as well to my cats. There is so much to learn from them, too.

\bigskip

\

{\em Cergy, January 2026}

\smallskip

{\em Teo Banica}

\baselineskip=15.95pt
\tableofcontents
\baselineskip=14pt

\part{Linear algebra}

\ \vskip50mm

\begin{center}
{\em So close, no matter how far

Couldn't be much more from the heart

Forever trusting who we are

And nothing else matters}
\end{center}

\chapter{Real matrices}

\section*{1a. Linear maps}

We are interested in what follows in symmetries, rotations, projections and other such basic transformations, in 2, 3 or even more dimensions. Such transformations appear a bit everywhere, in physics. To be more precise, each physical problem or equation has some ``symmetries'', and exploiting these symmetries is usually a useful thing.

\bigskip

Let us start with 2 dimensions, and leave 3 and more dimensions for later. The transformations of the plane $\mathbb R^2$ that we are interested in are as follows:

\index{linear map}
\index{affine map}

\begin{definition}
A map $f:\mathbb R^2\to\mathbb R^2$ is called affine when it maps lines to lines,
$$f(tx+(1-t)y)=tf(x)+(1-t)f(y)$$
for any $x,y\in\mathbb R^2$ and any $t\in\mathbb R$. If in addition $f(0)=0$, we call $f$ linear.
\end{definition}

As a first observation, our ``maps lines to lines'' interpretation of the equation in the statement assumes that the points are degenerate lines, and this in order for our interpretation to work when $x=y$, or when $f(x)=f(y)$. Also, what we call line is not exactly a set, but rather a dynamic object, think trajectory of a point on that line. We will be back to this later, once we will know more about such maps.

\bigskip

Here are some basic examples of symmetries, all being linear in the above sense:

\index{symmetry}

\begin{proposition}
The symmetries with respect to $Ox$ and $Oy$ are:
$$\binom{x}{y}\to\binom{x}{-y}\quad,\quad
\binom{x}{y}\to\binom{-x}{y}$$
The symmetries with respect to the $x=y$ and $x=-y$ diagonals are:
$$\binom{x}{y}\to\binom{y}{x}\quad,\quad 
\binom{x}{y}\to\binom{-y}{-x}$$
All these maps are linear, in the above sense.
\end{proposition}

\begin{proof}
The fact that all these maps are linear is clear, because they map lines to lines, in our sense, and they also map $0$ to $0$. As for the explicit formulae in the statement, these are clear as well, by drawing pictures for each of the maps involved.
\end{proof}

Here are now some basic examples of rotations, once again all being linear:

\index{rotation}

\begin{proposition}
The rotations of angle $0^\circ$ and of angle $90^\circ$ are:
$$\binom{x}{y}\to\binom{x}{y}\quad,\quad
\binom{x}{y}\to\binom{-y}{x}$$
The rotations of angle $180^\circ$ and of angle $270^\circ$ are:
$$\binom{x}{y}\to\binom{-x}{-y}\quad,\quad 
\binom{x}{y}\to\binom{y}{-x}$$
All these maps are linear, in the above sense.
\end{proposition}

\begin{proof}
As before, these rotations are all linear, for obvious reasons. As for the formulae in the statement, these are clear as well, by drawing pictures.
\end{proof}

Here are some basic examples of projections, once again all being linear:

\index{projection}

\begin{proposition}
The projections on $Ox$ and $Oy$ are:
$$\binom{x}{y}\to\binom{x}{0}\quad,\quad 
\binom{x}{y}\to\binom{0}{y}$$
The projections on the $x=y$ and $x=-y$ diagonals are:
$$\binom{x}{y}\to\frac{1}{2}\binom{x+y}{x+y}\quad,\quad 
\binom{x}{y}\to\frac{1}{2}\binom{x-y}{y-x}$$
All these maps are linear, in the above sense.
\end{proposition}

\begin{proof}
Again, these projections are all linear, and the formulae are clear as well, by drawing pictures, with only the last 2 formulae needing some explanations. In what regards the projection on the $x=y$ diagonal, the picture here is as follows:
\vskip-6mm
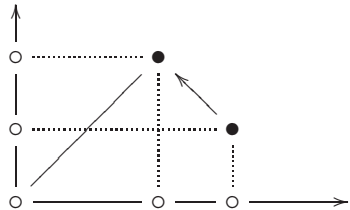
$$\xymatrix@R=15pt@C=15pt{
&\\
\circ\ar@{-}[d]\ar[u]\ar@{.}[rr]&&\bullet\ar@{.}[dd]\\
\circ\ar@{-}[d]\ar@{.}[rrr]&&&\bullet\ar[ul]\ar@{.}[d]\\
\circ\ar@{-}[uurr]\ar@{-}[rr]&\ &\circ\ar@{-}[r]&\circ\ar[rr]&&}$$

But this gives the result, since the $45^\circ$ triangle shows that this projection leaves invariant $x+y$, so we can only end up with the average $(x+y)/2$, as double coordinate. As for the projection on the $x=-y$ diagonal, the proof here is similar.
\end{proof}

Finally, we have the translations, which are as follows:

\index{translation}

\begin{proposition}
The translations are exactly the maps of the form
$$\binom{x}{y}\to\binom{x+p}{y+q}$$
with $p,q\in\mathbb R$, and these maps are all affine, in our sense.
\end{proposition}

\begin{proof}
A translation $f:\mathbb R^2\to\mathbb R^2$ is clearly affine, because it maps lines to lines. Also, such a translation is uniquely determined by the following vector:
$$f\binom{0}{0}=\binom{p}{q}$$

To be more precise, $f$ must be the map which takes a vector $\binom{x}{y}$, and adds this vector $\binom{p}{q}$ to it. But this gives the formula in the statement.
\end{proof}

Summarizing, we have many interesting examples of linear and affine maps. Let us develop now some general theory, for such maps. As a first result, we have:

\begin{theorem}
For a map $f:\mathbb R^2\to\mathbb R^2$, the following are equivalent:
\begin{enumerate}
\item $f$ is linear in our sense, mapping lines to lines, and $0$ to $0$.

\item $f$ maps sums to sums, $f(x+y)=f(x)+f(y)$, and satisfies $f(\lambda x)=\lambda f(x)$.
\end{enumerate}
\end{theorem}

\begin{proof}
This is something which comes from definitions, as follows:

\medskip

$(1)\implies(2)$ We know that $f$ satisfies the following equation, and $f(0)=0$:
$$f(tx+(1-t)y)=tf(x)+(1-t)f(y)$$

By setting $y=0$, and by using our assumption $f(0)=0$, we obtain, as desired:
$$f(tx)=tf(x)$$

As for the first condition, regarding sums, this can be established as follows:
\begin{eqnarray*}
f(x+y)
&=&f\left(2\cdot\frac{x+y}{2}\right)\\
&=&2f\left(\frac{x+y}{2}\right)\\
&=&2\cdot\frac{f(x)+f(y)}{2}\\
&=&f(x)+f(y)
\end{eqnarray*}

$(2)\implies(1)$ Conversely now, assuming that $f$ satisfies $f(x+y)=f(x)+f(y)$ and $f(\lambda x)=\lambda f(x)$, it follows that $f$ must map lines to lines, as shown by:
\begin{eqnarray*}
f(tx+(1-t)y)
&=&f(tx)+f((1-t)y)\\
&=&tf(x)+(1-t)f(y)
\end{eqnarray*} 

Also, we have $f(0)=f(2\cdot 0)=2f(0)$, which gives $f(0)=0$, as desired.
\end{proof}

The above result is very useful, and in practice, we will often use the condition (2) there, somewhat as a new definition for the linear maps. Let us record this as follows:

\begin{definition}[upgrade]
A map $f:\mathbb R^2\to\mathbb R^2$ is called:
\begin{enumerate}
\item Linear, when it satisfies $f(x+y)=f(x)+f(y)$ and $f(\lambda x)=\lambda f(x)$.

\item Affine, when it is of the form $f=g+x$, with $g$ linear, and $x\in\mathbb R^2$.
\end{enumerate}
\end{definition}

Before getting into the mathematics of linear maps, let us comment a bit more on the ``maps lines to lines'' feature of such maps. As mentioned after Definition 1.1, this requires thinking at lines as being ``dynamic'' objects, the point being that, when thinking at lines as being sets, this interpretation fails, as shown by the following map:
$$f\binom{x}{y}=\binom{x^3}{0}$$

However, in relation with all this we have the following useful result:

\begin{theorem}
For a continuous injective $f:\mathbb R^2\to\mathbb R^2$, the following are equivalent:
\begin{enumerate}
\item $f$ is affine in our sense, mapping lines to lines.

\item $f$ maps set-theoretical lines to set-theoretical lines.
\end{enumerate}
\end{theorem}

\begin{proof}
By composing $f$ with a translation, we can assume that we have $f(0)=0$. With this assumption made, the proof goes as follows:

\medskip

$(1)\implies(2)$ This is clear from definitions.

\medskip

$(2)\implies(1)$ Let us first prove that we have $f(x+y)=f(x)+f(y)$. We do this first in the case where our vectors are not proportional, $x\not\sim y$. In this case we have a proper parallelogram $(0,x,y,x+y)$, and since $f$ was assumed to be injective, it must map parallel lines to parallel lines, and so must map our parallelogram into a parallelogram $(0,f(x),f(y),f(x+y))$. But this latter parallelogram shows that we have:
$$f(x+y)=f(x)+f(y)$$

In the remaining case where our vectors are proportional, $x\sim y$, we can pick a sequence $x_n\to x$ satisfying $x_n\not\sim y$ for any $n$, and we obtain, as desired:
\begin{eqnarray*}
x_n\to x,x_n\not\sim y,\forall n
&\implies&f(x_n+y)=f(x_n)+f(y),\forall n\\
&\implies&f(x+y)=f(x)+f(y)
\end{eqnarray*}

Regarding now $f(\lambda x)=\lambda f(x)$, since $f$ maps lines to lines, it must map the line $0-x$ to the line $0-f(x)$, so we have a formula as follows, for any $\lambda,x$:
$$f(\lambda x)=\varphi_x(\lambda)f(x)$$

But since $f$ maps parallel lines to parallel lines, by Thales the function $\varphi_x:\mathbb R\to\mathbb R$ does not depend on $x$. Thus, we have a formula as follows, for any $\lambda,x$:
$$f(\lambda x)=\varphi(\lambda)f(x)$$

We know that we have $\varphi(0)=0$ and $\varphi(1)=1$, and we must prove that we have $\varphi(\lambda)=\lambda$ for any $\lambda$. For this purpose, we use a trick. On one hand, we have:
$$f((\lambda+\mu)x)=\varphi(\lambda+\mu)f(x)$$

On the other hand, since $f$ maps sums to sums, we have as well:
\begin{eqnarray*}
f((\lambda+\mu)x)
&=&f(\lambda x)+f(\mu x)\\
&=&\varphi(\lambda)f(x)+\varphi(\mu)f(x)\\
&=&(\varphi(\lambda)+\varphi(\mu))f(x)
\end{eqnarray*}

Thus our rescaling function $\varphi:\mathbb R\to\mathbb R$ satisfies the following conditions:
$$\varphi(0)=0\quad,\quad \varphi(1)=1\quad,\quad\varphi(\lambda+\mu)=\varphi(\lambda)+\varphi(\mu)$$

But with these conditions in hand, it is clear that we have $\varphi(\lambda)=\lambda$, first for all the inverses of integers, $\lambda=1/n$ with $n\in\mathbb N$, then for all rationals, $\lambda\in\mathbb Q$, and finally by continuity for all reals, $\lambda\in\mathbb R$. Thus, we have proved the following formula:
$$f(\lambda x)=\lambda f(x)$$

But this finishes the proof of $(2)\implies(1)$, and we are done.
\end{proof}

All this is nice, and there are some further things that can be said, but getting to business, Definition 1.7 is what we need. Indeed, we have the following powerful result, showing that the linear/affine maps $f:\mathbb R^2\to\mathbb R^2$ are fully described by $4/6$ parameters:

\index{linear map}
\index{affine map}

\begin{theorem}
The linear maps $f:\mathbb R^2\to\mathbb R^2$ are precisely the maps of type
$$f\binom{x}{y}=\binom{ax+by}{cx+dy}$$
and the affine maps $f:\mathbb R^2\to\mathbb R^2$ are precisely the maps of type
$$f\binom{x}{y}=\binom{ax+by}{cx+dy}+\binom{p}{q}$$
with the conventions from Definition 1.7 for such maps.
\end{theorem}

\begin{proof}
Assuming that $f$ is linear in the sense of Definition 1.7, we have:
\begin{eqnarray*}
f\binom{x}{y}
&=&f\left(\binom{x}{0}+\binom{0}{y}\right)\\
&=&f\binom{x}{0}+f\binom{0}{y}\\
&=&f\left(x\binom{1}{0}\right)+f\left(y\binom{0}{1}\right)\\
&=&xf\binom{1}{0}+yf\binom{0}{1}
\end{eqnarray*}

Thus, we obtain the formula in the statement, with $a,b,c,d\in\mathbb R$ being given by:
$$f\binom{1}{0}=\binom{a}{c}\quad,\quad f\binom{0}{1}=\binom{b}{d}$$

In the affine case now, we have as extra piece of data a vector, as follows:
$$f\binom{0}{0}=\binom{p}{q}$$

Indeed, if $f:\mathbb R^2\to\mathbb R^2$ is affine, then the following map must be linear:
$$f-\binom{p}{q}:\mathbb R^2\to\mathbb R^2$$

Thus, by using the formula in (1) we obtain the result.
\end{proof}

Moving ahead now, Theorem 1.9 is all that we need for doing some non-trivial mathematics, and so in practice, that will be our new definition for the linear and affine maps. In order to simplify now all that, which might be a bit complicated to memorize, the idea will be to put our parameters $a,b,c,d$ into a matrix, in the following way:

\index{matrix}
\index{matrix multiplication}

\begin{definition}
A matrix $A\in M_2(\mathbb R)$ is an array as follows:
$$A=\begin{pmatrix}a&b\\ c&d\end{pmatrix}$$
These matrices act on the vectors in the following way,
$$\begin{pmatrix}a&b\\ c&d\end{pmatrix}\binom{x}{y}=\binom{ax+by}{cx+dy}$$
the rule being ``multiply the rows of the matrix by the vector''.
\end{definition}

The above multiplication formula might seem a bit complicated, at a first glance, but it is not. Here is an example for it, quickly worked out:
$$\begin{pmatrix}1&2\\ 5&6\end{pmatrix}\binom{3}{1}=\binom{1\cdot 3+2\cdot 1}{5\cdot 3+6\cdot 1}=\binom{5}{21}$$

As already mentioned, all this comes from our findings from Theorem 1.9. Indeed, with the above multiplication convention for matrices and vectors, we can turn Theorem 1.9 into something much simpler, and better-looking, as follows:

\index{linear map}
\index{affine map}

\begin{theorem}
The linear maps $f:\mathbb R^2\to\mathbb R^2$ are precisely the maps of type
$$f(v)=Av$$
and the affine maps $f:\mathbb R^2\to\mathbb R^2$ are precisely the maps of type
$$f(v)=Av+w$$
with $A$ being a $2\times 2$ matrix, and with $v,w\in\mathbb R^2$ being vectors, written vertically.
\end{theorem}

\begin{proof}
With the above conventions, the formulae in Theorem 1.9 read:
$$f\binom{x}{y}=\begin{pmatrix}a&b\\ c&d\end{pmatrix}\binom{x}{y}$$
$$f\binom{x}{y}=\begin{pmatrix}a&b\\ c&d\end{pmatrix}\binom{x}{y}+\binom{p}{q}$$

Thus, we are led to the conclusions in the statement.
\end{proof}

Before going further, let us discuss some examples. First, we have:

\index{symmetry}

\begin{proposition}
The symmetries with respect to $Ox$ and $Oy$ are given by
$$\begin{pmatrix}1&0\\0&-1\end{pmatrix}\binom{x}{y}\quad,\quad
\begin{pmatrix}-1&0\\0&1\end{pmatrix}\binom{x}{y}$$
and the symmetries with respect to the $x=y$ and $x=-y$ diagonals are given by
$$\begin{pmatrix}0&1\\1&0\end{pmatrix}\binom{x}{y}\quad,\quad 
\begin{pmatrix}0&-1\\-1&0\end{pmatrix}\binom{x}{y}$$
with our conventions above for the matrix multiplication.
\end{proposition}

\begin{proof}
According to Proposition 1.2, the above transformations map $\binom{x}{y}$ to:
$$\binom{x}{-y}
\quad,\quad\binom{-x}{y}
\quad,\quad\binom{y}{x}
\quad,\quad\binom{-y}{-x}$$

But this gives the formulae in the statement, by guessing in each case the matrix which does the job, in the obvious way.
\end{proof}

Regarding now the basic rotations, we have here:

\index{rotation}

\begin{proposition}
The rotations of angle $0^\circ$ and of angle $90^\circ$ are given by
$$\begin{pmatrix}1&0\\0&1\end{pmatrix}\binom{x}{y}\quad,\quad 
\begin{pmatrix}0&-1\\1&0\end{pmatrix}\binom{x}{y}$$
and the rotations of angle $180^\circ$ and of angle $270^\circ$ are given by
$$\begin{pmatrix}-1&0\\0&-1\end{pmatrix}\binom{x}{y}\quad,\quad 
\begin{pmatrix}0&1\\-1&0\end{pmatrix}\binom{x}{y}$$
with our conventions above for the matrix multiplication.
\end{proposition}

\begin{proof}
As before, but by using Proposition 1.3, the vector $\binom{x}{y}$ maps to:
$$\binom{x}{y}
\quad,\quad\binom{-y}{x}
\quad,\quad\binom{-x}{-y}
\quad,\quad\binom{y}{-x}$$

But this gives the formulae in the statement, again by guessing the matrix.
\end{proof}

Finally, regarding the basic projections, we have here:

\index{projection}

\begin{proposition}
The projections on $Ox$ and $Oy$ are given by
$$\begin{pmatrix}1&0\\0&0\end{pmatrix}\binom{x}{y}\quad,\quad 
\begin{pmatrix}0&0\\0&1\end{pmatrix}\binom{x}{y}$$
and the projections on the $x=y$ and $x=-y$ diagonals are given by
$$\frac{1}{2}\begin{pmatrix}1&1\\1&1\end{pmatrix}\binom{x}{y}\quad,\quad 
\frac{1}{2}\begin{pmatrix}1&-1\\-1&1\end{pmatrix}\binom{x}{y}$$
with our conventions above for the matrix multiplication.
\end{proposition}

\begin{proof}
As before, but according now to Proposition 1.4, the vector $\binom{x}{y}$ maps to:
$$\binom{x}{0}
\quad,\quad\binom{0}{y}
\quad,\quad\frac{1}{2}\binom{x+y}{x+y}
\quad,\quad\frac{1}{2}\binom{x-y}{y-x}$$

But this gives the formulae in the statement, as usual by guessing the matrix.
\end{proof}

In addition to the above transformations, there are many other examples. We have for instance the null transformation, which is given by:
$$\begin{pmatrix}0&0\\0&0\end{pmatrix}\binom{x}{y}=\binom{0}{0}$$

Here is now a more bizarre map, but which can still be understood, however, as being the map which ``switches the coordinates, then kills the second one'':
$$\begin{pmatrix}0&1\\0&0\end{pmatrix}\binom{x}{y}=\binom{y}{0}$$

Even more bizarrely now, here is a certain linear map, whose interpretation is more complicated, and is left to you, reader:
$$\begin{pmatrix}1&1\\0&0\end{pmatrix}\binom{x}{y}=\binom{x+y}{0}$$

And here is another linear map, which once again, being something geometric, in 2 dimensions, can definitely be understood, at least in theory:
$$\begin{pmatrix}1&1\\0&1\end{pmatrix}\binom{x}{y}=\binom{x+y}{y}$$

Let us discuss now the computation of the arbitrary symmetries, rotations and projections. We begin with the rotations, whose formula is a must-know:

\index{rotation}

\begin{theorem}
The rotation of angle $t\in\mathbb R$ is given by the matrix
$$R_t=\begin{pmatrix}\cos t&-\sin t\\ \sin t&\cos t\end{pmatrix}$$
depending on $t\in\mathbb R$ taken modulo $2\pi$.
\end{theorem}

\begin{proof}
The rotation being linear, it must correspond to a certain matrix:
$$R_t=\begin{pmatrix}a&b\\ c&d\end{pmatrix}$$

We can guess this matrix, via its action on the basic coordinate vectors $\binom{1}{0}$ and $\binom{0}{1}$. Indeed, a quick picture shows that we must have:
$$\begin{pmatrix}a&b\\ c&d\end{pmatrix}\begin{pmatrix}1\\ 0\end{pmatrix}=
\begin{pmatrix}\cos t\\ \sin t\end{pmatrix}$$

Also, by paying attention to positives and negatives, we must have:
$$\begin{pmatrix}a&b\\ c&d\end{pmatrix}\begin{pmatrix}0\\ 1\end{pmatrix}=
\begin{pmatrix}-\sin t\\ \cos t\end{pmatrix}$$

Guessing now the matrix is not complicated, because the first equation gives us the first column, and the second equation gives us the second column:
$$\binom{a}{c}=\begin{pmatrix}\cos t\\ \sin t\end{pmatrix}\quad,\quad 
\binom{b}{d}=\begin{pmatrix}-\sin t\\ \cos t\end{pmatrix}$$

Thus, we can just put together these two vectors, and we obtain our matrix.
\end{proof}

Regarding now the symmetries, the formula here is as follows:

\index{symmetry}

\begin{theorem}
The symmetry with respect to the $Ox$ axis rotated by an angle $t/2\in\mathbb R$ is given by the matrix
$$S_t=\begin{pmatrix}\cos t&\sin t\\ \sin t&-\cos t\end{pmatrix}$$
depending on $t\in\mathbb R$ taken modulo $2\pi$.
\end{theorem}

\begin{proof}
As before, we can guess the matrix via its action on the basic coordinate vectors $\binom{1}{0}$ and $\binom{0}{1}$. A quick picture shows that we must have:
$$\begin{pmatrix}a&b\\ c&d\end{pmatrix}\begin{pmatrix}1\\ 0\end{pmatrix}=
\begin{pmatrix}\cos t\\ \sin t\end{pmatrix}$$

Also, by paying attention to positives and negatives, we must have:
$$\begin{pmatrix}a&b\\ c&d\end{pmatrix}\begin{pmatrix}0\\ 1\end{pmatrix}=
\begin{pmatrix}\sin t\\-\cos t\end{pmatrix}$$

Guessing now the matrix is not complicated, because we must have:
$$\binom{a}{c}=\begin{pmatrix}\cos t\\ \sin t\end{pmatrix}\quad,\quad
\binom{b}{d}=\begin{pmatrix}\sin t\\-\cos t\end{pmatrix}$$

Thus, we can just put together these two vectors, and we obtain our matrix.
\end{proof}

Finally, regarding the projections, the formula here is as follows:

\index{projection}

\begin{theorem}
The projection on the $Ox$ axis rotated by an angle $t/2\in\mathbb R$ is given by the matrix
$$P_t=\frac{1}{2}\begin{pmatrix}1+\cos t&\sin t\\ \sin t&1-\cos t\end{pmatrix}$$
depending on $t\in\mathbb R$ taken modulo $2\pi$.
\end{theorem}

\begin{proof}
We will need here some trigonometry, and more precisely the formulae for the duplication of the angles. Regarding the sine, the formula here is:
$$\sin(2t)=2\sin t\cos t$$

Regarding the cosine, we have here 3 equivalent formulae, as follows: 
\begin{eqnarray*}
\cos(2t)
&=&\cos^2t-\sin^2t\\
&=&2\cos^2t-1\\
&=&1-2\sin^2t
\end{eqnarray*}

Getting back now to our problem, some quick pictures, using similarity of triangles, and then the above trigonometry formulae, show that we must have:
$$P_t\begin{pmatrix}1\\ 0\end{pmatrix}
=\cos\frac{t}{2}\binom{\cos\frac{t}{2}}{\sin\frac{t}{2}}
=\frac{1}{2}\begin{pmatrix}1+\cos t\\ \sin t\end{pmatrix}$$
$$P_t\begin{pmatrix}0\\ 1\end{pmatrix}
=\sin\frac{t}{2}\binom{\cos\frac{t}{2}}{\sin\frac{t}{2}}
=\frac{1}{2}\begin{pmatrix}\sin t\\1-\cos t\end{pmatrix}$$

Now by putting together these two vectors, and we obtain our matrix.
\end{proof}

\section*{1b. Matrix calculus}

In order to formulate now our second theorem, dealing with compositions of maps, let us make the following multiplication convention, between matrices and matrices:
$$\begin{pmatrix}a&b\\ c&d\end{pmatrix}
\begin{pmatrix}p&q\\ r&s\end{pmatrix}
=\begin{pmatrix}ap+br&aq+bs\\ cp+dr&cq+ds\end{pmatrix}$$

This might look a bit complicated, but as before, in what was concerning multiplying matrices and vectors, the idea is very simple, namely ``multiply the rows of the first matrix by the columns of the second matrix''. With this convention, we have:

\index{matrix multiplication}
\index{composition of linear maps}

\begin{theorem}
If we denote by $f_A:\mathbb R^2\to\mathbb R^2$ the linear map associated to a matrix $A$, given by the formula 
$$f_A(v)=Av$$
then we have the following multiplication formula for such maps:
$$f_Af_B=f_{AB}$$
That is, the composition of linear maps corresponds to the multiplication of matrices.
\end{theorem}

\begin{proof}
We want to prove that we have the following formula, valid for any two matrices $A,B\in M_2(\mathbb R)$, and any vector $v\in\mathbb R^2$:
$$A(Bv)=(AB)v$$

For this purpose, let us write our matrices and vector as follows:
$$A=\begin{pmatrix}a&b\\ c&d\end{pmatrix}\quad,\quad B=\begin{pmatrix}p&q\\ r&s\end{pmatrix}\quad,\quad v=\binom{x}{y}$$

The formula that we want to prove becomes:
$$\begin{pmatrix}a&b\\ c&d\end{pmatrix}
\left[
\begin{pmatrix}p&q\\ r&s\end{pmatrix}
\binom{x}{y}
\right]=
\left[\begin{pmatrix}a&b\\ c&d\end{pmatrix}\begin{pmatrix}p&q\\ r&s\end{pmatrix}\right]
\binom{x}{y}$$

But this is the same as saying that:
$$\begin{pmatrix}a&b\\ c&d\end{pmatrix}
\binom{px+qy}{rx+sy}=
\begin{pmatrix}ap+br&aq+bs\\ cp+dr&cq+ds\end{pmatrix}
\binom{x}{y}$$

And this latter formula does hold indeed, because on both sides we get:
$$\binom{apx+aqy+brx+bsy}{cpx+cqy+drx+dsy}$$

Thus, we have proved the result.
\end{proof}

As a verification for the above result, let us compose two rotations. The computation here is as follows, yieding a rotation, as it should, and of the correct angle:
\begin{eqnarray*}
R_sR_t
&=&\begin{pmatrix}\cos s&-\sin s\\ \sin s&\cos s\end{pmatrix}\begin{pmatrix}\cos t&-\sin t\\ \sin t&\cos t\end{pmatrix}\\
&=&\begin{pmatrix}\cos s\cos t-\sin s\sin t&&-\cos s\sin t-\sin t\cos s\\
\sin s\cos t+\cos s\sin t&&-\sin s\sin t+\cos s\cos t\end{pmatrix}\\
&=&\begin{pmatrix}\cos(s+t)&-\sin(s+t)\\ \sin(s+t)&\cos(s+t)\end{pmatrix}\\
&=&R_{s+t}
\end{eqnarray*}

We are ready now to pass to 3 dimensions. The idea is to select from what we learned in 2 dimensions, nice results only, and generalize to 3 dimensions. We obtain:

\index{linear map}
\index{affine map}
\index{composition of linear maps}

\begin{theorem}
Consider a map $f:\mathbb R^3\to\mathbb R^3$.
\begin{enumerate}
\item $f$ is linear when it is of the form $f(v)=Av$, with $A\in M_3(\mathbb R)$.

\item $f$ is affine when $f(v)=Av+w$, with $A\in M_3(\mathbb R)$ and $w\in\mathbb R^3$.

\item We have the composition formula $f_Af_B=f_{AB}$, similar to the $2D$ one.
\end{enumerate}
\end{theorem}

\begin{proof}
Here (1,2) can be proved exactly as in the 2D case, with the multiplication convention being as usual, ``multiply the rows of the matrix by the vector'':
$$\begin{pmatrix}a&b&c\\ d&e&f\\ g&h&i\end{pmatrix}
\begin{pmatrix}x\\ y\\ z\end{pmatrix}
=\begin{pmatrix}ax+by+cz\\ dx+ey+fz\\ gx+hy+iz\end{pmatrix}$$

As for (3), once again the 2D idea applies, with the same product rule, ``multiply the rows of the first matrix by the columns of the second matrix'':
$$\begin{pmatrix}a&b&c\\ d&e&f\\ g&h&i\end{pmatrix}
\begin{pmatrix}p&q&r\\ s&t&u\\ v&w&x\end{pmatrix}\\
=\begin{pmatrix}
ap+bs+cv&aq+bt+cw&ar+bu+cx\\
dp+es+fv&dq+et+fw&dr+eu+fx\\
gp+hs+iv&gq+ht+iw&gr+hu+ix
\end{pmatrix}$$

Thus, we proved our theorem. Of course, we are going a bit fast here, but we will discuss all this in detail, right next, directly in arbitrary $N$ dimensions.
\end{proof}

We are now ready to discuss 4 and more dimensions. Before doing so, let us point out however that the maps of type $f:\mathbb R^3\to\mathbb R^2$, or $f:\mathbb R\to\mathbb R^2$, and so on, are not covered by our results. Since there are many interesting such maps, say obtained by projecting and then rotating, and so on, we will be interested here in the maps $f:\mathbb R^N\to\mathbb R^M$.

\bigskip

A bit of thinking suggests that such maps should come from the $M\times N$ matrices. Indeed, this is what happens at $M=N=2$ and $M=N=3$, of course. But this happens as well at $N=1$, because a linear map $f:\mathbb R\to\mathbb R^M$ can only be something of the form $f(\lambda)=\lambda v$, with $v\in\mathbb R^M$, and $v\in\mathbb R^M$ means that $v$ is a $M\times 1$ matrix. So, let us start with the product rule for the $M\times N$ matrices, which is as follows:

\index{matrix multiplication}
\index{product of matrices}
\index{rectangular matrix}

\begin{definition}
We can multiply the $M\times N$ matrices with $N\times K$ matrices,
$$\begin{pmatrix}
a_{11}&\ldots&a_{1N}\\
\vdots&&\vdots\\
a_{M1}&\ldots&a_{MN}
\end{pmatrix}
\begin{pmatrix}
b_{11}&\ldots&b_{1K}\\
\vdots&&\vdots\\
b_{N1}&\ldots&b_{NK}
\end{pmatrix}$$
the product being the $M\times K$ matrix given by the following formula,
$$\begin{pmatrix}
a_{11}b_{11}+\ldots+a_{1N}b_{N1}&\ldots\ldots&a_{11}b_{1K}+\ldots+a_{1N}b_{NK}\\
\vdots&&\vdots\\
\vdots&&\vdots\\
a_{M1}b_{11}+\ldots+a_{MN}b_{N1}&\ldots\ldots&a_{M1}b_{1K}+\ldots+a_{MN}b_{NK}
\end{pmatrix}$$
obtained via the usual rule ``multiply rows by columns''.
\end{definition}

Observe that this formula generalizes all the multiplication rules that we have been using so far, between various types of matrices and vectors. Thus, in practice, we can simply forget all the previous multiplication rules, and simply memorize this one.

\bigskip

In case the above formula looks hard to memorize, here is an alternative formulation of it, which is simpler and more powerful, by using the standard algebraic notation for the matrices, $A=(A_{ij})$, that we will heavily use, in what follows:

\begin{proposition}
The matrix multiplication is given by formula
$$(AB)_{ij}=\sum_kA_{ik}B_{kj}$$
with $A_{ij}$ standing for the entry of $A$ at row $i$ and column $j$.
\end{proposition}

\begin{proof}
This is indeed just a shorthand for the formula in Definition 1.20, by following the rule there, namely ``multiply the rows of $A$ by the columns of $B$''.  
\end{proof}

As an illustration for the power of the convention in Proposition 1.21, we have:

\begin{proposition}
We have the following formula, valid for any matrices $A,B,C$,
$$(AB)C=A(BC)$$
provided that the sizes of our matrices $A,B,C$ fit.
\end{proposition}

\begin{proof}
We have the following computation, using indices as above:
$$((AB)C)_{ij}
=\sum_k(AB)_{ik}C_{kj}
=\sum_{kl}A_{il}B_{lk}C_{kj}$$

On the other hand, we have as well the following computation:
$$(A(BC))_{ij}
=\sum_lA_{il}(BC)_{lj}
=\sum_{kl}A_{il}B_{lk}C_{kj}$$

Thus we have $(AB)C=A(BC)$, and we have proved our result.
\end{proof}

With this, we can now talk about linear maps between spaces of arbitrary dimension, generalizing what we have been doing so far. The main result here is as follows:

\index{linear map}
\index{affine map}
\index{composition of linear maps}

\begin{theorem}
Consider a map $f:\mathbb R^N\to\mathbb R^M$.
\begin{enumerate}
\item $f$ is linear when it is of the form $f(v)=Av$, with $A\in M_{M\times N}(\mathbb R)$.

\item $f$ is affine when $f(v)=Av+w$, with $A\in M_{M\times N}(\mathbb R)$ and $w\in\mathbb R^M$.

\item We have the composition formula $f_Af_B=f_{AB}$, whenever the sizes fit.
\end{enumerate}
\end{theorem}

\begin{proof}
We already know that this happens at $M=N=2$, and at $M=N=3$ as well. In general, the proof is similar, by doing some elementary computations.
\end{proof}

\index{identity matrix}
\index{null matrix}

As a first example here, we have the identity matrix, acting as the identity:
$$\begin{pmatrix}
1&&0\\
&\ddots\\
0&&1\end{pmatrix}
\begin{pmatrix}x_1\\ \vdots\\ x_N\end{pmatrix}
=\begin{pmatrix}x_1\\ \vdots\\ x_N\end{pmatrix}$$

Along the same lines, we have as well the null matrix $(0)_{ij}$, acting as the null map, $x\to0$. Here is now an important result, providing us with many examples:

\index{diagonal matrix}

\begin{proposition}
The diagonal matrices act as follows,
$$\begin{pmatrix}
\lambda_1&&0\\
&\ddots\\
0&&\lambda_N\end{pmatrix}
\begin{pmatrix}x_1\\ \vdots\\ x_N\end{pmatrix}
=\begin{pmatrix}\lambda_1x_1\\ \vdots\\ \lambda_Nx_N\end{pmatrix}$$
by multiplying each vector entry by a certain scalar.
\end{proposition}

\begin{proof}
This is clear, indeed, from definitions.
\end{proof}

As a more specialized example now, we have:

\index{flat matrix}
\index{all-one matrix}
\index{scalar product}
\index{orthogonal projection}
\index{all-one vector}

\begin{proposition}
The flat matrix, which is as follows,
$$\mathbb I_N=\begin{pmatrix}
1&\ldots&1\\
\vdots&&\vdots\\
1&\ldots&1\end{pmatrix}$$
acts via $N$ times the projection on the all-one vector.
\end{proposition}

\begin{proof}
The flat matrix acts in the following way:
$$\begin{pmatrix}
1&\ldots&1\\
\vdots&&\vdots\\
1&\ldots&1\end{pmatrix}
\begin{pmatrix}x_1\\ \vdots\\ x_N\end{pmatrix}
=\begin{pmatrix}x_1+\ldots+x_N\\ \vdots\\ x_1+\ldots+x_N\end{pmatrix}$$

Thus, in terms of the matrix $P=\mathbb I_N/N$, we have the following formula:
$$P\begin{pmatrix}x_1\\ \vdots\\ x_N\end{pmatrix}
=\frac{x_1+\ldots+x_N}{N}\begin{pmatrix}1\\ \vdots\\ 1\end{pmatrix}$$

Now since the linear map $f(x)=Px$ satisfies $f^2=f$, and since $Im(f)$ consists of the scalar multiples of the all-one vector $\xi\in\mathbb R^N$, we conclude that $f$ is a projection on $\mathbb R\xi$. Also, with the standard scalar product convention $<x,y>=\sum x_iy_i$, we have:
\begin{eqnarray*}
<f(x)-x,\xi>
&=&<f(x),\xi>-<x,\xi>\\
&=&\frac{\sum x_i}{N}\times N-\sum x_i\\
&=&0
\end{eqnarray*}

Thus, our projection is indeed orthogonal, and we are done. And more on this later in this chapter, when systematically discussing scalar products and orthogonality.
\end{proof}

\section*{1c. Diagonalization}

Let us develop now some general theory for the square matrices. We will need the following standard result, regarding the changes of coordinates in $\mathbb R^N$:

\index{basis}
\index{invertible matrix}
\index{bijective linear map}

\begin{theorem}
For a system $\{v_1,\ldots,v_N\}\subset\mathbb R^N$, the following are equivalent:
\begin{enumerate}
\item The vectors $v_i$ form a basis of $\mathbb R^N$, in the sense that each vector $x\in\mathbb R^N$ can be written in a unique way as a linear combination of these vectors:
$$x=\sum\lambda_iv_i$$

\item The following linear map associated to these vectors is bijective:
$$f:\mathbb R^N\to\mathbb R^N\quad,\quad 
\lambda\to\sum\lambda_iv_i$$

\item The matrix formed by these vectors, regarded as usual as column vectors,
$$P=[v_1,\ldots,v_N]\in M_N(\mathbb R)$$
is invertible, with respect to the usual multiplication of the matrices.
\end{enumerate}
\end{theorem}

\begin{proof}
Here the equivalence $(1)\iff(2)$ is clear from definitions, and the equivalence $(2)\iff(3)$ is clear as well, because we have $f(x)=Px$.
\end{proof}

Getting back now to the matrices, as an important definition, we have:

\index{eigenvalue}
\index{eigenvector}
\index{diagonalizable matrix}

\begin{definition}
Let $A\in M_N(\mathbb R)$ be a square matrix. We say that $v\in\mathbb R^N$ is an eigenvector of $A$, with corresponding eigenvalue $\lambda\in\mathbb R^N$, when:
$$Av=\lambda v$$
Also, we say that $A$ is diagonalizable when $\mathbb R^N$ has a basis formed by eigenvectors of $A$.
\end{definition}

We will see in a moment examples of eigenvectors and eigenvalues, and of diagonalizable matrices. However, even before seeing the examples, it is quite clear that these are key notions. Indeed, for a matrix $A\in M_N(\mathbb R)$, being diagonalizable is the best thing that can happen, because in this case, once the basis changed, $A$ becomes diagonal.

\bigskip

To be more precise here, we have the following result:

\index{change of basis}
\index{eigenvector basis}

\begin{proposition}
Assuming that $A\in M_N(\mathbb R)$ is diagonalizable, we have the formula
$$A=\begin{pmatrix}
\lambda_1\\
&\ddots\\
&&\lambda_N
\end{pmatrix}$$
with respect to the basis $\{v_1,\ldots,v_N\}$ of $\mathbb R^N$ consisting of eigenvectors of $A$.
\end{proposition}

\begin{proof}
This is clear from the definition of eigenvalues and eigenvectors, and from the formula of linear maps associated to diagonal matrices, from Proposition 1.24.
\end{proof}

Here is an equivalent form of the above result, which is often used in practice, when we prefer not to change the basis, and stay with the usual basis of $\mathbb R^N$:

\index{change of basis}
\index{diagonalization}
\index{passage matrix}
\index{diagonal form}

\begin{theorem}
Assuming that $A\in M_N(\mathbb R)$ is diagonalizable, with
$$v_1,\ldots,v_N\in\mathbb R^N\quad,\quad 
\lambda_1,\ldots,\lambda_N\in\mathbb R$$
as eigenvectors and corresponding eigenvalues, we have the formula
$$A=PDP^{-1}$$
with the matrices $P,D\in M_N(\mathbb R)$ being given by the formulae
$$P=[v_1,\ldots,v_N]\quad,\quad 
D=diag(\lambda_1,\ldots,\lambda_N)$$
and respectively called passage matrix, and diagonal form of $A$.
\end{theorem}

\begin{proof}
This can be viewed in two possible ways, as follows:

\medskip

(1) As already mentioned, with respect to the basis $v_1,\ldots,v_N\in\mathbb R^N$ formed by the eigenvectors, our matrix $A$ is given by:
$$A=\begin{pmatrix}
\lambda_1\\
&\ddots\\
&&\lambda_N
\end{pmatrix}$$

But this corresponds precisely to the formula $A=PDP^{-1}$ from the statement, with $P$ and its inverse appearing there due to our change of basis.

\medskip

(2) We can equally establish the formula in the statement by a direct computation. Indeed, we have $Pe_i=v_i$, where $\{e_1,\ldots,e_N\}$ is the standard basis of $\mathbb R^N$, and so:
$$APe_i
=Av_i
=\lambda_iv_i$$

On the other hand, once again by using $Pe_i=v_i$, we have as well:
$$PDe_i
=P\lambda_ie_i
=\lambda_iPe_i
=\lambda_iv_i$$

Thus we have $AP=PD$, and so $A=PDP^{-1}$, as claimed.
\end{proof}

Let us discuss now some basic examples, namely the rotations, symmetries and projections in 2 dimensions. The situation is very simple for the projections, as follows:

\index{projection}

\begin{proposition}
The projection on the $Ox$ axis rotated by an angle $t/2\in\mathbb R$,
$$P_t=\frac{1}{2}\begin{pmatrix}1+\cos t&\sin t\\ \sin t&1-\cos t\end{pmatrix}$$
is diagonalizable, its diagonal form being as follows,
$$P_t\sim\begin{pmatrix}1&0\\0&0\end{pmatrix}$$
and this regardless of the value of the angle $t/2$.
\end{proposition}

\begin{proof}
This is clear, because if we denote by $L$ the line where our projection projects, we can pick any vector $v\in L$, and this will be an eigenvector with eigenvalue 1, and then pick any vector $w\in L^\perp$, and this will be an eigenvector with eigenvalue 0. Thus, even without computations, we are led to the conclusion in the statement.
\end{proof}

The computation for the symmetries is similar, as follows:

\index{symmetry}

\begin{proposition}
The symmetry with respect to the $Ox$ axis rotated by $t/2\in\mathbb R$, 
$$S_t=\begin{pmatrix}\cos t&\sin t\\ \sin t&-\cos t\end{pmatrix}$$
is diagonalizable, its diagonal form being as follows,
$$S_t\sim\begin{pmatrix}1&0\\0&-1\end{pmatrix}$$
and this regardless of the value of the angle $t/2$.
\end{proposition}

\begin{proof}
This is again clear, because if we denote by $L$ the line with respect to which our symmetry symmetrizes, we can pick any vector $v\in L$, and this will be an eigenvector with eigenvalue 1, and then pick any vector $w\in L^\perp$, and this will be an eigenvector with eigenvalue $-1$. Thus, we are led to the conclusion in the statement.
\end{proof}

Regarding now the rotations, here the situation is different, as follows:

\index{rotation}

\begin{proposition}
The rotation of angle $t\in[0,2\pi)$, given by the formula
$$R_t=\begin{pmatrix}\cos t&-\sin t\\ \sin t&\cos t\end{pmatrix}$$
is diagonal at $t=0,\pi$, and is not diagonalizable at $t\neq0,\pi$.
\end{proposition}

\begin{proof}
The first assertion is clear, because at $t=0,\pi$ the rotations are:
$$R_0=\begin{pmatrix}1&0\\0&1\end{pmatrix}\quad,\quad 
R_\pi=\begin{pmatrix}-1&0\\0&-1\end{pmatrix}$$

As for the rotations of angle $t\neq0,\pi$, these clearly cannot have eigenvectors.
\end{proof}

Finally, here is one more example, which is the most important of them all:

\begin{theorem}
The following matrix is not diagonalizable,
$$J=\begin{pmatrix}0&1\\0&0\end{pmatrix}$$
because it has only $1$ eigenvector.
\end{theorem}

\begin{proof}
The above matrix, called $J$ en hommage to Jordan, acts as follows:
$$\begin{pmatrix}0&1\\0&0\end{pmatrix}\binom{x}{y}=\binom{y}{0}$$

Thus the eigenvector/eigenvalue equation $Jv=\lambda v$ reads:
$$\binom{y}{0}=\binom{\lambda x}{\lambda y}$$

We have then two cases, depending on $\lambda$, as follows, which give the result:

\medskip

(1) For $\lambda\neq0$ we must have $y=0$, coming from the second row, and so $x=0$ as well, coming from the first row, so we have no nontrivial eigenvectors. 

\medskip

(2) As for the case $\lambda=0$, here we must have $y=0$, coming from the first row, and so the eigenvectors here are the vectors of the form $\binom{x}{0}$. 
\end{proof}

\section*{1d. Scalar products}

In order to discuss some interesting examples of matrices, and their diagonalization, in arbitrary dimensions, we will need the following standard fact:

\index{scalar product}
\index{transpose matrix}

\begin{proposition}
Consider the scalar product on $\mathbb R^N$, given by:
$$<x,y>=\sum_ix_iy_i$$ 
We have then the following formula, valid for any vectors $x,y$ and any matrix $A$,
$$<Ax,y>=<x,A^ty>$$
with $A^t$ being the transpose matrix, $(A^t)_{ij}=A_{ji}$.
\end{proposition}

\begin{proof}
By linearity, it is enough to prove the above formula on the standard basis vectors $e_1,\ldots,e_N$ of $\mathbb R^N$. Thus, we want to prove that for any $i,j$ we have:
$$<Ae_j,e_i>=<e_j,A^te_i>$$

The scalar product being symmetric, this is the same as proving that:
$$<Ae_j,e_i>=<A^te_i,e_j>$$

On the other hand, for any matrix $M$ we have the following formula:
$$M_{ij}=<Me_j,e_i>$$

We conclude that the formula to be proved simply reads:
$$A_{ij}=(A^t)_{ji}$$

But this precisely the definition of $A^t$, and we are done.
\end{proof}

With this, we can develop some theory. We first have:

\index{orthogonal projection}

\begin{theorem}
The orthogonal projections are the matrices satisfying:
$$P^2=P^t=P$$
These projections are diagonalizable, with eigenvalues $0,1$.
\end{theorem}

\begin{proof}
It is obvious that a linear map $f(x)=Px$ is a projection precisely when:
$$P^2=P$$

In order now for this projection to be an orthogonal projection, the condition to be satisfied can be written and then processed as follows:
\begin{eqnarray*}
<Px-Py,Px-x>=0
&\iff&<x-y,P^tPx-P^tx>=0\\
&\iff&P^tPx-P^tx=0\\
&\iff&P^tP-P^t=0
\end{eqnarray*}

Thus we must have $P^t=P^tP$. Now observe that by transposing, we have as well:
$$P
=(P^tP)^t
=P^t(P^t)^t
=P^tP$$

Thus we must have $P=P^t$, as claimed. Finally, regarding the diagonalization assertion, this is clear by taking a basis of $Im(f)$, which consists of $1$-eigenvectors, and then completing with 0-eigenvectors, which can be found inside the orthogonal of $Im(f)$.
\end{proof}

Here is now a key computation of such projections:

\index{rank 1 projection}
\index{length of vector}

\begin{theorem}
The rank $1$ projections are given by the formula
$$P_x=\frac{1}{||x||^2}(x_ix_j)_{ij}$$
where the constant, $||x||=\sqrt{\sum_ix_i^2}$, is the length of the vector.
\end{theorem}

\begin{proof}
Consider a vector $y\in\mathbb R^N$. Its projection on $\mathbb Rx$ must be a certain multiple of $x$, and we are led in this way to the following formula:
$$P_xy
=\frac{<y,x>}{<x,x>}\,x
=\frac{1}{||x||^2}<y,x>x$$

With this in hand, we can now compute the entries of $P_x$, as follows:
\begin{eqnarray*}
(P_x)_{ij}
&=&<P_xe_j,e_i>\\
&=&\frac{1}{||x||^2}<e_j,x><x,e_i>\\
&=&\frac{x_jx_i}{||x||^2}
\end{eqnarray*}

Thus, we are led to the formula in the statement.
\end{proof}

As an application, we can recover a result that we already know, namely:

\index{projections}

\begin{proposition}
In $2$ dimensions, the rank $1$ projections, which are the projections on the $Ox$ axis rotated by an angle $t/2\in[0,\pi)$, are given by the following formula:
$$P_t=\frac{1}{2}\begin{pmatrix}1+\cos t&\sin t\\ \sin t&1-\cos t\end{pmatrix}$$
Together with the following two matrices, which are the rank $0$ and $2$ projections in $\mathbb R^2$,
$$0=\begin{pmatrix}0&0\\ 0&0\end{pmatrix}\quad,\quad
1=\begin{pmatrix}1&1\\ 1&1\end{pmatrix}$$
these are all the projections in $2$ dimensions.
\end{proposition}

\begin{proof}
The first assertion can be deduced from the general formula in Theorem 1.36, by plugging in the following vector, depending on a parameter $s\in[0,\pi)$:
$$x=\binom{\cos s}{\sin s}$$

Indeed, we obtain in this way the following matrix, which with $t=2s$ is the one in the statement, via the standard trigonometry formulae for the doubles of angles:
$$P_{2s}=\begin{pmatrix}\cos^2s&\cos s\sin s\\ \cos s\sin s&\sin^2 s\end{pmatrix}$$

As for the second assertion, this is clear from the first one, because outside rank 1 we can only have rank 0 or rank 2, corresponding to the matrices in the statement.
\end{proof}

Here is another interesting application, this time in $N$ dimensions:

\index{flat matrix}
\index{all-one matrix}
\index{all-one vector}

\begin{proposition}
The projection on the all-$1$ vector $\xi\in\mathbb R^N$ is
$$P_\xi=\frac{1}{N}\begin{pmatrix}
1&\ldots&1\\
\vdots&&\vdots\\
1&\ldots&1\end{pmatrix}$$
with the all-$1$ matrix on the right being called the flat matrix.
\end{proposition}

\begin{proof}
As already pointed out in the proof of Proposition 1.25, the matrix in the statement acts in the following way:
$$P_\xi
\begin{pmatrix}x_1\\ \vdots\\ x_N\end{pmatrix}
=\frac{x_1+\ldots+x_N}{N}\begin{pmatrix}1\\ \vdots\\ 1\end{pmatrix}$$

Thus $P_\xi$ is indeed a projection onto $\mathbb R\xi$, and the fact that this projection is indeed the orthogonal one follows either by a direct orthogonality computation, or by using the general formula in Theorem 1.36, by plugging in the all-1 vector $\xi$.
\end{proof}

Let us discuss now, as a final topic of this chapter, the isometries of $\mathbb R^N$. We have here the following general result:

\index{isometry}
\index{distance preservation}
\index{orthogonal matrix}
\index{polarization identity}

\begin{theorem}
The linear maps $f:\mathbb R^N\to\mathbb R^N$ which are isometries, in the sense that they preserve the distances, are those coming from the matrices satisfying:
$$U^t=U^{-1}$$ 
These latter matrices are called orthogonal, and they form a set $O_N\subset M_N(\mathbb R)$ which is stable under taking compositions, and inverses.
\end{theorem}

\begin{proof}
We have several things to be proved, the idea being as follows:

\medskip

(1) We recall that we can pass from scalar products to distances, as follows:
$$||x||=\sqrt{<x,x>}$$

Conversely, we can compute the scalar products in terms of distances, by using the polarization identity, which is as follows:
\begin{eqnarray*}
||x+y||^2-||x-y||^2
&=&||x||^2+||y||^2+2<x,y>-||x||^2-||y||^2+2<x,y>\\
&=&4<x,y>
\end{eqnarray*}

Now given a matrix $U\in M_N(\mathbb R)$, we have the following equivalences, with the first one coming from the above identities, and with the other ones being clear:
\begin{eqnarray*}
||Ux||=||x||
&\iff&<Ux,Uy>=<x,y>\\
&\iff&<x,U^tUy>=<x,y>\\
&\iff&U^tUy=y\\
&\iff&U^tU=1\\
&\iff&U^t=U^{-1}
\end{eqnarray*}

(2) The second assertion is clear from the definition of the isometries, and can be established as well by using matrices, and the $U^t=U^{-1}$ criterion.
\end{proof}

As a basic illustration here, we have:

\index{rotation}
\index{symmetry}
\index{orthogonal group}

\begin{theorem}
The rotations and symmetries in the plane, given by
$$R_t=\begin{pmatrix}\cos t&-\sin t\\ \sin t&\cos t\end{pmatrix}\quad,\quad 
S_t=\begin{pmatrix}\cos t&\sin t\\ \sin t&-\cos t\end{pmatrix}$$
are isometries. These are all the isometries in $2$ dimensions.
\end{theorem}

\begin{proof}
We already know that $R_t$ is the rotation of angle $t$. As for $S_t$, this is the symmetry with respect to the $Ox$ axis rotated by $t/2\in\mathbb R$. But this gives the result, since the isometries in 2 dimensions are obviously either rotations, or symmetries. 
\end{proof}

As a conclusion, the set $O_N$ from Theorem 1.39 is a quite fundamental object, with $O_2$ already consisting of some interesting $2\times2$ matrices, namely the matrices $R_t,S_t$. We will be back to $O_N$, which is a so-called group, and is actually one of the most important examples of groups, on several occasions, in what follows.

\section*{1e. Exercises}

The key thing in linear algebra is that of geometrically understanding the linear maps $x\to Ax$ associated to the matrices $A\in M_N(\mathbb R)$. Here is an exercise on this:

\begin{exercise}
Work out the geometric interpretation of the map $f(x)=Ax$, with
$$A\in M_2(\pm1)$$
and then discuss as well the diagonalization of these matrices.
\end{exercise}

To be more precise, there are $2^4=16$ matrices here, some of which were already discussed in the above. As a bonus exercise, you can try as well $A\in M_2(0,1)$, which is 16 more matrices. And for the black belt, try $A\in M_2(-1,0,1)$.

\begin{exercise}
Diagonalize explicitly the third flat matrix, namely
$$\mathbb I_3
=\begin{pmatrix}1&1&1\\1&1&1\\1&1&1\end{pmatrix}$$
and then study as well the general case, that of the matrix $\mathbb I_N$.
\end{exercise}

Here we already know from the above that the diagonal form is $D=(N,0,\ldots,0)$, and the problem is that of finding the passage matrix $P$, as to write the diagonalization formula $\mathbb I_N=PDP^{-1}$. The case to start with, as a warm-up for the exercise, is $N=2$, where $\mathbb I_2$ is twice the orthogonal projection on the $x=y$ diagonal, which was already discussed in the above. Then, go with $N=3$, and then with general $N\in\mathbb N$.

\begin{exercise}
Work out the trigonometry formulae
$$\sin(2t)=2\sin t\cos t\quad,\quad 
\cos(2t)=2\cos^2t-1$$ 
by using elementary methods, coming from plane geometry.
\end{exercise}

There are many ways of solving this exercise, and of course enjoy.

\begin{exercise}
Prove that the isometries in $2$ dimensions are either rotations, or symmetries, as to complete the proof of Theorem 1.40.
\end{exercise}

As before, there are many ways of dealing with this, all being nice geometry.

\begin{exercise}
Develop a theory of angles between the vectors $x,y\in\mathbb R^N$, by using the well-known formula
$$<x,y>=||x||\cdot||y||\cdot\cos t$$
that you should by the way fully understand first, in $N=2$ dimensions.
\end{exercise}

To be more precise, you must first make sure that the above formula holds indeed at $N=2$, as a theorem. Then, based on this, you can use this formula at $N\geq3$ too, but this time as a definition for the angle $t$ between $x,y$. There are many things that can be done here, and the more complete the theory that you develop, the better.

\chapter{The determinant}

\section*{2a. Matrix inversion}

We have seen in the previous chapter that most of the interesting maps $f:\mathbb R^N\to\mathbb R^N$ that we know, such as the rotations, symmetries and projections, are linear, and can be written in the following form, with $A\in M_N(\mathbb R)$ being a square matrix:
$$f(v)=Av$$

In this chapter we develop more general theory for such linear maps. We will be mostly motivated by the following fundamental result, which has countless concrete applications, and which is actually at the origin of the whole linear algebra theory:

\index{linear equations}
\index{matrix inversion}

\begin{theorem}
Any linear system of equations 
$$\begin{cases}
a_{11}x_1+a_{12}x_2+\ldots+a_{1N}x_N\!\!\!&=\ v_1\\
a_{21}x_1+a_{22}x_2+\ldots+a_{2N}x_N\!\!\!&=\ v_2\\
\ \ \vdots\\
a_{N1}x_1+a_{N2}x_2+\ldots+a_{NN}x_N\!\!\!&=\ v_N
\end{cases}$$
can be written in matrix form, as follows,
$$Ax=v$$
and when $A$ is invertible, its solution is given by $x=A^{-1}v$.
\end{theorem}

\begin{proof}
With linear algebra conventions, our system reads:
$$\begin{pmatrix}
a_{11}&a_{12}&\ldots&a_{1N}\\
a_{21}&a_{22}&\ldots&a_{2N}\\
\vdots&&&\vdots\\
a_{N1}&a_{N2}&\ldots&a_{NN}
\end{pmatrix}
\begin{pmatrix}
x_1\\
x_2\\
\vdots\\
x_N
\end{pmatrix}
=\begin{pmatrix}
v_1\\
v_2\\
\vdots\\
v_N
\end{pmatrix}$$

Thus, we are led to the conclusions in the statement.
\end{proof}

In practice, we are led to the question of inverting the matrices $A\in M_N(\mathbb R)$. And this is the same question as inverting the linear maps $f:\mathbb R^N\to\mathbb R^N$, due to:

\index{invertible matrix}

\begin{theorem}
A linear map $f:\mathbb R^N\to\mathbb R^N$, written as
$$f(v)=Av$$
is invertible precisely when $A$ is invertible, and in this case we have $f^{-1}(v)=A^{-1}v$.
\end{theorem}

\begin{proof}
This is something that we basically know, coming from the fact that, with the notation $f_A(v)=Av$, we have the following formula:
$$f_Af_B=f_{AB}$$

Thus, we are led to the conclusion in the statement.
\end{proof}

In order to study invertibility questions, for matrices and linear maps, let us begin with some examples. In the simplest case, in 2 dimensions, the result is as follows:

\index{inversion formula}

\begin{theorem}
We have the following inversion formula, for the $2\times2$ matrices:
$$\begin{pmatrix}a&b\\ c&d\end{pmatrix}^{-1}
=\frac{1}{ad-bc}\begin{pmatrix}d&-b\\ -c&a\end{pmatrix}$$
When $ad-bc=0$, the matrix is not invertible.
\end{theorem}

\begin{proof}
We have two assertions to be proved, the idea being as follows:

\medskip

(1) As a first observation, when $ad-bc=0$ we must have, for some $\lambda\in\mathbb R$:
$$b=\lambda a\quad,\quad 
d=\lambda c$$

Thus our matrix must be of the following special type:
$$\begin{pmatrix}a&b\\ c&d\end{pmatrix}=\begin{pmatrix}a&\lambda a\\ a&\lambda c\end{pmatrix}$$

But in this case the columns are proportional, so the linear map associated to the matrix is not invertible, and so the matrix itself is not invertible either.

\medskip

(2) When $ad-bc\neq 0$, let us look for an inversion formula of the following type:
$$\begin{pmatrix}a&b\\ c&d\end{pmatrix}^{-1}
=\frac{1}{ad-bc}\begin{pmatrix}*&*\\ *&*\end{pmatrix}$$

We must therefore solve the following equations:
$$\begin{pmatrix}a&b\\ c&d\end{pmatrix}
\begin{pmatrix}*&*\\ *&*\end{pmatrix}=
\begin{pmatrix}ad-bc&0\\ 0&ad-bc\end{pmatrix}$$

The obvious solution here is as follows:
$$\begin{pmatrix}a&b\\ c&d\end{pmatrix}
\begin{pmatrix}d&-b\\ -c&a\end{pmatrix}=
\begin{pmatrix}ad-bc&0\\ 0&ad-bc\end{pmatrix}$$

Thus, we are led to the formula in the statement.
\end{proof}

In order to deal now with the inversion problem in general, for the arbitrary matrices $A\in M_N(\mathbb R)$, we will use the same method as the one above, at $N=2$. Let us write indeed our matrix as follows, with $v_1,\ldots,v_N\in\mathbb R^N$ being its column vectors:
$$A=[v_1,\ldots,v_N]$$

We know from the general results from chapter 1 that, in order for $A$ to be invertible, the vectors $v_1,\ldots,v_N$ must be linearly independent. Thus, following the observations (1) from the above proof of Theorem 2.3, we are led into the question of understanding when a family of vectors $v_1,\ldots,v_N\in\mathbb R^N$ are linearly independent. 

\bigskip

In order to deal with this latter question, let us introduce the following notion:

\index{volume of parallelepiped}

\begin{definition}
Associated to any vectors $v_1,\ldots,v_N\in\mathbb R^N$ is the volume
$${\rm det}^+(v_1\ldots v_N)=vol<v_1,\ldots,v_N>$$
of the parallelepiped made by these vectors.
\end{definition}

Here the volume is taken in the standard $N$-dimensional sense. At $N=1$ this volume is a length, at $N=2$ this volume is an area, at $N=3$ this is the usual 3D volume, and so on. In general, the volume of a body $X\subset\mathbb R^N$ is by definition the number $vol(X)\in[0,\infty]$ of copies of the unit cube $C\subset\mathbb R^N$ which are needed for filling $X$, when allowing this unit cube to be divided into smaller cubes, for the needs of the filling operation.

\bigskip

In order to compute this volume we can use various geometric techniques, and we will see soon that, in what regards the case that we are interested in, namely that of the parallelepipeds $P\subset\mathbb R^N$, we can basically compute here everything, just by using very basic geometric techniques, essentially based on the Thales theorem.

\bigskip

In relation with our inversion problem, we have the following statement:

\index{invertible matrix}

\begin{theorem}
The quantity ${\rm det}^+$ that we constructed, regarded as a function of the corresponding square matrices, formed by column vectors,
$${\rm det}^+:M_N(\mathbb R)\to\mathbb R_+$$
has the property that a matrix $A\in M_N(\mathbb R)$ is invertible precisely when ${\rm det}^+(A)>0$.
\end{theorem}

\begin{proof}
This follows from Theorem 2.2, and from the general results from chapter 1, which tell us that a matrix $A\in M_N(\mathbb R)$ is invertible precisely when its column vectors $v_1,\ldots,v_N\in\mathbb R^N$ are linearly independent. But this latter condition is equivalent to the fact that we must have the following strict inequality: 
$$vol<v_1,\ldots,v_N>>0$$

Thus, we are led to the conclusion in the statement.
\end{proof}

Summarizing, all this leads us into the explicit computation of ${\rm det}^+$. As a first observation, in 1 dimension we obtain the absolute value of the real numbers:
$${\rm det}^+(a)=|a|$$

In 2 dimensions now, the computation is non-trivial, and we have the following result, making the link with our main result so far, namely Theorem 2.3:

\begin{theorem}
In $2$ dimensions we have the following formula,
$${\rm det}^+\begin{pmatrix}a&b\\ c&d\end{pmatrix}=|ad-bc|$$
with ${\rm det}^+:M_2(\mathbb R)\to\mathbb R_+$ being the function constructed above.
\end{theorem}

\begin{proof}
We must show that the area of the parallelogram formed by $\binom{a}{c},\binom{b}{d}$ equals $|ad-bc|$. We can assume $a,b,c,d>0$ for simplifying, the proof in general being similar. Moreover, by switching if needed the vectors $\binom{a}{c},\binom{b}{d}$, we can assume that we have:
$$\frac{a}{c}>\frac{b}{d}$$

According to these conventions, the picture of our parallelogram is as follows:
$$
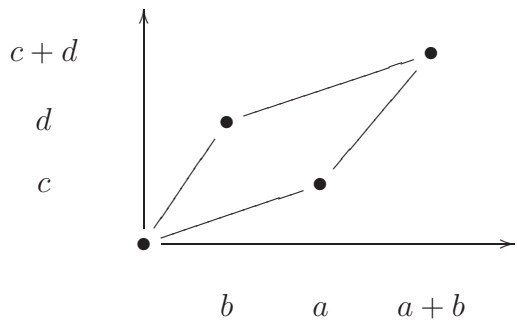
\xymatrix@R=10pt@C=15pt{
&&&\\
c+d&&&&\bullet\\
d&&\bullet\ar@{-}[urr]&&\\
c&&&\bullet\ar@{-}[uur]&\\
&\bullet\ar@{-}[urr]\ar[rrrr]\ar[uuuu]\ar@{-}[uur]&&&&\\
&&\ b\ &\ a\ &a+b}$$

Now let us slide the upper side downwards left, until we reach the $Oy$ axis. Our parallelogram, which has not changed its area in this process, becomes:
$$
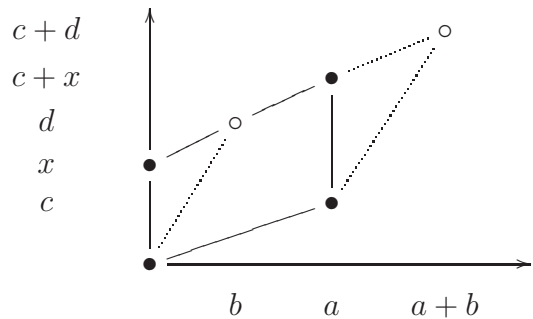
\xymatrix@R=2pt@C=16pt{
&&&\\
c+d&&&&\circ\\
c+x&&&\bullet\ar@{.}[ur]&\\
d&&\circ\ar@{-}[ur]&&&&&\\
x&\bullet\ar@{-}[ur]\ar[uuuu]&&&\\
c&&&\bullet\ar@{-}[uuu]\ar@{.}[uuuur]&\\
&&&&&&\\
&\bullet\ar@{-}[uurr]\ar@{-}[uuu]\ar[rrrr]\ar@{.}[uuuur]&&&&\\
&&\ b\ &\ a\ &a+b}$$

We can further modify this parallelogram, once again by not altering its area, by sliding the right side downwards, until we reach the $Ox$ axis:
$$
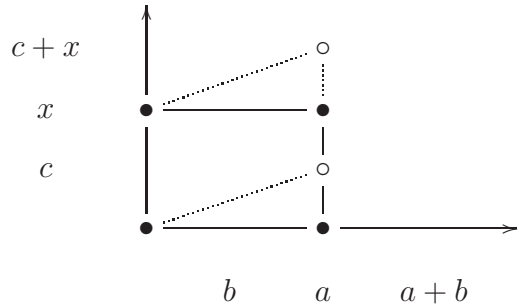
\xymatrix@R=10pt@C=15pt{
&&&\\
c+x&&&\circ&\\
x&\bullet\ar@{.}[urr]\ar[uu]\ar@{-}[rr]&&\bullet\ar@{.}[u]&\\
c&&&\circ\ar@{-}[u]&\\
&\bullet\ar@{.}[urr]\ar@{-}[uu]\ar@{-}[rr]&&\bullet\ar@{-}[u]\ar[rr]&&\\
&&\ b\ &\ a\ &a+b}$$

Let us compute now the area. Since our two sliding operations have not changed the area of the original parallelogram, this area is given by:
$$A=ax$$

In order to compute the quantity $x$, observe that in the context of the first move, we have two similar triangles, according to the following picture:
$$
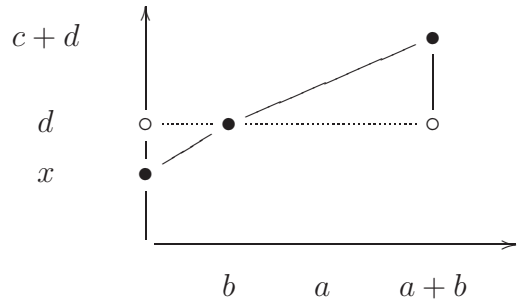
\xymatrix@R=5pt@C=15pt{
&&&&\\
c+d&&&&\bullet\\
&&&&&&\\
d&\circ\ar@{.}[r]\ar[uuu]&\bullet\ar@{.}[rr]\ar@{-}[uurr]&&\circ\ar@{-}[uu]\\
x&\bullet\ar@{-}[u]\ar@{-}[ur]&&&\\
&&&&\\
&\ar@{-}[uu]\ar[rrrr]&&&&\\
&&\ b\ &\ a\ &a+b}$$

Thus, we are led to the following equation for the number $x$:
$$\frac{d-x}{b}=\frac{c}{a}$$

By solving this equation, we obtain the following value for $x$:
$$x=d-\frac{bc}{a}$$

Thus the area of our parallelogram, or rather of the final rectangle obtained from it, which has the same area as the original parallelogram, is given by:
$$A=ax=ad-bc$$

Thus, we are led to the conclusion in the statement.
\end{proof}

\section*{2b. The determinant}

All the above is very nice, we obviously have a beginning of theory here. However, when looking carefully, we can see that our theory has a weakness, because:

\medskip

\begin{enumerate}
\item In 1 dimension the number $a$, which is the simplest function of $a$ itself, is certainly a better quantity than the number $|a|$.

\medskip

\item In 2 dimensions the number $ad-bc$, which is linear in $a,b,c,d$, is certainly a better quantity than the number $|ad-bc|$. 
\end{enumerate}

\medskip

So, let us upgrade now our theory, by constructing a better function, which does the same job, namely checking if the vectors are proportional, of the following type:
$$\det:M_N(\mathbb R)\to\mathbb R\quad,\quad 
\det=\pm {\rm det}^+$$

That is, we would like to have a clever, signed version of $\det^+$, satisfying:
$$\det(a)=a\quad,\quad 
\det\begin{pmatrix}a&b\\ c&d\end{pmatrix}=ad-bc$$

In order to do this, we must come up with a way of splitting the systems of vectors $v_1,\ldots,v_N\in\mathbb R^N$ into two classes, call them positive and negative. And here, the answer is quite clear, because a bit of thinking leads to the following definition:

\index{oriented system of vectors}
\index{unoriented system of vectors}
\index{sign of system of vectors}

\begin{definition}
A system of vectors $v_1,\ldots,v_N\in\mathbb R^N$ is called:
\begin{enumerate}
\item Oriented, if one can continuously pass from the standard basis to it.

\item Unoriented, otherwise.
\end{enumerate}
The associated sign is $+$ in the oriented case, and $-$ in the unoriented case. 
\end{definition}

As a first example, in 1 dimension the basis consists of the single vector $e=1$, which can be continuously deformed into any vector $a>0$. Thus, the sign is the usual one:
$$sgn(a)=
\begin{cases}
+&{\rm if}\ a>0\\
-&{\rm if}\ a<0
\end{cases}$$

Thus, in connection with our original question, we are definitely on the good track, because when multiplying $|a|$ by this sign we obtain $a$ itself, as desired:
$$a=sgn(a)|a|$$

In 2 dimensions now, the explicit formula of the sign is as follows:

\begin{proposition}
We have the following formula, valid for any $2$ vectors in $\mathbb R^2$,
$$sgn\left[\binom{a}{c},\binom{b}{d}\right]=sgn(ad-bc)$$
with the sign function on the right being the usual one, in $1$ dimension.
\end{proposition}

\begin{proof}
According to our conventions, the sign of $\binom{a}{c},\binom{b}{d}$ is as follows:

\medskip

(1) The sign is $+$ when these vectors come in this order with respect to the counterclockwise rotation in the plane, around 0.

\medskip

(2) The sign is $-$ otherwise, meaning when these vectors come in this order with respect to the clockwise rotation in the plane, around 0. 

\medskip

If we assume now $a,b,c,d>0$ for simplifying, we are left with comparing the angles having the numbers $c/a$ and $d/b$ as tangents, and we obtain in this way:
$$sgn\left[\binom{a}{c},\binom{b}{d}\right]=
\begin{cases}
+&{\rm if}\ \frac{c}{a}<\frac{d}{b}\\
-&{\rm if}\ \frac{c}{a}>\frac{d}{b}
\end{cases}$$

But this gives the formula in the statement. The proof in general is similar.
\end{proof}

Once again, in connection with our original question, we are on the good track, because when multiplying $|ad-bc|$ by this sign we obtain $ad-bc$ itself, as desired:
$$ad-bc=sgn(ad-bc)|ad-bc|$$

Let us look as well into the case $N=3$. Things here are more complicated, and we will discuss this later on. However, we have the following basic result:

\index{signature}

\begin{proposition}
Consider the standard basis of $\mathbb R^3$, namely:
$$e_1=\begin{pmatrix}1\\0\\0\end{pmatrix}
\qquad,\qquad 
e_2=\begin{pmatrix}0\\1\\0\end{pmatrix}
\qquad,\qquad 
e_3=\begin{pmatrix}0\\0\\1\end{pmatrix}$$
We have then the following sign computations:
\begin{enumerate}
\item $sgn(e_1,e_2,e_3)=+$.

\item $sgn(e_1,e_3,e_2)=-$.

\item $sgn(e_2,e_1,e_3)=-$.

\item $sgn(e_2,e_3,e_1)=+$.

\item $sgn(e_3,e_1,e_2)=+$.

\item $sgn(e_3,e_2,e_1)=-$.
\end{enumerate}
\end{proposition}

\begin{proof}
In each case the problem is whether one can continuously pass from $(e_1,e_2,e_3)$ to the basis in statement, and the computations can be done as follows:

\medskip

(1) In three of the cases under investigation, namely (2,3,6), one of the vectors is unchanged, and the other two are switched. Thus, we are more or less in 2 dimensions, and since the switch here clearly corresponds to $-$, the sign in these cases is $-$.

\medskip

(2) As for the remaining three cases, namely (1,4,5), here the sign can only be $+$, since things must be 50-50 between $+$ and $-$, say by symmetry reasons. And this is indeed the case, because what we have here are rotations of the standard basis.
\end{proof}

As already mentioned, we will be back to this later, with a general formula for the sign in 3 dimensions. This formula is quite complicated, the idea being that of making out of the $3\times3=9$ entries of our vectors a certain quantity, somewhat in the spirit of the one in Proposition 2.8, and then taking the sign of this quantity.

\bigskip

At the level of the general results now, we have:

\begin{proposition}
The orientation of a system of vectors changes as follows:
\begin{enumerate}
\item If we switch the sign of a vector, the associated sign switches.

\item If we permute two vectors, the associated sign switches as well.
\end{enumerate}
\end{proposition}

\begin{proof}
Both these assertions are clear from the definition of the sign, because the two operations in question change the orientation of the system of vectors.
\end{proof}

With the above notion in hand, we can now formulate:

\index{determinant}
\index{signed volume}

\begin{definition}
The determinant of $v_1,\ldots,v_N\in\mathbb R^N$ is the signed volume
$$\det(v_1\ldots v_N)=\pm vol<v_1,\ldots,v_N>$$
of the parallelepiped made by these vectors.
\end{definition}

In other words, we are upgrading here Definition 2.4, by adding a sign to the quantity ${\rm det}^+$ constructed there, as to potentially reach to good additivity properties:
$$\det(v_1\ldots v_N)=\pm {\rm det}^+(v_1\ldots v_N)$$

In relation with our original inversion problem for the square matrices, this upgrade does not change what we have so far, and we have the following statement:

\index{invertible matrix}

\begin{theorem}
The quantity $\det$ that we constructed, regarded as a function of the corresponding square matrices, formed by column vectors,
$$\det:M_N(\mathbb R)\to\mathbb R$$
has the property that a matrix $A\in M_N(\mathbb R)$ is invertible precisely when $\det(A)\neq 0$.
\end{theorem}

\begin{proof}
We know from Theorem 2.5 that a matrix $A\in M_N(\mathbb R)$ is invertible precisely when ${\rm det}^+(A)=|\det A|$ is strictly positive, and this gives the result.
\end{proof}

In the matrix context, we will often use the symbol $|\,.\,|$ instead of $\det$:
$$|A|=\det A$$

Let us try now to compute the determinant. In 1 dimension we have of course the formula $\det(a)=a$, because the absolute value fits, and so does the sign:
$$\det(a)
=sgn(a)\times|a|
=a$$

In 2 dimensions now, we have the following result:

\begin{theorem}
In $2$ dimensions we have the following formula,
$$\begin{vmatrix}a&b\\ c&d\end{vmatrix}=ad-bc$$
with $|\,.\,|=\det$ being the determinant function constructed above.
\end{theorem}

\begin{proof}
According to our definition, to the computation in Theorem 2.6, and to sign formula from Proposition 2.8, the determinant of a $2\times2$ matrix is given by:
\begin{eqnarray*}
\det\begin{pmatrix}a&b\\ c&d\end{pmatrix}
&=&sgn\left[\binom{a}{c},\binom{b}{d}\right]\times {\rm det}^+\begin{pmatrix}a&b\\ c&d\end{pmatrix}\\
&=&sgn\left[\binom{a}{c},\binom{b}{d}\right]\times|ad-bc|\\
&=&sgn(ad-bc)\times|ad-bc|\\
&=&ad-bc
\end{eqnarray*}

Thus, we have obtained the formula in the statement.
\end{proof}

\section*{2c. Basic properties}

In order to discuss now arbitrary dimensions, we will need a number of theoretical results. Here is a first series of formulae, coming straight from definitions:

\begin{theorem}
The determinant has the following properties:
\begin{enumerate}
\item When multiplying by scalars, the determinant gets multiplied as well:
$$\det(\lambda_1v_1,\ldots,\lambda_Nv_N)=\lambda_1\ldots\lambda_N\det(v_1,\ldots,v_N)$$

\item When permuting two columns, the determinant changes the sign:
$$\det(\ldots,u,\ldots,v,\ldots)=-\det(\ldots,v,\ldots,u,\ldots)$$

\item The determinant $\det(e_1,\ldots,e_N)$ of the standard basis of $\mathbb R^N$ is $1$.
\end{enumerate}
\end{theorem}

\begin{proof}
All this is clear from definitions, as follows:

\medskip

(1) This follows from definitions, and from Proposition 2.10 (1).

\medskip

(2) This follows as well from definitions, and from Proposition 2.10 (2).

\medskip

(3) This is clear from our definition of the determinant.
\end{proof}

As an application of the above result, we have:

\begin{theorem}
The determinant of a diagonal matrix is given by:
$$\begin{vmatrix}
\lambda_1\\ 
&\ddots\\
&&\lambda_N\end{vmatrix}=\lambda_1\ldots\lambda_N$$
That is, we obtain the product of diagonal entries, or of eigenvalues.
\end{theorem}

\begin{proof}
The formula in the statement is clear by using the rules (1) and (3) in Theorem 2.14, which in matrix terms give:
\begin{eqnarray*}
\begin{vmatrix}
\lambda_1\\ 
&\ddots\\
&&\lambda_N\end{vmatrix}
&=&\lambda_1\ldots\lambda_N
\begin{vmatrix}
1\\ 
&\ddots\\
&&1\end{vmatrix}\\
&=&\lambda_1\ldots\lambda_N
\end{eqnarray*}

As for the last assertion, this is rather a remark.
\end{proof}

The above result is very useful, and we will see in a moment that, more generally, the determinant of any diagonalizable matrix is the product of its eigenvalues.

\bigskip

In order to reach now to a more advanced theory, let us adopt the linear map point of view. In this setting, the definition of the determinant reformulates as follows: 

\begin{theorem}
Given a linear map, written as $f(v)=Av$, its ``inflation coefficient'', obtained as the signed volume of the image of the unit cube, is given by:
$$I_f=\det A$$
More generally, $I_f$ is the inflation ratio of any parallelepiped in $\mathbb R^N$, via the transformation $f$. In particular $f$ is invertible precisely when $\det A\neq0$.
\end{theorem}

\begin{proof}
The only non-trivial thing in all this is the fact that the inflation coefficient $I_f$, as defined above, is independent of the choice of the parallelepiped. But this is a generalization of the Thales theorem, which follows from the Thales theorem itself.
\end{proof}

As a first application of the above linear map viewpoint, we have:

\index{determinant of products}

\begin{theorem}
We have the following formula, valid for any matrices $A,B$:
$$\det(AB)=\det A\cdot\det B$$
In particular, we have $\det(AB)=\det(BA)$.
\end{theorem}

\begin{proof}
The decomposition formula in the statement follows by using the associated linear maps, which multiply as follows:
$$f_{AB}=f_Af_B$$

Indeed, when computing the determinant, by using the ``inflation coefficient'' viewpoint from Theorem 2.16, we obtain the same thing on both sides. As for the formula $\det(AB)=\det(BA)$, this is clear from the first formula, which is symmetric in $A,B$.
\end{proof}

Getting back now to explicit computations, we have the following key result:

\index{product of eigenvalues}

\begin{theorem}
The determinant of a diagonalizable matrix 
$$A\sim\begin{pmatrix}
\lambda_1\\ 
&\ddots\\
&&\lambda_N\end{pmatrix}$$
is the product of its eigenvalues, $\det A=\lambda_1\ldots\lambda_N$.
\end{theorem}

\begin{proof}
We know that a diagonalizable matrix can be written in the form $A=PDP^{-1}$, with $D=diag(\lambda_1,\ldots,\lambda_N)$. Now by using Theorem 2.17, we obtain:
\begin{eqnarray*}
\det A
&=&\det(PDP^{-1})\\
&=&\det(DP^{-1}P)\\
&=&\det D\\
&=&\lambda_1\ldots\lambda_N
\end{eqnarray*}

Thus, we are led to the formula in the statement.
\end{proof}

Here is another important result, which is very useful for diagonalization:

\index{characteristic polynomial}
\index{eigenvalue}
\index{roots of polynomial}

\begin{theorem}
The eigenvalues of a matrix $A\in M_N(\mathbb R)$ are the roots of
$$P(x)=\det(A-x1_N)$$
called characteristic polynomial of the matrix.
\end{theorem}

\begin{proof}
We have the following computation, using the fact that a linear map is bijective precisely when the determinant of the associated matrix is nonzero:
\begin{eqnarray*}
\exists v,Av=\lambda v
&\iff&\exists v,(A-\lambda 1_N)v=0\\
&\iff&\det(A-\lambda 1_N)=0
\end{eqnarray*}

Thus, we are led to the conclusion in the statement.
\end{proof}

Here are now some other computations, once again in arbitrary dimensions:

\begin{proposition}
We have the following results:
\begin{enumerate}
\item The determinant of an orthogonal matrix must be $\pm1$.

\item The determinant of a projection must be $0$ or $1$.
\end{enumerate}
\end{proposition}

\begin{proof}
These are elementary results, the idea being as follows:

\medskip

(1) Here the determinant must be indeed $\pm1$, because the orthogonal matrices map the unit cube to a copy of the unit cube.

\medskip

(2) Here the determinant is 0, because the projections flatten the unit cube, unless the projection in question is the identity, where the determinant is 1.
\end{proof}

In general now, at the theoretical level, we have the following key result:

\begin{theorem}
The determinant has the additivity property
$$\det(\ldots,u+v,\ldots)
=\det(\ldots,u,\ldots)
+\det(\ldots,v,\ldots)$$
valid for any choice of the vectors involved.
\end{theorem}

\begin{proof}
This follows by doing some elementary geometry, in the spirit of the computations in the proof of Theorem 2.6, as follows:

\medskip

(1) We can either use the Thales theorem, and then compute the volumes of all the parallelepipeds involved, by using basic algebraic formulae.

\medskip

(2) Or we can solve the problem in ``puzzle'' style, the idea being to cut the big parallelepiped, and then recover the small ones, after some manipulations.

\medskip

(3) We can do as well something hybrid, consisting in deforming the parallelepipeds involved, without changing their volumes, and then cutting and gluing. 
\end{proof}

As a basic application of the above result, we have:

\index{upper triangular matrix}
\index{lower triangular matrix}

\begin{theorem}
We have the following results:
\begin{enumerate}
\item The determinant of a diagonal matrix is the product of diagonal entries.

\item The same is true for the upper triangular matrices.

\item The same is true for the lower triangular matrices.
\end{enumerate}
\end{theorem}

\begin{proof}
All this can be deduced by using our various general formulae, as follows:

\medskip

(1) This is something that we already know, from Theorem 2.15.

\medskip

(2) This follows by using Theorem 2.14 and Theorem 2.21, then (1), as follows:
\begin{eqnarray*}
\begin{vmatrix}
\lambda_1&&&*\\ 
&\lambda_2\\
&&\ddots\\
0&&&\lambda_N\end{vmatrix}
&=&\begin{vmatrix}
\lambda_1&0&&*\\ 
&\lambda_2\\
&&\ddots\\
0&&&\lambda_N\end{vmatrix}\\
&&\vdots\\
&&\vdots\\
&=&\begin{vmatrix}
\lambda_1&&&0\\ 
&\lambda_2\\
&&\ddots\\
0&&&\lambda_N\end{vmatrix}\\
&=&\lambda_1\ldots\lambda_N
\end{eqnarray*}

(3) This follows as well from Theorem 2.14 and Theorem 2.21, then (1), by proceeding this time from right to left, from the last column towards the first column.
\end{proof}

We can see from the above that the rules in Theorem 2.14 and Theorem 2.21 are quite powerful, taken altogether. For future reference, let us record these rules:

\begin{theorem}
The determinant has the following properties:
\begin{enumerate}
\item When adding two columns, the determinants get added:
$$\det(\ldots,u+v,\ldots)
=\det(\ldots,u,\ldots)
+\det(\ldots,v,\ldots)$$

\item When multiplying columns by scalars, the determinant gets multiplied:
$$\det(\lambda_1v_1,\ldots,\lambda_Nv_N)=\lambda_1\ldots\lambda_N\det(v_1,\ldots,v_N)$$

\item When permuting two columns, the determinant changes the sign:
$$\det(\ldots,u,\ldots,v,\ldots)=-\det(\ldots,v,\ldots,u,\ldots)$$

\item The determinant $\det(e_1,\ldots,e_N)$ of the standard basis of $\mathbb R^N$ is $1$.
\end{enumerate}
\end{theorem}

\begin{proof}
This is something that we already know, which follows by putting together the various formulae from Theorem 2.14 and Theorem 2.21.
\end{proof}

As an important theoretical result now, which will ultimately lead to an algebraic reformulation of the whole determinant problematics, we have:

\begin{theorem}
The determinant of square matrices is the unique map
$$\det:M_N(\mathbb R)\to\mathbb R$$
satisfying the conditions in Theorem 2.23.
\end{theorem}

\begin{proof}
This can be done in two steps, as follows:

\medskip

(1) Our first claim is that any map $\det':M_N(\mathbb R)\to\mathbb R$ satisfying the conditions in Theorem 2.23 must coincide with $\det$ on the upper triangular matrices. But this is clear from the proof of Theorem 2.22, which only uses the rules in Theorem 2.23.

\medskip

(2) Our second claim is that we have $\det'=\det$, on all matrices. But this can be proved by putting the matrix in upper triangular form, by using operations on the columns, in the spirit of the manipulations from the proof of Theorem 2.22.
\end{proof}

Here is now another important theoretical result:

\index{row expansion}

\begin{theorem}
The determinant is subject to the row expansion formula
\begin{eqnarray*}
\begin{vmatrix}a_{11}&\ldots&a_{1N}\\
\vdots&&\vdots\\
a_{N1}&\ldots&a_{NN}\end{vmatrix}
&=&a_{11}\begin{vmatrix}a_{22}&\ldots&a_{2N}\\
\vdots&&\vdots\\
a_{N2}&\ldots&a_{NN}\end{vmatrix}
-a_{12}\begin{vmatrix}a_{21}&a_{23}&\ldots&a_{2N}\\
\vdots&\vdots&&\vdots\\
a_{N1}&a_{N3}&\ldots&a_{NN}\end{vmatrix}\\
&&+\ldots\ldots
+(-1)^{N+1}a_{1N}\begin{vmatrix}a_{21}&\ldots&a_{2,N-1}\\
\vdots&&\vdots\\
a_{N1}&\ldots&a_{N,N-1}\end{vmatrix}
\end{eqnarray*}
and this method fully computes it, by recurrence.
\end{theorem}

\begin{proof}
This follows from the fact that the formula in the statement produces a certain function $\det:M_N(\mathbb R)\to\mathbb R$, which has the 4 properties in Theorem 2.23.
\end{proof}

We can expand as well over the columns, as follows:

\index{column expansion}

\begin{theorem}
The determinant is subject to the column expansion formula
\begin{eqnarray*}
\begin{vmatrix}a_{11}&\ldots&a_{1N}\\
\vdots&&\vdots\\
a_{N1}&\ldots&a_{NN}\end{vmatrix}
&=&a_{11}\begin{vmatrix}a_{22}&\ldots&a_{2N}\\
\vdots&&\vdots\\
a_{N2}&\ldots&a_{NN}\end{vmatrix}
-a_{21}\begin{vmatrix}a_{12}&\ldots&a_{1N}\\
a_{32}&\ldots&a_{3N}\\
\vdots&&\vdots\\
a_{N2}&\ldots&a_{NN}\end{vmatrix}\\
&&+\ldots\ldots
+(-1)^{N+1}a_{N1}\begin{vmatrix}a_{12}&\ldots&a_{1N}\\
\vdots&&\vdots\\
a_{N-1,2}&\ldots&a_{N-1,N}\end{vmatrix}
\end{eqnarray*}
and this method fully computes it, by recurrence.
\end{theorem}

\begin{proof}
This follows by using the same argument as for the rows.
\end{proof}

We can now complement Theorem 2.23 with a similar result for the rows:

\begin{theorem}
The determinant has the following properties:
\begin{enumerate}
\item When adding two rows, the determinants get added:
$$\det\begin{pmatrix}\vdots\\ u+v\\ \vdots\end{pmatrix}
=\det\begin{pmatrix}\vdots\\ u\\ \vdots\end{pmatrix}
+\det\begin{pmatrix}\vdots \\ v\\ \vdots\end{pmatrix}$$

\item When multiplying row by scalars, the determinant gets multiplied:
$$\det\begin{pmatrix}\lambda_1v_1\\ \vdots\\ \lambda_Nv_N\end{pmatrix}
=\lambda_1\ldots\lambda_N\det\begin{pmatrix}v_1\\ \vdots\\ v_N\end{pmatrix}$$

\item When permuting two rows, the determinant changes the sign.
\end{enumerate}
\end{theorem}

\begin{proof}
This follows indeed by using the using various formulae established above, and is best seen by using the column expansion formula from Theorem 2.26.
\end{proof}

We can see from the above that the determinant is the subject to many interesting formulae, and that some of these formulae, when taken altogether, uniquely determine it. In all this, what is the most luminous is certainly the definition of the determinant as a volume. As for the second most luminous of our statements, this is Theorem 2.24, which is something a bit abstract, but both beautiful and useful. So, as a final theoretical statement now, here is an alternative reformulation of Theorem 2.24:

\index{multilinear form}

\begin{theorem}
The determinant of the systems of vectors
$$\det:\mathbb R^N\times\ldots\times\mathbb R^N\to\mathbb R$$
is multilinear, alternate and unital, and unique with these properties.
\end{theorem}

\begin{proof}
This is a fancy reformulation of Theorem 2.24, with the various properties of $\det$ from the statement being those from Theorem 2.23.
\end{proof}

As a conclusion to all this, we have now a full theory for the determinant, and we can freely use all the above results, definitions and theorems alike, and even start forgetting what is actually definition, and what is theorem.

\section*{2d. Sarrus and beyond}

As a first application of the above methods, we can now prove:

\index{Sarrus formula}

\begin{theorem}
The determinant of the $3\times3$ matrices is given by
$$\begin{vmatrix}a&b&c\\ d&e&f\\ g&h&i\end{vmatrix}=aei+bfg+cdh-ceg-bdi-afh$$
which can be memorized by using Sarrus' triangle method,
\begin{eqnarray*}
\det
&=&\begin{pmatrix}*&&\\ &*&\\ &&*\end{pmatrix}
+\begin{pmatrix}&*&\\ &&*\\ *&&\end{pmatrix}
+\begin{pmatrix}&&*\\ *&&\\ &*&\end{pmatrix}\\
&-&\begin{pmatrix}&&*\\ &*&\\ *&&\end{pmatrix}
+\begin{pmatrix}&*&\\ *&&\\ &&*\end{pmatrix}
+\begin{pmatrix}*&&\\ &&*\\ &*&\end{pmatrix}
\end{eqnarray*}
``triangles parallel to the diagonal, minus triangles parallel to the antidiagonal".
\end{theorem}

\begin{proof}
Here is the computation, using Theorem 2.25:
\begin{eqnarray*}
\begin{vmatrix}a&b&c\\ d&e&f\\ g&h&i\end{vmatrix}
&=&a\begin{vmatrix}e&f\\h&i\end{vmatrix}
-b\begin{vmatrix}d&f\\g&i\end{vmatrix}
+c\begin{vmatrix}d&e\\g&h\end{vmatrix}\\
&=&a(ei-fh)-b(di-fg)+c(dh-eg)\\
&=&aei-afh-bdi+bfg+cdh-ceg\\
&=&aei+bfg+cdh-ceg-bdi-afh
\end{eqnarray*}

Thus, we obtain the formula in the statement.
\end{proof}

As a first application, let us go back to the inversion problem for the $3\times3$ matrices, that we left open in the above. We can now solve this problem, as follows:

\index{matrix inversion}

\begin{theorem}
The inverses of the $3\times3$ matrices are given by
$$\begin{pmatrix}a&b&c\\ d&e&f\\ g&h&i\end{pmatrix}^{-1}
=\frac{1}{D}\begin{pmatrix}ei-fh&ch-bi&bf-ce\\ fg-di&ai-cg&cd-af\\ dh-eg&bg-ah&ae-bd\end{pmatrix}$$
with $D$ being the determimant. When $D=0$, the matrix is not invertible.
\end{theorem}

\begin{proof}
We can use here the same method as for the $2\times 2$ matrices. To be more precise, in order for the matrix to be invertible, we must have: 
$$D\neq0$$

The trick now is to look for solutions of the following problem:
$$\begin{pmatrix}a&b&c\\ d&e&f\\ g&h&i\end{pmatrix}
\begin{pmatrix}*&*&*\\ *&*&*\\ *&*&*\end{pmatrix}
=\begin{pmatrix}D&0&0\\ 0&D&0\\ 0&0&D\end{pmatrix}$$

We know from Theorem 2.29 that the determinant is given by:
$$D=aei+bfg+cdh-ceg-bdi-afh$$

But this leads, via some obvious choices, to the following solution:
$$\begin{pmatrix}*&*&*\\ *&*&*\\ *&*&*\end{pmatrix}
=\begin{pmatrix}ei-fh&ch-bi&bf-ce\\ fg-di&ai-cg&cd-af\\ dh-eg&bg-ah&ae-bd\end{pmatrix}$$

Thus, by rescaling, we obtain the formula in the statement.
\end{proof}

In fact, we can now fully solve the inversion problem, as follows:

\index{matrix inversion}
\index{checkered signs}

\begin{theorem}
The inverse of a square matrix, having nonzero determinant,
$$A=\begin{pmatrix}a_{11}&\ldots&a_{1N}\\
\vdots&&\vdots\\
a_{N1}&\ldots&a_{NN}\end{pmatrix}$$
is given by the following formula,
$$A^{-1}=\frac{1}{\det A}
\begin{pmatrix}
\det A^{(11)}&-\det A^{(21)}&\det A^{(31)}&\ldots\\
-\det A^{(12)}&\det A^{(22)}&-\det A^{(32)}&\ldots\\
\det A^{(13)}&-\det A^{(23)}&\det A^{(33)}&\ldots\\
\vdots&\vdots&\vdots&
\end{pmatrix}$$
where $A^{(ij)}$ is the matrix $A$, with the $i$-th row and $j$-th column removed.
\end{theorem}

\begin{proof}
This follows indeed by using the row expansion formula from Theorem 2.25, which in terms of the matrix $A^{-1}$ in the statement reads $AA^{-1}=1$.
\end{proof}

In practice, the above result leads to the following algorithm, which is quite easy to memorize, for computing the inverse:

\medskip

(1) Delete rows and columns, and compute the corresponding determinants.

\medskip

(2) Transpose, and add checkered signs.

\medskip

(3) Divide by the determinant.

\medskip

Observe that this generalizes our previous computations at $N=2,3$. As an illustration, consider an arbitrary $2\times2$ matrix, written as follows:
$$A=\begin{pmatrix}a&b\\ c&d\end{pmatrix}$$

By deleting rows and columns we obtain $1\times1$ matrices, and so the matrix formed by the determinants $\det(A^{(ij)})$ is as follows:
$$M=\begin{pmatrix}d&c\\ b&a\end{pmatrix}$$

Now by transposing, adding checkered signs and dividing by $\det A$, we obtain:
$$A^{-1}=\frac{1}{ad-bc}\begin{pmatrix}d&-b\\ -c&a\end{pmatrix}$$

Similarly, at $N=3$ what we obtain is the inversion formula from Theorem 2.30.

\bigskip

As a new application now, let us record the following result, at $N=4$:

\begin{theorem}
The determinant of the $4\times4$ matrices is given by
\begin{eqnarray*}
&&\begin{vmatrix}a_1&a_2&a_3&a_4\\ b_1&b_2&b_3&b_4\\ c_1&c_2&c_3&c_4\\ d_1&d_2&d_3&d_4\end{vmatrix}\\
&=&a_1b_2c_3d_4-a_1b_2c_4d_3-a_1b_3c_2d_4+a_1b_3c_4d_2+a_1b_4c_2d_3-a_1b_4c_3d_2\\
&-&a_2b_1c_3d_4+a_2b_1c_4d_3+a_2b_3c_1d_4-a_2b_3c_4d_1-a_2b_4c_1d_3+a_2b_4c_3d_1\\
&+&a_3b_1c_2d_4+a_3b_1c_4d_2-a_3b_2c_1d_4+a_3b_2c_4d_1+a_3b_4c_1d_2-a_3b_4c_2d_1\\
&-&a_4b_1c_2d_3+a_4b_1c_3d_2-a_4b_2c_1d_3-a_4b_2c_3d_1-a_4b_3c_1d_2+a_4b_3c_2d_1
\end{eqnarray*}
and the formula of the inverse is as follows, involving $16$ Sarrus determinants,
$$A^{-1}=\frac{1}{\det A}
\begin{pmatrix}
\det A^{(11)}&-\det A^{(21)}&\det A^{(31)}&-\det A^{(41)}\\
-\det A^{(12)}&\det A^{(22)}&-\det A^{(32)}&\det A^{(42)}\\
\det A^{(13)}&-\det A^{(23)}&\det A^{(33)}&-\det A^{(43)}\\
-\det A^{(14)}&\det A^{(24)}&-\det A^{(34)}&\det A^{(44)}
\end{pmatrix}$$
where $A^{(ij)}$ is the matrix $A$, with the $i$-th row and $j$-th column removed.
\end{theorem}

\begin{proof}
The formula for the determinant follows by developing over the first row, then by using the Sarrus formula, for each of the 4 smaller determinants which appear:
\begin{eqnarray*}
\begin{vmatrix}a_1&a_2&a_3&a_4\\ b_1&b_2&b_3&b_4\\ c_1&c_2&c_3&c_4\\ d_1&d_2&d_3&d_4\end{vmatrix}
&=&a_1\begin{vmatrix}b_2&b_3&b_4\\ c_2&c_3&c_4\\ d_2&d_3&d_4\end{vmatrix}
-a_2\begin{vmatrix}b_1&b_3&b_4\\ c_1&c_3&c_4\\ d_1&d_3&d_4\end{vmatrix}\\
&+&a_3\begin{vmatrix}b_1&b_2&b_4\\ c_1&c_2&c_4\\ d_1&d_2&d_4\end{vmatrix}
-a_4\begin{vmatrix}b_1&b_2&b_3\\ c_1&c_2&c_3\\ d_1&d_2&d_3\end{vmatrix}
\end{eqnarray*}

As for the formula of the inverse, this is something that we already know.
\end{proof}

Let us discuss now the general formula of the determinant, at arbitrary values $N\in\mathbb N$ of the matrix size, generalizing those that we have at $N=2,3,4$. We will need:

\index{permutation}

\begin{definition}
A permutation of $\{1,\ldots,N\}$ is a bijection, as follows:
$$\sigma:\{1,\ldots,N\}\to\{1,\ldots,N\}$$
The set of such permutations is denoted $S_N$.
\end{definition}

There are many possible notations for the permutations, the basic one consisting in writing the numbers $1,\ldots,N$, and below them, their permuted versions:
$$\sigma=\begin{pmatrix}
1&2&3&4&5\\
2&1&4&5&3
\end{pmatrix}$$

Another method, which is faster, is by using diagrams, acting from top to bottom:
$$\xymatrix@R=3mm@C=3.5mm{
&\ar@{-}[ddr]&\ar@{-}[ddl]&\ar@{-}[ddrr]&\ar@{-}[ddl]&\ar@{-}[ddl]\\
\sigma=\\
&&&&&}$$

Here are some basic properties of the permutations:

\begin{theorem}
The permutations have the following properties:
\begin{enumerate}
\item There are $N!$ of them.

\item They are stable by composition, and inversion.
\end{enumerate}
\end{theorem}

\begin{proof}
In order to construct a permutation $\sigma\in S_N$, we have:

\smallskip

-- $N$ choices for the value of $\sigma(N)$.

-- $(N-1)$ choices for the value of $\sigma(N-1)$.

-- $(N-2)$ choices for the value of $\sigma(N-2)$.

$\vdots$

-- and so on, up to 1 choice for the value of $\sigma(1)$.

\smallskip

Thus, we have $N!$ choices, as claimed. As for the second assertion, this is clear.
\end{proof}

We will need the following key result:

\index{signature}
\index{number of inversions}
\index{transpositions}
\index{odd cycles}
\index{crossings}

\begin{theorem}
The permutations have a signature function
$$\varepsilon:S_N\to\{\pm1\}$$
which can be defined in the following equivalent ways:
\begin{enumerate}
\item As $(-1)^c$, where $c$ is the number of inversions.

\item As $(-1)^t$, where $t$ is the number of transpositions.

\item As $(-1)^o$, where $o$ is the number of odd cycles.

\item As $(-1)^x$, where $x$ is the number of crossings.

\item As the sign of the corresponding permuted basis of $\mathbb R^N$.
\end{enumerate}
\end{theorem}

\begin{proof}
This is something important, and quite subtle, to be systematically used in what follows. As a first observation, we can see right away a relation with the determinant, coming from (5). Thus, we already have some knowledge here, for instance coming from Proposition 2.9, which computes the signature of the permutations $\sigma\in S_3$.

\medskip

In practice now, we have explain what the numbers $c,t,o,x$ appearing in (1-4) above exactly are, then why they are well-defined modulo 2, then why they are equal to each other, and finally why the constructions (1-4) yield the same sign as (5).

\medskip

Let us begin with the first two steps, namely precise definition of $c,t,o,x$, and fact that these numbers are well-defined modulo 2:

\medskip

(1) The idea here is that given any two numbers $i<j$ among $1,\ldots,N$, the permutation  can either keep them in the same order, $\sigma(i)<\sigma(j)$, or invert them:
$$\sigma(j)>\sigma(i)$$

Now by making $i<j$ vary over all pairs of numbers in $1,\ldots,N$, we can count the number of inversions, and call it $c$. This is an integer, $c\in\mathbb N$, which is well-defined.

\medskip

(2) Here the idea, which is something quite intuitive, is that any permutation appears as a product of switches, also called transpositions: 
$$i\leftrightarrow j$$

The decomposition as a product of transpositions is not unique, but the number $t$ of the needed transpositions is unique, when considered modulo 2. This follows for instance from the equivalence of (2) with (1,3,4,5), explained below.

\medskip

(3) Here the point is that any permutation decomposes, in a unique way, as a product of cycles, which are by definition permutations of the following type:
$$i_1\to i_2\to i_3\to\ldots\ldots\to i_k\to i_1$$

Some of these cycles have even length, and some others have odd length. By counting those having odd length, we obtain a well-defined number $o\in\mathbb N$.

\medskip

(4) Here the method is that of drawing the permutation, as we usually do, and by avoiding triple crossings, and then counting the number of crossings. This number $x$ depends on the way we draw the permutations, but modulo 2, we always get the same number. Indeed, this follows from the fact that we can continuously pass from a drawing to each other, and that when doing so, the number of crossings can only jump by $\pm2$.

\medskip

Summarizing, we have 4 different definitions for the signature of the permutations, which all make sense, constructed according to (1-4) above. Regarding now the fact that we always obtain the same number, this can be established as follows:

\medskip

(1)=(2) This is clear, because any transposition inverts once, modulo 2.

\medskip

(1)=(3) This is clear as well, because the odd cycles invert once, modulo 2.

\medskip

(1)=(4) This comes from the fact that the crossings correspond to inversions.

\medskip

(2)=(3) This follows by decomposing the cycles into transpositions.

\medskip

(2)=(4) This comes from the fact that the crossings correspond to transpositions.

\medskip

(3)=(4) This follows by drawing a product of cycles, and counting the crossings.

\medskip

Finally, in what regards the equivalence of all these constructions with (5), here simplest is to use (2). Indeed, we already know that the sign of a system of vectors switches when interchanging two vectors, and so the equivalence between (2,5) is clear. 
\end{proof}

We can now formulate a key result, as follows:

\index{determinant formula}
\index{Sarrus formula}

\begin{theorem}
We have the following formula for the determinant,
$$\det A=\sum_{\sigma\in S_N}\varepsilon(\sigma)A_{1\sigma(1)}\ldots A_{N\sigma(N)}$$
with the signature function being the one introduced above.
\end{theorem}

\begin{proof}
This follows by recurrence over $N\in\mathbb N$, as follows:

\medskip

(1) When developing the determinant over the first column, we obtain a signed sum of $N$ determinants of size $(N-1)\times(N-1)$. But each of these determinants can be computed by developing over the first column too, and so on, and we are led to the conclusion that we have a formula as in the statement, with $\varepsilon(\sigma)\in\{-1,1\}$ being certain coefficients.

\medskip

(2) But these latter coefficients $\varepsilon(\sigma)\in\{-1,1\}$ can only be the signatures of the corresponding permutations $\sigma\in S_N$, with this being something that can be viewed again by recurrence, with either of the definitions (1-5) in Theorem 2.35 for the signature.
\end{proof}

The above result is something quite tricky, and in order to get familiar with it, there is nothing better than doing some computations. As a first, basic example, in 2 dimensions we recover the usual formula of the determinant, the details being as follows:
\begin{eqnarray*}
\begin{vmatrix}a&b\\ c&d\end{vmatrix}
&=&\varepsilon(|\,|)\cdot ad+\varepsilon(\slash\hskip-2mm\backslash)\cdot cb\\
&=&1\cdot ad+(-1)\cdot cb\\
&=&ad-bc
\end{eqnarray*}

In 3 dimensions now, we recover the Sarrus formula:
$$\begin{vmatrix}a&b&c\\ d&e&f\\ g&h&i\end{vmatrix}=aei+bfg+cdh-ceg-bdi-afh$$

Observe that the triangles in the Sarrus formula correspond to the permutations of $\{1,2,3\}$, and their signs correspond to the signatures of these permutations:
\begin{eqnarray*}
\det
&=&\begin{pmatrix}*&&\\ &*&\\ &&*\end{pmatrix}
+\begin{pmatrix}&*&\\ &&*\\ *&&\end{pmatrix}
+\begin{pmatrix}&&*\\ *&&\\ &*&\end{pmatrix}\\
&-&\begin{pmatrix}&&*\\ &*&\\ *&&\end{pmatrix}
+\begin{pmatrix}&*&\\ *&&\\ &&*\end{pmatrix}
+\begin{pmatrix}*&&\\ &&*\\ &*&\end{pmatrix}
\end{eqnarray*}

Also, in 4 dimensions, we recover the formula that we already know, as follows:

\begin{theorem}
The determinant of the $4\times4$ matrices is given by
\begin{eqnarray*}
&&\begin{vmatrix}a_1&a_2&a_3&a_4\\ b_1&b_2&b_3&b_4\\ c_1&c_2&c_3&c_4\\ d_1&d_2&d_3&d_4\end{vmatrix}\\
&=&a_1b_2c_3d_4-a_1b_2c_4d_3-a_1b_3c_2d_4+a_1b_3c_4d_2+a_1b_4c_2d_3-a_1b_4c_3d_2\\
&-&a_2b_1c_3d_4+a_2b_1c_4d_3+a_2b_3c_1d_4-a_2b_3c_4d_1-a_2b_4c_1d_3+a_2b_4c_3d_1\\
&+&a_3b_1c_2d_4+a_3b_1c_4d_2-a_3b_2c_1d_4+a_3b_2c_4d_1+a_3b_4c_1d_2-a_3b_4c_2d_1\\
&-&a_4b_1c_2d_3+a_4b_1c_3d_2-a_4b_2c_1d_3-a_4b_2c_3d_1-a_4b_3c_1d_2+a_4b_3c_2d_1
\end{eqnarray*}
with the generic term being of the following form, with $\sigma\in S_4$,
$$\pm a_{\sigma(1)}b_{\sigma(2)}c_{\sigma(3)}d_{\sigma(4)}$$
and with the sign being $\varepsilon(\sigma)$, computable by using Theorem 2.35.
\end{theorem}

\begin{proof}
We can indeed recover this formula as well as a particular case of Theorem 2.36. To be more precise, the permutations in the statement are listed according to the lexicographic order, and the computation of the corresponding signatures is something elementary, by using the various rules from Theorem 2.35.
\end{proof}

As another application, we have the following key result:

\index{transpose matrix}

\begin{theorem}
We have the formula
$$\det A=\det A^t$$
valid for any square matrix $A$.
\end{theorem}

\begin{proof}
This follows from the formula in Theorem 2.36. Indeed, we have:
\begin{eqnarray*}
\det A^t
&=&\sum_{\sigma\in S_N}\varepsilon(\sigma)(A^t)_{1\sigma(1)}\ldots(A^t)_{N\sigma(N)}\\
&=&\sum_{\sigma\in S_N}\varepsilon(\sigma)A_{\sigma(1)1}\ldots A_{\sigma(N)N}\\
&=&\sum_{\sigma\in S_N}\varepsilon(\sigma)A_{1\sigma^{-1}(1)}\ldots A_{N\sigma^{-1}(N)}\\
&=&\sum_{\sigma\in S_N}\varepsilon(\sigma^{-1})A_{1\sigma^{-1}(1)}\ldots A_{N\sigma^{-1}(N)}\\
&=&\sum_{\sigma\in S_N}\varepsilon(\sigma)A_{1\sigma(1)}\ldots A_{N\sigma(N)}\\
&=&\det A
\end{eqnarray*}

Thus, we are led to the formula in the statement.
\end{proof}

Good news, this is the end of the general theory that we wanted to develop. We have now in our bag all the needed techniques for computing the determinant.

\bigskip

Here is however a nice and important example of a determinant, whose computation uses some interesting new techniques, going beyond what has been said above:

\index{Vandermonde formula}

\begin{theorem}
We have the Vandermonde determinant formula
$$\begin{vmatrix}
1&1&1&\ldots\ldots&1\\
x_1&x_2&x_3&\ldots\ldots&x_N\\
x_1^2&x_2^2&x_3^2&\ldots\ldots&x_N\\
\vdots&\vdots&\vdots&&\vdots\\
\vdots&\vdots&\vdots&&\vdots\\
x_1^{N-1}&x_2^{N-1}&x_3^{N-1}&\ldots\ldots&x_N^{N-1}
\end{vmatrix}
=\prod_{i>j}(x_i-x_j)$$
valid for any $x_1,\ldots,x_N\in\mathbb R$.
\end{theorem}

\begin{proof}
Let us first do some checks. At $N=2$ the formula holds indeed:
$$\begin{vmatrix}
1&1\\
a&b\end{vmatrix}=b-a$$

At $N=3$ now, the Vandermonde formula holds too, as shown by:
\begin{eqnarray*}
\begin{vmatrix}
1&1&1\\
a&b&c\\
a^2&b^2&c^2\end{vmatrix}
&=&bc^2+ab^2+a^2c-a^2b-b^2c-ac^2\\
&=&(bc^2-ac^2)+(ab^2-a^2b)+(a^2c-b^2c)\\
&=&(b-a)(c^2+ab-ac-bc)\\
&=&(b-a)(c-a)(c-b)
\end{eqnarray*}

In general, by expanding over the columns, we can see that the determinant in question, say $D$, is a polynomial in the variables $x_1,\ldots,x_N$, having degree $N-1$ in each variable. Now observe that when setting $x_i=x_j$, for some indices $i\neq j$, our matrix will have two identical columns, and so its determinant $D$ will vanish:
$$x_i=x_j\implies D=0$$

But this gives us the key to the computation of $D$. Indeed, $D$ must be divisible by $x_i-x_j$ for any $i\neq j$, and so we must have a formula of the following type:
$$D=c\prod_{i>j}(x_i-x_j)$$

Moreover, since the product on the right is, exactly as $D$ itself, a polynomial in the variables $x_1,\ldots,x_N$, having degree $N-1$ in each variable, we conclude that the quantity $c$ must be a constant, not depending on any of the variables $x_1,\ldots,x_N$:
$$c\in\mathbb R$$

In order to finish the computation, it remains to find the value of this constant $c$. But this can be done for instance by recurrence, and we obtain $c=1$, as desired.
\end{proof}

Summarizing, we are now experts in the computation of the determinant, and moving on, we should investigate the next problem, namely the diagonalization one.

\bigskip

But here, as a key input, we know from Theorem 2.19 that the eigenvalues of a matrix $A\in M_N(\mathbb R)$ appear as roots of the characteristic polynomial:
$$P(x)=\det(A-x1_N)$$

Thus, with the determinant theory developed above, we can in principle compute these eigenvalues, and solve the diagonalization problem afterwards. 

\bigskip

The problem, however, is that certain real matrices can have characteristic polynomials of type $P(x)=x^2+1$, and this suggests that these matrices might be not diagonalizable over $\mathbb R$, but be diagonalizable over $\mathbb C$ instead. And so, before getting into diagonalization problems, we must upgrade our theory, and talk about complex matrices. We will do this in the next chapter, and afterwards, we will go back to the diagonalization problem.

\section*{2e. Exercises}

There has been a lot of exciting theory in this chapter, with some details sometimes missing, and our exercises will be mainly about this. First, we have:

\begin{exercise}
Fill in all the geometric details in the basic theory of the determinant, by using the same type of arguments as those in the proof of
$$\det\begin{pmatrix}a&b\\c&d\end{pmatrix}=ad-bc$$
which was fully proved in the above, namely geometric manipulations, and Thales.
\end{exercise}

To be more precise here, passed some issues with the sign and orientation, which are all elementary, the above $2\times2$ determinant formula was subject of Theorem 2.6, coming with a full and honest proof. The problem is that of using the same arguments, namely basic geometry, as to have a full proof of Theorem 2.16 and Theorem 2.21 as well.

\begin{exercise}
Prove with full details, based on the above, that the determinant of the systems of vectors
$$\det:\mathbb R^N\times\ldots\times\mathbb R^N\to\mathbb R$$
is multilinear, alternate and unital, and unique with these properties. Then try to prove as well this directly, without any reference to geometry.
\end{exercise}

To be more precise, in what regards the first question, this is something that we already discussed in the above, with only a few details missing, and the problem is that of recovering these details. As for the second question, this is something more tricky, and there are several possible approaches here, all being interesting and enjoyable.

\begin{exercise}
Work out, with full details, the theory of the signature map
$$\varepsilon:S_N\to\{\pm1\}$$
as outlined in Theorem 2.35 and its proof. 
\end{exercise}

As before, these are things that we already discussed, with a few details missing.

\begin{exercise}
Prove that for a matrix $H\in M_N(\pm1)$, we have
$$|\det H|\leq N^{N/2}$$
and then find the maximizers of $|\det H|$, at small values of $N$.
\end{exercise}

Here the first question is theoretical, and its proof should not be difficult. As for the second question, which is quite tricky, the higher the $N\in\mathbb N$ you get to, the better.

\chapter{Complex matrices}

\section*{3a. Complex numbers}

We have seen that the study of the real matrices $A\in M_N(\mathbb R)$ suggests the use of the complex numbers. Indeed, even simple matrices like the $2\times2$ ones can, at least in a formal sense, have complex eigenvalues. In what follows we discuss the complex matrices $A\in M_N(\mathbb C)$. We will see that the theory here is much more complete than in the real case. As an application, we will solve in this way problems left open in the real case.

\bigskip

Let us begin with the complex numbers. There is a lot of magic here, and we will carefully explain this material. Their definition is as follows:

\index{complex numbers}

\begin{definition}
The complex numbers are variables of the form
$$x=a+ib$$
which add in the obvious way, and multiply according to the rule $i^2=-1$.
\end{definition}

In other words, we consider variables as above, without bothering for the moment with their precise meaning. Now consider two such complex numbers:
$$x=a+ib\quad,\quad 
y=c+id$$

The formula for the sum is then the obvious one, as follows:
$$x+y=(a+c)+i(b+d)$$

As for the formula of the product, by using the rule $i^2=-1$, we obtain:
\begin{eqnarray*}
xy
&=&(a+ib)(c+id)\\
&=&ac+iad+ibc+i^2bd\\
&=&ac+iad+ibc-bd\\
&=&(ac-bd)+i(ad+bc)
\end{eqnarray*}

Thus, the complex numbers as introduced above are well-defined. The multiplication formula is of course quite tricky, and hard to memorize, but we will see later some alternative ways, which are more conceptual, for performing the multiplication. 

\bigskip

The advantage of using the complex numbers comes from the fact that the equation $x^2=1$ has now a solution, $x=i$. In fact, this equation has two solutions, namely:
$$x=\pm i$$

This is of course very good news. More generally, we have the following result:

\begin{theorem}
The complex solutions of $ax^2+bx+c=0$ with $a,b,c\in\mathbb R$ are
$$x_{1,2}=\frac{-b\pm\sqrt{b^2-4ac}}{2a}$$
with the square root of negative real numbers being defined as $\sqrt{-m}=\pm i\sqrt{m}$.
\end{theorem}

\begin{proof}
We can write our equation in the following way:
\begin{eqnarray*}
ax^2+bx+c=0
&\iff&x^2+\frac{b}{a}x+\frac{c}{a}=0\\
&\iff&\left(x+\frac{b}{2a}\right)^2-\frac{b^2}{4a^2}+\frac{c}{a}=0\\
&\iff&\left(x+\frac{b}{2a}\right)^2=\frac{b^2-4ac}{4a^2}\\
&\iff&x+\frac{b}{2a}=\pm\frac{\sqrt{b^2-4ac}}{2a}
\end{eqnarray*}

Thus, we are led to the conclusion in the statement.
\end{proof}

We will be back later to this, with generalizations. Getting back now to Definition 3.1 as it is, we can represent the complex numbers in the plane, as follows:

\index{complex numbers}

\begin{proposition}
The complex numbers, written as usual
$$x=a+ib$$
can be represented in the plane, according to the following identification:
$$x=\binom{a}{b}$$
With this convention, the sum of complex numbers is the usual sum of vectors. 
\end{proposition}

\begin{proof}
Consider indeed two arbitrary complex numbers:
$$x=a+ib\quad,\quad 
y=c+id$$

Their sum is then by definition the following complex number:
$$x+y=(a+c)+i(b+d)$$

Now let us represent $x,y$ in the plane, as in the statement:
$$x=\binom{a}{b}\quad,\quad y=\binom{c}{d}$$

In this picture, their sum is given by the following formula:
$$x+y=\binom{a+c}{b+d}$$

But this is indeed the vector corresponding to $x+y$, so we are done.
\end{proof}

Observe that in the above picture, the real numbers correspond to the numbers on the $Ox$ axis. As for the purely imaginary numbers, these lie on the $Oy$ axis, with:
$$i=\binom{0}{1}$$

All this is very nice, but in order to understand now the multiplication, we must do something more complicated, namely using polar coordinates. Let us start with:

\index{polar coordinates}
\index{modulus of complex number}
\index{argument of complex number}

\begin{definition}
The complex numbers $x=a+ib$ can be written in polar coordinates,
$$x=r(\cos t+i\sin t)$$
with the connecting formulae being
$$a=r\cos t\quad,\quad 
b=r\sin t$$
and in the other sense being
$$r=\sqrt{a^2+b^2}\quad,\quad 
\tan t=b/a$$
and with $r,t$ being called modulus, and argument.
\end{definition}

There is a clear relation here with the vector notation from Proposition 3.3, because $r$ is the length of the vector, and $t$ is the angle made by the vector with the $Ox$ axis. As a basic example here, the number $i$ takes the following form:
$$i=\cos\left(\frac{\pi}{2}\right)+i\sin\left(\frac{\pi}{2}\right)$$

The point now is that in polar coordinates, the multiplication formula for the complex numbers, which was so far something quite opaque, takes a very simple form:

\begin{theorem}
Two complex numbers written in polar coordinates,
$$x=r(\cos s+i\sin s)\quad,\quad 
y=p(\cos t+i\sin t)$$
multiply according to the following formula:
$$xy=rp(\cos(s+t)+i\sin(s+t))$$
In other words, the moduli multiply, and the arguments sum up.
\end{theorem}

\begin{proof}
This can be proved by doing some trigonometry, as follows:

\medskip

(1) Recall first the definition of $\sin,\cos$, as being the sides of a right triangle having angle $t$. Our first claim is that we have the Pythagoras' theorem, namely:
$$\sin^2t+\cos^2t=1$$

But this comes from the following well-known, remarkable picture, with the edges of the outer and inner square being respectively $\sin t+\cos t$ and 1:
$$
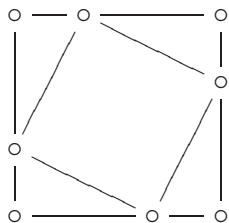
\xymatrix@R=13pt@C=13pt{
\circ\ar@{-}[r]\ar@{-}[dd]&\circ\ar@{-}[rr]\ar@{-}[drr]&&\circ\ar@{-}[d]\\
&&&\circ\ar@{-}[dd]^{\sin t}&\\
\circ\ar@{-}[d]\ar@{-}[uur]\ar@{-}[drr]&&&&\\
\circ\ar@{-}[rr]&&\circ\ar@{-}[r]_{\cos t}\ar@{-}[uur]^1&\circ
}$$

Indeed, when computing the area of the outer square, in two ways, we obtain:
$$(\sin t+\cos t)^2=1+4\times\frac{\sin t\cos t}{2}$$

Now when expanding we obtain $\sin^2t+\cos^2t=1$, as claimed.

\medskip

(2) Next in line, our claim is that we have the following formulae:
$$\sin(s+t)=\cos s\sin t+\sin s\cos t$$
$$\cos(s+t)=\cos s\cos t-\sin s\sin t$$

To be more precise, let us first establish this formula. In order to do so, consider the following picture, consisting of a length 1 line segment, with angles $s,t$ drawn on each side, and with everything being completed, and lengths computed, as indicated:
$$
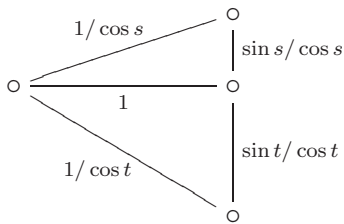
\xymatrix@R=15pt@C=70pt{
&\circ\ar@{-}[d]^{\sin s/\cos s}\\
\circ\ar@{-}[ur]^{1/\cos s}\ar@{-}[r]_1\ar@{-}[ddr]_{1/\cos t}&\circ\ar@{-}[dd]^{\sin t/\cos t}\\
\\
&\circ
}$$

Now let us compute the area of the big triangle, or rather the double of that area. We can do this in two ways, either directly, with a formula involving $\sin(s+t)$, or by using the two small triangles, involving functions of $s,t$. We obtain in this way:
$$\frac{1}{\cos s}\cdot\frac{1}{\cos t}\cdot\sin(s+t)=\frac{\sin s}{\cos s}\cdot 1+\frac{\sin t}{\cos t}\cdot 1$$

But this gives the formula for $\sin(s+t)$ claimed above. Now by using this formula for $\sin(s+t)$ we can deduce as well the formula for $\cos(s+t)$, as follows:
\begin{eqnarray*}
\cos(s+t)
&=&\sin\left(\frac{\pi}{2}-s-t\right)\\
&=&\sin\left[\left(\frac{\pi}{2}-s\right)+(-t)\right]\\
&=&\sin\left(\frac{\pi}{2}-s\right)\cos(-t)+\cos\left(\frac{\pi}{2}-s\right)\sin(-t)\\
&=&\cos s\cos t-\sin s\sin t
\end{eqnarray*}

(3) Now back to complex numbers, we want to prove that $x=r(\cos s+i\sin s)$ and $y=p(\cos t+i\sin t)$ multiply according to the following formula:
$$xy=rp(\cos(s+t)+i\sin(s+t))$$

We can assume that we have $r=p=1$, by dividing everything by these numbers. Now with this assumption made, we have the following computation:
\begin{eqnarray*}
xy
&=&(\cos s+i\sin s)(\cos t+i\sin t)\\
&=&(\cos s\cos t-\sin s\sin t)+i(\cos s\sin t+\sin s\cos t)\\
&=&\cos(s+t)+i\sin(s+t)
\end{eqnarray*}

Thus, we are led to the conclusion in the statement.
\end{proof}

The above result, which was based on some non-trivial trigonometry, is quite powerful. As a basic application of it, we can now compute powers, as follows:

\index{powers of complex number}

\begin{theorem}
The powers of a complex number, written in polar form,
$$x=r(\cos t+i\sin t)$$
are given by the following formula, valid for any exponent $k\in\mathbb N$:
$$x^k=r^k(\cos kt+i\sin kt)$$
Moreover, this formula holds in fact for any $k\in\mathbb Z$, and even for any $k\in\mathbb Q$.
\end{theorem}

\begin{proof}
Given a complex number $x$, written in polar form as above, and an exponent $k\in\mathbb N$, we have indeed the following computation, with $k$ terms everywhere:
\begin{eqnarray*}
x^k
&=&x\ldots x\\
&=&r(\cos t+i\sin t)\ldots r(\cos t+i\sin t)\\
&=&r\ldots r(\cos(t+\ldots+t)+i\sin(t+\ldots+t))\\
&=&r^k(\cos kt+i\sin kt)
\end{eqnarray*}

Thus, we are done with the case $k\in\mathbb N$. Regarding now the generalization to the case $k\in\mathbb Z$, it is enough here to do the verification for $k=-1$, where the formula is:
$$x^{-1}=r^{-1}(\cos(-t)+i\sin(-t))$$

But this number $x^{-1}$ is indeed the inverse of $x$, because:
\begin{eqnarray*}
xx^{-1}
&=&r(\cos t+i\sin t)\cdot r^{-1}(\cos(-t)+i\sin(-t))\\
&=&\cos(t-t)+i\sin(t-t)\\
&=&\cos 0+i\sin 0\\
&=&1
\end{eqnarray*}

Finally, regarding the generalization to the case $k\in\mathbb Q$, it is enough to do the verification for exponents of type $k=1/n$, with $n\in\mathbb N$. The claim here is that:
$$x^{1/n}=r^{1/n}\left[\cos\left(\frac{t}{n}\right)+i\sin\left(\frac{t}{n}\right)\right]$$

In order to prove this, let us compute the $n$-th power of this number. We can use the power formula for the exponent $n\in\mathbb N$, that we already established, and we obtain:
\begin{eqnarray*}
(x^{1/n})^n
&=&(r^{1/n})^n\left[\cos\left(n\cdot\frac{t}{n}\right)+i\sin\left(n\cdot\frac{t}{n}\right)\right]\\
&=&r(\cos t+i\sin t)\\
&=&x
\end{eqnarray*}

Thus, we have indeed a $n$-th root of $x$, and our proof is now complete.
\end{proof}

We should mention that there is a bit of ambiguity in the above, in the case of the exponents $k\in\mathbb Q$, due to the fact that the square roots, and the higher roots as well, can take multiple values, in the complex number setting. We will be back to this.

\section*{3b. Euler formula}

We would like to discuss now the final and most convenient writing of the complex numbers, which is a well-known variation on the polar writing, as follows:
$$x=re^{it}$$

In what follows we will not really need the true power of this formula, which is of analytic nature, due to occurrence of the number $e$. However, we would like to use the notation $x=re^{it}$, as everyone does, among others because it simplifies the writing. The point indeed with the above formula comes from the following deep result:

\index{polar writing}

\begin{theorem}
We have the following formula, valid for any $t\in\mathbb R$,
$$e^{it}=\cos t+i\sin t$$
where $e=2.7182\ldots$ is the usual constant from analysis.
\end{theorem}

\begin{proof}
This is something quite tricky, the idea being as follows:

\medskip

(1) As a first question, what is $e$? In answer, there are two equivalent definitions of it, one as a limit, and the other one as the sum of a series, as follows:
$$e=\lim_{n\to\infty}\left(1+\frac{1}{n}\right)^n=\sum_{k=0}^\infty\frac{1}{k!}$$

Next, what is the exponential function? Again, we have two equivalent definitions here, which can be deduced from the above two formulae, as follows:
$$e^x=\lim_{n\to\infty}\left(1+\frac{x}{n}\right)^n=\sum_{k=0}^\infty\frac{x^k}{k!}$$

(2) Next, can we really apply this exponential function to complex numbers? And the answer here is yes, due to the following estimate, based on the series approach:
$$|e^x|
=\left|\sum_{k=0}^\infty\frac{x^k}{k!}\right|
\leq\sum_{k=0}^\infty\left|\frac{x^k}{k!}\right|
=\sum_{k=0}^\infty\frac{|x|^k}{k!}
=e^{|x|}<\infty$$

Now with this done, what can we say about $e^x$? And as a basic fact here, we have:
\begin{eqnarray*}
e^{x+y}
&=&\sum_{k=0}^\infty\frac{(x+y)^k}{k!}\\
&=&\sum_{k=0}^\infty\sum_{s=0}^k\binom{k}{s}\cdot\frac{x^sy^{k-s}}{k!}\\
&=&\sum_{k=0}^\infty\sum_{s=0}^k\frac{x^sy^{k-s}}{s!(k-s)!}\\
&=&e^xe^y
\end{eqnarray*}

(3) Our next claim is that $e^x$ is continuous. Indeed, at $x=0$ this comes from:
$$|e^t-1|
=\left|\sum_{k=1}^\infty\frac{t^k}{k!}\right|
\leq\sum_{k=1}^\infty\left|\frac{t^k}{k!}\right|
=\sum_{k=1}^\infty\frac{|t|^k}{k!}
=e^{|t|}-1$$

As for the continuity of $x\to e^x$ in general, this can be deduced as follows:
$$\lim_{t\to0}e^{x+t}
=\lim_{t\to0}e^xe^t
=e^x\lim_{t\to0}e^t
=e^x\cdot 1
=e^x$$

(4) Getting now towards what we want to do, our first claim is that for $t\in\mathbb R$ we have $e^{it}\in\mathbb T$, unit circle. In order to prove this, observe that we have, for any $x\in\mathbb C$:
$$e^{\bar{x}}=\sum_{k=0}^\infty\frac{\bar{x}^k}{k!}=\overline{\sum_{k=0}^\infty\frac{x^k}{k!}}=\overline{e^x}$$

Also, we have as well the following computation, again for any $x\in\mathbb C$:
$$e^xe^{-x}=e^{x-x}=e^0=1\implies (e^x)^{-1}=e^{-x}$$

But with these two formulae in hand, we can prove our claim. Indeed, the above two formulae, applied with $x=it$, with $t\in\mathbb R$, give the following equalities:
$$e^{-it}=\overline{e^{it}}\quad,\quad
(e^{it})^{-1}=e^{-it}$$

Thus the number $z=e^{it}$ has the property $z^{-1}=\bar{z}$, and so $z\in\mathbb T$, as claimed.

\medskip

(5) Time now for the proof of $e^{it}=\cos t+i\sin t$. We know that the operation $t\to e^{it}$ is continuous, and maps sums in $\mathbb R$ to products in $\mathbb T$. But in view of this, skipping some details, that we will leave as an exercise, we can conclude that this operation must appear by ``wrapping". That is, we must have a formula as follows, for a certain $\alpha\in\mathbb R$:
$$e^{it}=\cos(\alpha t)+i\sin(\alpha t)$$

In order now to find the parameter $\alpha\in\mathbb R$, let us look at what happens around $t=0$. And here, we have the following elementary estimate, obtained by truncating $\exp$:
$$e^{it}\simeq 1+it$$

On the other hand, according to some basic trigonometry for $\sin$, $\cos$, done in the old way, on the unit circle, we have as well the following estimate, again around $t=0$:
$$\cos(\alpha t)+i\sin(\alpha t)\simeq 1+i\alpha t$$

Thus, we must have $\alpha=1$, which gives the Euler formula, as desired.

\medskip

(6) As an alternative proof for the Euler formula, which is certainly quicker, but unfortunately hides what is going on, geometrically, we can kill the problem with calculus. Indeed, we have the following formulae, with the first one being clear, and the other two being obtained from the usual formulae of $\sin(x+t)$ and $\cos(x+t)$, with $t\simeq 0$:
$$(e^x)'=e^x\quad,\quad(\sin x)'=\cos x\quad,\quad(\cos x)'=-\sin x$$

In order to prove the Euler formula, consider the following function $f:\mathbb R\to\mathbb C$:
$$f(t)=\frac{\cos t+i\sin t}{e^{it}}$$

By using standard calculus rules, the derivative of this function is given by:
\begin{eqnarray*}
f'(t)
&=&(e^{-it}(\cos t+i\sin t))'\\
&=&-ie^{-it}(\cos t+i\sin t)+e^{-it}(-\sin t+i\cos t)\\
&=&e^{-it}(-i\cos t+\sin t)+e^{-it}(-\sin t+i\cos t)\\
&=&0
\end{eqnarray*}

Thus $f$ is constant, equalling $f(0)=1$, and we have proved the Euler formula.

\medskip

(7) Finally, no discussion about the Euler formula would be complete without performing the following computation, based on the definition of the exponential:
\begin{eqnarray*}
e^{it}
&=&\sum_{k=0}^\infty\frac{(it)^k}{k!}\\
&=&\sum_{k=2l}\frac{(it)^k}{k!}+\sum_{k=2l+1}\frac{(it)^k}{k!}\\
&=&\sum_{l=0}^\infty(-1)^l\frac{t^{2l}}{(2l)!}+i\sum_{l=0}^\infty(-1)^l\frac{t^{2l+1}}{(2l+1)!}
\end{eqnarray*}

Indeed, we obtain in this way, via Euler, the following formulae for $\cos$ and $\sin$:
$$\cos t=\sum_{l=0}^\infty(-1)^l\frac{t^{2l}}{(2l)!}\quad,\quad 
\sin t=\sum_{l=0}^\infty(-1)^l\frac{t^{2l+1}}{(2l+1)!}$$

Which is nice, these being the Taylor series of $\cos$ and $\sin$, coming from the formulae $\sin'=\cos$ and $\cos'=-\sin$, discussed in (6). However, and here comes the point, the fact that we have equalities $=$ as above, instead of just $\simeq$, and with these equalities being valid at any $t\in\mathbb R$, is something well beyond the theory of real Taylor series, coming from the Euler formula, proved as in (5), or as in (6). And, good to know, all this.
\end{proof}

As a first interesting consequence of the Euler formula, we have:

\begin{theorem}
We have the following formula,
$$e^{\pi i}=-1$$
and we have $E=mc^2$ as well.
\end{theorem}

\begin{proof}
We have two assertions here, the idea being as follows:

\medskip

(1) The first formula, $e^{\pi i}=-1$, which is actually the main formula in mathematics, comes from Theorem 3.7, by setting $t=\pi$. Indeed, we obtain:
\begin{eqnarray*}
e^{\pi i}
&=&\cos\pi+i\sin\pi\\
&=&-1+i\cdot 0\\
&=&-1
\end{eqnarray*}

(2) As for $E=mc^2$, which is the main formula in physics, this is something deep as well. Although we will not really need it here, we recommend learning it too, for symmetry reasons between math and physics, say from Feynman \cite{fe1}, \cite{fe2}, \cite{fe3}.
\end{proof}

Now back to our $x=re^{it}$ objectives, with the above theory in hand we can indeed use from now on this notation, the complete statement being as follows:

\index{polar coordinates}

\begin{theorem}
The complex numbers $x=a+ib$ can be written in polar coordinates,
$$x=re^{it}$$
with the connecting formulae being
$$a=r\cos t\quad,\quad 
b=r\sin t$$
and in the other sense being
$$r=\sqrt{a^2+b^2}\quad,\quad 
\tan t=b/a$$
and with $r,t$ being called modulus, and argument.
\end{theorem}

\begin{proof}
This is just a reformulation of Definition 3.4, by using the formula $e^{it}=\cos t+i\sin t$ from Theorem 3.7, and multiplying everything by $r$.
\end{proof}

We can now go back to the basics, and we have the following result:

\index{multiplication of complex numbers}

\begin{theorem}
In polar coordinates, the complex numbers multiply as
$$re^{is}\cdot pe^{it}=rp\,e^{i(s+t)}$$
with the arguments $s,t$ being taken modulo $2\pi$.
\end{theorem}

\begin{proof}
This is something that know from Theorem 3.5, reformulated by using the notations from Theorem 3.9. Observe that this follows as well from $e^{x+y}=e^xe^y$.
\end{proof}

We can now investigate more complicated operations, as follows:

\begin{theorem}
We have the following operations on the complex numbers:
\begin{enumerate}
\item Inversion: $(re^{it})^{-1}=r^{-1}e^{-it}$.

\item Square roots: $\sqrt{re^{it}}=\pm\sqrt{r}e^{it/2}$.

\item Powers: $(re^{it})^a=r^ae^{ita}$.
\end{enumerate}
\end{theorem}

\begin{proof}
This is something that we already know, from Theorem 3.6, but we can now discuss all this, from a more conceptual viewpoint, the idea being as follows:

\medskip

(1) We have indeed the following computation, using Theorem 3.10:
$$(re^{it})(r^{-1}e^{-it})
=rr^{-1}\cdot e^{i(t-t)}
=1$$

(2) Once again by using Theorem 3.10, we have:
$$(\pm\sqrt{r}e^{it/2})^2
=(\sqrt{r})^2e^{i(t/2+t/2)}
=re^{it}$$

(3) Given an arbitrary number $a\in\mathbb R$, we can define, as stated:
$$(re^{it})^a=r^ae^{ita}$$

And, due to Theorem 3.10, this operation $x\to x^a$ is indeed the correct one.
\end{proof}

We can now go back to the degree 2 equations, and we have:

\index{degree 2 equation}

\begin{theorem}
The complex solutions of $ax^2+bx+c=0$ with $a,b,c\in\mathbb C$ are
$$x_{1,2}=\frac{-b\pm\sqrt{b^2-4ac}}{2a}$$
with the square root of complex numbers being defined as above.
\end{theorem}

\begin{proof}
This is clear, the computations being the same as in the real case. To be more precise, our degree 2 equation can be written as follows:
$$\left(x+\frac{b}{2a}\right)^2=\frac{b^2-4ac}{4a^2}$$

Now since we know from Theorem 3.11 (2) that any complex number has a square root, we are led to the conclusion in the statement.
\end{proof}

More generally now, we have the following key result, in arbitrary degree:

\index{roots of polynomials}
\index{complex roots}

\begin{theorem}
Any polynomial $P\in\mathbb C[X]$ decomposes as
$$P=c(X-a_1)\ldots (X-a_N)$$
with $c\in\mathbb C$ and with $a_1,\ldots,a_N\in\mathbb C$.
\end{theorem}

\begin{proof}
The problem is that of proving that our polynomial has at least one root, because afterwards we can proceed by recurrence. We prove this by contradiction. So, assume that $P$ has no roots, and pick a number $z\in\mathbb C$ where $|P|$ attains its minimum:
$$|P(z)|=\min_{x\in\mathbb C}|P(x)|>0$$ 

Since $Q(t)=P(z+t)-P(z)$ is a polynomial which vanishes at $t=0$, this polynomial must be of the form $ct^k$ + higher terms, with $c\neq0$, and with $k\geq1$ being an integer. We obtain from this that, with $t\in\mathbb C$ small, we have the following estimate:
$$P(z+t)\simeq P(z)+ct^k$$

Now let us write $t=rw$, with $r>0$ small, and with $|w|=1$. Our estimate becomes:
$$P(z+rw)\simeq P(z)+cr^kw^k$$

Now recall that we have assumed $P(z)\neq0$. We can therefore choose $w\in\mathbb T$ such that $cw^k$ points in the opposite direction to that of $P(z)$, and we obtain in this way:
\begin{eqnarray*}
|P(z+rw)|
&\simeq&|P(z)+cr^kw^k|\\
&=&|P(z)|(1-|c|r^k)
\end{eqnarray*}

Now by choosing $r>0$ small enough, as for the error in the first estimate to be small, and overcame by the negative quantity $-|c|r^k$, we obtain from this:
$$|P(z+rw)|<|P(z)|$$

But this contradicts our definition of $z\in\mathbb C$, as a point where $|P|$ attains its minimum. Thus $P$ has a root, and by recurrence it has $N$ roots, as stated.
\end{proof}

All this is very nice, and we will see applications in a moment. As a last topic now regarding the complex numbers, we have the roots of unity:

\index{roots of unity}

\begin{theorem}
The equation $x^N=1$ has $N$ complex solutions, namely
$$\left\{w^k\Big|k=0,1,\ldots,N-1\right\}\quad,\quad w=e^{2\pi i/N}$$
which are called roots of unity of order $N$.
\end{theorem}

\begin{proof}
This follows from Theorem 3.10. Indeed, with $x=re^{it}$ our equation reads:
$$r^Ne^{itN}=1$$

Thus $r=1$, and $t\in[0,2\pi)$ must be a multiple of $2\pi/N$, as stated.
\end{proof}

As an illustration here, the roots of unity of small order, along with some of their basic properties, which are very useful for computations, are as follows:

\medskip

\underline{$N=1$}. Here the unique root of unity is 1.

\medskip

\underline{$N=2$}. Here we have two roots of unity, namely 1 and $-1$.

\medskip

\underline{$N=3$}. Here we have 1, then $w=e^{2\pi i/3}$, and then $w^2=\bar{w}=e^{4\pi i/3}$.

\medskip

\underline{$N=4$}. Here the roots of unity, read as usual counterclockwise, are $1,i,-1,-i$.

\medskip

\underline{$N=5$}. Here, with $w=e^{2\pi i/5}$, the roots of unity are $1,w,w^2,w^3,w^4$.

\medskip

\underline{$N=6$}. Here a useful alternative writing is $\{\pm1,\pm w,\pm w^2\}$, with $w=e^{2\pi i/3}$.

\medskip

The roots of unity are very useful variables, and have many interesting properties. As a first application, we can now solve the ambiguity questions related to the extraction of $N$-th roots, from Theorem 3.6 and Theorem 3.11, the statement being as follows:

\index{roots of unity}

\begin{theorem}
Any nonzero $x=re^{it}$ has exactly $N$ roots of order $N$, namely
$$y=r^{1/N}e^{it/N}$$
multiplied by the $N$ roots of unity of order $N$.
\end{theorem}

\begin{proof}
We must solve the equation $z^N=x$, over the complex numbers. Since the number $y$ in the statement clearly satisfies $y^N=x$, our equation reformulates as: 
$$z^N=x\iff z^N=y^N\iff \left(\frac{z}{y}\right)^N=1$$

Thus, we are led to the conclusion in the statement.
\end{proof}

The roots of unity appear in connection with many other questions, and there are many useful formulae relating them, which are good to know, as for instance:

\index{roots of unity}
\index{barycenter}

\begin{theorem}
The roots of unity, $\{w^k\}$ with $w=e^{2\pi i/N}$, have the property
$$\sum_{k=0}^{N-1}(w^k)^s=N\delta_{N|s}$$
for any exponent $s\in\mathbb N$, where on the right we have a Kronecker symbol.
\end{theorem}

\begin{proof}
The numbers in the statement, when written more conveniently as $(w^s)^k$ with $k=0,\ldots,N-1$, form a certain regular polygon in the plane $P_s$. Thus, if we denote by $C_s$ the barycenter of this polygon, we have the following formula:
$$\frac{1}{N}\sum_{k=0}^{N-1}w^{ks}=C_s$$

Now observe that in the case $N\slash\hskip-1.6mm|\,s$ our polygon $P_s$ is non-degenerate, circling around the unit circle, and having center $C_s=0$. As for the case $N|s$, here the polygon is degenerate, lying at 1, and having center $C_s=1$. Thus, we have the following formula:
$$C_s=\delta_{N|s}$$

Thus, we obtain the formula in the statement.
\end{proof}

\section*{3c. Complex matrices}

Back now to linear algebra, our first task will be that of extending the results that we know, from the real case, to the complex case. We first have:

\index{linear map}
\index{rectangular matrix}

\begin{theorem}
The linear maps $f:\mathbb C^N\to\mathbb C^M$ are the maps of the form
$$f(x)=Ax$$
with $A$ being a rectangular matrix, $A\in M_{M\times N}(\mathbb C)$.
\end{theorem}

\begin{proof}
This follows as in the real case. Indeed, $f:\mathbb C^N\to\mathbb C^M$ must send a vector $x\in\mathbb C^N$ to a certain vector $f(x)\in\mathbb C^M$, all whose components are linear combinations of the components of $x$. Thus, we can write, for certain complex numbers $a_{ij}\in\mathbb C$:
$$f\begin{pmatrix}
x_1\\
\vdots\\
\vdots\\
x_N
\end{pmatrix}
=\begin{pmatrix}
a_{11}x_1+\ldots+a_{1N}x_N\\
\vdots\\
\vdots\\
a_{M1}x_1+\ldots+a_{MN}x_N
\end{pmatrix}$$

But the parameters $a_{ij}\in\mathbb C$ can be regarded as being the entries of a matrix:
$$A=(a_{ij})\in M_{M\times N}(\mathbb C)$$

Now with the usual convention for the rectangular matrix multiplication, exactly as in the real case, the above formula is precisely the one in the statement.
\end{proof}

We have as well the following result, again inspired from the real case:

\index{invertible matrix}

\begin{theorem}
A linear map $f:\mathbb C^N\to\mathbb C^M$, written as 
$$f(v)=Av$$
is invertible precisely when $A$ is invertible, and in this case we have $f^{-1}(v)=A^{-1}v$.
\end{theorem}

\begin{proof}
As in the real case, with the convention $f_A(v)=Av$, we have the following multiplication formula for such linear maps:
$$f_Af_B(v)=f_{AB}(v)$$

But this shows that $f_Af_B=1$ is equivalent to $AB=1$, as desired. 
\end{proof}

With respect to the real case, some subtleties appear at the level of the scalar products, isometries and projections. The basic theory here is as follows:

\index{scalar product}
\index{adjoint matrix}
\index{isometry}
\index{unitary matrix}
\index{projection}
\index{orthogonal projection}
\index{rank 1 projection}

\begin{theorem}
Consider the usual scalar product $<x,y>=\sum_ix_i\bar{y}_i$ on $\mathbb C^N$.
\begin{enumerate}
\item We have the following formula, with $(A^*)_{ij}=\bar{A}_{ji}$ being the adjoint matrix:
$$<Ax,y>=<x,A^*y>$$

\item A linear map $f:\mathbb C^N\to\mathbb C^N$, written as $f(x)=Ux$ with $U\in M_N(\mathbb C)$, is an isometry precisely when $U$ is unitary, in the sense that:
$$U^*=U^{-1}$$

\item  A linear map $f:\mathbb C^N\to\mathbb C^N$, written as $f(x)=Px$ with $P\in M_N(\mathbb C)$, is a porojection precisely when $P$ is projection, in the sense that: 
$$P^2=P^*=P$$

\item Also, the formula for the rank $1$ projections is $P_x=\frac{1}{||x||^2}(x_i\bar{x}_j)_{ij}$.
\end{enumerate}
\end{theorem}

\begin{proof}
This follows as in the real case, with modifications where needed:

\medskip

(1) By using the standard basis of $\mathbb C^N$, we want to prove that for any $i,j$ we have:
$$<Ae_j,e_i>=<e_j,A^*e_i>$$

The scalar product being now antisymmetric, this is the same as proving that:
$$<Ae_j,e_i>=\overline{<A^*e_i,e_j>}$$

On the other hand, for any matrix $M$ we have the following formula:
$$M_{ij}=<Me_j,e_i>$$

Thus, the formula to be proved simply reads $A_{ij}=\overline{(A^*)_{ji}}$, as desired.

\medskip

(2) Let first recall that we can pass from scalar products to distances, as follows:
$$||x||=\sqrt{<x,x>}$$

Conversely, we can compute the scalar products in terms of distances, by using the complex polarization identity, which is as follows:
\begin{eqnarray*}
&&||x+y||^2-||x-y||^2+i||x+iy||^2-i||x-iy||^2\\
&=&||x||^2+||y||^2-||x||^2-||y||^2+i||x||^2+i||y||^2-i||x||^2-i||y||^2\\
&&+2Re(<x,y>)+2Re(<x,y>)+2iIm(<x,y>)+2iIm(<x,y>)\\
&=&4<x,y>
\end{eqnarray*}

\index{polarization identity}

Now given a matrix $U\in M_N(\mathbb C)$, we have the following equivalences, with the first one coming from the above identities, and with the other ones being clear:
\begin{eqnarray*}
||Ux||=||x||
&\iff&<Ux,Uy>=<x,y>\\
&\iff&<x,U^*Uy>=<x,y>\\
&\iff&U^*Uy=y\\
&\iff&U^*U=1\\
&\iff&U^*=U^{-1}
\end{eqnarray*}

(3) As in the real case, $P$ is an abstract projection, not necessarily orthogonal, when $P^2=P$. The point now is that this projection is orthogonal when:
\begin{eqnarray*}
<Px-Py,Px-x>=0
&\iff&<x-y,P^*Px-P^*x>=0\\
&\iff&P^*Px-P^*x=0\\
&\iff&P^*P-P^*=0
\end{eqnarray*}

Thus we must have $P^*=P^*P$. Now observe that by conjugating, we obtain:
$$P
=(P^*P)^*
=P^*(P^*)^*
=P^*P$$

Now by comparing with the original relation, $P^*=P^*P$, we conclude that $P=P^*$. Thus, we have shown that any orthogonal projection must satisfy, as claimed:
$$P^2=P^*=P$$

Conversely, if this condition is satisfied, $P^2=P$ shows that $P$ is a projection, and $P=P^*$ shows via the above computation that $P$ is indeed orthogonal.

\medskip

(4) Once again in analogy with the real case, we have the following formula:
$$P_xy
=\frac{<y,x>}{<x,x>}\,x
=\frac{1}{||x||^2}<y,x>x$$

With this in hand, we can now compute the entries of $P_x$, as follows:
$$(P_x)_{ij}
=<P_xe_j,e_i>
=\frac{1}{||x||^2}<e_j,x><x,e_i>
=\frac{\bar{x}_jx_i}{||x||^2}$$

Thus, we are led to the formula in the statement.
\end{proof}

We can talk as well about eigenvalues and eigenvectors, as in the real case:

\index{eigenvalue}
\index{eigenvector}
\index{diagonalizable matrix}

\begin{definition}
Let $A\in M_N(\mathbb C)$ be a square matrix. When $Av=\lambda v$ we say that:
\begin{enumerate}
\item $v$ is an eigenvector of $A$.

\item $\lambda$ is an eigenvalue of $A$.
\end{enumerate}
We say that $A$ is diagonalizable when $\mathbb C^N$ has a basis of eigenvectors of $A$.
\end{definition}

When $A$ is diagonalizable, in that basis of eigenvectors we can write:
$$A=\begin{pmatrix}
\lambda_1\\
&\ddots\\
&&\lambda_N
\end{pmatrix}$$

In general, this means that we have a formula as follows, with $D$ diagonal:
$$A=PDP^{-1}$$

Indeed, we can take $P$ to be the matrix formed by the eigenvectors:
$$P=[v_1\ldots v_N]$$

As a first interesting result now, regarding the real matrices, we have:

\index{characteristic polynomial}

\begin{theorem}
The eigenvalues of a real matrix $A\in M_N(\mathbb R)$ are the roots of
$$P(x)=\det(A-x1_N)$$
and in particular, any such matrix $A\in M_N(\mathbb R)$ has at least $1$ complex eigenvalue.
\end{theorem}

\begin{proof}
The first assertion is something that we already know, coming from:
\begin{eqnarray*}
\exists v,Av=\lambda v
&\iff&\exists v,(A-\lambda 1_N)v=0\\
&\iff&\det(A-\lambda 1_N)=0
\end{eqnarray*}

As for the second assertion, this follows from the first assertion, and from Theorem 3.13, which shows in particular that $P$ has at least $1$ complex root.
\end{proof} 

It is possible to further build on these results, but this is quite long, and we will rather do this in the next chapter. For the moment, let us just keep in mind the conclusion that a real matrix $A\in M_N(\mathbb R)$ has substantially more chances of being diagonalizable over the complex numbers, than over the real numbers. As an illustration for this principle, and as a first concrete result, which is of true complex nature, we have:

\index{rotation}

\begin{theorem}
The rotation of angle $t\in\mathbb R$ in the real plane, namely
$$R_t=\begin{pmatrix}\cos t&-\sin t\\ \sin t&\cos t\end{pmatrix}$$
can be diagonalized over the complex numbers, as follows:
$$R_t=\frac{1}{2}\begin{pmatrix}1&1\\i&-i\end{pmatrix}
\begin{pmatrix}e^{-it}&0\\0&e^{it}\end{pmatrix}
\begin{pmatrix}1&-i\\1&i\end{pmatrix}$$
Over the real numbers this is impossible, unless $t=0,\pi$.
\end{theorem}

\begin{proof}
The last assertion is something clear, that we already know, coming from the fact that at $t\neq0,\pi$ our rotation is a ``true'' rotation, having no eigenvectors in the plane. Regarding the first assertion, the point is that we have the following computation:
\begin{eqnarray*}
R_t\binom{1}{i}
&=&\begin{pmatrix}\cos t&-\sin t\\ \sin t&\cos t\end{pmatrix}\binom{1}{i}\\
&=&\binom{\cos t-i\sin t}{i\cos t+\sin t}\\
&=&e^{-it}\binom{1}{i}
\end{eqnarray*}

We have as well a second eigenvector, as follows:
\begin{eqnarray*}
R_t\binom{1}{-i}
&=&\begin{pmatrix}\cos t&-\sin t\\ \sin t&\cos t\end{pmatrix}\binom{1}{-i}\\
&=&\binom{\cos t+i\sin t}{-i\cos t+\sin t}\\
&=&e^{it}\binom{1}{-i}
\end{eqnarray*}

Thus our matrix $R_t$ is diagonalizable over $\mathbb C$, with the diagonal form being:
$$R_t\sim\begin{pmatrix}e^{-it}&0\\0&e^{it}\end{pmatrix}$$

As for the passage matrix, obtained by putting together the eigenvectors, this is:
$$P=\begin{pmatrix}1&1\\i&-i\end{pmatrix}$$

In order to invert now $P$, we can use the standard inversion formula for the $2\times2$ complex matrices, which is similar to the one in the real case, and gives:
$$P^{-1}
=\frac{1}{-2i}\begin{pmatrix}-i&-1\\-i&1\end{pmatrix}
=\frac{1}{2}\begin{pmatrix}1&-i\\1&i\end{pmatrix}$$

Our diagonalization formula is therefore as follows:
$$R_t=\frac{1}{2}\begin{pmatrix}1&1\\i&-i\end{pmatrix}
\begin{pmatrix}e^{-it}&0\\0&e^{it}\end{pmatrix}
\begin{pmatrix}1&-i\\1&i\end{pmatrix}$$

Thus, we are led to the conclusion in the statement.
\end{proof}

\section*{3d. The determinant}

Regarding now the determinant, for the complex matrices it is more convenient to use an abstract approach, and this due to our lack of geometric intuition with the space $\mathbb C^N$, at $N\geq2$, and with the ``complex volumes'' of the bodies there. So, let us formulate:

\index{determinant}
\index{determinant formula}

\begin{definition}
The determinant of a complex matrix $A\in M_N(\mathbb C)$ is given by
$$\det A=\sum_{\sigma\in S_N}\varepsilon(\sigma)A_{1\sigma(1)}\ldots A_{N\sigma(N)}$$
with $\varepsilon=\pm1$ being the signature of the permutations.
\end{definition}

Generally speaking, the theory of the determinant from the real case extends well. To be more precise, we first have the following result, summarizing the basics:

\begin{theorem}
The determinant has the following properties:
\begin{enumerate}
\item When adding two columns, the determinants get added:
$$\det(\ldots,u+v,\ldots)
=\det(\ldots,u,\ldots)
+\det(\ldots,v,\ldots)$$

\item When multiplying columns by scalars, the determinant gets multiplied:
$$\det(\lambda v_1,\ldots,\lambda_Nv_N)=\lambda_1\ldots\lambda_N\det(v_1,\ldots,v_N)$$

\item When permuting two columns, the determinant changes the sign:
$$\det(\ldots,v,\ldots,w,\ldots)=-\det(\ldots,w,\ldots,v,\ldots)$$

\item Also, the determinant of the identity matrix is $1$. 
\end{enumerate}
\end{theorem}

\begin{proof}
This follows indeed by doing some elementary algebraic computations with permutations, which are similar to those that we did before in the real case, but done now backwards, based on the formula of the determinant from Definition 3.23.
\end{proof}

We have as well a similar result for the rows, which is equally useful, as follows:

\begin{theorem}
The determinant has the following properties:
\begin{enumerate}
\item When adding two rows, the determinants get added:
$$\det\begin{pmatrix}\vdots\\ u+v\\ \vdots\end{pmatrix}
=\det\begin{pmatrix}\vdots\\ u\\ \vdots\end{pmatrix}
+\det\begin{pmatrix}\vdots\\ v\\ \vdots\end{pmatrix}$$

\item When multiplying rows by scalars, the determinant gets multiplied:
$$\det\begin{pmatrix}\lambda_1v_1\\ \vdots\\ \lambda_Nv_N\end{pmatrix}
=\lambda_1\ldots\lambda_N\det\begin{pmatrix}v_1\\ \vdots\\ v_N\end{pmatrix}$$

\item When permuting two rows, the determinant changes the sign.
\end{enumerate}
\end{theorem}

\begin{proof}
This follows once again by doing some algebraic computations with permutations, based on the formula of the determinant from Definition 3.23.
\end{proof}

Next in line, we have the following result, which is very useful in practice:

\index{row expansion}

\begin{theorem}
The determinant is subject to the row expansion formula
\begin{eqnarray*}
\begin{vmatrix}a_{11}&\ldots&a_{1N}\\
\vdots&&\vdots\\
a_{N1}&\ldots&a_{NN}\end{vmatrix}
&=&a_{11}\begin{vmatrix}a_{22}&\ldots&a_{2N}\\
\vdots&&\vdots\\
a_{N2}&\ldots&a_{NN}\end{vmatrix}\\
&-&a_{12}\begin{vmatrix}a_{21}&a_{23}&\ldots&a_{2N}\\
\vdots&\vdots&&\vdots\\
a_{N1}&a_{N3}&\ldots&a_{NN}\end{vmatrix}\\
&\vdots&\\
&\vdots&\\
&+&(-1)^{N+1}a_{1N}\begin{vmatrix}a_{21}&\ldots&a_{2,N-1}\\
\vdots&&\vdots\\
a_{N1}&\ldots&a_{N,N-1}\end{vmatrix}
\end{eqnarray*}
and this method fully computes it, by recurrence.
\end{theorem}

\begin{proof}
This follows indeed by doing some elementary algebraic computations.
\end{proof}

We can expand as well over the columns, as follows:

\index{column expansion}

\begin{theorem}
The determinant is subject to the column expansion formula
\begin{eqnarray*}
\begin{vmatrix}a_{11}&\ldots&a_{1N}\\
\vdots&&\vdots\\
a_{N1}&\ldots&a_{NN}\end{vmatrix}
&=&a_{11}\begin{vmatrix}a_{22}&\ldots&a_{2N}\\
\vdots&&\vdots\\
a_{N2}&\ldots&a_{NN}\end{vmatrix}\\
&-&a_{21}\begin{vmatrix}a_{12}&\ldots&a_{1N}\\
a_{32}&\ldots&a_{3N}\\
\vdots&&\vdots\\
a_{N2}&\ldots&a_{NN}\end{vmatrix}\\
&\vdots&\\
&\vdots&\\
&+&(-1)^{N+1}a_{N1}\begin{vmatrix}a_{12}&\ldots&a_{1N}\\
\vdots&&\vdots\\
a_{N-1,2}&\ldots&a_{N-1,N}\end{vmatrix}
\end{eqnarray*}
and this method fully computes it, by recurrence.
\end{theorem}

\begin{proof}
Once again, this follows by doing some algebraic computations.
\end{proof}

Still in analogy with the real case, we have the following result:

\index{multilinear form}

\begin{theorem}
The determinant of the systems of vectors
$$\det:\mathbb C^N\times\ldots\times\mathbb C^N\to\mathbb C$$
is multilinear, alternate and unital, and unique with these properties.
\end{theorem}

\begin{proof}
This is something that we know in the real case, and the proof in the complex case is similar, with the conditions in the statement corresponding to those in Theorem 3.24. It is possible to prove this result as well directly, by doing some abstract algebra.
\end{proof}

Finally, once again at the general level, let us record the following result:

\begin{theorem}
We have the following formulae,
$$\det\bar{A}=\overline{\det A}\quad,\quad 
\det A^t=\det A\quad,\quad 
\det A^*=\overline{\det A}$$
valid for any square matrix $A\in M_N(\mathbb C)$.
\end{theorem}

\begin{proof}
The first formula is clear from Definition 3.23, because when conjugating the entries of $A$, the determinant will get conjugated:
$$\det\bar{A}=\sum_{\sigma\in S_N}\varepsilon(\sigma)\overline{A_{1\sigma(1)}\ldots A_{N\sigma(N)}}$$

The second formula follows as in the real case, as follows:
\begin{eqnarray*}
\det A^t
&=&\sum_{\sigma\in S_N}\varepsilon(\sigma)(A^t)_{1\sigma(1)}\ldots(A^t)_{N\sigma(N)}\\
&=&\sum_{\sigma\in S_N}\varepsilon(\sigma)A_{\sigma(1)1}\ldots A_{\sigma(N)N}\\
&=&\sum_{\sigma\in S_N}\varepsilon(\sigma)A_{1\sigma^{-1}(1)}\ldots A_{N\sigma^{-1}(N)}\\
&=&\sum_{\sigma\in S_N}\varepsilon(\sigma^{-1})A_{1\sigma^{-1}(1)}\ldots A_{N\sigma^{-1}(N)}\\
&=&\sum_{\sigma\in S_N}\varepsilon(\sigma)A_{1\sigma(1)}\ldots A_{N\sigma(N)}\\
&=&\det A
\end{eqnarray*}

As for the third formula, this follows from the first two formulae, by using:
$$\det A^*=\det \bar{A}^t$$

Thus, we are led to the conclusions in the statement.
\end{proof}

Summarizing, the theory from the real case extends well, and we have complex analogues of all results. As in the real case, as a main application of all this, we have:

\index{matrix inversion}

\begin{theorem}
The inverse of a square matrix, having nonzero determinant,
$$A=\begin{pmatrix}a_{11}&\ldots&a_{1N}\\
\vdots&&\vdots\\
a_{N1}&\ldots&a_{NN}\end{pmatrix}$$
is given by the following formula,
$$A^{-1}=\frac{1}{\det A}
\begin{pmatrix}
\det A^{(11)}&-\det A^{(21)}&\det A^{(31)}&\ldots\\
-\det A^{(12)}&\det A^{(22)}&-\det A^{(32)}&\ldots\\
\det A^{(13)}&-\det A^{(23)}&\det A^{(33)}&\ldots\\
\vdots&\vdots&\vdots&
\end{pmatrix}$$
where $A^{(ij)}$ is the matrix $A$, with the $i$-th row and $j$-th column removed.
\end{theorem}

\begin{proof}
This follows indeed by using the row expansion formula from Theorem 3.26, which in terms of the matrix $A^{-1}$ in the statement reads $AA^{-1}=1$.
\end{proof}

As a final topic now, regarding the complex matrices, let us discuss some interesting examples of such matrices, which definitely do not exist in the real setting, and which are very useful, even in connection with real matrix questions. Let us start with:

\begin{definition}
The Fourier matrix is as follows, 
$$F_N=(w^{ij})_{ij}$$
with $w=e^{2\pi i/N}$, and with the convention that the indices are
$$i,j\in\{0,1,\ldots,N-1\}$$ 
and are taken modulo $N$.
\end{definition}

Here the conventions regarding the indices are standard, and are there for various reasons, as for instance for having the first row and column consisting of 1 entries. Indeed, in standard matrix form, and with the above conventions for the indices, we have:
$$F_N=
\begin{pmatrix}
1&1&1&\ldots&1\\
1&w&w^2&\ldots&w^{N-1}\\
1&w^2&w^4&\ldots&w^{2(N-1)}\\
\vdots&\vdots&\vdots&&\vdots\\
1&w^{N-1}&w^{2(N-1)}&\ldots&w^{(N-1)^2}
\end{pmatrix}$$

Thus, what we have here is a Vandermonde matrix, in the sense of chapter 2, of very special type. Let us record as well the first few values of these matrices:

\index{Fourier matrix}

\begin{proposition}
The second Fourier matrix is as follows,
$$F_2=\begin{pmatrix}1&1\\1&-1\end{pmatrix}$$
the third Fourier matrix is as follows, with $w=e^{2\pi i/3}$,
$$F_3=\begin{pmatrix}1&1&1\\1&w&w^2\\1&w^2&w\end{pmatrix}$$
and the fourth Fourier matrix is as follows, with $i^2=-1$ as usual,
$$F_4=\begin{pmatrix}1&1&1&1\\1&i&-1&-i\\1&-1&1&-1\\1&-i&-1&i\end{pmatrix}$$
with the above convention for the indices, $i,j\in\{0,1,\ldots,N-1\}$.
\end{proposition}

\begin{proof}
All these formulae are clear from definitions, with our usual convention for the indices of the Fourier matrices, from Definition 3.31.
\end{proof}

Our claim now is that the Fourier matrix can be used in order to solve a variety of linear algebra questions, a bit in a same way as the Fourier transform can be used in order to solve analysis questions. Before discussing all this, however, let us analyze the Fourier matrix $F_N$, from a linear algebra perspective. We have the following result:

\begin{theorem}
The Fourier matrix $F_N$ has the following properties:
\begin{enumerate}
\item It is symmetric, $F_N^t=F_N$.

\item The matrix $F_N/\sqrt{N}$ is unitary.

\item Its inverse is the matrix $F_N^*/N$.
\end{enumerate}
\end{theorem}

\begin{proof}
This is a collection of elementary results, the idea being as follows:

\medskip

(1) This is indeed clear from definitions.

\medskip

(2) The row vectors $R_0,\ldots,R_{N-1}$ of the rescaled matrix $F_N/\sqrt{N}$ have all length 1, and by using the barycenter formula in Theorem 3.16, we have, for any $i\neq j$:
$$<R_i,R_j>
=\frac{1}{N}\sum_kw^{ik}w^{-jk}
=\frac{1}{N}\sum_k(w^{i-j})^k
=0$$

Thus, $R_0,\ldots,R_{N-1}$ are pairwise orthogonal, and so $F_N/\sqrt{N}$ is unitary, as claimed.

\medskip

(3) This follows from (1) and (2), because for a symmetric matrix, the adjoint is the conjugate, and in the unitary case, this is the inverse.
\end{proof}

Now back to our motivations, we were saying before that the Fourier matrix is to linear algebra what the Fourier transform is to analysis, namely advanced technology. In order to discuss now an illustrating application of the theory developed above, let us go back to our favorite example of a $N\times N$ matrix, namely the flat matrix:
$$\mathbb I_N=\begin{pmatrix}
1&\ldots&1\\
\vdots&&\vdots\\
1&\ldots&1\end{pmatrix}$$

This is a real matrix, and we know that we have $\mathbb I_N=NP_N$, with $P_N$ being the projection on the all-1 vector $\xi=(1)_i\in\mathbb R^N$. Thus, $\mathbb I_N$ diagonalizes over $\mathbb R$:
$$\mathbb I_N\sim
\begin{pmatrix}
N\\
&0\\
&&\ddots\\
&&&0
\end{pmatrix}$$

The problem, however, is that when looking for 0-eigenvectors, in order to have an explicit diagonalization formula, we must solve the following equation:
$$x_1+\ldots+x_N=0$$

And this is not an easy task, if our objective is that of finding a nice, explicit basis for the space of solutions. To be more precise, if we want linearly independent vectors $v_1,\ldots,v_{N-1}\in\mathbb R^N$, each with components summing up to 0, and which are given by simple formulae, of type $(v_i)_j=$ explicit function of $i,j$, we are in trouble.

\bigskip

Fortunately, the complex numbers come to the rescue, and we have:

\index{flat matrix}
\index{Fourier matrix}

\begin{theorem}
The flat matrix of size $N$, namely
$$\mathbb I_N=\begin{pmatrix}
1&\ldots&1\\
\vdots&&\vdots\\
1&\ldots&1\end{pmatrix}$$
has the following explicit diagonalization, over the complex numbers,
$$\mathbb I_N=\frac{1}{N}F_NQF_N^*$$
with $F_N=(w^{ij})_{ij}$ being the Fourier matrix, and with $Q=diag(N,0,\ldots,0)$.
\end{theorem}

\begin{proof}
Indeed, the 0-eigenvector problem discussed above can be solved explicitly over the complex numbers, by using the formula in Theorem 3.16, with the solution $(v_i)_j=w^{ij}$, with $w=e^{2\pi i/N}$. Thus, we are led to the conclusion in the statement.
\end{proof}

There are many other uses of the Fourier matrix $F_N$, along the same lines. We will be back to all this in chapter 7 below, with a complete discussion of the Fourier matrices, and of their natural generalizations, called complex Hadamard matrices.

\section*{3e. Exercises}

As a first exercise, in relation with the complex numbers, we have:

\begin{exercise}
Try to use a complex number type idea in order to multiply the vectors of $\mathbb R^3$, and then $\mathbb R^4$, and report on what you found.
\end{exercise}

This is something quite tricky, and a piece of hint, do not worry if you find nothing interesting at $N=3$. However, the $N=4$ case is definitely worth some study.

\begin{exercise}
Can you use complex numbers in order to explicitly find the roots of arbitrary degree $3$ polynomials, a bit in the same way as in degree $2$?
\end{exercise}

This is actually something quite tricky, and if stuck, look up on the internet, or in a good calculus book of your choice, ``Cardano formula'', which is the keyword for this.

\begin{exercise}
Write down a complete proof for the Euler formula
$$e^{it}=\cos t+i\sin t$$
using any method of your choice. 
\end{exercise}

This is something that we discussed in the above, but with our proofs however still missing a few details, regarding the basic properties of the function $e^x$. Thus, you can either try to recover these details, or go with some other idea, of your choice.

\begin{exercise}
Find a geometric interpretation of the formula
$$\begin{pmatrix}\cos t&-\sin t\\ \sin t&\cos t\end{pmatrix}
=\frac{1}{2}\begin{pmatrix}1&1\\i&-i\end{pmatrix}
\begin{pmatrix}e^{-it}&0\\0&e^{it}\end{pmatrix}
\begin{pmatrix}1&-i\\1&i\end{pmatrix}$$
which diagonalizes the rotation of angle $t\in\mathbb R$ in the real plane.
\end{exercise}

This is something quite tricky, and of course, enjoy.

\begin{exercise}
Develop a complete theory of diagonalization for the $2\times2$ matrices, notably by deciding when exactly such a matrix is diagonalizable.
\end{exercise}

This is quite non-trivial, but all the needed ingredients are in the above.

\begin{exercise}
Work out all the details of the diagonalization formula
$$\mathbb I_N=\frac{1}{N}F_NQF_N^*$$
with $Q=diag(N,0,\ldots,0)$, and then try formulating a generalization of this.
\end{exercise}

Here the first question is standard, amounting in completing the proof that was given in the above. As for the second question, this is something more tricky.

\chapter{Diagonalization}

\section*{4a. Diagonalization}

In this chapter we discuss the diagonalization question, with a number of advanced results, for the complex matrices $A\in M_N(\mathbb C)$. Our techniques will apply of course to the real case too, $A\in M_N(\mathbb R)$, and we will obtain in this way a number of non-trivial results regarding the diagonalization of such matrices, over the complex numbers.

\bigskip

Let us begin with a reminder of the basic diagonalization theory, that we already know. The basic theory that we have so far can be summarized as follows:

\index{diagonalizable matrix}

\begin{theorem}
Assuming that a matrix $A\in M_N(\mathbb C)$ is diagonalizable, in the sense that $\mathbb C^N$ has a basis formed by eigenvectors of $A$, we have
$$A=PDP^{-1}$$
where $P=[v_1\ldots v_N]$ is the square matrix formed by the eigenvectors of $A$, and $D=diag(\lambda_1,\ldots,\lambda_N)$ is the diagonal matrix formed by the corresponding eigenvalues.
\end{theorem}

\begin{proof}
This is something that we already know, coming by changing the basis. We can prove this by direct computation as well, because we have $Pe_i=v_i$, and so the matrices $A$ and $PDP^{-1}$ follow to act in the same way on the basis vectors $v_i$:
\begin{eqnarray*}
PDP^{-1}v_i
&=&PDe_i\\
&=&P\lambda_ie_i\\
&=&\lambda_iPe_i\\
&=&\lambda_iv_i
\end{eqnarray*}

Thus, the matrices $A$ and $PDP^{-1}$ coincide, as stated.
\end{proof}

In general, in order to study the diagonalization problem, the idea is that the eigenvectors can be grouped into linear spaces, called eigenspaces:

\index{eigenspaces}

\begin{definition}
Given $A\in M_N(\mathbb C)$, for any eigenvalue $\lambda\in\mathbb C$ we let
$$E_\lambda=\left\{v\in\mathbb C^N\Big|Av=\lambda v\right\}$$
be the vector space formed by the corresponding eigenvectors.
\end{definition}

As an illustration for this, consider a diagonalizable matrix $A\in M_N(\mathbb C)$, with the diagonalization chosen as for the eigenvalues to appear grouped, as follows:
$$A\sim\begin{pmatrix}
\lambda_1\\
&\ddots\\
&&\lambda_1\\
&&&\ddots\\
&&&&\lambda_k\\
&&&&&\ddots\\
&&&&&&\lambda_k
\end{pmatrix}$$

The corresponding eigenspaces are then as follows, in an obvious direct sum position, with $d_1,\ldots,d_k$ being the multiplicities of the eigenvalues $\lambda_1,\ldots,\lambda_k$:
$$E_{\lambda_1}=\left\{\begin{pmatrix}
x_1\\
\vdots\\
x_{d_1}\\
0\\
\vdots\\
0\\
0\\
\vdots\\
0
\end{pmatrix},x_i\in\mathbb C\right\}
\qquad\ldots\qquad
E_{\lambda_k}=\left\{\begin{pmatrix}
0\\
\vdots\\
0\\
0\\
\vdots\\
0\\
x_1\\
\vdots\\
x_{d_k}
\end{pmatrix},x_i\in\mathbb C\right\}$$

In general, we have the following result, which is valid for any matrix:

\index{dimension inequality}

\begin{theorem}
The eigenspaces of an arbitrary matrix $A\in M_N(\mathbb C)$, given by
$$E_\lambda=\left\{v\in\mathbb C^N\Big|Av=\lambda v\right\}$$
are in a direct sum position, in the sense that given vectors $v_1\in E_{\lambda_1},\ldots,v_k\in E_{\lambda_k}$ corresponding to different eigenvalues $\lambda_1,\ldots,\lambda_k$, we have:
$$\sum_ic_iv_i=0\implies c_i=0$$
In particular, we have the following dimension inequality, with the sum being over all the eigenvalues $\lambda\in\mathbb C$ of our matrix $A$,
$$\sum_\lambda\dim(E_\lambda)\leq N$$
and our matrix is diagonalizable precisely when we have equality.
\end{theorem}

\begin{proof}
We prove the first assertion by recurrence on $k\in\mathbb N$. Assume by contradiction that we have a formula as follows, with the scalars $c_1,\ldots,c_k$ being not all zero:
$$c_1v_1+\ldots+c_kv_k=0$$

By dividing by one of these scalars, we can assume that our formula is:
$$v_k=c_1v_1+\ldots+c_{k-1}v_{k-1}$$

Now let us apply $A$ to this vector. On the left we obtain:
$$Av_k
=\lambda_kc_1v_1+\ldots+\lambda_kc_{k-1}v_{k-1}$$

On the right we obtain something different, as follows:
\begin{eqnarray*}
A(c_1v_1+\ldots+c_{k-1}v_{k-1})
&=&c_1Av_1+\ldots+c_{k-1}Av_{k-1}\\
&=&c_1\lambda_1v_1+\ldots+c_{k-1}\lambda_{k-1}v_{k-1}
\end{eqnarray*}

We conclude from this that the following equality must hold:
$$\lambda_kc_1v_1+\ldots+\lambda_kc_{k-1}v_{k-1}=c_1\lambda_1v_1+\ldots+c_{k-1}\lambda_{k-1}v_{k-1}$$

On the other hand, we know by recurrence that the vectors $v_1,\ldots,v_{k-1}$ must be linearly independent. Thus, the coefficients must be equal, at right and at left:
$$\lambda_kc_1=c_1\lambda_1$$
$$\vdots$$
$$\lambda_kc_{k-1}=c_{k-1}\lambda_{k-1}$$

Now since at least one $c_i$ must be nonzero, from the corresponding equality $\lambda_kc_i=c_i\lambda_i$ we obtain $\lambda_k=\lambda_i$, which is a contradiction. Thus our proof by recurrence of the first assertion is complete. As for the second assertion, this follows from the first one.
\end{proof}

The above result is something quite intuitive, and in the case of a diagonalizable matrix, this comes from the discussion before the statement. As a second illustration, let us see as well what happens for the simplest non-diagonalizable matrix, namely:
$$J=\begin{pmatrix}0&1\\0&0\end{pmatrix}$$

As observed in chapter 1, this matrix has $\lambda=0$ as unique eigenvalue, with:
$$E_0=\left\{\binom{x}{0},x\in\mathbb C\right\}$$

Thus, the diagonalization condition in Theorem 4.3 is not satisfied indeed.

\bigskip

In order to reach now to more advanced results, we can use the characteristic polynomial. Here is a result summarizing and improving our knowledge of the subject:

\index{characteristic polynomial}

\begin{theorem}
Given a matrix $A\in M_N(\mathbb C)$, consider its characteristic polynomial:
$$P(x)=\det(A-x1_N)$$
The eigenvalues of $A$ are then the roots of $P$. Also, we have the inequality
$$\dim(E_\lambda)\leq m_\lambda$$
where $m_\lambda$ is the multiplicity of $\lambda$, as root of $P$.
\end{theorem}

\begin{proof}
The first assertion follows from the following computation, using the fact that a linear map is bijective when the determinant of the associated matrix is nonzero:
\begin{eqnarray*}
\exists v,Av=\lambda v
&\iff&\exists v,(A-\lambda 1_N)v=0\\
&\iff&\det(A-\lambda 1_N)=0
\end{eqnarray*}

Regarding now the second assertion, given an eigenvalue $\lambda$ of our matrix $A$, consider the dimension of the corresponding eigenspace:
$$d_\lambda=\dim(E_\lambda)$$

By changing the basis of $\mathbb C^N$, as for the eigenspace $E_\lambda$ to be spanned by the first $d_\lambda$ basis elements, our matrix becomes as follows, with $B$ being a certain smaller matrix:
$$A\sim\begin{pmatrix}\lambda 1_{d_\lambda}&0\\0&B\end{pmatrix}$$

We conclude that the characteristic polynomial of $A$ is of the following form:
$$P_A
=P_{\lambda 1_{d_\lambda}}P_B
=(\lambda-x)^{d_\lambda}P_B$$

Thus we have $m_\lambda\geq d_\lambda$, which leads to the conclusion in the statement.
\end{proof}

We can put together Theorem 4.3 and Theorem 4.4, and by using as well the fact that any complex polynomial of degree $N$ has exactly $N$ complex roots, when counted with multiplicities, that we know from chapter 3, we obtain the following result:

\index{diagonalization}
\index{characteristic polynomial}
\index{eigenspaces}

\begin{theorem}
Given a matrix $A\in M_N(\mathbb C)$, consider its characteristic polynomial
$$P(X)=\det(A-X1_N)$$ 
then factorize this polynomial, by computing the complex roots, with multiplicities,
$$P(X)=(-1)^N(X-\lambda_1)^{n_1}\ldots(X-\lambda_k)^{n_k}$$
and finally compute the corresponding eigenspaces, for each eigenvalue found:
$$E_i=\left\{v\in\mathbb C^N\Big|Av=\lambda_iv\right\}$$
The dimensions of these eigenspaces satisfy then the following inequalities,
$$\dim(E_i)\leq n_i$$
and $A$ is diagonalizable precisely when we have equality for any $i$.
\end{theorem}

\begin{proof}
This follows by combining the above results. Indeed, by summing the inequalities $\dim(E_\lambda)\leq m_\lambda$ from Theorem 4.4, we obtain an inequality as follows:
$$\sum_\lambda\dim(E_\lambda)\leq\sum_\lambda m_\lambda\leq N$$

On the other hand, we know from Theorem 4.3 that our matrix is diagonalizable when we have global equality. Thus, we are led to the conclusion in the statement.
\end{proof}

This was for the main result of linear algebra. There are countless applications of this, and generally speaking, advanced linear algebra consists in further building on Theorem 4.5. Let us record as well a useful algorithmic version of the above result:

\begin{theorem}
The square matrices $A\in M_N(\mathbb C)$ can be diagonalized as follows:
\begin{enumerate}
\item Compute the characteristic polynomial.

\item Factorize the characteristic polynomial.

\item Compute the eigenvectors, for each eigenvalue found.

\item If there are no $N$ eigenvectors, $A$ is not diagonalizable.

\item Otherwise, $A$ is diagonalizable, $A=PDP^{-1}$.
\end{enumerate}
\end{theorem}

\begin{proof}
This is an informal reformulation of Theorem 4.5, with (4) referring to the total number of linearly independent eigenvectors found in (3), and with $A=PDP^{-1}$ in (5) being the usual diagonalization formula, with $P,D$ being as before.
\end{proof}

As a remark here, in step (3) it is always better to start with the eigenvalues having big multiplicity. Indeed, a multiplicity 1 eigenvalue, for instance, can never lead to the end of the computation, via (4), simply because the eigenvectors always exist.

\bigskip

As a key consequence of Theorem 4.5, which is very useful in practice, we have:

\begin{theorem}
If a matrix $A\in M_N(\mathbb C)$ has distinct eigenvalues, then it is diagonalizable. Moreover, this is indeed the case, for the generic matrices.
\end{theorem}

\begin{proof}
The first assertion is clear from Theorem 4.3, because the criterion there for diagonalization is trivially satisfied when the eigenvalues are different, as follows:
$$\sum_\lambda\dim(E_\lambda)
=\sum_\lambda1
=N$$

As for the second assertion, this is something quite intuitive, coming from the fact that $N$ numbers $\lambda_1,\ldots,\lambda_N\in\mathbb C$ picked at random must be distinct. Of course, this does not stand as a formal proof, but we will come back to this in a moment, with a proof.
\end{proof}

Getting back now to Theorem 4.5, or rather to Theorem 4.6, the main problem raised by the diagonalization procedure is the computation of the roots of characteristic polynomials. As a first observation here, in degree 2 we have the following trick:

\begin{proposition}
The roots of a degree $2$ polynomial of the form
$$P=X^2-aX+b$$
are precisely the numbers $r,s$ satisfying $r+s=a$, $rs=b$.
\end{proposition}

\begin{proof}
This is indeed something trivial, coming from $P=(X-r)(X-s)$.
\end{proof}

In the matrix setting now, the result coming from this is as follows:

\begin{theorem}
Consider an arbitrary $2\times2$ matrix, written as follows:
$$A=\begin{pmatrix}
a&b\\
c&d\end{pmatrix}$$
\begin{enumerate}
\item The eigenvalues are the numbers $r,s$ satisfying $r+s=a+d$, $rs=ad-bc$.

\item When $(a+d)^2\neq4(ad-bc)$ we have $r\neq s$, and $A$ is diagonalizable.

\item Otherwise, $r=s$, and $A$ is diagonalizable precisely when $A=\binom{r\ 0}{0\ r}$.
\end{enumerate}
\end{theorem}

\begin{proof}
This is something straightforward, coming from Proposition 4.8:

\medskip

(1) We have indeed the following computation, which gives the result:
$$\det(A-X1_2)=\begin{vmatrix}a-X&b\\ c&d-X\end{vmatrix}=X^2-(a+d)X+(ad-bc)$$

(2) Here the first assertion comes from $\Delta=(a+d)^2-4(ad-bc)$ for the degree 2 polynomial found above, and the second assertion comes from Theorem 4.7.

\medskip

(3) Assuming $\Delta=0$ we have indeed $r=s$, and then, according to Theorem 4.5, the diagonalization condition reads $E_r=\mathbb C^2$, so $Ax=rx$ for any $x$, and so $A=\binom{r\ 0}{0\ r}$.
\end{proof}

In higher dimensions things certainly get more complicated, but we have:

\index{sum of eigenvalues}
\index{product of eigenvalues}

\begin{theorem}
The complex eigenvalues of a matrix $A\in M_N(\mathbb C)$, counted with multiplicities, have the following properties:
\begin{enumerate}
\item Their sum is the trace.

\item Their product is the determinant.
\end{enumerate}
\end{theorem}

\begin{proof}
Consider indeed the characteristic polynomial $P$ of the matrix:
\begin{eqnarray*}
P(X)
&=&\det(A-X1_N)\\
&=&(-1)^NX^N+(-1)^{N-1}Tr(A)X^{N-1}+\ldots+\det(A)
\end{eqnarray*}

We can factorize this polynomial, by using its $N$ complex roots, and we obtain:
\begin{eqnarray*}
P(X)
&=&(-1)^N(X-\lambda_1)\ldots(X-\lambda_N)\\
&=&(-1)^NX^N+(-1)^{N-1}\left(\sum_i\lambda_i\right)X^{N-1}+\ldots+\prod_i\lambda_i
\end{eqnarray*}

Thus, we are led to the conclusion in the statement.
\end{proof}

Regarding now the intermediate terms, we have here the following result:

\index{symmetric functions}
\index{characteristic polynomial}

\begin{theorem}
Assume that $A\in M_N(\mathbb C)$ has eigenvalues $\lambda_1,\ldots,\lambda_N\in\mathbb C$, counted with multiplicities. The basic symmetric functions of these eigenvalues, namely
$$c_k=\sum_{i_1<\ldots<i_k}\lambda_{i_1}\ldots\lambda_{i_k}$$
are then given by the fact that the characteristic polynomial of the matrix is:
$$P(X)=(-1)^N\sum_{k=0}^N(-1)^kc_kX^k$$
Moreover, all symmetric functions of the eigenvalues, such as the sums of powers
$$d_s=\lambda_1^s+\ldots+\lambda_N^s$$
appear as polynomials in these characteristic polynomial coefficients $c_k$.
\end{theorem}

\begin{proof}
These results can be proved by doing some algebra, as follows:

\medskip

(1) Consider indeed the characteristic polynomial $P$ of the matrix, factorized by using its $N$ complex roots, taken with multiplicities. By expanding, we obtain:
\begin{eqnarray*}
P(X)
&=&(-1)^N(X-\lambda_1)\ldots(X-\lambda_N)\\
&=&(-1)^NX^N+(-1)^{N-1}\left(\sum_i\lambda_i\right)X^{N-1}+\ldots+\prod_i\lambda_i\\
&=&(-1)^NX^N+(-1)^{N-1}c_1X^{N-1}+\ldots+(-1)^0c_N\\
&=&(-1)^N\left(X^N-c_1X^{N-1}+\ldots+(-1)^Nc_N\right)
\end{eqnarray*}

With the convention $c_0=1$, we are led to the conclusion in the statement.

\medskip

(2) This is something standard, coming by doing some abstract algebra. Working out the formulae for the sums of powers $d_s=\sum_i\lambda_i^s$, at small values of the exponent $s\in\mathbb N$, is an excellent exercise, which shows how to proceed in general, by recurrence.
\end{proof}

Finally, getting back to the general factorization question for polynomials, we have the following result, which can be useful as well, in the linear algebra context:

\begin{theorem}
Assume that we have a polynomial as follows, with integer coefficients, and with the leading term being $1$:
$$P=X^N+a_{N-1}X^{N-1}+\ldots+a_1X+a_0$$
The integer roots of $P$ must then divide the last coefficient $a_0$.
\end{theorem}

\begin{proof}
This is clear, because any integer root $c\in\mathbb Z$ of our polynomial must satisfy:
$$c^N+a_{N-1}c^{N-1}+\ldots+a_1c+a_0=0$$

But modulo $c$, this equation simply reads $a_0=0$, as desired.
\end{proof}

\section*{4b. Density tricks}

Let us go back now to Theorem 4.7, and more specifically, to the last assertion there, which was quite a strong statement. In order to discuss this, we first have:

\index{simple roots}

\begin{theorem}
For a matrix $A\in M_N(\mathbb C)$ the following conditions are equivalent,
\begin{enumerate}
\item The eigenvalues are different, $\lambda_i\neq\lambda_j$,

\item The characteristic polynomial $P$ has simple roots,

\item The characteristic polynomial satisfies $(P,P')=1$,
\end{enumerate}
and in this case, the matrix is diagonalizable.
\end{theorem}

\begin{proof}
The equivalences in the statement are clear, the idea being as follows:

\medskip

$(1)\iff(2)$ This follows indeed from Theorem 4.5.

\medskip

$(2)\iff(3)$ This is standard, the double roots of $P$ being roots of $P'$.

\medskip

As for the last assertion, this is something that we know, from Theorem 4.7.
\end{proof}

As an important comment, the assumptions of Theorem 4.13 can be effectively verified in practice, without the need for factorizing polynomials, the idea here being that of using the condition (3) there. In order to discuss this, let us start with:

\index{resultant}
\index{common roots}

\begin{theorem}
Given two polynomials $P,Q\in\mathbb C[X]$, written as follows,
$$P=c(X-a_1)\ldots(X-a_k)\quad,\quad 
Q=d(X-b_1)\ldots(X-b_l)$$
the following quantity, which is called resultant of $P,Q$,
$$R(P,Q)=c^ld^k\prod_{ij}(a_i-b_j)$$
is a polynomial in the coefficients of $P,Q$, with integer coefficients, and we have
$$R(P,Q)=0$$
precisely when $P,Q$ have a common root.
\end{theorem}

\begin{proof}
This is something quite tricky, the idea being as follows:

\medskip

(1) Given two polynomials $P,Q\in\mathbb C[X]$, we can certainly construct the quantity $R(P,Q)$ in the statement, with the role of the normalization factor $c^ld^k$ to become clear later on, and then we have $R(P,Q)=0$ precisely when $P,Q$ have a common root:
$$R(P,Q)=0\iff \exists i,j,a_i=b_j$$

(2) As bad news, however, this quantity $R(P,Q)$, defined in this way, is a priori not very useful in practice, because it depends on the roots $a_i,b_j$ of our polynomials $P,Q$, that we cannot compute in general. However, and here comes our point, as we will prove below, it turns out that $R(P,Q)$ is in fact a polynomial in the coefficients of $P,Q$, with integer coefficients, and this is where the power of $R(P,Q)$ comes from.

\medskip

(3) You might perhaps say, nice, but why not doing things the other way around, that is, formulating our theorem with the explicit formula of $R(P,Q)$, in terms of the coefficients of $P,Q$, and then proving that we have $R(P,Q)=0$, via roots and everything. Good point, but this is not exactly obvious, the formula of $R(P,Q)$ in terms of the coefficients of $P,Q$ being something quite complicated. In short, trust me, let us prove our theorem as stated, and for alternative formulae of $R(P,Q)$, we will see later. 

\medskip

(4) Getting started now, let us expand the formula of $R(P,Q)$, by making all the multiplications there, abstractly, in our head. Everything being symmetric in $a_1,\ldots,a_k$, we obtain in this way certain symmetric functions in these variables, which will be therefore certain polynomials in the coefficients of $P$. Moreover, due to our normalization factor $c^l$, these polynomials in the coefficients of $P$ will have integer coefficients.

\medskip

(5) With this done, let us look now what happens with respect to the remaining variables $b_1,\ldots,b_l$, which are the roots of $Q$. Once again what we have here are certain symmetric functions in these variables $b_1,\ldots,b_l$, and these symmetric functions must be certain polynomials in the coefficients of $Q$. Moreover, due to our normalization factor $d^k$, these polynomials in the coefficients of $Q$ will have integer coefficients.

\medskip

(6) Thus, we are led to the conclusion in the statement, that $R(P,Q)$ is a polynomial in the coefficients of $P,Q$, with integer coefficients, and with the remark that the $c^ld^k$ factor is there for these latter coefficients to be indeed integers, instead of rationals.
\end{proof}

All this might seem a bit complicated, so as an illustration, let us work out an example. Consider the case of a polynomial of degree 2, and a polynomial of degree 1:
$$P=ax^2+bx+c\quad,\quad 
Q=dx+e$$

In order to compute the resultant, let us factorize our polynomials:
$$P=a(x-p)(x-q)\quad,\quad 
Q=d(x-r)$$

The resultant can be then computed as follows, by using the method above:
\begin{eqnarray*}
R(P,Q)
&=&ad^2(p-r)(q-r)\\
&=&ad^2(pq-(p+q)r+r^2)\\
&=&cd^2+bd^2r+ad^2r^2\\
&=&cd^2-bde+ae^2
\end{eqnarray*}

Finally, observe that $R(P,Q)=0$ corresponds indeed to the fact that $P,Q$ have a common root. Indeed, the root of $Q$ is $r=-e/d$, and we have:
$$P(r)
=\frac{ae^2}{d^2}-\frac{be}{d}+c
=\frac{R(P,Q)}{d^2}$$

Thus we have $P(r)=0$ precisely when $R(P,Q)=0$, as predicted by Theorem 4.14. 

\bigskip

Regarding now the explicit formula of the resultant $R(P,Q)$, this is something quite complicated, and there are several methods for dealing with this problem. There is a slight similarity between Theorem 4.14 and the Vandermonde determinants discussed in chapter 2, and we have in fact the following formula for $R(P,Q)$:

\index{resultant}

\begin{theorem}
The resultant of two polynomials, written as 
$$P=p_kX^k+\ldots+p_1X+p_0\quad,\quad 
Q=q_lX^l+\ldots+q_1X+q_0$$
appears as the determinant of an associated matrix, as follows,
$$R(P,Q)=
\begin{vmatrix}
p_k&&&q_l\\
\vdots&\ddots&&\vdots&\ddots\\
p_0&&p_k&q_0&&q_l\\
&\ddots&\vdots&&\ddots&\vdots\\
&&p_0&&&q_0
\end{vmatrix}
$$
with the matrix having size $k+l$, and having $0$ coefficients at the blank spaces.
\end{theorem}

\begin{proof}
This is something quite clever, due to Sylvester, as follows:

\medskip

(1) Consider the vector space $\mathbb C_k[X]$ formed by the polynomials of degree $<k$:
$$\mathbb C_k[X]=\left\{P\in\mathbb C[X]\Big|\deg P<k\right\}$$

This is a vector space of dimension $k$, having as basis the monomials $1,X,\ldots,X^{k-1}$. Now given polynomials $P,Q$ as in the statement, consider the following linear map:
$$\Phi:\mathbb C_l[X]\times\mathbb C_k[X]\to\mathbb C_{k+l}[X]\quad,\quad (A,B)\to AP+BQ$$

\medskip

(2) Our first claim is that with respect to the standard bases for all the vector spaces involved, namely those consisting of the monomials $1,X,X^2,\ldots$, the matrix of $\Phi$ is the matrix in the statement. But this is something which is clear from definitions.

\medskip

(3) Our second claim is that $\det\Phi=0$ happens precisely when $P,Q$ have a common root. Indeed, our polynomials $P,Q$ having a common root means that we can find $A,B$ such that $AP+BQ=0$, and so that $(A,B)\in\ker\Phi$, which reads $\det\Phi=0$.

\medskip

(4) Finally, our claim is that we have $\det\Phi=R(P,Q)$. But this follows from the uniqueness of the resultant, up to a scalar, and with this uniqueness property being elementary to establish, along the lines of the proof of Theorem 4.14.
\end{proof}

As an illustration, consider our favorite polynomials, as before:
$$P=ax^2+bx+c\quad,\quad 
Q=dx+e$$

According to the above result, the resultant should be then, as it should:
$$R(P,Q)
=\begin{vmatrix}
a&d&0\\
b&e&d\\
c&0&e
\end{vmatrix}
=ae^2-bde+cd^2$$

Now back to our diagonalization questions, we want to compute $R(P,P')$, where $P$ is the characteristic polynomial. So, we need one more piece of theory, as follows:

\index{discriminant}
\index{double root}
\index{single roots}

\begin{theorem}
Given a polynomial $P\in\mathbb C[X]$, written as
$$P(X)=cX^N+dX^{N-1}+\ldots$$
its discriminant, defined as being the following quantity,
$$\Delta(P)=\frac{(-1)^{\binom{N}{2}}}{c}R(P,P')$$
is a polynomial in the coefficients of $P$, with integer coefficients, and
$$\Delta(P)=0$$
happens precisely when $P$ has a double root.
\end{theorem}

\begin{proof}
The fact that the discriminant $\Delta(P)$ is a polynomial in the coefficients of $P$, with integer coefficients, comes from Theorem 4.14, coupled with the fact that the division by the leading coefficient $a$ is indeed possible, under $\mathbb Z$, as being shown by:
$$R(P,P')=
\begin{vmatrix}
a&&&Na\\
\vdots&\ddots&&\vdots&\ddots\\
z&&a&y&&Na\\
&\ddots&\vdots&&\ddots&\vdots\\
&&z&&&y
\end{vmatrix}
$$

Also, the fact that we have $\Delta(P)=0$ precisely when $P$ has a double root is clear from Theorem 4.14. Finally, let us mention that the sign $(-1)^{\binom{N}{2}}$ is there for various reasons, including the compatibility with the formula $\Delta(P)=b^2-4ac$ in degree 2.
\end{proof}

As an illustration, let us see what happens in degree 2. Here we have:
$$P=aX^2+bX+c\quad,\quad 
P'=2aX+b$$

Thus, the resultant is given by the following formula:
\begin{eqnarray*}
R(P,P')
&=&ab^2-b(2a)b+c(2a)^2\\
&=&4a^2c-ab^2\\
&=&-a(b^2-4ac)
\end{eqnarray*}

It follows that the discriminant of our polynomial is, as it should:
$$\Delta(P)=b^2-4ac$$

Alternatively, we can use the formula in Theorem 4.15, and we obtain:
$$\Delta(P)=
-\frac{1}{a}\begin{vmatrix}
a&2a&\\
b&b&2a\\
c&&b
\end{vmatrix}
=b^2-4ac$$

At the theoretical level now, we have the following result, which is not trivial:

\begin{theorem}
The discriminant of a polynomial $P$ is given by the formula
$$\Delta(P)=a^{2N-2}\prod_{i<j}(r_i-r_j)^2$$
where $a$ is the leading coefficient, and $r_1,\ldots,r_N$ are the roots.
\end{theorem}

\begin{proof}
This is something quite tricky, the idea being as follows:

\medskip

(1) The first thought goes to the formula in Theorem 4.14, so let us see what that formula teaches us, in the case $Q=P'$. Let us write $P,P'$ as follows:
$$P=a(x-r_1)\ldots(x-r_N)$$ 
$$P'=Na(x-p_1)\ldots(x-p_{N-1})$$

According to Theorem 4.14, the resultant of $P,P'$ is then given by:
$$R(P,P')=a^{N-1}(Na)^N\prod_{ij}(r_i-p_j)$$

And bad news, this is not exactly what we wished for, namely the formula in the statement. That is, we are on the good way, but certainly have to work some more.

\medskip

(2) Obviously, we must get rid of the roots $p_1,\ldots,p_{N-1}$ of the polynomial $P'$. In order to do this, let us rewrite the formula that we found in (1) in the following way:
\begin{eqnarray*}
R(P,P')
&=&N^Na^{2N-1}\prod_i\left(\prod_j(r_i-p_j)\right)\\
&=&N^Na^{2N-1}\prod_i\frac{P'(r_i)}{Na}\\
&=&a^{N-1}\prod_iP'(r_i)
\end{eqnarray*}

(3) In order to compute now $P'$, and more specifically the values $P'(r_i)$ that we are interested in, we can use the Leibnitz rule. So, consider our polynomial:
$$P(x)=a(x-r_1)\ldots(x-r_N)$$ 

The Leibnitz rule for derivatives tells us that $(fg)'=f'g+fg'$, but then also that $(fgh)'=f'gh+fg'h+fgh'$, and so on. Thus, for our polynomial, we obtain:
$$P'(x)=a\sum_i(x-r_1)\ldots\underbrace{(x-r_i)}_{missing}\ldots(x-r_N)$$ 

Now when applying this formula to one of the roots $r_i$, we obtain:
$$P'(r_i)=a(r_i-r_1)\ldots\underbrace{(r_i-r_i)}_{missing}\ldots(r_i-r_N)$$ 

By making now the product over all indices $i$, this gives the following formula:
$$\prod_iP'(r_i)=a^N\prod_{i\neq j}(r_i-r_j)$$

(4) Time now to put everything together. By taking the formula in (2), making the normalizations in Theorem 4.16, and then using the formula found in (3), we obtain:
\begin{eqnarray*}
\Delta(P)
&=&(-1)^{\binom{N}{2}}a^{N-2}\prod_iP'(r_i)\\
&=&(-1)^{\binom{N}{2}}a^{2N-2}\prod_{i\neq j}(r_i-r_j)\\
&=&a^{2N-2}\prod_{i<j}(r_i-r_j)^2
\end{eqnarray*}

Thus, we are led to the conclusion in the statement.
\end{proof}

Now back to our linear algebra questions, we can upgrade Theorem 4.13, as follows:

\begin{theorem}
For a matrix $A\in M_N(\mathbb C)$ the following conditions are equivalent,
\begin{enumerate}
\item The eigenvalues are different, $\lambda_i\neq\lambda_j$,

\item The characteristic polynomial $P$ has simple roots,

\item The discriminant of $P$ is nonzero, $\Delta(P)\neq0$,
\end{enumerate}
and in this case, the matrix is diagonalizable.
\end{theorem}

\begin{proof}
This is indeed an upgrade of Theorem 4.13, by replacing the condition (3) there with the condition $\Delta(P)\neq0$, which is something better, computational.
\end{proof}

As mentioned before, in Theorem 4.7,  one can prove that the matrices having distinct eigenvalues are ``generic'', so the above result basically captures the whole situation. We have in fact the following collection of density results, all being very useful:

\index{density}
\index{matrices with distinct eigenvalues}

\begin{theorem}
The following happen, inside $M_N(\mathbb C)$:
\begin{enumerate}
\item The invertible matrices are dense.

\item The matrices having distinct eigenvalues are dense.

\item The diagonalizable matrices are dense.
\end{enumerate}
\end{theorem}

\begin{proof}
These are quite advanced linear algebra results, which can be proved as follows, with the technology that we have so far:

\medskip

(1) This is clear, intuitively speaking, because the invertible matrices are given by the condition $\det A\neq 0$. Thus, the set formed by these matrices appears as the complement of the hypersurface $\det A=0$, and so must be dense inside $M_N(\mathbb C)$, as claimed. 

\medskip

(2) Here we can use a similar argument, this time by saying that the set formed by the matrices having distinct eigenvalues appears as the complement of the hypersurface given by $\Delta(P_A)=0$, and so must be dense inside $M_N(\mathbb C)$, as claimed. 

\medskip

(3) This follows from (2), via the fact that the matrices having distinct eigenvalues are diagonalizable, that we know from Theorem 4.18. There are of course some other proofs as well, for instance by putting the matrix in Jordan form.
\end{proof}

As an application of the above results, and of our methods in general, we can now establish a number of useful and interesting linear algebra results, as follows:

\index{products of matrices}
\index{functions of matrices}

\begin{theorem}
The following happen:
\begin{enumerate}
\item We have $P_{AB}=P_{BA}$, for any two matrices $A,B\in M_N(\mathbb C)$.

\item $AB,BA$ have the same eigenvalues, with the same multiplicities.

\item If $A$ has eigenvalues $\lambda_1,\ldots,\lambda_N$, then $f(A)$ has eigenvalues $f(\lambda_1),\ldots,f(\lambda_N)$.
\end{enumerate}
\end{theorem}

\begin{proof}
These results can be deduced by using Theorem 4.19, as follows:

\medskip

(1) It follows from definitions that the characteristic polynomial of a matrix is invariant under conjugation, in the sense that we have the following formula:
$$P_C=P_{ACA^{-1}}$$

Now observe that, when assuming that $A$ is invertible, we have:
$$AB=A(BA)A^{-1}$$

Thus, we have the result when $A$ is invertible. By using now Theorem 4.19 (1), we conclude that this formula holds for any matrix $A$, by continuity. 

\medskip

(2) This is a reformulation of (1) above, via the fact that $P$ encodes the eigenvalues, with multiplicities, which is hard to prove with bare hands.

\medskip

(3) This is something more informal, the idea being that this is clear for the diagonal matrices $D$, then for the diagonalizable matrices $PDP^{-1}$, and finally for all matrices, by using Theorem 4.19 (3), provided that $f$ has suitable regularity properties. 
\end{proof}

The last assertion in the above theorem remains of course to be clarified, and we will be back to this in chapter 8 below, with details, when doing spectral theory. 

\section*{4c. Spectral theorems}

Let us go back now to the diagonalization question. Here is a key result:

\index{spectral theorem}
\index{self-adjoint matrix}

\begin{theorem}
Any matrix $A\in M_N(\mathbb C)$ which is self-adjoint, $A=A^*$, is diagonalizable, with the diagonalization being of the following type,
$$A=UDU^*$$
with $U\in U_N$, and with $D\in M_N(\mathbb R)$ diagonal. The converse holds too.
\end{theorem}

\begin{proof}
As a first remark, the converse trivially holds, because if we take a matrix of the form $A=UDU^*$, with $U$ unitary and $D$ diagonal and real, then we have:
$$A^*
=(UDU^*)^*
=UD^*U^*
=UDU^*
=A$$

In the other sense now, assume that $A$ is self-adjoint, $A=A^*$.  Our first claim is that the eigenvalues are real. Indeed, assuming $Av=\lambda v$, we have:
\begin{eqnarray*}
\lambda<v,v>
&=&<Av,v>\\
&=&<v,Av>\\
&=&<v,\lambda v>\\
&=&\bar{\lambda}<v,v>
\end{eqnarray*}

Thus we obtain $\lambda\in\mathbb R$, as claimed. Our next claim now is that the eigenspaces corresponding to different eigenvalues are pairwise orthogonal. Assume indeed that:
$$Av=\lambda v\quad,\quad 
Aw=\mu w$$

We have then the following computation, using $\lambda,\mu\in\mathbb R$:
\begin{eqnarray*}
\lambda<v,w>
&=&<Av,w>\\
&=&<v,Aw>\\
&=&<v,\mu w>\\
&=&\mu<v,w>
\end{eqnarray*}

Thus $\lambda\neq\mu$ implies $v\perp w$, as claimed. In order now to finish, it remains to prove that the eigenspaces span $\mathbb C^N$. For this purpose, we will use a recurrence method. Let us pick an eigenvector, $Av=\lambda v$. Assuming $v\perp w$, we have:
\begin{eqnarray*}
<Aw,v>
&=&<w,Av>\\
&=&<w,\lambda v>\\
&=&\lambda<w,v>\\
&=&0
\end{eqnarray*}

Thus, if $v$ is an eigenvector, then the vector space $v^\perp$ is invariant under $A$. In order to do the recurrence, it still remains to prove that the restriction of $A$ to the vector space $v^\perp$ is self-adjoint. But this comes from a general property of the self-adjoint matrices, that we will explain now. Our claim is that an arbitary square matrix $A$ is self-adjoint precisely when the following happens, for any vector $v$:
$$<Av,v>\in\mathbb R$$

Indeed, the fact that the above scalar product is real is equivalent to:
$$<(A-A^*)v,v>=0$$

But this is equivalent, by developing the scalar product, to $A=A^*$, so our claim is proved. Now back to our questions, it is clear from our self-adjointness criterion above that the restriction of $A$ to any invariant subspace, and in particular to the subspace $v^\perp$, is self-adjoint. Thus, we can proceed by recurrence, and we obtain the result.
\end{proof}

Let us record as well the real version of the above result:

\index{spectral theorem}
\index{symmetric matrix}

\begin{theorem}
Any matrix $A\in M_N(\mathbb R)$ which is symmetric, in the sense that
$$A=A^t$$
is diagonalizable, with the diagonalization being of the following type,
$$A=UDU^t$$
with $U\in O_N$, and with $D\in M_N(\mathbb R)$ diagonal. The converse holds too.
\end{theorem}

\begin{proof}
As before, the converse trivially holds, because if we take a matrix of the form $A=UDU^t$, with $U$ orthogonal and $D$ diagonal and real, then we have $A^t=A$. In the other sense now, this follows from Theorem 4.21, and its proof.
\end{proof}

As basic examples of self-adjoint matrices, we have the orthogonal projections. The diagonalization result regarding them is as follows:

\index{projection}

\begin{proposition}
The matrices $P\in M_N(\mathbb C)$ which are projections, $P^2=P^*=P$, 
are precisely those which diagonalize as follows,
$$P=UDU^*$$
with $U\in U_N$, and with $D\in M_N(0,1)$ being diagonal.
\end{proposition}

\begin{proof}
This is clear, geometrically, with the diagonalization being as follows, with the 1-eigenspace being the image of $P$, and the 0-eigenspace being the kernel:
$$P\sim\begin{pmatrix}
1\\
&\ddots\\
&&1\\
&&&0\\
&&&&\ddots\\
&&&&&0
\end{pmatrix}$$

Alternatively, we can get this algebrically, from $P^2=P^*=P$. Indeed, $P^*=P$ shows that the eigenvalues are real, and then, assuming $Pv=\lambda v$, by using $P^2=P$ we get:
\begin{eqnarray*}
\lambda<v,v>
&=&<Pv,v>\\
&=&<P^2v,v>\\
&=&<Pv,Pv>\\
&=&<\lambda v,\lambda v>\\
&=&\lambda^2<v,v>
\end{eqnarray*}

We therefore have $\lambda\in\{0,1\}$, and the rest comes from Theorem 4.21.
\end{proof}

In the real case, the result regarding the projections is as follows:

\begin{proposition}
The matrices $P\in M_N(\mathbb R)$ which are projections,
$$P^2=P^t=P$$
are precisely those which diagonalize as follows,
$$P=UDU^t$$
with $U\in O_N$, and with $D\in M_N(0,1)$ being diagonal.
\end{proposition}

\begin{proof}
This follows indeed from Proposition 4.23, and its proof. 
\end{proof}

An important class of self-adjoint matrices, that we will discuss now, which includes all projections, are the positive matrices. The general theory here is as follows:

\index{positive matrix}

\begin{theorem}
For a matrix $A\in M_N(\mathbb C)$ the following conditions are equivalent, and if they are satisfied, we say that $A$ is positive, and write $A\geq0$:
\begin{enumerate}
\item $A=B^2$, with $B=B^*$.

\item $A=CC^*$, for some $C\in M_N(\mathbb C)$.

\item $<Ax,x>\geq0$, for any vector $x\in\mathbb C^N$.

\item $A=A^*$, and the eigenvalues are positive, $\lambda_i\geq0$.

\item $A=UDU^*$, with $U\in U_N$ and with $D\in M_N(\mathbb R_+)$ diagonal.
\end{enumerate}
\end{theorem}

\begin{proof}
The idea is that the equivalences in the statement basically follow from some elementary computations, with only Theorem 4.21 needed, at some point:

\medskip

$(1)\implies(2)$ This is clear, because we can take $C=B$.

\medskip

$(2)\implies(3)$ This comes indeed from the following computation:
$$<Ax,x>
=<CC^*x,x>
=<C^*x,C^*x>
\geq0$$

$(3)\implies(4)$ By using the fact that $<Ax,x>$ is real, we have:
$$<Ax,x>
=<x,A^*x>
=<A^*x,x>$$

Thus we have $A=A^*$, and the remaining assertion, regarding the eigenvalues, follows from the following computation, assuming $Ax=\lambda x$:
$$<Ax,x>
=<\lambda x,x>
=\lambda<x,x>
\geq0$$

$(4)\implies(5)$ This follows indeed by using Theorem 4.21.

\medskip

$(5)\implies(1)$ Assuming $A=UDU^*$ as in the statement, we can set $B=U\sqrt{D}U^*$. Then this matrix $B$ is self-adjoint, and its square is given by:
\begin{eqnarray*}
B^2
&=&U\sqrt{D}U^*\cdot U\sqrt{D}U^*\\
&=&UDU^*\\
&=&A
\end{eqnarray*}

Thus, we are led to the conclusion in the statement.
\end{proof}

Let us record as well the following technical version of the above result:

\index{strictly positive matrix}

\begin{theorem}
For a matrix $A\in M_N(\mathbb C)$ the following conditions are equivalent, and if they are satisfied, we say that $A$ is strictly positive, and write $A>0$:
\begin{enumerate}
\item $A=B^2$, with $B=B^*$, invertible.

\item $A=CC^*$, for some $C\in M_N(\mathbb C)$ invertible.

\item $<Ax,x>>0$, for any nonzero vector $x\in\mathbb C^N$.

\item $A=A^*$, and the eigenvalues are strictly positive, $\lambda_i>0$.

\item $A=UDU^*$, with $U\in U_N$ and with $D\in M_N(\mathbb R_+^*)$ diagonal.
\end{enumerate}
\end{theorem}

\begin{proof}
This follows either from Theorem 4.25, by adding the various extra assumptions in the statement, or from the proof of Theorem 4.25, by modifying where needed.
\end{proof}

The positive matrices are quite important, for a number of reasons. On one hand, these are the matrices $A\in M_N(\mathbb C)$ having a square root $\sqrt{A}\in M_N(\mathbb C)$, as shown by our positivity condition (1). On the other hand, any matrix $A\in M_N(\mathbb C)$ produces the positive matrix $A^*A\in M_N(\mathbb C)$, as shown by our positivity condition (2). We can combine these two observations, and we are led to the following construction, for any $A\in M_N(\mathbb C)$:
$$A\to\sqrt{A^*A}$$

Which is something quite interesting, because at $N=1$ what we have here is the construction of the absolute value of complex numbers, $|z|=\sqrt{z\bar{z}}$. This suggests using the notation $|A|=\sqrt{A^*A}$, and then looking for a decomposition result of type:
$$A=U|A|$$

We will be back to this type of decomposition later, called polar decomposition, at the end of the present chapter, after developing some more general theory.

\bigskip

Let us discuss now the case of the unitary matrices. We have here:

\index{spectral theorem}
\index{unitary matrix}

\begin{theorem}
Any matrix $U\in M_N(\mathbb C)$ which is unitary, $U^*=U^{-1}$, 
is diagonalizable, with the eigenvalues on $\mathbb T$. More precisely we have
$$U=VDV^*$$
with $V\in U_N$, and with $D\in M_N(\mathbb T)$ diagonal. The converse holds too.
\end{theorem}

\begin{proof}
As a first remark, the converse trivially holds, because given a matrix of type $U=VDV^*$, with $V\in U_N$, and with $D\in M_N(\mathbb T)$ being diagonal, we have:
\begin{eqnarray*}
U^*
&=&(VDV^*)^*\\
&=&VD^*V^*\\
&=&VD^{-1}V^{-1}\\
&=&(V^*)^{-1}D^{-1}V^{-1}\\
&=&(VDV^*)^{-1}\\
&=&U^{-1}
\end{eqnarray*}

Let us prove now the first assertion, stating that the eigenvalues of a unitary matrix $U\in U_N$ belong to $\mathbb T$. Indeed, assuming $Uv=\lambda v$, we have:
\begin{eqnarray*}
<v,v>
&=&<U^*Uv,v>\\
&=&<Uv,Uv>\\
&=&<\lambda v,\lambda v>\\
&=&|\lambda|^2<v,v>
\end{eqnarray*}

Thus we obtain $\lambda\in\mathbb T$, as claimed. Our next claim now is that the eigenspaces corresponding to different eigenvalues are pairwise orthogonal. Assume indeed that:
$$Uv=\lambda v\quad,\quad 
Uw=\mu w$$

We have then the following computation, using $U^*=U^{-1}$ and $\lambda,\mu\in\mathbb T$:
\begin{eqnarray*}
\lambda<v,w>
&=&<\lambda v,w>\\
&=&<Uv,w>\\
&=&<v,U^*w>\\
&=&<v,U^{-1}w>\\
&=&<v,\mu^{-1}w>\\
&=&\mu<v,w>
\end{eqnarray*}

Thus $\lambda\neq\mu$ implies $v\perp w$, as claimed. In order now to finish, it remains to prove that the eigenspaces span $\mathbb C^N$. For this purpose, we will use a recurrence method. Let us pick an eigenvector, $Uv=\lambda v$. Assuming $v\perp w$, we have:
\begin{eqnarray*}
<Uw,v>
&=&<w,U^*v>\\
&=&<w,U^{-1}v>\\
&=&<w,\lambda^{-1}v>\\
&=&\lambda<w,v>\\
&=&0
\end{eqnarray*}

Thus, if $v$ is an eigenvector, then the vector space $v^\perp$ is invariant under $U$. Now since $U$ is an isometry, so is its restriction to this space $v^\perp$. Thus this restriction is a unitary, and so we can proceed by recurrence, and we obtain the result.
\end{proof}

Let us record as well the real version of the above result, in a weak form:

\index{spectral theorem}
\index{orthogonal matrix}

\begin{theorem}
Any matrix $U\in M_N(\mathbb R)$ which is orthogonal, $U^t=U^{-1}$,
is diagonalizable, with the eigenvalues on $\mathbb T$. More precisely we have
$$U=VDV^*$$
with $V\in U_N$, and with $D\in M_N(\mathbb T)$ being diagonal.
\end{theorem}

\begin{proof}
This follows indeed from Theorem 4.27.
\end{proof}

Observe that the above result does not provide us with a complete characterization of the matrices $U\in M_N(\mathbb R)$ which are orthogonal. To be more precise, the question left is that of understanding when the matrices of type $U=VDV^*$, with $V\in U_N$, and with $D\in M_N(\mathbb T)$ being diagonal, are real, and this is something non-trivial.

\bigskip

As an illustration, for the simplest unitaries that we know, namely the rotations in the real plane, we have the following formula, that we know well from chapter 3:
$$\begin{pmatrix}\cos t&-\sin t\\ \sin t&\cos t\end{pmatrix}
=\frac{1}{2}\begin{pmatrix}1&1\\i&-i\end{pmatrix}
\begin{pmatrix}e^{-it}&0\\0&e^{it}\end{pmatrix}
\begin{pmatrix}1&-i\\1&i\end{pmatrix}$$

We will be back to such questions later, when discussing the orthogonal groups.

\section*{4d. Normal matrices}

Back to generalities, the self-adjoint matrices and the unitary matrices are particular cases of the general notion of a ``normal matrix'', and we have here:

\index{spectral theorem}
\index{normal matrix}

\begin{theorem}
Any matrix $A\in M_N(\mathbb C)$ which is normal, $AA^*=A^*A$, is diagonalizable, with the diagonalization being of the following type,
$$A=UDU^*$$
with $U\in U_N$, and with $D\in M_N(\mathbb C)$ diagonal. The converse holds too.
\end{theorem}

\begin{proof}
As a first remark, the converse trivially holds, because if we take a matrix of the form $A=UDU^*$, with $U$ unitary and $D$ diagonal, then we have:
\begin{eqnarray*}
AA^*
&=&UDU^*\cdot UD^*U^*\\
&=&UDD^*U^*\\
&=&UD^*DU^*\\
&=&UD^*U^*\cdot UDU^*\\
&=&A^*A
\end{eqnarray*}

In the other sense now, this is something more technical. Our first claim is that a matrix $A$ is normal precisely when the following happens, for any vector $v$:
$$||Av||=||A^*v||$$

Indeed, the above equality can be written as follows:
$$<AA^*v,v>=<A^*Av,v>$$

But this is equivalent to $AA^*=A^*A$, by using the polarization identity. Our claim now is that $A,A^*$ have the same eigenvectors, with conjugate eigenvalues:
$$Av=\lambda v\implies A^*v=\bar{\lambda}v$$

Indeed, this follows from the following computation, and from the trivial fact that if $A$ is normal, then so is any matrix of type $A-\lambda 1_N$:
\begin{eqnarray*}
||(A^*-\bar{\lambda}1_N)v||
&=&||(A-\lambda 1_N)^*v||\\
&=&||(A-\lambda 1_N)v||\\
&=&0
\end{eqnarray*}

Let us prove now, by using this, that the eigenspaces of $A$ are pairwise orthogonal. Assuming $Av=\lambda v$ and $Aw=\mu w$ with $\lambda\neq\mu$, we have:
\begin{eqnarray*}
\lambda<v,w>
&=&<\lambda v,w>\\
&=&<Av,w>\\
&=&<v,A^*w>\\
&=&<v,\bar{\mu}w>\\
&=&\mu<v,w>
\end{eqnarray*}

Thus $\lambda\neq\mu$ implies $v\perp w$, as claimed. In order to finish now the proof, it remains to prove that the eigenspaces of $A$ span the whole $\mathbb C^N$. This is something that we have already seen for the self-adjoint matrices, and for the unitaries, and we will use here these results, in order to deal with the general normal case. As a first observation, given an arbitrary matrix $A$, the matrix $AA^*$ is self-adjoint:
$$(AA^*)^*=AA^*$$

Thus, we can diagonalize this matrix $AA^*$, as follows, with the passage matrix being a unitary, $V\in U_N$, and with the diagonal form being real, $E\in M_N(\mathbb R)$:
$$AA^*=VEV^*$$

Now observe that, for matrices of type $A=UDU^*$, which are those that we supposed to deal with, we have $V=U,E=D\bar{D}$. In particular, $A$ and $AA^*$ have the same eigenspaces. So, this will be our idea, proving that the eigenspaces of $AA^*$ are eigenspaces of $A$. In order to do so, let us pick two eigenvectors $v,w$ of the matrix $AA^*$, corresponding to different eigenvalues, $\lambda\neq\mu$. The eigenvalue equations are then as follows:
$$AA^*v=\lambda v\quad,\quad 
AA^*w=\mu w$$

We have the following computation, using the normality condition $AA^*=A^*A$, and the fact that the eigenvalues of $AA^*$, and in particular $\mu$, are real:
\begin{eqnarray*}
\lambda<Av,w>
&=&<\lambda Av,w>\\
&=&<A\lambda v,w>\\
&=&<AAA^*v,w>\\
&=&<AA^*Av,w>\\
&=&<Av,AA^*w>\\
&=&<Av,\mu w>\\
&=&\mu<Av,w>
\end{eqnarray*}

We conclude that we have $<Av,w>=0$. But this reformulates as follows:
$$\lambda\neq\mu\implies A(E_\lambda)\perp E_\mu$$

Now since the eigenspaces of $AA^*$ are pairwise orthogonal, and span the whole $\mathbb C^N$, we deduce from this that these eigenspaces are invariant under $A$:
$$A(E_\lambda)\subset E_\lambda$$

But with this result in hand, we can finish. Indeed, we can decompose the problem, and the matrix $A$ itself, following these eigenspaces of $AA^*$, which in practice amounts in saying that we can assume that we only have 1 eigenspace. But by rescaling, this is the same as assuming that we have $AA^*=1$, and with this done, we are now into the unitary case, that we know how to solve, as explained in Theorem 4.27.
\end{proof}

Let us discuss now the polar decomposition. We first have the following result:

\index{absolute value}
\index{modulus}

\begin{theorem}
Given a matrix $A\in M_N(\mathbb C)$, we can construct a matrix $|A|$ as follows, by using the fact that $A^*A$ is diagonalizable, with positive eigenvalues:
$$|A|=\sqrt{A^*A}$$
This matrix $|A|$ is then positive, and its square is $|A|^2=A$. In the case $N=1$, we obtain in this way the usual absolute value of the complex numbers.
\end{theorem}

\begin{proof}
Consider indeed the matrix $A^*A$, which is normal. According to Theorem 4.29, we can diagonalize this matrix as follows, with $U\in U_N$, and with $D$ diagonal:
$$A=UDU^*$$

Since we have $A^*A\geq0$, it follows that we have $D\geq0$, which means that the entries of $D$ are real, and positive. Thus we can extract the square root $\sqrt{D}$, and then set:
$$\sqrt{A^*A}=U\sqrt{D}U^*$$

Now if we call this latter matrix $|A|$, we are led to the conclusions in the statement, namely $|A|\geq0$, and $|A|^2=A$. Finally, the last assertion is clear from definitions.
\end{proof}

We can now formulate a first polar decomposition result, as follows:

\index{polar decomposition}

\begin{theorem}
Any invertible matrix $A\in M_N(\mathbb C)$ decomposes as
$$A=U|A|$$
with $U\in U_N$, and with $|A|=\sqrt{A^*A}$ as above.
\end{theorem}

\begin{proof}
According to our definition of the modulus, $|A|=\sqrt{A^*A}$, we have:
\begin{eqnarray*}
<|A|x,|A|y>
&=&<x,|A|^2y>\\
&=&<x,A^*Ay>\\
&=&<Ax,Ay>
\end{eqnarray*}

Thus we can define a unitary matrix $U\in U_N$ by the following formula:
$$U(|A|x)=Ax$$

But this formula shows that we have $A=U|A|$, as desired.
\end{proof}

Observe that at $N=1$ we obtain in this way the usual polar decomposition of the nonzero complex numbers. More generally now, we have the following result:

\index{polar decomposition}
\index{partial isometry}

\begin{theorem}
Any square matrix $A\in M_N(\mathbb C)$ decomposes as
$$A=U|A|$$
with $U$ being a partial isometry, and with $|A|=\sqrt{A^*A}$ as above.
\end{theorem}

\begin{proof}
Once again, this follows by comparing the actions of $A,|A|$ on the vectors $v\in\mathbb C^N$, and deducing from this the existence of a partial isometry $U$ as above. Alternatively, we can get this from Theorem 4.31, applied on the complement of the 0-eigenvectors. 
\end{proof}

And with this, good news, done with linear algebra. We have learned many things in the past 100 pages, and our knowledge of the subject is quite decent, and we will stop here. In the remainder of this book we will be rather looking into applications.

\section*{4e. Exercises}

Things have been quite dense in this chapter, which was our last one on basic linear algebra, with some details missing. As a first exercise, in relation with abstract vector calculus, that we somehow assumed to be reasonably known, we have:

\begin{exercise}
Clarify the theory of linear spaces $V\subset\mathbb C^N$, notably with:
\begin{enumerate}
\item A standard discussion regarding generating sets, linear independence, bases.

\item Injectivity, surjectivity and bijectivity of the linear maps $f:\mathbb C^N\to\mathbb C^N$.

\item More generally, $\dim(\ker f)+\dim(Im f)=N$, for such maps $f:\mathbb C^N\to\mathbb C^N$.
\end{enumerate}
Then, extend this into a theory of linear spaces $V$, not necessarily subspaces of $\mathbb C^N$.
\end{exercise}

Here the first question is something quite standard, by using our linear algebra knowledge. As for the second question, things are a bit more tricky here, because once the abstract linear spaces $V$ are defined, the only available tool is recurrence.

\begin{exercise}
Work out what happens to the main diagonalization theorem for the matrices $A\in M_N(\mathbb C)$, in the cases $A\in M_2(\mathbb C)$, $A\in M_N(\mathbb R)$, and $A\in M_2(\mathbb R)$.
\end{exercise}

As before, this is a rather theoretical exercise, the point being that of carefully reviewing all the material above, in the 3 particular cases which are indicated.

\begin{exercise}
Clarify which functions can be applied to which matrices, as to have results stating that the eigenvalues of $f(A)$ are $f(\lambda_1),\ldots,f(\lambda_N)$.
\end{exercise}

This exercise is actually quite difficult, with various technical assumptions being needed on both $f$ and $A$, as for everything to work fine. We will be back to this.

\begin{exercise}
Work out specialized spectral theorems for the orthogonal matrices $U\in O_N$, going beyond what has been said in the above.
\end{exercise}

To be more precise here, we have proved many spectral theorems in the above, but the case $U\in O_N$, where our statement here was something quite weak, coming without a converse, is obviously still in need of discussion. Again, this is something non-trivial.

\begin{exercise}
Prove that any matrix can be put in Jordan form, 
$$A\sim\begin{pmatrix}
J_1\\
&\ddots\\
&&J_k
\end{pmatrix}\qquad,\qquad J_i=\begin{pmatrix}
\lambda_i&1\\
&\ddots&\ddots\\
&&\lambda_i&1\\
&&&\lambda_i
\end{pmatrix}$$
with the size of each Jordan block $J_i$ being the multiplicity of $\lambda_i$.
\end{exercise}

This is something useful, because it applies to any matrix $A\in M_N(\mathbb C)$, without assumptions, and is somewhat the ``nuclear option'' in linear algebra.

\part{Matrix analysis}

\ \vskip50mm

\begin{center}
{\em Everything dies, baby, that's a fact

But maybe everything that dies some day comes back

Put your makeup on, fix your hair up pretty

And meet me tonight in Atlantic City}
\end{center}

\chapter{Basic calculus}

\section*{5a. Real analysis}

We discuss in what follows some applications of the theory that we developed above, to basic  questions in analysis. The idea will be that the functions of several variables $f:\mathbb R^N\to\mathbb R^M$ can be locally approximated by linear maps, in the same way as the functions $f:\mathbb R\to\mathbb R$ can be locally approximated by using derivatives:
$$f(x+t)\simeq f(x)+f'(x)t\quad,\quad f'(x)\in M_{M\times N}(\mathbb R)$$

There are many things that can be said here, and at order 2 too, and we will be quite brief. Getting started now, let us first discuss the simplest case, $f:\mathbb R\to\mathbb R$. Here we have the following result, which is the starting point for everything in analysis:

\index{derivative}

\begin{theorem}
Any function $f:\mathbb R\to\mathbb R$ is approximately locally affine,
$$f(x+t)\simeq f(x)+f'(x)t$$
with $f'(x)\in\mathbb R$ being the derivative of $f$ at the point $x$, given by
$$f'(x)=\lim_{t\to0}\frac{f(x+t)-f(x)}{t}$$
provided that this latter limit converges indeed.
\end{theorem}

\begin{proof}
This is something trivial, because if the limit in the statement converges, by multiplying by $t$ we obtain the above estimate for $f(x+t)$. Observe also that, by drawing the graph of $f$, we can see that $f'(x)$ compute the slope, at the given point $x$. Finally, as a basic counterexample, observe that $f(x)=|x|$ is not differentiable at $x=0$.
\end{proof}

As a first illustration, the derivatives of power functions are as follows:

\begin{proposition}
We have the differentiation formula
$$(x^p)'=px^{p-1}$$
valid for any exponent $p\in\mathbb R$.
\end{proposition}

\begin{proof}
In the case $p\in\mathbb N$ we can use the binomial formula, which gives:
$$(x+t)^p
=x^p+px^{p-1}t+\ldots+t^p
\simeq x^p+px^{p-1}t$$

Next, for $p\in\mathbb Q$, we can write $p=m/n$, with $m\in\mathbb N$ and $n\in\mathbb Z$, and we have: 
\begin{eqnarray*}
(x+t)^{m/n}-x^{m/n}
&=&\frac{(x+t)^m-x^m}{(x+t)^{m(n-1)/n}+\ldots+x^{m(n-1)/n}}\\
&\simeq&\frac{mx^{m-1}t}{nx^{m(n-1)/n}}\\
&=&\frac{m}{n}\cdot x^{m/n-1}\cdot t
\end{eqnarray*}

But then, the general case, $p\in\mathbb R$, follows too, via a continuity argument.
\end{proof}

There are many other computations that can be done, and we will be back to this later. Now back to the general level, let us record here the following key result:

\index{Leinbitz rule}
\index{chain rule}

\begin{theorem}
The derivatives are subject to the following rules:
\begin{enumerate}
\item Leibnitz rule: $(fg)'=f'g+fg'$.

\item Chain rule: $(f\circ g)'=f'(g)g'$.
\end{enumerate}
\end{theorem}

\begin{proof}
Both formulae follow from the definition of the derivative, as follows:

\medskip

(1) Regarding products, we have the following computation:
\begin{eqnarray*}
(fg)(x+t)
&=&f(x+t)g(x+t)\\
&\simeq&(f(x)+f'(x)t)(g(x)+g'(x)t)\\
&\simeq&f(x)g(x)+(f'(x)g(x)+f(x)g'(x))t
\end{eqnarray*}

(2) Regarding compositions, we have the following computation:
\begin{eqnarray*}
(f\circ g)(x+t)
&=&f(g(x+t))\\
&\simeq&f(g(x)+g'(x)t)\\
&\simeq&f(g(x))+f'(g(x))g'(x)t
\end{eqnarray*}

Thus, we are led to the conclusions in the statement.
\end{proof}

There are many applications of the derivative, summarized as follows:

\index{local minimum}
\index{local maximum}
\index{Rolle theorem}
\index{mean value theorem}
\index{vanishing derivative}

\begin{theorem}
Given a differentiable function $f:[a,b]\to\mathbb R$, we have:
\begin{enumerate}
\item The local minima and maxima of $f$ appear at the points where $f'(x)=0$.

\item Rolle theorem: if $f(a)=f(b)$, we must have $f'(c)=0$, for some $c\in(a,b)$. 

\item Mean value theorem: $\frac{f(b)-f(a)}{b-a}=f'(c)$, for some $c\in(a,b)$.

\item Main theorem: if $f'=0$ then $f$ must be constant.
\end{enumerate}
\end{theorem}

\begin{proof}
Here (1) is clear from $f(x+t)\simeq f(x)+f'(x)t$, then $(1)\implies(2)$ is clear too, and then $(3)$ comes from (2), applied to the following function:
$$g(x)=f(x)-\frac{f(b)-f(a)}{b-a}\cdot x$$

As for (4), which is extremely useful in practice, this follows from (3).
\end{proof}

At a more advanced level now, we can talk about second derivatives, and we have:

\index{Taylor formula}

\begin{theorem}
Any twice differentiable $f:\mathbb R\to\mathbb R$ is approximately locally quadratic,
$$f(x+t)\simeq f(x)+f'(x)t+\frac{f''(x)}{2}\,t^2$$
with $f''(x)$ being the derivative of the function $f':\mathbb R\to\mathbb R$ at the point $x$.
\end{theorem}

\begin{proof}
This is something quite intuitive, when thinking geometrically. In practice, we can use L'H\^opital's rule, stating that the $0/0$ type limits can be computed as:
$$\frac{f(x)}{g(x)}\simeq\frac{f'(x)}{g'(x)}$$

Observe that this formula holds indeed, as an application of Theorem 5.1. Now by using this, if we denote by $\varphi(t)\simeq P(t)$ the formula to be proved, we have:
\begin{eqnarray*}
\frac{\varphi(t)-P(t)}{t^2}
&\simeq&\frac{\varphi'(t)-P'(t)}{2t}\\
&\simeq&\frac{\varphi''(t)-P''(t)}{2}\\
&=&\frac{f''(x)-f''(x)}{2}\\
&=&0
\end{eqnarray*}

Thus, we are led to the conclusion in the statement.
\end{proof}

The above result substantially improves Theorem 5.1, and there are many applications of it. We can improve for instance Theorem 5.4 (1), as follows:

\index{local minimum}
\index{local maximum}

\begin{theorem}
The local extrema of a twice differentiable function $f:\mathbb R\to\mathbb R$ appear at the points $x\in\mathbb R$ where $f'(x)=0$, as follows:
\begin{enumerate}
\item If $f''(x)>0$ we have a local minimum.

\item If $f''(x)<0$ we have a local maximum.

\item If $f''(x)=0$ things are undetermined.
\end{enumerate}
\end{theorem}

\begin{proof}
The first assertion is something that we already know. As for the second assertion, we can use the formula in Theorem 5.5, which in the case $f'(x)=0$ reads:
$$f(x+t)\simeq f(x)+\frac{f''(x)}{2}\,t^2$$

Indeed, assuming $f''(x)\neq 0$, it is clear that the condition $f''(x)>0$ will produce a local minimum, and that the condition $f''(x)<0$ will produce a local maximum.
\end{proof}

We can further develop the above method, at order 3, at order 4, and so on, the ultimate result on the subject, called Taylor formula, being as follows:

\index{Taylor formula}

\begin{theorem}
Assuming that $f:\mathbb R\to\mathbb R$ is $n$ times differentiable, we have
$$f(x+t)\simeq\sum_{k=0}^n\frac{f^{(k)}(x)}{k!}\,t^k$$
where $f^{(k)}(x)$ are the higher derivatives of $f$ at the point $x$.
\end{theorem}

\begin{proof}
We use the same method as in the proof of Theorem 5.5. Indeed, if we denote by $\varphi(t)\simeq P(t)$ the approximation to be proved, we have:
\begin{eqnarray*}
\frac{\varphi(t)-P(t)}{t^n}
&\simeq&\frac{\varphi'(t)-P'(t)}{nt^{n-1}}\\
&\simeq&\frac{\varphi''(t)-P''(t)}{n(n-1)t^{n-2}}\\
&\vdots&\\
&\simeq&\frac{\varphi^{(n)}(t)-P^{(n)}(t)}{n!}\\
&=&0
\end{eqnarray*}

Thus, we are led to the conclusion in the statement.
\end{proof}

As a basic application of derivatives and the Taylor formula, we have:

\index{Taylor formula}
\index{trigonometric functions}
\index{exp and log}

\begin{theorem}
We have the following formulae,
$$\sin t=\sum_{l=0}^\infty(-1)^l\frac{t^{2l+1}}{(2l+1)!}\quad,\quad 
\cos t=\sum_{l=0}^\infty(-1)^l\frac{t^{2l}}{(2l)!}$$
as well as the following formulae,
$$e^t=\sum_{k=0}^\infty\frac{t^k}{k!}\quad,\quad 
\log(1+t)=\sum_{k=0}^\infty(-1)^{k+1}\frac{t^k}{k}$$
as Taylor series, and in general as well, with $|t|<1$ needed for $\log$.
\end{theorem}

\begin{proof}
There are several statements here, the proofs being as follows:

\medskip

(1) Regarding $\sin$ and $\cos$, we can use here the following well-known formulae:
$$\sin(x+t)=\sin x\cos t+\cos x\sin t$$
$$\cos(x+t)=\cos x\cos t-\sin x\sin t$$

With these formulae in hand we can appproximate both $\sin$ and $\cos$, and we get:
$$(\sin x)'=\cos x\quad,\quad 
(\cos x)'=-\sin x$$

Thus, we can differentiate $\sin$ and $\cos$ as many times as we want to, and so we can compute the corresponding Taylor series, and we obtain the formulae in the statement.

\medskip

(2) Regarding $\exp$ and $\log$, here the needed formulae, which lead to the formulae in the statement for the corresponding Taylor series, are as follows:
$$(e^x)'=e^x\quad,\quad 
(\log x)'=x^{-1}\quad,\quad
(x^p)'=px^{p-1}$$

(3) Finally, the fact that the Taylor formulae in the statement are exact, and extend beyond the small $t$ setting, is something standard too. Indeed, for $\exp$ this is clear, for $\sin$, $\cos$ this is something that we know from chapter 3, coming from the Euler formula, and for $\log$ this is something which follows from some standard computations.
\end{proof}

As another basic application of derivatives and the Taylor formula, we have:

\index{binomial formula}
\index{generalized binomial formula}

\begin{theorem}
We have the generalized binomial formula
$$(1+t)^p=\sum_{k=0}^\infty\binom{p}{k}t^k$$
with the generalized binomial coefficients being given by
$$\binom{p}{k}=\frac{p(p-1)\ldots(p-k+1)}{k!}$$
for any $p\in\mathbb R$, and any $|t|<1$. With $p\in\mathbb N$, we recover the usual binomial formula.
\end{theorem}

\begin{proof}
As before with the various functions in Theorem 5.8, the Taylor series assertion is clear. Regarding now the fact that the formula is indeed exact, and extends beyond the small $t$ setting, if $f$ is the series in the statement, we have:
$$(1+t)f'(t)=pf(t)$$

Now by using this formula, we have the following computation:
$$\left((1+t)^{-p}f(t)\right)'=-p(1+t)^{-p-1}f(t)+(1+t)^{-p}f'(t)=0$$

Thus we have $f(t)=c(1+t)^p$, with $c=f(0)=1$, as desired.
\end{proof}

As a main application of the above formula, we can now extract square roots:

\index{square root}
\index{Catalan numbers}
\index{central binomial coefficients}

\begin{theorem}
We have the following formula,
$$\sqrt{1+t}=1-2\sum_{k=1}^\infty C_{k-1}\left(\frac{-t}{4}\right)^k$$
with $C_k=\frac{1}{k+1}\binom{2k}{k}$ being the Catalan numbers. Also, we have
$$\frac{1}{\sqrt{1+t}}=\sum_{k=0}^\infty D_k\left(\frac{-t}{4}\right)^k$$
with $D_k=\binom{2k}{k}$ being the central binomial coefficients.
\end{theorem}

\begin{proof}
At $p=1/2$, the generalized binomial coefficients are:
\begin{eqnarray*}
\binom{1/2}{k}
&=&\frac{1/2(-1/2)\ldots(3/2-k)}{k!}\\
&=&(-1)^{k-1}\frac{1\cdot 3\cdot 5\ldots(2k-3)}{2^kk!}\\
&=&(-1)^{k-1}\frac{(2k-2)!}{2^{k-1}(k-1)!2^kk!}\\
&=&-2\left(\frac{-1}{4}\right)^kC_{k-1}
\end{eqnarray*}

At $p=-1/2$, the generalized binomial coefficients are:
\begin{eqnarray*}
\binom{-1/2}{k}
&=&\frac{-1/2(-3/2)\ldots(1/2-k)}{k!}\\
&=&(-1)^k\frac{1\cdot 3\cdot 5\ldots(2k-1)}{2^kk!}\\
&=&(-1)^k\frac{(2k)!}{2^kk!2^kk!}\\
&=&\left(\frac{-1}{4}\right)^kD_k
\end{eqnarray*}

Thus, we obtain the formulae in the statement.
\end{proof}

Let us discuss as well the basics of integration theory. We first have:

\index{Riemann integration}

\begin{definition}
We have the Riemann integration formula,
$$\int_a^bf(x)dx=\lim_{N\to\infty}\sum_{k=1}^N\frac{b-a}{N}\times f\left(a+\frac{b-a}{N}\cdot k\right)$$
which can serve as a formal definition for the integral.
\end{definition}

To be more precise, given a continuous function $f:[a,b]\to\mathbb R$, we can try to compute the signed area below its graph, called integral and denoted $\int_a^bf(x)dx$, and by approximating with rectangles, in the obvious way, we are led to the Riemann formula. 

\bigskip

As an illustration for this, with some arithmetic know-how, for the computation of sums of type $1^p+2^p+\ldots+N^p$, we have the following formula, for $p\in\mathbb N$:
$$\int_0^1x^pdx=
\lim_{N\to\infty}\frac{1^p+2^p+\ldots+N^p}{N^{p+1}}=\frac{1}{p+1}$$

However, such things remain a bit amateurish. At the more advanced level, the point is that the derivatives and integrals are related in several subtle ways, as follows:

\index{change of variable}
\index{partial integration}

\begin{theorem}
We have the following formulae, called fundamental theorem of calculus, integration by parts formula, and change of variable formula,
$$\int_a^bF'(x)dx=\Big[F\Big]_a^b$$
$$\int_a^b(f'g+fg')(x)dx=\Big[fg\Big]_a^b$$
$$\int_a^bf(x)dx=\int_{\varphi^{-1}(a)}^{\varphi^{-1}(b)}f(\varphi(t))\varphi'(t)dt$$
with the convention $[F]_a^b=F(b)-F(a)$, for the first two formulae.
\end{theorem}

\begin{proof}
To start with, given a continuous function $f:[a,b]\to\mathbb R$, by integrating $\min f\leq f\leq \max f$ we obtain the following formula, called mean value property:
$$\exists c\in[a,b]\quad,\quad \int_a^bf(x)dx=(b-a)f(c)$$

Next, this mean value property shows that we have the following implication:
$$I(x)=\int_a^xf(s)ds\implies I'=f$$

Now given $F:\mathbb R\to\mathbb R$ as in the statement, by using this with $f=F'$, we obtain $I'=F'$. Since $I(a)=0$, this reads $F(x)=I(x)+F(a)$, and with $x=b$ we get:
$$F(b)=\int_a^bF'(x)dx+F(a)$$

Thus, first formula proved, and the second and third formulae follow as well.
\end{proof}

\section*{5b. Several variables}

Let us discuss now what happens in several variables. At order 1, we haves:

\index{partial derivatives}
\index{continuously differentiable}

\begin{theorem}
A function $f:\mathbb R^N\to\mathbb R^M$ is continuously differentiable,
$$f(x+t)\simeq f(x)+f'(x)t$$
with $f'(x)$ linear, and $x\to f'(x)$ continuous, precisely when it has partial derivatives,
$$\frac{df_i}{dx_j}(x)=\lim_{t\to 0}\frac{f_i(x+te_j)-f_i(x)}{t}$$
which depend continuously on $x$. In this case the derivative is
$$f'(x)=\left(\frac{df_i}{dx_j}(x)\right)_{ij}\in M_{M\times N}(\mathbb R)$$ 
acting on the vectors $t\in\mathbb R^N$ by usual multiplication.
\end{theorem}

\begin{proof}
The formula in the statement makes sense indeed, as follows:
$$f\begin{pmatrix}x_1+t_1\\ \vdots\\ x_N+t_N\end{pmatrix}
\simeq f\begin{pmatrix}x_1\\ \vdots\\ x_N\end{pmatrix}
+\begin{pmatrix}
\frac{df_1}{dx_1}(x)&\ldots&\frac{df_1}{dx_N}(x)\\
\vdots&&\vdots\\
\frac{df_M}{dx_1}(x)&\ldots&\frac{df_M}{dx_N}(x)
\end{pmatrix}\begin{pmatrix}t_1\\ \vdots\\ t_N\end{pmatrix}$$

Getting now to the proof of this formula, this goes as follows:

\medskip

(1) First of all, at $N=M=1$ what we have is a usual 1-variable function $f:\mathbb R\to\mathbb R$, and the formula in the statement is something that we know well, namely:
$$f(x+t)\simeq f(x)+f'(x)t$$

(2) Let us discuss now the case $N=2,M=1$. Here what we have is a function $f:\mathbb R^2\to\mathbb R$, and by using twice the basic approximation result from (1), we obtain:
\begin{eqnarray*}
f\binom{x_1+t_1}{x_2+t_2}
&\simeq&f\binom{x_1+t_1}{x_2}+\frac{df}{dx_2}(x)t_2\\
&\simeq&f\binom{x_1}{x_2}+\frac{df}{dx_1}(x)t_1+\frac{df}{dx_2}(x)t_2\\
&=&f\binom{x_1}{x_2}+\begin{pmatrix}\frac{df}{dx_1}(x)&\frac{df}{dx_2}(x)\end{pmatrix}\binom{t_1}{t_2}
\end{eqnarray*}

(3) More generally, we can deal in this way with the general case $M=1$, with the formula here, obtained via a straightforward recurrence, being as follows:
\begin{eqnarray*}
f\begin{pmatrix}x_1+t_1\\ \vdots\\ x_N+t_N\end{pmatrix}
&\simeq&f\begin{pmatrix}x_1\\ \vdots\\ x_N\end{pmatrix}+\frac{df}{dx_1}(x)t_1+\ldots+\frac{df}{dx_N}(x)t_N\\
&=&f\begin{pmatrix}x_1\\ \vdots\\ x_N\end{pmatrix}+
\begin{pmatrix}\frac{df}{dx_1}(x)&\ldots&\frac{df}{dx_N}(x)\end{pmatrix}
\begin{pmatrix}t_1\\ \vdots\\ t_N\end{pmatrix}
\end{eqnarray*}

(4) But this gives the result in the case where both $N,M\in\mathbb N$ are arbitrary too. Indeed, consider a function $f:\mathbb R^N\to\mathbb R^M$, and let us write it as follows:
$$f=\begin{pmatrix}f_1\\ \vdots\\ f_M\end{pmatrix}$$

We can apply (3) to each of the components $f_i:\mathbb R^N\to\mathbb R$, and we get:
$$f_i\begin{pmatrix}x_1+t_1\\ \vdots\\ x_N+t_N\end{pmatrix}
\simeq f_i\begin{pmatrix}x_1\\ \vdots\\ x_N\end{pmatrix}+
\begin{pmatrix}\frac{df_i}{dx_1}(x)&\ldots&\frac{df_i}{dx_N}(x)\end{pmatrix}
\begin{pmatrix}t_1\\ \vdots\\ t_N\end{pmatrix}$$

(5) But this collection of $M$ formulae tells us precisely that the following happens, as an equality, or rather approximation, of vectors in $\mathbb R^M$:
$$f\begin{pmatrix}x_1+t_1\\ \vdots\\ x_N+t_N\end{pmatrix}
\simeq f\begin{pmatrix}x_1\\ \vdots\\ x_N\end{pmatrix}
+\begin{pmatrix}
\frac{df_1}{dx_1}(x)&\ldots&\frac{df_1}{dx_N}(x)\\
\vdots&&\vdots\\
\frac{df_M}{dx_1}(x)&\ldots&\frac{df_M}{dx_N}(x)
\end{pmatrix}\begin{pmatrix}t_1\\ \vdots\\ t_N\end{pmatrix}$$

Thus, we are led to the conclusion in the statement.
\end{proof}

Generally speaking, Theorem 5.13 is what we need to know for upgrading from calculus to multivariable calculus. As a standard result here, we have:

\index{chain rule}

\begin{theorem}
We have the chain derivative formula
$$(f\circ g)'(x)=f'(g(x))\cdot g'(x)$$
as an equality of matrices.
\end{theorem}

\begin{proof}
Consider indeed a composition of functions, as follows:
$$f:\mathbb R^N\to\mathbb R^M\quad,\quad 
g:\mathbb R^K\to\mathbb R^N\quad,\quad 
f\circ g:\mathbb R^K\to\mathbb R^M$$

According to Theorem 5.13, the derivatives of these functions are certain linear maps, corresponding to certain rectangular matrices, as follows:
$$f'(g(x))\in M_{M\times N}(\mathbb R)\quad,\quad 
g'(x)\in M_{N\times K}(\mathbb R)\quad\quad
(f\circ g)'(x)\in M_{M\times K}(\mathbb R)$$

Thus, our formula makes sense indeed. As for proof, this comes from:
\begin{eqnarray*}
(f\circ g)(x+t)
&=&f(g(x+t))\\
&\simeq&f(g(x)+g'(x)t)\\
&\simeq&f(g(x))+f'(g(x))g'(x)t
\end{eqnarray*}

Thus, we are led to the conclusion in the statement.
\end{proof}

Next, we can talk about higher derivatives, in the obvious way, simply by performing the operation of taking derivatives recursively. To be more precise, we have:

\index{higher derivative}
\index{Clairaut formula}

\begin{theorem}
Given $f:\mathbb R^N\to\mathbb R$, we can talk about its higher derivatives
$$\frac{d^kf}{dx_{i_1}\ldots dx_{i_k}}=\frac{d}{dx_{i_1}}\cdots\frac{d}{dx_{i_k}}(f)$$
provided that these derivatives exist indeed. Moreover, due to the Clairaut formula,
$$\frac{d^2f}{dx_idx_j}=\frac{d^2f}{dx_jdx_i}$$
the order in which these higher derivatives are computed is irrelevant.
\end{theorem}

\begin{proof}
There are several things going on here, the idea being as follows:

\medskip

(1) First of all, we can talk about the quantities in the statement, with the remark of course that at each step of our recursion, the corresponding partial derivative can exist of not. We will say in what follows that our function is $n$ times differentiable if the quantities in the statement exist at any $k\leq n$, and smooth, if this works with $n=\infty$.

\medskip

(2) Regarding the second assertion, this is self-explanatory, based on the Clairaut formula, which is something elementary, coming from the mean value theorem.

\medskip

(3) In practice now, we can permute the order of our partial derivative computations, and a standard way of doing this is by differentiating first with respect to $x_1$, as many times as needed, then with respect to $x_2$, and so on. Thus, the collection of partial derivatives can be written, in a more convenient form, as follows:
$$\frac{d^kf}{dx_1^{k_1}\ldots dx_N^{k_N}}=\frac{d^{k_1}}{dx_1^{k_1}}\cdots\frac{d^{k_N}}{dx_N^{k_N}}(f)$$

(4) To be more precise, here $k\in\mathbb N$ is as usual the global order of our derivatives, the exponents $k_1,\ldots,k_N\in\mathbb N$ are subject to the condition $k_1+\ldots+k_N=k$, and the operations on the right are the familiar one-variable higher derivative operations.
\end{proof}

Regarding now the Taylor formula, in several variables, at order 2, we have:

\index{local minimum}
\index{local maximum}
\index{Hessian matrix}

\begin{theorem}
Given a function $f:\mathbb R^N\to\mathbb R$, construct its Hessian, as being:
$$f''(x)=\left(\frac{d^2f}{dx_idx_j}(x)\right)_{ij}$$ 
We have then the following order $2$ approximation of $f$ around a given $x\in\mathbb R^N$,
$$f(x+t)\simeq f(x)+f'(x)t+\frac{<f''(x)t,t>}{2}$$
relating the positivity properties of $f''$ to the local minima and maxima of $f$.
\end{theorem}

\begin{proof}
This is something very standard, the idea being as follows:

\medskip

(1) At $N=1$ the Hessian matrix is the $1\times1$ matrix having as entry the usual $f''(x)$, and the formula in the statement is something that we know well, namely:
$$f(x+t)\simeq f(x)+f'(x)t+\frac{f''(x)t^2}{2}$$

(2) In general, our claim is that the formula in the statement follows from the one-variable formula above, applied to the restriction of $f$ to the following segment in $\mathbb R^N$:
$$I=[x,x+t]$$

To be more precise, let $y\in\mathbb R^N$, and consider the following function, with $r\in\mathbb R$:
$$g(r)=f(x+ry)$$

We know from (1) that the Taylor formula for $g$, at the point $r=0$, reads:
$$g(r)\simeq g(0)+g'(0)r+\frac{g''(0)r^2}{2}$$

And our claim is that, with $t=ry$, this is precisely the formula in the statement.

\medskip

(3) So, let us see if our claim is correct. By using the chain rule, we have:
$$g'(r)=f'(x+ry)\cdot y$$

By using again the chain rule, we can compute the second derivative as well:
\begin{eqnarray*}
g''(r)
&=&(f'(x+ry)\cdot y)'\\
&=&\left(\sum_i\frac{df}{dx_i}(x+ry)\cdot y_i\right)'\\
&=&\sum_i\sum_j\frac{d^2f}{dx_idx_j}(x+ry)\cdot\frac{d(x+ry)_j}{dr}\cdot y_i\\
&=&\sum_i\sum_j\frac{d^2f}{dx_idx_j}(x+ry)\cdot y_iy_j\\
&=&<f''(x+ry)y,y>
\end{eqnarray*}

(4) Time now to conclude. We know that we have $g(r)=f(x+ry)$, and according to our various computations above, we have the following formulae:
$$g(0)=f(x)\quad,\quad 
g'(0)=f'(x)\quad,\quad 
g''(0)=<f''(x)y,y>$$

Buit with this data in hand, the usual Taylor formula for our one variable function $g$, at order 2, at the point $r=0$, takes the following form, with $t=ry$:
\begin{eqnarray*}
f(x+ry)
&\simeq&f(x)+f'(x)ry+\frac{<f''(x)y,y>r^2}{2}\\
&=&f(x)+f'(x)t+\frac{<f''(x)t,t>}{2}
\end{eqnarray*}

Thus, we have obtained the formula in the statement. Finally, the last assertion, regarding the local extrema, is something standard, as in the one-variable case.
\end{proof}

As a complement to Theorem 5.16, very useful in practice, let us record:

\index{Hessian eigenvalues}

\begin{theorem}
Given a twice differentiable function $f:\mathbb R^N\to\mathbb R$, assume that $f'(x)=0$, and let $\lambda_1,\ldots,\lambda_N$ be the eigenvalues of $f''(x)$. Then:
\begin{enumerate}
\item $\lambda_i\geq0$ is needed for $x$ to be a local minimum.

\item $\lambda_i>0$ guarantees that $x$ is a local minimum.

\item $\lambda_i\leq0$ is needed for $x$ to be a local maximum.

\item $\lambda_i<0$ guarantees that $x$ is a local maximum.
\end{enumerate}
\end{theorem}

\begin{proof}
This comes from Theorem 5.16 and from linear algebra, as follows:

\medskip

(1) We know from chapter 4 that the Hessian matrix $f''(x)$, which is symmetric, is diagonalized by a certain matrix $U\in O_N$. But with this in hand, we can change the basis of $\mathbb R^N$, with the help of this matrix $U\in O_N$, and the Taylor formula becomes:
$$f(x+t)\simeq f(x)+\sum_{i=1}^N\lambda_it_i^2$$

And this latter formula, obviously, gives all the assertions in the statement.

\medskip

(2) This was for the theory, but in practice, there are some other things that can be useful. Consider for instance a function $f:\mathbb R^2\to\mathbb R$, whose Hessian looks as follows:
$$f''(x)=\begin{pmatrix}a&b\\ c&d\end{pmatrix}$$

The eigenvalues are then given by the following trace and determinant equations:
$$\lambda_1+\lambda_2=a+d\quad,\quad\lambda_1\lambda_2=ad-bc$$

Thus, without even computing the eigenvalues, we can say right away, depending on the signs of $a+d$, $ad-bc$, if we are in one of the situations (1,2,3,4) in the statement.

\medskip

(3) In more dimensions things are more complicated, but there are still tricks, that can help, and the more you learn and know here, the better your analysis will be.
\end{proof}

\section*{5c. Multiple integrals}

Getting now to integration matters, in several variables, we certainly have an analogue of Definition 5.11, and we can usually compute the multiple integrals by iterating one-variable integrals. At the theoretical level, as a key result here, we have:

\index{change of variable}
\index{Jacobian}
\index{multiple integral}

\begin{theorem}
Given a transformation $\varphi=(\varphi_1,\ldots,\varphi_N)$, we have
$$\int_Ef(x)dx=\int_{\varphi^{-1}(E)}f(\varphi(t))|J_\varphi(t)|dt$$
with the $J_\varphi$ quantity, called Jacobian, being given by
$$J_\varphi(t)=\det\left[\left(\frac{d\varphi_i}{dx_j}(t)\right)_{ij}\right]$$ 
and with this generalizing the $1$-variable formula that we know well.
\end{theorem}

\begin{proof}
This is something quite tricky, the idea being as follows:

\medskip

(1) Observe first that this generalizes indeed the change of variable formula in 1 dimension, from Theorem 5.12, the point here being that the absolute value on the derivative appears as to compensate for the lack of explicit bounds for the integral.

\medskip 

(2) In general now, we can first argue that, the formula in the statement being linear in $f$, we can assume $f=1$. Thus we want to prove $vol(E)=\int_{\varphi^{-1}(E)}|J_\varphi(t)|dt$, and with $D={\varphi^{-1}(E)}$, this amounts in proving $vol(\varphi(D))=\int_D|J_\varphi(t)|dt$.

\medskip

(3) Now since this latter formula is additive with respect to $D$, it is enough to prove that $vol(\varphi(D))=\int_D J_\varphi(t)dt$, for small cubes $D$, and assuming $J_\varphi>0$. But for $\varphi$ linear this follows by using the definition of the determinant as a volume, as in chapter 2.

\medskip

(4) In order to prove now the theorem, as stated, let us rather focus on the transformations used $\varphi$, instead of the functions to be integrated $f$. Our first claim is that the validity of the theorem is stable under taking compositions of such transformations $\varphi$.

\medskip

(5) In order to prove this claim, consider a composition, as follows:
$$\varphi:E\to F\quad,\quad 
\psi:D\to E\quad,\quad 
\varphi\circ\psi:D\to F$$

Assuming that the theorem holds for $\varphi,\psi$, we have the following computation:
\begin{eqnarray*}
\int_Ff(x)dx
&=&\int_Ef(\varphi(s))|J_\varphi(s)|ds\\
&=&\int_Df(\varphi\circ\psi(t))|J_\varphi(\psi(t))|\cdot|J_\psi(t)|dt\\
&=&\int_Df(\varphi\circ\psi(t))|J_{\varphi\circ\psi}(t)|dt
\end{eqnarray*}

Thus, our theorem holds as well for $\varphi\circ\psi$, and we have proved our claim.

\medskip

(6) Next, as a key ingredient, let us examine the case where we are in $N=2$ dimensions, and our transformation $\varphi$ has one of the following special forms:
$$\varphi(x,y)=(\psi(x,y),y)\quad,\quad\varphi(x,y)=(x,\psi(x,y))$$

By symmetry, it is enough to deal with the first case. Here the Jacobian is $d\psi/dx$, and by replacing if needed $\psi\to-\psi$, we can assume that this Jacobian is positive, $d\psi/dx>0$. Now by assuming as before that $D=\varphi^{-1}(E)$ is a rectangle, $D=[a,b]\times[c,d]$, we can prove our formula by using the change of variables in 1 dimension, as follows:
\begin{eqnarray*}
\int_Ef(s)ds
&=&\int_{\varphi(D)}f(x,y)dxdy\\
&=&\int_c^d\int_{\psi(a,y)}^{\psi(b,y)}f(x,y)dxdy\\
&=&\int_c^d\int_a^bf(\psi(x,y),y)\frac{d\psi}{dx}\,dxdy\\
&=&\int_Df(\varphi(t))J_\varphi(t)dt
\end{eqnarray*}

(7) But with this, we can now prove the theorem, in $N=2$ dimensions. Indeed, given a transformation $\varphi=(\varphi_1,\varphi_2)$, consider the following two transformations:
$$\phi(x,y)=(\varphi_1(x,y),y)\quad,\quad \psi(x,y)=(x,\varphi_2\circ\phi^{-1}(x,y))$$

We have then $\varphi=\psi\circ\phi$, and by using (6) for $\psi,\phi$, which are of the special form there, and then (5) for composing, we conclude that the theorem holds for $\varphi$, as desired.

\medskip

(8) Thus, theorem proved in $N=2$ dimensions, at least in the generic situation, and we will leave the remaining details as an exercise. And the extension of the above proof to arbitrary $N$ dimensions is straightforward, that we will leave as an exercise too.
\end{proof}

We can discuss now some more advanced questions, related to the computation of volumes of the spheres, and to the integration over spheres. Let us start with:

\index{polar coordinates}
\index{Jacobian}

\begin{theorem}
We have polar coordinates in $2$ dimensions,
$$\begin{cases}
x\!\!\!&=\ r\cos t\\
y\!\!\!&=\ r\sin t
\end{cases}$$
the corresponding Jacobian being $J=r$.
\end{theorem}

\begin{proof}
This is something elementary, the Jacobian being given by:
\begin{eqnarray*}
J
&=&\begin{vmatrix}
\cos t&-r\sin t\\
\sin t&r\cos t
\end{vmatrix}\\
&=&r\cos^2t+r\sin^2t\\
&=&r
\end{eqnarray*}

Thus, we have indeed the formula in the statement.
\end{proof}

We can now compute the Gauss integral, which is the best calculus formula ever:

\index{Gauss integral}

\begin{theorem}
We have the following formula,
$$\int_\mathbb Re^{-x^2}dx=\sqrt{\pi}$$
called Gauss integral formula.
\end{theorem}

\begin{proof}
This is something truly magic, the idea being as follows:

\medskip

(1) To start with, we can certainly integrate $e^{-x^2}$ by using the formula of the exponential series, and the primitive which is worth 0 at $x=0$ is given by:
$$\int e^{-x^2}=\sum_{k=0}^\infty(-1)^k\frac{x^{2k+1}}{(2k+1)k!}$$

However, this series is not computable, in terms of the known, familiar series.

\medskip

(2) Next, we can still ask for the computation of $\int_\mathbb Re^{-x^2}dx$, who knows. And here, another surprise awaits us, this is undoable, with bare hands. However, and here comes the magic, the Gauss integral can be computed by using two dimensions, as follows:
\begin{eqnarray*}
\left(\int_\mathbb Re^{-x^2}dx\right)^2
&=&\int_\mathbb R\int_\mathbb Re^{-x^2-y^2}dxdy\\
&=&\int_0^{2\pi}\int_0^\infty e^{-r^2}rdrdt\\
&=&2\pi\int_0^\infty\left(-\frac{e^{-r^2}}{2}\right)'dr\\
&=&2\pi\left[0-\left(-\frac{1}{2}\right)\right]\\
&=&\pi
\end{eqnarray*}

(3) Amazing, all this. We will heavily use the Gauss integral, in what follows.
\end{proof}

Getting now to 3 dimensions, we have here the following result:

\index{spherical coordinates}

\begin{theorem}
We have spherical coordinates in $3$ dimensions,
$$\begin{cases}
x\!\!\!&=\ r\cos s\\
y\!\!\!&=\ r\sin s\cos t\\
z\!\!\!&=\ r\sin s\sin t
\end{cases}$$
the corresponding Jacobian being $J(r,s,t)=r^2\sin s$.
\end{theorem}

\begin{proof}
The fact that we have indeed spherical coordinates is clear. Regarding now the Jacobian, this is given by the following formula:
\begin{eqnarray*}
J(r,s,t)
&=&\begin{vmatrix}
\cos s&-r\sin s&0\\
\sin s\cos t&r\cos s\cos t&-r\sin s\sin t\\
\sin s\sin t&r\cos s\sin t&r\sin s\cos t
\end{vmatrix}\\
&=&r^2\sin s\sin t
\begin{vmatrix}\cos s&-r\sin s\\ \sin s\sin t&r\cos s\sin t\end{vmatrix}
+r\sin s\cos t\begin{vmatrix}\cos s&-r\sin s\\ \sin s\cos t&r\cos s\cos t\end{vmatrix}\\
&=&r\sin s\sin^2 t
\begin{vmatrix}\cos s&-r\sin s\\ \sin s&r\cos s\end{vmatrix}
+r\sin s\cos^2 t\begin{vmatrix}\cos s&-r\sin s\\ \sin s&r\cos s\end{vmatrix}\\
&=&r\sin s(\sin^2t+\cos^2t)\begin{vmatrix}\cos s&-r\sin s\\ \sin s&r\cos s\end{vmatrix}\\
&=&r\sin s\times 1\times r\\
&=&r^2\sin s
\end{eqnarray*}

Thus, we have indeed the formula in the statement.
\end{proof}

Let us work out now the spherical coordinate formula in $N$ dimensions. The result here, which generalizes those at $N=2,3$, is as follows:

\index{spherical coordinates}
\index{Jacobian}

\begin{theorem}
We have spherical coordinates in $N$ dimensions,
$$\begin{cases}
x_1\!\!\!&=\ r\cos t_1\\
x_2\!\!\!&=\ r\sin t_1\cos t_2\\
\vdots\\
x_{N-1}\!\!\!&=\ r\sin t_1\sin t_2\ldots\sin t_{N-2}\cos t_{N-1}\\
x_N\!\!\!&=\ r\sin t_1\sin t_2\ldots\sin t_{N-2}\sin t_{N-1}
\end{cases}$$
the Jacobian being $J(r,t)=r^{N-1}\sin^{N-2}t_1\sin^{N-3}t_2\,\ldots\,\sin^2t_{N-3}\sin t_{N-2}$.
\end{theorem}

\begin{proof}
As before, the fact that we have spherical coordinates is clear. Regarding now the Jacobian, also as before, by developing over the last column, we have:
\begin{eqnarray*}
J_N
&=&r\sin t_1\ldots\sin t_{N-2}\sin t_{N-1}\times \sin t_{N-1}J_{N-1}\\
&+&r\sin t_1\ldots \sin t_{N-2}\cos t_{N-1}\times\cos t_{N-1}J_{N-1}\\
&=&r\sin t_1\ldots\sin t_{N-2}(\sin^2 t_{N-1}+\cos^2 t_{N-1})J_{N-1}\\
&=&r\sin t_1\ldots\sin t_{N-2}J_{N-1}
\end{eqnarray*}

Thus, we obtain the formula in the statement, by recurrence.
\end{proof}

As an application, let us compute now the volumes of spheres. For this purpose, we must understand how the products of coordinates integrate over spheres. Let us start with the case $N=2$. Here the sphere is the unit circle $\mathbb T$, and with $z=e^{it}$ the coordinates are $\cos t,\sin t$. We can first integrate arbitrary powers of these coordinates, as follows:

\index{double factorials}
\index{trigonometric integral}
\index{Wallis formula}

\begin{proposition}
We have the following formulae,
$$\int_0^{\pi/2}\cos^pt\,dt=\int_0^{\pi/2}\sin^pt\,dt=\left(\frac{\pi}{2}\right)^{\varepsilon(p)}\frac{p!!}{(p+1)!!}$$
where $\varepsilon(p)=1$ if $p$ is even, and $\varepsilon(p)=0$ if $p$ is odd, and where
$$m!!=(m-1)(m-3)(m-5)\ldots$$
with the product ending at $2$ if $m$ is odd, and ending at $1$ if $m$ is even.
\end{proposition}

\begin{proof}
Let us first compute the integral on the left $I_p$. We have:
\begin{eqnarray*}
(\cos^pt\sin t)'
&=&p\cos^{p-1}t(-\sin t)\sin t+\cos^pt\cos t\\
&=&p\cos^{p+1}t-p\cos^{p-1}t+\cos^{p+1}t\\
&=&(p+1)\cos^{p+1}t-p\cos^{p-1}t
\end{eqnarray*}

By integrating between $0$ and $\pi/2$, we obtain the following formula:
$$(p+1)I_{p+1}=pI_{p-1}$$

Thus we can compute $I_p$ by recurrence, and we obtain:
\begin{eqnarray*}
I_p
&=&\frac{p-1}{p}\,I_{p-2}\\
&=&\frac{p-1}{p}\cdot\frac{p-3}{p-2}\,I_{p-4}\\
&=&\frac{p-1}{p}\cdot\frac{p-3}{p-2}\cdot\frac{p-5}{p-4}\,I_{p-6}\\
&&\vdots\\
&=&\frac{p!!}{(p+1)!!}\,I_{1-\varepsilon(p)}
\end{eqnarray*}

Thus, we obtain the result, by recurrence. As for the second formula, regarding $\sin t$, this follows from the first formula, with the change of variables $t=\frac{\pi}{2}-s$.
\end{proof}

We can now compute the volumes of the spheres, as follows:

\index{volume of sphere}
\index{double factorial}

\begin{theorem}
The volume of the unit sphere in $\mathbb R^N$ is given by
$$V=\left(\frac{\pi}{2}\right)^{[N/2]}\frac{2^N}{(N+1)!!}$$
with the convention
$$N!!=(N-1)(N-3)(N-5)\ldots$$
with the product ending at $2$ if $N$ is odd, and ending at $1$ if $N$ is even.
\end{theorem}

\begin{proof}
If we denote by $B^+$ the positive part of the unit sphere, we have:
\begin{eqnarray*}
V^+
&=&\int_{B^+}1\\
&=&\int_0^1\int_0^{\pi/2}\ldots\int_0^{\pi/2}r^{N-1}\sin^{N-2}t_1\ldots\sin t_{N-2}\,drdt_1\ldots dt_{N-1}\\
&=&\int_0^1r^{N-1}\,dr\int_0^{\pi/2}\sin^{N-2}t_1\,dt_1\ldots\int_0^{\pi/2}\sin t_{N-2}dt_{N-2}\int_0^{\pi/2}1dt_{N-1}\\
&=&\frac{1}{N}\times\left(\frac{\pi}{2}\right)^{[N/2]}\times\frac{(N-2)!!}{(N-1)!!}\cdot\frac{(N-3)!!}{(N-2)!!}\ldots\frac{2!!}{3!!}\cdot\frac{1!!}{2!!}\cdot1\\
&=&\frac{1}{N}\times\left(\frac{\pi}{2}\right)^{[N/2]}\times\frac{1}{(N-1)!!}\\
&=&\left(\frac{\pi}{2}\right)^{[N/2]}\frac{1}{(N+1)!!}
\end{eqnarray*}

Thus, we are led to the formula in the statement.
\end{proof}

As main particular cases of the above formula, we have:

\index{volume of sphere}

\begin{proposition}
The volumes of the low-dimensional spheres are as follows:
\begin{enumerate}
\item At $N=1$, the length of the unit interval is $V=2$.

\item At $N=2$, the area of the unit disk is $V=\pi$.

\item At $N=3$, the volume of the unit sphere is $V=\frac{4\pi}{3}$

\item At $N=4$, the volume of the corresponding unit sphere is $V=\frac{\pi^2}{2}$.
\end{enumerate}
\end{proposition}

\begin{proof}
These are all particular cases of the formula in Theorem 5.24.
\end{proof}

\section*{5d. Stirling estimates}

The formula in Theorem 5.24 is certainly nice, but in practice, we would like to have estimates for that sphere volumes too. For this purpose, we will need:

\index{Stirling formula}
\index{Riemann sum}

\begin{theorem}
We have the Stirling formula
$$N!\simeq\left(\frac{N}{e}\right)^N\sqrt{2\pi N}$$
valid in the $N\to\infty$ limit.
\end{theorem}

\begin{proof}
This is something quite tricky, the idea being as follows:

\medskip

(1) Let us first see what we can get with Riemann sums. We have:
$$\log(N!)
=\sum_{k=1}^N\log k
\approx\int_1^N\log x\,dx
=N\log N-N+1$$

By exponentiating, this gives the following estimate, which is not bad:
$$N!\approx\left(\frac{N}{e}\right)^N\cdot e$$

(2) We can improve our estimate by replacing the rectangles from the Riemann sum approach to the integrals by trapezoids. In practice, this gives the following estimate:
$$\log(N!)
\approx\int_1^N\log x\,dx+\frac{\log 1+\log N}{2}
=N\log N-N+1+\frac{\log N}{2}$$

By exponentiating, this gives the following estimate, which gets us closer:
$$N!\approx\left(\frac{N}{e}\right)^N\cdot e\cdot\sqrt{N}$$

(3) In order to conclude, we must take some kind of mathematical magnifier, and carefully estimate the error made in (2). Fortunately, this mathematical magnifier exists, called Euler-Maclaurin formula, and after some computations, this leads to:
$$N!\simeq\left(\frac{N}{e}\right)^N\sqrt{2\pi N}$$

(4) However, all this remains a bit complicated, so we would like to present now an alternative approach to (3), which also misses some details, but better does the job, explaining where the $\sqrt{2\pi}$ factor comes from. First, by partial integration we have:
$$N!=\int_0^\infty x^Ne^{-x}dx$$

Since the integrand is sharply peaked at $x=N$, as you can see by computing the derivative of $\log(x^Ne^{-x})$, this suggests writing $x=N+y$, and we obtain:
\begin{eqnarray*}
\log(x^Ne^{-x})
&=&N\log x-x\\
&=&N\log(N+y)-(N+y)\\
&=&N\log N+N\log\left(1+\frac{y}{N}\right)-(N+y)\\
&\simeq&N\log N+N\left(\frac{y}{N}-\frac{y^2}{2N^2}\right)-(N+y)\\
&=&N\log N-N-\frac{y^2}{2N}
\end{eqnarray*}

By exponentiating, we obtain from this the following estimate:
$$x^Ne^{-x}\simeq\left(\frac{N}{e}\right)^Ne^{-y^2/2N}$$

(5) Now by integrating, and using the Gauss formula, we obtain from this:
\begin{eqnarray*}
N!
&=&\int_0^\infty x^Ne^{-x}dx\\
&\simeq&\int_{-N}^N\left(\frac{N}{e}\right)^Ne^{-y^2/2N}\,dy\\
&\simeq&\left(\frac{N}{e}\right)^N\int_\mathbb Re^{-y^2/2N}\,dy\\
&=&\left(\frac{N}{e}\right)^N\sqrt{2\pi N}
\end{eqnarray*}

Thus, we have proved the Stirling formula, as formulated in the statement.
\end{proof}

We can now estimate the volumes of the spheres, as follows:

\begin{theorem}
The volume of the unit sphere in $\mathbb R^N$ is given by
$$V\simeq\left(\frac{2\pi e}{N}\right)^{N/2}\frac{1}{\sqrt{\pi N}}$$
in the $N\to\infty$ limit.
\end{theorem}

\begin{proof}
This is very standard, using the formula in Theorem 5.24, as follows:

\medskip

(1) The double factorials can be estimated by using the Stirling formula. Indeed, in the case where $N=2K$ is even, we have the following computation:
\begin{eqnarray*}
(N+1)!!
&=&2^KK!\\
&\simeq&\left(\frac{2K}{e}\right)^K\sqrt{2\pi K}\\
&=&\left(\frac{N}{e}\right)^{N/2}\sqrt{\pi N}
\end{eqnarray*}

(2) As for the case where $N=2K-1$ is odd, here the estimate goes as follows:
\begin{eqnarray*}
(N+1)!!
&=&\frac{(2K)!}{2^KK!}\\
&\simeq&\frac{1}{2^K}\left(\frac{2K}{e}\right)^{2K}\sqrt{4\pi K}\left(\frac{e}{K}\right)^K\frac{1}{\sqrt{2\pi K}}\\
&=&\left(\frac{2K}{e}\right)^K\sqrt{2}\\
&=&\left(\frac{N+1}{e}\right)^{(N+1)/2}\sqrt{2}\\
&=&\left(\frac{N}{e}\right)^{N/2}\left(\frac{N+1}{N}\right)^{N/2}\sqrt{\frac{N+1}{e}}\cdot\sqrt{2}\\
&\simeq&\left(\frac{N}{e}\right)^{N/2}\sqrt{e}\cdot\sqrt{\frac{N}{e}}\cdot\sqrt{2}\\
&=&\left(\frac{N}{e}\right)^{N/2}\sqrt{2N}
\end{eqnarray*}

(3) Now back to the spheres, when $N$ is even, the estimate goes as follows:
\begin{eqnarray*}
V
&=&\left(\frac{\pi}{2}\right)^{N/2}\frac{2^N}{(N+1)!!}\\
&\simeq&\left(\frac{\pi}{2}\right)^{N/2}2^N\left(\frac{e}{N}\right)^{N/2}\frac{1}{\sqrt{\pi N}}\\
&=&\left(\frac{2\pi e}{N}\right)^{N/2}\frac{1}{\sqrt{\pi N}}
\end{eqnarray*}

(4) As for the case where $N$ is odd, here the estimate goes as follows:
\begin{eqnarray*}
V
&=&\left(\frac{\pi}{2}\right)^{(N-1)/2}\frac{2^N}{(N+1)!!}\\
&\simeq&\left(\frac{\pi}{2}\right)^{(N-1)/2}2^N\left(\frac{e}{N}\right)^{N/2}\frac{1}{\sqrt{2N}}\\
&=&\sqrt{\frac{2}{\pi}}\left(\frac{2\pi e}{N}\right)^{N/2}\frac{1}{\sqrt{2N}}\\
&=&\left(\frac{2\pi e}{N}\right)^{N/2}\frac{1}{\sqrt{\pi N}}
\end{eqnarray*}

Thus, we are led to the uniform formula in the statement.
\end{proof}

Good to have the above estimates, and in what regards their practical use, more later. By the way, no discussion here would be complete without a word on the gamma function, and we will certainly have an exercise about this, at the end of this chapter.

\bigskip

Getting back now to our main result so far, Theorem 5.24, we can compute in the same way the area of the sphere, the result being as follows:

\index{area of sphere}

\begin{theorem}
The area of the unit sphere in $\mathbb R^N$ is given by
$$A=\left(\frac{\pi}{2}\right)^{[N/2]}\frac{2^N}{(N-1)!!}$$
with the our usual convention for double factorials, namely:
$$N!!=(N-1)(N-3)(N-5)\ldots$$
In particular, at $N=2,3,4$ we obtain respectively $A=2\pi,4\pi,2\pi^2$.
\end{theorem}

\begin{proof}
Regarding the first assertion, we can use here the standard fact, which is elementary, that the area and volume of the sphere in $\mathbb R^N$ are related by the following formula, which together with Theorem 5.24 gives the result:
$$A=N\cdot V$$ 

Alternatively, we can of course redo the computations in the proof of Theorem 5.24, and we obtain the result. As for the last assertion, this can be either worked out directly, or deduced from the results for volumes that we have so far, by multiplying by $N$.
\end{proof}

So long for high dimensional spheres and their volumes. All this is very useful when dealing with Fourier analysis, harmonic functions are related equations, such as the wave and heat ones, and exercise of course for you, to learn more about all this.

\section*{5e. Exercises}

There has been a lot of material in this chapter. In what regards the functions of one variable, and more specifically the second derivative, the standard exercise here is:

\begin{exercise}
Given a convex function $f:\mathbb R\to\mathbb R$, prove that we have the following Jensen inequality, for any $x_1,\ldots,x_N\in\mathbb R$, and any $\lambda_1,\ldots,\lambda_N>0$ summing up to $1$,
$$f(\lambda_1x_1+\ldots+\lambda_Nx_N)\leq\lambda_1f(x_1)+\ldots+\lambda_Nx_N$$
with equality when $x_1=\ldots=x_N$. In particular, by taking the weights $\lambda_i$ to be all equal, we obtain the following Jensen inequality, valid for any $x_1,\ldots,x_N\in\mathbb R$,
$$f\left(\frac{x_1+\ldots+x_N}{N}\right)\leq\frac{f(x_1)+\ldots+f(x_N)}{N}$$
and once again with equality when $x_1=\ldots=x_N$. Prove also that a similar statement holds for the concave functions, with all the inequalities being reversed.
\end{exercise}

This is something very classical, enjoy. For a bonus point, try the functions of several variables as well, and comment on the condition $f''\geq0$ in this case.

\begin{exercise}
Prove that for $p\in(1,\infty)$ we have the following inequality,
$$\left|\frac{x_1+\ldots+x_N}{N}\right|^p\leq\frac{|x_1|^p+\ldots+|x_N|^p}{N}$$
and that for $p\in(0,1)$ we have the following reverse inequality
$$\left|\frac{x_1+\ldots+x_N}{N}\right|^p\geq\frac{|x_1|^p+\ldots+|x_N|^p}{N}$$
with in both cases equality precisely when $|x_1|=\ldots=|x_N|$.
\end{exercise}

As a bonus exercise here, try as well, directly, the case $p=2$.

\begin{exercise}
Develop the theory of the gamma function, defined as
$$\Gamma(s)=\int_0^\infty x^{s-1}e^{-x}dx$$
notably by establishing the following formula, for any $N\in\mathbb N$,
$$\Gamma(N)=(N-1)!$$
and then comment on the formulae for the volumes and areas of spheres.
\end{exercise}

To be more precise, the first question is that of establishing the well-known formula $\Gamma(s+1)=s\Gamma(s)$. The next step is that of computing $\Gamma(s)$ for $s\in\mathbb N/2$, with the above formula in the case $s\in\mathbb N$. And then, the problem is that of deciding if all this can be useful in connection with the formulae for the volumes and areas of spheres.

\chapter{Normal laws}

\section*{6a. Random variables}

In this chapter we discuss the basics of probability theory, as an application of the methods developed in chapter 5. With the idea in mind of doing things a bit abstractly, remember after all that we are algebraists, in this book, as a starting point, we have:

\index{moments}
\index{law}
\index{distribution}

\begin{definition}
Let $X$ be a probability space, that is, a space with a probability measure, and with the corresponding integration denoted $E$, and called expectation.
\begin{enumerate}
\item The random variables are the real functions $f\in L^\infty(X)$.

\item The moments of such a variable are the numbers $M_k(f)=E(f^k)$.

\item The law of such a variable is the measure given by $M_k(f)=\int_\mathbb Rx^kd\mu_f(x)$.
\end{enumerate}
Also, we call mean and variance of $f$ the numbers $E=M_1$ and $V=M_2-M_1^2$.
\end{definition}

All this is self-explanatory, save for the existence of the law $\mu_f$, which is not exactly trivial. But we can do this by looking at formulae of the following type:
$$E(\varphi(f))=\int_\mathbb R\varphi(x)d\mu_f(x)$$

Indeed, having this for monomials $\varphi(x)=x^n$, as above, is the same as having it for polynomials $\varphi\in\mathbb R[X]$, which in turn is the same as having it for the characteristic functions $\varphi=\chi_I$ of measurable sets $I\subset\mathbb R$. Thus, in the end, what we need is:
$$P(f\in I)=\mu_f(I)$$

But this latter formula can serve as a definition for $\mu_f$, and we are done. Next, regarding the key notion of independence, we can formulate here:

\index{independence}

\begin{definition}
Two variables $f,g\in L^\infty(X)$ are called independent when
$$E(f^kg^l)=E(f^k)\,E(g^l)$$
happens, for any $k,l\in\mathbb N$.
\end{definition}

Again, this definition, which was quick, hides some non-trivial things. The idea is a bit as before, namely that of looking at formulae of the following type:
$$E[\varphi(f)\psi(g)]=E[\varphi(f)]\,E[\psi(g)]$$

To be more precise, passing as before from monomials to polynomials, then to characteristic functions, we are led to the usual definition of independence, namely:
$$P(f\in I,g\in J)=P(f\in I)\,P(g\in J)$$

As a first result now, in order to deal with independence, we have:

\index{convolution}

\begin{theorem}
Assuming that $f,g\in L^\infty(X)$ are independent, we have
$$\mu_{f+g}=\mu_f*\mu_g$$
where $*$ is the convolution of real probability measures.
\end{theorem}

\begin{proof}
We have the following computation, using the independence of $f,g$:
$$\int_\mathbb Rx^kd\mu_{f+g}(x)
=E((f+g)^k)
=\sum_r\binom{k}{r}M_r(f)M_{k-r}(g)$$

On the other hand, we have as well the following computation:
\begin{eqnarray*}
\int_\mathbb Rx^kd(\mu_f*\mu_g)(x)
&=&\int_{\mathbb R\times\mathbb R}(x+y)^kd\mu_f(x)d\mu_g(y)\\
&=&\sum_r\binom{k}{r}M_r(f)M_{k-r}(g)
\end{eqnarray*}

Thus $\mu_{f+g}$ and $\mu_f*\mu_g$ have the same moments, so they coincide, as claimed.
\end{proof}

As a second result on independence, which is more advanced, we have:

\index{independence}
\index{Fourier transform}

\begin{theorem}
Assuming that $f,g\in L^\infty(X)$ are independent, we have
$$F_{f+g}=F_fF_g$$
where $F_f(x)=E(e^{ixf})$ is the Fourier transform.
\end{theorem}

\begin{proof}
This is something very standard, based on Theorem 6.3, as follows:
\begin{eqnarray*}
F_{f+g}(x)
&=&\int_\mathbb Re^{ixz}d(\mu_f*\mu_g)(z)\\
&=&\int_{\mathbb R\times\mathbb R}e^{ix(z+t)}d\mu_f(z)d\mu_g(t)\\
&=&\int_\mathbb Re^{ixz}d\mu_f(z)\int_\mathbb Re^{ixt}d\mu_g(t)\\
&=&F_f(x)F_g(x)
\end{eqnarray*}

Thus, we are led to the conclusion in the statement.
\end{proof}

All the above is very nice, we have some interesting theory going on. Let us discuss now some illustrations. We will first talk about discrete probability. First, we have:

\index{Bernoulli law}

\begin{definition}
The Bernoulli law of parameter $x\in[0,1]$ is the law
$$\rho_x=(1-x)\delta_0+x\delta_1$$
appearing when flipping a biased coin, $P({\rm heads})=x$, $P({\rm tails})=1-x$.
\end{definition}

To be more precise, when flipping a biased coin as above, and betting heads, your winning law is $\rho_x$. Next, let us flip the biased coin several times in a row. This leads to:

\index{binomial law}
\index{convolution}
\index{independence}

\begin{theorem}
When flipping a $x$-biased coin $n$ times in a row, the law is
$$\rho_{xn}=\sum_{k=0}^n\binom{n}{k}x^k(1-x)^{n-k}\delta_k$$
called binomial law of parameters $x\in[0,1]$ and $n\in\mathbb N$.
\end{theorem}

\begin{proof}
This is something very standard, the idea being as follows:

\medskip

(1) Observe first that at $n=1$ we have indeed the Bernoulli law $\rho_x$. 

\medskip

(2) In general, we can argue that when flipping the coin $n$ times in a row, and betting heads, the probability of winning $k$ times, among our $n$ attempts, is given by:
$$P(k\ {\rm wins})=\binom{n}{k}P({\rm heads})^kP({\rm tails})^{n-k}=\binom{n}{k}x^k(1-x)^{n-k}$$

Thus, we are led to the formula of $\rho_{xn}$ in the statement.

\medskip

(3) Alternatively, and being a bit more formal, since our $n$ coin tosses are independent, and independence corresponds to convolution, at the level of laws, we have:
\begin{eqnarray*}
\rho_{xn}
&=&\rho_x^{*n}\\
&=&\Big[(1-x)\delta_0+x\delta_1\Big]^{*n}\\
&=&\sum_{k=0}^n\binom{n}{k}x^k(1-x)^{n-k}\,\delta_1^{*k}*\delta_0^{*n-k}\\
&=&\sum_{k=0}^n\binom{n}{k}x^k(1-x)^{n-k}\delta_k
\end{eqnarray*}

(4) Thus, one way or another, we are led to the formula in the statement.
\end{proof}

Getting now to the study of the binomial laws, we have here:

\index{mean}
\index{variance}
\index{moments}

\begin{theorem}
The binomial law $\rho_{xn}$ has the following properties:
\begin{enumerate}
\item The mean is $E=nx$.

\item The variance is $V=nx(1-x)$.
\end{enumerate}
\end{theorem}

\begin{proof}
In what regards the mean, the computation is as follows:
\begin{eqnarray*}
E
&=&\sum_{k=1}^nk\binom{n}{k}x^k(1-x)^{n-k}\\
&=&\sum_{k=1}^n\frac{n!}{(k-1)!(n-k)!}x^k(1-x)^{n-k}\\
&=&nx\sum_{k=1}^n\frac{(n-1)!}{(k-1)!(n-k)!}x^{k-1}(1-x)^{n-k}\\
&=&nx\sum_{t=0}^{n-1}\binom{n-1}{t}x^t(1-x)^{n-t-1}\\
&=&nx(x+1-x)^{n-1}\\
&=&nx
\end{eqnarray*}

With the same trick, we can compute the difference of the first two moments:
\begin{eqnarray*}
M_2-M_1
&=&\sum_{k=2}^n(k^2-k)\binom{n}{k}x^k(1-x)^{n-k}\\
&=&\sum_{k=2}^n\frac{n!}{(k-2)!(n-k)!}x^k(1-x)^{n-k}\\
&=&n(n-1)x^2\sum_{k=2}^n\frac{(n-2)!}{(k-2)!(n-k)!}x^{k-2}(1-x)^{n-k}\\
&=&n(n-1)x^2\sum_{t=0}^{n-2}\binom{n-2}{t}x^t(1-x)^{n-t-2}\\
&=&n(n-1)x^2(x+1-x)^{n-2}\\
&=&n(n-1)x^2
\end{eqnarray*}

We conclude that the second moment is given by the following formula:
$$M_2=n(n-1)x^2+nx=nx((n-1)x+1)$$

As for the variance $V=M_2-M_1^2$, this is given by the following formula:
$$V=nx((n-1)x+1)-(nx)^2=nx(1-x)$$

Thus, we are led to the conclusions in the statement.
\end{proof}

Many other things can be said about the binomial laws, and we will be back to this. Moving on, the central objects in discrete probability theory are the Poisson laws:

\index{Poisson law}

\begin{definition}
The Poisson law of parameter $1$ is the measure
$$p_t=\frac{1}{e}\sum_{k\in\mathbb N}\frac{\delta_k}{k!}$$
and more generally, the Poisson law of parameter $t>0$ is the measure
$$p_t=e^{-t}\sum_{k\in\mathbb N}\frac{t^k}{k!}\,\delta_k$$
with the letter ``p'' standing for Poisson.
\end{definition}

Observe that $p_t$ has indeed mass 1, with this coming from $e^t=\sum_kt^k/k!$. Regarding the mean and variance, these are as follows, and more on this in a moment:
$$E=V=t$$

Many interesting things can be said about the Poisson laws. Going now directly for the kill, Fourier transform computation, we have here the following result:

\index{Fourier transform}
\index{convolution semigroup}

\begin{theorem}
The Fourier transform of $p_t$ is given by:
$$F_{p_t}(y)=\exp\left((e^{iy}-1)t\right)$$
In particular we have $p_s*p_t=p_{s+t}$, called convolution semigroup property.
\end{theorem}

\begin{proof}
We have indeed the following computation, for the Fourier transform:
\begin{eqnarray*}
F_{p_t}(y)
&=&e^{-t}\sum_k\frac{t^k}{k!}F_{\delta_k}(y)\\
&=&e^{-t}\sum_k\frac{t^k}{k!}\,e^{iky}\\
&=&e^{-t}\sum_k\frac{(e^{iy}t)^k}{k!}\\
&=&\exp\left((e^{iy}-1)t\right)
\end{eqnarray*}

As for the second assertion, this follows from the fact that $\log F_{p_t}$ is linear in $t$, via the linearization property for the convolution from Theorem 6.4.
\end{proof}

We can now establish the Poisson Limit Theorem, as follows:

\index{PLT}
\index{Poisson Limit Theorem}
\index{Bernoulli laws}

\begin{theorem}[PLT]
We have the following convergence, in moments,
$$\left(\left(1-\frac{t}{n}\right)\delta_0+\frac{t}{n}\delta_1\right)^{*n}\to p_t$$
for any $t>0$.
\end{theorem}

\begin{proof}
If we denote by $\nu_n$ the measure under the convolution sign, we have the following computation, for the Fourier transform of the limit: 
\begin{eqnarray*}
F_{\delta_r}(y)=e^{iry}
&\implies&F_{\nu_n}(y)=\left(1-\frac{t}{n}\right)+\frac{t}{n}e^{iy}\\
&\implies&F_{\nu_n^{*n}}(y)=\left(\left(1-\frac{t}{n}\right)+\frac{t}{n}e^{iy}\right)^n\\
&\implies&F_{\nu_n^{*n}}(y)=\left(1+\frac{(e^{iy}-1)t}{n}\right)^n\\
&\implies&F(y)=\exp\left((e^{iy}-1)t\right)
\end{eqnarray*}

Thus, we obtain indeed the Fourier transform of $p_t$, as desired.
\end{proof}

At the level of the moments now, the result is quite interesting, as follows:

\index{Bell numbers}
\index{partitions}

\begin{theorem}
The moments of $p_1$ are the Bell numbers,
$$M_k(p_1)=|P(k)|$$
where $P(k)$ is the set of partitions of $\{1,\ldots,k\}$. More generally, we have
$$M_k(p_t)=\sum_{\pi\in P(k)}t^{|\pi|}$$
for any $t>0$, where $|.|$ is the number of blocks. In particular, $E=V=t$.
\end{theorem}

\begin{proof}
We know that the moments of $p_1$ are given by the following formula:
$$M_k=\frac{1}{e}\sum_r\frac{r^k}{r!}$$

We therefore have the following recurrence formula for these moments:
\begin{eqnarray*}
M_{k+1}
&=&\frac{1}{e}\sum_r\frac{r^k}{r!}\left(1+\frac{1}{r}\right)^k\\
&=&\frac{1}{e}\sum_r\frac{r^k}{r!}\sum_s\binom{k}{s}r^{-s}\\
&=&\sum_s\binom{k}{s}M_{k-s}
\end{eqnarray*}

But the Bell numbers $B_k=|P(k)|$ satisfy the same recurrence, trivially, so we have $M_k=B_k$, as claimed. As for the proof of the formula at $t>0$ arbitrary, this is similar. Finally, regarding the mean and variance, $E=t$ is clear, and $V=(t^2+t)-t^2=t$.
\end{proof}

All the above was of course quite quick, but we will be back to this, in chapter 11.

\section*{6b. Central limits}

Getting now to the continuous case, as a key application of the Gauss integral formula, established in chapter 5, we can introduce the normal laws, as follows:

\begin{definition}
The normal law of parameter $1$ is the following measure:
$$g_1=\frac{1}{\sqrt{2\pi}}e^{-x^2/2}dx$$
More generally, the normal law of parameter $t>0$ is the following measure:
$$g_t=\frac{1}{\sqrt{2\pi t}}e^{-x^2/2t}dx$$
These are also called Gaussian distributions, with ``g'' standing for Gauss.
\end{definition}

Observe that the above laws have indeed mass 1, as they should. This follows indeed from the Gauss formula, which gives, with $x=\sqrt{2t}\,y$:
\begin{eqnarray*}
\int_\mathbb R e^{-x^2/2t}dx
&=&\int_\mathbb R e^{-y^2}\sqrt{2t}\,dy\\
&=&\sqrt{2t}\int_\mathbb R e^{-y^2}dy\\
&=&\sqrt{2\pi t}
\end{eqnarray*}

Generally speaking, the normal laws appear as bit everywhere, in real life. The reasons behind this phenomenon come from the Central Limit Theorem (CLT), that we will explain in a moment, after developing some general theory. As a first result, we have:

\begin{proposition}
We have the variance formula
$$V(g_t)=t$$
valid for any $t>0$.
\end{proposition}

\begin{proof}
The first moment is 0, because our normal law $g_t$ is centered. As for the second moment, this can be computed as follows:
\begin{eqnarray*}
M_2
&=&\frac{1}{\sqrt{2\pi t}}\int_\mathbb Rx^2e^{-x^2/2t}dx\\
&=&\frac{1}{\sqrt{2\pi t}}\int_\mathbb R(tx)\left(-e^{-x^2/2t}\right)'dx\\
&=&\frac{1}{\sqrt{2\pi t}}\int_\mathbb Rte^{-x^2/2t}dx\\
&=&t
\end{eqnarray*}

We conclude from this that the variance is $V=M_2=t$.
\end{proof}

Here is another result, which is the key one for the study of the normal laws:

\begin{theorem}
We have the following formula, valid for any $t>0$:
$$F_{g_t}(x)=e^{-tx^2/2}$$
In particular, the normal laws satisfy $g_s*g_t=g_{s+t}$, for any $s,t>0$.
\end{theorem}

\begin{proof}
The Fourier transform formula can be established as follows:
\begin{eqnarray*}
F_{g_t}(x)
&=&\frac{1}{\sqrt{2\pi t}}\int_\mathbb Re^{-y^2/2t+ixy}dy\\
&=&\frac{1}{\sqrt{2\pi t}}\int_\mathbb Re^{-(y/\sqrt{2t}-\sqrt{t/2}ix)^2-tx^2/2}dy\\
&=&\frac{1}{\sqrt{2\pi t}}\int_\mathbb Re^{-z^2-tx^2/2}\sqrt{2t}dz\\
&=&\frac{1}{\sqrt{\pi}}e^{-tx^2/2}\int_\mathbb Re^{-z^2}dz\\
&=&\frac{1}{\sqrt{\pi}}e^{-tx^2/2}\cdot\sqrt{\pi}\\
&=&e^{-tx^2/2}
\end{eqnarray*}

As for the last assertion, this follows from the fact that $\log F_{g_t}$ is linear in $t$, via the linearization property for the convolution from Theorem 6.4.
\end{proof}

We are now ready to state and prove the CLT, as follows:

\begin{theorem}[CLT]
Given random variables $f_1,f_2,f_3,\ldots\in L^\infty(X)$ which are i.i.d., centered, and with variance $t>0$, we have, with $n\to\infty$, in moments,
$$\frac{1}{\sqrt{n}}\sum_{i=1}^nf_i\sim g_t$$
where $g_t$ is the Gaussian law of parameter $t$, having as density $\frac{1}{\sqrt{2\pi t}}e^{-y^2/2t}dy$.
\end{theorem}

\begin{proof}
In terms of moments, the Fourier transform is given by:
\begin{eqnarray*}
F_f(x)
&=&E\left(\sum_{k=0}^\infty\frac{(ixf)^k}{k!}\right)\\
&=&\sum_{k=0}^\infty\frac{(ix)^kE(f^k)}{k!}\\
&=&\sum_{k=0}^\infty\frac{i^kM_k(f)}{k!}\,x^k
\end{eqnarray*}

We conclude that the Fourier transform of the variable in the statement is:
\begin{eqnarray*}
F(x)
&=&\left[F_f\left(\frac{x}{\sqrt{n}}\right)\right]^n\\
&=&\left[1-\frac{tx^2}{2n}+O(n^{-2})\right]^n\\
&\simeq&\left[1-\frac{tx^2}{2n}\right]^n\\
&\simeq&e^{-tx^2/2}
\end{eqnarray*}

But this latter function being the Fourier transform of $g_t$, we obtain the result.
\end{proof}

Let us discuss now some further properties of the normal law. We first have:

\begin{proposition}
The even moments of the normal law are the numbers
$$M_k(g_t)=t^{k/2}\times k!!$$
where $k!!=(k-1)(k-3)(k-5)\ldots\,$, and the odd moments vanish. 
\end{proposition}

\begin{proof}
We have the following computation, valid for any integer $k\in\mathbb N$:
\begin{eqnarray*}
M_k
&=&\frac{1}{\sqrt{2\pi t}}\int_\mathbb Ry^ke^{-y^2/2t}dy\\
&=&\frac{1}{\sqrt{2\pi t}}\int_\mathbb R(ty^{k-1})\left(-e^{-y^2/2t}\right)'dy\\
&=&\frac{1}{\sqrt{2\pi t}}\int_\mathbb Rt(k-1)y^{k-2}e^{-y^2/2t}dy\\
&=&t(k-1)\times\frac{1}{\sqrt{2\pi t}}\int_\mathbb Ry^{k-2}e^{-y^2/2t}dy\\
&=&t(k-1)M_{k-2}
\end{eqnarray*}

Thus by recurrence, we are led to the formula in the statement.
\end{proof}

We have the following alternative formulation of the above result:

\begin{proposition}
The moments of the normal law are the numbers
$$M_k(g_t)=t^{k/2}|P_2(k)|$$
where $P_2(k)$ is the set of pairings of $\{1,\ldots,k\}$.
\end{proposition}

\begin{proof}
Let us count the pairings of $\{1,\ldots,k\}$. In order to have such a pairing, we must pair $1$ with one of the numbers $2,\ldots,k$, and then use a pairing of the remaining $k-2$ numbers. Thus, we have the following recurrence formula:
$$|P_2(k)|=(k-1)|P_2(k-2)|$$

As for the initial data, this is $P_1=0$, $P_2=1$. Thus, we are led to the result.
\end{proof}

We are not done yet, and here is one more improvement of the above:

\begin{theorem}
The moments of the normal law are the numbers
$$M_k(g_t)=\sum_{\pi\in P_2(k)}t^{|\pi|}$$
where $P_2(k)$ is the set of pairings of $\{1,\ldots,k\}$, and $|.|$ is the number of blocks.
\end{theorem}

\begin{proof}
This follows indeed from Proposition 6.17, because the number of blocks of a pairing of $\{1,\ldots,k\}$ is trivially $k/2$, independently of the pairing.
\end{proof}

Observe the similarity with Theorem 6.11, regarding the moments of the Poisson laws. We will see later that many other interesting probability distributions are subject to similar formulae regarding their moments, involving partitions, and a lot of exciting combinatorics. Discussing this will be in fact a main theme of the present book.

\section*{6c. Spherical integrals}

Let us discuss now the computation of the arbitrary integrals over the sphere, and their asymptotics, which will lead us into some key examples of normal variables. We will need a technical result extending the trigonometric formulae from chapter 5, namely:

\index{trigonometric integral}

\begin{theorem}
We have the following formula,
$$\int_0^{\pi/2}\cos^pt\sin^qt\,dt=\left(\frac{\pi}{2}\right)^{\varepsilon(p)\varepsilon(q)}\frac{p!!q!!}{(p+q+1)!!}$$
where $\varepsilon(p)=1$ if $p$ is even, and $\varepsilon(p)=0$ if $p$ is odd, and where
$$m!!=(m-1)(m-3)(m-5)\ldots$$
with the product ending at $2$ if $m$ is odd, and ending at $1$ if $m$ is even.
\end{theorem}

\begin{proof}
Let $I_{pq}$ be the integral in the statement. In order to do the partial integration, a bit as we previously did at $p=0$ or $q=0$, in chapter 5, observe that we have:
\begin{eqnarray*}
(\cos^pt\sin^qt)'
&=&p\cos^{p-1}t(-\sin t)\sin^qt\\
&+&\cos^pt\cdot q\sin^{q-1}t\cos t\\
&=&-p\cos^{p-1}t\sin^{q+1}t+q\cos^{p+1}t\sin^{q-1}t
\end{eqnarray*}

By integrating between $0$ and $\pi/2$, we obtain, for $p,q>0$:
$$pI_{p-1,q+1}=qI_{p+1,q-1}$$

Thus, we can compute $I_{pq}$ by recurrence. When $q$ is even we have:
\begin{eqnarray*}
I_{pq}
&=&\frac{q-1}{p+1}\,I_{p+2,q-2}\\
&=&\frac{q-1}{p+1}\cdot\frac{q-3}{p+3}\,I_{p+4,q-4}\\
&=&\frac{q-1}{p+1}\cdot\frac{q-3}{p+3}\cdot\frac{q-5}{p+5}\,I_{p+6,q-6}\\
&=&\vdots\\
&=&\frac{p!!q!!}{(p+q)!!}\,I_{p+q}
\end{eqnarray*}

But the last term comes from the formulae in chapter 5, and we obtain the result:
\begin{eqnarray*}
I_{pq}
&=&\frac{p!!q!!}{(p+q)!!}\,I_{p+q}\\
&=&\frac{p!!q!!}{(p+q)!!}\left(\frac{\pi}{2}\right)^{\varepsilon(p+q)}\frac{(p+q)!!}{(p+q+1)!!}\\
&=&\left(\frac{\pi}{2}\right)^{\varepsilon(p)\varepsilon(q)}\frac{p!!q!!}{(p+q+1)!!}
\end{eqnarray*}

Observe that this gives the result for $p$ even as well, by symmetry. Indeed, we have $I_{pq}=I_{qp}$, by using the following change of variables:
$$t=\frac{\pi}{2}-s$$

In the remaining case now, where both $p,q$ are odd, we can use once again the formula $pI_{p-1,q+1}=qI_{p+1,q-1}$ established above, and the recurrence goes as follows:
\begin{eqnarray*}
I_{pq}
&=&\frac{q-1}{p+1}\,I_{p+2,q-2}\\
&=&\frac{q-1}{p+1}\cdot\frac{q-3}{p+3}\,I_{p+4,q-4}\\
&=&\frac{q-1}{p+1}\cdot\frac{q-3}{p+3}\cdot\frac{q-5}{p+5}\,I_{p+6,q-6}\\
&=&\vdots\\
&=&\frac{p!!q!!}{(p+q-1)!!}\,I_{p+q-1,1}
\end{eqnarray*}

In order to compute the last term, observe that we have:
\begin{eqnarray*}
I_{p1}
&=&\int_0^{\pi/2}\cos^pt\sin t\,dt\\
&=&-\frac{1}{p+1}\int_0^{\pi/2}(\cos^{p+1}t)'\,dt\\
&=&\frac{1}{p+1}
\end{eqnarray*}

Thus, we can finish our computation in the case $p,q$ odd, as follows:
\begin{eqnarray*}
I_{pq}
&=&\frac{p!!q!!}{(p+q-1)!!}\,I_{p+q-1,1}\\
&=&\frac{p!!q!!}{(p+q-1)!!}\cdot\frac{1}{p+q}\\
&=&\frac{p!!q!!}{(p+q+1)!!}
\end{eqnarray*}

Thus, we obtain the formula in the statement, the exponent of $\pi/2$ appearing there being $\varepsilon(p)\varepsilon(q)=0\cdot 0=0$ in the present case, and this finishes the proof.
\end{proof}

We can now integrate over the spheres, as follows:

\index{spherical integral}
\index{double factorials}

\begin{theorem}
The polynomial integrals over the unit sphere $S^{N-1}_\mathbb R\subset\mathbb R^N$, with respect to the normalized, mass $1$ measure, are given by the following formula,
$$\int_{S^{N-1}_\mathbb R}x_1^{k_1}\ldots x_N^{k_N}\,dx=\frac{(N-1)!!k_1!!\ldots k_N!!}{(N+\Sigma k_i-1)!!}$$
valid when all exponents $k_i$ are even. If an exponent is odd, the integral vanishes.
\end{theorem}

\begin{proof}
Assume first that one of the exponents $k_i$ is odd. We can make then the following change of variables, which shows that the integral in the statement vanishes:
$$x_i\to-x_i$$

Assume now that all the exponents $k_i$ are even. As a first observation, the result holds indeed at $N=2$, due to the formula from Theorem 6.19, which reads:
\begin{eqnarray*}
\int_0^{\pi/2}\cos^pt\sin^qt\,dt
&=&\left(\frac{\pi}{2}\right)^{\varepsilon(p)\varepsilon(q)}\frac{p!!q!!}{(p+q+1)!!}\\\
&=&\frac{p!!q!!}{(p+q+1)!!}
\end{eqnarray*}

Indeed, this formula computes the integral in the statement over the first quadrant. But since the exponents $p,q\in\mathbb N$ are assumed to be even, the integrals over the other quadrants are given by the same formula, so when averaging we obtain the result. 

\medskip

In the general case now, where the dimension $N\in\mathbb N$ is arbitrary, the integral in the statement can be written in spherical coordinates, as follows:
$$I=\frac{2^N}{A}\int_0^{\pi/2}\ldots\int_0^{\pi/2}x_1^{k_1}\ldots x_N^{k_N}J\,dt_1\ldots dt_{N-1}$$

Here $A$ is the area of the sphere, $J$ is the Jacobian, and the $2^N$ factor comes from the restriction to the $1/2^N$ part of the sphere where all the coordinates are positive. According to our formulae in chapter 5, the normalization constant in front of the integral is:
$$\frac{2^N}{A}=\left(\frac{2}{\pi}\right)^{[N/2]}(N-1)!!$$

As for the unnormalized integral, by using the various formulae from chapter 5, for the spherical coordinates and their Jacobian, this is given by:
\begin{eqnarray*}
I'=\int_0^{\pi/2}\ldots\int_0^{\pi/2}
&&(\cos t_1)^{k_1}
(\sin t_1\cos t_2)^{k_2}\\
&&\vdots\\
&&(\sin t_1\sin t_2\ldots\sin t_{N-2}\cos t_{N-1})^{k_{N-1}}\\
&&(\sin t_1\sin t_2\ldots\sin t_{N-2}\sin t_{N-1})^{k_N}\\
&&\sin^{N-2}t_1\sin^{N-3}t_2\ldots\sin^2t_{N-3}\sin t_{N-2}\\
&&dt_1\ldots dt_{N-1}
\end{eqnarray*}

By rearranging the terms, we obtain the following formula:
\begin{eqnarray*}
I'
&=&\int_0^{\pi/2}\cos^{k_1}t_1\sin^{k_2+\ldots+k_N+N-2}t_1\,dt_1\\
&&\int_0^{\pi/2}\cos^{k_2}t_2\sin^{k_3+\ldots+k_N+N-3}t_2\,dt_2\\
&&\vdots\\
&&\int_0^{\pi/2}\cos^{k_{N-2}}t_{N-2}\sin^{k_{N-1}+k_N+1}t_{N-2}\,dt_{N-2}\\
&&\int_0^{\pi/2}\cos^{k_{N-1}}t_{N-1}\sin^{k_N}t_{N-1}\,dt_{N-1}
\end{eqnarray*}

Now by using the above-mentioned formula at $N=2$, this gives:
\begin{eqnarray*}
I'
&=&\frac{k_1!!(k_2+\ldots+k_N+N-2)!!}{(k_1+\ldots+k_N+N-1)!!}\left(\frac{\pi}{2}\right)^{\varepsilon(N-2)}\\
&&\frac{k_2!!(k_3+\ldots+k_N+N-3)!!}{(k_2+\ldots+k_N+N-2)!!}\left(\frac{\pi}{2}\right)^{\varepsilon(N-3)}\\
&&\vdots\\
&&\frac{k_{N-2}!!(k_{N-1}+k_N+1)!!}{(k_{N-2}+k_{N-1}+l_N+2)!!}\left(\frac{\pi}{2}\right)^{\varepsilon(1)}\\
&&\frac{k_{N-1}!!k_N!!}{(k_{N-1}+k_N+1)!!}\left(\frac{\pi}{2}\right)^{\varepsilon(0)}
\end{eqnarray*}

Now let $F$ be the part involving the double factorials, and $P$ be the part involving the powers of $\pi/2$, so that $I'=F\cdot P$. Regarding $F$, by cancelling terms we have:
$$F=\frac{k_1!!\ldots k_N!!}{(\Sigma k_i+N-1)!!}$$

As in what regards $P$, by summing the exponents, we obtain $P=\left(\frac{\pi}{2}\right)^{[N/2]}$. We can now put everything together, and we obtain:
\begin{eqnarray*}
I
&=&\frac{2^N}{A}\times F\times P\\
&=&\left(\frac{2}{\pi}\right)^{[N/2]}(N-1)!!\times\frac{k_1!!\ldots k_N!!}{(\Sigma k_i+N-1)!!}\times\left(\frac{\pi}{2}\right)^{[N/2]}\\
&=&\frac{(N-1)!!k_1!!\ldots k_N!!}{(\Sigma k_i+N-1)!!}
\end{eqnarray*}

Thus, we are led to the conclusion in the statement.
\end{proof}

We have the following useful generalization of the above formula:

\index{spherical integral}

\begin{theorem}
We have the following integration formula over the sphere $S^{N-1}_\mathbb R\subset\mathbb R^N$, with respect to the normalized measure, valid for any exponents $k_i\in\mathbb N$,
$$\int_{S^{N-1}_\mathbb R}|x_1^{k_1}\ldots x_N^{k_N}|\,dx=\left(\frac{2}{\pi}\right)^{\Sigma(k_1,\ldots,k_N)}\frac{(N-1)!!k_1!!\ldots k_N!!}{(N+\Sigma k_i-1)!!}$$
with $\Sigma=[odds/2]$ if $N$ is odd and $\Sigma=[(odds+1)/2]$ if $N$ is even, where ``odds'' denotes the number of odd numbers in the sequence $k_1,\ldots,k_N$.
\end{theorem}

\begin{proof}
As before, the formula holds at $N=2$, due to Theorem 6.19. In general, the integral in the statement can be written in spherical coordinates, as follows:
$$I=\frac{2^N}{A}\int_0^{\pi/2}\ldots\int_0^{\pi/2}x_1^{k_1}\ldots x_N^{k_N}J\,dt_1\ldots dt_{N-1}$$

Here $A$ is the area of the sphere, $J$ is the Jacobian, and the $2^N$ factor comes from the restriction to the $1/2^N$ part of the sphere where all the coordinates are positive. The normalization constant in front of the integral is, as before:
$$\frac{2^N}{A}=\left(\frac{2}{\pi}\right)^{[N/2]}(N-1)!!$$

As for the unnormalized integral, this can be written as before, as follows:
\begin{eqnarray*}
I'
&=&\int_0^{\pi/2}\cos^{k_1}t_1\sin^{k_2+\ldots+k_N+N-2}t_1\,dt_1\\
&&\int_0^{\pi/2}\cos^{k_2}t_2\sin^{k_3+\ldots+k_N+N-3}t_2\,dt_2\\
&&\vdots\\
&&\int_0^{\pi/2}\cos^{k_{N-2}}t_{N-2}\sin^{k_{N-1}+k_N+1}t_{N-2}\,dt_{N-2}\\
&&\int_0^{\pi/2}\cos^{k_{N-1}}t_{N-1}\sin^{k_N}t_{N-1}\,dt_{N-1}
\end{eqnarray*}

Now by using the formula at $N=2$, we get:
\begin{eqnarray*}
I'
&=&\frac{\pi}{2}\cdot\frac{k_1!!(k_2+\ldots+k_N+N-2)!!}{(k_1+\ldots+k_N+N-1)!!}\left(\frac{2}{\pi}\right)^{\delta(k_1,k_2+\ldots+k_N+N-2)}\\
&&\frac{\pi}{2}\cdot\frac{k_2!!(k_3+\ldots+k_N+N-3)!!}{(k_2+\ldots+k_N+N-2)!!}\left(\frac{2}{\pi}\right)^{\delta(k_2,k_3+\ldots+k_N+N-3)}\\
&&\vdots\\
&&\frac{\pi}{2}\cdot\frac{k_{N-2}!!(k_{N-1}+k_N+1)!!}{(k_{N-2}+k_{N-1}+k_N+2)!!}\left(\frac{2}{\pi}\right)^{\delta(k_{N-2},k_{N-1}+k_N+1)}\\
&&\frac{\pi}{2}\cdot\frac{k_{N-1}!!k_N!!}{(k_{N-1}+k_N+1)!!}\left(\frac{2}{\pi}\right)^{\delta(k_{N-1},k_N)}
\end{eqnarray*}

In order to compute this quantity, let us denote by $F$ the part involving the double factorials, and by $P$ the part involving the powers of $\pi/2$, so that we have:
$$I'=F\cdot P$$

Regarding $F$, there are many cancellations there, and we end up with:
$$F=\frac{k_1!!\ldots k_N!!}{(\Sigma k_i+N-1)!!}$$

As in what regards $P$, the $\delta$ exponents on the right sum up to the following number:
$$\Delta(k_1,\ldots,k_N)=\sum_{i=1}^{N-1}\delta(k_i,k_{i+1}+\ldots+k_N+N-i-1)$$

In other words, with this notation, the above formula reads:
\begin{eqnarray*}
I'
&=&\left(\frac{\pi}{2}\right)^{N-1}\frac{k_1!!k_2!!\ldots k_N!!}{(k_1+\ldots+k_N+N-1)!!}\left(\frac{2}{\pi}\right)^{\Delta(k_1,\ldots,k_N)}\\
&=&\left(\frac{2}{\pi}\right)^{\Delta(k_1,\ldots,k_N)-N+1}\frac{k_1!!k_2!!\ldots k_N!!}{(k_1+\ldots+k_N+N-1)!!}\\
&=&\left(\frac{2}{\pi}\right)^{\Sigma(k_1,\ldots,k_N)-[N/2]}\frac{k_1!!k_2!!\ldots k_N!!}{(k_1+\ldots+k_N+N-1)!!}
\end{eqnarray*}

To be more precise, the formula relating $\Delta$ to $\Sigma$ follows from a number of simple observations, the first of which being the fact that, due to obvious parity reasons, the sequence of $\delta$ numbers appearing in the definition of $\Delta$ cannot contain two consecutive zeroes. Now together with $I=(2^N/V)I'$, this gives the formula in the statement.
\end{proof}

Summarizing, we have complete results for the integration over the spheres, with the answers involving various multinomial type coefficients, defined in terms of factorials, or of double factorials. All these formulae are of course very useful, in practice.

\bigskip

As a basic application of all this, we have the following result:

\index{hyperspherical law}
\index{normal law}

\begin{theorem}
The moments of the hyperspherical variables are
$$\int_{S^{N-1}_\mathbb R}x_i^kdx=\frac{(N-1)!!k!!}{(N+k-1)!!}$$
and the normalized hyperspherical variables 
$$y_i=\frac{x_i}{\sqrt{N}}$$
become normal and independent with $N\to\infty$.
\end{theorem}

\begin{proof}
We have two things to be proved, the idea being as follows:

\medskip

(1) The formula in the statement follows from the general integration formula over the sphere, from Theorem 6.20. Indeed, that formula gives:
$$\int_{S^{N-1}_\mathbb R}x_i^kdx=\frac{(N-1)!!k!!}{(N+k-1)!!}$$

Now observe that with $N\to\infty$ we have the following estimate:
\begin{eqnarray*}
\int_{S^{N-1}_\mathbb R}x_i^kdx
&=&\frac{(N-1)!!}{(N+k-1)!!}\times k!!\\
&\simeq&N^{k/2}k!!\\
&=&N^{k/2}M_k(g_1)
\end{eqnarray*}

Thus, the variables $y_i=\frac{x_i}{\sqrt{N}}$ become normal with $N\to\infty$.

\medskip

(2) As for the asymptotic independence result, this is standard as well, once again by using Theorem 6.20, for computing mixed moments, and taking the $N\to\infty$ limit.
\end{proof}

As a comment here, all this might seem quite specialized. However, we will see later on that all this is related to linear algebra, and more specifically to the fine study of the group $O_N$ formed by the orthogonal matrices. But more on this later.

\section*{6d. Complex spheres}

Let us discuss now the complex analogues of all the above. We must first introduce the complex analogues of the normal laws, and this can be done as follows:

\index{complex normal law}
\index{complex Gaussian law}

\begin{definition}
The complex Gaussian law of parameter $t>0$ is
$$G_t=law\left(\frac{1}{\sqrt{2}}(a+ib)\right)$$
where $a,b$ are independent, each following the law $g_t$.
\end{definition}

The combinatorics of these laws is a bit more complicated than in the real case, and we will be back to this in a moment. But to start with, we have:

\index{convolution}

\begin{theorem}
The complex Gaussian laws have the property
$$G_s*G_t=G_{s+t}$$
for any $s,t>0$, and so they form a convolution semigroup.
\end{theorem}

\begin{proof}
This follows indeed from the real result, for the usual Gaussian laws, established in above, by taking real and imaginary parts.
\end{proof}

We have as well the following complex analogue of the CLT:

\index{complex CLT}

\begin{theorem}[CCLT]
Given complex random variables $f_1,f_2,f_3,\ldots\in L^\infty(X)$, which are i.i.d., centered, and with variance $t>0$, we have, with $n\to\infty$, in moments,
$$\frac{1}{\sqrt{n}}\sum_{i=1}^nf_i\sim G_t$$
where $G_t$ is the complex Gaussian law of parameter $t$.
\end{theorem}

\begin{proof}
This follows indeed from the real CLT, established above, simply by taking the real and imaginary parts of all the variables involved.
\end{proof}

Regarding now the moments, things are a bit more complicated than before, because our variables are now complex instead of real. In order to deal with this issue, we will use ``colored moments'', which are the expectations of the ``colored powers'', with these latter powers being defined by the following formulae, and multiplicativity:
$$f^\emptyset=1\quad,\quad f^\circ=f\quad,\quad f^\bullet=\bar{f}$$

With these conventions made, the result is as follows, with a pairing of a colored integer $k=\circ\bullet\bullet\circ\ldots$ being called matching when it pairs $\circ$ symbols with $\bullet$ symbols:

\index{moments}
\index{matching pairings}

\begin{theorem}
The moments of the complex normal law are the numbers
$$M_k(G_t)=\sum_{\pi\in\mathcal P_2(k)}t^{|\pi|}$$
where $\mathcal P_2(k)$ are the matching pairings of $\{1,\ldots,k\}$, and $|.|$ is the number of blocks.
\end{theorem}

\begin{proof}
This can be done in several steps, as follows:

\medskip

(1) We recall from the above that the moments of the real Gaussian law $g_1$, with respect to integer exponents $k\in\mathbb N$, are the following numbers:
$$m_k=|P_2(k)|$$

(2) We will show here that in what concerns the complex Gaussian law $G_1$, a similar result holds. Numerically, we will prove that we have the following formula, where a colored integer $k=\circ\bullet\bullet\circ\ldots$ is called uniform when it contains the same number of $\circ$ and $\bullet$\,, and where $|k|\in\mathbb N$ is the length of such a colored integer:
$$M_k=\begin{cases}
(|k|/2)!&(k\ {\rm uniform})\\
0&(k\ {\rm not\ uniform})
\end{cases}$$

Now since the matching partitions $\pi\in\mathcal P_2(k)$ are counted by exactly the same numbers, and this for trivial reasons, we will obtain the formula in the statement, namely:
$$M_k=|\mathcal P_2(k)|$$

(3) This was for the plan. In practice now, we must compute the moments, with respect to colored integer exponents $k=\circ\bullet\bullet\circ\ldots$\,, of the variable in Definition 6.23:
$$c=\frac{1}{\sqrt{2}}(a+ib)$$

As a first observation, in the case where such an exponent $k=\circ\bullet\bullet\circ\ldots$ is not uniform in $\circ,\bullet$\,, a rotation argument shows that the corresponding moment of $c$ vanishes. To be more precise, the variable $c'=wc$ can be shown to be complex Gaussian too, for any $w\in\mathbb C$, and from $M_k(c)=M_k(c')$ we obtain $M_k(c)=0$, in this case.

\medskip

(4) In the uniform case now, where $k=\circ\bullet\bullet\circ\ldots$ consists of $p$ copies of $\circ$ and $p$ copies of $\bullet$\,, the corresponding moment can be computed as follows:
\begin{eqnarray*}
M_k
&=&\int(c\bar{c})^p\\
&=&\frac{1}{2^p}\int(a^2+b^2)^p\\
&=&\frac{1}{2^p}\sum_s\binom{p}{s}\int a^{2s}\int b^{2p-2s}\\
&=&\frac{1}{2^p}\sum_s\binom{p}{s}(2s)!!(2p-2s)!!\\
&=&\frac{1}{2^p}\sum_s\frac{p!}{s!(p-s)!}\cdot\frac{(2s)!}{2^ss!}\cdot\frac{(2p-2s)!}{2^{p-s}(p-s)!}\\
&=&\frac{p!}{4^p}\sum_s\binom{2s}{s}\binom{2p-2s}{p-s}
\end{eqnarray*}

(5) In order to finish now the computation, let us recall that we have the following formula, coming from the generalized binomial formula, or from the Taylor formula:
$$\frac{1}{\sqrt{1+t}}=\sum_{k=0}^\infty\binom{2k}{k}\left(\frac{-t}{4}\right)^k$$

By taking the square of this series, we obtain the following formula:
\begin{eqnarray*}
\frac{1}{1+t}
&=&\sum_{ks}\binom{2k}{k}\binom{2s}{s}\left(\frac{-t}{4}\right)^{k+s}\\
&=&\sum_p\left(\frac{-t}{4}\right)^p\sum_s\binom{2s}{s}\binom{2p-2s}{p-s}
\end{eqnarray*}

Now by looking at the coefficient of $t^p$ on both sides, we conclude that the sum on the right equals $4^p$. Thus, we can finish the moment computation in (4), as follows:
$$M_p=\frac{p!}{4^p}\times 4^p=p!$$

(6) As a conclusion, if we denote by $|k|$ the length of a colored integer $k=\circ\bullet\bullet\circ\ldots$\,, the moments of the variable $c$ in the statement are given by:
$$M_k=\begin{cases}
(|k|/2)!&(k\ {\rm uniform})\\
0&(k\ {\rm not\ uniform})
\end{cases}$$

On the other hand, the numbers $|\mathcal P_2(k)|$ in the statement are given by exactly the same formula. Indeed, in order to have matching pairings of $k$, our exponent $k=\circ\bullet\bullet\circ\ldots$ must be uniform, consisting of $p$ copies of $\circ$ and $p$ copies of $\bullet$, with:
$$p=\frac{|k|}{2}$$

But then the matching pairings of $k$ correspond to the permutations of the $\bullet$ symbols, as to be matched with $\circ$ symbols, and so we have $p!$ such matching pairings. Thus, we have exactly the same formula as for the moments of $c$, and this finishes the proof.
\end{proof}

There are of course many other possible proofs for the above result, which are all instructive, and some further theory as well, that can be developed for the complex normal variables, which is very interesting too. We refer here to Feller \cite{fel}, or Durrett \cite{dur}. We will be back to this, on several occasions, in what follows.

\bigskip 

In practice, we also need to know how to compute joint moments of independent normal variables. We have here the following result, to be used later on:

\index{complex normal law}
\index{colored moments}
\index{matching pairings}
\index{Wick formula}

\begin{theorem}[Wick formula]
Given independent variables $f_i$, each following the complex normal law $G_t$, with $t>0$ being a fixed parameter, we have the formula
$$E\left(f_{i_1}^{k_1}\ldots f_{i_s}^{k_s}\right)=t^{s/2}\#\left\{\pi\in\mathcal P_2(k)\Big|\pi\leq\ker i\right\}$$
where $k=k_1\ldots k_s$ and $i=i_1\ldots i_s$, for the joint moments of these variables.
\end{theorem}

\begin{proof}
This is something well-known, and the basis for all possible computations with complex normal variables, which can be proved in two steps, as follows:

\medskip

(1) Let us first discuss the case where we have a single variable $f$, which amounts in taking $f_i=f$ for any $i$ in the formula in the statement. What we have to compute here are the moments of $f$, with respect to colored integer exponents $k=\circ\bullet\bullet\circ\ldots\,$, and the formula in the statement tells us that these moments must be:
$$E(f^k)=t^{|k|/2}|\mathcal P_2(k)|$$

But this is the formula in Theorem 6.26, so we are done with this case.

\medskip

(2) In general now, when expanding the product $f_{i_1}^{k_1}\ldots f_{i_s}^{k_s}$ and rearranging the terms, we are left with doing a number of computations as in (1), and then making the product of the expectations that we found. But this amounts in counting the partitions in the statement, with the condition $\pi\leq\ker i$ there standing for the fact that we are doing the various type (1) computations independently, and then making the product.
\end{proof}

The above statement is one of the possible formulations of the Wick formula, and there are in fact many more formulations, which are all useful. Here is an alternative such formulation, which is quite popular, and that we will also use in what follows:

\begin{theorem}[Wick formula 2]
Given independent variables $f_i$, each following the complex normal law $G_t$, with $t>0$ being a fixed parameter, we have the formula
$$E\left(f_{i_1}\ldots f_{i_k}f_{j_1}^*\ldots f_{j_k}^*\right)=t^k\#\left\{\pi\in S_k\Big|i_{\pi(r)}=j_r,\forall r\right\}$$
for the non-vanishing joint moments of these variables.
\end{theorem}

\begin{proof}
This follows from the usual Wick formula, from Theorem 6.27. With some changes in the indices and notations, the formula there reads:
$$E\left(f_{I_1}^{K_1}\ldots f_{I_s}^{K_s}\right)=t^{s/2}\#\left\{\sigma\in\mathcal P_2(K)\Big|\sigma\leq\ker I\right\}$$

Now observe that we have $\mathcal P_2(K)=\emptyset$, unless the colored integer $K=K_1\ldots K_s$ is uniform, in the sense that it contains the same number of $\circ$ and $\bullet$ symbols. Up to permutations, the non-trivial case, where the moment is non-vanishing, is the case where the colored integer $K=K_1\ldots K_s$ is of the following special form:
$$K=\underbrace{\circ\circ\ldots\circ}_k\ \underbrace{\bullet\bullet\ldots\bullet}_k$$

So, let us focus on this case, which is the non-trivial one. Here we have $s=2k$, and we can write the multi-index $I=I_1\ldots I_s$ in the following way:
$$I=i_1\ldots i_k\ j_1\ldots j_k$$

With these changes made, the above usual Wick formula reads:
$$E\left(f_{i_1}\ldots f_{i_k}f_{j_1}^*\ldots f_{j_k}^*\right)=t^k\#\left\{\sigma\in\mathcal P_2(K)\Big|\sigma\leq\ker(ij)\right\}$$

The point now is that the matching pairings $\sigma\in\mathcal P_2(K)$, with $K=\circ\ldots\circ\bullet\ldots\bullet\,$, of length $2k$, as above, correspond to the permutations $\pi\in S_k$, in the obvious way. With this identification made, the above modified usual Wick formula becomes:
$$E\left(f_{i_1}\ldots f_{i_k}f_{j_1}^*\ldots f_{j_k}^*\right)=t^k\#\left\{\pi\in S_k\Big|i_{\pi(r)}=j_r,\forall r\right\}$$

Thus, we have reached to the formula in the statement, and we are done.
\end{proof}

Finally, here is one more formulation of the Wick formula, which is useful as well:

\begin{theorem}[Wick formula 3]
Given independent variables $f_i$, each following the complex normal law $G_t$, with $t>0$ being a fixed parameter, we have the formula
$$E\left(f_{i_1}f_{j_1}^*\ldots f_{i_k}f_{j_k}^*\right)=t^k\#\left\{\pi\in S_k\Big|i_{\pi(r)}=j_r,\forall r\right\}$$
for the non-vanishing joint moments of these variables.
\end{theorem}

\begin{proof}
This follows from our second Wick formula, from Theorem 6.28, simply by permuting the terms, as to have an alternating sequence of plain and conjugate variables. Alternatively, we can start with Theorem 6.27, and then perform the same manipulations as in the proof of Theorem 6.28, but with the exponent being this time as follows: 
$$K=\underbrace{\circ\bullet\circ\bullet\ldots\ldots\circ\bullet}_{2k}$$

Thus, we are led to the conclusion in the statement.
\end{proof}

In relation now with the spheres, we first have the following variation of the integration formula in Theorem 6.20, dealing this time with integrals over the complex sphere:

\index{spherical integral}

\begin{theorem}
We have the following integration formula over the complex sphere $S^{N-1}_\mathbb C\subset\mathbb R^N$, with respect to the normalized measure, 
$$\int_{S^{N-1}_\mathbb C}|z_1|^{2l_1}\ldots|z_N|^{2l_N}\,dz=4^{\sum l_i}\frac{(2N-1)!l_1!\ldots l_n!}{(2N+\sum l_i-1)!}$$
valid for any exponents $l_i\in\mathbb N$. As for the other polynomial integrals in $z_1,\ldots,z_N$ and their conjugates $\bar{z}_1,\ldots,\bar{z}_N$, these all vanish.
\end{theorem}

\begin{proof}
Consider an arbitrary polynomial integral over $S^{N-1}_\mathbb C$, written as follows:
$$I=\int_{S^{N-1}_\mathbb C}z_{i_1}\bar{z}_{i_2}\ldots z_{i_{2l-1}}\bar{z}_{i_{2l}}\,dz$$

(1) By using transformations of type $p\to\lambda p$ with $|\lambda|=1$, we see that $I$ vanishes, unless each $z_a$ appears as many times as $\bar{z}_a$ does, and this gives the last assertion.

\medskip

(2) Assume now that we are in the non-vanishing case. Then the $l_a$ copies of $z_a$ and the $l_a$ copies of $\bar{z}_a$ produce by multiplication a factor $|z_a|^{2l_a}$, so we have:
$$I=\int_{S^{N-1}_\mathbb C}|z_1|^{2l_1}\ldots|z_N|^{2l_N}\,dz$$

Now by using the standard identification $S^{N-1}_\mathbb C\simeq S^{2N-1}_\mathbb R$, we obtain:
\begin{eqnarray*}
I
&=&\int_{S^{2N-1}_\mathbb R}(x_1^2+y_1^2)^{l_1}\ldots(x_N^2+y_N^2)^{l_N}\,d(x,y)\\
&=&\sum_{r_1\ldots r_N}\binom{l_1}{r_1}\ldots\binom{l_N}{r_N}\int_{S^{2N-1}_\mathbb R}x_1^{2l_1-2r_1}y_1^{2r_1}\ldots x_N^{2l_N-2r_N}y_N^{2r_N}\,d(x,y)
\end{eqnarray*}

(3) By using the formula in Theorem 6.20, we obtain:
\begin{eqnarray*}
&&I\\
&=&\sum_{r_1\ldots r_N}\binom{l_1}{r_1}\ldots\binom{l_N}{r_N}\frac{(2N-1)!!(2r_1)!!\ldots(2r_N)!!(2l_1-2r_1)!!\ldots (2l_N-2r_N)!!}{(2N+2\sum l_i-1)!!}\\
&=&\sum_{r_1\ldots r_N}\binom{l_1}{r_1}\ldots\binom{l_N}{r_N}
\frac{(2N-1)!(2r_1)!\ldots (2r_N)!(2l_1-2r_1)!\ldots (2l_N-2r_N)!}{(2N+\sum l_i-1)!r_1!\ldots r_N!(l_1-r_1)!\ldots (l_N-r_N)!}
\end{eqnarray*}

(4) We can rewrite the sum on the right in the following way:
\begin{eqnarray*}
&&I\\
&=&\sum_{r_1\ldots r_N}\frac{l_1!\ldots l_N!(2N-1)!(2r_1)!\ldots (2r_N)!(2l_1-2r_1)!\ldots (2l_N-2r_N)!}{(2N+\sum l_i-1)!(r_1!\ldots r_N!(l_1-r_1)!\ldots (l_N-r_N)!)^2}\\
&=&\sum_{r_1}\binom{2r_1}{r_1}\binom{2l_1-2r_1}{l_1-r_1}\ldots\sum_{r_N}\binom{2r_N}{r_N}\binom{2l_N-2r_N}{l_N-r_N}\frac{(2N-1)!l_1!\ldots l_N!}{(2N+\sum l_i-1)!}\\
&=&4^{l_1}\times\ldots\times 4^{l_N}\times\frac{(2N-1)!l_1!\ldots l_N!}{(2N+\sum l_i-1)!}
\end{eqnarray*}

Thus, we obtain the formula in the statement.
\end{proof}

Regarding now the hyperspherical variables, investigated in the above in the real case, we have similar results for the complex spheres, as follows:

\index{hyperspherical law}

\begin{theorem}
The rescaled coordinates on the complex sphere $S^{N-1}_\mathbb C$, 
$$w_i=\frac{z_i}{\sqrt{N}}$$
become complex Gaussian and independent with $N\to\infty$.
\end{theorem}

\begin{proof}
We have two assertions to be proved, the idea being as follows:

\medskip

(1) The assertion about the laws follows exactly as in the real case, by using this time Theorem 6.30 as a main technical ingredient.

\medskip

(2) As for the independence result, this follows as well as in the real case, by using this time the Wick formula as a main technical ingredient.
\end{proof}

As a conclusion to all this, we have now a good level in linear algebra, and also in probability. And this can only open up a whole new set of perspectives, on what further books can be read, in relation with geometry, analysis, and physics. 

\bigskip

As for algebra and probability, stay with us. The story is far from being over with what we learned, and dozens of further interesting things to follow. We still have 250 more pages, and there will be algebra and probability in them, that is promised.

\section*{6e. Exercises}

We have learned many interesting things in this chapter, and there are many possible exercises about this. First, in connection with the CLT, we have:

\begin{exercise}
Work out the precise convergence conclusions in the CLT, 
$$\frac{1}{\sqrt{n}}\sum_{i=1}^nf_i\sim g_t$$
going beyond the convergence in moments, which was established in the above.
\end{exercise}

This is a bit vague, but at this stage, learning more theory would be a good thing. Of course, in case all this looks a bit complicated, don't hesitate to look it up. As already mentioned, some good references for probability are Durrett \cite{dur} and Feller \cite{fel}.

\begin{exercise}
Find an alternative proof for the moment formula 
$$M_k(G_t)=\sum_{\pi\in\mathcal P_2(k)}t^{|\pi|}$$
using a method of your choice.
\end{exercise}

Again, this is a bit vague, and many things that you can try. As before, in case you lack a new idea here, don't hesitate to look it up, and report on what you learned.

\begin{exercise}
Find a probability measure $\nu$ whose moments are given by
$$M_k(\nu)=|NC_2(k)|$$
then find as well a probability measure $\eta$ whose moments are given by
$$M_k(\eta)=|NC(k)|$$
where $NC$ stands for ``noncrossing''. Then try as well the parametric case.
\end{exercise}

These latter exercises are actually quite difficult, but still doable, with some patience, and you will learn many interesting things in this way, notably in relation with the moment problem, which is a key topic in advanced probability. By the way, for a bonus point, try to solve as well the question left, regarding the noncrossing matching pairings. With this latter question being also quite difficult, but definitely worth studying.

\begin{exercise}
Compute the density of the hyperspherical law at $N=4$, that is, the law of one of the coordinates over the unit sphere $S^3_\mathbb R\subset\mathbb R^4$.
\end{exercise}

This might look a bit specialized, but trust me, it is a must-do exercise, and if you find something quite interesting, as an answer here, do not be surprised. After all, $S^3_\mathbb R$ is the sphere of space-time, having its own magic. We will be back to this.

\chapter{Special matrices}

\section*{7a. Fourier matrices}

In this chapter we go back to basic linear algebra questions. We will be interested in various classes of ``special matrices'', and in the tools for dealing with them. As a first and central example here, which is obviously special, we have the flat matrix:

\index{flat matrix}
\index{all-one matrix}
\index{all-one vector}

\begin{definition}
The flat matrix $\mathbb I_N$ is the all-one $N\times N$ matrix:
$$\mathbb I_N=\begin{pmatrix}
1&\ldots&1\\
\vdots&&\vdots\\
1&\ldots&1\end{pmatrix}$$
Equivalently, $\mathbb I_N/N$ is the orthogonal projection on the all-one vector $\xi\in\mathbb C^N$.
\end{definition}

Observe that $\mathbb I_N$ has a lot of interesting properties, such as being circulant, and bistochastic. The idea will be that many techniques that can be applied to $\mathbb I_N$, with quite trivial results, apply to such special classes of matrices, with non-trivial consequences.

\bigskip

A first interesting question regarding $\mathbb I_N$ concerns its diagonalization. Since $\mathbb I_N$ is a multiple of a rank 1 projection, we have right away the following result:

\begin{proposition}
The flat matrix diagonalizes as follows,
$$\mathbb I_N=P
\begin{pmatrix}
N\\
&0\\
&&\ddots\\
&&&0\end{pmatrix}P^{-1}$$
where $P\in M_N(\mathbb C)$ can be any matrix formed by the all one-vector $\xi$, followed by $N-1$ linearly independent solutions $x\in\mathbb C^N$ of the equation $x_1+\ldots+x_N=0$.
\end{proposition}

\begin{proof}
This follows indeed from our linear algebra knowledge from chapters 1-4, by using the fact that $\mathbb I_N/N$ is the orthogonal projection onto $\mathbb C\xi$.
\end{proof}

In practice now, the problem which is left is that of finding an explicit matrix $P\in M_N(\mathbb C)$, as above. To be more precise, there are plently of solutions here, some of them being even real, $P\in M_N(\mathbb R)$, and the problem is that of finding a ``nice'' such solution, say having the property that $P_{ij}$ appears as an explicit function of $i,j$.

\bigskip

Long story short, we are led to the question of solving, in a somewhat canonical and elegant way, the following equation, over the real or the complex numbers:
$$x_1+\ldots+x_N=0$$

And this question is more tricky than it seems. To be more precise, there is no hope of doing this over the real numbers. As in what regards the complex numbers, there is a ray of light here coming from the roots of unity. So, let us formulate:

\index{Fourier matrix}

\begin{definition}
The Fourier matrix $F_N$ is the following matrix, with $w=e^{2\pi i/N}$:
$$F_N=
\begin{pmatrix}
1&1&1&\ldots&1\\
1&w&w^2&\ldots&w^{N-1}\\
1&w^2&w^4&\ldots&w^{2(N-1)}\\
\vdots&\vdots&\vdots&&\vdots\\
1&w^{N-1}&w^{2(N-1)}&\ldots&w^{(N-1)^2}
\end{pmatrix}$$
That is, $F_N=(w^{ij})_{ij}$, with indices $i,j\in\{0,1,\ldots,N-1\}$, taken modulo $N$.
\end{definition}

Before getting further, observe that this matrix $F_N$ is ``special'' too, but in a different sense, its main properties being the fact that it is a Vandermonde matrix, and also, a rescaled unitary. We will axiomatize later the matrices of this type.

\bigskip

Getting back now to the diagonalization problem for the flat matrix $\mathbb I_N$, this can be solved by using the Fourier matrix $F_N$, in the following elegant way:

\index{flat matrix}
\index{Fourier matrix}

\begin{theorem}
The flat matrix diagonalizes as follows,
$$\mathbb I_N=\frac{1}{N}\,F_N
\begin{pmatrix}
N\\
&0\\
&&\ddots\\
&&&0\end{pmatrix}F_N^*$$
with $F_N=(w^{ij})_{ij}$ being the Fourier matrix.
\end{theorem}

\begin{proof}
According to Proposition 7.2, and with indices $i,j\in\{0,1,\ldots,N-1\}$, we are left with finding the 0-eigenvectors of $\mathbb I_N$, which amounts in solving:
$$x_0+\ldots+x_{N-1}=0$$

But for this purpose, we use the root of unity $w=e^{2\pi i/N}$, and more specifically, the following standard formula, that we know from chapter 3:
$$\sum_{i=0}^{N-1}w^{ij}=N\delta_{j0}$$

Indeed, this formula shows that for $j=1,\ldots,N-1$, the vector $v_j=(w^{ij})_i$ is a 0-eigenvector. Moreover, these vectors are pairwise orthogonal, because we have:
$$<v_j,v_k>
=\sum_iw^{ij-ik}
=N\delta_{jk}$$

Thus, we have our basis $\{v_1,\ldots,v_{N-1}\}$ of 0-eigenvectors, and since the $N$-eigenvector is $\xi=v_0$, the passage matrix $P$ that we are looking is given by:
$$P=\begin{bmatrix}v_0&v_1&\ldots&v_{N-1}\end{bmatrix}$$

But this is precisely the Fourier matrix, $P=F_N$. In order to finish now, observe that the above computation of $<v_i,v_j>$ shows that $F_N/\sqrt{N}$ is unitary, and so:
$$F_N^{-1}=\frac{1}{N}\,F_N^*$$

Thus, we are led to the diagonalization formula in the statement.
\end{proof}

Generally speaking, the above result will be the template for what we will be doing here. On one hand we will have special matrices to be studied, of $\mathbb I_N$ type, and on the other hand we will have special matrices that can be used as tools, of $F_N$ type. Let us begin with a discussion of the ``tools''. Inspired by $F_N$, let us formulate:

\index{Hadamard matrix}

\begin{definition}
A complex Hadamard matrix is a square matrix
$$H\in M_N(\mathbb T)$$
where $\mathbb T$ is the unit circle, satisfying the following equivalent conditions:
\begin{enumerate}
\item The rows are pairwise orthogonal.

\item The columns are pairwise orthogonal.

\item The rescaled matrix $H/\sqrt{N}$ is unitary.

\item The rescaled matrix $H^t/\sqrt{N}$ is unitary.
\end{enumerate}
\end{definition}

Here the fact that the above conditions are indeed equivalent comes from basic linear algebra, and more specifically from the fact that a matrix $U\in M_N(\mathbb C)$ is a unitary precisely when the rows, or columns, have norm 1, and are pairwise orthogonal.

\bigskip

We already know, from the proof of Theorem 7.4, that the Fourier matrix $F_N$ is a complex Hadamard matrix. There are many other examples of complex Hadamard matrices, and the basic theory of such matrices can be summarized as follows:

\begin{proposition}
The class of $N\times N$ complex Hadamard matrices is as follows:
\begin{enumerate}
\item It contains the Fourier matrix $F_N$.

\item It is stable under taking tensor products.

\item It is stable under taking transposes, conjugates and adjoints.

\item It is stable under permuting rows, or permuting columns.

\item It is stable under multiplying rows or columns by numbers in $\mathbb T$.
\end{enumerate}
\end{proposition}

\begin{proof}
All this is elementary, the idea being as follows:

\medskip

(1) This is something that we already know, from the proof of Theorem 7.4.

\medskip

(2) Assume that $H\in M_M(\mathbb T)$ and $K\in M_N(\mathbb T)$ are Hadamard matrices, and consider their tensor product, which in double index notation is as follows:
$$(H\otimes K)_{ia,jb}=H_{ij}K_{ab}$$

We have then $H\otimes K\in M_{MN}(\mathbb T)$, and the rows $R_{ia}$ of this matrix are pairwise orthogonal, as shown by the following computation:
\begin{eqnarray*}
<R_{ia},R_{kc}>
&=&\sum_{jb}H_{ij}K_{ab}\cdot\bar{H}_{kj}\bar{K}_{cb}\\
&=&\sum_jH_{ij}\bar{H}_{kj}\sum_bK_{ab}\bar{K}_{cb}\\
&=&MN\delta_{ik}\delta_{ac}
\end{eqnarray*}

(3) We know that the set formed by the $N\times N$ complex Hadamard matrices appears as follows, with the intersection being taken inside $M_N(\mathbb C)$:
$$X_N=M_N(\mathbb T)\cap\sqrt{N}U_N$$

The set $M_N(\mathbb T)$ is stable under the operations in the statement. As for the set $\sqrt{N}U_N$, here we can use the well-known fact that if a matrix is unitary, $U\in U_N$, then so is its complex conjugate $\bar{U}=(\bar{U}_{ij})$, the inversion formulae being as follows:
$$U^*=U^{-1}\quad,\quad 
U^t=\bar{U}^{-1}$$

Thus the unitary group $U_N$ is stable under the following operations:
$$U\to U^t\quad,\quad 
U\to\bar{U}\quad,\quad 
U\to U^*$$

It follows that the above set $X_N$ is stable as well under these operations, as desired.

\medskip

(4-5) These assertions are clear from definitions, because permuting rows or columns, or multiplying them by numbers in $\mathbb T$, leaves invariant both $M_N(\mathbb T)$ and $\sqrt{N}U_N$.
\end{proof}

In the above result, the assertions (1,2) are really important, and (3,4,5) are rather technical remarks. As a consequence, coming from (1,2), let us formulate:

\index{generalized Fourier matrix}
\index{Walsh matrix}
\index{Hadamard matrix}

\begin{theorem}
The following matrices, called generalized Fourier matrices,
$$F_{N_1,\ldots,N_k}=F_{N_1}\otimes\ldots\otimes F_{N_k}$$
are Hadamard, for any choice of $N_1,\ldots,N_k$. In particular the following matrices,
$$W_N=\begin{pmatrix}1&1\\1&-1\end{pmatrix}^{\otimes k}$$
having size $N=2^k$, and called Walsh matrices, are all Hadamard.
\end{theorem}

\begin{proof}
The first assertion comes from Proposition 7.6. As for the second assertion, this comes from this, by taking $N_1=\ldots=N_k=2$. Indeed, the matrix that we get is:
$$F_{2,\ldots,2}
=F_2^{\otimes k}
=\begin{pmatrix}1&1\\1&-1\end{pmatrix}^{\otimes k}$$

Thus, we are led to the conclusion in the statement.
\end{proof}

As an illustration for the above result, the second Walsh matrix, which is an Hadamard matrix having real entries, as is the case with all the Walsh matrices, is as follows:
$$W_4=\begin{pmatrix}1&1&1&1\\ 1&-1&1&-1\\ 1&1&-1&-1\\ 1&-1&-1&1\end{pmatrix}$$

In order to work out now some classification results, let us formulate:

\index{equivalent Hadamard matrices}

\begin{definition}
Two complex Hadamard matrices are called equivalent, and we write $H\sim K$, when it is possible to pass from $H$ to $K$ via the following operations:
\begin{enumerate}
\item Permuting the rows, or permuting the columns.

\item Multiplying the rows or columns by numbers in $\mathbb T$.
\end{enumerate}
\end{definition}

To be more precise, this is based on Proposition 7.6. Also, we have not taken into account all the results there, because the operations $H\to H^t,\bar{H},H^*$ are far more subtle than those in (1,2) above, and can complicate things, if included in the equivalence. Now with this notion of equivalence in hand, we first have the following result:

\begin{theorem}
The Hadamard matrices at $N=2,3,4$ are up to equivalence
$$F_2=\begin{pmatrix}1&1\\ 1&-1\end{pmatrix}\quad,\quad
F_3=\begin{pmatrix}1&1&1\\ 1&w&w^2\\ 1&w^2&w\end{pmatrix}\quad,\quad 
F_4^q=\begin{pmatrix}
1&1&1&1\\
1&-1&1&-1\\
1&q&-1&-q\\ 
1&-q&-1&q
\end{pmatrix}$$
with $w=e^{2\pi i/3}$, and with  $q\in\mathbb T$.
\end{theorem}

\begin{proof}
This is something elementary, the idea being as follows:

\medskip

(1) At $N=2$ the result is clear, because up to equivalence we can put our matrix in the following form, and then the Hadamard condition gives $x=-1$:
$$H=\begin{pmatrix}1&1\\ 1&x\end{pmatrix}$$

(2) At $N=3$ now, again up to equivalence, we can assume that our matrix is:
$$H=\begin{pmatrix}1&1&1\\ 1&x&y\\ 1&z&t\end{pmatrix}$$

The orthogonality conditions between the rows of this matrix read:
$$x+y=-1\quad,\quad 
z+t=-1\quad,\quad 
x\bar{z}+y\bar{t}=-1$$

In order to process these conditions, consider an equation of the following type:
$$p+q=-1\quad,\quad p,q\in\mathbb T$$

Now observe that this equation tells us that the triangle having vertices at $1,p,q$ must be equilateral, and so, that we must have $\{p,q\}=\{w,w^2\}$, with $w=e^{2\pi i/3}$. By using this fact, for the first two equations, we conclude that we must have:
$$\{x,y\}=\{w,w^2\}\quad,\quad 
\{z,t\}=\{w,w^2\}$$

As for the third equation, this gives $x\neq z$. Thus, $H$ is either the Fourier matrix $F_3$, or the matrix obtained from $F_3$ by permuting the last two columns, and we are done.

\medskip

(3) As for the proof at $N=4$, where what we get are certain deformations of $F_4$, covering for instance $W_4$, this is similar, and we will leave this as an exercise.
\end{proof}

At $N=5$ things get more complicated, and following Haagerup \cite{haa}, we have:

\begin{theorem}
The only Hadamard matrix at $N=5$ is the Fourier matrix,
$$F_5=\begin{pmatrix}
1&1&1&1&1\\
1&w&w^2&w^3&w^4\\
1&w^2&w^4&w&w^3\\
1&w^3&w&w^4&w^2\\
1&w^4&w^3&w^2&w
\end{pmatrix}$$
with $w=e^{2\pi i/5}$, up to the standard equivalence relation for such matrices. 
\end{theorem}

\begin{proof}
This is something quite technical, the idea being as follows:

\medskip

(1) Consider an Hadamard matrix $H\in M_5(\mathbb T)$, chosen dephased, as follows:
$$H=\begin{pmatrix}
1&1&1&1&1\\
1&a&x&*&*\\
1&y&b&*&*\\
1&*&*&*&*\\
1&*&*&*&*
\end{pmatrix}$$

By using the orthogonality of rows and columns, and doing some computations, we eventually conclude that the numbers $a,b,x,y$ must satisfy the following equations:
$$(a-b)(a-xy)(b-xy)=0$$
$$(x-y)(x-ab)(y-ab)=0$$

(2) Our claim now is that, by doing some combinatorics, we can actually obtain from this $a=b$ and $x=y$, up to the equivalence relation for the Hadamard matrices: 
$$H\sim\begin{pmatrix}
1&1&1&1&1\\
1&a&x&*&*\\
1&x&a&*&*\\
1&*&*&*&*\\
1&*&*&*&*
\end{pmatrix}$$

Indeed, the above two equations lead to 9 possible cases, the first of which is, as desired, $a=b$ and $x=y$. As for the remaining 8 cases, here again things are determined by 2 parameters, and in practice, we can always permute the first 3 rows and 3 columns, and then dephase our matrix, as for our matrix to take the above special form.

\medskip

(3) But with this in hand, the combinatorics of the scalar products between the first 3 rows, and between the first 3 columns as well, becomes something which is quite simple to investigate. By doing a routine study here, and then completing it with a study of the lower right $2\times2$ corner as well, we are led to 2 possible cases, as follows:
$$H\sim\begin{pmatrix}
1&1&1&1&1\\
1&a&b&c&d\\
1&b&a&d&c\\
1&c&d&a&b\\
1&d&c&b&a
\end{pmatrix}\quad,\quad 
H\sim\begin{pmatrix}
1&1&1&1&1\\
1&a&b&c&d\\
1&b&a&d&c\\
1&c&d&b&a\\
1&d&c&a&b
\end{pmatrix}$$

(4) Next, a routine study shows that the first case is in fact not possible. Regarding now the second case, the orthogonality equations there are as follows:
\begin{eqnarray*}
a+b+c+d&=&-1\\
2Re(a\bar{b})+2Re(c\bar{d})&=&-1\\
a\bar{c}+c\bar{b}+b\bar{d}+d\bar{a}&=&-1
\end{eqnarray*}

Now observe that the third equation can be written in the following form:
\begin{eqnarray*}
Re[(a+b)(\bar{c}+\bar{d})]&=&-1\\
Im[(a-b)(\bar{c}-\bar{d})]&=&0
\end{eqnarray*}

By using now $a,b,c,d\in\mathbb T$, we conclude that we can find $s,t\in\mathbb R$ such that:
$$a+b=is(a-b)\quad,\quad 
c+d=it(c-d)$$

By plugging in these values, our system of equations simplifies, as follows:
\begin{eqnarray*}
(a+b)+(c+d)&=&-1\\
|a+b|^2+|c+d|^2&=&3\\
(a+b)(\bar{c}+\bar{d})&=&-1
\end{eqnarray*}

(5) Now observe that the last equation implies in particular that we have:
$$|a+b|^2\cdot|c+d|^2=1$$

Thus $|a+b|^2,|c+d|^2$ must be roots of $X^2-3X+1=0$, and this gives:
$$\Big\{|a+b|\,,\,|c+d|\Big\}=\left\{\frac{\sqrt{5}+1}{2}\,,\,\frac{\sqrt{5}-1}{2}\right\}$$

Which is very good news, because, obviously, we are now into 5-th roots of unity. 

\medskip

(6) Next, we have 2 cases to be considered. The first one is as follows, with $z\in\mathbb T$:
$$a+b=\frac{\sqrt{5}+1}{2}\,z\quad,\quad 
c+d=-\frac{\sqrt{5}-1}{2}\,z$$

But from $a+b+c+d=-1$ we obtain $z=-1$, and by using this we conclude that we have $b=\bar{a}$, $d=\bar{c}$. Thus we have the following formulae:
$$Re(a)=\cos(2\pi/5)\quad,\quad 
Re(c)=\cos(\pi/5)$$

We conclude that we have an equivalence $H\sim F_5$, as claimed. As for the second case, with the variables $a,b$ and $c,d$ interchanged, this leads to $H\sim F_5$ as well.
\end{proof}

At $N=6$ now, things explode, and we have here all sorts of matrices, related or not to $F_6$, and not classified yet. As an example here, we have the following matrix of Bj\"orck and Fr\"oberg, with $a\in\mathbb T$ being one of the roots of $a^2+(\sqrt{3}-1)a+1=0$:
$$BF_6=\begin{pmatrix}
1&ia&-a&-i&-\bar{a}&i\bar{a}\\
i\bar{a}&1&ia&-a&-i&-\bar{a}\\
-\bar{a}&i\bar{a}&1&ia&-a&-i\\
-i&-\bar{a}&i\bar{a}&1&ia&-a\\
-a&-i&-\bar{a}&i\bar{a}&1&ia\\
ia&-a&-i&-\bar{a}&i\bar{a}&1
\end{pmatrix}$$

Finally, let us mention that the generalized Fourier matrices, and the Hadamard matrices in general, have many applications, to questions in coding, radio transmissions, quantum physics, and many more. We refer here for instance to the book of Bengtsson-\.Zyczkowski \cite{bzy}, and to the papers of Bj\"orck \cite{bjo}, Haagerup \cite{haa}, Idel-Wolf \cite{iwo}, Jones \cite{jo4}, Sylvester \cite{syl}. We will be back to these matrices later, on several occasions.

\section*{7b. Circulant matrices}

Let us go back now to the general linear algebra considerations from the beginning of this chapter. We have seen that $F_N$ diagonalizes in an elegant way the flat matrix $\mathbb I_N$, and the idea in what follows will be that of $F_N$, or other real or complex Hadamard matrices, can be used in order to deal with other matrices, of $\mathbb I_N$ type.

\bigskip

A first feature of the flat matrix $\mathbb I_N$ is that it is circulant, in the following sense:

\index{circulant matrix}

\begin{definition}
A real or complex matrix $M$ is called circulant if
$$M_{ij}=\xi_{j-i}$$
for a certain vector $\xi$, with the indices taken modulo $N$.
\end{definition}

The circulant matrices are beautiful mathematical objects, which appear of course in many serious problems as well. As an example, at $N=4$, we must have:
$$M=\begin{pmatrix}
a&b&c&d\\
d&a&b&c\\
c&d&a&b\\
b&c&d&a
\end{pmatrix}$$

The point now is that, while certainly gently looking, these matrices can be quite diabolic, when it comes to diagonalization, and other problems. For instance, when $M$ is real, the computations with $M$ are usually very complicated over the real numbers. Fortunately the complex numbers and the Fourier matrices are there, and we have:

\index{Fourier-diagonal matrix}
\index{discrete Fourier transform}

\begin{theorem}
For a matrix $M\in M_N(\mathbb C)$, the following are equivalent:
\begin{enumerate}
\item $M$ is circulant, $M_{ij}=\xi_{j-i}$, for a certain vector $\xi\in\mathbb C^N$.

\item $M$ is Fourier-diagonal, $M=F_NQF_N^*$, for a certain diagonal matrix $Q$.
\end{enumerate}
If so, $\xi=F_N^*q$, where $q\in\mathbb C^N$ is the column vector formed by the diagonal entries of $Q$.
\end{theorem}

\begin{proof}
This follows indeed from some basic computations with roots of unity:

\medskip

$(1)\implies(2)$ Assuming $M_{ij}=\xi_{j-i}$, the matrix $Q=F_N^*MF_N$ is diagonal, due to:
\begin{eqnarray*}
Q_{ij}
&=&\sum_{kl}w^{-ik}M_{kl}w^{lj}\\
&=&\sum_{kl}w^{jl-ik}\xi_{l-k}\\
&=&\sum_{kr}w^{j(k+r)-ik}\xi_r\\
&=&\sum_rw^{jr}\xi_r\sum_kw^{(j-i)k}\\
&=&N\delta_{ij}\sum_rw^{jr}\xi_r
\end{eqnarray*}

$(2)\implies(1)$ Assuming now $Q=diag(q_1,\ldots,q_N)$, the matrix $M=F_NQF_N^*$ is circulant, as shown by the following computation:
$$M_{ij}
=\sum_kw^{ik}Q_{kk}w^{-jk}
=\sum_kw^{(i-j)k}q_k$$

To be more precise, in this formula the last term depends only on $j-i$, and so shows that we have $M_{ij}=\xi_{j-i}$, with $\xi$ being the following vector:
$$\xi_i
=\sum_kw^{-ik}q_k
=(F_N^*q)_i$$

Thus, we are led to the conclusions in the statement.
\end{proof}

As a basic illustration for the above result, for the circulant matrix $M=\mathbb I_N$ we recover in this way the diagonalization result from Theorem 7.4, namely:
$$\mathbb I_N=\frac{1}{N}\,F_N
\begin{pmatrix}
N\\
&0\\
&&\ddots\\
&&&0\end{pmatrix}F_N^*$$

The above result is something quite powerful, and very useful, and suggests doing everything in Fourier, when dealing with circulant matrices. And we can use here:

\index{discrete Fourier transform}

\begin{theorem}
The various basic sets of $N\times N$ circulant matrices are as follows, with the convention that associated to any $q\in\mathbb C^N$ is the matrix  $Q=diag(q_1,\ldots,q_N)$:
\begin{enumerate}
\item The set of all circulant matrices is:
$$M_N(\mathbb C)^{circ}=\left\{F_NQF_N^*\Big|q\in\mathbb C^N\right\}$$

\item The set of all circulant unitary matrices is:
$$U_N^{circ}=\left\{\frac{1}{N}F_NQF_N^*\Big|q\in\mathbb T^N\right\}$$

\item The set of all circulant orthogonal matrices is:
$$O_N^{circ}=\left\{\frac{1}{N}F_NQF_N^*\Big|q\in\mathbb T^N,\bar{q}_i=q_{-i},\forall i\right\}$$
\end{enumerate}
In addition, in this picture, the first row vector of $F_NQF_N^*$ is given by $\xi=F_N^*q$.
\end{theorem}

\begin{proof}
All this follows from Theorem 7.12, as follows:

\medskip

(1) This assertion, along with the last one, is Theorem 7.12 itself.

\medskip

(2) This is clear from (1), and from the fact that the rescaled matrix $F_N/\sqrt{N}$ is unitary, because the eigenvalues of a unitary matrix must be on the unit circle $\mathbb T$.

\medskip

(3) This follows from (2), because the matrix is real when $\xi_i=\bar{\xi}_i$, and in Fourier transform, $\xi=F_N^*q$, this corresponds to the condition $\bar{q}_i=q_{-i}$.
\end{proof}

As a last topic regarding the circulant matrices, which is somehow one level above the considerations above, let us discuss the circulant Hadamard matrices. We first have:

\begin{proposition}
The following are circulant and symmetric Hadamard matrices,
$$F_2'=\begin{pmatrix}i&1\\1&i\end{pmatrix}\quad,\quad
F_3'=\begin{pmatrix}w&1&1\\1&w&1\\1&1&w\end{pmatrix}\quad,\quad 
F_4''=\begin{pmatrix}-1&\nu&1&\nu\\\nu&-1&\nu&1\\1&\nu&-1&\nu\\ \nu&1&\nu&-1\end{pmatrix}$$
where $w=e^{2\pi i/3},\nu=e^{\pi i/4}$, equivalent to the Fourier matrices $F_2,F_3,F_4$.
\end{proposition}

\begin{proof}
The orthogonality between rows being clear, we have here complex Hadamard matrices. The fact that we have an equivalence $F_2\sim F_2'$ follows from:
$$\begin{pmatrix}1&1\\1&-1\end{pmatrix}
\sim\begin{pmatrix}i&i\\1&-1\end{pmatrix}
\sim\begin{pmatrix}i&1\\1&i\end{pmatrix}$$

At $N=3$ now, the equivalence $F_3\sim F_3'$ can be constructed as follows:
$$\begin{pmatrix}1&1&1\\1&w&w^2\\1&w^2&w\end{pmatrix}
\sim\begin{pmatrix}1&1&w\\1&w&1\\w&1&1\end{pmatrix}
\sim\begin{pmatrix}w&1&1\\1&w&1\\1&1&w\end{pmatrix}$$

As for the case $N=4$, here the equivalence $F_4\sim F_4''$ can be constructed as follows, where we use the logarithmic notation $[k]_s=e^{2\pi ki/s}$, with respect to $s=8$:
$$\begin{bmatrix}0&0&0&0\\0&2&4&6\\0&4&0&4\\0&6&4&2\end{bmatrix}_8
\sim\begin{bmatrix}0&1&4&1\\1&4&1&0\\4&1&0&1\\1&0&1&4\end{bmatrix}_8
\sim\begin{bmatrix}4&1&0&1\\1&4&1&0\\0&1&4&1\\1&0&1&4\end{bmatrix}_8
$$

Thus, the Fourier matrices $F_2,F_3,F_4$ can be put indeed in circulant form.
\end{proof}

In order to discuss now the general case, we will use a technical method for dealing with the circulant matrices, namely Bj\"orck's cyclic root formalism \cite{bjo}, as follows:

\index{circulant matrix}
\index{cyclic root}

\begin{theorem}
Assume that a matrix $H\in M_N(\mathbb T)$ is circulant, $H_{ij}=\gamma_{j-i}$. Then $H$ is is a complex Hadamard matrix if and only if the vector 
$$z=(z_0,z_1,\ldots,z_{N-1})$$
given by $z_i=\gamma_i/\gamma_{i-1}$ satisfies the following equations:
\begin{eqnarray*}
z_0+z_1+\ldots+z_{N-1}&=&0\\
z_0z_1+z_1z_2+\ldots+z_{N-1}z_0&=&0\\
\ldots\\
z_0z_1\ldots z_{N-2}+\ldots+z_{N-1}z_0\ldots z_{N-3}&=&0\\
z_0z_1\ldots z_{N-1}&=&1
\end{eqnarray*}
If so is the case, we say that $z=(z_0,\ldots,z_{N-1})$ is a cyclic $N$-root.
\end{theorem}

\begin{proof}
This follows indeed from a direct computation, the idea being that, with $H_{ij}=\gamma_{j-i}$ as above, the orthogonality conditions between the rows are best written in terms of the variables $z_i=\gamma_i/\gamma_{i-1}$, and correspond to the equations in the statement.
\end{proof}

Now back to the Fourier matrices, we have the following result:

\begin{theorem}
Given $N\in\mathbb N$, construct the following complex numbers:
$$\nu=e^{\pi i/N}\quad,\quad 
q=\nu^{N-1}\quad,\quad 
w=\nu^2$$
We have then a cyclic $N$-root, given by the following formula,
$$(q,qw,qw^2,\ldots,qw^{N-1})$$
and the corresponding complex Hadamard matrix $F_N'$ is circulant and symmetric, and equivalent to the Fourier matrix $F_N$.
\end{theorem}

\begin{proof}
Given two numbers $q,w\in\mathbb T$, let us find out when $(q,qw,qw^2,\ldots,qw^{N-1})$ is a cyclic root. We have two conditions to be verified, as follows:

\medskip

(1) In order for the $=0$ equations in Theorem 7.15 to be satisfied, the value of $q$ is irrelevant, and $w$ must be a primitive $N$-root of unity. 

\medskip

(2) As for the $=1$ equation in Theorem 7.15, this states that we must have: 
$$q^Nw^{\frac{N(N-1)}{2}}=1$$

Thus, we must have $q^N=(-1)^{N-1}$, so with the values of $q,w\in\mathbb T$ in the statement, we have a cyclic $N$-root. Now construct $H_{ij}=\gamma_{j-i}$ as in Theorem 7.15. We have:
\begin{eqnarray*}
\gamma_k=\gamma_{-k}
&\iff&q^{k+1}w^{\frac{k(k+1)}{2}}=q^{-k+1}w^{\frac{k(k-1)}{2}}\\
&\iff&q^{2k}w^k=1\\
&\iff&q^2=w^{-1}
\end{eqnarray*}

But this latter condition holds indeed, because we have:
$$q^2
=\nu^{2N-2}
=\nu^{-2}
=w^{-1}$$

We conclude that our circulant matrix $H$ is symmetric as well, as claimed. It remains to construct an equivalence $H\sim F_N$. In order to do this, observe that, due to our conventions $q=\nu^{N-1},w=\nu^2$, the first row vector of $H$ is given by:
\begin{eqnarray*}
\gamma_k
&=&q^{k+1}w^{\frac{k(k+1)}{2}}\\
&=&\nu^{(N-1)(k+1)}\nu^{k(k+1)}\\
&=&\nu^{(N+k-1)(k+1)}
\end{eqnarray*}

Thus, the entries of $H$ are given by the following formula:
\begin{eqnarray*}
H_{-i,j}
&=&H_{0,i+j}\\
&=&\nu^{(N+i+j-1)(i+j+1)}\\
&=&\nu^{i^2+j^2+2ij+Ni+Nj+N-1}\\
&=&\nu^{N-1}\cdot\nu^{i^2+Ni}\cdot\nu^{j^2+Nj}\cdot\nu^{2ij}
\end{eqnarray*}

We conclude that the matrix $H=(H_{ij})$ is equivalent to the following matrix:
$$H'=(H_{-i,j})$$

Now regarding this latter matrix $H'$, observe that in the above formula, the factors $\nu^{N-1}$, $\nu^{i^2+Ni}$, $\nu^{j^2+Nj}$ correspond respectively to a global multiplication by a scalar, and to row and column multiplications by scalars. Thus $H'$ is equivalent to the matrix $H''$ obtained from it by deleting these factors. But this latter matrix, given by $H''_{ij}=\nu^{2ij}$ with $\nu=e^{\pi i/N}$, is precisely the Fourier matrix $F_N$, and we are done.
\end{proof}

As an illustration, at $N=2,3$ we obtain the old matrices $F_2',F_3'$. As for the case $N=4$, here we obtain the following matrix, with $\nu=e^{\pi i/4}$:
$$F_4'=\begin{pmatrix}
\nu^3&1&\nu^7&1\\
1&\nu^3&1&\nu^7\\
\nu^7&1&\nu^3&1\\
1&\nu^7&1&\nu^3
\end{pmatrix}$$

This matrix is equivalent to the matrix $F_4''$ from Proposition 7.14, with the equivalence $F_4'\sim F_4''$ being obtained by multiplying everything by the number $\nu=e^{\pi i/4}$. 

\bigskip

There are many other things that can be said about the circulant Hadamard matrices, and about the Fourier matrices, and we refer here to Bj\"orck \cite{bjo} and Haagerup \cite{haa}.

\section*{7c. Bistochastic matrices}

Getting back now to the main idea behind what we are doing, namely building on the relation between $\mathbb I_N$ and $F_N$, let us study now the class of bistochastic matrices:

\index{row-stochastic matrix}
\index{column-stochastic matrix}
\index{bistochastic matrix}

\begin{definition}
A square matrix $M\in M_N(\mathbb C)$ is called bistochastic if each row and each column sum up to the same number:
$$\begin{matrix}
M_{11}&\ldots&M_{1N}&\to&\lambda\\
\vdots&&\vdots\\
M_{N1}&\ldots&M_{NN}&\to&\lambda\\
\downarrow&&\downarrow\\
\lambda&&\lambda
\end{matrix}$$
If this happens only for the rows, or only for the columns, the matrix is called row-stochastic, respectively column-stochastic.
\end{definition}

As a basic example of a bistochastic matrix, we have of course the flat matrix $\mathbb I_N$. In fact, the various above notions of stochasticity are closely related to $\mathbb I_N$, or rather to the all-one vector $\xi$ that the matrix $\mathbb I_N/N$ projects on, in the following way:

\index{all-one vector}

\begin{proposition}
Let $M\in M_N(\mathbb C)$ be a square matrix.
\begin{enumerate}
\item $M$ is row stochastic, with sums $\lambda$, when $M\xi=\lambda\xi$.

\item $M$ is column stochastic, with sums $\lambda$, when $M^t\xi=\lambda\xi$.

\item $M$ is bistochastic, with sums $\lambda$, when $M\xi=M^t\xi=\lambda\xi$.
\end{enumerate}
\end{proposition}

\begin{proof}
All these assertions are clear from definitions, because when multiplying a matrix by $\xi$, we obtain the vector formed by the row sums.
\end{proof}

As an observation here, we can reformulate if we want the above statement in a purely matrix-theoretic form, by using the flat matrix $\mathbb I_N$, as follows:

\begin{proposition}
Let $M\in M_N(\mathbb C)$ be a square matrix.
\begin{enumerate}
\item $M$ is row stochastic, with sums $\lambda$, when $M\mathbb I_N=\lambda\mathbb I_N$.

\item $M$ is column stochastic, with sums $\lambda$, when $\mathbb I_NM=\lambda\mathbb I_N$.

\item $M$ is bistochastic, with sums $\lambda$, when $M\mathbb I_N=\mathbb I_NM=\lambda\mathbb I_N$.
\end{enumerate}
\end{proposition}

\begin{proof}
This follows from Proposition 7.18, and from the fact that both the rows and the columns of the flat matrix $\mathbb I_N$ are copies of the all-one vector $\xi$.
\end{proof}

In what follows we will be mainly interested in the unitary bistochastic matrices, which are quite interesting objects. As a first result, regarding such matrices, we have:

\begin{theorem}
For a unitary matrix $U\in U_N$, the following conditions are equivalent:
\begin{enumerate}
\item $H$ is bistochastic, with sums $\lambda$.

\item $H$ is row stochastic, with sums $\lambda$, and $|\lambda|=1$.

\item $H$ is column stochastic, with sums $\lambda$, and $|\lambda|=1$.
\end{enumerate}
\end{theorem}

\begin{proof}
By using a symmetry argument we just need to prove $(1)\iff(2)$, and both the implications are elementary, as follows:

\medskip

$(1)\implies(2)$ If we denote by $U_1,\ldots,U_N\in\mathbb C^N$ the rows of $U$, we have indeed:
\begin{eqnarray*}
1
&=&\sum_i<U_1,U_i>\\
&=&\sum_jU_{1j}\sum_i\bar{U}_{ij}\\
&=&\sum_jU_{1j}\cdot\bar{\lambda}\\
&=&|\lambda|^2
\end{eqnarray*}

$(2)\implies(1)$ Consider the all-one vector $\xi=(1)_i\in\mathbb C^N$. The fact that $U$ is row-stochastic with sums $\lambda$ reads:
\begin{eqnarray*}
\sum_jU_{ij}=\lambda,\forall i
&\iff&\sum_jU_{ij}\xi_j=\lambda\xi_i,\forall i\\
&\iff&U\xi=\lambda\xi
\end{eqnarray*}

Also, the fact that $U$ is column-stochastic with sums $\lambda$ reads:
\begin{eqnarray*}
\sum_iU_{ij}=\lambda,\forall j
&\iff&\sum_jU_{ij}\xi_i=\lambda\xi_j,\forall j\\
&\iff&U^t\xi=\lambda\xi
\end{eqnarray*}

We must prove that the first condition implies the second one, provided that the row sum $\lambda$ satisfies $|\lambda|=1$. But this follows from the following computation:
\begin{eqnarray*}
U\xi=\lambda\xi
&\implies&U^*U\xi=\lambda U^*\xi\\
&\implies&\xi=\lambda U^*\xi\\
&\implies&\xi=\bar{\lambda}U^t\xi\\
&\implies&U^t\xi=\lambda\xi
\end{eqnarray*}

Thus, we have proved both the implications, and we are done.
\end{proof}

The unitary bistochastic matrices are stable under a number of operations, and in particular under taking products, and we have the following result:

\index{bistochastic group}
\index{discrete Fourier transform}

\begin{theorem}
The real and complex bistochastic groups, which are the sets
$$B_N\subset O_N\quad,\quad 
C_N\subset U_N$$
consisting of matrices which are bistochastic, are isomorphic to $O_{N-1}$, $U_{N-1}$.
\end{theorem}

\begin{proof}
Let us pick a unitary matrix $F\in U_N$ satisfying the following condition, where $e_0,\ldots,e_{N-1}$ is the standard basis of $\mathbb C^N$, and where $\xi$ is the all-one vector:
$$Fe_0=\frac{1}{\sqrt{N}}\xi$$ 

Observe that such matrices $F\in U_N$ exist indeed, the basic example being the normalized Fourier matrix $F_N/\sqrt{N}$. We have then, by using the above property of $F$:
\begin{eqnarray*}
u\xi=\xi
&\iff&uFe_0=Fe_0\\
&\iff&F^*uFe_0=e_0\\
&\iff&F^*uF=diag(1,w)
\end{eqnarray*}

Thus we have isomorphisms as in the statement, given by $w_{ij}\to(F^*uF)_{ij}$.
\end{proof}

We will be back to $B_N,C_N$ later in this book, when doing group theory. In relation now with the Hadamard matrices, as a first remark, the first Walsh matrix $W_2$ looks better in complex bistochastic form, modulo the standard equivalence relation:
$$\begin{pmatrix}1&1\\1&-1\end{pmatrix}
\sim\begin{pmatrix}i&i\\1&-1\end{pmatrix}
\sim\begin{pmatrix}i&1\\1&i\end{pmatrix}$$

The second Walsh matrix $W_4=W_2\otimes W_2$ can be put as well in complex bistochastic form, as follows, and also looks better in bistochastic form:
$$\ \begin{pmatrix}1&1&1&1\\ 1&-1&1&-1\\ 1&1&-1&-1\\ 1&-1&-1&1\end{pmatrix}
\sim\begin{pmatrix}
-1&1&1&1\\
1&-1&1&1\\
1&1&-1&1\\
1&1&1&-1
\end{pmatrix}$$

In fact, by using the above formulae, we are led to the following statement:

\index{Walsh matrix}
\index{bistochastic matrix}

\begin{proposition}
All the Walsh matrices, $W_N=W_2^{\otimes n}$ with $N=2^n$, can be put in bistochastic form, up to the standard equivalence relation, as follows:
\begin{enumerate}
\item The matrices $W_N$ with $N=4^n$ admit a real bistochastic form, namely:
$$W_N\sim\begin{pmatrix}
-1&1&1&1\\
1&-1&1&1\\
1&1&-1&1\\
1&1&1&-1
\end{pmatrix}^{\otimes n}$$

\item The matrices $W_N$ with $N=2\times4^n$ admit a complex bistochastic form, namely:
$$W_N\sim\begin{pmatrix}i&1\\1&i\end{pmatrix}\otimes\begin{pmatrix}
-1&1&1&1\\
1&-1&1&1\\
1&1&-1&1\\
1&1&1&-1
\end{pmatrix}^{\otimes n}$$
\end{enumerate}
\end{proposition}

\begin{proof}
This follows indeed from the above discussion.
\end{proof}

Regarding now the question of putting the general Hadamard matrices, real or complex, in complex bistochastic form, things here are tricky. We first have:

\index{bistochastic Hadamard matrix}

\begin{theorem}
The class of the bistochastic complex Hadamard matrices has the following properties:
\begin{enumerate}
\item It contains the circulant symmetric forms $F_N'$ of the Fourier matrices $F_N$.

\item It is stable under permuting rows and columns.

\item It is stable under taking tensor products.
\end{enumerate}
In particular, any generalized Fourier matrix $F_{N_1,\ldots,N_k}=F_{N_1}\otimes\ldots\otimes F_{N_k}$ can be put in bistochastic and symmetric form, up to the equivalence relation.
\end{theorem}

\begin{proof}
We have several things to be proved, the idea being as follows:

\medskip

(1) We know from the above that any Fourier matrix $F_N$ has a circulant and symmetric form $F_N'$. But since circulant implies bistochastic, this gives the result.

\medskip

(2) The claim regarding permuting rows and columns is clear. 

\medskip

(3) Assuming that $H,K$ are bistochastic, with sums $\lambda,\mu$, we have:
\begin{eqnarray*}
\sum_{ia}(H\otimes K)_{ia,jb}
&=&\sum_{ia}H_{ij}K_{ab}\\
&=&\sum_iH_{ij}\sum_aK_{ab}\\
&=&\lambda\mu
\end{eqnarray*}

We have as well the following computation:
\begin{eqnarray*}
\sum_{jb}(H\otimes K)_{ia,jb}
&=&\sum_{jb}H_{ij}K_{ab}\\
&=&\sum_jH_{ij}\sum_bK_{ab}\\
&=&\lambda\mu
\end{eqnarray*}

Thus, the matrix $H\otimes K$ is bistochastic as well. 

\medskip

(4) As for the last assertion, this follows from (1,2,3).
\end{proof}

In general now, putting an arbitrary complex Hadamard matrix in bistochastic form can be theoretically done, according to a general theorem of Idel-Wolf \cite{iwo}. The proof of this latter theorem is however based on a quite advanced, and non-explicit argument, coming from symplectic geometry, and there are many interesting open questions here.

\section*{7d. Hadamard conjecture}

As a final topic for this chapter, let us discuss now the real Hadamard matrices. The definition here, going back to 19th century work of Sylvester \cite{syl}, is as follows:

\index{Hadamard matrix}

\begin{definition}
A real Hadamard matrix is a square binary matrix, 
$$H\in M_N(\pm1)$$
whose rows are pairwise orthogonal, with respect to the scalar product on $\mathbb R^N$. 
\end{definition}

Observe that we do not really need real numbers in order to talk about the Hadamard matrices, because the orthogonality condition tells us that, when comparing two rows, the number of matchings should equal the number of mismatchings. 

\bigskip

As a first result regarding such matrices, we have:

\begin{proposition}
For a square matrix $H\in M_N(\pm1)$, the following are equivalent:
\begin{enumerate}
\item The rows of $H$ are pairwise orthogonal, and so $H$ is Hadamard.

\item The columns of $H$ are pairwise orthogonal, and so $H^t$ is Hadamard.

\item The rescaled matrix $U=H/\sqrt{N}$ is orthogonal, $U\in O_N$.
\end{enumerate}
\end{proposition}

\begin{proof}
This is something that we already know for the complex Hadamard matrices, with the orthogonal group $O_N$ being replaced by the unitary group $U_N$. In the real case the proof is similar, with everything coming from definitions, and linear algebra. 
\end{proof}

As an abstract consequence of the above result, let us record:

\begin{theorem}
The set of the $N\times N$ Hadamard matrices is
$$Y_N=M_N(\pm 1)\cap\sqrt{N}O_N$$
where $O_N$ is the orthogonal group, the intersection being taken inside $M_N(\mathbb R)$.
\end{theorem}

\begin{proof}
This follows from Proposition 7.25, which tells us that an arbitrary matrix $H\in M_N(\pm1)$ belongs to $Y_N$ if and only if it belongs to $\sqrt{N}O_N$.
\end{proof}

As a conclusion here, the set $Y_N$ that we are interested in appears as a kind of set of ``special rational points'' of the real algebraic manifold $\sqrt{N}O_N$. Moving now forward, as before in the complex matrix case, it is convenient to introduce:

\index{Hadamard equivalence}

\begin{definition}
Two real Hadamard matrices are called equivalent, and we write $H\sim K$, when it is possible to pass from $H$ to $K$ via the following operations:
\begin{enumerate}
\item Permuting the rows, or the columns.

\item Multiplying the rows or columns by $-1$.
\end{enumerate}
\end{definition}

Observe that we do not include the transposition operation $H\to H^t$ in our list of allowed operations. This is because Proposition 7.25, while looking quite elementary, rests however on a deep linear algebra fact, namely that the transpose of an orthogonal matrix is orthogonal as well, and this can produce complications later on.

\bigskip

Let us do now some classification work. Here is the result at $N=4$:

\begin{proposition}
There is only one Hadamard matrix at $N=4$, namely
$$W_4=W_2\otimes W_2$$
up to the standard equivalence relation for such matrices.
\end{proposition}

\begin{proof}
Consider an Hadamard matrix $H\in M_4(\pm1)$, assumed to be dephased:
$$H=\begin{pmatrix}1&1&1&1\\ 1&a&b&c\\ 1&d&e&f\\ 1&g&h&i\end{pmatrix}$$

By orthogonality of the first 2 rows we must have $\{a,b,c\}=\{-1,-1,1\}$, and so by permuting the last 3 columns, we can further assume that our matrix is as follows:
$$H=\begin{pmatrix}1&1&1&1\\ 1&-1&1&-1\\ 1&m&n&o\\ 1&p&q&r\end{pmatrix}$$

By orthogonality of the first 2 columns we must have $\{m,p\}=\{-1,1\}$, and so by permuting the last 2 rows, we can further assume that our matrix is as follows:
$$H=\begin{pmatrix}1&1&1&1\\ 1&-1&1&-1\\ 1&1&x&y\\ 1&-1&z&t\end{pmatrix}$$

Now from the orthogonality of the rows and columns we obtain 
$x=y=-1$, and then $z=-1,t=1$. Thus, up to equivalence we have $H=W_4$, as claimed.
\end{proof}

The case $N=5$ is excluded, because the orthogonality condition forces $N\in 2\mathbb N$. The point now is that the case $N=6$ is excluded as well, because we have:

\index{size of Hadamard matrix}

\begin{proposition}
The size of an Hadamard matrix must be 
$$N\in\{2\}\cup 4\mathbb N$$
with this coming from the orthogonality condition between the first $3$ rows.
\end{proposition}

\begin{proof}
By permuting the rows and columns or by multiplying them by $-1$, as to rearrange the first 3 rows, we can always assume that our matrix looks as follows:
$$H=\begin{pmatrix}
1\ldots\ldots 1&1\ldots\ldots 1&1\ldots\ldots 1&1\ldots\ldots 1\\
1\ldots\ldots 1&1\ldots\ldots 1&-1\ldots -1&-1\ldots -1\\
1\ldots\ldots 1&-1\ldots -1&1\ldots\ldots 1&-1\ldots -1\\
\underbrace{\ldots\ldots\ldots}_x&\underbrace{\ldots\ldots\ldots}_y&\underbrace{\ldots\ldots\ldots}_z&\underbrace{\ldots\ldots\ldots}_t
\end{pmatrix}$$

Now if we denote by $x,y,z,t$ the sizes of the 4 block columns, as indicated, the orthogonality conditions between the first 3 rows give the following system of equations:
$$(1\perp 2)\quad:\quad x+y=z+t$$
$$(1\perp 3)\quad:\quad x+z=y+t$$
$$(2\perp 3)\quad:\quad x+t=y+z$$

The numbers $x,y,z,t$ being such that the average of any two equals the average of the other two, and so equals the global average, the solution of our system is:
$$x=y=z=t$$

Thus the matrix size $N=x+y+z+t$ must be a multiple of 4, as claimed.
\end{proof}

The above result, and various other findings, suggest the following conjecture:

\index{Hadamard conjecture}
\index{HC}
\index{number of the beast}

\begin{conjecture}[Hadamard Conjecture (HC)]
There is at least one Hadamard matrix 
$$H\in M_N(\pm1)$$
for any integer $N\in 4\mathbb N$.
\end{conjecture}

This conjecture, going back to the 19th century, is one of the most beautiful statements in combinatorics, linear algebra, and mathematics in general. Quite remarkably, the numeric verification so far goes up to the number of the beast:
$$\mathfrak N=666$$

Our purpose now will be that of gathering some evidence for this conjecture. At $N=4,8$ we have the Walsh matrices $W_4,W_8$. Thus, the next existence problem comes at $N=12$. And here, we can use the following key construction, due to Paley:

\index{Paley matrices}

\begin{theorem}
Let $q=p^r$ be an odd prime power, define 
$$\chi:\mathbb F_q\to\{-1,0,1\}$$
by $\chi(0)=0$, $\chi(a)=1$ if $a=b^2$ for some $b\neq0$, and $\chi(a)=-1$ otherwise, and finally set 
$$Q_{ab}=\chi(a-b)$$
We have then constructions of Hadamard matrices, as follows:
\begin{enumerate}
\item Paley $1$: if $q=3(4)$ we have a matrix of size $N=q+1$, as follows:
$$P_N^1=1+\begin{pmatrix}
0&1&\ldots&1\\
-1\\
\vdots&&Q\\
-1
\end{pmatrix}$$

\item Paley $2$: if $q=1(4)$ we have a matrix of size $N=2q+2$, as follows:
$$P_N^2=\begin{pmatrix}
0&1&\ldots&1\\
1\\
\vdots&&Q\\
1
\end{pmatrix}\quad:\quad 0\to\begin{pmatrix}1&-1\\ -1&-1\end{pmatrix}\quad,\quad\pm1\to\pm\begin{pmatrix}1&1\\1&-1\end{pmatrix}$$
\end{enumerate}
These matrices are skew-symmetric $(H+H^t=2)$, respectively symmetric $(H=H^t)$.
\end{theorem}

\begin{proof}
In order to simplify the presentation, we will denote by $1$ all the identity matrices, of any size, and by $\mathbb I$ all the rectangular all-one matrices, of any size as well. It is elementary to check that the matrix $Q_{ab}=\chi(a-b)$ has the following properties:
$$QQ^t=q1-\mathbb I\quad,\quad 
Q\mathbb I=\mathbb IQ=0$$

In addition, we have the following formulae, which are elementary as well, coming from the fact that $-1$ is a square in $\mathbb F_q$ precisely when $q=1(4)$:
$$q=1(4)\implies Q=Q^t\ \ \,$$
$$q=3(4)\implies Q=-Q^t$$

With these observations in hand, the proof goes as follows:

\medskip

(1) With our conventions for the symbols $1$ and $\mathbb I$, the matrix in the statement is:
$$P_N^1=\begin{pmatrix}1&\mathbb I\\ -\mathbb I&1+Q\end{pmatrix}$$
 
With this formula in hand, the Hadamard matrix condition follows from:
\begin{eqnarray*}
P_N^1(P_N^1)^t
&=&\begin{pmatrix}1&\mathbb I\\ -\mathbb I&1+Q\end{pmatrix}\begin{pmatrix}1&-\mathbb I\\ \mathbb I&1-Q\end{pmatrix}\\
&=&\begin{pmatrix}N&0\\ 0&\mathbb I+1-Q^2\end{pmatrix}\\
&=&\begin{pmatrix}N&0\\ 0&N\end{pmatrix}
\end{eqnarray*}

(2) If we denote by $G,F$ the matrices in the statement, which replace respectively the $0,1$ entries, then we have the following formula for our matrix:
$$P_N^2=\begin{pmatrix}0&\mathbb I\\ \mathbb I&Q\end{pmatrix}\otimes F+1\otimes G$$

With this formula in hand, the Hadamard matrix condition follows from:
\begin{eqnarray*}
(P_N^2)^2
&=&\begin{pmatrix}0&\mathbb I\\ \mathbb I&Q\end{pmatrix}^2\otimes F^2+\begin{pmatrix}1&0\\ 0&1\end{pmatrix}\otimes G^2+\begin{pmatrix}0&\mathbb I\\ \mathbb I&Q\end{pmatrix}\otimes(FG+GF)\\
&=&\begin{pmatrix}q&0\\ 0&q\end{pmatrix}\otimes 2+\begin{pmatrix}1&0\\ 0&1\end{pmatrix}\otimes 2+\begin{pmatrix}0&\mathbb I\\ \mathbb I&Q\end{pmatrix}\otimes0\\
&=&\begin{pmatrix}N&0\\ 0&N\end{pmatrix}
\end{eqnarray*}

Finally, the last assertion is clear, from the above formulae relating $Q,Q^t$.
\end{proof}

The above constructions allow us to get well beyond the Walsh matrix level:

\index{HC}
\index{Paley matrices}

\begin{theorem}
The HC is verified at least up to $N=88$, as follows:
\begin{enumerate}
\item At $N=4,8,16,32,64$ we have Walsh matrices.

\item At $N=12,20,24,28,44,48,60,68,72,80,84,88$ we have Paley $1$ matrices.

\item At $N=36,52,76$ we have Paley $2$ matrices.

\item At $N=40,56$ we have Paley $1$ matrices tensored with $W_2$.
\end{enumerate}
\end{theorem}

\begin{proof}
First of all, the numbers in (1-4) are indeed all the multiples of 4, up to 88. As for the various assertions, the proof here goes as follows:

\medskip

(1) This is clear from the definition of the Walsh matrices.

\medskip

(2) Since $N-1$ takes the values $q=11,19,23,27,43,47,59,67,71,79,83,87$, all prime powers, we can indeed apply the Paley 1 construction, in all these cases.

\medskip

(3) Since $N=4(8)$ here, and $N/2-1$ takes the values $q=17,25,37$, all prime powers, we can indeed apply the Paley 2 construction, in these cases.

\medskip

(4) At $N=40$ we have indeed $P_{20}^1\otimes W_2$, and at $N=56$ we have $P_{28}^1\otimes W_2$.
\end{proof}

As a continuation of all this, at $N=92$ we have $92-1=7\times13$, so the Paley 1 construction does not work, and $92/2=46$, so the Paley 2 construction, or tensoring with $W_2$, does not work either. However, we can use here the following result:

\begin{theorem}
Assuming that $A,B,C,D\in M_K(\pm1)$ are circulant, symmetric, pairwise commute and satisfy the condition
$$A^2+B^2+C^2+D^2=4K$$
the following $4K\times4K$ matrix is Hadamard, called of Williamson type:
$$H=\begin{pmatrix}
A&B&C&D\\
-B&A&-D&C\\
-C&D&A&-B\\
-D&-C&B&A
\end{pmatrix}$$
Moreover, matrices $A,B,C,D$ as above exist at $K=23$, where $4K=92$.
\end{theorem}

\begin{proof}
Consider the quaternion units $1,i,j,k\in M_4(0,1)$, which describe the positions of the $A,B,C,D$ entries in the matrix $H$ from the statement. We have then:
$$H=A\otimes 1+B\otimes i+C\otimes j+D\otimes k$$

Assuming now that $A,B,C,D$ are symmetric, we have:
\begin{eqnarray*}
HH^t
&=&(A\otimes 1+B\otimes i+C\otimes j+D\otimes k)\\
&&(A\otimes 1-B\otimes i-C\otimes j-D\otimes k)\\
&=&(A^2+B^2+C^2+D^2)\otimes 1-([A,B]-[C,D])\otimes i\\
&&-([A,C]-[B,D])\otimes j-([A,D]-[B,C])\otimes k
\end{eqnarray*}

Now assume that our matrices $A,B,C,D$ pairwise commute, and satisfy the condition in the statement. In this case, it follows from the above formula that we have:
$$HH^t=4K$$

Thus, we obtain indeed an Hadamard matrix, as claimed. However, finding such matrices is in general a difficult task, and this is where Williamson's extra assumption in the statement, that $A,B,C,D$ should be taken circulant, comes from. Finally, regarding the $K=23$ and $N=92$ example, this comes via a computer search. 
\end{proof}

Things get even worse at higher values of $N$, where more and more complicated constructions are needed. The whole subject is quite technical, and, as already mentioned, human knowledge here stops so far at the number of the beast, namely:
$$\mathfrak N=666$$

Switching topics now, another well-known open question concerns the circulant case. Given a binary vector $\gamma\in(\pm 1)^N$, one can ask whether the matrix $H\in M_N(\pm 1)$ defined by $H_{ij}=\gamma_{j-i}$ is Hadamard or not. Here is a solution to the problem:
$$K_4=\begin{pmatrix}-1&1&1&1\\ 1&-1&1&1\\ 1&1&-1&1\\ 1&1&1&-1\end{pmatrix}$$

More generally, any vector $\gamma\in(\pm 1)^4$ satisfying $\sum\gamma_i=\pm 1$ is a solution to the problem. The following conjecture, from the 50s, states that there are no other solutions:

\index{Circulant Hadamard conjecture}
\index{CHC}

\begin{conjecture}[Circulant Hadamard Conjecture (CHC)]
The only Hadamard matrices which are circulant are
$$K_4=\begin{pmatrix}-1&1&1&1\\ 1&-1&1&1\\ 1&1&-1&1\\ 1&1&1&-1\end{pmatrix}$$
and its conjugates, regardless of the value of $N\in\mathbb N$.
\end{conjecture}

The fact that such a simple-looking problem is still open might seem quite surprising. Indeed, if we denote by $S\subset\{1,\ldots,N\}$ the set of positions of the $-1$ entries of $\gamma$, the Hadamard matrix condition is simply, for any $k\neq 0$, taken modulo $N$:
$$|S\cap(S+k)|=|S|-N/4$$ 

Thus, the above conjecture simply states that at $N\neq 4$, such a set $S$ cannot exist. This is a well-known problem in combinatorics, raised by Ryser a long time ago.

\bigskip

Summarizing, we have many interesting questions in the real case. The situation is quite different from the one in complex case, where at any $N\in\mathbb N$ we have the Fourier matrix $F_N$, which makes the HC problematics dissapear. Since $F_N$ can be put in circulant form, the CHC dissapears as well. There are however many interesting questions in the complex case, for the most in relation with questions in quantum physics.

\section*{7e. Exercises}

We have learned many interesting things in this chapter, and our exercises will focus on the complex Hadamard matrices, which were the central objects, in all this. First, we have the following standard fact, dealing with deformations of such matrices:

\begin{exercise}
If $H\in M_M(\mathbb T)$ and $K\in M_N(\mathbb T)$ are Hadamard matrices, so is 
$$H\otimes_QK\in M_{MN}(\mathbb T)$$
given by the following formula, with $Q\in M_{M\times N}(\mathbb T)$,
$$(H\otimes_QK)_{ia,jb}=Q_{ib}H_{ij}K_{ab}$$
called Di\c t\u a deformation of $H\otimes K$, with parameter $Q$.
\end{exercise}

Normally this is just a quick, standard verification. More difficult, however, is the question of explicitly writing down the matrices that can be constructed in this way, because this requires things like struggling with double indices. Good luck here.

\begin{exercise}
Prove that the only complex Hadamard matrices at $N=4$ are, up to the standard equivalence relation, the matrices 
$$F_4^q=\begin{pmatrix}
1&1&1&1\\
1&-1&1&-1\\
1&q&-1&-q\\ 
1&-q&-1&q
\end{pmatrix}$$
with $q\in\mathbb T$, which appear as Di\c t\u a deformations of $W_4=F_2\otimes F_2$. 
\end{exercise}

Here the first question is quite standard, in the spirit of the computations at $N=3$, mentioned before. As for the second question, good luck here with the double indices.

\begin{exercise}
Given an Hadamard matrix $H\in M_5(\mathbb T)$, chosen dephased,
$$H=\begin{pmatrix}
1&1&1&1&1\\
1&a&x&*&*\\
1&y&b&*&*\\
1&*&*&*&*\\
1&*&*&*&*
\end{pmatrix}$$
prove that the numbers $a,b,x,y$ must satisfy $(x-y)(x-ab)(y-ab)=0$.
\end{exercise}

This is something quite tricky, called Haagerup lemma, and in case you're stuck with this, you can of course take a look at Haagerup's paper \cite{haa}. As bonus exercise, using this lemma, work out the full details of the classification at $N=5$.

\chapter{Infinite dimensions}

\section*{8a. Hilbert spaces}

We have seen so far the basics of linear algebra, concerning linear maps and matrices, the determinant, the diagonalization procedure, and some applications. In this chapter, motivated by quantum mechanics, we discuss what happens in infinite dimensions.

\bigskip

To be more precise, among the main discoveries of the 1920s, due to Heisenberg, Schr\"odinger and others was the fact that small particles like electrons cannot really be described by their position vectors $v\in\mathbb R^3$, and instead we must use their so-called wave functions $\psi:\mathbb R^3\to\mathbb C$. Thus, the natural space for quantum mechanics, or at least for the quantum mechanics of the 1920s, is not our usual $V=\mathbb R^3$, but rather the infinite dimensional space $H=L^2(\mathbb R^3)$ of such wave functions $\psi$. And more recent versions of quantum mechanics are built on the same idea, namely infinite dimensional spaces.

\bigskip

Getting started now, we would like to look at linear algebra over infinite dimensional spaces. However, this is not very interesting, due to a number of technical reasons, the idea being that the infinite dimensionality prevents us from doing many basic things, to the point that we cannot even have things started. So, the idea will be that of using infinite dimensional vector spaces with some extra structure, as follows: 

\index{scalar product}

\begin{definition}
A scalar product on a complex vector space $H$ is an operation 
$$H\times H\to\mathbb C$$
denoted $(x,y)\to <x,y>$, satisfying the following conditions:
\begin{enumerate}
\item $<x,y>$ is linear in $x$, and antilinear in $y$.

\item $\overline{<x,y>}=<y,x>$, for any $x,y$.

\item $<x,x>>0$, for any $x\neq0$.
\end{enumerate}
\end{definition}

As a basic example here, we have the finite dimensional vector space $H=\mathbb C^N$, with its usual scalar product, which is as follows:
$$<x,y>=\sum_ix_i\bar{y}_i$$ 

There are many other examples, and notably various spaces of $L^2$ functions, which naturally appear in problems coming from physics. We will discuss them later.

\bigskip

In order to study the scalar products, let us formulate the following definition:

\index{norm of vector}

\begin{definition}
The norm of a vector $x\in H$ is the following quantity:
$$||x||=\sqrt{<x,x>}$$
We also call this number length of $x$, or distance from $x$ to the origin.
\end{definition}

In analogy with what happens in finite dimensions, we have two important results regarding the norms. First is the Cauchy-Schwarz inequality, as follows:

\index{Cauchy-Schwarz inequality}

\begin{theorem}
We have the Cauchy-Schwarz inequality
$$|<x,y>|\leq||x||\cdot||y||$$
and the equality case holds precisely when $x,y$ are proportional.
\end{theorem}

\begin{proof}
Consider the following quantity, depending on a real variable $t\in\mathbb R$, and on a variable on the unit circle, $w\in\mathbb T$:
$$f(t)=||twx+y||^2$$

By developing $f$, we can see that this is a degree 2 polynomial in $t$:
\begin{eqnarray*}
f(t)
&=&<twx+y,twx+y>\\
&=&t^2<x,x>+tw<x,y>+t\bar{w}<y,x>+<y,y>\\
&=&t^2||x||^2+2tRe(w<x,y>)+||y||^2
\end{eqnarray*}

Since $f$ is obviously positive, its discriminant must be negative:
$$4Re(w<x,y>)^2-4||x||^2\cdot||y||^2\leq0$$

But this is equivalent to the following condition:
$$|Re(w<x,y>)|\leq||x||\cdot||y||$$

Now the point is that we can arrange for the number $w\in\mathbb T$ to be such that the quantity $w<x,y>$ is real. Thus, we obtain the Cauchy-Schwarz inequality:
$$|<x,y>|\leq||x||\cdot||y||$$

Finally, the study of the equality case is straightforward, by using the fact that the discriminant of $f$ vanishes precisely when we have a root. But this leads to the conclusion in the statement, namely that the vectors $x,y$ must be proportional. 
\end{proof}
 
As a second main result now, we have the Minkowski inequality:

\index{Minkovski inequality}

\begin{theorem}
We have the Minkowski inequality
$$||x+y||\leq||x||+||y||$$
and the equality case holds precisely when $x,y$ are proportional.
\end{theorem}

\begin{proof}
This follows indeed from the Cauchy-Schwarz inequality, as follows:
\begin{eqnarray*}
&&||x+y||\leq||x||+||y||\\
&\iff&||x+y||^2\leq(||x||+||y||)^2\\
&\iff&||x||^2+||y||^2+2Re<x,y>\leq||x||^2+||y||^2+2||x||\cdot||y||\\
&\iff&Re<x,y>\leq||x||\cdot||y||
\end{eqnarray*}

As for the equality case, this is clear from Cauchy-Schwarz as well.
\end{proof}
  
As a consequence of this, we have the following result:

\index{distance}
\index{metric space}

\begin{theorem}
The following function is a distance on $H$,
$$d(x,y)=||x-y||$$
in the usual sense, that of the abstract metric spaces.
\end{theorem}

\begin{proof}
This follows indeed from the Minkowski inequality, which corresponds to the triangle inequality, the other two axioms for a distance being trivially satisfied.
\end{proof}

The above result is quite important, because it shows that we can do geometry in our present setting, a bit as in the finite dimensional case. Still in connection with this, doing geometry, we have the following key technical result, which can be very useful:

\index{polarization identity}

\begin{proposition}
The scalar products can be recovered from distances, via the formula
$$4<x,y>
=||x+y||^2-||x-y||^2+i||x+iy||^2-i||x-iy||^2$$
called complex polarization identity.
\end{proposition}

\begin{proof}
This is something that we already met before, in finite dimensions. In arbitrary dimensions the proof is similar, as follows:
\begin{eqnarray*}
&&||x+y||^2-||x-y||^2+i||x+iy||^2-i||x-iy||^2\\
&=&||x||^2+||y||^2-||x||^2-||y||^2+i||x||^2+i||y||^2-i||x||^2-i||y||^2\\
&&+2Re(<x,y>)+2Re(<x,y>)+2iIm(<x,y>)+2iIm(<x,y>)\\
&=&4<x,y>
\end{eqnarray*}

Thus, we are led to the conclusion in the statement.
\end{proof}

Let us discuss now some more advanced aspects. In order to do analysis on our spaces, we need the Cauchy sequences that we construct to converge. This is something which is automatic in finite dimensions, but in arbitrary dimensions, this can fail.

\bigskip

Thus, we must add an extra axiom, stating that our vector space $H$ is complete with respect to the norm. It is convenient here to formulate a detailed new definition, as follows, which will be the starting point for our various considerations to follow:

\index{Hilbert space}
\index{scalar product}

\begin{definition}
A Hilbert space is a complex vector space $H$ given with a scalar product $<x,y>$, satisfying the following conditions:
\begin{enumerate}
\item $<x,y>$ is linear in $x$, and antilinear in $y$.

\item $\overline{<x,y>}=<y,x>$, for any $x,y$.

\item $<x,x>>0$, for any $x\neq0$.

\item $H$ is complete with respect to the norm $||x||=\sqrt{<x,x>}$.
\end{enumerate}
\end{definition}

In other words, we have taken here Definition 8.1, and added the condition that $H$ must be complete with respect to the norm $||x||=\sqrt{<x,x>}$, that we know indeed to be a norm, according to the Minkowski inequality proved above. 

\bigskip

As a basic example, since in finite dimensions the completness axiom is automatically satisfied, we have as before the space $H=\mathbb C^N$, with its usual scalar product:
$$<x,y>=\sum_ix_i\bar{y}_i$$

More generally now, we have the following construction of Hilbert spaces:

\index{square-summable}

\begin{proposition}
The sequences of numbers $x=(x_i)$ which are square-summable,
$$\sum_i|x_i|^2<\infty$$
form a Hilbert space, denoted $l^2(\mathbb N)$, with the following scalar product:
$$<x,y>=\sum_ix_i\bar{y}_i$$
In fact, given any index set $I$, we can construct a Hilbert space $l^2(I)$, in this way.
\end{proposition}

\begin{proof}
The fact that we have indeed a complex vector space with a scalar product is elementary, and the fact that this space is indeed complete is very standard too. We will leave all the verifications here, which are straightforward, as an exercise.
\end{proof}

On the other hand, we can talk as well about spaces of functions, as follows:

\begin{proposition}
Given an interval $X\subset\mathbb R$, the quantity 
$$<f,g>=\int_Xf(x)\overline{g(x)}dx$$
is a scalar product, making $H=L^2(X)$ a Hilbert space.
\end{proposition}

\begin{proof}
Once again this is routine, coming this time from basic measure theory, with  $H=L^2(X)$ being the space of square-integrable functions $f:X\to\mathbb C$, with the convention that two such functions are identified when they coincide almost everywhere. 
\end{proof}

The point now is that we can unify the above two constructions, as follows:

\begin{theorem}
Given a measured space $X$, the quantity 
$$<f,g>=\int_Xf(x)\overline{g(x)}dx$$
is a scalar product, making $H=L^2(X)$ a Hilbert space.
\end{theorem}

\begin{proof}
Here the first assertion is clear, and the fact that the Cauchy sequences converge is clear as well, by taking the pointwise limit, and using a standard argument. As before with our previous such results, we will leave the verifications here as an exercise.
\end{proof}

Observe that with $X=\{1,\ldots,N\}$  we obtain the space $H=\mathbb C^N$. Also, with $X=\mathbb N$, with the counting measure, we obtain the space $H=l^2(\mathbb N)$. In fact, with an arbitrary set $I$, once again with the counting mesure, we obtain the space $H=l^2(I)$. Thus, the construction in Theorem 8.10 unifies all the Hilbert space constructions that we have.

\bigskip

Quite remarkably, the converse of this holds, in the sense that any Hilbert space must be of the form $L^2(X)$. This follows indeed from the following key result, which tells us that, in addition to this, we can always assume that $X=I$ is a discrete space:

\index{algebraic basis}
\index{orthonormal basis}
\index{Gram-Schmidt}

\begin{theorem}
Let $H$ be a Hilbert space.
\begin{enumerate}
\item Any algebraic basis of this space $\{f_i\}_{i\in I}$ can be turned into an orthonormal basis $\{e_i\}_{i\in I}$, by using the Gram-Schmidt procedure.

\item Thus, $H$ has an orthonormal basis, and so we have $H\simeq l^2(I)$, with $I$ being the indexing set for this orthonormal basis.
\end{enumerate}
\end{theorem}

\begin{proof}
There are several things going on here, the idea being as follows:

\medskip

(1) In finite dimensions, we can turn any vector space basis $\{f_i\}_{i\in I}$ into an orthogonal basis $\{e_i\}_{i\in I}$, by using the Gram-Schmidt procedure, as follows, with $\alpha_i,\beta_i,\gamma_i,\ldots$ being uniquely determined by the fact at each step, $e_k$ must be orthogonal to $f_1,\ldots,f_{k-1}$:
$$e_1=f_1$$
$$e_2=f_2+\alpha_1f_1$$
$$e_3=f_3+\beta_1f_1+\beta_2f_2$$
$$e_4=f_4+\gamma_1f_1+\gamma_2f_2+\gamma_3f_3$$
$$\vdots$$

And then, by replacing $e_i\to e_i/||e_i||$, we have our orthonormal basis, as desired.

\medskip

(2) In general, the same method works, namely Gram-Schmidt, with a subtlety coming from the fact that the basis $\{e_i\}_{i\in I}$ will not span in general the whole $H$, but just a dense subspace of it, as it is in fact obvious by looking at the standard basis of $l^2(\mathbb N)$. 

\medskip

(3) And there is a second subtlety as well, coming from the fact that the recurrence procedure needed for Gram-Schmidt must be replaced by some sort of ``transfinite recurrence'', using standard tools from logic, and more specifically the Zorn lemma.
\end{proof}

We have the following definition, based on the above:

\index{separable Hilbert space}

\begin{definition}
A Hilbert space $H$ is called separable when the following equivalent conditions are satisfied:
\begin{enumerate}
\item $H$ has a countable algebraic basis $\{f_i\}_{i\in\mathbb N}$.

\item $H$ has a countable orthonormal basis $\{e_i\}_{i\in\mathbb N}$.

\item We have $H\simeq l^2(\mathbb N)$, isomorphism of Hilbert spaces.
\end{enumerate}
\end{definition}

As a main question now, are the Hilbert spaces coming from quantum mechanics, such as the Schr\"odinger space $H=L^2(\mathbb R^3)$ of wave functions of the electron, separable? In answer, up to some simple operations, involving tensor products and stretching, we must solve the question for $H=L^2[0,1]$. And here, following Weierstrass, we have:

\index{orthogonal polynomials}

\begin{theorem}
The following happen, regarding the functions $f:[0,1]\to\mathbb C$:
\begin{enumerate}
\item Any continuous function $f:[0,1]\to\mathbb C$ can be uniformly approximated by polynomials. Thus, $\{x^n\}_{n\in\mathbb N}$ is an algebraic basis of the space $L^2[0,1]$.

\item By applying Gram-Schmidt we obtain certain polynomials $\{L_n\}_{n\in\mathbb N}$, the modified Legendre polynomials, which give an explicit isomorphism $L^2[0,1]\simeq l^2(\mathbb N)$.
\end{enumerate}
\end{theorem}

\begin{proof}
This is something very classical, the idea being as follows:

\medskip

(1) Consider the following polynomials, called Bernstein polynomials:
$$b_{kn}(x)=\binom{n}{k}x^k(1-x)^{n-k}$$

Then, given $f:[0,1]\to\mathbb R$ continuous, consider the following polynomials:
$$f_n(x)=\sum_{k=0}^nf\left(\frac{k}{n}\right)b_{kn}(x)$$

Our claim is that we have $f_n\to_uf$, uniform convergence on $[0,1]$.

\medskip

(2) In order to prove this, observe that the polynomials $b_{kn}$ encode the densities of the binomial laws $\rho_{xn}$. Thus, we have the following formulae, with the first one corresponding to the fact that $\rho_{xn}$ is indeed a probability measure, and with the second and third formulae coming from our mean and variance computations from chapter 6:
$$\sum_{k=0}^nb_{kn}(x)=1$$
$$\sum_{k=0}^n\frac{k}{n}\cdot b_{kn}(x)=x$$
$$\sum_{k=0}^n\left(x-\frac{k}{n}\right)^2b_{kn}(x)=\frac{x(1-x)}{n}$$

(3) In order to estimate now the error $|f_n-f|$, we can use the uniform continuity property of $f$. So, pick $\varepsilon>0$, and then $\delta>0$ such that the following happens:
$$|x-y|<\delta\implies|f(x)-f(y)|<\varepsilon$$

(4) We have then the following estimate, using this, and with $M=\sup|f|$:
\begin{eqnarray*}
|f_n(x)-f(x)|
&=&\left|\sum_{k=0}^nf\left(\frac{k}{n}\right)b_{kn}(x)-\sum_{k=0}^nf(x)b_{kn}(x)\right|\\
&\leq&\sum_{k=0}^n\left|f\left(\frac{k}{n}\right)-f(x)\right|b_{kn}(x)\\
&=&\sum_{\left|x-\frac{k}{n}\right|<\delta}\left|f\left(\frac{k}{n}\right)-f(x)\right|b_{kn}(x)
+\sum_{\left|x-\frac{k}{n}\right|\geq\delta}\left|f\left(\frac{k}{n}\right)-f(x)\right|b_{kn}(x)\\
&\leq&\varepsilon+M\sum_{\left|x-\frac{k}{n}\right|\geq\delta}b_{kn}(x)
\end{eqnarray*}

(5) In order to deal with the sum on the right, we will need some standard estimates. Let us first recall the Markov inequality, which is something trivial, as follows:
$$P\big(|\varphi|\geq b\big)\leq\frac{E(\varphi)}{b}$$

By using this with $\varphi=(\psi-E)^2$, with $E=E(\psi)$, we obtain the Chebycheff inequality:
$$P\big(|\psi-E|\geq a\big)
\leq\frac{E((\psi-E)^2)}{a^2}
=\frac{V}{a^2}$$

(6) The point now is that this latter inequality applies to the last sum in (4), with $\psi$ being a variable following the binomial law $\rho_{xn}$, rescaled to $[0,1]$, and gives:
\begin{eqnarray*}
\sum_{\left|x-\frac{k}{n}\right|\geq\delta}b_{kn}(x)
&\leq&\sum_{k=0}^n\delta^{-2}\left(x-\frac{k}{n}\right)^2b_{kn}(x)\\
&=&\delta^{-2}\,\frac{x(1-x)}{n}\\
&\leq&\frac{\delta^{-2}}{4n}
\end{eqnarray*}

(7) Now by putting everything together, we obtain the following estimate:
$$|f_n(x)-f(x)|\leq \varepsilon+\frac{\delta^{-2}M}{4n}$$

Thus we have indeed $|f_n-f|\to0$, uniform convergence, as desired. Finally, in what regards orthogonalization, we will leave some learning here as an exercise.
\end{proof}

As a conclusion to all this, we are interested in 1 space, namely the unique separable Hilbert space $H$, but due to various technical reasons, it is often better to forget that we have $H=l^2(\mathbb N)$, and say instead that we have $H=L^2(X)$, with $X$ being a separable measured space, or simply say that $H$ is an abstract separable Hilbert space.

\section*{8b. Linear operators}

Let us get now into the study of linear operators $T:H\to H$, which will eventually lead us into the correct infinite dimensional version of linear algebra. We first have:

\begin{proposition}
For a linear operator $T:H\to H$, the following are equivalent:
\begin{enumerate}
\item $T$ is continuous.

\item $T$ is continuous at $0$.

\item $T(B)\subset cB$ for some $c<\infty$, where $B\subset H$ is the unit ball.

\item $T$ is bounded, in the sense that $||T||=\sup_{||x||\leq1}||Tx||$ satisfies $||T||<\infty$.
\end{enumerate}
\end{proposition}

\begin{proof}
This is something elementary, the idea being as follows:

\medskip

$(1)\iff(2)$ This is indeed clear from the linearity of $T$.

\medskip

$(2)\iff(3)$ This is again something clear, coming from definitions.

\medskip

$(3)\iff(4)$ Again, this is clear, with the number $||T||$ appearing in (4) being the infimum of the numbers $c$ making the condition (3) work.

\medskip

$(4)\iff(1)$ This is something clear too, coming from the definiton of continuity.
\end{proof}

Regarding now the bounded operators, we have the following result, about them:

\index{operator algebra}
\index{complex algebra}
\index{Banach algebra}

\begin{theorem}
The linear operators $T:H\to H$ which are bounded,
$$||T||=\sup_{||x||\leq1}||Tx||<\infty$$
form a complex algebra with unit $B(H)$, having the property 
$$||ST||\leq||S||\cdot||T||$$
and which is complete with respect to the norm.
\end{theorem}

\begin{proof}
The fact that we have indeed an algebra, satisfying the product condition in the statement, follows from the following estimates, which are all elementary:
$$||S+T||\leq||S||+||T||\quad,\quad 
||\lambda T||=|\lambda|\cdot||T||\quad,\quad
||ST||\leq||S||\cdot||T||$$

Summarizing, we have indeed an algebra, satisfying the product condition in the statement. Regarding now the last assertion, if $\{T_n\}\subset B(H)$ is Cauchy then $\{T_nx\}$ is Cauchy for any $x\in H$, so we can define the limit $T=\lim_{n\to\infty}T_n$ by setting:
$$Tx=\lim_{n\to\infty}T_nx$$

Let us first check that the application $x\to Tx$ is linear. We have:
\begin{eqnarray*}
T(x+y)
&=&\lim_{n\to\infty}T_n(x+y)\\
&=&\lim_{n\to\infty}T_n(x)+T_n(y)\\
&=&\lim_{n\to\infty}T_n(x)+\lim_{n\to\infty}T_n(y)\\
&=&T(x)+T(y)
\end{eqnarray*}

Similarly, we have as well the following computation:
\begin{eqnarray*}
T(\lambda x)
&=&\lim_{n\to\infty}T_n(\lambda x)\\
&=&\lambda\lim_{n\to\infty}T_n(x)\\
&=&\lambda T(x)
\end{eqnarray*}

Thus we have a linear map $T:A\to A$. It remains to prove that we have $T\in B(H)$, and that we have $T_n\to T$ in norm. For this purpose, observe that we have:
\begin{eqnarray*}
&&||T_n-T_m||\leq\varepsilon\ ,\ \forall n,m\geq N\\
&\implies&||T_nx-T_mx||\leq\varepsilon\ ,\ \forall||x||=1\ ,\ \forall n,m\geq N\\
&\implies&||T_nx-Tx||\leq\varepsilon\ ,\ \forall||x||=1\ ,\ \forall n\geq N\\
&\implies&||T_Nx-Tx||\leq\varepsilon\ ,\ \forall||x||=1\\
&\implies&||T_N-T||\leq\varepsilon
\end{eqnarray*}

As a first consequence, we obtain $T\in B(H)$, because we have:
\begin{eqnarray*}
||T||
&=&||T_N+(T-T_N)||\\
&\leq&||T_N||+||T-T_N||\\
&\leq&||T_N||+\varepsilon\\
&<&\infty
\end{eqnarray*}

As a second consequence, we obtain $T_N\to T$ in norm, and we are done.
\end{proof}

As a useful complement to the above result, in the presence of a basis, we have:

\index{linear operator}
\index{infinite matrix}

\begin{theorem}
Let $H$ be a Hilbert space, with orthonormal basis $\{e_i\}_{i\in I}$. The bounded operators $T\in B(H)$ can be then identified with matrices $M\in M_I(\mathbb C)$ via
$$Tx=Mx\quad,\quad M_{ij}=<Te_j,e_i>$$
and we obtain in this way an embedding as follows, which is multiplicative:
$$B(H)\subset M_I(\mathbb C)$$
In the case $H=\mathbb C^N$ we obtain in this way the usual isomorphism $B(H)\simeq M_N(\mathbb C)$. In the separable case we obtain in this way a proper embedding $B(H)\subset M_\infty(\mathbb C)$.
\end{theorem}

\begin{proof}
We have several assertions to be proved, the idea being as follows:

\medskip

(1) Regarding the first assertion, given a bounded operator $T:H\to H$, let us associate to it a matrix $M\in M_I(\mathbb C)$ as in the statement, by the following formula:
$$M_{ij}=<Te_j,e_i>$$

It is clear that this correspondence $T\to M$ is linear, and also that its kernel is $\{0\}$. Thus, we have an embedding of linear spaces $B(H)\subset M_I(\mathbb C)$.

\medskip

(2) Our claim now is that this embedding is multiplicative. But this is clear too, because if we denote by $T\to M_T$ our correspondence, we have:
\begin{eqnarray*}
(M_{ST})_{ij}
&=&<STe_j,e_i>\\
&=&\left<S\sum_k<Te_j,e_k>e_k,e_i\right>\\
&=&\sum_k<Se_k,e_i><Te_j,e_k>\\
&=&\sum_k(M_S)_{ik}(M_T)_{kj}\\
&=&(M_SM_T)_{ij}
\end{eqnarray*}

(3) Finally, we must prove that the original operator $T:H\to H$ can be recovered from its matrix $M\in M_I(\mathbb C)$ via the formula in the statement, namely $Tx=Mx$. But this latter formula holds for the vectors of the basis, $x=e_j$, because we have:
$$(Te_j)_i
=<Te_j,e_i>
=M_{ij}
=(Me_j)_i$$

Now by linearity we obtain from this that the formula $Tx=Mx$ holds everywhere, on any vector $x\in H$, and this finishes the proof of the first assertion.

\medskip

(4) In finite dimensions we obtain of course an isomorphism, and this because any usual matrix $M\in M_N(\mathbb C)$ determines a linear operator $T:\mathbb C^N\to\mathbb C^N$, according to the formula $<Te_j,e_i>=M_{ij}$. In infinite dimensions, however, we do not have an isomorphism. For instance on $H=l^2(\mathbb N)$ the following matrix does not define a linear operator:
$$M=\begin{pmatrix}
1&1&1&\ldots\\
1&1&1&\ldots\\
1&1&1&\ldots\\
\vdots&\vdots&\vdots
\end{pmatrix}$$

Thus, we are led to the conclusions in the statement.
\end{proof}

As a third and last main result about the bounded operators, we have:

\index{adjoint operator}

\begin{theorem}
The normed algebra $B(H)$ has an involution $T\to T^*$, given by 
$$<Tx,y>=<x,T^*y>$$
which is antilinear, antimultiplicative, and is an isometry, in the sense that:
$$||T||=||T^*||$$
Moreover, the norm the involution are related as well by $||TT^*||=||T||^2$.
\end{theorem}

\begin{proof}
We have several things to be proved, the idea being as follows:

\medskip

(1) As a preliminary fact, that we will need in what follows, our claim is that any linear form $\varphi:H\to\mathbb C$ must be of the following type, for a certain vector $z\in H$:
$$\varphi(x)=<x,z>$$

Indeed, this is something clear for any Hilbert space of type $H=l^2(I)$. But, by using a basis, any Hilbert space is of this form, and so we have proved our claim.

\medskip

(2) The existence of the adjoint operator $T^*$, given by the formula in the statement, comes from the fact that the function $\varphi(x)=<Tx,y>$ being a linear map $H\to\mathbb C$, we must have a formula as follows, for a certain vector $T^*y\in H$:
$$\varphi(x)=<x,T^*y>$$

Moreover, since this vector is unique, $T^*$ is unique too, and we have as well:
$$(S+T)^*=S^*+T^*\quad,\quad 
(\lambda T)^*=\bar{\lambda}T^*$$
$$(ST)^*=T^*S^*\quad,\quad
(T^*)^*=T$$

Observe also that we have indeed $T^*\in B(H)$, because:
\begin{eqnarray*}
||T||
&=&\sup_{||x||=1}\sup_{||y||=1}<Tx,y>\\
&=&\sup_{||y||=1}\sup_{||x||=1}<x,T^*y>\\
&=&||T^*||
\end{eqnarray*}

(3) Regarding now the last assertion, observe that we have:
$$||TT^*||
\leq||T||\cdot||T^*||
=||T||^2$$

On the other hand, we have as well the following estimate:
\begin{eqnarray*}
||T||^2
&=&\sup_{||x||=1}|<Tx,Tx>|\\
&=&\sup_{||x||=1}|<x,T^*Tx>|\\
&\leq&||T^*T||
\end{eqnarray*}

By replacing $T\to T^*$ we obtain from this that we have as well $||T||^2\leq||TT^*||$. Thus, we have obtained the needed inequality, and we are done.
\end{proof}

As an observation here, in the context of the construction $T\to M$ from Theorem 8.16, the adjoint operation $T\to T^*$ takes a very simple form, namely:
$$(M^*)_{ij}=\overline{M}_{ji}$$ 

However, this is a bit theoretical, because for spaces like $L^2[0,1]$, which do not have a simple orthonormal basis, the embedding $B(H)\subset M_I(\mathbb C)$ that we have is not very concrete. Thus, while the bounded operators $T:H\to H$ are basically some infinite matrices, it is better to think of these operators as being objects on their own.

\section*{8c. Spectral theory}

We will be interested in what follows in the algebra $B(H)$, and its closed subalgebras $A\subset B(H)$. It is convenient to formulate the following definition:

\index{Banach algebra}

\begin{definition}
A Banach algebra is a complex algebra with unit $A$, having a vector space norm $||.||$ satisfying
$$||ab||\leq||a||\cdot||b||$$
and which makes it a Banach space, in the sense that the Cauchy sequences converge.
\end{definition}

As said above, the basic examples of Banach algebras, or at least the basic examples that we will be interested in here, are the operator algebra $B(H)$, and its norm closed subalgebras $A\subset B(H)$, such as the algebras $A=<T>$ generated by a single operator $T\in B(H)$. There are many other examples, and more on this later.

\bigskip

Generally speaking, the elements $a\in A$ of a Banach algebra can be thought of as being bounded operators on some Hilbert space, which is not present. With this idea in mind, we can emulate spectral theory in our setting, the starting point being:

\index{spectrum}

\begin{definition}
The spectrum of an element $a\in A$ is the set
$$\sigma(a)=\left\{\lambda\in\mathbb C\Big|a-\lambda\not\in A^{-1}\right\}$$
where $A^{-1}\subset A$ is the set of invertible elements.
\end{definition}

As a basic example, the spectrum of a usual matrix $M\in M_N(\mathbb C)$ is the collection of its eigenvalues, taken of course without multiplicities. In the case of the trivial algebra $A=\mathbb C$, appearing at $N=1$, the spectrum of an element is the element itself.

\bigskip

As a first, basic result regarding spectra, we have:

\index{shift}

\begin{proposition}
We have the following formula, valid for any $a,b\in A$:
$$\sigma(ab)\cup\{0\}=\sigma(ba)\cup\{0\}$$
Also, there are examples where $\sigma(ab)\neq\sigma(ba)$.
\end{proposition}

\begin{proof}
We will first prove that we have the following implication:
$$1\notin\sigma(ab)\implies1\notin\sigma(ba)$$

For this purpose, assume that $1-ab$ is invertible, with inverse denoted $c$:
$$c=(1-ab)^{-1}$$

We have then the following formulae, relating our variables $a,b,c$:
$$abc=cab=c-1$$

By using these formulae, we obtain the following equality:
\begin{eqnarray*}
(1+bca)(1-ba)
&=&1+bca-ba-bcaba\\
&=&1+bca-ba-bca+ba\\
&=&1
\end{eqnarray*}

On the other hand, a similar computation shows that we have as well:
$$(1-ba)(1+bca)=1$$

Thus $1-ba$ is invertible, with inverse $1+bca$, which proves our claim. Now by multiplying by scalars, we deduce from this that for any $\lambda\in\mathbb C-\{0\}$ we have:
$$\lambda\notin\sigma(ab)\implies\lambda\notin\sigma(ba)$$ 

But this leads to the conclusion in the statement, namely:
$$\sigma(ab)\cup\{0\}=\sigma(ba)\cup\{0\}$$

Regarding now the last claim, we know from linear algebra that $\sigma(ab)=\sigma(ba)$ holds for the usual matrices, for instance because of the above, and because $ab$ is invertible if any only if $ba$ is. However, this latter fact fails for general operators on Hilbert spaces. Indeed, we can take our operator $a$ to be the shift on the space $l^2(\mathbb N)$, given by:
$$S(e_i)=e_{i+1}$$

As for $b$, we can take the adjoint of $S$, which is the following operator:
$$S^*(e_i)=\begin{cases}
e_{i-1}&{\rm if}\ i>0\\
0&{\rm if}\ i=0
\end{cases}$$

Let us compose now these two operators. In one sense, we have:
$$S^*S=1\implies 0\notin\sigma(SS^*)$$

In the other sense, however, the situation is different, as follows:
$$SS^*=Proj(e_0^\perp)\implies 0\in\sigma(SS^*)$$

Thus, the spectra do not match on $0$, and we have our counterexample, as desired.
\end{proof}

Let us discuss now a second basic result about spectra, which is something very useful. Given an arbitrary Banach algebra element $a\in A$, and a rational function $f=P/Q$ having poles outside the spectrum $\sigma(a)$, we can construct the following element:
$$f(a)=P(a)Q(a)^{-1}$$

For simplicity, and due to the fact that the elements $P(a),Q(a)$ commute, so that the order is irrelevant, we write this element as a usual fraction, as follows:
$$f(a)=\frac{P(a)}{Q(a)}$$

With this convention, we have the following result:

\index{rational calculus}
\index{rational function}

\begin{theorem}
We have the ``rational functional calculus'' formula
$$\sigma(f(a))=f(\sigma(a))$$
valid for any rational function $f\in\mathbb C(X)$ having poles outside $\sigma(a)$.
\end{theorem}

\begin{proof}
In order to prove this result, we can proceed in two steps, as follows:

\medskip

(1) Assume first that we are in the polynomial function case, $f\in\mathbb C[X]$. We pick a scalar $\lambda\in\mathbb C$, and we decompose the polynomial $f-\lambda$ into factors:
$$f(X)-\lambda=c(X-r_1)\ldots(X-r_n)$$

By using this formula, we have then, as desired:
\begin{eqnarray*}
\lambda\notin\sigma(f(a))
&\iff&f(a)-\lambda\in A^{-1}\\
&\iff&c(a-r_1)\ldots(a-r_n)\in A^{-1}\\
&\iff&a-r_1,\ldots,a-r_n\in A^{-1}\\
&\iff&r_1,\ldots,r_n\notin\sigma(a)\\
&\iff&\lambda\notin f(\sigma(a))
\end{eqnarray*}

(2) Assume now that we are in the general rational function case, $f\in\mathbb C(X)$. We pick a scalar $\lambda\in\mathbb C$, we write $f=P/Q$, and we set:
$$F=P-\lambda Q$$

By using now what we found in (1), for this polynomial, we obtain:
\begin{eqnarray*}
\lambda\in\sigma(f(a))
&\iff&F(a)\notin A^{-1}\\
&\iff&0\in\sigma(F(a))\\
&\iff&0\in F(\sigma(a))\\
&\iff&\exists\mu\in\sigma(a),F(\mu)=0\\
&\iff&\lambda\in f(\sigma(a))
\end{eqnarray*}

Thus, we have obtained the formula in the statement.
\end{proof}

Summarizing, we have a beginning of theory. In order to advance, we will need:

\begin{proposition}
Let $A$ be a Banach algebra.
\begin{enumerate}
\item $||a||<1\implies(1-a)^{-1}=1+a+a^2+\ldots$

\item The set $A^{-1}$ is open.

\item The map $a\to a^{-1}$ is differentiable.
\end{enumerate}
\end{proposition}

\begin{proof}
All these assertions are elementary, as follows:

\medskip

(1) This follows as in the scalar case, the computation being as follows, provided that everything converges under the norm, which amounts in saying that $||a||<1$:
\begin{eqnarray*}
(1-a)(1+a+a^2+\ldots)
&=&1-a+a-a^2+a^2-a^3+\ldots\\
&=&1
\end{eqnarray*}

(2) Assuming $a\in A^{-1}$, let us pick $b\in A$ such that we have:
$$||a-b||<\frac{1}{||a^{-1}||}$$

By using this, we have then the following norm estimate:
\begin{eqnarray*}
||1-a^{-1}b||
&=&||a^{-1}(a-b)||\\
&\leq&||a^{-1}||\cdot||a-b||\\
&<&1
\end{eqnarray*}

Thus by (1) we obtain $a^{-1}b\in A^{-1}$, and so $b\in A^{-1}$, as desired.

\medskip

(3) This follows as in the scalar case, where the derivative of $f(t)=t^{-1}$ is:
$$f'(t)=-t^{-2}$$

To be more precise, in the present Banach algebra setting the derivative is no longer a number, but rather a linear transformation. But this linear transformation can be found by developing the function $f(a)=a^{-1}$ at order 1, as follows:
\begin{eqnarray*}
(a+h)^{-1}
&=&((1+ha^{-1})a)^{-1}\\
&=&a^{-1}(1+ha^{-1})^{-1}\\
&=&a^{-1}(1-ha^{-1}+(ha^{-1})^2-\ldots)\\
&\simeq&a^{-1}(1-ha^{-1})\\
&=&a^{-1}-a^{-1}ha^{-1}
\end{eqnarray*}

We conclude that the derivative that we are looking for is:
$$f'(a)h=-a^{-1}ha^{-1}$$

Thus, we are led to the conclusion in the statement.
\end{proof}

We can now formulate a key theorem about the Banach algebras, as follows:

\index{spectrum}

\begin{theorem}
The spectrum of any Banach algebra element $\sigma(a)\subset\mathbb C$ is:
\begin{enumerate}
\item Compact.

\item Contained in the disc $D_0(||a||)$.

\item Non-empty.
\end{enumerate}
\end{theorem}

\begin{proof}
This can be proved by using the above results, as follows:

\medskip

(1) In view of (2) below, it is enough to prove that $\sigma(a)$ is closed. But this follows from the following computation, with $|\varepsilon|$ being small:
\begin{eqnarray*}
\lambda\notin\sigma(a)
&\implies&a-\lambda\in A^{-1}\\
&\implies&a-\lambda-\varepsilon\in A^{-1}\\
&\implies&\lambda+\varepsilon\notin\sigma(a)
\end{eqnarray*}

(2) This follows indeed from the following computation:
\begin{eqnarray*}
\lambda>||a||
&\implies&\Big|\Big|\frac{a}{\lambda}\Big|\Big|<1\\
&\implies&1-\frac{a}{\lambda}\in A^{-1}\\
&\implies&\lambda-a\in A^{-1}\\
&\implies&\lambda\notin\sigma(a)
\end{eqnarray*}

(3) Assume by contradiction $\sigma(a)=\emptyset$. Given a linear form $f\in A^*$, consider the following map, which is well-defined, due to our assumption $\sigma(a)=\emptyset$:
$$\varphi:\mathbb C\to\mathbb C\quad,\quad 
\lambda\to f((a-\lambda)^{-1})$$

By using Proposition 8.22 this map is differentiable, and so is a power series:
$$\varphi(\lambda)=\sum_{k=0}^\infty c_k\lambda^k$$

On the other hand, we have the following estimate, coming from definitions:
\begin{eqnarray*}
\lambda\to\infty
&\implies&a-\lambda\to\infty\\
&\implies&(a-\lambda)^{-1}\to0\\
&\implies&\varphi(\lambda)\to0
\end{eqnarray*}

Thus by the Liouville theorem from complex analysis we obtain $\varphi=0$, and since $f\in A^*$ was arbitrary, this gives $(a-\lambda)^{-1}=0$. But this is a contradiction, as desired.
\end{proof}

This was for the basic spectral theory in Banach algebras, which notably applies to the case $A=B(H)$. It is possible to go beyond the above, for instance with a holomorphic function extension of the rational functional calculus formula $\sigma(f(a))=f(\sigma(a))$ from Theorem 8.21. Also, in the case of the algebras of operators, more can be said.

\section*{8d. Operator algebras}

Let us get back now to the operator algebra $B(H)$. We know from Theorem 8.17 that this algebra has an involution $T\to T^*$, and this suggests formulating:

\index{operator algebra}

\begin{definition}
A $C^*$-algebra is a complex algebra with unit $A$, having:
\begin{enumerate}
\item A norm $a\to||a||$, making it a Banach algebra.

\item An involution $a\to a^*$, which satisfies $||aa^*||=||a||^2$, for any $a\in A$.
\end{enumerate}
\end{definition}

At the level of the basic examples, we know from Theorem 8.17 that the full operator algebra $B(H)$ is a $C^*$-algebra, in the above sense. More generally, any closed $*$-subalgebra $A\subset B(H)$ is a $C^*$-algebra. We will see later on that any $C^*$-algebra appears in fact in this way, as a closed $*$-subalgebra $A\subset B(H)$, for a certain Hilbert space $H$.

\bigskip

For the moment, we are interested in developing the theory of $C^*$-algebras, without reference to operators, or Hilbert spaces. As a first observation, we have:

\begin{proposition}
If $X$ is an abstract compact space,  the algebra $C(X)$ of continuous functions $f:X\to\mathbb C$ is a $C^*$-algebra, with structure as follows:
\begin{enumerate}
\item The norm is the usual sup norm of the functions, given by:
$$||f||=\sup_{x\in X}|f(x)|$$

\item The involution is the usual involution of the functions, given by:
$$f^*(x)=\overline{f(x)}$$
\end{enumerate}
This algebra is commutative, in the sense that $fg=gf$, for any $f,g$.
\end{proposition}

\begin{proof}
Almost everything here is trivial. Observe that we have indeed:
\begin{eqnarray*}
||ff^*||
&=&\sup_{x\in X}|f(x)\overline{f(x)}|\\
&=&\sup_{x\in X}|f(x)|^2\\
&=&||f||^2
\end{eqnarray*}

Thus, the axioms are satisfied, and finally $fg=gf$ is clear.
\end{proof}

Our claim now is that any commutative $C^*$-algebra appears as above. This is something non-trivial, which requires a number of preliminaries. We will need:

\index{spectral radius}

\begin{definition}
Given an element $a\in A$, its spectral radius 
$$\rho (a)\in(0,||a||)$$
is the radius of the smallest disk centered at $0$ containing $\sigma(a)$. 
\end{definition}

Here we have included a number of results that we already know, from Theorem 8.23, namely the fact that the spectrum is nonzero, and contained in the disk $D_0(||a||)$.

\bigskip

We have the following key result, extending our spectral theory knowledge, from the general Banach algebra setting, to the present $C^*$-algebra setting:

\index{unitary operator}
\index{self-adjoint operator}
\index{normal operator}
\index{spectral radius}

\begin{theorem}
Let $A$ be a $C^*$-algebra.
\begin{enumerate}
\item The spectrum of a unitary element $(a^*=a^{-1}$) is on the unit circle. 

\item The spectrum of a self-adjoint element ($a=a^*$) consists of real numbers. 

\item The spectral radius of a normal element ($aa^*=a^*a$) is equal to its norm.
\end{enumerate}
\end{theorem}

\begin{proof}
We use the various results established above, and notably the rational calculus formula from Theorem 8.21, and the various results from Theorem 8.23:

\medskip

(1) Assuming $a^*=a^{-1}$, we have the following norm computations:
$$||a||=\sqrt{||aa^*||}=\sqrt{1}=1$$
$$||a^{-1}||=||a^*||=||a||=1$$

Now if we denote by $D$ the unit disk, we obtain from this:
$$||a||=1\implies\sigma(a)\subset D$$
$$||a^{-1}||=1\implies\sigma(a^{-1})\subset D$$

On the other hand, by using the rational function $f(z)=z^{-1}$, we have:
$$\sigma(a^{-1})\subset D\implies \sigma(a)\subset D^{-1}$$

Now by putting everything together we obtain, as desired:
$$\sigma(a)\subset D\cap D^{-1}=\mathbb T$$

(2) This follows by using the result (1), just established above, and Theorem 8.21, with the following rational function, depending on a parameter $t\in\mathbb R$:
$$f(z)=\frac{z+it}{z-it}$$

Indeed, for $t>>0$ the element $f(a)$ is well-defined, and we have:
\begin{eqnarray*}
\left(\frac{a+it}{a-it}\right)^*
&=&\frac{(a+it)^*}{(a-it)^*}\\
&=&\frac{a-it}{a+it}\\
&=&\left(\frac{a+it}{a-it}\right)^{-1}
\end{eqnarray*}

Thus the element $f(a)$ is a unitary, and by using (1) its spectrum is contained in $\mathbb T$. We conclude from this that we have the following inclusion:
$$f(\sigma(a))=\sigma(f(a))\subset\mathbb T$$

But this shows, by applying the inverse of $f$, that we have, as desired:
$$\sigma(a)\subset f^{-1}(\mathbb T)=\mathbb R$$

(3) We already know that we have the inequality in one sense, $\rho(a)\leq ||a||$, and this for any $a\in A$. For the reverse inequality, when $a$ is normal, we fix a number as follows:
$$\rho>\rho(a)$$

We have then the following computation, with the convention that the integration over the circle $|z|=\rho$ is normalized, as for the integral of the 1 function to be 1:
\begin{eqnarray*}
\int_{|z|=\rho}\frac{z^n}{z -a}\,dz
&=&\int_{|z|=\rho}\sum_{k=0}^\infty z^{n-k-1}a^k\,dz\\
&=&\sum_{k=0}^\infty\left(\int_{|z|=\rho}z^{n-k-1}dz\right)a^k\\
&=&\sum_{k=0}^\infty\delta_{n,k+1}a^k\\
&=&a^{n-1}
\end{eqnarray*}

Here we have used the following formula, with $m\in\mathbb Z$, whose proof is elementary:
$$\int_{|z|=\rho}z^m\,dz=\delta_{m0}$$

By applying now the norm and taking $n$-th roots we obtain from the above formula, modulo some elementary manipulations, the following estimate: 
$$\rho\geq\lim_{n\to\infty}||a^n||^{1/n}$$

Now recall that $\rho$ was by definiton an arbitrary number satisfying $\rho>\rho(a)$. Thus, we have obtained the following estimate, valid for any $a\in A$:
$$\rho(a)\geq\lim_{n\to\infty}||a^n||^{1/n}$$

In order to finish, we must prove that when $a$ is normal, this estimate implies the missing estimate, namely $\rho(a)\geq||a||$. We can proceed in two steps, as follows:

\medskip

\underline{Step 1}. In the case $a=a^*$ we have $||a^n||=||a||^n$ for any exponent of the form $n=2^k$, by using the $C^*$-algebra condition $||aa^*||=||a||^2$, and by taking $n$-th roots we get:
$$\rho(a)\geq||a||$$

Thus, we are done with the self-adjoint case, with the result $\rho(a)=||a||$.

\medskip

\underline{Step 2}. In the general normal case $aa^*=a^*a$ we have $a^n(a^n)^*=(aa^*)^n$, and by using this, along with the result from Step 1, applied to $aa^*$, we obtain:
\begin{eqnarray*}
\rho(a)
&\geq&\lim_{n\to\infty}||a^n||^{1/n}\\
&=&\sqrt{\lim_{n\to\infty}||a^n(a^n)^*||^{1/n}}\\
&=&\sqrt{\lim_{n\to\infty}||(aa^*)^n||^{1/n}}\\
&=&\sqrt{\rho(aa^*)}\\
&=&\sqrt{||a||^2}\\
&=&||a||
\end{eqnarray*}

Thus, we are led to the conclusion in the statement.
\end{proof}

As a first comment, the spectral radius formula $\rho(a)=||a||$ does not hold in general, the simplest counterexample being the following non-normal matrix:
$$M=\begin{pmatrix}0&1\\0&0\end{pmatrix}$$

As another comment, we can combine the formula $\rho(a)=||a||$ for normal elements with the formula $||aa^*||=||a||^2$, and we are led to the following statement:

\index{operator algebra}

\begin{proposition}
In a $C^*$-algebra, the norm is given by
$$||a||=\sqrt{\sup\left\{\lambda\in\mathbb C\Big| aa^*-\lambda\notin A^{-1}\right\}}$$
and so is an algebraic quantity.
\end{proposition}

\begin{proof}
We have the following computation, using the condition $||aa^*||=||a||^2$, then the spectral radius formula for $aa^*$, and finally the definition of the spectral radius:
\begin{eqnarray*}
||a||
&=&\sqrt{||aa^*||}\\
&=&\sqrt{\rho(aa^*)}\\
&=&\sqrt{\sup\left\{\lambda\in\mathbb C\Big| \lambda\in\sigma(aa^*)\right\}}\\
&=&\sqrt{\sup\left\{\lambda\in\mathbb C\Big| aa^*-\lambda\notin A^{-1}\right\}}
\end{eqnarray*}

Thus, we are led to the conclusion in the statement.
\end{proof}

The above result is quite interesting, because it raises the possibility of axiomatizing the $C^*$-algebras as being the Banach $*$-algebras having the property that the formula in Proposition 8.28 defines a norm, which must satisfy the usual $C^*$-algebra conditions. However, this is something rather philosophical, and we will not follow this path.

\bigskip

Good news, we are now in position of proving a key result, namely:

\index{Gelfand theorem}

\begin{theorem}[Gelfand]
Any commutative $C^*$-algebra is the form 
$$A=C(X)$$
with the compact space $X$, called spectrum of $A$, and denoted
$$X=Spec(A)$$
appearing as the space of Banach algebra characters $\chi :A\to\mathbb C$.
\end{theorem}

\begin{proof}
This can be deduced from our spectral theory results, as follows:

\medskip

(1) Given a commutative $C^*$-algebra $A$, we can define indeed $X$ to be the set of characters $\chi :A\to\mathbb C$, with the topology making continuous all the evaluation maps:
$$ev_a:\chi\to\chi(a)$$

Then $X$ is a compact space, and $a\to ev_a$ is a morphism of algebras: 
$$ev:A\to C(X)$$

(2) We first prove that $ev$ is involutive. We use the following formula:
$$a=\frac{a+a^*}{2}-i\cdot\frac{i(a-a^*)}{2}$$

Thus it is enough to prove the following equality, for self-adjoint elements $a$:
$$ev_{a^*}=ev_a^*$$

But this is the same as proving that $a=a^*$ implies that $ev_a$ is a real function, which is in turn true, because $ev_a(\chi)=\chi(a)$ is an element of $\sigma(a)$, contained in $\mathbb R$.

\medskip

(3) Since $A$ is commutative, each element is normal, so $ev$ is isometric:
$$||ev_a||
=\rho(a)
=||a||$$

(4) It remains to prove that $ev$ is surjective. But this follows from the Stone-Weierstrass theorem, because $ev(A)$ is a closed subalgebra of $C(X)$, which separates the points.
\end{proof}

As a first consequence of the Gelfand theorem, we can extend the rational calculus formula from Theorem 8.21, to the case of the normal elements, as follows:

\index{continuous functional calculus}
\index{normal operator}

\begin{theorem}
We have the ``continuous functional calculus'' formula
$$\sigma(f(a))=f(\sigma(a))$$
valid for any normal element $a\in A$, and any continuous function $f\in C(\sigma(a))$.
\end{theorem}

\begin{proof}
Since our element $a$ is normal, the $C^*$-algebra $<a>$ that is generates is commutative, and the Gelfand theorem gives an identification as follows:
$$<a>=C(X)$$

In order to compute $X$, observe that the map $X\to\sigma(a)$ given by evaluation at $a$ is bijective. Thus, we have an identification of compact spaces, as follows:
$$X=\sigma(a)$$

As a conclusion, the Gelfand theorem provides us with an identification as follows:
$$<a>=C(\sigma(a))$$

Now given $f\in C(\sigma(a))$, we can define indeed an element $f(a)\in A$, with $f\to f(a)$ being a morphism of $C^*$-algebras, and we have $\sigma(f(a))=f(\sigma(a))$, as claimed. 
\end{proof}

The above result adds to a series of similar statements, namely Theorem 8.21, dealing with rational calculus, and the known holomorphic calculus in Banach algebras, briefly mentioned after Theorem 8.23. However, the story is not over here, because in certain special $C^*$-algebras, such as the matrix algebras $M_N(\mathbb C)$, or more generally the so-called von Neumann algebras, we can apply if we want arbitrary measurable functions to the normal elements, and we still have $\sigma(f(a))=f(\sigma(a))$. We will not get here into this.

\bigskip

As another important remark, the above result, or rather the formula $<a>=C(\sigma(a))$ from its proof, when applied to the normal operators $T\in B(H)$, is more of less the spectral theorem for such operators. Once again, we will not get here into this.

\bigskip

As a last topic, let us discuss now the GNS representation theorem, providing us with embeddings $A\subset B(H)$. We will need some more spectral theory, as follows:

\index{positive operator}
\index{square root}

\begin{proposition}
For a normal element $a\in A$, the following are equivalent:
\begin{enumerate}
\item $a$ is positive, in the sense that $\sigma(a)\subset[0,\infty)$.

\item $a=b^2$, for some $b\in A$ satisfying $b=b^*$.

\item $a=cc^*$, for some $c\in A$.
\end{enumerate}
\end{proposition}

\begin{proof}
This is something very standard, as follows:

\medskip

$(1)\implies(2)$ Since $a$ is normal, we can use Theorem 8.30, and set $b=\sqrt{a}$. 

\medskip

$(2)\implies(3)$ This is trivial, because we can set $c=b$. 

\medskip

$(3)\implies(1)$ We proceed by contradiction. By multiplying $c$ by a suitable element of $<cc^*>$, we are led to the existence of an element $d\neq0$ satisfying $-dd^*\geq0$. By writing now $d=x+iy$ with $x=x^*,y=y^*$ we have:
$$dd^*+d^*d=2(x^2+y^2)\geq0$$

Thus $d^*d\geq0$. But this contradicts the elementary fact that $\sigma(dd^*),\sigma(d^*d)$ must coincide outside $\{0\}$, that we know from Proposition 8.20.
\end{proof}

Here is now the GNS representation theorem for the $C^*$-algebras, due to Gelfand, Naimark and Segal, along with the idea of the proof:

\index{operator algebra}
\index{GNS theorem}

\begin{theorem}[GNS theorem]
Let $A$ be a $C^*$-algebra.
\begin{enumerate}
\item $A$ appears as a closed $*$-subalgebra $A\subset B(H)$, for some Hilbert space $H$. 

\item When $A$ is separable (usually the case), $H$ can be chosen to be separable.

\item When $A$ is finite dimensional, $H$ can be chosen to be finite dimensional. 
\end{enumerate}
\end{theorem}

\begin{proof}
This is something quite tricky, the idea being as follows:

\medskip

(1) Let us first discuss the commutative case, $A=C(X)$. Our claim here is that if we pick a probability measure on $X$, we have an embedding as follows:
$$C(X)\subset B(L^2(X))\quad,\quad 
f\to(g\to fg)$$

Indeed, given a function $f\in C(X)$, consider the operator $T_f(g)=fg$, acting on $H=L^2(X)$. Observe that $T_f$ is indeed well-defined, and bounded as well, because:
$$||fg||_2
=\sqrt{\int_X|f(x)|^2|g(x)|^2dx}
\leq||f||_\infty||g||_2$$

The application $f\to T_f$ being linear, involutive, continuous, and injective as well, we obtain in this way a $C^*$-algebra embedding $C(X)\subset B(H)$, as claimed.

\medskip

(2) In general, we can use a similar idea, with the positivity issues being taken care of by Proposition 8.31. Indeed, assuming that a linear form $\varphi:A\to\mathbb C$ has suitable positivity properties, making it analogous to the integration functionals $\int_X:A\to\mathbb C$ from the commutative case, we can define a scalar product on $A$, by the following formula:
$$<a,b>=\varphi(ab^*)$$

By completing we obtain a Hilbert space $H$, and we have an embedding as follows:
$$A\subset B(H)\quad,\quad 
a\to(b\to ab)$$

Thus we obtain the assertion (1), and a careful examination of the construction $A\to H$, outlined above, shows that the assertions (2,3) are in fact proved as well.
\end{proof}

There are of course many other things that can be said about bounded operators and operator algebras, but for our purposes here, the above material, and especially the Gelfand theorem, will be basically all that we will need, in what follows. For more on all this, we refer as usual to our favorite analysis authors, namely Rudin \cite{ru2} and Lax \cite{la2}. And for even more, this time in relation with physics, go with Connes \cite{con}.

\section*{8e. Exercises}

The present chapter was an introduction to linear algebra in infinite dimensions, and most of our exercises here will be about continuations of this. We first have:

\begin{exercise}
Find an explicit orthonormal basis of the Hilbert space $H=L^2[0,1]$, by applying the Gram-Schmidt procedure to the polynomials $f_n=x^n$, with $n\in\mathbb N$.
\end{exercise}

This is something both fundamental and a bit scary, and the answer can be found by doing an internet search with the keyword ``orthogonal polynomials''.

\begin{exercise}
Develop a theory of projections, isometries and symmetries inside $B(H)$, notably by examining the validity of the formula
$$\lim_{n\to\infty}(PQ)^n=P\wedge Q$$
when talking about projections, and also by taking into account the fact that
$$UU^*=1\iff U^*U=1$$
does not necessarily hold in infinite dimensions, when talking about isometries.
\end{exercise}

There are countless possible things to be done here, with all this being very useful, leading you to a much better understanding of the linear operators. Enjoy.

\begin{exercise}
Prove that for the usual matrices $A,B\in M_N(\mathbb C)$ we have
$$\sigma^+(AB)=\sigma^+(BA)$$
where $\sigma^+$ denotes the set of eigenvalues, taken with multiplicities.
\end{exercise}

As a remark, we have seen that $\sigma(AB)=\sigma(BA)$ holds outside $\{0\}$, and the equality on $\{0\}$ holds as well, because $AB$ is invertible if and only if $BA$ is invertible. However, in what regards the eigenvalues taken with multiplicities, things are more tricky.

\begin{exercise}
Clarify, with examples and counterexamples, the relation between the eigenvalues of an operator $T\in B(H)$, and its spectrum $\sigma(T)\subset\mathbb C$. 
\end{exercise}

Here, as usual, the counterexamples could only come from the shift operator $S$, on the space $H=l^2(\mathbb N)$. As a bonus exercise here, try computing the spectrum of $S$.

\begin{exercise}
Develop a theory of noncommutative geometry, by formally writing any $C^*$-algebra, not necessarily commutative, as
$$A=C(X)$$
with $X$ being a ``compact quantum space'', and report on what you found.
\end{exercise}

This is of course a very broad question, and countless things can be done here, all interesting and beautiful. We will be actually back to this, later in this book.

\part{Group theory}

\ \vskip50mm

\begin{center}
{\em Castles out of fairy tales

Timbers shivered where once there sailed

The lovesick men who caught her eye

And no one knew but Lorelei}
\end{center}

\chapter{Finite groups}

\section*{9a. Groups, examples}

We have seen so far the basics of linear algebra, with the conclusion that the theory is very useful, and quickly becomes non-trivial. We have seen as well some abstract applications, to questions in analysis and combinatorics, and with some results in the infinite dimensional case as well. All this is of course very useful in physics.

\bigskip

In this second half of this book we discuss a related topic, which is of key interest, namely the matrix groups. The theory here is once again very useful in connection with various questions in physics, the general idea being that any physical system $S$ has a group of symmetries $G(S)$, whose study can lead to concrete results about $S$.

\bigskip

Let us begin with some abstract aspects. A group is something very simple, namely a set, with a composition operation, which must satisfy what we should expect from a ``multiplication''. The precise definition of the groups is as follows: 

\index{multiplication}
\index{group}
\index{abelian group}
\index{associativity}

\begin{definition}
A group is a set $G$ with a multiplication operation 
$$(g,h)\to gh$$
which must satisfy the following conditions:
\begin{enumerate}
\item Associativity: we have $(gh)k=g(hk)$, for any $g,h,k\in G$.

\item Unit: there is an element $1\in G$ such that $g1=1g=g$, for any $g\in G$.

\item Inverses: for any $g\in G$ there is $g^{-1}\in G$ such that $gg^{-1}=g^{-1}g=1$.
\end{enumerate}
\end{definition}

The multiplication law is not necessarily commutative. In the case where it is, in the sense that $gh=hg$, for any $g,h\in G$, we call $G$ abelian, en hommage to Abel, and we usually denote its multiplication, unit and inverse operation as follows:
$$(g,h)\to g+h\quad,\quad
0\in G\quad,\quad
g\to-g$$

However, this is not a general rule, and rather the converse is true, in the sense that if a group is denoted as above, this means that the group must be abelian.

\bigskip

At the level of examples, we have for instance the symmetric group $S_N$. There are many other examples, with typically the basic systems of numbers that we know being abelian groups, and the basic sets of matrices being non-abelian groups. Once again, this is of course not a general rule. Here are some basic examples and counterexamples:

\index{groups of numbers}

\begin{proposition}
We have the following groups, and non-groups:
\begin{enumerate}
\item $(\mathbb Z,+)$ is a group.

\item $(\mathbb Q,+)$, $(\mathbb R,+)$, $(\mathbb C,+)$ are groups as well.

\item $(\mathbb N,+)$ is not a group.

\item $(\mathbb Q^*,\cdot\,)$ is a group.

\item $(\mathbb R^*,\cdot\,)$, $(\mathbb C^*,\cdot\,)$ are groups as well.

\item $(\mathbb N^*,\cdot\,)$, $(\mathbb Z^*,\cdot\,)$ are not groups.
\end{enumerate}
\end{proposition}

\begin{proof}
All this is clear from the definition of the groups, as follows:

\medskip

(1) The group axioms are indeed satisfied for $\mathbb Z$, with the sum $g+h$ being the usual sum, 0 being the usual 0, and $-g$ being the usual $-g$.

\medskip

(2) Once again, the axioms are satisfied for $\mathbb Q,\mathbb R,\mathbb C$, with the remark that for $\mathbb Q$ we are using here the fact that the sum of two rational numbers is rational, coming from:
$$\frac{a}{b}+\frac{c}{d}=\frac{ad+bc}{bd}$$

(3) In $\mathbb N$ we do not have inverses, so we do not have a group:
$$-1\notin\mathbb N$$

(4) The group axioms are indeed satisfied for $\mathbb Q^*$, with the product $gh$ being the usual product, 1 being the usual 1, and $g^{-1}$ being the usual $g^{-1}$. Observe that we must remove indeed the element $0\in\mathbb Q$, because in a group, any element must be invertible.

\medskip

(5) Once again, the axioms are satisfied for $\mathbb R^*,\mathbb C^*$, with the remark that for $\mathbb C$ we are using here the fact that the nonzero complex numbers can be inverted, coming from:
$$\frac{1}{a+ib}=\frac{a-ib}{a^2+b^2}$$

(6) Here in $\mathbb N^*,\mathbb Z^*$ we do not have inverses, so we do not have groups, as claimed.
\end{proof}

There are many interesting groups coming from linear algebra, as follows:

\index{groups of matrices}
\index{special linear group}
\index{general linear group}
\index{orthogonal group}
\index{unitary group}
\index{special orthogonal group}
\index{special unitary group}

\begin{theorem}
We have the following groups:
\begin{enumerate}
\item $(\mathbb R^N,+)$ and $(\mathbb C^N,+)$.

\item $(M_N(\mathbb R),+)$ and $(M_N(\mathbb C),+)$.

\item $(GL_N(\mathbb R),\cdot\,)$ and $(GL_N(\mathbb C),\cdot\,)$, the invertible matrices.

\item $(SL_N(\mathbb R),\cdot\,)$ and $(SL_N(\mathbb C),\cdot\,)$, with S standing for ``special'', meaning $\det=1$.

\item $(O_N,\cdot\,)$ and $(U_N,\cdot\,)$, the orthogonal and unitary matrices.

\item $(SO_N,\cdot\,)$ and $(SU_N,\cdot\,)$, with S standing as above for $\det=1$.
\end{enumerate}
\end{theorem}

\begin{proof}
All this is clear from definitions, and from our linear algebra knowledge:

\medskip

(1) The axioms are indeed clearly satisfied for $\mathbb R^N,\mathbb C^N$, with the sum being the usual sum of vectors, $-v$ being the usual $-v$, and the null vector $0$ being the unit. 

\medskip

(2) Once again, the axioms are clearly satisfied for $M_N(\mathbb R),M_N(\mathbb C)$, with the sum being the usual sum of matrices, $-M$ being the usual $-M$, and the null matrix $0$ being the unit. Observe that what we have here is in fact a particular case of (1), because any $N\times N$ matrix can be regarded as a $N^2\times1$ vector, and so at the group level we have:
$$(M_N(\mathbb R),+)\simeq(\mathbb R^{N^2},+)\quad,\quad 
(M_N(\mathbb C),+)\simeq(\mathbb C^{N^2},+)$$

(3) Regarding now $GL_N(\mathbb R),GL_N(\mathbb C)$, these are groups because the product of invertible matrices is invertible, according to the following formula:
$$(AB)^{-1}=B^{-1}A^{-1}$$

Observe that at $N=1$ we obtain the groups $(\mathbb R^*,\cdot),(\mathbb C^*,\cdot)$. At $N\geq2$ the groups $GL_N(\mathbb R),GL_N(\mathbb C)$ are not abelian, because we do not have $AB=BA$ in general.

\medskip

(4) The sets $SL_N(\mathbb R),SL_N(\mathbb C)$ formed by the real and complex matrices of determinant 1 are subgroups of the groups in (3), because of the following formula, which shows that the matrices satisfying $\det A=1$ are stable under multiplication:
$$\det(AB)=\det(A)\det(B)$$

(5) Regarding now $O_N,U_N$, here the group property is clear too from definitions, and is best seen by using the associated linear maps, because the composition of two isometries is an isometry. Equivalently, assuming $U^*=U^{-1}$ and $V^*=V^{-1}$, we have:
$$(UV)^*
=V^*U^*
=V^{-1}U^{-1}
=(UV)^{-1}$$

(6) The sets of matrices $SO_N,SU_N$ in the statement are obtained by intersecting the groups in (4) and (5), and so they are groups indeed:
$$SO_N=O_N\cap SL_N(\mathbb R)\quad,\quad 
SU_N=U_N\cap SL_N(\mathbb C)$$

Thus, all the sets in the statement are indeed groups, as claimed.
\end{proof}

Let us focus now on the finite case. The simplest finite group is the cyclic group:

\index{cyclic group}
\index{roots of unity}

\begin{definition}
The cyclic group $\mathbb Z_N$ is defined as follows:
\begin{enumerate}
\item As the additive group of remainders modulo $N$.

\item As the multiplicative group of the $N$-th roots of unity.
\end{enumerate}
\end{definition}

Observe that (1,2) are indeed equivalent, because if we set $w=e^{2\pi i/N}$, then any remainder modulo $N$ defines a $N$-th root of unity, according to the following formula:
$$k\to w^k$$

We obtain in this way all the $N$-roots of unity, so our correspondence is bijective. Moreover, our correspondence transforms the sum of remainders modulo $N$ into the multiplication of the $N$-th roots of unity, due to the following formula:
$$w^kw^l=w^{k+l}$$

Thus, the groups defined in (1,2) are isomorphic, via $k\to w^k$, and we agree to denote by $\mathbb  Z_N$ the corresponding group, and call it cyclic group. With the following comment:

\begin{comment}
Both the above conventions for $\mathbb Z_N$ are useful. The additive one
$$\mathbb Z_N=\{0,1,2,\ldots,N-1\}$$
is good for doing quick algebra, while the multiplicative one, with $\mathbb Z_N$ being
$$\xymatrix@R=20pt@C=0pt{
&w^2\ar@{.}@/^/[rr]\ar@{.}@/_/[dl]&&w\ar@{.}@/^/[dr]\\
w^3&&\ \ \ast&&1\\
&\ar@{.}@/_/[rr]\ar@{.}@/^/[ul]&&w^{N-1}\ar@{.}@/_/[ur]}$$
with $w=e^{2\pi i/N}$, is obviously ``cyclic'', and brings geometric understanding.
\end{comment}

Observe now that the cyclic groups $\mathbb Z_N$ are by definition abelian. We can construct further abelian groups by taking products of such cyclic groups, as follows:

\begin{theorem}
The following groups are all finite, and abelian,
$$G=\mathbb Z_{N_1}\times\ldots\times\mathbb Z_{N_k}$$
for any choice of the numbers $N_1,\ldots,N_k\in\mathbb N$.
\end{theorem}

\begin{proof}
This is something trivial, coming from the fact that a product of abelian groups must be abelian too. We will see later, at the end of this chapter, that any finite abelian group must appear as above, as a product of cyclic groups.
\end{proof}

Moving on, another interesting example of finite group, which is more advanced, and non-abelian this time, is the dihedral group $D_N$, which appears as follows:

\index{dihedral group}
\index{regular polygon}

\begin{definition}
The dihedral group $D_N$ is the symmetry group of 
$$\xymatrix@R=12pt@C=13pt{
&\bullet\ar@{-}[r]\ar@{-}[dl]&\bullet\ar@{-}[dr]\\
\bullet\ar@{-}[d]&&&\bullet\ar@{-}[d]\\
\bullet\ar@{-}[dr]&&&\bullet\ar@{-}[dl]\\
&\bullet\ar@{-}[r]&\bullet}$$
that is, of the regular polygon having $N$ vertices.
\end{definition}

Here are some basic examples of regular $N$-gons, at small values of the parameter $N\in\mathbb N$, and of their symmetry groups:

\medskip

\underline{$N=2$}. Here the $N$-gon is just a segment, and its symmetries are the identity $id$ and the obvious symmetry $\tau$. Thus $D_2=\{id,\tau\}$, and in group theory terms, $D_2=\mathbb Z_2$.

\medskip

\underline{$N=3$}. Here the $N$-gon is an equilateral triangle, and the symmetries are the $3!=6$ possible permutations of the vertices. Thus we have $D_3=S_3$.

\medskip

\underline{$N=4$}. Here the $N$-gon is a square, and as symmetries we have 4 rotations, of angles $0^\circ,90^\circ,180^\circ,270^\circ$, as well as 4 symmetries, with respect to the 4 symmetry axes, which are the 2 diagonals, and the 2 segments joining the midpoints of opposite sides.

\medskip

\underline{$N=5$}. Here the $N$-gon is a regular pentagon, and as symmetries we have 5 rotations, of angles $0^\circ,72^\circ,144^\circ,216^\circ,288^\circ$, as well as 5 symmetries, with respect to the 5 symmetry axes, which join the vertices to the midpoints of the opposite sides.

\medskip

\underline{$N=6$}. Here the $N$-gon is a regular hexagon, and we have 6 rotations, of angles $0^\circ,60^\circ,120^\circ,180^\circ,240^\circ,300^\circ$, and 6 symmetries, with respect to the 6 symmetry axes, which are the 3 diagonals, and the 3 segments joining the midpoints of opposite sides.

\medskip

We can see from the above that the various dihedral groups $D_N$ have many common features, and that there are some differences as well. In general, we have:

\begin{proposition}
The dihedral group $D_N$ has $2N$ elements, as follows:
\begin{enumerate}
\item We have $N$ rotations $R_1,\ldots,R_N$, with $R_k$ being the rotation of angle $2k\pi/N$. When labeling the vertices $1,\ldots,N$, the rotation formula is $R_k:i\to k+i$.

\item We have $N$ symmetries $S_1,\ldots,S_N$, with $S_k$ being the symmetry with respect to the $Ox$ axis rotated by $k\pi/N$. The symmetry formula is $S_k:i\to k-i$.
\end{enumerate}
\end{proposition}

\begin{proof}
This is clear, indeed. To be more precise, $D_N$ consists of:

\medskip

(1) The $N$ rotations, of angles $2k\pi/N$ with $k=1,\ldots,N$.

\medskip

(2) The $N$ symmetries with respect to the $N$ possible symmetry axes, which are the $N$ medians of the $N$-gon when $N$ is odd, and are the $N/2$ diagonals plus the $N/2$ lines connecting the midpoints of opposite edges, when $N$ is even.
\end{proof}

With the above description of $D_N$ in hand, we can forget if we want about geometry and the regular $N$-gon, and talk about $D_N$ abstractly, as follows:

\begin{theorem}
The dihedral group $D_N$ is the group having $2N$ elements, $R_1,\ldots,R_N$ and $S_1,\ldots,S_N$, called rotations and symmetries, which multiply as follows,
$$R_kR_l=R_{k+l}\quad,\quad 
R_kS_l=S_{k+l}$$
$$S_kR_l=S_{k-l}\quad,\quad 
S_kS_l=R_{k-l}$$
with all indices being taken modulo $N$.
\end{theorem}

\begin{proof}
With notations from Proposition 9.8, the various compositions between rotations and symmetries can be computed as follows:
$$R_kR_l\ :\ i\to l+i\to k+l+i$$
$$R_kS_l\ :\ i\to l-i\to k+l-i$$
$$S_kR_l\ :\ i\to l+i\to k-l-i$$
$$S_kS_l\ :\ i\to l-i\to k-l+i$$

But these are exactly the formulae for $R_{k+l},S_{k+l},S_{k-l},R_{k-l}$, as stated. Now since a group is uniquely determined by its multiplication rules, this gives the result.
\end{proof}

Observe that $D_N$ has the same cardinality as $E_N=\mathbb Z_N\times\mathbb Z_2$. We obviously don't have $D_N\simeq E_N$, because $D_N$ is not abelian, while $E_N$ is. So, our next goal will be that of proving that $D_N$ appears by ``twisting'' $E_N$. In order to do this, let us start with:

\begin{proposition}
The group $E_N=\mathbb Z_N\times\mathbb Z_2$ is the group having $2N$ elements, $r_1,\ldots,r_N$ and $s_1,\ldots,s_N$, which multiply according to the following rules,
$$r_kr_l=r_{k+l}\quad,\quad 
r_ks_l=s_{k+l}$$
$$s_kr_l=s_{k+l}\quad,\quad
s_ks_l=r_{k+l}$$
with all the indices being taken modulo $N$.
\end{proposition}

\begin{proof}
With the notation $\mathbb Z_2=\{1,\tau\}$, the elements of the product group $E_N=\mathbb Z_N\times\mathbb Z_2$ can be labeled $r_1,\ldots,r_N$ and $s_1,\ldots,s_N$, as follows:
$$r_k=(k,1)\quad,\quad
s_k=(k,\tau)$$

These elements multiply then according to the formulae in the statement. Now since a group is uniquely determined by its multiplication rules, this gives the result.
\end{proof}

Let us compare now Theorem 9.9 and Proposition 9.10. In order to formally obtain $D_N$ from $E_N$, we must twist some of the multiplication rules of $E_N$, namely:
$$s_kr_l=s_{k+l}\to s_{k-l}\quad,\quad 
s_ks_l=r_{k+l}\to r_{k-l}$$

Informally, this amounts in following the rule ``$\tau$ switches the sign of what comes afterwards", and we are led in this way to the following definition:

\index{crossed product}

\begin{definition}
Given groups $H,K$, with an action $K\curvearrowright H$, the crossed product
$$G=H\rtimes K$$
is the set $H\times K$, with multiplication $(g,s)(h,t)=(gh^s,st)$.
\end{definition}

It is routine to check that $G$ is indeed a group. Observe that when the action is trivial, $h^s=h$ for any $h\in H$ and $s\in K$, we obtain the usual product $H\times K$.  

\bigskip

Now with this technology in hand, by getting back to the dihedral group $D_N$, we can improve Theorem 9.9, into a final result on the subject, as follows:

\index{dihedral group}
\index{crossed product decomposition}

\begin{theorem}
We have a crossed product decomposition as follows,
$$D_N=\mathbb Z_N\rtimes\mathbb Z_2$$
with $\mathbb Z_2=\{1,\tau\}$ acting on $\mathbb Z_N$ via switching signs, $k^\tau=-k$.
\end{theorem}

\begin{proof}
We have an action $\mathbb Z_2\curvearrowright\mathbb Z_N$ given by the formula in the statement, namely $k^\tau=-k$, so we can consider the corresponding crossed product group:
$$L_N=\mathbb Z_N\rtimes\mathbb Z_2$$

In order to understand the structure of $L_N$, we follow Proposition 9.10. The elements of $L_N$ can indeed be labeled $\rho_1,\ldots,\rho_N$ and $\sigma_1,\ldots,\sigma_N$, as follows:
$$\rho_k=(k,1)\quad,\quad 
\sigma_k=(k,\tau)$$

Now when computing the products of such elements, we basically obtain the formulae in Proposition 9.10, perturbed as in Definition 9.11. To be more precise, we have:
$$\rho_k\rho_l=\rho_{k+l}\quad,\quad 
\rho_k\sigma_l=\sigma_{k+l}$$
$$\sigma_k\rho_l=\sigma_{k+l}\quad,\quad 
\sigma_k\sigma_l=\rho_{k+l}$$

But these are exactly the multiplication formulae for $D_N$, from Theorem 9.9. Thus, we have an isomorphism $D_N\simeq L_N$ given by $R_k\to\rho_k$ and $S_k\to\sigma_k$, as desired.
\end{proof}

As a third basic example of a finite group, we have the symmetric group $S_N$. This is a group that we already met, when talking about the determinant, and we have:

\index{symmetric group}
\index{permutation group}
\index{signature}

\begin{theorem}
The permutations of $\{1,\ldots,N\}$ form a group, denoted $S_N$, and called symmetric group. This group has $N!$ elements. The signature map
$$\varepsilon:S_N\to\mathbb Z_2$$
can be regarded as being a group morphism, with values in $\mathbb Z_2=\{\pm1\}$, and 
$$A_N=\left\{\sigma\in S_N\Big|\varepsilon(\sigma)=1\right\}$$
is a subgroup having $N!/2$ elements, called alternating group.
\end{theorem}

\begin{proof}
As explained in chapter 2, the group property is clear, and the count is clear as well. As for the last assertion, recall the following formula, also from chapter 2:
$$\varepsilon(\sigma\tau)=\varepsilon(\sigma)\varepsilon(\tau)$$

But this tells us precisely that $\varepsilon$ is a group morphism, and we can see as well from this that $A_N\subset S_N$ is indeed a subgroup. Finally, with $\tau\in S_N$ being any transposition we have $S_N=A_N\sqcup\tau A_N$, and it follows that we have $|A_N|=N!/2$, as claimed.
\end{proof}

We will be back to $S_N$ on many occasions, in what follows. At an even more advanced level now, we have the hyperoctahedral group $H_N$, which appears as follows:

\index{hyperoctahedral group}
\index{hypercube}

\begin{definition}
The hyperoctahedral group $H_N\subset O_N$ is the group formed by the symmetries of the unit cube in $\mathbb R^N$,
$$\xymatrix@R=18pt@C=20pt{
&\bullet\ar@{-}[rr]&&\bullet\\
\bullet\ar@{-}[rr]\ar@{-}[ur]&&\bullet\ar@{-}[ur]\\
&\bullet\ar@{-}[rr]\ar@{-}[uu]&&\bullet\ar@{-}[uu]\\
\bullet\ar@{-}[uu]\ar@{-}[ur]\ar@{-}[rr]&&\bullet\ar@{-}[uu]\ar@{-}[ur]
}$$
viewed as a graph, or equivalently, as a metric space.
\end{definition}

Here the equivalence at the end is clear from definitions, because any symmetry of the cube graph must preserve the lengths of the edges, and so we have:
$$G(\square_{graph})=G(\square_{metric})$$

The hyperoctahedral group is a quite interesting group, whose definition, as a symmetry group, reminds that of the dihedral group $D_N$. So, let us start our study in the same way as we did for $D_N$, with a discussion at small values of $N\in\mathbb N$:

\medskip

\underline{$N=1$}. Here the 1-cube is the segment, whose symmetries are the identity $id$ and the flip $\tau$. Thus, we obtain the group with 2 elements, which is a very familiar object: 
$$H_1=D_2=S_2=\mathbb Z_2$$

\underline{$N=2$}. Here the 2-cube is the square, and so the corresponding symmetry group is the dihedral group $D_4$, which is a group that we know well:
$$H_2=D_4=\mathbb Z_4\rtimes\mathbb Z_2$$

\underline{$N=3$}. Here the 3-cube is the usual cube, and the situation is considerably more complicated, because this usual cube has no less than 48 symmetries.

\medskip

All this looks quite complicated, but fortunately we can count $H_N$, as follows:

\begin{theorem}
We have the cardinality formula
$$|H_N|=2^NN!$$
coming from the fact that $H_N$ is the symmetry group of the coordinate axes of $\mathbb R^N$.
\end{theorem}

\begin{proof}
This follows from some geometric thinking, as follows:

\medskip

(1) Consider the standard cube in $\mathbb R^N$, centered at 0, and having as vertices the points having coordinates $\pm1$. With this picture in hand, it is clear that the symmetries of the cube coincide with the symmetries of the $N$ coordinate axes of $\mathbb R^N$.

\medskip

(2) In order to count now these latter symmetries, a bit as we did for the dihedral group, observe first that we have $N!$ permutations of these $N$ coordinate axes. 

\medskip

(3) But each of these permutations of the coordinate axes $\sigma\in S_N$ can be further ``decorated'' by a sign vector $e\in\{\pm1\}^N$, consisting of the possible $\pm1$ flips which can be applied to each coordinate axis, at the arrival. Thus, we have:
$$|H_N|
=|S_N|\cdot|\mathbb Z_2^N|
=N!\cdot2^N$$

Thus, we are led to the conclusions in the statement. 
\end{proof}

As in the dihedral group case, it is possible to go beyond this, as follows:

\index{crossed product}
\index{wreath product}
\index{hyperoctahedral group}

\begin{theorem}
We have a wreath product decomposition $H_N=\mathbb Z_2\wr S_N$, which means by definition that we have a crossed product decomposition
$$H_N=\mathbb Z_2^N\rtimes S_N$$
with the permutations $\sigma\in S_N$ acting on the elements $e\in\mathbb Z_2^N$ as follows:
$$\sigma(e_1,\ldots,e_N)=(e_{\sigma(1)},\ldots,e_{\sigma(N)})$$
In particular we have, as found before, the cardinality formula $|H_N|=2^NN!$.
\end{theorem}

\begin{proof}
As explained in the proof of Theorem 9.15, the elements of $H_N$ can be identified with the pairs $g=(e,\sigma)$ consisting of a permutation $\sigma\in S_N$, and a sign vector $e\in\mathbb Z_2^N$, so that at the level of the cardinalities, we have:
$$|H_N|=|\mathbb Z_2^N\times S_N|$$

To be more precise, given an element $g\in H_N$, the element $\sigma\in S_N$ is the corresponding permutation of the $N$ coordinate axes, regarded as unoriented lines in $\mathbb R^N$, and $e\in\mathbb Z_2^N$ is the vector collecting the possible flips of these coordinate axes, at the arrival. Now observe that the product formula for two such pairs $g=(e,\sigma)$ is as follows, with the permutations $\sigma\in S_N$ acting on the elements $f\in\mathbb Z_2^N$ as in the statement:
$$(e,\sigma)(f,\tau)=(ef^\sigma,\sigma\tau)$$

Thus, we are precisely in the framework of Definition 9.11, and we conclude that we have a crossed product decomposition, as follows:
$$H_N=\mathbb Z_2^N\rtimes S_N$$
 
Thus, we are led to the conclusion in the statement, with the formula $H_N=\mathbb Z_2\wr S_N$ being just a shorthand for the decomposition $H_N=\mathbb Z_2^N\rtimes S_N$ that we found.
\end{proof}

\section*{9b. Cayley theorem}

At the level of the general theory now, we have the following fundamental result regarding the finite groups, due to Cayley:

\index{Cayley embedding}

\begin{theorem}
Given a finite group $G$, we have an embedding as follows,
$$G\subset S_N\quad,\quad g\to(h\to gh)$$
with $N=|G|$. Thus, any finite group is a permutation group. 
\end{theorem}

\begin{proof}
Given a group element $g\in G$, we can associate to it the following map:
$$\sigma_g:G\to G\quad,\quad 
h\to gh$$

Since $gh=gh'$ implies $h=h'$, this map is bijective, and so is a permutation of $G$, viewed as a set. Thus, with $N=|G|$, we can view this map as a usual permutation, $\sigma_G\in S_N$. Summarizing, we have constructed so far a map as follows:
$$G\to S_N\quad,\quad 
g\to\sigma_g$$

Our first claim is that this is a group morphism. Indeed, this follows from:
$$\sigma_g\sigma_h(k)
=\sigma_g(hk)
=ghk
=\sigma_{gh}(k)$$

It remains to prove that this group morphism is injective. But this follows from:
\begin{eqnarray*}
g\neq h
&\implies&\sigma_g(1)\neq\sigma_h(1)\\
&\implies&\sigma_g\neq\sigma_h
\end{eqnarray*}

Thus, we are led to the conclusion in the statement.
\end{proof}

Observe that in the above statement the embedding $G\subset S_N$ that we constructed depends on a particular writing $G=\{g_1,\ldots,g_N\}$, which is needed in order to identify the permutations of $G$ with the elements of the symmetric group $S_N$. This is not very good, in practice, and as an illustration, for the basic examples of groups that we know, the Cayley theorem provides us with embeddings as follows:
$$\mathbb Z_N\subset S_N\quad,\quad 
D_N\subset S_{2N}\quad,\quad 
S_N\subset S_{N!}\quad,\quad 
H_N\subset S_{2^NN!}$$

And here the first embedding is the good one, the second one is not the best possible one, but can be useful, and the third and fourth embeddings are useless. Thus, as a conclusion, the Cayley theorem remains something quite theoretical. We will be back to this later on, with a systematic study of the ``representation'' problem.

\bigskip

Getting back now to our main series of finite groups, $\mathbb Z_N\subset D_N\subset S_N\subset H_N$, these are of course permutation groups, according to the above. However, and perhaps even more interestingly, these are as well subgroups of the orthogonal group $O_N$:
$$\mathbb Z_N\subset D_N\subset S_N\subset H_N\subset O_N$$

Indeed, we have $H_N\subset O_N$, because any transformation of the unit cube in $\mathbb R^N$ must extend into an isometry of the whole $\mathbb R^N$, in the obvious way. Now in view of this, it makes sense to look at the finite subgroups $G\subset O_N$. With two remarks, namely: 

\bigskip

(1) Although we do not have examples yet, following our general ``complex is better than real'' philosophy, it is better to look at the general subgroups $G\subset U_N$.

\bigskip

(2) Also, it is better to upgrade our study to the case where $G$ is compact, and this in order to cover some interesting continuous groups, such as $O_N,U_N,SO_N,SU_N$. 

\bigskip

Long story short, we are led in this way to the study of the closed subgroups $G\subset U_N$. Let us start our discussion here with the following simple fact:

\index{closed subgroup}

\begin{proposition}
The closed subgroups $G\subset U_N$ are precisely the closed sets of matrices $G\subset U_N$ satisfying the following conditions:
\begin{enumerate}
\item $U,V\in G\implies UV\in G$.

\item $1\in G$.

\item $U\in G\implies U^{-1}\in G$.
\end{enumerate}
\end{proposition}

\begin{proof}
This is clear from definitions, the only point with this statement being the fact that a subset $G\subset U_N$ can be a group or not, as indicated above.
\end{proof}

As a second result now regarding the closed subgroups $G\subset U_N$, let us prove that any finite group $G$ appears in this way. This is something more or less clear from what we have, but let us make this precise. We first have the following key result:

\index{permutation matrix}

\begin{theorem}
We have a group embedding as follows, obtained by regarding $S_N$ as the permutation group of the $N$ coordinate axes of $\mathbb R^N$,
$$S_N\subset O_N$$
which makes $\sigma\in S_N$ correspond to the matrix having $1$ on row $\sigma(j)$ and column $j$, for any $j$, and having $0$ entries elsewhere.
\end{theorem}

\begin{proof}
This is something quite fundamental, the idea being as follows:

\medskip

(1) To start with, we can certainly regard $S_N$ as being the permutation group of the $N$ coordinate axes of $\mathbb R^N$. Now since these permutations of the $N$ coordinate axes of $\mathbb R^N$ are isometries, this provides us with a group embedding $S_N\subset O_N$, as stated.

\medskip

(2) Regarding now the formula of this embedding, we have by definition:
$$\sigma(e_j)=e_{\sigma(j)}$$

Thus, the permutation matrix corresponding to $\sigma$ is given by:
$$\sigma_{ij}=
\begin{cases}
1&{\rm if}\ \sigma(j)=i\\
0&{\rm otherwise}
\end{cases}$$

We are theferore led to the conclusion in the statement.
\end{proof}

We can combine the above result with the Cayley theorem, and we obtain the following result, which is something very nice, having theoretical importance:

\index{finite group}
\index{Cayley embedding}
\index{permutation group}

\begin{theorem}
Given a finite group $G$, we have an embedding as follows,
$$G\subset O_N\quad,\quad g\to(e_h\to e_{gh})$$
with $N=|G|$. Thus, any finite group is an orthogonal matrix group.
\end{theorem}

\begin{proof}
The Cayley theorem gives an embedding as follows:
$$G\subset S_N\quad,\quad g\to(h\to gh)$$

On the other hand, Theorem 9.19 provides us with an embedding as follows:
$$S_N\subset O_N\quad,\quad 
\sigma\to(e_i\to e_{\sigma(i)})$$

Thus, we are led to the conclusion in the statement.
\end{proof}

The same remarks as for the Cayley theorem apply. First, the embedding $G\subset O_N$ that we constructed depends on a particular writing $G=\{g_1,\ldots,g_N\}$. And also, for the basic examples of groups that we know, the embeddings that we obtain are as follows:
$$\mathbb Z_N\subset O_N\quad,\quad 
D_N\subset O_{2N}\quad,\quad 
S_N\subset O_{N!}\quad,\quad 
H_N\subset O_{2^NN!}$$

As before, here the first embedding is the good one, the second one is not the best possible one, but can be useful, and the third and fourth embeddings are useless. 

\bigskip

Summarizing, in order to advance, it is better to forget about the Cayley theorem, and build on Theorem 9.19 instead. In relation with the basic groups, we have:

\begin{theorem}
We have the following finite groups of matrices:
\begin{enumerate}
\item $\mathbb Z_N\subset O_N$, the cyclic permutation matrices.

\item $D_N\subset O_N$, the dihedral permutation matrices.

\item $S_N\subset O_N$, the permutation matrices.

\item $H_N\subset O_N$, the signed permutation matrices.
\end{enumerate}
\end{theorem}

\begin{proof}
This is something self-explanatory, the idea being that Theorem 9.19 provides us with embeddings as follows, given by the permutation matrices:
$$\mathbb Z_N\subset D_N\subset S_N\subset O_N$$

In addition, looking back at the definition of $H_N$, this group inserts into the embedding on the right, $S_N\subset H_N\subset O_N$. Thus, we are led to the conclusion that all our 4 groups appear as groups of suitable ``permutation type matrices''. To be more precise:

\medskip

(1) The cyclic permutation matrices are by definition the matrices as follows, with 0 entries elsewhere, and form a group, which is isomorphic to the cyclic group $\mathbb Z_N$:
$$U=\begin{pmatrix}
&&&1\\
&&&&\ddots\\
&&&&&1\\
1\\
&\ddots\\
&&1
\end{pmatrix}$$

(2) The dihedral matrices are the above cyclic permutation matrices, plus some suitable symmetry permutation matrices, and form a group which is isomorphic to $D_N$.

\medskip

(3) The permutation matrices, which by Theorem 9.19 form a group which is isomorphic to $S_N$, are the $0-1$ matrices having exactly one 1 on each row and column.

\medskip

(4) Finally, regarding the signed permutation matrices, these are by definition the $(-1)-0-1$ matrices having exactly one nonzero entry on each row and column, and by Theorem 9.15 these matrices form a group, which is isomorphic to $H_N$.
\end{proof}

Finally, let us mention that when looking, more generally, at the finite subgroups of the unitary groups, we have many interesting examples too. More on these later. 

\section*{9c. General theory}

Let us go back now to the abstract groups, as defined in the beginning of this chapter, and develop some theory, without relation to linear algebra. We first have:

\begin{theorem}
Given a finite group $G$ and a subgroup $H\subset G$, the sets
$$G/H=\{gH\Big|g\in G\}\quad,\quad H\backslash G=\{Hg\Big|g\in G\}$$
both consist of partitions of $G$ into subsets of size $H$, and we have the formula
$$|G|=|H|\cdot|G/H|=|H|\cdot|H\backslash G|$$
which shows that the order of the subgroup divides the order of the group:
$$|H|\ \big|\ |G|$$
When $H\subset G$ is normal, $gH=Hg$ for any $g\in G$, the space $G/H=H\backslash G$ is a group.
\end{theorem}

\begin{proof}
There are several assertions here, which are in fact all trivial, when deduced in the precise order indicated in the statement. To be more precise, the partition claim for $G/H$ can be deduced as follows, and the proof for $H\backslash G$ is similar:
$$gH\cap kH\neq\emptyset\iff g^{-1}k\in H\iff gH=kH$$

With this in hand, the cardinality formulae are all clear, and it remains to prove the last assertion. But here, the point is that when $H\subset G$ is normal, we have:
$$gH=kH,sH=tH\implies gsH=gtH=gHt=kHt=ktH$$

Thus $G/H=H\backslash G$ is a indeed group, with multiplication $(gH)(sH)=gsH$.
\end{proof}

As a main consequence of the above result, which is equally useful, we have:

\begin{theorem}
Given a finite group $G$, any $g\in G$ generates a cyclic subgroup
$$<g>=\{1,g,g^2,\ldots,g^{k-1}\}$$
with $k=ord(g)$ being the smallest number $k\in\mathbb N$ satisfying $g^k=1$. Also, we have
$$ord(g)\ \big|\ |G|$$
that is, the order of any group element divides the order of the group.
\end{theorem}

\begin{proof}
As before with Theorem 9.22, we have opted here for a long collection of statements, which are all trivial, when deduced in the above precise order. To be more precise, consider the semigroup $<g>\subset G$ formed by the sequence of powers of $g$:
$$<g>=\{1,g,g^2,g^3,\ldots\}\subset G$$

Since $G$ was assumed to be finite, the sequence of powers must cycle, $g^n=g^m$ for some $n<m$, and so we have $g^k=1$, with $k=m-n$. Thus, we have in fact:
$$<g>=\{1,g,g^2,\ldots,g^{k-1}\}$$

Moreover, we can choose $k\in\mathbb N$ to be minimal with this property, and with this choice, we have a set without repetitions. Thus $<g>\subset G$ is indeed a group, and more specifically a cyclic group, of order $k=ord(g)$. Finally, $ord(g)\,|\,|G|$ follows from Theorem 9.22.
\end{proof}

More concretely now, groups are meant to act on sets, and we have here:

\begin{proposition}
Given an action $G\curvearrowright X$ and a point $x\in X$, we have
$$|G(x)|=|G|/|G_x|$$
where $G_x=\{g\in G|g(x)=x\}$. In particular, the cardinality of orbits divides $|G|$.
\end{proposition}

\begin{proof}
In order to prove this, we will construct a bijection, as follows:
$$\varphi:G/G_x\to G(x)$$

But the formula of $\varphi$ can only be something straightforward, as follows:
$$\varphi(gG_x)=g(x)$$

So, let us see if this works. To start with, $\varphi$ is well-defined and injective, due to:
\begin{eqnarray*}
gG_x=hG_x
&\iff&g^{-1}h\in G_x\\
&\iff&g^{-1}h(x)=x\\
&\iff&g(x)=h(x)
\end{eqnarray*}

But $\varphi$ is clearly surjective too, and we therefore obtain the result.
\end{proof}

As an application of the above technology, we have the following key result:

\index{Cauchy theorem}

\begin{theorem}[Cauchy]
Given a finite group $G$, and a prime number satisfying
$$p\ |\ |G|$$
$G$ has an element of order $p$. Equivalently, $G$ has a subgroup of order $p$.
\end{theorem}

\begin{proof}
We must find $g\neq1$ with $g^p=1$. In order to do so, let us set:
$$X=\left\{(g_1,\ldots,g_p)\in G^p\Big|g_1\ldots g_p=1\right\}$$

We have then an obvious action $\mathbb Z_p\curvearrowright X$, by rotation, as follows:
$$k(g_1,\ldots,g_p)=(g_{k+1},\ldots,g_{k+p})$$

Now let us decompose $X$ into orbits. This gives the following formula, with $F\subset X$ being the fixed points, and with the sum being over the non-trivial orbits $O$:
$$|X|=|F|+\sum_{|O|\geq2}|O|$$

Next, let us look at this equality modulo $p$. To start with, we have:
$$|X|=|G|^{p-1}=0(p)$$

Also, in what regards the fixed points, we can say here that we have:
$$(1,\ldots,1)\in F\implies |F|\geq1$$

Finally, by Proposition 9.24 the size of any orbit must divide $|\mathbb Z_p|=p$, and so:
$$|O|\geq2\implies |O|=p$$

Now by putting everything together, modulo our $p\geq2$, we conclude that:
$$|F|\geq2$$

But this is exactly what we need, because the fixed points are precisely the elements $(g,\ldots,g)\in G^p$ with $g^p=1$. Thus, we have found $g\neq1$ with $g^p=1$, as desired.
\end{proof}

Moving on, this time with some inspiration from linear algebra, let us call unitary representation of $G$ any group morphism $u:G\to U_N$. This is a key notion, and of particular interest is the case $N=1$, where we have the following result:

\index{representation}
\index{character}
\index{dual group}
\index{group of characters}
\index{cyclic group}
\index{product of cyclic groups}
\index{self-dual group}

\begin{theorem}
Given a finite group $G$, the group morphisms $\chi:G\to\mathbb T$, called characters of $G$, form a finite abelian group $\widehat{G}$, called Pontrjagin dual of $G$. We have:
\begin{enumerate}
\item The dual of a cyclic group is the group itself, $\widehat{\mathbb Z}_N=\mathbb Z_N$.

\item The dual of a product is the product of duals, $\widehat{G\times H}=\widehat{G}\times\widehat{H}$.

\item Any product of cyclic groups $G=\mathbb Z_{N_1}\times\ldots\times\mathbb Z_{N_k}$ is self-dual, $G=\widehat{G}$.
\end{enumerate}
\end{theorem}

\begin{proof}
We have several assertions here, the idea being as follows:

\medskip

(1) Our first claim is that $\widehat{G}$ is a group, with the pointwise multiplication, namely:
$$(\chi\rho)(g)=\chi(g)\rho(g)$$

Indeed, if $\chi,\rho$ are characters, so is $\chi\rho$, and so the multiplication is well-defined on $\widehat{G}$. Regarding the unit, this is the trivial character, constructed as follows:
$$1:G\to\mathbb T\quad,\quad 
g\to1$$ 

Finally, we have inverses, with the inverse of $\chi:G\to\mathbb T$ being its conjugate:
$$\bar{\chi}:G\to\mathbb T\quad,\quad 
g\to\overline{\chi(g)}$$

Next, our claim is that the group $\widehat{G}$ is finite. Indeed, assuming that we have a character $\chi:G\to\mathbb T$, we have the following formula, for any group element $g\in G$:
$$g^k=1\implies\chi(g)^k=1$$

Thus $\chi(g)$ must be one of the $k$-th roots of unity, and in particular there are finitely many choices for $\chi(g)$. Finally, the fact that $\widehat{G}$ is abelian follows from definitions.

\medskip

(2) Next, in the cyclic group case, a character $\chi:\mathbb Z_N\to\mathbb T$ is uniquely determined by its value $z=\chi(g)$ on the standard generator $g\in\mathbb Z_N$. But this value must satisfy:
$$z^N=1$$

We conclude that we must have $z\in\mathbb Z_N$. Conversely, any $N$-th root of unity $z\in\mathbb Z_N$ defines a certain character $\chi:\mathbb Z_N\to\mathbb T$, by setting, for any $r\in\mathbb N$:
$$\chi(g^r)=z^r$$

Summarizing, we have indeed an identification $\widehat{\mathbb Z}_N=\mathbb Z_N$, as claimed.

\medskip

(3) Regarding now products of groups, a character $\chi:G\times H\to\mathbb T$ must satisfy:
$$\chi(g,h)=\chi\left[(g,1)(1,h)\right]=\chi(g,1)\chi(1,h)$$

Thus $\chi$ must appear as the product of its restrictions $\chi_{|G},\chi_{|H}$, which must be both characters, and this gives $\chi\in\widehat{G}\times\widehat{H}$, as desired. Finally, the last assertion is clear.
\end{proof}

As a continuation, we can get some further insight into duality by using the spectral theory methods developed in chapter 8, and we have the following result:

\index{finite abelian group}
\index{Fourier transform}
\index{dual group}

\begin{theorem}
Given a finite abelian group $G$, we have an isomorphism of commutative $C^*$-algebras as follows, obtained by linearizing/delinearizing the characters:
$$\mathbb C[G]\simeq C(\widehat{G})$$
Also, the Pontrjagin duality is indeed a duality, in the sense that we have $G=\widehat{\widehat{G}}$.
\end{theorem}

\begin{proof}
We have several assertions here, the idea being as follows:

\medskip

(1) Given a finite abelian group $G$, consider indeed the group algebra $\mathbb C[G]$, having as elements the formal combinations of elements of $G$, and with involution given by:
$$g^*=g^{-1}$$

This $*$-algebra is then a $C^*$-algebra, with norm coming by making act $\mathbb C[G]$ on itself, so by the Gelfand theorem we obtain an isomorphism as follows:
$$\mathbb C[G]=C(X)$$

To be more precise, $X$ is the space of the $*$-algebra characters as follows:
$$\chi:\mathbb C[G]\to\mathbb C$$

The point now is that by delinearizing, such a $*$-algebra character must come from a usual group character of $G$, obtained by restricting to $G$, as follows:
$$\chi:G\to\mathbb T$$

Thus we have $X=\widehat{G}$, and we are led to the isomorphism in the statement, namely:
$$\mathbb C[G]\simeq C(\widehat{G})$$

(2) In order to prove now the second assertion, consider the following group morphism, which is available for any finite group $G$, not necessarily abelian:
$$G\to\widehat{\widehat{G}}\quad,\quad 
g\to(\chi\to\chi(g))$$

Our claim is that in the case where $G$ is abelian, this is an isomorphism. As a first observation, we only need to prove that this morphism is injective or surjective, because the cardinalities match, according to the following formula, coming from (1):
$$|G|=\dim\mathbb C[G]=\dim C(\widehat{G})=|\widehat{G}|$$

(3) We will prove that the above morphism is injective. For this purpose, let us compute its kernel. We know that $g\in G$ is in the kernel when the following happens:
$$\chi(g)=1\quad,\quad\forall\chi\in\widehat{G}$$

But this means precisely that $g\in\mathbb C[G]$ is mapped, via the isomorphism $\mathbb C[G]\simeq C(\widehat{G})$ constructed in (1), to the constant function $1\in C(\widehat{G})$, and now by getting back to $\mathbb C[G]$ via our isomorphism, this shows that we have indeed $g=1$, which ends the proof.
\end{proof}

\section*{9d. Abelian groups}

Let us go back now to the finite abelian groups, with the aim of proving that these are exactly the products of cyclic groups. Let us start with a basic result, as follows:

\begin{proposition}
Given a finite abelian group $G$, and $p|\,|G|$, the set
$$G_p=\left\{g\in G\Big|\exists k\in\mathbb N,g^{p^k}=1\right\}$$
is a subgroup, having as order the biggest power of $p$ dividing $|G|$.
\end{proposition}

\begin{proof}
This is something elementary, the idea being as follows:

\medskip

(1) To start with, the fact that the set in the statement $G_p\subset G$ is a subgroup is clear, coming from the following computation, valid inside any abelian group:
$$g^a=1,h^b=1\implies (gh)^{ab}=g^ah^b=1$$

Indeed, given two elements $g,h\in G$, having as orders powers of $p$, this computation shows that  $gh\in G$ has as order a certain power of $p$ too, as desired.

\medskip

(2) Next, assuming $|G|=p^kn$ with $(n,p)=1$, we must show that we have $|G_p|=p^k$. But this is best seen by contradiction. Indeed, assuming $p|\,|G/G_p|$, by Cauchy we would have a certain non-trivial element $hG_p\in G/G_p$ of order $p$. But this means $h\notin G_p$, $h^p\in G_p$, which in turn reads $h\notin G_p$, $h\in G_p$, which is contradictory.
\end{proof}

As a continuation of this, we have the following key result:

\begin{theorem}
Given a finite abelian group $G$, we have
$$G=\prod_pG_p$$
with $G_p\subset G$ with $p$ prime being the subgroups constructed above.
\end{theorem}

\begin{proof}
By using the fact that our group $G$ is abelian, we have a group morphism as follows, with the order of the factors when computing $\prod_pg_p$ being irrelevant:
$$\prod_pG_p\to G\quad,\quad (g_p)\to\prod_pg_p$$

(1) Our first claim is that this morphism is injective. Indeed, let us consider an element in its kernel, which amounts in having an equation of the following type:
$$g_1\ldots g_k=1$$

Now since the elements $g_1$ and $g_2\ldots g_k$, which are inverse to each other, must have the same order, and the order of $g_1$ is a certain prime power, and that of $g_2\ldots g_k$ is not divisible by that prime, we conclude that the kernel is trivial, as claimed. 

\medskip

(2) It remains to prove that our morphism is surjective. But this can be done in the pedestrian way, by picking $g\in G$, writing its order as $ord(g)=p_1^{a_1}\ldots p_k^{a_k}$, and doing some arithmetic in order to reach to a writing of type $g=g_1\ldots g_k$, with $g_i\in G_{p_i}$.
\end{proof}

Getting now to what we wanted to do, structure theorem for the abelian groups, Theorem 9.29 does half of the job. For the other half, we must decompose the components $G_p$. With the convention that $p$-group means $|G|=p^k$, for some $k\in\mathbb N$, we have:

\index{abelian p-group}

\begin{theorem}
The abelian $p$-groups decompose as follows:
$$G=\mathbb Z_{p^{r_1}}\times\ldots\times\mathbb Z_{p^{r_s}}$$
That is, the abelian $p$-groups are the products of cyclic $p$-groups.
\end{theorem}

\begin{proof}
We can do this by recurrence on $|G|$, as follows:

\medskip

(1) Let us pick $g\in G$ of maximal order, say $ord(g)=p^k$, and consider the subgroup $H=<g>$ that it generates, inside $G$. By recurrence, the quotient group $G/H$ must decompose as follows, with the components $C_i$ being cyclic groups:
$$G/H=C_1\times\ldots\times C_n$$

Our goal will be that of producing, out of this, an isomorphism as follows:
$$G=H\times C_1\times\ldots\times C_n$$

(2) Let us start by fixing some notation. The subgroups $C_i\subset G/H$ appearing above being cyclic, we can denote them as $C_i=\{z_i^aH\}$, with $z_iH\in C_i$ being some chosen generators for them. And with this, the isomorphism that we have is:
$$\varphi:C_1\times\ldots\times C_n\to G/H\quad,\quad 
(z_1^{a_1}H,\ldots,z_n^{a_n}H)\to z_1^{a_1}\ldots z_n^{a_n}H$$

Our more precise claim now, which will prove the result, is that, with a suitable choice of the generators $z_iH\in C_i$, we can lift this into an isomorphism as follows:
$$\psi:H\times C_1\times\ldots\times C_n\to G\quad,\quad 
(g^a,z_1^{a_1}H,\ldots,z_n^{a_n}H)\to g^az_1^{a_1}\ldots z_n^{a_n}$$

(3) In order to do this, let us look at one of the components, $C=C_i$. If we pick an arbitrary generator $zH\in C$, with $z\in G$, the following happens, trivially:
$$ord(zH)|ord(z)$$

And our claim now, which will provide us with what is needed in (2), is that we can always arrange for our generator $zH\in C$, with $z\in G$, as to have equality:
$$ord(zH)=ord(z)$$

(4) Summarizing, we have eventually found something concrete to prove, in relation with what we want to do, so let us prove this. Let us start with an arbitrary generator $xH\in C$, with $x\in G$. Consider the two orders mentioned in (3), namely:
$$p^r=ord(xH)\quad,\quad p^s=ord(x)\quad,\quad r\leq s$$

Our goal will be that of suitably modifying our generator $xH$, as to have $r=s$.

\medskip

(5) In order to do so, let us look at the following group element $y\in G$:
$$y=x^{p^r}\quad,\quad ord(y)=p^{s-r}$$

Since $ord(xH)=p^r$ we have $ord(yH)=1$, which means $y\in H$. Now since $H=<g>$ was the group generated by $g$, we can write $y$ as follows, with $(n,p)=1$:
$$y=g^{np^t}$$

Now recall that $g\in G$ was chosen of maximal order $p^k$. Thus, we have:
$$ord(y)=p^{k-t}$$

We conclude that we have $s-r=k-t$. Now consider the following element:
$$z=xg^{-np^{t-r}}$$

Our claim is that this is the element $z\in G$ that we were looking for, in (3).

\medskip

(6) Indeed, we first have the following computation, which gives $ord(z)\leq p^r$:
$$z^{p^r}=x^{p^r}g^{-np^t}=y\cdot y^{-1}=1$$

Also, $zH=xH=C$, and so $ord(zH)=|C|=p^r$. Thus we have, as desired:
$$ord(zH)=ord(z)=p^r$$

(7) Time for the endgame. Let us go back to the isomorphism in (2), which was as follows, and with the generators $z_iH\in C_i$ with $z_i\in G$ being chosen as above:
$$\varphi:C_1\times\ldots\times C_n\to G/H\quad,\quad 
(z_1^{a_1}H,\ldots,z_n^{a_n}H)\to z_1^{a_1}\ldots z_n^{a_n}H$$

Our claim is that this lifts into an isomorphism as follows:
$$\psi:H\times C_1\times\ldots\times C_n\to G\quad,\quad 
(g^a,z_1^{a_1}H,\ldots,z_n^{a_n}H)\to g^az_1^{a_1}\ldots z_n^{a_n}$$

(8) Indeed, this latter map is well-defined, due to $ord(z_iH)=ord(z_i)$. It is also clear that $\psi$ is a group morphism. Also, since $\varphi$ is surjective, so must be $\psi$. Finally, since the cardinalities of the domain and range match, $\psi$ must be an isomorphism, as desired.
\end{proof}

Time now to put everything together. We obtain the following remarkable result:

\index{finite abelian group}

\begin{theorem}
The finite abelian groups are the products of cyclic groups:
$$G=\mathbb Z_{N_1}\times\ldots\times\mathbb Z_{N_k}$$
Moreover, we can choose the numbers $N_i$ to be prime powers.
\end{theorem}

\begin{proof}
This follows indeed by putting together all the above, and more specifically, by combining Theorem 9.29 and Theorem 9.30. As further remarks on this:

\medskip

(1) In fact, what you need to know is just the first part of the present theorem, because the second part is easy to recover, thanks to the following elementary isomorphisms:
$$\mathbb Z_{p_1^{a_1}\ldots p_k^{a_k}}=\mathbb Z_{p_1^{a_1}}\times\ldots\times\mathbb Z_{p_k^{a_k}}$$

(2) There is a uniqueness assertion too, which is elementary, stating that with $G$ fully split, with $N_i$ prime powers, the components will be unique, up to permutation.
\end{proof}

As an application of the above, and in relation with characters, let us go back to the generalized Fourier matrices, from chapter 7. We have here the following result:

\index{Fourier transform}
\index{Fourier matrix}
\index{generalized Fourier matrix}

\begin{theorem}
Given a finite abelian group $G$, with dual group $\widehat{G}=\{\chi:G\to\mathbb T\}$, consider the corresponding Fourier coupling, namely:
$$\mathcal F_G:G\times\widehat{G}\to\mathbb T\quad,\quad 
(i,\chi)\to\chi(i)$$
\begin{enumerate}
\item Via the standard isomorphism $G\simeq\widehat{G}$, this Fourier coupling can be regarded as a square matrix, $F_G\in M_G(\mathbb T)$, which is a complex Hadamard matrix.

\item In the case of the cyclic group $G=\mathbb Z_N$ we obtain in this way, via the standard identification $\mathbb Z_N=\{1,\ldots,N\}$, the Fourier matrix $F_N$.

\item In general, when using a decomposition $G=\mathbb Z_{N_1}\times\ldots\times\mathbb Z_{N_k}$, the corresponding Fourier matrix is given by $F_G=F_{N_1}\otimes\ldots\otimes F_{N_k}$.
\end{enumerate}
\end{theorem}

\begin{proof}
This follows indeed by using the above finite abelian group theory:

\medskip

(1) With the identification $G\simeq\widehat{G}$ made our matrix is given by $(F_G)_{i\chi}=\chi(i)$, and the scalar products between the rows are computed as follows:
$$<R_i,R_j>
=\sum_\chi\chi(i)\overline{\chi(j)}
=\sum_\chi\chi(i-j)
=|G|\cdot\delta_{ij}$$

Thus, we obtain indeed a complex Hadamard matrix.

\medskip

(2) This follows from the well-known and elementary fact that, via the identifications $\mathbb Z_N=\widehat{\mathbb Z_N}=\{1,\ldots,N\}$, the Fourier coupling here is as follows, with $w=e^{2\pi i/N}$:
$$(i,j)\to w^{ij}$$

(3) We use here the following formula that we know, for the duals of products: 
$$\widehat{H\times K}=\widehat{H}\times\widehat{K}$$

At the level of the corresponding Fourier couplings, we obtain from this:
$$F_{H\times K}=F_H\otimes F_K$$

Now by decomposing $G$ into cyclic groups, as in the statement, and by using (2) for the cyclic components, we obtain the formula in the statement.
\end{proof}

As a nice application of the above result, we have:

\index{Walsh matrix}

\begin{theorem}
The Walsh matrix, $W_N$ with $N=2^n$, which is given by
$$W_N=\begin{pmatrix}1&1\\1&-1\end{pmatrix}^{\otimes n}$$
is the Fourier matrix of the finite abelian group $K_N=\mathbb Z_2^n$.
\end{theorem}

\begin{proof}
We know that the first Walsh matrix is a Fourier matrix:
$$W_2=F_2=F_{K_2}$$

Now by taking tensor powers we obtain from this that we have, for any $N=2^n$:
$$W_N
=W_2^{\otimes n}
=F_{K_2}^{\otimes n}
=F_{K_2^n}
=F_{K_N}$$

Thus, we are led to the conclusion in the statement.
\end{proof}

Summarizing, we have now a better understanding of the generalized Fourier matrices, and of the complex Hadamard matrices in general, and also a new and fresh point of view on the various discrete Fourier analysis considerations from chapter 7.

\bigskip

All this is quite interesting, suggesting among others that we should have a deeper relation between group theory and Fourier analysis. In answer, this is indeed the case, with the ultimate result here stating that associated to any locally compact abelian group $G$ is a Fourier transform, which can be useful for many purposes. Good to know.

\section*{9e. Exercises}

There are many things that can be said about groups, especially in the matrix case, $G\subset U_N$, and we will discuss this later in this book. Our exercises here will rather focus on the abstract groups, as in the end of the present chapter, and we first have:

\begin{exercise}
Given a locally compact abelian group $G$, prove that its group characters, which must be by definition continuous,
$$\chi:G\to\mathbb T$$
form a locally compact abelian group, denoted $\widehat{G}$, and called dual of $G$.
\end{exercise}

Here locally compact means that any group element $g\in G$ has a neighborhood which is compact, a bit in analogy with what happens for the real numbers $r\in\mathbb R$.

\begin{exercise}
Prove that the integers are dual to the unit circle, and vice versa:
$$\widehat{\mathbb Z}=\mathbb T\quad,\quad 
\widehat{\mathbb T}=\mathbb Z$$
Also, prove that the group of real numbers is self-dual, $\widehat{\mathbb R}=\mathbb R$.
\end{exercise}

To be more precise, we already know from the above that we have $\widehat{\mathbb Z}_N=\mathbb Z_N$, for any $N\in\mathbb N$, and the first question, regarding $\mathbb Z$ and $\mathbb T$, is a kind of ``$N=\infty$'' version of this. As for the second question, regarding $\mathbb R$, this is related to all this as well. 

\begin{exercise}
Prove that the finitely generated abelian groups are
$$G=\mathbb Z_{N_1}\times\ldots\times\mathbb Z_{N_k}$$
with the convention $\mathbb Z_\infty=\mathbb Z$, and that the compact matrix abelian groups are
$$H=\mathbb Z_{N_1}\times\ldots\times\mathbb Z_{N_k}$$
with this time the convention $\mathbb Z_\infty=\mathbb T$. Also, prove that $G=\widehat{H}$ and $H=\widehat{G}$.
\end{exercise}

This exercise, generalizing everything that we know, or almost, is actually something quite tricky, requiring a good knowledge of both algebra and analysis.

\begin{exercise}
Clarify the relation between the dualities
$$\widehat{\mathbb Z}_N=\mathbb Z_N\quad,\quad 
\widehat{\mathbb Z}=\mathbb T\quad,\quad 
\widehat{\mathbb T}=\mathbb Z\quad,\quad 
\widehat{\mathbb R}=\mathbb R$$
and the various types of Fourier transforms available.
\end{exercise}

To be more precise here, the problem is that of understanding why the above 3 dualities correspond to the main 3 types of known Fourier transforms, namely the discrete Fourier transforms, the usual Fourier series, and the usual Fourier transforms. And with the remark that this is something that we already know, for the first duality.

\chapter{Rotation groups}

\section*{10a. Rotation groups}

We have seen that there are many interesting examples of finite groups $G$, which usually appear as groups of orthogonal matrices $G\subset O_N$, or unitary matrices $G\subset U_N$. In this chapter we have a closer look at the subgroups $G\subset U_N$. We have:

\begin{question}
What are the subgroups of the $4$ main rotation groups,
$$\xymatrix@R=50pt@C=50pt{
SU_N\ar[r]&U_N\\
SO_N\ar[u]\ar[r]&O_N\ar[u]}$$
in low dimensions, $N=2,3,\ldots$? What about generic dimensions $N\in\mathbb N$?
\end{question}

Let us start with the following result, regarding the 4 main rotation groups themselves, which is something very useful, that we will use many times, in what follows:

\begin{proposition}
The following happen, regarding the main rotation groups:
\begin{enumerate}
\item $U\in O_N\implies\det U=\pm1$.

\item $O_N=SO_N\sqcup(-SO_N)$, when $N$ is odd.

\item $U\in U_N\implies|\det U|=1$.

\item $U_N=\bigcup_{w\in\mathbb T}wSU_N$, for any $N$.
\end{enumerate}
\end{proposition}

\begin{proof}
This is something elementary, coming from definitions, as follows:

\medskip

(1) This comes indeed from the following computation:
\begin{eqnarray*}
U\in O_N
&\implies&U^t=U^{-1}\\
&\implies&\det(U^t)=\det(U^{-1})\\
&\implies&\det U=(\det U)^{-1}\\
&\implies&\det U=\pm1
\end{eqnarray*}

(2) According to (1) we have the following decomposition formula, with $\overline{SO}_N\subset O_N$ standing for the set of orthogonal matrices having determinant $-1$:
$$O_N=SO_N\sqcup\overline{SO}_N$$

Now the point is that when $N$ is odd we have $\det(-U)=-\det U$, for any matrix $U\in M_N(\mathbb R)$, and by using this, we can see right away that we have:
$$\overline{SO}_N=-SO_N$$

Thus, we are led to the decomposition formula in the statement, namely:
$$O_N=SO_N\sqcup(-SO_N)$$

By the way, observe that this fails when $N$ is even, and in a quite drastic way, for instance because at $N=2$ the group $SO_2$ consists of the rotations of the plane, while the other component $\overline{SO}_2$ consists of the symmetries of the plane. More on this later.

\medskip

(3) This follows from the following computation, similar to the one in (1):
\begin{eqnarray*}
U\in U_N
&\implies&U^*=U^{-1}\\
&\implies&\det(U^*)=\det(U^{-1})\\
&\implies&\overline{\det U}=(\det U)^{-1}\\
&\implies&|\det U|=1
\end{eqnarray*}

(4) According to (3) we have the following decomposition formula, with $SU_N^{(z)}\subset U_N$ standing for the set of unitary matrices having determinant $z\in\mathbb T$, and coming with the warning that, contrary to the decomposition in (2), this is not a decomposition into connected components, due to the continuous nature of the parameter $z\in\mathbb T$:
$$U_N=\bigsqcup_{z\in\mathbb T}SU_N^{(z)}$$

Still following (2), let us try now to relate the components $SU_N^{(z)}$ to the main component, $SU_N=SU_N^{(1)}$. But this is an easy task in the present complex case, because we can extract $N$-th roots of any complex number. Indeed, let $w\in\mathbb T$ be such that:
$$w^N=z$$

Now given an arbitrary matrix $U\in SU_N^{(z)}$, the rescaled matrix $V=U/w$ is unitary, $V\in U_N$. As for the determinant of this latter matrix, this is given by:
\begin{eqnarray*}
\det(V)
&=&\det(U/w)\\
&=&\det U/w^N\\
&=&z/z\\
&=&1
\end{eqnarray*}

Thus we have $V\in SU_N$, and so $U\in wSU_N$, and with this in hand, our previous decomposition of $U_N$ takes the following form, which is the one in the statement:
$$U_N=\bigcup_{w\in\mathbb T}wSU_N$$

(5) Finally, observe that this latter decomposition is no longer a disjoint union, due to the choice needed in the above, when solving $w^N=z$. As yet another remark, getting back now to (2), all this suggests some complex number trickery, based on $i^2=-1$, in order to deal with $O_N$ when $N$ is even. We will leave some exploration here as an interesting exercise, and with the remark however that the $N=2$ case, discussed in (2), shows that we cannot really expect very concrete things to arise, in this way.
\end{proof}

With this discussed, time for some classification work, at small values of $N$. To start with, at $N=1$ all our matrices are just numbers, and the main rotation groups are:
$$\xymatrix@R=17pt@C=18pt{
SU_1\ar[rr]&&U_1&&\{1\}\ar[rr]&&\mathbb T\\
&&&=\\
SO_1\ar[uu]\ar[rr]&&O_1\ar[uu]&&\{1\}\ar[uu]\ar[rr]&&\{\pm1\}\ar[uu]}$$

Equivalently, with $\mathbb Z_s$ standing as usual for the group of $s$-th roots of unity, and with the extra convention $\mathbb Z_\infty=\mathbb T$, that we already used in chapter 9, the diagram is:
$$\xymatrix@R=17pt@C=18pt{
\mathbb Z_1\ar[rr]&&\mathbb Z_\infty\\
\\
\mathbb Z_1\ar[uu]\ar[rr]&&\mathbb Z_2\ar[uu]}$$

Now the point is that, with the finite subgroups of the cyclic groups being cyclic, we are led to the following result, answering Question 10.1 at $N=1$:

\begin{theorem}
The finite subgroups of the basic continuous groups at $N=1$ are:
$$\xymatrix@R=17pt@C=18pt{
SU_1\ar[rr]&&U_1&&\mathbb Z_1\ar[rr]&&\{\mathbb Z_n|n\in\mathbb N\}\\
&&&:\\
SO_1\ar[uu]\ar[rr]&&O_1\ar[uu]&&\mathbb Z_1\ar[uu]\ar[rr]&&\{\mathbb Z_1,\mathbb Z_2\}\ar[uu]}$$
That is, all the finite rotation groups at $N=1$ are cyclic.
\end{theorem}

\begin{proof}
This is certainly something trivial, with only some explanations regarding the subgroups of $U_1=\mathbb T$ being needed, with the situation here being as follows:

\medskip

(1) To start with, the unit circle $\mathbb T$ has many subgroups, as you can see by picking some random numbers $\{z_i\}\subset\mathbb T$, finitely many, or countably many, or even uncountably many, and looking at the group $G=<z_i>$ that they generate, which can vary a lot. 

\medskip

(2) However, when looking at the finite subgroups $G\subset\mathbb T$, things are easy, due to:
\begin{eqnarray*}
|G|=m
&\implies&g^m=1, \forall g\in G\\
&\implies&g\in\mathbb Z_m,\ \forall g\in G\\
&\implies&G\subset\mathbb Z_m\\
&\implies&G=\mathbb Z_n,\ n|m
\end{eqnarray*}

Thus, end of the story, and we are led to the conclusion in the statement. 

\medskip

(3) Finally, let us mention that in what regards the infinite subgroups $G\subset\mathbb T$, when restricting the attention to those which are closed, we only have one solution, namely $G=\mathbb T$ itself. Thus, as a generalization of the present result, we can say that all closed rotation groups at $N=1$, finite or not, are cyclic, with our usual convention $\mathbb Z_\infty=\mathbb T$.
\end{proof}
 
At $N=2$ now, let us first study $SO_2$, $O_2$ are their subgroups. In what regards the groups $SO_2$, $O_2$ themselves, these are groups that we know well, and this since chapter 1, but always good to talk about them again. Their basic theory is as follows:

\begin{theorem}
We have the following results:
\begin{enumerate}
\item $SO_2$ is the group of usual rotations in the plane, which are given by:
$$R_t=\begin{pmatrix}\cos t&-\sin t\\ \sin t&\cos t\end{pmatrix}$$

\item $O_2$ consists in addition of the usual symmetries in the plane, given by:
$$S_t=\begin{pmatrix}\cos t&\sin t\\ \sin t&-\cos t\end{pmatrix}$$

\item Abstractly speaking, we have group isomorphisms as follows,
$$SO_2\simeq\mathbb T\quad,\quad O_2=\mathbb T\rtimes\mathbb Z_2$$
with the second one coming from $O_2=SO_2\rtimes<S_t>$, for any symmetry $S_t$.
\end{enumerate}
\end{theorem}

\begin{proof}
These are basically things that we know, as follows:

\medskip

(1) This is clear, because the only isometries of the plane which preserve the orientation are the usual rotations. As for the formula of $R_t$, rotation of angle $t$, this is something that we know well from chapter 1, obtained by computing $R_t\binom{1}{0}$ and $R_t\binom{0}{1}$.

\medskip

(2) This is clear too, because rotations left aside, we are left with the symmetries of the plane, in the usual sense. As for formula of $S_t$, symmetry with respect to $Ox$ rotated by $t/2$, this is something that we know too, obtained by computing $S_t\binom{1}{0}$ and $S_t\binom{0}{1}$.

\medskip

(3) The first assertion is clear, because the angles $t\in\mathbb R$, taken as usual modulo $2\pi$, form the group $\mathbb T$. As for the second assertion, the proof here is similar to the proof of the crossed product decomposition $D_n=\mathbb Z_n\rtimes\mathbb Z_2$ for the dihedral groups.
\end{proof}

Getting now to the subgroups of $SO_2,O_2$, we have the following result:

\begin{theorem}
The finite subgroups of $SO_2,O_2$ are as follows:
\begin{enumerate}
\item The finite subgroups of $SO_2$ are the cyclic groups $\mathbb Z_n$.

\item For $O_2$, we obtain in addition the dihedral groups $D_n$.
\end{enumerate}
\end{theorem}

\begin{proof}
This is again something elementary, as follows:

\medskip

(1) This is indeed something clear, geometrically, which formally comes from $SO_2\simeq\mathbb T$, via the discussion from Theorem 10.4, regarding the same group there, $U_1\simeq\mathbb T$. 

\medskip

(2) In order to prove this, consider a finite subgroup as follows:
$$G\subset O_2\quad,\quad G\not\subset SO_2$$

According to (1), we have a formula as follows, for a certain $n\in\mathbb N$:
$$G\cap SO_2=\mathbb Z_n$$

Now let us pick $S\in G-SO_2$. Since products of symmetries are rotations, any other element $T\in G-SO_2$ must satisfy $ST\in\mathbb Z_n$, and so $T\in S\mathbb Z_n$. We conclude that our group $G$ must appear as follows, coming from a subgroup $\mathbb Z_n\subset\mathbb T$, and a symmetry $S\in O_2$:
$$G=\mathbb Z_n\sqcup S\mathbb Z_n$$

But this latter group must have the same multiplication table as the dihedral group $D_n$, and  conclude that we have an isomorphism $G\simeq D_n$, as desired.
\end{proof}

Quite nice the above, and in fact we can do better, as follows:

\begin{theorem}
The finite rotation groups in $2$ dimensions appear as the symmetry groups of the regular polygons,
$$\xymatrix@R=12pt@C=13pt{
&\bullet\ar@{-}[r]\ar@{-}[dl]&\bullet\ar@{-}[dr]\\
\bullet\ar@{-}[d]&&&\bullet\ar@{-}[d]\\
\bullet\ar@{-}[dr]&&&\bullet\ar@{-}[dl]\\
&\bullet\ar@{-}[r]&\bullet}$$
with these polygons being taken unoriented as above, or oriented.
\end{theorem}

\begin{proof}
This is indeed self-explanatory, based on Theorem 10.5, and with the remark that in what regards $SO_2$, looking at the symmetries of an oriented polygon, or at the orientation-preserving symmetries of an unoriented polygon, is the same thing.
\end{proof}

The above result looks quite exciting, and it is tempting at this point to forget our next task, namely understanding what happens in 2 complex dimensions, and move instead to 3 real dimensions, with the following interesting question in mind:

\begin{question}
Can we have a $3D$ analogue of Theorem 10.6 going, with regular polygons replaced by regular polyhedra?
\end{question}

And good question this is. We will see in the next section, following Plato, and then Euler, and Klein and others, that the answer to this question is remarkably ``yes'', and with this solving our group theory problems, in 3 real dimensions. 

\bigskip

As for the 2 complex dimensions, these will be not forgotten either, and we will see later, following again Euler, Klein and others, including this time Rodrigues, Hamilton, and also Pauli, Dirac and other physicists, that things are quite interesting here too.

\section*{10b. Klein subgroups}

With Question 10.7 in mind, let us get now into 3D geometry, and symmetries. At the start of everything, we have the following remarkable result, going back to Plato:

\begin{theorem}
There are $5$ regular polyhedra, called Platonic solids, namely:
\begin{enumerate}
\item Tetrahedron, having $4$ vertices and $4$ faces.

\item Octahedron, having $6$ vertices and $8$ faces.

\item Cube, having $8$ vertices and $6$ faces.

\item Icosahedron, having $12$ vertices and $20$ faces.

\item Dodecahedron, having $20$ vertices and $12$ faces.
\end{enumerate}
\end{theorem} 

\begin{proof}
Many things can be said here, the idea being as follows:

\medskip

(1) Let us try to figure out how a regular polyhedron looks like. There are a number of faces meeting at each vertex, $\geq3$ faces to be more precise, and when flattening the polyhedron there, we can see appear an angle $t$, called angle defect at that vertex:
$$\xymatrix@R=15pt@C=12pt{
&\ar@{-}[ddr]&&\\
&{\rm face}&&t\ar@.[dl]\\
\ar@{-}[rr]&&\ast\ar@{-}[rr]&&\\
&{\rm face}&{\rm face}&{\rm face}\\
&\ar@{-}[uur]&&\ar@{-}[uul]
}$$

Now since hexagons and higher have angles $\geq120^\circ$, these cannot be used for constructing polyhedra, due to $t>0$. In fact, still due to $t>0$, we are left with 5 cases:

\medskip

-- Polyhedron made of triangles, with 3 or 4 or 5 faces meeting at each vertex.

\medskip

-- Polyhedron made of squares, with 3 faces meeting at each vertex.

\medskip

-- Polyhedron made of pengatons, with 3 faces meeting at each vertex.

\medskip

(2) Now let us try to construct the solutions. In the first case, polyhedron made of triangles, with 3 faces meeting at each vertex, we obtain the tetrahedron:
$$\xymatrix@R=15pt@C=20pt{
&\bullet\ar@{-}[ddl]\ar@{-}[ddrr]\ar@{-}[dddr]\\
\\
\bullet\ar@{-}[drr]\ar@{.}[rrr]&&&\bullet\ar@{-}[dl]\\
&&\bullet}$$

(3) Two other obvious solutions, corresponding to the second and fourth cases above, triangles meeting $\times4$, and squares meeting $\times3$, are the octahedron and the cube:
$$\xymatrix@R=5pt@C=2pt{
&&&\bullet\ar@{.}[ddl]\ar@{-}[dddlll]\ar@{-}[ddrrr]\ar@{-}[dddr]\\
\\
&&\bullet\ar@{.}[rrrr]&&&&\bullet\\
\bullet\ar@{-}[rrrr]\ar@{.}[urr]&&&&\bullet\ar@{-}[urr]\\
\\
&&&\bullet\ar@{-}[uur]\ar@{-}[uuurrr]\ar@{-}[uulll]\ar@{.}[uuul]}
\qquad\qquad\qquad
\xymatrix@R=13pt@C=13pt{
&\bullet\ar@{-}[rr]&&\bullet\\
\bullet\ar@{-}[rr]\ar@{-}[ur]&&\bullet\ar@{-}[ur]\\
&\bullet\ar@{.}[rr]\ar@{.}[uu]&&\bullet\ar@{-}[uu]\\
\bullet\ar@{-}[uu]\ar@{.}[ur]\ar@{-}[rr]&&\bullet\ar@{-}[uu]\ar@{-}[ur]
}$$

Before going further, observe that there is a relation between these two polyhedra, with the vertices of the octahedron appearing at the middle of the faces of the cube, and vice versa. Due to this, we say that the octahedron and cube are dual, and with this explaining why their number of vertices and faces are interchanged, as follows:
$$(6,8)\leftrightarrow(8,6)$$

By the way, observe that the tetrahedron is self-dual, $(4,4)\leftrightarrow(4,4)$. These dualities will be quite important to us later, when looking at the symmetries of our polyhedra.

\medskip

(4) Back to constructing solutions, we are left with studying the third and fifth cases in (1), namely triangles meeting $\times5$, and pentagons meeting $\times3$. And here, by some kind of miracle, we have indeed solutions, namely the icosahedron and dodecahedron, which look as follows, with in each case half of the faces, those facing us, represented:
$$\xymatrix@R=2.6pt@C=10pt{
&&\bullet\ar@{-}[dddl]\ar@{-}[dddll]\ar@{-}[dddr]\ar@{-}[dddrr]\\
\\
\\
\bullet\ar@{-}[r]\ar@{-}[dd]&\bullet\ar@{-}[rr]\ar@{-}[ddl]\ar@{-}[ddr]&&\bullet\ar@{-}[r]\ar@{-}[ddl]\ar@{-}[ddr]&\bullet\ar@{-}[dd]\\
\\
\bullet\ar@{-}[rr]\ar@{-}[dddrr]&&\bullet\ar@{-}[rr]\ar@{-}[ddd]&&\bullet\ar@{-}[dddll]\\
\\
\\
&&\bullet}\qquad\qquad\qquad
\xymatrix@R=1pt@C=7pt{
&&&\bullet\ar@{-}[dll]\ar@{-}[drr]\\
&\bullet\ar@{-}[dr]\ar@{-}[ddl]&&&&\bullet\ar@{-}[dl]\ar@{-}[ddr]\\
&&\bullet\ar@{-}[rr]\ar@{-}[dd]&&\bullet\ar@{-}[dd]\\
\bullet\ar@{-}[dddr]&&&&&&\bullet\ar@{-}[dddl]\\
&&\bullet\ar@{-}[dr]\ar@{-}[ddl]&&\bullet\ar@{-}[dl]\ar@{-}[ddr]\\
&&&\bullet\ar@{-}[dd]\\
&\bullet\ar@{-}[dr]&&&&\bullet\ar@{-}[dl]\\
&&\bullet\ar@{-}[r]&\bullet&\bullet\ar@{-}[l]}$$

As before with the octahedron and cube, these two latter polyhedra are dual, with this intechanging their number of vertices and faces, $(12,20)\leftrightarrow(20,12)$.
\end{proof}

Getting back now to Question 10.7, we would like to compute the symmetry groups $G\subset O_3$ and $SG\subset SO_3$ of the various Platonic solids that we found, and then try to prove that these are basically all the finite subgroups of $O_3$ and $SO_3$.

\bigskip

In order to do so, let us begin with some generalities regarding $O_3,SO_3$ and their subgroups. We have here the following elementary result, further building on what we know from Proposition 10.2, regarding the groups $O_N,SO_N$ with $N$ odd:

\begin{proposition}
The following happen, regarding $O_3,SO_3$ and their subgroups:
\begin{enumerate}
\item The central symmetry $-1\in O_3$ is not orientation-preserving, $-1\notin SO_3$.

\item We have a disjoint union decomposition $O_3=SO_3\sqcup(-SO_3)$.

\item This decomposition gives an identification $O_3=SO_3\times\mathbb Z_2$.

\item More generally, assuming $G\subset O_3$, $-1\in G$, we have $G=SG\times\mathbb Z_2$.
\end{enumerate}
\end{proposition}

\begin{proof}
This is something elementary, as follows:

\medskip

(1) This is best viewed by using the determinant, $\det(-1)=-1$.

\medskip

(2) This follows indeed from $\det U=\pm1$ for $U\in O_3$, and from (1).

\medskip

(3) This is the group-theoretical reformulation of the decomposition in (2).

\medskip

(4) This is similar, based on $G=SG\sqcup(-SG)$, coming from $-1\in G$.
\end{proof}

Getting now to the symmetry groups that we are interested in, those of the Platonic solids found in the previous section, we have the following result, about them:

\begin{theorem}
The symmetry groups $G\subset O_3$ and the orientation-preserving symmetry groups $SG\subset SO_3$ of the Platonic solids are as follows:
\begin{enumerate}
\item Tetrahedron: $G=S_4$, $SG=A_4$.

\item Octahedron and cube: $G=S_4\times\mathbb Z_2$, $SG=S_4$.

\item Icosahedron and dodecahedron: $G=A_5\times\mathbb Z_2$, $SG=A_5$.
\end{enumerate}
\end{theorem} 

\begin{proof}
This basically comes from our experience from chapter 9, with some extra work needed for the  icosahedron and dodecahedron, the idea being as follows:

\medskip

(1) In what regards the tetrahedron, we certainly have $G=S_4$, and then $SG=A_4$, and with this latter group being usually called tetrahedral group. Observe that, contrary to what happens for the other polyhedra, the central symmetry $-1\in O_3$ is not a symmetry of the tetrahedron, so Proposition 10.9 (4) does not apply. In fact, we have $S_4\neq A_4\times\mathbb Z_2$ as abstract groups, because none of the transpositions $\tau\in S_4$ is central.

\medskip

(2) Regarding now the cube, here we have $G=H_3$, and $SG=SH_3=S_4$, as we know well since chapter 9, and with this latter $S_4$ being best understood as acting on the diagonals of the cube. Then, due to $-1\in G$, Proposition 10.9 (4) applies, and gives:
$$G=S_4\times\mathbb Z_2$$

(3) As for the ocathedron, this being dual to the cube, the symmetry groups are the same. Let us mention also that $SG=S_4$ is called octahedral group, and with this explaining why $G=H_3$, which is twice as big, is called hyperoctahedral group.

\medskip

(4) In what regards now the icosahedron and dodecahedron, these are dual too, so they have the same symmetry groups. In order to compute these common symmetry groups, let us look at the dodecahedron, whose picture, facing us, was as follows:
$$\xymatrix@R=1pt@C=7pt{
&&&\bullet\ar@{-}[dll]\ar@{-}[drr]\\
&\bullet\ar@{-}[dr]\ar@{-}[ddl]&&&&\bullet\ar@{-}[dl]\ar@{-}[ddr]\\
&&\bullet\ar@{-}[rr]\ar@{-}[dd]&&\bullet\ar@{-}[dd]\\
\bullet\ar@{-}[dddr]&&&&&&\bullet\ar@{-}[dddl]\\
&&\bullet\ar@{-}[dr]\ar@{-}[ddl]&&\bullet\ar@{-}[dl]\ar@{-}[ddr]\\
&&&\circ\ar@{-}[dd]\\
&\bullet\ar@{-}[dr]&&&&\bullet\ar@{-}[dl]\\
&&\bullet\ar@{-}[r]&\bullet&\bullet\ar@{-}[l]}$$

Now let us pick a vertex, say the one marked $\circ$, and look at the 3 faces meeting at this vertex. A symmetry $g\in SG$ must then send this vertex $\circ$ to one of the 20 available vertices $*$ of the dodecahedron, and then there is an extra $\times3$ choice, coming from the permutation of the 3 faces, at the arrival, around $*$. Thus, we conclude that we have:
$$|SG|=20\times3=60$$

(5) Before further commenting on the dodecahedron, it is worth noticing that our method above applies to any regular polyhedron $P$. Indeed, if we denote by $v$ the number of vertices, and by $m$ the number of faces meeting at any vertex, we obtain:
$$|SG|=vm=2e$$

In addition, $|SG|=2e$ can be seen as well directly, because any symmetry $g\in SG$ is uniquely determined by its action on a given edge, up to a $\times2$ choice at the arrival. Needless to say, all this fits with the data for our various polyhedra, as follows:
$$e=6,12,30\implies|SG|=12,24,60$$

(6) Getting back now to the dodecahedron, as a conclusion to the above discussion, we have two ways at looking at the corresponding group $SG$, coming from:
$$|SG|=20\times3=30\times2$$

Observe also that we have an embedding $\mathbb Z_5\subset SG$, obtained by rotating any given face of the dodecahedron. Now by putting everything together, this shows, via some routine abstract algebra that we will leave as an exercise, that we have, as claimed:
$$SG=A_5$$

(7) But you might wonder if there is a simpler proof for this, using a clever embedding $SG\subset S_5$, say a bit as before with $SH_3=S_4$ acting on the diagonals of the cube. In answer, yes, but with this being a bit neuron-burning, the idea being that we have exactly 5 cubes having vertices among the 20 vertices of the dodecahedron, and the symmetries $g\in SG$ come from permutations of these 5 cubes, which must be alternating.

\medskip

(8) Finally, still talking dodecahedron and isocahedron, these have central symmetry $-1\in G$, so by Proposition 10.9 (4) we obtain $G=A_5\times\mathbb Z_2$, as claimed.
\end{proof}

Good work that we did, and time now to answer Question 10.7, regarding the classification of finite groups of 3D rotations. In order to deal with this, we will need:

\index{rotation axis}

\begin{theorem}[Euler]
Any usual rotation in $3D$ space
$$U\in SO_3$$
has a rotation axis. 
\end{theorem}

\begin{proof}
We have the following computation, using some linear algebra magic:
\begin{eqnarray*}
\det(U-1)
&=&\det(U^t-1)\\
&=&\det(U^t(1-U))\\
&=&\det(U^t)\det(1-U)\\
&=&\det(1-U)
\end{eqnarray*}

Thus $\det(U-1)=0$, which tells us that $U$ must have a $1$-eigenvector:
$$U\xi=\xi$$

Thus, we got our rotation axis for our abstract rotation $U\in SO_3$, as desired. 
\end{proof}

We can now answer Question 10.7 positively, as follows:

\begin{theorem}[Klein]
The finite subgroups of $SO_3$ are as follows,
\begin{enumerate}
\item Cyclic, $\mathbb Z_n$.

\item Dihedral, $D_n$.

\item Tetrahedral, $A_4$.

\item Octahedral, $S_4$.

\item Icosahedral, $A_5$.
\end{enumerate}
all appearing as symmetry groups of regular polygons and polyhedra.
\end{theorem}

\begin{proof}
This is something truly remarbable, the idea being as follows:

\medskip

(1) To start with, we certainly have as examples the groups in the statement. Indeed, those in (1,2) come from Theorem 10.5, via the following standard embedding:
$$O_2\subset SO_3\quad,\quad U\to 
\begin{pmatrix}U&0\\0&\det U\end{pmatrix}$$ 

As for those in (3,4,5), these are the groups that we found in Theorem 10.10. 

\medskip

(2) Regarding now the converse, assume that $G\subset SO_3$ is finite. Given $g\in G-\{1\}$, consider its rotation axis coming from Theorem 10.11, and then the two points $\pm x$ where this axis intersects the unit sphere $S^2\subset\mathbb R^3$, called poles of $g$. We can consider then the set $X\subset S^2$ of all poles of all elements $g\in G-\{1\}$, and we have an action as follows:
$$G\curvearrowright X$$

(3) In order to exploit this latter action, we can use the following counting trick, due to Burnside, which is valid for any finite group action on a finite set, $G\curvearrowright X$:
\begin{eqnarray*}
\sum_{g\in G}|X^g|
&=&\sum_{x\in X}|G_x|\\
&=&|G|\sum_{x\in X}\frac{1}{|Gx|}\\
&=&|G|\sum_{O\in X/G}|O|\cdot\frac{1}{|O|}\\
&=&|G|\cdot|X/G|
\end{eqnarray*}

To be more precise, here $X^g\subset X$ is the set of fixed points by $g\in G$, and $G_x\subset G$ is the stabilizer of $x\in X$, and we have used the general theory from chapter 9.

\medskip

(4) Now let us see what the Burnside formula gives, for the action in (2). If we denote by $N$ the number of orbits of our action $G\curvearrowright X$, this formula reads:
$$|X|+2(|G|-1)=N|G|$$

Now observe that this latter formula can be further processed in the following way, with $\{x_1,\ldots,x_N\}\subset X$ being a set of representatives for the orbits of $G\curvearrowright X$:
\begin{eqnarray*}
2\left(1-\frac{1}{|G|}\right)
&=&N-\frac{|X|}{|G|}\\
&=&N-\frac{1}{|G|}\sum_{i=1}^N[G:G_{x_i}]\\
&=&\sum_{i=1}^N1-\frac{1}{|G_{x_i}|}
\end{eqnarray*}

(5) And the point is that this latter formula is exactly what we need. Indeed, observe that the left term and the right components are subject to the following estimates:
$$2\left(1-\frac{1}{|G|}\right)<2\quad,\quad1-\frac{1}{|G_{x_i}|}\geq\frac{1}{2}$$

We conclude that we must have $N=2,3$, which is a big win, we are almost there.

\medskip

(6) In practice now, in the case $N=2$, the formula that we found in (4) reads:
$$\frac{2}{|G|}=\frac{1}{|G_x|}+\frac{1}{|G_y|}$$

But a quick study shows that the solution here is $G=\mathbb Z_n$, corresponding to:
$$\frac{2}{n}=\frac{1}{n}+\frac{1}{n}$$

(7) Regarding now the case $N=3$, here the formula found in (4) reads:
$$1+\frac{2}{|G|}=\frac{1}{|G_x|}+\frac{1}{|G_y|}+\frac{1}{|G_z|}$$

But here we have 4 possible cases, corresponding to the following solutions of this:
$$1+\frac{2}{2n}=\frac{1}{2}+\frac{1}{2}+\frac{1}{n}\qquad,\qquad
1+\frac{2}{12}=\frac{1}{2}+\frac{1}{3}+\frac{1}{3}$$
$$1+\frac{2}{24}=\frac{1}{2}+\frac{1}{3}+\frac{1}{4}\qquad,\qquad 
1+\frac{2}{60}=\frac{1}{2}+\frac{1}{3}+\frac{1}{5}$$

And a study of these cases, that we will leave as an instructive exercise, leads to the other solutions in the statement, namely $G=D_n$, $G=A_4$, $G=S_4$, $G=A_5$.
\end{proof}

Very nice all this. We should mention that, with a bit more work, based on the above, the finite subgroups of $O_3$ can be classified too, using Proposition 10.9, and with this being something quite straightforward. We will leave this, again, as an instructive exercise.

\section*{10c. Euler-Rodrigues}

Moving forward, let us go back now to $N=2$ dimensions, but with a study in the complex case. We first have here the following result, which is elementary:

\index{unit sphere}

\begin{proposition}
We have the following formula,
$$SU_2=\left\{\begin{pmatrix}a&b\\ -\bar{b}&\bar{a}\end{pmatrix}\ \Big|\ |a|^2+|b|^2=1\right\}$$
which makes $SU_2$ isomorphic to the unit complex sphere $S^1_\mathbb C\subset\mathbb C^2$.
\end{proposition}

\begin{proof}
Indeed, according to the usual matrix rules, for a matrix $U=\binom{a\ b}{c\ d}$ having determinant 1, the complex isometry condition $U^*=U^{-1}$ reads:
$$\begin{pmatrix}\bar{a}&\bar{c}\\\bar{b}&\bar{d}\end{pmatrix}
=\begin{pmatrix}d&-b\\-c&a\end{pmatrix}$$

Thus $U$ is as in the statement, and with $|a|^2+|b|^2=1$ coming from $\det U=1$.
\end{proof}

Here is a useful reformulation of the above result, using real numbers:

\begin{proposition}
We have the formula
$$SU_2=\left\{\begin{pmatrix}x+iy&z+it\\ -z+it&x-iy\end{pmatrix}\ \Big|\ x^2+y^2+z^2+t^2=1\right\}$$
which makes $SU_2$ isomorphic to the unit real sphere $S^3_\mathbb R\subset\mathbb R^3$.
\end{proposition}

\begin{proof}
This is indeed self-explanatory, coming from Proposition 10.13.
\end{proof}

At a more advanced level now, here is yet another reformulation of what we have:

\index{Pauli matrices}
\index{quaternions}
\index{spin matrices}

\begin{theorem}
We have the following formula,
$$SU_2=\left\{xc_1+yc_2+zc_3+tc_4\ \Big|\ x^2+y^2+z^2+t^2=1\right\}$$
where $c_1,c_2,c_3,c_4$ are matrices given by
$$c_1=\begin{pmatrix}1&0\\ 0&1\end{pmatrix}\quad,\quad
c_2=\begin{pmatrix}i&0\\ 0&-i\end{pmatrix}\quad,\quad 
c_3=\begin{pmatrix}0&1\\ -1&0\end{pmatrix}\quad,\quad 
c_4=\begin{pmatrix}0&i\\ i&0\end{pmatrix}$$
called Pauli spin matrices.
\end{theorem}

\begin{proof}
According to Proposition 10.14 the elements $U\in SU_2$ are the matrices as follows, depending on parameters $x,y,z,t\in\mathbb R$ satisfying $x^2+y^2+z^2+t^2=1$:
$$U=x\begin{pmatrix}1&0\\ 0&1\end{pmatrix}
+y\begin{pmatrix}i&0\\ 0&-i\end{pmatrix}
+z\begin{pmatrix}0&1\\ -1&0\end{pmatrix}+
t\begin{pmatrix}0&i\\ i&0\end{pmatrix}$$

Thus, we are led to the conclusions in the statement.
\end{proof}

The above result is often the most convenient one, when dealing with $SU_2$. This is because the Pauli matrices have a number of remarkable properties, which are very useful when doing computations. These properties can be summarized as follows:

\begin{proposition}
The Pauli matrices multiply according to the formulae
$$c_2^2=c_3^2=c_4^2=-1$$
$$c_2c_3=-c_3c_2=c_4$$
$$c_3c_4=-c_4c_3=c_2$$
$$c_4c_2=-c_2c_4=c_3$$
they conjugate according to the following rules,
$$c_1^*=c_1\ ,\ c_2^*=-c_2\ ,\ c_3^*=-c_3\ ,\ c_4^*=-c_4$$
and they form an orthonormal basis of $M_2(\mathbb C)$, with respect to the scalar product
$$<a,b>=tr(ab^*)$$
with $tr:M_2(\mathbb C)\to\mathbb C$ being the normalized trace of $2\times 2$ matrices, $tr=Tr/2$.
\end{proposition}

\begin{proof}
The first two assertions, regarding the multiplication and conjugation rules for the Pauli matrices, follow from some elementary computations. As for the last assertion, this follows by using these rules. Indeed, the fact that the Pauli matrices are pairwise orthogonal follows from computations of the following type, for $i\neq j$:
$$<c_i,c_j>
=tr(c_ic_j^*)
=tr(\pm c_ic_j)
=tr(\pm c_k)
=0$$

As for the fact that the Pauli matrices have norm 1, this follows from:
$$<c_i,c_i>
=tr(c_ic_i^*)
=tr(\pm c_i^2)
=tr(c_1)
=1$$

Thus, we are led to the conclusion in the statement.
\end{proof}

Moving on, we would like to discuss now a key relation between $SU_2$ and $SO_3$. Let us start with the following construction, whose goal will become clear in a moment:

\index{adjoint action}

\begin{proposition}
The adjoint action $SU_2\curvearrowright M_2(\mathbb C)$, given by
$$T_U(M)=UMU^*$$
leaves invariant the following real vector subspace of $M_2(\mathbb C)$,
$$E=span_\mathbb R(c_1,c_2,c_3,c_4)$$
and we obtain in this way a group morphism $SU_2\to GL_4(\mathbb R)$.
\end{proposition}

\begin{proof}
We have two assertions to be proved, as follows:

\medskip

(1) We must first prove that, with $E\subset M_2(\mathbb C)$ being the real vector space in the statement, we have the following implication:
$$U\in SU_2,M\in E\implies UMU^*\in E$$

But this is clear from the multiplication rules for the Pauli matrices, from Proposition 10.16. Indeed, let us write our matrices $U,M$ as follows:
$$U=xc_1+yc_2+zc_3+tc_4$$
$$M=ac_1+bc_2+cc_3+dc_4$$

We know that the coefficients $x,y,z,t$ and $a,b,c,d$ are real, due to $U\in SU_2$ and $M\in E$. The point now is that when computing $UMU^*$, by using the various rules from Proposition 10.16, we obtain a matrix of the same type, namely a combination of $c_1,c_2,c_3,c_4$, with real coefficients. Thus, we have $UMU^*\in E$, as desired.

\medskip

(2) In order to conclude, let us identify $E\simeq\mathbb R^4$, by using the basis $c_1,c_2,c_3,c_4$. The result found in (1) shows that we have a correspondence as follows:
$$SU_2\to M_4(\mathbb R)\quad,\quad U\to (T_U)_{|E}$$

Now observe that for any $U\in SU_2$ and any $M\in M_2(\mathbb C)$ we have:
$$T_{U^*}T_U(M)=U^*UMU^*U=M$$

Thus $T_{U^*}=T_U^{-1}$, and so the correspondence that we found can be written as:
$$SU_2\to GL_4(\mathbb R)\quad,\quad U\to (T_U)_{|E}$$

But this a group morphism, due to the following computation:
$$T_UT_V(M)=UVMV^*U^*=T_{UV}(M)$$

Thus, we are led to the conclusion in the statement.
\end{proof}

The point now, which makes the link with $SO_3$, and which will ultimately elucidate the structure of $SO_3$, is that Proposition 10.17 can be improved as follows:

\begin{theorem}
The adjoint action $SU_2\curvearrowright M_2(\mathbb C)$ leaves invariant the space
$$F=span_\mathbb R(c_2,c_3,c_4)$$
and we obtain in this way a group morphism $SU_2\to SO_3$.
\end{theorem}

\begin{proof}
We can do this in several steps, as follows:

\medskip

(1) Our first claim is that the group morphism $SU_2\to GL_4(\mathbb R)$ constructed in Proposition 10.17 is in fact a morphism $SU_2\to O_4$. In order to prove this, recall the following formula, valid for any $U\in SU_2$, from the proof of Proposition 10.17:
$$T_{U^*}=T_U^{-1}$$

We want to prove that the matrices $T_U\in GL_4(\mathbb R)$ are orthogonal, and in view of the above formula, it is enough to prove that we have:
$$T_U^*=(T_U)^t$$

So, let us prove this. For any two matrices $M,N\in E$, we have:
\begin{eqnarray*}
<T_{U^*}(M),N>
&=&<U^*MU,N>\\
&=&tr(U^*MUN)\\
&=&tr(MUNU^*)
\end{eqnarray*}

On the other hand, we have as well the following formula:
\begin{eqnarray*}
<(T_U)^t(M),N>
&=&<M,T_U(N)>\\
&=&<M,UNU^*>\\
&=&tr(MUNU^*)
\end{eqnarray*}

Thus we have indeed $T_U^*=(T_U)^t$, which proves our $SU_2\to O_4$ claim.

\medskip

(2) In order now to finish, recall that we have by definition $c_1=1$, as a matrix. Thus, the action of $SU_2$ on the vector $c_1\in E$ is given by:
$$T_U(c_1)=Uc_1U^*=UU^*=1=c_1$$ 

We conclude that $c_1\in E$ is invariant under $SU_2$, and by orthogonality the following subspace of $E$ must be invariant as well under the action of $SU_2$:
$$e_1^\perp=span_\mathbb R(c_2,c_3,c_4)$$

Now if we call this subspace $F$, and we identify $F\simeq\mathbb R^3$ by using the basis $c_2,c_3,c_4$, we obtain by restriction to $F$ a morphism of groups as follows:
$$SU_2\to O_3$$

But since this morphism is continuous and $SU_2$ is connected, its image must be connected too. Now since the target group decomposes as $O_3=SO_3\sqcup(-SO_3)$, and $1\in SU_2$ gets mapped to $1\in SO_3$, the whole image must lie inside $SO_3$, and we are done.
\end{proof}

We can now formulate a key result, due to Euler-Rodrigues, as follows:

\begin{theorem}
We have a double cover map, obtained via the adjoint representation,
$$SU_2\to SO_3$$
and this map produces the Euler-Rodrigues formula
$$U=\begin{pmatrix}
x^2+y^2-z^2-t^2&2(yz-xt)&2(xz+yt)\\
2(xt+yz)&x^2+z^2-y^2-t^2&2(zt-xy)\\
2(yt-xz)&2(xy+zt)&x^2+t^2-y^2-z^2
\end{pmatrix}$$
for the generic elements of $SO_3$.
\end{theorem}

\begin{proof}
We have several things to be proved here, the idea being as follows:

\medskip

(1) Our first claim is that, with respect to the standard basis $c_1,c_2,c_3,c_4$ of the vector space $\mathbb R^4=span(c_1,c_2,c_3,c_4)$, the morphism $T:SU_2\to GL_4(\mathbb R)$ is given by:
$$T_U=\begin{pmatrix}
1&0&0&0\\
0&x^2+y^2-z^2-t^2&2(yz-xt)&2(xz+yt)\\
0&2(xt+yz)&x^2+z^2-y^2-t^2&2(zt-xy)\\
0&2(yt-xz)&2(xy+zt)&x^2+t^2-y^2-z^2
\end{pmatrix}$$

(2) Indeed, with notations from Proposition 10.17 and its proof, let us first look at the action $L:SU_2\curvearrowright\mathbb R^4$ by left multiplication, which is by definition given by:
$$L_U(M)=UM$$

In order to compute the matrix of this action, let us write, as usual:
$$U=xc_1+yc_2+zc_3+tc_4$$
$$M=ac_1+bc_2+cc_3+dc_4$$

By using the multiplication formulae in Proposition 10.16, we obtain:
\begin{eqnarray*}
UM
&=&(xc_1+yc_2+zc_3+tc_4)(ac_1+bc_2+cc_3+dc_4)\\
&=&(xa-yb-zc-td)c_1\\
&+&(xb+ya+zd-tc)c_2\\
&+&(xc-yd+za+tb)c_3\\
&+&(xd+yc-zb+ta)c_4
\end{eqnarray*}

We conclude that the matrix of the left action considered above is:
$$L_U=\begin{pmatrix}
x&-y&-z&-t\\
y&x&-t&z\\
z&t&x&-y\\
t&-z&y&x
\end{pmatrix}$$

(3) Similarly, let us look now at the action $R:SU_2\curvearrowright\mathbb R^4$ by right multiplication, which is by definition given by the following formula:
$$R_U(M)=MU^*$$

In order to compute the matrix of this action, let us write, as before:
$$U=xc_1+yc_2+zc_3+tc_4$$
$$M=ac_1+bc_2+cc_3+dc_4$$

By using the multiplication formulae in Proposition 10.16, we obtain:
\begin{eqnarray*}
MU^*
&=&(ac_1+bc_2+cc_3+dc_4)(xc_1-yc_2-zc_3-tc_4)\\
&=&(ax+by+cz+dt)c_1\\
&+&(-ay+bx-ct+dz)c_2\\
&+&(-az+bt+cx-dy)c_3\\
&+&(-at-bz+cy+dx)c_4
\end{eqnarray*}

We conclude that the matrix of the right action considered above is:
$$R_U=\begin{pmatrix}
x&y&z&t\\
-y&x&-t&z\\
-z&t&x&-y\\
-t&-z&y&x
\end{pmatrix}$$

(4) Now by composing, the matrix of the adjoint matrix in the statement is:
\begin{eqnarray*}
T_U
&=&R_UL_U\\
&=&\begin{pmatrix}
x&y&z&t\\
-y&x&-t&z\\
-z&t&x&-y\\
-t&-z&y&x
\end{pmatrix}
\begin{pmatrix}
x&-y&-z&-t\\
y&x&-t&z\\
z&t&x&-y\\
t&-z&y&x
\end{pmatrix}\\
&=&\begin{pmatrix}
1&0&0&0\\
0&x^2+y^2-z^2-t^2&2(yz-xt)&2(xz+yt)\\
0&2(xt+yz)&x^2+z^2-y^2-t^2&2(zt-xy)\\
0&2(yt-xz)&2(xy+zt)&x^2+t^2-y^2-z^2
\end{pmatrix}
\end{eqnarray*}

(5) Summarizing, we have proved our claim in (1). We conclude that, when looking at $T:SU_2\to GL_4(\mathbb R)$ as a group morphism $SU_2\to O_4$, what we have in fact is a group morphism $SU_2\to O_3$, and even $SU_2\to SO_3$, given by the Euler-Rodrigues formula.

\medskip

(6) Next, the kernel of this morphism is elementary to compute, as follows:
\begin{eqnarray*}
\ker(SU_2\to SO_3)
&=&\left\{U\in SU_2\Big|T_U(M)=M,\forall M\in E\right\}\\
&=&\left\{U\in SU_2\Big|Uc_i=c_iU,\forall i\right\}\\
&=&\{\pm1\}
\end{eqnarray*}

(7) Finally, in what regards the surjectivity, we can argue here for instance that since each rotation $U\in SO_3$ is uniquely determined by its rotation axis, plus its rotation angle $t\in[0,2\pi)$, we are led to the conclusion that $U$ is uniquely determined by an element of $SU_2/\{\pm 1\}$, and so appears indeed via the Euler-Rodrigues formula, as stated.
\end{proof}

Getting back now to our finite subgroup questions, we have:

\begin{theorem}[Klein]
The subgroups of $SU_2$ are as follows:
\begin{enumerate}
\item Cyclic, $\mathbb Z_n$.

\item Dicyclic, $DC_n$.

\item Binary tetrahedral, lifting $A_4$.

\item Binary octahedral, lifting $S_4$.

\item Binary icosahedral, lifting $A_5$.
\end{enumerate}
\end{theorem}

\begin{proof}
This is indeed something quite standard, from what we have, the idea being that the various groups in (2-5) appear as lifts via $SU_2\to SO_3$ of the groups in Theorem 10.12 (2-5). We will leave some further learning here as an instructive exercise.
\end{proof}

Good work that we did, but the story is not over with this, because we can talk about $SU_3$ as well. As usual, exercise for you, to learn more about all this.

\section*{10d. Symplectic groups}

We have learned many interesting things in small dimensions, and time now to discuss the high dimensions as well. We will be interested in finding uniform families of subgroups $G_N\subset O_N$ or $G_N\subset U_N$, either finite or continuous. Let us start our study with:

\index{bistochastic matrix}

\begin{definition}
A square matrix $M\in M_N(\mathbb C)$ is called bistochastic if each row and each column sum up to the same number:
$$\begin{matrix}
M_{11}&\ldots&M_{1N}&\to&\lambda\\
\vdots&&\vdots\\
M_{N1}&\ldots&M_{NN}&\to&\lambda\\
\downarrow&&\downarrow\\
\lambda&&\lambda
\end{matrix}$$
If this happens only for the rows, or only for the columns, the matrix is called row-stochastic, respectively column-stochastic.
\end{definition}

In what follows we will be interested in the unitary bistochastic matrices, which are quite interesting objects. As a first result, regarding such matrices, we have:

\begin{proposition}
For a unitary matrix $U\in U_N$, the following are equivalent:
\begin{enumerate}
\item $H$ is bistochastic, with sums $\lambda$.

\item $H$ is row stochastic, with sums $\lambda$, and $|\lambda|=1$.

\item $H$ is column stochastic, with sums $\lambda$, and $|\lambda|=1$.
\end{enumerate}
\end{proposition}

\begin{proof}
This is something that we know from chapter 7, with $(1)\iff(2)$ being elementary, and with the further equivalence with (3) coming by symmetry.
\end{proof}

The unitary bistochastic matrices are stable under a number of operations, and in particular under taking products. Thus, these matrices form a group. We have:

\index{bistochastic group}
\index{discrete Fourier transform}

\begin{theorem}
The real and complex bistochastic groups, which are the sets
$$B_N\subset O_N\quad,\quad 
C_N\subset U_N$$
consisting of matrices which are bistochastic, are isomorphic to $O_{N-1}$, $U_{N-1}$.
\end{theorem}

\begin{proof}
This is something that we know too from chapter 7. To be more precise, let us pick a matrix $F\in U_N$, such as the Fourier matrix $F_N$, satisfying the following condition, where $e_0,\ldots,e_{N-1}$ is the standard basis of $\mathbb C^N$, and where $\xi$ is the all-one vector:
$$Fe_0=\frac{1}{\sqrt{N}}\xi$$ 

We have then, by using the above property of $F$:
\begin{eqnarray*}
u\xi=\xi
&\iff&uFe_0=Fe_0\\
&\iff&F^*uFe_0=e_0\\
&\iff&F^*uF=diag(1,w)
\end{eqnarray*}

Thus we have isomorphisms as in the statement, given by $w_{ij}\to(F^*uF)_{ij}$.
\end{proof}

We will be back to $B_N,C_N$ later. Moving ahead now, as yet another basic example of a continuous group, we have the symplectic group $Sp_N$. Let us begin with:

\index{super-space}
\index{super-identity}

\begin{definition}
The ``super-space'' $\bar{\mathbb C}^N$ is the usual space $\mathbb C^N$, with its standard basis $\{e_1,\ldots,e_N\}$, with a chosen sign $\varepsilon=\pm 1$, and a chosen involution on the indices:
$$i\to\bar{i}$$
The ``super-identity'' matrix is $J_{ij}=\delta_{i\bar{j}}$ for $i\leq j$ and $J_{ij}=\varepsilon\delta_{i\bar{j}}$ for $i\geq j$.
\end{definition}

Up to a permutation of the indices, we have a decomposition $N=2p+q$, such that the involution is, in standard permutation notation:
$$(12)\ldots (2p-1,2p)(2p+1)\ldots (q)$$

Thus, up to a base change, the super-identity is as follows, where $N=2p+q$ and $\varepsilon=\pm 1$, with the $1_q$ block at right disappearing if $\varepsilon=-1$:
$$J=\begin{pmatrix}
0&1\ \ \ \\
\varepsilon 1&0_{(0)}\\
&&\ddots\\
&&&0&1\ \ \ \\
&&&\varepsilon 1&0_{(p)}\\
&&&&&1_{(1)}\\
&&&&&&\ddots\\
&&&&&&&1_{(q)}
\end{pmatrix}$$

In the case $\varepsilon=1$, the super-identity is the following matrix:
$$J_+(p,q)=\begin{pmatrix}
0&1\ \ \ \\
1&0_{(1)}\\
&&\ddots\\
&&&0&1\ \ \ \\
&&&1&0_{(p)}\\
&&&&&1_{(1)}\\
&&&&&&\ddots\\
&&&&&&&1_{(q)}
\end{pmatrix}$$

In the case $\varepsilon=-1$ now, the diagonal terms vanish, and the super-identity is:
$$J_-(p,0)=\begin{pmatrix}
0&1\ \ \ \\
-1&0_{(1)}\\
&&\ddots\\
&&&0&1\ \ \ \\
&&&-1&0_{(p)}
\end{pmatrix}$$

With the above notions in hand, we have the following result:

\index{super-orthogonal group}
\index{symplectic group}

\begin{theorem}
The super-orthogonal group, which is by definition
$$\bar{O}_N=\left\{U\in U_N\Big|U=J\bar{U}J^{-1}\right\}$$
with $J$ being the super-identity matrix, is as follows:
\begin{enumerate}
\item At $\varepsilon=1$ we have $\bar{O}_N=O_N$.

\item At $\varepsilon=-1$ we have $\bar{O}_N=Sp_N$.
\end{enumerate}
\end{theorem}

\begin{proof}
These is something quite tricky, the idea being as follows:

\medskip

(1) At $\varepsilon=1$, consider the root of unity $w=e^{\pi i/4}$, and let us set:
$$K=\frac{1}{\sqrt{2}}\begin{pmatrix}w&w^7\\ w^3&w^5\end{pmatrix}$$

This matrix $K$ is then unitary, and we have the following formula:
$$K\begin{pmatrix}0&1\\1&0\end{pmatrix}K^t=1$$

Thus the following matrix is unitary as well, and satisfies $CJC^t=1$:
$$C=\begin{pmatrix}K^{(1)}\\&\ddots\\&&K^{(p)}\\&&&1_q\end{pmatrix}$$

Now in terms of $V=CUC^*$, the relations $U=J\bar{U}J^{-1}=$ unitary simply read:
$$V=\bar{V}={\rm unitary}$$

We conclude that we have an isomorphism $\bar{O}_N=O_N$ as in the statement. 

\medskip

(2) At $\varepsilon=-1$, this depends a bit on what you call symplectic group $Sp_N$, and for our purposes here, we will take the above formula $Sp_N=\bar{O}_N$ as a definition for it.
\end{proof}

We can say more about the symplectic group $Sp_N$, as follows:

\begin{theorem}
The symplectic group $Sp_N\subset U_N$, which is by definition
$$Sp_N=\left\{U\in U_N\Big|U=J\bar{U}J^{-1}\right\}$$
with $J$ being as above, consists of the $SU_2$ patterned matrices,
$$U=\begin{pmatrix}
a&b&\ldots\\
-\bar{b}&\bar{a}\\
\vdots&&\ddots
\end{pmatrix}$$
which are unitary, $U\in U_N$. In particular, we have $Sp_2=SU_2$.
\end{theorem}

\begin{proof}
At $N=2$, to start with, given a matrix $U=\binom{a\ b}{c\ d}$, the condition $U=J\bar{U}J^{-1}$ reformulates as follows, which gives $d=\bar{a}$ and $c=-\bar{b}$, as desired:
$$\begin{pmatrix}a&b\\c&d\end{pmatrix}
\begin{pmatrix}0&1\\-1&0\end{pmatrix}
=\begin{pmatrix}0&1\\-1&0\end{pmatrix}
\begin{pmatrix}\bar{a}&\bar{b}\\\bar{c}&\bar{d}\end{pmatrix}
\iff
\begin{pmatrix}-b&a\\-d&c\end{pmatrix}
=\begin{pmatrix}\bar{c}&\bar{d}\\-\bar{a}&-\bar{b}\end{pmatrix}$$

In the general case, $N\in2\mathbb N$, the proof is similar, with the condition $U=J\bar{U}J^{-1}$ corresponding precisely to the fact that $U$ must be $SU_2$-patterned, as stated.
\end{proof}

We will be back later to the symplectic groups, towards the end of the present book, with more results about them. In the meantime, have a look at the mechanics book of Arnold \cite{ar2}, which explains what the symplectic groups and geometry are good for.

\bigskip

As a last topic of discussion, now that we have a decent understanding of the main continuous groups of unitary matrices $G\subset U_N$, let us go back to the finite groups from the beginning of this chapter, and make a link with the material there. We first have:

\begin{theorem}
The full complex reflection group $K_N\subset U_N$, given by
$$K_N=M_N(\mathbb T\cup\{0\})\cap U_N$$
decomposes as $K_N=\mathbb T\wr S_N$, with $S_N$ acting on $\mathbb T^N$ by permuting the factors.
\end{theorem}

\begin{proof}
This is something quite similar to what we know from chapter 9 regarding the hyperoctahedral group $H_N\subset O_N$, and we will leave the various details here as an exercise. With the comment that we will be back to this later, in chapter 12. 
\end{proof}

Next, we can talk about the reflection subgroup of any subgroup $G\subset U_N$, as follows: 

\begin{definition}
Given $G\subset U_N$, we can define its reflection subgroup to be
$$K=G\cap K_N$$
with the intersection taken inside $U_N$.
\end{definition}

Many things can be said in relation with this, but let us not stop here. Indeed, given an intermediate subgroup $H_N\subset G\subset U_N$, we can view it as follows:
$$\xymatrix@R=20pt@C=20pt{
K_N\ar[rr]&&U_N\\
&G\ar@.[ur]&\\
H_N\ar[rr]\ar[uu]\ar@.[ur]&&O_N\ar[uu]}$$

Thus, we have some sort of 2D orientation for the subgroups $H_N\subset G\subset U_N$, and this suggests extending the construction in Definition 10.28, in the following way:

\begin{definition}
Associated to any intermediate compact group $H_N\subset G\subset U_N$ are its discrete, real, complex and smooth versions, given by the formulae
$$G^d=G\cap K_N\quad,\quad 
G^r=G\cap O_N$$
$$G^c=<G,K_N>\quad,\quad 
G^s=<G,O_N>$$
with $<\,,>$ being the topological generation operation, involving taking a closure.
\end{definition}

But with this in hand, it is natural now to formulate the following definition:

\index{oriented group}

\begin{definition}
A compact group $H_N\subset G\subset U_N$ is called oriented if
$$\xymatrix@R=38pt@C=43pt{
K_N\ar[r]&G^c\ar[r]&U_N\\
G^d\ar[u]\ar[r]&G\ar[r]\ar[u]&G^s\ar[u]\\
H_N\ar[r]\ar[u]&G^r\ar[u]\ar[r]&O_N\ar[u]}$$
is an intersection and generation diagram, in the sense that any of its square subdiagrams $A\subset B,C\subset D$ satisfies $A=B\cap C$ and $D=<B,C>$.
\end{definition}

And this notion is quite interesting, because most of the basic examples of closed subgroups $G\subset U_N$, finite or continuous, are oriented. In fact, we have:

\begin{question}
What are the oriented groups $H_N\subset G\subset U_N$? What about the oriented groups coming in families, $G=(G_N)$, with $N\in\mathbb N$?
\end{question}

And we will stop here our discussion, sometimes a good question is better as hunting trophy than a final theorem, or at least that's what my cats say. We will be back to this in Part IV  below, under a number of supplementary assumptions on the groups $G$ that we consider, which will allow us to derive a number of classification results.

\section*{10e. Exercises}

There has been a lot of theory in this chapter, and this is just the tip of the iceberg, on what can be said about the rotation groups. As a first exercise, we have:

\begin{exercise}
Prove that for a convex polyhedron we have the Euler formula
$$v+f=e+2$$
with $v,e,f$ being the number of vertices, edges and faces. 
\end{exercise}

This is normally not very difficult, by recurrence, and as a bonus exercise, reprove the Plato theorem by using this, with the data for regular polyhedra being as follows:
$$\begin{array}{cccccccccccc}
&T&O&C&I&D\\
v&4&6&8&12&20\\
f&4&8&6&20&12\\
e&6&12&12&30&30
\end{array}$$

As second bonus exercise, learn also about the Euler formula for planar graphs, and for higher genus graphs. There are many interesting things here to be learned.

\begin{exercise}
Work out all the details of the Euler-Rodrigues formula, by using the fact that any rotation in $\mathbb R^3$ has a rotation axis.
\end{exercise}

Here the problem, once the rotation axis found, is that of drawing the picture, identifying the relevant angles, and then doing the math in terms of these angles.

\begin{exercise}
Work out the theory of the subgroups of $O_N,U_N$ constructed via
$$(\det U)^d=1$$
with $d\in\mathbb N\cup\{\infty\}$, which generalize both $O_N,U_N$ and $SO_N,SU_N$.
\end{exercise}

There are many things that can be done here, and the more, the better.

\begin{exercise}
Look up the literature, and find the relevance of the symplectic groups, and of symplectic geometry in general, to questions in physics.
\end{exercise}

As before with the previous exercise, many things that can be learned and done here, especially from classical mechanics books, and the more you learn, the better.

\begin{exercise}
Find and then write down a brief account of the Shephard-Todd theorem, stating that the irreducible complex reflection groups are
$$H_N^{sd}=\left\{U\in M_N(\mathbb Z_s\cup\{0\})\cap U_N\Big|(\det U)^d=1\right\}$$
along with a number of exceptional examples, more precisely $34$ of them.
\end{exercise}

As before with the previous exercises, the more you learn here, the better.

\chapter{Symmetric groups}

\section*{11a. Character laws}

We would like to develop in what follows some general theory for the compact subgroups $G\subset U_N$, usually taken finite, with our main example being the symmetric group $S_N\subset O_N$. Let us start with a notion that we already met in chapter 9, namely:

\index{representation}
\index{character}

\begin{definition}
A representation of a finite group $G$ is a group morphism
$$u:G\to U_N$$
into a unitary group. The character of such a representation is the function
$$\chi:G\to\mathbb C\quad,\quad 
g\to Tr(u_g)$$
where $Tr$ is the usual, unnormalized trace of the $N\times N$ matrices.
\end{definition}

As explained in chapter 9, the simplest case of all this, namely $N=1$, is of particular importance. Here the representations coincide with their characters, and are by definition the group morphisms as follows, called characters of the group:
$$\chi:G\to\mathbb T$$

These characters from an abelian group $\widehat{G}$, and when $G$ itself is abelian, the correspondence $G\to\widehat{G}$ is a duality, in the sense that it maps $\widehat{G}\to G$ as well. Moreover, a more detailed study shows that we have in fact an isomorphism $G\simeq\widehat{G}$, with this being something quite subtle, related at the same time to the structure theorem for the finite abelian groups, $G\simeq\mathbb Z_{N_1}\times\ldots\times\mathbb Z_{N_k}$, and to the Fourier transforms over such groups.

\bigskip

Let us summarize this discussion, along with a little more, as follows:

\begin{theorem}
The characters of a finite group $\chi:G\to\mathbb T$ factorize as
$$\chi:G\to G_{ab}\to\mathbb T$$
with $G_{ab}$ being the abelianization of $G$, given by the formula
$$G_{ab}=G\big/\big<gh=hg\big>$$
and so correspond to the elements of the dual $\widehat{G}_{ab}\simeq G_{ab}$ of this abelianization.
\end{theorem}

\begin{proof}
Here the fact that the characters factorize indeed as $\chi:G\to G_{ab}\to\mathbb T$ is clear from definitions, and the last assertion comes from the discussion above.
\end{proof}

In what follows we will be interested in the general case, $N\in\mathbb N$. It is technically convenient to assume that the representation $u:G\to U_N$ is faithful, by replacing if necessary $G$ with its image. Thus, we are led to the following definition:

\index{main character}

\begin{definition}
The main character of a compact group $G\subset U_N$ is the map
$$\chi:G\to\mathbb C\quad,\quad 
g\to Tr(g)$$
which associates to the group elements, viewed as unitary matrices, their trace.
\end{definition}

We will see in a moment some motivations for the study of these characters. From a naive viewpoint, which is ours at the present stage, we want to do some linear algebra with our group elements $g\in U_N$, and we have several choices here, as follows:

\medskip

(1) A first idea would be to look at the determinant, $\det g\in\mathbb T$. However, this is usually not a very interesting quantity, for instance because $g\in O_N$ implies $\det g=\pm1$. Also, for groups like $SO_N,SU_N$, this determinant is by definition 1.

\medskip

(2) A second idea would be to try to compute eigenvalues and eigenvectors for the group elements $g\in G$, and then solve diagonalization questions for these elements. However, all this is quite complicated, so this idea is not good either.

\medskip

(3) Thus, we are left with looking at the trace, $Tr(g)\in\mathbb C$. We will see soon that this is a very reasonable choice, with the mathematics being at the same time non-trivial, doable, and also interesting, for a whole number of reasons.

\medskip

Before starting our study, let us mention as well the more advanced reasons leading to the study of characters. The idea here is that a given finite or compact group $G$ can have several representations $\pi:G\to U_N$, and these representations can be studied via their characters $\chi_\pi:G\to\mathbb C$, with a well-known and deep theorem basically stating that $\pi$ can be recovered from its character $\chi_\pi$. We will be back to this later.

\bigskip

As a basic result now regarding the characters, we have:

\index{central function}
\index{conjugacy classes}

\begin{theorem}
Given a compact group $G\subset U_N$, its main character $\chi:G\to\mathbb C$ is a central function, in the sense that it satisfies the following condition:
$$\chi(gh)=\chi(hg)$$
Equivalently, $\chi$ is constant on the conjugacy classes of $G$.
\end{theorem}

\begin{proof}
This is clear from the fact that the trace of matrices satisfies:
$$Tr(AB)=Tr(BA)$$

Thus, we are led to the conclusion in the statement.
\end{proof}

As before, there is some interesting mathematics behind all this. We will prove later, when doing representation theory, that any central function $f:G\to\mathbb C$ appears as a linear combination of characters $\chi_\pi:G\to\mathbb C$ of representations $\pi:G\to U_N$.

\bigskip

In order to work out now some examples, let us get back now to our main examples of finite groups, constructed in chapter 9, which were as follows:
$$\mathbb Z_N\subset D_N\subset S_N\subset H_N$$

We will do in what follows some character computations for these groups. Let us start with the following result, which covers $\mathbb Z_N\subset D_N\subset S_N$, or rather tells us what is to be done with these groups, in relation with their main characters:

\index{fixed points}
\index{symmetric group}
\index{main character}

\begin{proposition}
For the symmetric group, regarded as group of permutation matrices, $S_N\subset O_N$, the main character counts the number of fixed points:
$$\chi(g)=\#\left\{i\in\{1,\ldots,N\}\Big|\sigma(i)=i\right\}$$
The same goes for any $G\subset S_N$, regarded as a matrix group via $G\subset S_N\subset O_N$.
\end{proposition}

\begin{proof}
This is indeed clear from definitions, because the diagonal entries of the permutation matrices correspond to the fixed points of the permutation.
\end{proof}

Summarizing, we are left with counting fixed points. For the simplest possible group, namely the cyclic group $\mathbb Z_N\subset S_N$, the computation is as follows:

\index{cyclic group}

\begin{proposition}
The main character of $\mathbb Z_N\subset O_N$ is given by:
$$\chi(g)=\begin{cases}
0&{\rm if}\ g\neq1\\
N&{\rm if}\ g=1
\end{cases}$$
Thus, at the probabilistic level, we have the following formula,
$$law(\chi)=\left(1-\frac{1}{N}\right)\delta_0+\frac{1}{N}\delta_N$$
telling us that the main character $\chi$ follows a Bernoulli law.
\end{proposition}

\begin{proof}
The first formula is clear, because the cyclic permutation matrices have 0 on the diagonal, and so 0 as trace, unless the matrix is the identity, having trace $N$. As for the second formula, this is a probabilistic reformulation of the first one.
\end{proof}

For the dihedral group now, which is the next one in our hierarchy, the computation is more interesting, and the final answer is no longer uniform in $N$, as follows:

\index{dihedral group}

\begin{proposition}
For the dihedral group $D_N\subset S_N$ we have
$$law(\chi)=\begin{cases}
\left(\frac{3}{4}-\frac{1}{2N}\right)\delta_0+\frac{1}{4}\delta_2+\frac{1}{2N}\delta_N&(N\ even)\\
&\\
\left(\frac{1}{2}-\frac{1}{2N}\right)\delta_0+\frac{1}{2}\delta_1+\frac{1}{2N}\delta_N&(N\ odd)
\end{cases}$$
with this law being no longer uniform in $N$.
\end{proposition}

\begin{proof}
The dihedral group $D_N$ consists indeed of:

\medskip

--  $N$ symmetries, having each $1$ fixed point when $N$ is odd, and having 0 or 2 fixed points, distributed $50-50$, when $N$ is even.

\medskip

-- $N$ rotations, each having $0$ fixed points, except for the identity, which is technically a rotation too, and which has $N$ fixed points.

\medskip

Thus, we are led to the formulae in the statement.
\end{proof}

Regarding now the symmetric group $S_N$ itself, the permutations having no fixed points at all are called derangements, and the first question which appears, which is a classical question in combinatorics, is that of counting these derangements. We will need:

\begin{proposition}
We have the following formula,
$$\left|\left(\bigcup_iA_i\right)^c\right|=|A|-\sum_i|A_i|+\sum_{i<j}|A_i\cap A_j|-\sum_{i<j<k}|A_i\cap A_j\cap A_k|+\ldots$$
called inclusion-exclusion principle.
\end{proposition}

\begin{proof}
This is indeed quite clear, by thinking a bit, as before, as follows:

\medskip

(1) In order to count $(\cup_iA_i)^c$, we certainly have to start with $|A|$.

\medskip

(2) Then, we obviously have to remove each $|A_i|$, and so remove $\sum_i|A_i|$.

\medskip

(3) But then, we have to put back each $|A_i\cap A_j|$, and so put back $\sum_{i<j}|A_i\cap A_j|$.

\medskip

$\vdots$

\medskip

(4) And so on, which leads to the formula in the statement.
\end{proof}

We can now do the computation for $S_N$, leading to the following remarkable result:

\index{random permutation}
\index{derangement}
\index{fixed points}

\begin{theorem}
The probability for a random $\sigma\in S_N$ to be a derangement is:
$$P=1-\frac{1}{1!}+\frac{1}{2!}-\ldots+(-1)^{N-1}\frac{1}{(N-1)!}+(-1)^N\frac{1}{N!}$$
Thus, we have the following asymptotic formula, in the $N\to\infty$ limit,
$$P\simeq\frac{1}{e}$$
where $e=2.7182\ldots$ is the usual constant from analysis.
\end{theorem}

\begin{proof}
This is something very classical, which is best viewed by using the inclusion-exclusion principle. Consider indeed the following sets of permutations:
$$S_N^i=\left\{\sigma\in S_N\Big|\sigma(i)=i\right\}$$

The set of permutations having no fixed points, or derangements, is then:
$$X_N=\left(\bigcup_i S_N^i\right)^c$$

In order to compute now the cardinality $|X_N|$, consider as well the following sets, depending on indices $i_1<\ldots<i_k$, obtained by taking intersections:
$$S_N^{i_1\ldots i_k}=S_N^{i_1}\cap\ldots\cap S_N^{i_k}$$

In other words, these latter sets are given by the following formula:
$$S_N^{i_1\ldots i_k}=\left\{\sigma\in S_N\Big|\sigma(i_1)=i_1,\ldots,\sigma(i_k)=i_k\right\}$$

The inclusion-exclusion principle tells us that we have:
$$|X_N|=|S_N|-\sum_i|S_N^i|+\sum_{i<j}|S_N^{ij}|-\ldots+(-1)^N\sum_{i_1<\ldots<i_N}|S_N^{i_1\ldots i_N}|$$

Thus, the probability that we are interested in is given by:
\begin{eqnarray*}
P
&=&\frac{1}{N!}\left(|S_N|-\sum_i|S_N^i|+\sum_{i<j}|S_N^{ij}|-\ldots+(-1)^N\sum_{i_1<\ldots<i_N}|S_N^{i_1\ldots i_N}|\right)\\
&=&\frac{1}{N!}\sum_{k=0}^N(-1)^k\sum_{i_1<\ldots<i_k}|S_N^{i_1\ldots i_k}|\\
&=&\frac{1}{N!}\sum_{k=0}^N(-1)^k\sum_{i_1<\ldots<i_k}(N-k)!\\
&=&\frac{1}{N!}\sum_{k=0}^N(-1)^k\binom{N}{k}(N-k)!\\
&=&\sum_{k=0}^N\frac{(-1)^k}{k!}
\end{eqnarray*}

Since on the right we have the expansion of $1/e$, we obtain the result.
\end{proof}

The above result is something remarkable, and there are many versions and generalizations of it. We will discuss this gradually, in what follows, all this being key material. To start with, in terms of characters, the above result reformulates as follows:

\index{main character}

\begin{proposition}
For the symmetric group $S_N\subset O_N$, the probability for main character $\chi:S_N\to\mathbb N$ to vanish is given by the following formula:
$$P(\chi=0)=1-\frac{1}{1!}+\frac{1}{2!}-\ldots+(-1)^{N-1}\frac{1}{(N-1)!}+(-1)^N\frac{1}{N!}$$
Thus we have the formula $P(\chi=0)\simeq1/e$, in the $N\to\infty$ limit.
\end{proposition}

\begin{proof}
This follows indeed by combining Proposition 11.5, which tells us that $\chi$ counts the number of fixed points, with Theorem 11.9.
\end{proof}

Let us discuss now, more generally, what happens when counting permutations having exactly $k$ fixed points. The result here, extending Theorem 11.9, is as follows:

\begin{theorem}
The probability for a random permutation $\sigma\in S_N$ to have exactly $k$ fixed points is given by the following formula:
$$P=\frac{1}{k!}\left(1-\frac{1}{1!}+\frac{1}{2!}-\ldots+(-1)^{N-1}\frac{1}{(N-1)!}+(-1)^N\frac{1}{N!}\right)$$
Thus we have the formula $P\simeq1/(ek!)$, in the $N\to\infty$ limit.
\end{theorem}

\begin{proof}
We already know, from Theorem 11.9, that this formula holds at $k=0$. In the general case now, we have to count the permutations $\sigma\in S_N$ having exactly $k$ points. Since having such a permutation amounts in choosing $k$ points among $1,\ldots,N$, and then permuting the $N-k$ points left, without fixed points allowed, we have:
\begin{eqnarray*}
\#\left\{\sigma\in S_N\Big|\chi(\sigma)=k\right\}
&=&\binom{N}{k}\#\left\{\sigma\in S_{N-k}\Big|\chi(\sigma)=0\right\}\\
&=&\frac{N!}{k!(N-k)!}\#\left\{\sigma\in S_{N-k}\Big|\chi(\sigma)=0\right\}\\
&=&N!\times\frac{1}{k!}\times\frac{\#\left\{\sigma\in S_{N-k}\Big|\chi(\sigma)=0\right\}}{(N-k)!}
\end{eqnarray*}

Now by dividing everything by $N!$, we obtain from this the following formula:
$$\frac{\#\left\{\sigma\in S_N\Big|\chi(\sigma)=k\right\}}{N!}=\frac{1}{k!}\times\frac{\#\left\{\sigma\in S_{N-k}\Big|\chi(\sigma)=0\right\}}{(N-k)!}$$

By using now the computation at $k=0$, that we already have, from Theorem 11.9, it follows that with $N\to\infty$ we have the following estimate:
\begin{eqnarray*}
P(\chi=k)
&\simeq&\frac{1}{k!}\cdot P(\chi=0)\\
&\simeq&\frac{1}{k!}\cdot\frac{1}{e}
\end{eqnarray*}

Thus, we are led to the conclusion in the statement.
\end{proof}

As before, in regards with derangements, we can reformulate what we found in terms of the main character, and we obtain in this way the following statement:

\begin{theorem}
For the symmetric group $S_N\subset O_N$, the distribution of the main character $\chi:S_N\to\mathbb N$ is given by the following formula:
$$P(\chi=k)=\frac{1}{k!}\left(1-\frac{1}{1!}+\frac{1}{2!}-\ldots+(-1)^{N-1}\frac{1}{(N-1)!}+(-1)^N\frac{1}{N!}\right)$$
Thus we have the following asymptotic formula, in the $N\to\infty$ limit,
$$P(\chi=k)\simeq\frac{1}{ek!}$$
with $e=2.7182\ldots$ being the usual constant from analysis.
\end{theorem}

\begin{proof}
This follows indeed by combining Proposition 11.5, which tells us that $\chi$ counts the number of fixed points, with Theorem 11.11.
\end{proof}

\section*{11b. Poisson limits}

In order to best interpret the above results, we will need some probability theory. We already met the Poisson laws in chapter 6, but the discussion there was quite brief, and time now to review all this in detail. We first have the following definition:

\index{Poisson law}

\begin{definition}
The Poisson law of parameter $1$ is the following measure,
$$p_1=\frac{1}{e}\sum_{k\in\mathbb N}\frac{\delta_k}{k!}$$
and the Poisson law of parameter $t>0$ is the following measure,
$$p_t=e^{-t}\sum_{k\in\mathbb N}\frac{t^k}{k!}\,\delta_k$$
with the letter ``p'' standing for Poisson.
\end{definition}

Observe that these laws have indeed mass 1, as they should, and this due to the following well-known formula, which is the foundational formula of calculus:
$$e^t=\sum_k\frac{t^k}{k!}$$

We will see in the moment why these measures appear a bit everywhere, in discrete contexts, the reasons behind this coming from the Poisson Limit Theorem (PLT). Let us first develop some general theory. We first have the following result:

\begin{proposition}
The mean and variance of the Poisson law $p_t$ are 
$$E=V=t$$
for any $t>0$. In particular, at $t=1$ we have $E=V=1$.
\end{proposition}

\begin{proof}
In what regards the mean of the Poisson law $p_t$, this is given by:
$$E=e^{-t}\sum_{k\geq1}\frac{t^kk}{k!}=e^{-t}\sum_{k\geq1}\frac{t^k}{(k-1)!}=e^{-t}\times te^t=t$$

Let us compute now the second moment. This can be done as follows:
\begin{eqnarray*}
M_2
&=&e^{-t}\sum_{k\geq1}\frac{t^kk^2}{k!}\\
&=&e^{-t}\left(\sum_{k\geq1}\frac{t^k(k-1)}{(k-1)!}+\sum_{k\geq1}\frac{t^k}{(k-1)!}\right)\\
&=&e^{-t}(t^2e^t+te^t)\\
&=&t^2+t
\end{eqnarray*}

Thus, the variance is $V=(t^2+t)-t^2=t$, as claimed.
\end{proof}

At a more advanced level now, we first have the following result:

\index{convolution}

\begin{theorem}
We have the following formula, for any $s,t>0$,
$$p_s*p_t=p_{s+t}$$
so the Poisson laws form a convolution semigroup.
\end{theorem}

\begin{proof}
We know that the convolution of Dirac masses is given by $\delta_k*\delta_l=\delta_{k+l}$, and by using this formula and the binomial formula, we obtain:
\begin{eqnarray*}
p_s*p_t
&=&e^{-s}\sum_k\frac{s^k}{k!}\,\delta_k*e^{-t}\sum_l\frac{t^l}{l!}\,\delta_l\\
&=&e^{-s-t}\sum_{kl}\frac{s^kt^l}{k!l!}\delta_{k+l}\\
&=&e^{-s-t}\sum_n\delta_n\sum_{k+l=n}\frac{s^kt^l}{k!l!}\\
&=&e^{-s-t}\sum_n\frac{\delta_n}{n!}\sum_{k+l=n}\frac{n!}{k!l!}s^kt^l\\\
&=&e^{-s-t}\sum_n\frac{(s+t)^n}{n!}\,\delta_n\\
&=&p_{s+t}
\end{eqnarray*}

Thus, we are led to the conclusion in the statement.
\end{proof}

Along the same lines, we have as well the following result:

\index{formal exponential}
\index{convolution}

\begin{theorem}
The Poisson laws appear as formal exponentials
$$p_t=\sum_k\frac{t^k(\delta_1-\delta_0)^{*k}}{k!}$$
with respect to the convolution of measures $*$.
\end{theorem}

\begin{proof}
By using the binomial formula, the measure at right is:
\begin{eqnarray*}
\mu
&=&\sum_k\frac{t^k}{k!}\sum_{p+q=k}(-1)^q\frac{k!}{p!q!}\delta_p\\
&=&\sum_kt^k\sum_{p+q=k}(-1)^q\frac{\delta_p}{p!q!}\\
&=&\sum_p\frac{t^p\delta_p}{p!}\sum_q\frac{(-1)^q}{q!}\\
&=&\frac{1}{e}\sum_p\frac{t^p\delta_p}{p!}\\
&=&p_t
\end{eqnarray*}

Thus, we are led to the conclusion in the statement.
\end{proof}

As in the continuous case, for the normal laws, our main tool for dealing with the Poisson laws will be the Fourier transform. The formula here is as follows:

\index{Poisson law}
\index{Fourier transform}

\begin{theorem}
The Fourier transform of $p_t$ is given by
$$F_{p_t}(x)=\exp\left((e^{ix}-1)t\right)$$
for any $t>0$.
\end{theorem}

\begin{proof}
We know that the Fourier transform of a variable $f$ is given, by definition, by the formula $F_f(x)=E(e^{ixf})$. We therefore obtain the following formula:
\begin{eqnarray*}
F_{p_t}(x)
&=&e^{-t}\sum_k\frac{t^k}{k!}F_{\delta_k}(x)\\
&=&e^{-t}\sum_k\frac{t^k}{k!}\,e^{ikx}\\
&=&e^{-t}\sum_k\frac{(e^{ix}t)^k}{k!}\\
&=&\exp(-t)\exp(e^{ix}t)\\
&=&\exp\left((e^{ix}-1)t\right)
\end{eqnarray*}

Thus, we have reached to the formula in the statement.
\end{proof}

Observe that we obtain in this way another proof for the convolution semigroup property of the Poisson laws, that we established above, by using the fact, that we know from chapter 6, that the logarithm of the Fourier transform linearizes the convolution.

\bigskip

We can now establish the Poisson Limit Theorem (PLT), as follows:

\index{Poisson limit theorem}
\index{PLT}

\begin{theorem}
We have the following convergence, in moments,
$$\left(\left(1-\frac{t}{n}\right)\delta_0+\frac{t}{n}\delta_1\right)^{*n}\to p_t$$
for any $t>0$.
\end{theorem}

\begin{proof}
Let us denote by $\mu_n$ the measure under the convolution sign:
$$\mu_n=\left(1-\frac{t}{n}\right)\delta_0+\frac{t}{n}\delta_1$$

We have then the following computation, for the law in the statement: 
\begin{eqnarray*}
F_{\delta_r}(x)=e^{irx}
&\implies&F_{\mu_n}(x)=\left(1-\frac{t}{n}\right)+\frac{t}{n}e^{ix}\\
&\implies&F_{\mu_n^{*n}}(x)=\left(\left(1-\frac{t}{n}\right)+\frac{t}{n}e^{ix}\right)^n\\
&\implies&F_{\mu_n^{*n}}(x)=\left(1+\frac{(e^{ix}-1)t}{n}\right)^n\\
&\implies&F(x)=\exp\left((e^{ix}-1)t\right)
\end{eqnarray*}

Thus, we obtain the Fourier transform of $p_t$, as desired.
\end{proof}

There are of course many other things that can be said about the PLT, including examples and illustrations, and more technical results regarding the convergence, and we refer here to any standard probability book, such as Feller \cite{fel} or Durrett \cite{dur}. In what follows, we will be rather doing more combinatorics. To start with, we have:

\index{Bell numbers}

\begin{theorem}
The moments of $p_1$ are the Bell numbers,
$$M_k(p_1)=|P(k)|$$
where $P(k)$ is the set of partitions of $\{1,\ldots,k\}$.
\end{theorem}

\begin{proof}
The moments of $p_1$ are given by the following formula:
$$M_k=\frac{1}{e}\sum_{n\geq1}\frac{n^k}{n!}$$

We therefore have the following recurrence formula, for these moments:
\begin{eqnarray*}
M_{k+1}
&=&\frac{1}{e}\sum_{n\geq1}\frac{n^{k+1}}{n!}\\
&=&\frac{1}{e}\sum_{m\geq0}\frac{(m+1)^k}{m!}\\
&=&\frac{1}{e}\sum_{m\geq0}\frac{m^k}{m!}\left(1+\frac{1}{m}\right)^k\\
&=&\frac{1}{e}\sum_{m\geq0}\frac{m^k}{m!}\sum_{s=0}^k\binom{k}{s}m^{-s}\\
&=&\sum_{s=0}^k\binom{k}{s}\cdot\frac{1}{e}\sum_{m\geq0}\frac{m^{k-s}}{m!}\\
&=&\sum_{s=0}^k\binom{k}{s}M_{k-s}
\end{eqnarray*}

Next, let us try now to find a recurrence for the Bell numbers:
$$B_k=|P(k)|$$

A partition of $\{1,\ldots,k+1\}$ appears by choosing $s$ neighbors for $1$, among the $k$ numbers available, and then partitioning the $k-s$ elements left. Thus, we have:
$$B_{k+1}=\sum_{s=0}^k\binom{k}{s}B_{k-s}$$

Thus, the numbers $M_k$ satisfy the same recurrence as the numbers $B_k$. Regarding now the initial values, for the moments of $p_1$, according to Proposition 11.14, these are:
$$M_0=1\quad,\quad 
M_1=1$$

Now by using the above recurrence for the moments, we obtain from this:
$$M_2
=\sum_s\binom{1}{s}M_{k-s}
=1+1
=2$$

Thus, we can say that the initial values for the moments of $p_1$ are:
$$M_1=1\quad,\quad 
M_2=2$$

As for the Bell numbers, here the initial values are as follows:
$$B_1=1\quad,\quad 
B_2=2$$

Thus the initial values coincide, and so these numbers are equal, as stated.
\end{proof}

More generally, we have the following result, regarding $p_t$ with $t>0$:

\begin{theorem}
The moments of $p_t$ are given by
$$M_k(p_t)=\sum_{\pi\in P(k)}t^{|\pi|}$$
where $|.|$ is the number of blocks.
\end{theorem}

\begin{proof}
Observe first that the formula in the statement generalizes the one in Theorem 11.19, because at $t=1$ we obtain, as we should:
$$M_k(p_1)
=\sum_{\pi\in P(k)}1^{|\pi|}
=|P(k)|
=B_k$$

In general now, the moments of $p_t$ with $t>0$ are given by:
$$M_k=e^{-t}\sum_{n\geq1}\frac{t^nn^k}{n!}$$

We therefore have the following recurrence formula, for these moments:
\begin{eqnarray*}
M_{k+1}
&=&e^{-t}\sum_{n\geq1}\frac{t^nn^{k+1}}{n!}\\
&=&e^{-t}\sum_{m\geq0}\frac{t^{m+1}(m+1)^k}{m!}\\
&=&e^{-t}\sum_{m\geq0}\frac{t^{m+1}m^k}{m!}\left(1+\frac{1}{m}\right)^k\\
&=&e^{-t}\sum_{m\geq0}\frac{t^{m+1}m^k}{m!}\sum_{s=0}^k\binom{k}{s}m^{-s}\\
&=&\sum_{s=0}^k\binom{k}{s}\cdot e^{-t}\sum_{m\geq0}\frac{t^{m+1}m^{k-s}}{m!}\\
&=&t\sum_{s=0}^k\binom{k}{s}M_{k-s}
\end{eqnarray*}

As for the initial values, according to Proposition 11.14, these are as follows:
$$M_1=t\quad,\quad 
M_2=t+t^2$$

On the other hand, consider the numbers in the statement, namely:
$$S_k=\sum_{\pi\in P(k)}t^{|\pi|}$$

Since a partition of $\{1,\ldots,k+1\}$ appears by choosing $s$ neighbors for $1$, among the $k$ numbers available, and then partitioning the $k-s$ elements left, we have:
$$S_{k+1}=t\sum_{s=0}^k\binom{k}{s}S_{k-s}$$

As for the initial values of these numbers, these are as follows:
$$S_1=t\quad,\quad 
S_2=t+t^2$$

Thus the initial values coincide, so these numbers are the moments, as stated.
\end{proof}

Observe the analogy with the moment formulae for $g_t$ and $G_t$, from chapter 6. To be more precise, the moments of the main laws come from partitions, as follows:

\begin{theorem}
The moments of the Poisson laws $p_t$, normal laws $g_t$ and complex normal laws $G_t$ are given by the same formula, namely
$$M_k=\sum_{\pi\in D(k)}t^{|\pi|}$$
with $|.|$ being the number of blocks, which at $t=1$ simplifies into
$$M_k=|D(k)|$$
with $D$ being respectively the partitions $P$, the pairings $P_2$, and the matching pairings $\mathcal P_2$.
\end{theorem}

\begin{proof}
This follows indeed by putting togeter the results from chapter 6 regarding the normal laws $g_t,G_t$, and the results here regarding the Poisson laws $p_t$.
\end{proof}

We will be back later with some more conceptual explanations for this result.

\section*{11c. Truncated characters}

With the above probabilistic preliminaries done, let us get back now to finite groups, and compute laws of characters. As a first piece of good news, our main result so far, namely Theorem 11.12, reformulates into something very simple, as follows:

\begin{theorem}
For the symmetric group $S_N\subset O_N$ we have
$$\chi\sim p_1$$
in the $N\to\infty$ limit.
\end{theorem}

\begin{proof}
This is indeed a reformulation of Theorem 11.12, which tells us that with $N\to\infty$ we have the following estimate:
$$P(\chi=k)\simeq\frac{1}{ek!}$$

But, according to our definition of the Poisson laws, this tells us precisely that the asymptotic law of the main character $\chi$ is Poisson (1), as stated.
\end{proof}

An interesting question now is that of recovering all the Poisson laws $p_t$, by using group theory. In order to do this, let us formulate the following definition:

\index{truncated character}

\begin{definition}
Given a closed subgroup $G\subset U_N$, the function
$$\chi:G\to\mathbb C\quad,\quad 
\chi_t(g)=\sum_{i=1}^{[tN]}g_{ii}$$
is called main truncated character of $G$, of parameter $t\in(0,1]$.
\end{definition}

As before with the plain characters, there is some general theory behind this definition, and we will discuss this later on, more systematically, in Part IV. 

\bigskip

Getting back now to the symmetric groups, we first have the following result:

\begin{proposition}
For the symmetric group $S_N\subset O_N$ the coordinate functions are
$$g_{ij}=\chi\left(\sigma\in S_N\Big|\sigma(j)=i\right)$$
and in this picture, the truncated characters count the number of partial fixed points
$$\chi_t(\sigma)=\#\left\{ i\in\{1,\ldots,[tN]\}\Big|\sigma(i)=i\right\}$$
with respect to the truncation parameter $t\in(0,1]$.
\end{proposition}

\begin{proof}
All this is clear from definitions, with the formula for the coordinates being clear from the definition of the embedding $S_N\subset O_N$, and with the character formulae following from it, by summing over $i=j$. To be more precise, we have:
\begin{eqnarray*}
\chi_t(\sigma)
&=&\sum_{i=1}^{[tN]}\sigma_{ii}\\
&=&\sum_{i=1}^{[tN]}\delta_{\sigma(i)i}\\
&=&\#\left\{i\in\{1,\ldots,[tN]\}\Big|\sigma(i)=i\right\}
\end{eqnarray*}

Thus, we are led to the conclusions in the statement.
\end{proof}

Regarding now the asymptotic laws of the truncated characters, the result here, generalizing everything that we have so far, is as follows:

\index{symmetric group}
\index{Poisson law}

\begin{theorem}
For the symmetric group $S_N\subset O_N$ we have
$$\chi_t\sim p_t$$
in the $N\to\infty$ limit, for any $t\in(0,1]$.
\end{theorem}

\begin{proof}
We already know from Theorem 11.22 that the result holds at $t=1$. In general, the proof is similar, the idea being as follows:

\medskip

(1) Consider indeed the following sets, as in the proof of Theorem 11.22, or rather as in the proof of Theorem 11.9, leaading to Theorem 11.22:
$$S_N^i=\left\{\sigma\in S_N\Big|\sigma(i)=i\right\}$$

The set of permutations having no fixed points among $1,\ldots,[tN]$ is then:
$$X_N=\left(\bigcup_{i\leq[tN]}S_N^i\right)^c$$

In order to compute now the cardinality $|X_N|$, consider as well the following sets, depending on indices $i_1<\ldots<i_k$, obtained by taking intersections:
$$S_N^{i_1\ldots i_k}=S_N^{i_1}\cap\ldots\cap S_N^{i_k}$$

As before in the proof of Theorem 11.9, we obtain by inclusion-exclusion that:
\begin{eqnarray*}
P(\chi_t=0)
&=&\frac{1}{N!}\sum_{k=0}^{[tN]}(-1)^k\sum_{i_1<\ldots<i_k<[tN]}|S_N^{i_1\ldots i_k}|\\
&=&\frac{1}{N!}\sum_{k=0}^{[tN]}(-1)^k\sum_{i_1<\ldots<i_k<[tN]}(N-k)!\\
&=&\frac{1}{N!}\sum_{k=0}^{[tN]}(-1)^k\binom{[tN]}{k}(N-k)!\\
&=&\sum_{k=0}^{[tN]}\frac{(-1)^k}{k!}\cdot\frac{[tN]!(N-k)!}{N!([tN]-k)!}
\end{eqnarray*}

With $N\to\infty$, we obtain from this the following estimate:
\begin{eqnarray*}
P(\chi_t=0)
&\simeq&\sum_{k=0}^{[tN]}\frac{(-1)^k}{k!}\cdot t^k\\
&=&\sum_{k=0}^{[tN]}\frac{(-t)^k}{k!}\\
&\simeq&e^{-t}
\end{eqnarray*}

(2) More generally now, by counting the permutations $\sigma\in S_N$ having exactly $k$ fixed points among $1,\ldots,[tN]$, as in the proof of Theorem 11.11, our claim is that we get:
$$P(\chi_t=k)\simeq\frac{t^k}{k!e^t}$$

We already know from (1) that this formula holds at $k=0$. In the general case now, we have to count the permutations $\sigma\in S_N$ having exactly $k$ fixed points among $1,\ldots,[tN]$. Since having such a permutation amounts in choosing $k$ points among $1,\ldots,[tN]$, and then permuting the $N-k$ points left, without fixed points among $1,\ldots,[tN]$ allowed, we obtain the following formula, where $s\in(0,1]$ is such that $[s(N-k)]=[tN]-k$:
\begin{eqnarray*}
\#\left\{\sigma\in S_N\Big|\chi_t(\sigma)=k\right\}
&=&\binom{[tN]}{k}\#\left\{\sigma\in S_{N-k}\Big|\chi_s(\sigma)=0\right\}\\
&=&\frac{[tN]!}{k!([tN]-k)!}\#\left\{\sigma\in S_{N-k}\Big|\chi_s(\sigma)=0\right\}\\
&=&\frac{1}{k!}\times\frac{[tN]!(N-k)!}{([tN]-k)!}\times\frac{\#\left\{\sigma\in S_{N-k}\Big|\chi_s(\sigma)=0\right\}}{(N-k)!}
\end{eqnarray*}

Now by dividing everything by $N!$, we obtain from this the following formula:
$$\frac{\#\left\{\sigma\in S_N\Big|\chi_t(\sigma)=k\right\}}{N!}=\frac{1}{k!}\times\frac{[tN]!(N-k)!}{N!([tN]-k)!}\times\frac{\#\left\{\sigma\in S_{N-k}\Big|\chi_s(\sigma)=0\right\}}{(N-k)!}$$

By using now the computation at $k=0$, that we already have, from (1) above, it follows that with $N\to\infty$ we have the following estimate:
\begin{eqnarray*}
P(\chi_t=k)
&\simeq&\frac{1}{k!}\times\frac{[tN]!(N-k)!}{N!([tN]-k)!}\cdot P(\chi_s=0)\\
&\simeq&\frac{t^k}{k!}\cdot P(\chi_s=0)\\
&\simeq&\frac{t^k}{k!}\cdot\frac{1}{e^s}
\end{eqnarray*}

Now recall that the parameter $s\in(0,1]$ was chosen in the above such that:
$$[s(N-k)]=[tN]-k$$

Thus in the $N\to\infty$ limit we have $s=t$, and so we obtain, as claimed:
$$P(\chi_t=k)\simeq\frac{t^k}{k!}\cdot\frac{1}{e^t}$$

It follows that we obtain in the limit a Poisson law of parameter $t$, as stated.
\end{proof}

\section*{11d. Further results}

All the above is quite interesting, and is at the core of the theory that we want to develop, so let us further build on all this, with a number of more specialized results on the subject, which will be sometimes research-grade. We will be following \cite{bco}.

\bigskip

To start with, let us first present a new, instructive proof for the above character results. The point indeed is that we can approach the problems as well directly, by integrating over $S_N$, and in order to do so, we can use the following result:

\index{polynomial integrals}

\begin{theorem}
Consider the symmetric group $S_N$, with its standard coordinates: 
$$g_{ij}=\chi\left(\sigma\in S_N\Big|\sigma(j)=i\right)$$
The products of these coordinates span the algebra $C(S_N)$, and the arbitrary integrals over $S_N$ are given, modulo linearity, by the formula
$$\int_{S_N}g_{i_1j_1}\ldots g_{i_kj_k}=\begin{cases}
\frac{(N-|\ker i|)!}{N!}&{\rm if}\ \ker i=\ker j\\
0&{\rm otherwise}
\end{cases}$$
where $\ker i$ denotes as usual the partition of $\{1,\ldots,k\}$ whose blocks collect the equal indices of $i$, and where $|.|$ denotes the number of blocks.
\end{theorem}

\begin{proof}
The first assertion follows from the Stone-Weierstrass theorem, because the standard coordinates $g_{ij}$ separate the points of $S_N$, and so the algebra $<g_{ij}>$ that they generate must be equal to the whole function algebra $C(S_N)$:
$$<g_{ij}>=C(S_N)$$

Regarding now the second assertion, according to the definition of the matrix coordinates $g_{ij}$, the integrals in the statement are given by:
$$\int_{S_N}g_{i_1j_1}\ldots g_{i_kj_k}=\frac{1}{N!}\#\left\{\sigma\in S_N\Big|\sigma(j_1)=i_1,\ldots,\sigma(j_k)=i_k\right\}$$

Now observe that the existence of $\sigma\in S_N$ as above requires:
$$i_m=i_n\iff j_m=j_n$$

Thus, the above integral vanishes when the following condition is satisfied:
$$\ker i\neq\ker j$$

Regarding now the case $\ker i=\ker j$, if we denote by $b\in\{1,\ldots,k\}$ the number of blocks of this partition $\ker i=\ker j$, we have $N-b$ points to be sent bijectively to $N-b$ points, and so $(N-b)!$ solutions, and the integral is $\frac{(N-b)!}{N!}$, as claimed.
\end{proof}

As an illustration for the above formula, we can recover the computation of the asymptotic laws of the truncated characters $\chi_t$. We have indeed:

\begin{theorem}
For the symmetric group $S_N\subset O_N$, regarded as a compact group of matrices, $S_N\subset O_N$, via the standard permutation matrices, the truncated character
$$\chi_t(g)=\sum_{i=1}^{[tN]}g_{ii}$$
counts the number of fixed points among $\{1,\ldots,[tN]\}$, and its law with respect to the counting measure becomes, with $N\to\infty$, a Poisson law of parameter $t$. 
\end{theorem}

\begin{proof}
The first assertion comes from the following formula:
$$g_{ij}=\chi\left(\sigma\Big|\sigma(j)=i\right)$$

Regarding now the second assertion, we can use here the integration formula in Theorem 11.26. With $S_{kb}$ being the Stirling numbers, counting the partitions of $\{1,\ldots,k\}$ having exactly $b$ blocks, we have indeed the following formula:
\begin{eqnarray*}
\int_{S_N}\chi_t^k
&=&\sum_{i_1,\ldots,i_k=1}^{[tN]}\int_{S_N}g_{i_1i_1}\ldots g_{i_ki_k}\\
&=&\sum_{\pi\in P(k)}\frac{[tN]!}{([tN]-|\pi|!)}\cdot\frac{(N-|\pi|!)}{N!}\\
&=&\sum_{b=1}^{[tN]}\frac{[tN]!}{([tN]-b)!}\cdot\frac{(N-b)!}{N!}\cdot S_{kb}
\end{eqnarray*}

In particular with $N\to\infty$ we obtain the following formula:
$$\lim_{N\to\infty}\int_{S_N}\chi_t^k=\sum_{b=1}^kS_{kb}t^b$$

But this is the $k$-th moment of the Poisson law $p_t$, and so we are done.
\end{proof}

Summarizing, we have a good understanding of our main result so far, involving the characters of the symmetric group $S_N$ and the Poisson laws of parameter $t\in(0,1]$, by using 2 different methods. We will see in a moment a third proof as well, and we will be actually back to this in Part IV too, with a fourth method too.

\bigskip

As another result now regarding $S_N$, here is a useful related formula:

\begin{theorem}
We have the law formula
$${\rm law}(g_{11}+\ldots +g_{ss})=\frac{s!}{N!}\sum_{p=0}^s\frac{(N-p)!}{(s-p)!}
\cdot\frac{\left(\delta_1-\delta_0\right)^{*p}}{p!}$$ 
where $g_{ij}$ are the standard coordinates of $S_N\subset O_N$.
\end{theorem}

\begin{proof}
We have the following moment formula, where $m_f$ is the number of permutations of $\{1,\ldots ,N\}$ having exactly $f$ fixed points in the set $\{1,\ldots ,s\}$: 
$$\int_{S_N}(u_{11}+\ldots +u_{ss})^k=\frac{1}{N!}\sum_{f=0}^sm_ff^k$$

Thus the law in the statement, say $\nu_{sN}$, is the following average
of Dirac masses:
$$\nu_{sN}=\frac{1}{N!}\sum_{f=0}^s m_{f}\,\delta_f$$

Now observe that the permutations contributing to $m_f$ are obtained by choosing $f$
points in the set $\{1,\ldots ,s\}$, then by permuting the remaining $N-f$ points in $\{1,\ldots ,n\}$ in such a way that there is no fixed point in $\{1,\ldots,s\}$. But these latter permutations are counted as follows: we start with all permutations, we substract those having one fixed point, we add those having two fixed points, and so on. We obtain in this way:
\begin{eqnarray*}
\nu_{sN}
&=&\frac{1}{N!}\sum_{f=0}^s\begin{pmatrix}s\\
f\end{pmatrix}\left(\sum_{k=0}^{s-f}(-1)^k
\begin{pmatrix}s-f\\ k\end{pmatrix}(N-f-k)!\right)\,\delta_f\\
&=&\sum_{f=0}^s\sum_{k=0}^{s-f}(-1)^k\frac{1}{N!}\cdot
\frac{s!}{f!(s-f)!}\cdot\frac{(s-f)!(N-f-k)!}{k!(s-f-k)!}\,\delta_f\\
&=&\frac{s!}{N!}\sum_{f=0}^s\sum_{k=0}^{s-f}\frac{(-1)^k(N-f-k)!}{f!k!(s-f-k)!}\,\delta_f
\end{eqnarray*}

We can proceed as follows, by using the new index $p=f+k$:
\begin{eqnarray*}
\nu_{sN}
&=&\frac{s!}{N!}\sum_{p=0}^s\sum_{k=0}^{p}\frac{(-1)^k
(N-p)!}{(p-k)!k!(s-p)!}\,\delta_{p-k}\\
&=&\frac{s!}{N!}\sum_{p=0}^s\frac{(N-p)!}{(s-p)!p!}
\sum_{k=0}^{p}(-1)^k\begin{pmatrix}p\\
k\end{pmatrix}\,\delta_{p-k}\\
&=&\frac{s!}{N!}\sum_{p=0}^s\frac{(N-p)!}{(s-p)!}\cdot
\frac{\left(\delta_1-\delta_0\right)^{*p}}{p!}
\end{eqnarray*}

Here $*$ is convolution of real measures, and the assertion follows.
\end{proof}

Observe that the above formula is finer than most of our previous formulae regarding truncated characters, which were asymptotic, because it is valid at any $N\in\mathbb N$. 

\bigskip

We can use the above formula as follows, in order to get yet another proof of our main result so far, regarding the Poisson laws, along with a bit more:

\begin{theorem}
Let $g_{ij}$ be the standard coordinates of $C(S_N)$.
\begin{enumerate}
\item $u_{11}+\ldots +u_{ss}$ with $s=o(N)$ is a projection of trace $s/N$. 

\item $u_{11}+\ldots +u_{ss}$ with $s=tN+o(N)$ is Poisson of parameter $t$.
\end{enumerate}
\end{theorem}

\begin{proof}
We can use indeed the formula in Theorem 11.28, as follows:

\medskip

(1) With $s$ fixed and $N\to\infty$ we have the following estimate:
\begin{eqnarray*}
&&{\rm law}(u_{11}+\ldots +u_{ss})\\
&=&\sum_{p=0}^s\frac{(N-p)!}{N!}\cdot\frac{s!}{(s-p)!}
\cdot\frac{\left(\delta_1-\delta_0\right)^{*p}}{p!}\\
&=&\delta_0+\frac{s}{N}\,(\delta_1-\delta_0)+O(N^{-2})
\end{eqnarray*}

But the law on the right is that of a projection of trace $s/N$, as desired.

\medskip

(2) We have a law formula of the following type:
$${\rm law}(u_{11}+\ldots +u_{ss})=
\sum_{p=0}^sc_p\cdot\frac{(\delta_1-\delta_0)^{*p}}{p!}$$

The coefficients $c_p$ can be estimated by using the Stirling formula, as follows:
\begin{eqnarray*}
c_p
&=&\frac{(tN)!}{N!}\cdot\frac{(N-p)!}{(tN-p)!}\\
&\simeq&\frac{(tN)^{tN}}{N^N}\cdot\frac{(N-p)^{N-p}}{(tN-p)^{tN-p}}\\
&=&\left(\frac{tN}{tN-p}\right)^{tN-p} \left(
\frac{N-p}{N}\right)^{N-p}\left( \frac{tN}{N}\right)^p
\end{eqnarray*}

But the last expression can be estimated by using the definition of the exponentials, and we obtain in this way the following estimate:
$$c_p
\simeq e^pe^{-p}t^p
=t^p$$

We can now compute the Fourier transform with respect to a variable $y$:
\begin{eqnarray*}
{\mathcal F}\left( {\rm law}(u_{11}+\ldots +u_{ss})\right)
&\simeq&\sum_{p=0}^st^p\cdot\frac{(e^y-1)^p}{p!}\\
&=&e^{t(e^y-1)}
\end{eqnarray*}

But this is precisely the Fourier transform of the Poisson law $p_t$, as computed in Theorem 11.17, and this gives the second assertion.
\end{proof}

Let us discuss now, as an instructive variation of the above, the computation for the alternating group $A_N\subset S_N$. We will see that with $N\to\infty$ nothing changes, and with this being part of a more general phenomenon, regarding more general types of reflection groups and subgroups, that we will further discuss in the next chapter.

\bigskip

Let us start with some algebraic considerations. We first have:

\begin{proposition}
For the symmetric group, regarded as group of permutations of the $N$ coordinate axes of $\mathbb R^N$, and so as group of permutation matrices, 
$$S_N\subset O_N$$
the determinant is the signature. The subgroup $A_N\subset S_N$ given by 
$$A_N=S_N\cap SO_N$$
and called alternating group, consists of the even permutations.
\end{proposition}

\begin{proof}
In this statement the first assertion is clear from the definition of the determinant, and of the permutation matrices, and all the rest is standard.
\end{proof}

Regarding now character computations, the best here is to use an analogue of Theorem 11.26. To be more precise, we have here the following result:

\begin{theorem}
Consider the alternating group $A_N$, regarded as group of permutation matrices, with its standard coordinates: 
$$g_{ij}=\chi\left(\sigma\in A_N\Big|\sigma(j)=i\right)$$
The products of these coordinates span the algebra $C(A_N)$, and the arbitrary integrals over $A_N$ are given, modulo linearity, by the formula
$$\int_{A_N}g_{i_1j_1}\ldots g_{i_kj_k}\simeq\begin{cases}
\frac{(N-|\ker i|)!}{N!}&{\rm if}\ \ker i=\ker j\\
0&{\rm otherwise}
\end{cases}$$
with $N\to\infty$, where $\ker i$ denotes as usual the partition of $\{1,\ldots,k\}$ whose blocks collect the equal indices of $i$, and where $|.|$ denotes the number of blocks.
\end{theorem}

\begin{proof}
The first assertion follows from the Stone-Weierstrass theorem, because the standard coordinates $g_{ij}$ separate the points of $A_N$, and so we have:
$$<g_{ij}>=C(A_N)$$

Regarding now the second assertion, according to the definition of the standard coordinates $g_{ij}$, the integrals in the statement are given by:
$$\int_{A_N}g_{i_1j_1}\ldots g_{i_kj_k}=\frac{1}{N!/2}\#\left\{\sigma\in A_N\Big|\sigma(j_1)=i_1,\ldots,\sigma(j_k)=i_k\right\}$$

Now observe that the existence of $\sigma\in A_N$ as above requires:
$$i_m=i_n\iff j_m=j_n$$

Thus, the above integral vanishes when the following holds:
$$\ker i\neq\ker j$$

Regarding now the case $\ker i=\ker j$, if we denote by $b\in\{1,\ldots,k\}$ the number of blocks of this partition $\ker i=\ker j$, we have $N-b$ points to be sent bijectively to $N-b$ points. But when assuming $N>>0$, and more specifically $N>k$, half of these bijections will be alternating, and so we have  $(N-b)!/2$ solutions. Thus, the integral is:
\begin{eqnarray*}
\int_{A_N}g_{i_1j_1}\ldots g_{i_kj_k}
&=&\frac{1}{N!/2}\#\left\{\sigma\in A_N\Big|\sigma(j_1)=i_1,\ldots,\sigma(j_k)=i_k\right\}\\
&=&\frac{(N-b)!/2}{N!/2}\\
&=&\frac{(N-b)!}{N!}
\end{eqnarray*}

Thus, we are led to the conclusion in the statement.
\end{proof}

As an application of the above formula, we can now compute the asymptotic laws of the truncated characters $\chi_t$, for the alternating group. We have indeed:

\begin{theorem}
For the alternating group $A_N\subset O_N$, regarded as a compact group of matrices, $A_N\subset O_N$, via the standard permutation matrices, the truncated character
$$\chi_t(g)=\sum_{i=1}^{[tN]}g_{ii}$$
counts the number of fixed points among $\{1,\ldots,[tN]\}$, and its law with respect to the counting measure becomes, with $N\to\infty$, a Poisson law of parameter $t$. 
\end{theorem}

\begin{proof}
The first assertion comes from the following formula:
$$g_{ij}=\chi\left(\sigma\Big|\sigma(j)=i\right)$$

Regarding now the second assertion, we can use here the integration formula in Theorem 11.31. With $S_{kb}$ being the Stirling numbers, counting the partitions of $\{1,\ldots,k\}$ having exactly $b$ blocks, we have the following formula:
\begin{eqnarray*}
\int_{A_N}\chi_t^k
&=&\sum_{i_1\ldots i_k=1}^{[tN]}\int_{A_N}g_{i_1i_1}\ldots g_{i_ki_k}\\
&\simeq&\sum_{\pi\in P(k)}\frac{[tN]!}{([tN]-|\pi|!)}\cdot\frac{(N-|\pi|!)}{N!}\\
&=&\sum_{b=1}^{[tN]}\frac{[tN]!}{([tN]-b)!}\cdot\frac{(N-b)!}{N!}\cdot S_{kb}
\end{eqnarray*}

In particular with $N\to\infty$ we obtain the following formula:
$$\lim_{N\to\infty}\int_{A_N}\chi_t^k=\sum_{b=1}^kS_{kb}t^b$$

But this is the $k$-th moment of the Poisson law $p_t$, and so we are done.
\end{proof}

Summarizing, when passing from the symmetric group $S_N$ to its subgroup $A_N\subset S_N$, in what concerns character computations, with $N\to\infty$ nothing changes. This is actually part of a more general phenomenon, regarding more general types of reflection groups and subgroups, that we will further discuss in the next chapter.

\bigskip

As a conclusion now to all this, we have seen that the truncated characters $\chi_t$ of the symmetric group $S_N$ have the Poisson laws $p_t$ as limiting distributions, with $N\to\infty$. Moreover, we have seen several proofs for this fundamental fact, using inclusion-exclusion, direct integration, and convolution exponentials and Fourier transforms as well.

\bigskip

We will keep building on all this in the next chapter, by stating and proving similar results for more general reflection groups $G\subset U_N$. Also, we will be back to the symmetric group $S_N$ and to the Poisson laws in Part IV, with a fourth proof for our results, using representation theory, and a property of $S_N$ called easiness. More on this later.

\bigskip

Finally, as an important theoretical remark, in relation with all this, recall from the beginning of this chapter that for the cyclic group $\mathbb Z_N\subset O_N$ the computation was not very interesting, leading to a Bernoulli law having trivial asymptotics, while for the dihedral group $D_N\subset O_N$ the law of the main character, not that interesting either, was not even uniform in $N$. You might probably ask then, what is wrong with $\mathbb Z_N$ and $D_N$? In answer, these groups are not ``easy'', and more on easiness, later in this book.

\section*{11e. Exercises}

There are many interesting possible exercises in connection with the above. First, in relation with derangements and fixed points, we have:

\begin{exercise}
Compute the number of derangements in $S_4$, by explicitly listing them, and then comment on the estimate of 
$$e=2.7182\ldots$$
that you obtain in this way.
\end{exercise}

Here the first question is of course elementary, but the problem is that of finding out what the best notation for permutations is, in order to solve this problem quickly. As for the second question, that you can investigate at higher $N$ too, based on the various formulae established in this chapter, this is something quite instructive too.

\begin{exercise}
Show that the probability for a length $1$ needle to intersect, when thrown, a $1$-spaced grid is $2/\pi$, and then comment on the estimate on
$$\pi=3.1415\ldots$$
that you obtain in this way.
\end{exercise}

Here the first question is quite tricky, because there are several possible ways of modelling the problem, but only one of them gives the correct, real-life answer. As for the second question, this is a good introduction to applied mathematics too.

\begin{exercise}
Find some formulae for the Bell numbers $B_k$, or rather for their generating series, or suitable transforms of that series, and the more the better.
\end{exercise}

There is a lot of interesting mathematics here, and after solving the exercise, you can check the internet, and complete your knowledge with more things.

\begin{exercise}
Show that the truncated characters of $S_N$, suitably moved over the diagonal, as to not overlap, become independent with $N\to\infty$.
\end{exercise}

Here the formulation is of course a bit loose, but this is intentional, and finding the precise formulation is part of the exercise. As for the proof, this can only come by using the various integration formulae over $S_N$ established in the above.

\begin{exercise}
Find some alternative proofs for the fact, that we already know, that the truncated charcters for $A_N\subset O_N$ become Poisson, with $N\to\infty$.
\end{exercise}

This is a bit technical, the problem being that of picking the best alternative proof for $S_N$, from the above, and then extending it to $A_N$. As a bonus exercise, you can work out as well independence aspects for $A_N$, in the spirit of the previous exercise.

\chapter{Reflection groups}

\section*{12a. Real reflections}

We have seen in the previous chapter that some interesting phenomena, in relation with the law of the main character, appear for the symmetric group $S_N$, in the $N\to\infty$ limit. All this suggests looking at more general reflection groups. Let us begin by discussing the hyperoctahedral group $H_N$. We recall from chapter 9 that we have:

\index{hyperoctahedral group}

\begin{theorem}
Consider the hyperoctahedral group $H_N$, which appears as the symmetry group of the $N$-cube, or the symmetry group of the $N$ coordinate axes of $\mathbb R^N$:
$$S_N\subset H_N\subset O_N$$
In matrix terms, $H_N$ consists of the permutation-type matrices having $\pm1$ as nonzero entries, and we have a wreath product decomposition as follows:
$$H_N=\mathbb Z_2\wr S_N$$
In this picture, the main character counts the signed number of fixed points, among the coordinate axes, and its truncations count the truncations of such numbers.
\end{theorem}

\begin{proof}
This is something that we discussed before, the idea being that the first assertions are clear, and that the wreath product decomposition in the statement comes from a crossed product decomposition $H_N=\mathbb Z_2^N\rtimes S_N$. As for the assertions regarding the main character and its truncations, once again these are clear, as for $S_N$.
\end{proof}

Regarding now the character laws, we can compute them by using the same method as for the symmetric group $S_N$, namely inclusion-exclusion, and we have:

\index{main character}

\begin{theorem}
For the hyperoctahedral group $H_N\subset O_N$, the law of the variable
$$\chi_t=\sum_{i=1}^{[tN]}g_{ii}$$
becomes with $N\to\infty$ the following measure
$$b_t=e^{-t}\sum_{k=-\infty}^\infty\delta_k\sum_{p=0}^\infty \frac{(t/2)^{|k|+2p}}{(|k|+p)!p!}$$ 
where $\delta_k$ is the Dirac mass at $k\in\mathbb Z$.
\end{theorem}

\begin{proof}
We follow \cite{bbc}. We regard $H_N$ as being the symmetry group of the graph $I_N=\{I^1,\ldots ,I^N\}$ formed by $N$ segments. The diagonal coefficients are given by:
$$u_{ii}(g)=\begin{cases}
\ 0\ \mbox{ if $g$ moves $I^i$}\\
\ 1\ \mbox{ if $g$ fixes $I^i$}\\
-1\mbox{ if $g$ returns $I^i$}
\end{cases}$$

We denote by $\uparrow g,\downarrow g$ the number of segments among $\{I^1,\ldots ,I^s\}$ which are fixed, respectively returned by an element $g\in H_N$. With this notation, we have:
$$u_{11}+\ldots+u_{ss}=\uparrow g-\downarrow g$$

Let us denote by $P_N$ probabilities computed over the group $H_N$. The density of the law of $u_{11}+\ldots+u_{ss}$ at a point $k\geq 0$ is then given by the following formula:
\begin{eqnarray*}
D(k)
&=&P_N(\uparrow g-\downarrow g=k)\\
&=&\sum_{p=0}^\infty P_N(\uparrow g=k+p, \downarrow g=p)
\end{eqnarray*}

Assume first that we have $t=1$. We use the fact, that we know well from chapter 11, that the probability of $\sigma\in S_N$ to have no fixed points is asymptotically given by:
$$P_0=\frac{1}{e}$$

Thus the probability of $\sigma\in S_N$ to have $m$ fixed points is asymptotically given by:
$$P_m=\frac{1}{em!}$$

In terms of probabilities over $H_N$, we obtain from this, as desired:
\begin{eqnarray*}
\lim_{N\to\infty}D(k)
&=&\lim_{N\to\infty}\sum_{p=0}^\infty(1/2)^{k+2p}\begin{pmatrix}k+2p\\ k+p\end{pmatrix} P_N(\uparrow g+\downarrow g=k+2p)\\ 
&=&\sum_{p=0}^\infty(1/2)^{k+2p}\begin{pmatrix}k+2p\\
k+p\end{pmatrix}\frac{1}{e(k+2p)!}\\
&=&\frac{1}{e}\sum_{p=0}^\infty \frac{(1/2)^{k+2p}}{(k+p)!p!}
\end{eqnarray*}

As for the general case $t\in(0,1]$, here the result follows by performing some modifications in the above computation. The asymptotic density is computed as follows:
\begin{eqnarray*}
\lim_{N\to\infty}D(k)
&=&\lim_{N\to\infty}\sum_{p=0}^\infty(1/2)^{k+2p}\begin{pmatrix}k+2p\\ k+p\end{pmatrix} P_N(\uparrow g+\downarrow g=k+2p)\\
&=&\sum_{p=0}^\infty(1/2)^{k+2p}\begin{pmatrix}k+2p\\
k+p\end{pmatrix}\frac{t^{k+2p}}{e^t(k+2p)!}\\
&=&e^{-t}\sum_{p=0}^\infty \frac{(t/2)^{k+2p}}{(k+p)!p!}
\end{eqnarray*}

Together with $D(-k)=D(k)$, this gives the formula in the statement.
\end{proof}

The above result is quite interesting, because the densities there are the Bessel functions of the first kind. Due to this fact, the limiting measures are called Bessel laws:

\index{Bessel function}
\index{Bessel law}

\begin{definition}
The Bessel law of parameter $t>0$ is the measure
$$b_t=e^{-t}\sum_{k=-\infty}^\infty\delta_k\,f_k(t/2)$$
with the density being the following function,
$$f_k(t)=\sum_{p=0}^\infty \frac{t^{|k|+2p}}{(|k|+p)!p!}$$
called Bessel function of the first kind.
\end{definition}

Let us study now these Bessel laws, in analogy with what we know from chapter 11, regarding the Poisson laws. We first have the following result:

\index{convolution}

\begin{theorem}
The Bessel laws $b_t$ have the property
$$b_s*b_t=b_{s+t}$$
so they form a truncated one-parameter semigroup
with respect to convolution.
\end{theorem}

\begin{proof}
Again, we follow \cite{bbc}. We use the formula in Definition 12.3, namely:
$$b_t=e^{-t}\sum_{k=-\infty}^\infty\delta_k\,f_k(t/2)$$

The Fourier transform of this measure is given by the following formula:
$$Fb_t(y)=e^{-t}\sum_{k=-\infty}^\infty e^{ky}\,f_k(t/2)$$

We compute now the derivative with respect to $t$:
$$Fb_t(y)'=-Fb_t(y)+\frac{e^{-t}}{2}\sum_{k=-\infty}^\infty e^{ky}\,f_k'(t/2)$$

On the other hand, the derivative of $f_k$ with $k\geq 1$ is given by:
\begin{eqnarray*}
f_k'(t)
&=&\sum_{p=0}^\infty \frac{(k+2p)t^{k+2p-1}}{(k+p)!p!}\\
&=&\sum_{p=0}^\infty \frac{(k+p)t^{k+2p-1}}{(k+p)!p!}+\sum_{p=0}^\infty\frac{p\,t^{k+2p-1}}{(k+p)!p!}\\
&=&\sum_{p=0}^\infty \frac{t^{k+2p-1}}{(k+p-1)!p!}+\sum_{p=1}^\infty\frac{t^{k+2p-1}}{(k+p)!(p-1)!}\\
&=&\sum_{p=0}^\infty \frac{t^{(k-1)+2p}}{((k-1)+p)!p!}+\sum_{p=1}^\infty\frac{t^{(k+1)+2(p-1)}}{((k+1)+(p-1))!(p-1)!}\\
&=&f_{k-1}(t)+f_{k+1}(t)
\end{eqnarray*}

This computation works in fact for any $k$, so we get:
\begin{eqnarray*}
Fb_t(y)'
&=&-Fb_t(y)+\frac{e^{-t}}{2}
\sum_{k=-\infty}^\infty e^{ky} (f_{k-1}(t/2)+f_{k+1}(t/2))\\
&=&-Fb_t(y)+\frac{e^{-t}}{2} \sum_{k=-\infty}^\infty
e^{(k+1)y}f_{k}(t/2)+e^{(k-1)y}f_{k}(t/2)\\
&=&-Fb_t(y)+\frac{e^{y}+e^{-y}}{2}\,Fb_t(y)\\
&=&\left(\frac{e^{y}+e^{-y}}{2}-1\right)Fb_t(y)
\end{eqnarray*}

Thus the log of the Fourier transform is linear in $t$, and we get
the assertion.
\end{proof}

In order to further discuss all this, we will need a number of probabilistic preliminaries. We recall that, conceptually speaking, the Poisson laws are the laws appearing via the Poisson Limit Theorem (PLT), stating that we have the following convergence:
$$\left(\left(1-\frac{t}{n}\right)\delta_0+\frac{t}{n}\delta_1\right)^{*n}\to p_t$$ 

In order to generalize this construction, as to cover the Bessel laws found above, in connection with the hyperoctahedral group $H_N$, we have the following notion:

\index{compound Poisson law}

\begin{definition}
Associated to any compactly supported positive measure $\nu$ on $\mathbb C$ is the probability measure
$$p_\nu=\lim_{n\to\infty}\left(\left(1-\frac{c}{n}\right)\delta_0+\frac{1}{n}\nu\right)^{*n}$$
where $c=mass(\nu)$, called compound Poisson law.
\end{definition}

In other words, what we are doing here is to generalize the construction in the Poisson Limit Theorem, by allowing the only parameter there, which was the positive real number $t>0$, to be replaced by a certain probability measure $\nu$, of arbitrary mass $c>0$.

\bigskip

In what follows we will be mainly interested in the case where $\nu$ is discrete, as is for instance the measure $\nu=t\delta_1$ with $t>0$, which produces via the above limiting procedure the Poisson laws. In fact, we will be mainly interested in the case where $\nu$ is a multiple of the uniform measure on the $s$-th roots of unity, and more on this later.

\bigskip

The following result allows us to detect compound Poisson laws:

\index{Fourier transform}

\begin{proposition}
For a discrete measure, $\nu=\sum_{i=1}^sc_i\delta_{z_i}$ with $c_i>0$ and $z_i\in\mathbb C$, we have the formula
$$F_{p_\nu}(y)=\exp\left(\sum_{i=1}^sc_i(e^{iyz_i}-1)\right)$$
where $F$ denotes as usual the Fourier transform.
\end{proposition}

\begin{proof}
Let $\mu_n$ be the measure appearing in Definition 12.5, namely:
$$\mu_n=\left(1-\frac{c}{n}\right)\delta_0+\frac{1}{n}\nu$$

We have the following computation, in the context of Definition 12.5:
\begin{eqnarray*}
&&F_{\mu_n}(y)=\left(1-\frac{c}{n}\right)+\frac{1}{n}\sum_{i=1}^sc_ie^{iyz_i}\\
&\implies&F_{\mu_n^{*n}}(y)=\left(\left(1-\frac{c}{n}\right)+\frac{1}{n}\sum_{i=1}^sc_ie^{iyz_i}\right)^n\\
&\implies&F_{p_\nu}(y)=\exp\left(\sum_{i=1}^sc_i(e^{iyz_i}-1)\right)
\end{eqnarray*}

Thus, we have obtained the formula in the statement.
\end{proof}

We have as well the following result, providing an alternative to Definition 12.5, and which will be our formulation of the Compound Poisson Limit Theorem (CPLT):

\index{compound Poisson Limit theorem}
\index{CPLT}

\begin{theorem}
For a discrete measure, $\nu=\sum_{i=1}^sc_i\delta_{z_i}$ with $c_i>0$ and $z_i\in\mathbb C$, we have the formula
$$p_\nu={\rm law}\left(\sum_{i=1}^sz_i\alpha_i\right)$$
where the variables $\alpha_i$ are Poisson $(c_i)$, independent.
\end{theorem}

\begin{proof}
Let $\alpha$ be the sum of Poisson variables in the statement:
$$\alpha=\sum_{i=1}^sz_i\alpha_i$$

By using some well-known Fourier transform formulae, we have:
\begin{eqnarray*}
F_{\alpha_i}(y)=\exp(c_i(e^{iy}-1))
&\implies&F_{z_i\alpha_i}(y)=\exp(c_i(e^{iyz_i}-1))\\
&\implies&F_\alpha(y)=\exp\left(\sum_{i=1}^sc_i(e^{iyz_i}-1)\right)
\end{eqnarray*}

Thus we have the same formula as in Proposition 12.6, as desired.
\end{proof}

Getting back now to the Bessel laws, we have the following result:

\index{Bessel law}

\begin{theorem}
The Bessel laws $b_t$ are compound Poisson laws, given by
$$b_t=p_{t\varepsilon}$$
where $\varepsilon=\frac{1}{2}(\delta_{-1}+\delta_1)$ is the uniform measure on $\mathbb Z_2$.
\end{theorem}

\begin{proof}
This follows indeed by comparing the formula of the Fourier transform of $b_t$, from the proof of Theorem 12.4, with the formula in Proposition 12.6.
\end{proof}

As a conclusion to this, when discussing the asymptotic character law for the basic finite subgroups $G\subset U_N$, such as $G=S_N,H_N$, it is all about compound Poisson laws.

\section*{12b. Complex reflections}

Our next task will be that of unifying and generalizing the results that we have for $S_N,H_N$. For this purpose, consider the following remarkable family of groups:

\index{complex reflection group}

\begin{definition}
The complex reflection group $H_N^s\subset U_N$, depending on parameters
$$N\in\mathbb N\quad,\quad s\in\mathbb N\cup\{\infty\}$$
is the group of permutation-type matrices with $s$-th roots of unity as entries,
$$H_N^s=M_N(\mathbb Z_s\cup\{0\})\cap U_N$$
with the convention $\mathbb Z_\infty=\mathbb T$, at $s=\infty$.
\end{definition}

This construction is something quite tricky, that will keep us busy, for the remainder of this section. As a first observation, at $s=1,2$ we obtain the following groups:
$$H_N^1=S_N\quad,\quad 
H_N^2=H_N$$

Another important particular case of the above construction is $s=\infty$, where we obtain a group which is actually not finite, but is still compact, denoted as follows:
$$K_N\subset U_N$$

This latter group $K_N$ is called full complex reflection group, and will appear many times, in what follows. In view of this, let us highlight its definition, as follows:

\index{full reflection group}

\begin{definition}
The full complex reflection group is given by:
$$K_N=M_N(\mathbb T\cup\{0\})\cap U_N$$
That is, $K_N$ is the group of permutation-type matrices with entries from $\mathbb T$.
\end{definition}

In fact, we already met $K_N$ at the end of chapter 10, when talking about the reflection subgroup of an arbitrary group $G\subset U_N$, which was constructed as follows:
$$K=G\cap K_N$$

Summarizing, $K_N$ seems to be a quite interesting object, with its precise potential remaining to be determined. So, let us first have a look at it at small values of $N$:

\medskip

\underline{$N=1$}.  What we have is the unit circle, $K_1=\mathbb T$.

\medskip 

\underline{$N=2$}. Here $K_2$ consists of the matrices as follows, with nonzero entries in $\mathbb T$:
$$\begin{pmatrix}x&0\\0&y\end{pmatrix}\qquad,\qquad
\begin{pmatrix}0&x\\y&0\end{pmatrix}$$

\underline{$N=3$}. Here $K_3$ consists of the matrices as follows, with nonzero entries in $\mathbb T$:
$$\begin{pmatrix}x&0&0\\0&y&0\\0&0&z\end{pmatrix}\quad,\quad
\begin{pmatrix}0&x&0\\y&0&0\\0&0&z\end{pmatrix}\quad,\quad
\begin{pmatrix}x&0&0\\0&0&y\\0&z&0\end{pmatrix}$$
$$\begin{pmatrix}0&0&x\\0&y&0\\z&0&0\end{pmatrix}\quad,\quad
\begin{pmatrix}0&0&x\\y&0&0\\0&z&0\end{pmatrix}\quad,\quad
\begin{pmatrix}0&x&0\\0&0&y\\z&0&0\end{pmatrix}$$

\underline{$N\geq4$}. And so on, you get the point, what we have is a bit like before for $H_N$, permutation matrices, but this time decorated by numbers in $\mathbb T$.

\medskip

Generally speaking, $K_N$ contains all the interesting finite groups $G\subset U_N$ that we know, including $S_N,H_N$, and more generally the groups $H_N^s$ from Definition 12.9. Quite remarkably, the dihedral group $D_N$ can be viewed as well as a subgroup, as follows:

\begin{theorem}
We have an embedding $D_N\subset K_2$, coming as follows,
$$D_N=\left\{\begin{pmatrix}x&0\\0&y\end{pmatrix}\ ,\ 
\begin{pmatrix}0&x\\y&0\end{pmatrix}\Big|x=y^{-1}\in\mathbb Z_N\right\}\subset K_2$$
obtained by augmenting the standard copy $\mathbb Z_N\subset K_2$ with a twisted copy of it.
\end{theorem}

\begin{proof}
The matrices patterned as in the statement form indeed a group, and when adding the extra condition $xy=1$, this remains a group. In order now to establish the isomorphism with $D_N$, let us label our group elements as follows, with $xy=zt=1$:
$$R_x=\begin{pmatrix}x&0\\0&y\end{pmatrix}\quad,\quad 
S_z=\begin{pmatrix}0&z\\t&0\end{pmatrix}$$

We have then the following computations, for the products of these elements:
$$R_xR_z=\begin{pmatrix}x&0\\0&y\end{pmatrix}
\begin{pmatrix}z&0\\0&t\end{pmatrix}
=\begin{pmatrix}xz&0\\0&yt\end{pmatrix}
=R_{xz}$$
$$R_xS_z=\begin{pmatrix}x&0\\0&y\end{pmatrix}
\begin{pmatrix}0&z\\t&0\end{pmatrix}
=\begin{pmatrix}0&xz\\yt&0\end{pmatrix}
=S_{xz}$$
$$S_xR_z=
\begin{pmatrix}0&x\\y&0\end{pmatrix}
\begin{pmatrix}z&0\\0&t\end{pmatrix}\
=\begin{pmatrix}0&xt\\yz&0\end{pmatrix}
=S_{xz^{-1}}$$
$$S_xS_z=\begin{pmatrix}0&x\\y&0\end{pmatrix}
\begin{pmatrix}0&z\\t&0\end{pmatrix}
=\begin{pmatrix}xt&0\\0&yz\end{pmatrix}
=R_{xz^{-1}}$$

But, we recognize here the table of multiplication of $D_N$, as desired.
\end{proof}

Summarizing, good idea to pass to complex numbers, and the complex reflection groups $H_N^s\subset U_N$ from Definition 12.9, with special attention to the group $H_N^\infty=K_N$ from Definition 12.10, which contains them all, will be our new objects of interest.

\bigskip

Let us start our study by summarizing some basic observations, as follows:

\begin{proposition}
The complex reflection groups $H_N^s\subset U_N$ are as follows:
\begin{enumerate}
\item At $s=1$ we have $H_N^1=S_N$, having cardinality $|S_N|=N!$.

\item At $s=2$ we have $H_N^2=H_N$, having cardinality $|H_N|=2^NN!$.

\item At $s=\infty$ we have $H_N^\infty=K_N$, having cardinality $|K_N|=\infty$.
\end{enumerate}
\end{proposition}

\begin{proof}
This is clear indeed from the discussion made after Definition 12.9, and with the cardinality results at $s=1$ and $s=2$ being something that we know well.
\end{proof}

Let us record as well the following result, which is something elementary too:

\begin{proposition}
We have inclusions as follows, for any $r,s\in\mathbb N\cup\{\infty\}$:
$$r|s\implies H_r\subset H_s$$
In particular, we have inclusions $S_N\subset H_N^s\subset K_N$, for any $s\in\mathbb N\cup\{\infty\}$.
\end{proposition}

\begin{proof}
With the cyclic group $\mathbb Z_s$ being viewed as usual, as being the group of the $s$-th roots of unity in the complex plane, we have inclusions as follows:
$$r|s\implies \mathbb Z_r\subset \mathbb Z_s$$

Thus, with the group $H_N^s$ constructed as in Definition 12.9, for $r|s$ we have:
\begin{eqnarray*}
H_N^r
&=&M_N(\mathbb Z_r\cup\{0\})\cap U_N\\
&\subset&M_N(\mathbb Z_s\cup\{0\})\cap U_N\\
&=&H_N^s
\end{eqnarray*} 

Finally, the last assertion is clear, and comes also from this, via $1|s|\infty$, for any $s$.
\end{proof}

Coming next, in analogy with what we know about $S_N,H_N$, we first have:

\begin{proposition}
The number of elements of $H_N^s$ with $s\in\mathbb N$ is:
$$|H_N^s|=s^NN!$$
At $s=\infty$, the group $K_N=H_N^\infty$ that we obtain is infinite.
\end{proposition}

\begin{proof}
This is indeed clear from our definition of $H_N^s$, as a matrix group, because there are $N!$ choices for a permutation-type matrix, and then $s^N$ choices for the corresponding $s$-roots of unity, which must decorate the $N$ nonzero entries.
\end{proof}

Once again in analogy with what we know at $s=1,2$, we have as well:

\index{wreath product}

\begin{theorem}
We have a wreath product decomposition
$$H_N^s=\mathbb Z_s^N\rtimes S_N=\mathbb Z_s\wr S_N$$
with the permutations $\sigma\in S_N$ acting on the elements $e\in\mathbb Z_s^N$ as follows:
$$\sigma(e_1,\ldots,e_N)=(e_{\sigma(1)},\ldots,e_{\sigma(N)})$$
In particular we have, as found before, the cardinality formula $|H_N^s|=s^NN!$.
\end{theorem}

\begin{proof}
As explained in the proof of Proposition 12.14, the elements of $H_N^s$ can be identified with the pairs $g=(e,\sigma)$ consisting of a permutation $\sigma\in S_N$, and a decorating vector $e\in\mathbb Z_s^N$, so that at the level of the cardinalities, we have:
$$|H_N|=|\mathbb Z_s^N\times S_N|$$

Now observe that the product formula for two such pairs $g=(e,\sigma)$ is as follows, with the permutations $\sigma\in S_N$ acting on the elements $f\in\mathbb Z_s^N$ as in the statement:
$$(e,\sigma)(f,\tau)=(ef^\sigma,\sigma\tau)$$

Thus, we are in the framework of the crossed products, and we obtain $H_N^s=\mathbb Z_s^N\rtimes S_N$. But this can be written, by definition, as $H_N^s=\mathbb Z_s\wr S_N$, and we are done.
\end{proof}

Finally, in relation with geometric aspects, the above groups appear as follows:

\begin{theorem}
The complex reflection group $H_N^s$ appears as a symmetry group,
$$H_N^s=G(C_s\ldots C_s)$$
with $C_s\ldots C_s$ consisting of $N$ disjoint copies of the oriented cycle $C_s$.
\end{theorem}

\begin{proof}
This is something elementary, the idea being as follows:

\medskip

(1) Consider first the oriented cycle $C_s$, which looks as follows:
$$\xymatrix@R=15pt@C=16pt{
&\bullet\ar[r]&\bullet\ar[dr]\\
\bullet\ar[ur]&&&\bullet\ar[d]\\
\bullet\ar[u]&&&\bullet\ar[dl]\\
&\bullet\ar[ul]&\bullet\ar[l]}$$

It is then clear that the symmetry group of this graph is the cyclic group $\mathbb Z_s$.

\medskip

(2) In the general case now, where we have $N\in\mathbb N$ disjoint copies of the above cycle $C_s$, we must suitably combine the corresponding $N$ copies of the cyclic group $\mathbb Z_s$. But this leads to the wreath product group $H_N^s=\mathbb Z_s\wr S_N$, as stated.
\end{proof}

Moving on, the story with the complex reflection groups is not over with the groups $H_N^s$ constructed in Definition 12.9, because we can do more generally, as follows:

\begin{theorem}
We have subgroups of the basic complex reflection groups,
$$H_N^{sd}=\left\{U\in H_N^s\Big|(\square\, U)^d=1\right\}$$
with $\square$ being the product of nonzero entries, covering all examples of reflection groups.
\end{theorem}

\begin{proof}
This is something very standard, the idea as follows:

\medskip

(1) To start with, with $\square$ being as above, we have a group morphism as follows:
$$\square:H_N^s\to\mathbb Z_s$$

Thus, for any $d|s$, we can define a subgroup $H_N^{sd}\subset H_N^s$ as in the statement.

\medskip

(2) At the level of basic examples now, we certainly have the groups $H_N^s=H_N^{ss}$. Also, recall from Theorem 12.11 that we have an identification as follows: 
$$D_N=\left\{\begin{pmatrix}x&0\\0&y\end{pmatrix}\ ,\ 
\begin{pmatrix}0&x\\y&0\end{pmatrix}\Big|x=y^{-1}\in\mathbb Z_N\right\}\subset K_2$$

But this translates into $D_N=H_2^{N1}$, so the dihedral group $D_N$ is covered too.
\end{proof}

As a conclusion to all this, good work that we did, and we will stop here with our construction of complex reflection groups, due to a famous classification result of Shephard and Todd, that we would like to explain now. To start with, we can talk about complex reflections and about complex reflection groups abstractly, as follows:

\begin{definition}
We can talk about reflections and reflection groups, as follows:
\begin{enumerate}
\item A reflection is a symmetry $S\in U_N$ with respect to a hyperplane $P\subset\mathbb C^N$.

\item A reflection group is a group $G\subset U_N$ generated by reflections, $G=<S_i>$.

\item Such a reflection group is called irreducible when it has no invariant subspaces.
\end{enumerate}
\end{definition}

Observe that we have not assumed $G$ to be finite, in the above, and with this making the above formalism quite broad, for instance with many continuous groups $G\subset U_N$ being reflection groups, in the above sense. Still in this setting, with no finiteness assumption on $G$, these reflection groups are best investigated by writing them as follows:
$$G=\left<S_1,\ldots,S_n\Big|(S_iS_j)^{m_{ij}}=1\right>$$ 

And there has been a lot of work here, by Coxeter and others. Getting now to the finite group case, any reflection group appears as product of irreducible reflection groups, and in what regards these latter groups, we have the following classification result:

\index{complex reflection group}
\index{Shephard-Todd}

\begin{theorem}
The irreducible complex reflection groups are
$$H_N^{sd}=\left\{U\in H_N^s\Big|(\square\, U)^d=1\right\}$$
along with $34$ exceptional examples.
\end{theorem}

\begin{proof}
This is something quite advanced, that we will not attempt to prove here, or even explain in detail, with the list of 34 exceptional cases, and we refer here to the paper of Shephard and Todd \cite{sto}, and to the subsequent literature on the subject. 
\end{proof}

\section*{12c. Bessel laws}

Back now to probability, in order to do the character computations for $H_N^s$, and why not for $H_N^{sd}$ too, we will need a number of further preliminaries. Let us start with:

\index{generalized Bessel laws}
\index{complex Bessel laws}

\begin{definition}
The Bessel law of level $s\in\mathbb N\cup\{\infty\}$ and parameter $t>0$ is
$$b_t^s=p_{t\varepsilon_s}$$
with $\varepsilon_s$ being the uniform measure on the $s$-th roots of unity.
\end{definition}

Observe that at $s=1,2$ we obtain the Poisson and real Bessel laws:
$$b^1_t=p_t\quad,\quad 
b^2_t=b_t$$

Another important particular case is $s=\infty$, where we obtain a measure which is actually not discrete, that we will denote as follows:
$$b^\infty_t=B_t$$

As a basic result on these laws, generalizing those before about $p_t,b_t$, we have:

\begin{theorem}
The generalized Bessel laws $b^s_t$ have the property
$$b^s_t*b^s_{t'}=b^s_{t+t'}$$
so they form a truncated one-parameter semigroup
with respect to convolution.
\end{theorem}

\begin{proof}
This follows indeed from the Fourier transform formula from Proposition 12.6, because for the Bessel laws, the log of this Fourier transform is linear in $t$.
\end{proof}

Regarding now the moments, the result here is as follows:

\begin{theorem}
The moments of the Bessel law $b^s_t$ are the numbers
$$M_k=|P^s(k)|$$
where $P^s(k)$ is the set of partitions of $\{1,\ldots,k\}$ satisfying
$$\#\circ=\#\bullet(s)$$
as a weighted sum, in each block.
\end{theorem}

\begin{proof}
This is something more technical, the idea being as follows:

\medskip

(1) We know that the formula in the statement holds at $s=1$, where $b^1_t=p_t$ is the Poisson law of parameter $t>0$, and $P^1=P$ is the set of all partitions. 

\medskip

(2) The formula in the statement holds also at $s=2$, where $b^2_t=b_t$ is the real Bessel law of parameter $t>0$, and $P^2=P_{even}$ is the set of partitions with even blocks. 

\medskip

(3) Next, at $s=\infty$ the measure in the statement is the complex Bessel law $b^\infty_t=B_t$, the set of partitions is $P^\infty=\mathcal P_{even}$, and the result can be proved, in a similar way.

\medskip

(4) Finally, with the cases $s=1,2,\infty$ understood, the generalization to the case $s\in\mathbb N\cup\{\infty\}$ is quite straightforward, by doing some combinatorics. See \cite{bb+}.
\end{proof}

Getting back now to the reflection groups, we have the following result:

\begin{theorem}
For the complex reflection group $H_N^s=\mathbb Z_s\wr S_N$ we have
$$\chi_t\sim b^s_t$$
with $N\to\infty$, where $b_t^s=p_{t\varepsilon_s}$ is the Bessel law constructed above.
\end{theorem}

\begin{proof}
The best here is to proceed in two steps, as follows:

\medskip

(1) Let us first work out the case $t=1$. Since the limit probability for a random permutation to have exactly $k$ fixed points is $e^{-1}/k!$, we get:
$$\lim_{N\to\infty}law(\chi_1)=e^{-1}\sum_{k=0}^\infty \frac{1}{k!}\,\varepsilon_s^{*k}$$

On the other hand, we get from the definition of the Bessel law $b^s_1$:
\begin{eqnarray*}
b^s_1
&=&\lim_{N\to\infty}\left(\left(1-\frac{1}{N}\right)\delta_0+\frac{1}{N}\,\varepsilon_s\right)^{*N}\\
&=&\lim_{N\to\infty}\sum_{k=0}^N\begin{pmatrix}N\\ k\end{pmatrix}\left(1-\frac{1}{N}\right)^{N-k}\frac{1}{N^k}\,\varepsilon_s^{*k}\\
&=&e^{-1}\sum_{k=0}^\infty\frac{1}{k!}\,\varepsilon_s^{*k}
\end{eqnarray*}

But this gives the assertion for $t=1$, as desired.

\medskip

(2) Now in the case where $t>0$ is arbitrary, we can use the same method, by performing the following modifications to the above computation:
\begin{eqnarray*}
\lim_{N\to\infty}law(\chi_t)
&=&e^{-t}\sum_{k=0}^\infty\frac{t^k}{k!}\,\varepsilon_s^{*k}\\
&=&\lim_{N\to\infty}\left(\left(1-\frac{1}{N}\right)\delta_0+\frac{1}{N}\,\varepsilon_s\right)^{*[tN]}\\
&=&b^s_t
\end{eqnarray*}

Thus, we are led to the conclusion in the statement.
\end{proof}

Let us develop now some more theory for the Bessel laws, following \cite{bb+}. According to our various results above, these Bessel laws appear in practice as follows:

\begin{theorem}
The Bessel laws are given by the formula
$$b^s_t={\rm law}\left(\sum_{k=1}^sw^ka_k\right)$$
with $a_1,\ldots,a_s$ being Poisson $(t/s)$ and independent, and $w=e^{2\pi i/s}$.
\end{theorem}

\begin{proof}
This comes indeed from our general formula from Theorem 12.7.
\end{proof}

We will need in our computations the level $s$ exponential function, given by:
$$\exp_sz=\sum_{k=0}^\infty\frac{z^{sk}}{(sk)!}=\frac{1}{s}\sum_{k=1}^s\exp(w^kz)$$

Observe also that at $s=1,2$ we have the following formulae:
$$\exp_1=\exp\quad,\quad 
\exp_2=\cosh$$

We have the following result, regarding the Fourier transform of the Bessel laws:

\begin{theorem}
The Fourier transform of $b^s_t$ is given by
$$\log F^s_t(z)=t\left(\exp_sz-1\right)$$
so in particular the measures $b^s_t$ are additive with respect to $t$.
\end{theorem}

\begin{proof}
Consider, as in Theorem 12.24, the following variable:
$$a=\sum_{k=1}^sw^ka_k$$

We have the following computation, for the corresponding Fourier transform:
\begin{eqnarray*}
\log F_a(z)
&=&\sum_{k=1}^s\log F_{a_k}(w^kz)\\
&=&\sum_{k=1}^s\frac{t}{s}\left(\exp(w^kz)-1\right)
\end{eqnarray*}

But this gives the following formula, in terms of the above function $\exp_s$:
\begin{eqnarray*}
\log F_a(z)
&=&t\left(\left(\frac{1}{s}\sum_{k=1}^s\exp(w^kz)\right)-1\right)\\
&=&t\left(\exp_s(z)-1\right)
\end{eqnarray*}

Now since $b^s_t$ is the law of $a$, this gives the formula in the statement.
\end{proof}

Regarding now the densities of the Bessel laws, these are as follows:

\begin{theorem}
We have the following formula,
$$b^s_t=e^{-t}\sum_{p_1=0}^\infty\ldots\sum_{p_s=0}^\infty\frac{1}{p_1!\ldots p_s!}\,\left(\frac{t}{s}\right)^{p_1+\ldots+p_s}\delta\left(\sum_{k=1}^sw^kp_k\right)$$
where $w=e^{2\pi i/s}$, and the $\delta$ symbol is a Dirac mass.
\end{theorem}

\begin{proof}
The Fourier transform of the measure on the right is given by:
\begin{eqnarray*}
F(z)
&=&e^{-t}\sum_{p_1=0}^\infty\ldots\sum_{p_s=0}^\infty\frac{1}{p_1!\ldots p_s!}\left(\frac{t}{s}\right)^{p_1+\ldots+p_s}F\delta\left(\sum_{k=1}^sw^kp_k\right)(z)\\
&=&e^{-t}\sum_{p_1=0}^\infty\ldots\sum_{p_s=0}^\infty\frac{1}{p_1!\ldots p_s!}\left(\frac{t}{s}\right)^{p_1+\ldots+p_s}\exp\left(\sum_{k=1}^sw^kp_kz\right)\\
&=&e^{-t}\sum_{r=0}^\infty\left(\frac{t}{s}\right)^r\sum_{\Sigma p_i=r}\frac{\exp\left(\sum_{k=1}^sw^kp_kz\right)}{p_1!\ldots p_s!}
\end{eqnarray*}

We multiply by $e^t$, and we compute the derivative with respect to $t$:
\begin{eqnarray*}
(e^tF(z))'
&=&\sum_{r=1}^\infty\frac{r}{s}\left(\frac{t}{s}\right)^{r-1}\sum_{\Sigma p_i=r}\frac{\exp\left(\sum_{k=1}^sw^kp_kz\right)}{p_1!\ldots p_s!}\\
&=&\frac{1}{s}\sum_{r=1}^\infty\left(\frac{t}{s}\right)^{r-1}\sum_{\Sigma p_i=r}\left(\sum_{l=1}^sp_l\right)\frac{\exp\left(\sum_{k=1}^sw^kp_kz\right)}{p_1!\ldots p_s!}\\
&=&\frac{1}{s}\sum_{r=1}^\infty\left(\frac{t}{s}\right)^{r-1}\sum_{\Sigma p_i=r}\sum_{l=1}^s\frac{\exp\left(\sum_{k=1}^sw^kp_kz\right)}{p_1!\ldots p_{l-1}!(p_l-1)!p_{l+1}!\ldots p_s!}\\
\end{eqnarray*}

By using the variable $u=r-1$, we get:
\begin{eqnarray*}
(e^tF(z))'
&=&\frac{1}{s}\sum_{u=0}^\infty\left(\frac{t}{s}\right)^u\sum_{\Sigma q_i=u}\sum_{l=1}^s\frac{\exp\left(w^lz+\sum_{k=1}^sw^kq_kz\right)}{q_1!\ldots q_s!}\\
&=&\left(\frac{1}{s}\sum_{l=1}^s\exp(w^lz)\right)\left(\sum_{u=0}^\infty\left(\frac{t}{s}\right)^u\sum_{\Sigma q_i=u}\frac{\exp\left(\sum_{k=1}^sw^kq_kz\right)}{q_1!\ldots q_s!}\right)\\
&=&(\exp_sz)(e^tF(z))
\end{eqnarray*}

On the other hand, consider the following function:
$$\Phi(t)=\exp(t\exp_sz)$$

This function satisfies as well the equation found above, namely:
$$\Phi'(t)=(\exp_sz)\Phi(t)$$

We conclude from this that we have the following equality of functions:
$$e^tF(z)=\Phi(t)$$

But this gives the following formula, for the logarithm of the Fourier transform:
\begin{eqnarray*}
\log F
&=&\log(e^{-t}\exp(t\exp_sz))\\
&=&\log(\exp(t(\exp_sz-1)))\\
&=&t(\exp_sz-1)
\end{eqnarray*}

Thus, we are led to the formulae in the statement.
\end{proof}

\section*{12d. Wigner laws}

In the continuous group case now, as a continuation of the above investigations, an interesting input comes from the various computations done some time ago in chapter 6. In order to discuss all this, let us first recall some useful formulae from chapter 6. One of the key results there, which is very useful in practice, was as follows:

\begin{theorem}
The polynomial integrals over the unit sphere $S^{N-1}_\mathbb R\subset\mathbb R^N$, with respect to the normalized, mass $1$ measure, are given by the following formula,
$$\int_{S^{N-1}_\mathbb R}x_1^{k_1}\ldots x_N^{k_N}\,dx=\frac{(N-1)!!k_1!!\ldots k_N!!}{(N+\Sigma k_i-1)!!}$$
valid when all exponents $k_i$ are even. If an exponent is odd, the integral vanishes.
\end{theorem}

\begin{proof}
This is something that we know from chapter 6, the idea being that the $N=2$ case is solved by the Wallis formula, and that the general case, $N\in\mathbb N$, follows from this, by using spherical coordinates and the Fubini theorem.
\end{proof}

As an application of the above formula, also following chapter 6, we have:

\index{hyperspherical law}

\begin{theorem}
The moments of the hyperspherical variables are
$$\int_{S^{N-1}_\mathbb R}x_i^kdx=\frac{(N-1)!!k!!}{(N+k-1)!!}$$
and the rescalings $y_i=x_i/\sqrt{N}$ become normal and independent with $N\to\infty$.
\end{theorem}

\begin{proof}
This is something that we know from chapter 6, coming from:
\begin{eqnarray*}
\int_{S^{N-1}_\mathbb R}x_i^kdx
&=&\frac{(N-1)!!k!!}{(N+k-1)!!}\\
&\simeq&N^{k/2}k!!\\
&=&N^{k/2}M_k(g_1)
\end{eqnarray*}

As for the asymptotic independence result, this is standard as well, once again by using Theorem 12.27, for computing mixed moments, and taking the $N\to\infty$ limit.
\end{proof}

Now back to groups, we can talk as well about rotation groups, as follows:

\begin{theorem}
We have the integration formula
$$\int_{O_N}U_{ij}^kdU=\frac{(N-1)!!k!!}{(N+k-1)!!}$$
and the rescalings $V_{ij}=U_{ij}/\sqrt{N}$ become normal and independent with $N\to\infty$.
\end{theorem}

\begin{proof}
We use the well-known fact that we have an embedding as follows, for any $i$, which makes correspond the respective integration functionals:
$$C(S^{N-1}_\mathbb R)\subset C(O_N)\quad,\quad 
x_i\to U_{1i}$$

With this identification made, the result follows from Theorem 12.28.
\end{proof}

We have similar results in the unitary case. First, we have:

\begin{theorem}
We have the following integration formula over the complex sphere $S^{N-1}_\mathbb C\subset\mathbb R^N$, with respect to the normalized measure, 
$$\int_{S^{N-1}_\mathbb C}|z_1|^{2l_1}\ldots|z_N|^{2l_N}\,dz=4^{\sum l_i}\frac{(2N-1)!l_1!\ldots l_n!}{(2N+\sum l_i-1)!}$$
valid for any exponents $l_i\in\mathbb N$. As for the other polynomial integrals in $z_1,\ldots,z_N$ and their conjugates $\bar{z}_1,\ldots,\bar{z}_N$, these all vanish.
\end{theorem}

\begin{proof}
As before, this is something that we know from chapter 6, and which can be proved either directly, or by using the formula in Theorem 12.27.
\end{proof}

We can talk about complex hyperspherical laws, and we have:

\begin{theorem}
The rescaled coordinates on the complex sphere $S^{N-1}_\mathbb C$, 
$$w_i=\frac{z_i}{\sqrt{N}}$$
become complex Gaussian and independent with $N\to\infty$.
\end{theorem}

\begin{proof}
This follows as in the proof of Theorem 12.28, by using Theorem 12.30.
\end{proof}

In relation now with rotation groups, the result that we obtain is as follows:

\begin{theorem}
For the unitary group $U_N$, the normalized coordinates
$$V_{ij}=\frac{U_{ij}}{\sqrt{N}}$$
become complex Gaussian and independent with $N\to\infty$.
\end{theorem}

\begin{proof}
We use the well-known fact that we have an embedding as follows, for any $i$, which makes correspond the respective integration functionals:
$$C(S^{N-1}_\mathbb C)\subset C(U_N)\quad,\quad 
x_i\to U_{1i}$$

With this identification made, the result follows from Theorem 12.31.
\end{proof}

Our claim now is that the above results can be reformulated in terms of the truncated characters introduced in chapter 11. Let us recall indeed from there that we have:

\begin{definition}
Given a closed subgroup $G\subset U_N$, the function
$$\chi:G\to\mathbb C\quad,\quad 
\chi_t(g)=\sum_{i=1}^{[tN]}g_{ii}$$
is called main truncated character of $G$, of parameter $t\in(0,1]$.
\end{definition}

In connection now with the present considerations, the point is that with the above notion in hand, our results above reformulate as follows:

\begin{theorem}
For the orthogonal and unitary groups $O_N,U_N$, the rescalings
$$\chi=\frac{\chi_{1/N}}{\sqrt{N}}$$
become respectively real and complex Gaussian, in the $N\to\infty$ limit.
\end{theorem}

\begin{proof}
According to our conventions, given a closed subgroup $G\subset U_N$, the main character truncated at $t=1/N$ is simply the first coordinate:
$$\chi_{1/N}(g)=g_{11}$$

With this remark made, the conclusions from the statement follow from the computations performed above, for the laws of coordinates on $O_N,U_N$.
\end{proof}

It is possible to get beyond such results, by using advanced representation theory methods, with full results about all the truncated characters, and in particular about the main characters. We will be back to this in Part IV below.

\bigskip

As a last topic now for this chapter, let us discuss the case where $N$ is fixed. Things are quite complicated here, and as a main goal, we would like to find the law of the main character for our favorite rotation groups, namely $SU_2$ and $SO_3$. 

\bigskip

In order to do so, we will need some combinatorial preliminaries. We first have the following well-known result, which is the cornerstone of all modern combinatorics:

\index{Catalan numbers}
\index{noncrossing pairings}

\begin{theorem}
The Catalan numbers, which are by definition given by
$$C_k=|NC_2(2k)|$$
satisfy the following recurrence formula,
$$C_{k+1}=\sum_{a+b=k}C_aC_b$$ 
and their generating series, given by definition by
$$f(z)=\sum_{k\geq0}C_kz^k$$
satisfies the following degree $2$ equation,
$$zf^2-f+1=0$$
and we have the following explicit formula for these numbers:
$$C_k=\frac{1}{k+1}\binom{2k}{k}$$
Numerically, these numbers are $1,1,2,5,14,42,132,429,1430,4862,16796,\ldots$
\end{theorem}

\begin{proof}
We must count the noncrossing pairings of $\{1,\ldots,2k\}$. But such a pairing appears by pairing 1 to an odd number, $2a+1$, and then inserting a noncrossing pairing of $\{2,\ldots,2a\}$, and a noncrossing pairing of $\{2a+2,\ldots,2l\}$. We conclude from this that we have the following recurrence formula for the Catalan numbers:
$$C_k=\sum_{a+b=k-1}C_aC_b$$ 

In terms of the generating series $f$, the above recurrence gives:
\begin{eqnarray*}
zf^2
&=&\sum_{a,b\geq0}C_aC_bz^{a+b+1}\\
&=&\sum_{k\geq1}\sum_{a+b=k-1}C_aC_bz^k\\
&=&\sum_{k\geq1}C_kz^k\\
&=&f-1
\end{eqnarray*}

Thus the generating series $f$ satisfies the following degree 2 equation:
$$zf^2-f+1=0$$

By choosing the solution which is bounded at $z=0$, we obtain:
$$f(z)=\frac{1-\sqrt{1-4z}}{2z}$$ 

By using now the Taylor formula for $\sqrt{x}$, we obtain the following formula:
$$f(z)=\sum_{k\geq0}\frac{1}{k+1}\binom{2k}{k}z^k$$

It follows that the Catalan numbers are given by the formula the statement.
\end{proof}

The Catalan numbers are central objects in probability as well, and we have the following key result here, complementing the formulae from Theorem 12.35:

\index{semicircle law}
\index{Wigner law}
\index{Catalan numbers}

\begin{theorem}
The normalized Wigner semicircle law, which is by definition
$$\gamma_1=\frac{1}{2\pi}\sqrt{4-x^2}dx$$
has the Catalan numbers as even moments. As for the odd moments, these all vanish. 
\end{theorem}

\begin{proof}
The even moments of the Wigner law can be computed with the change of variable $x=2\cos t$, and we are led to the following formula:
\begin{eqnarray*}
M_{2k}
&=&\frac{1}{\pi}\int_0^2\sqrt{4-x^2}x^{2k}dx\\
&=&\frac{1}{\pi}\int_0^{\pi/2}\sqrt{4-4\cos^2t}\,(2\cos t)^{2k}2\sin t\,dt\\
&=&\frac{4^{k+1}}{\pi}\int_0^{\pi/2}\cos^{2k}t\sin^2t\,dt\\
&=&\frac{4^{k+1}}{\pi}\cdot\frac{\pi}{2}\cdot\frac{(2k)!!2!!}{(2k+3)!!}\\
&=&2\cdot 4^k\cdot\frac{(2k)!/2^kk!}{2^{k+1}(k+1)!}\\
&=&C_k
\end{eqnarray*}

As for the odd moments, these all vanish, because the density of $\gamma_1$ is an even function. Thus, we are led to the conclusion in the statement.
\end{proof}

We can now formulate our result regarding $SU_2$, as follows: 

\begin{theorem}
The main character of $SU_2$, given by
$$\chi\begin{pmatrix}a&b\\ -\bar{b}&\bar{a}\end{pmatrix}=2Re(a)$$
follows a Wigner semicircle law $\gamma_1$.
\end{theorem}

\begin{proof}
The idea is that this follows by identifying $SU_2$ with the sphere $S^3_\mathbb R\subset\mathbb R^4$, and the uniform measure on $SU_2$ with the uniform measure on this sphere. Indeed, in terms of the standard parametrization of $SU_2$, from chapter 10, written in real form, we have the following formula, for the main character of $SU_2$:
$$\chi\begin{pmatrix}x+iy&z+it\\ -z+it&x-iy\end{pmatrix}=2x$$

We are therefore left with computing the law of the following variable:
$$x\in C(S^3_\mathbb R)$$

But for this purpose, we can use moments. Indeed, Theorem 12.27 gives:
\begin{eqnarray*}
\int_{S^3_\mathbb R}x^{2k}
&=&\frac{3!!(2k)!!}{(2k+3)!!}\\
&=&2\cdot\frac{3\cdot5\cdot7\ldots (2k-1)}{2\cdot4\cdot6\ldots (2k+2)}\\
&=&2\cdot\frac{(2k)!}{2^kk!2^{k+1}(k+1)!}\\
&=&\frac{1}{4^k}\cdot\frac{1}{k+1}\binom{2k}{k}\\
&=&\frac{C_k}{4^k}
\end{eqnarray*}

Thus the variable $2x\in C(S^3_\mathbb R)$ has the Catalan numbers as even moments, and so by Theorem 12.36 its distribution is the Wigner semicircle law $\gamma_1$, as claimed.
\end{proof}

In order to do the computation for $SO_3$, we will need some more probabilistic preliminaries, which are standard random matrix theory material. Let us start with:

\index{noncrossing partitions}
\index{noncrossing pairings}
\index{fattening partitions}
\index{shrinking partitions}

\begin{proposition}
We have a bijection $NC(k)\simeq NC_2(2k)$, as follows:
\begin{enumerate}
\item The application $NC(k)\to NC_2(2k)$ is the ``fattening'' one, obtained by doubling all the legs, and doubling all the strings too.

\item Its inverse $NC_2(2k)\to NC(k)$ is the ``shrinking'' application, obtained by collapsing pairs of consecutive neighbors.
\end{enumerate}
\end{proposition}

\begin{proof}
This is something self-explanatory, and in order to see how this works, let us discuss an example. Consider a noncrossing partition, say the following one:
$$\xymatrix@R=10pt@C=10pt{
\ar@{-}[rrrrrrr]&&&&&&&\\
&\ar@{-}[rr]&&&&\ar@{-}[r]&&\\
1\ar@{-}[uu]&2\ar@{-}[u]&3\ar@{-}[u]&4\ar@{-}[u]&5\ar@{-}[uu]&6\ar@{-}[u]&7\ar@{-}[u]&8\ar@{-}[uu]
}$$

Now let us ``fatten'' this partition, by doubling everything, as follows:
$$\xymatrix@R=10pt@C=10pt{
\ar@{=}[rrrrrrr]&&&&&&&\\
&\ar@{=}[rr]&&&&\ar@{=}[r]&&\\
11'\ar@{=}[uu]&22'\ar@{=}[u]&33'\ar@{=}[u]&44'\ar@{=}[u]&55'\ar@{=}[uu]&66'\ar@{=}[u]&77'\ar@{=}[u]&88'\ar@{=}[uu]
}$$

Now by relabeling the points $1,\ldots,16$, what we have is indeed a noncrossing pairing. As for the reverse operation, that is obviously obtained by ``shrinking'' our pairing, by collapsing pairs of consecutive neighbors, that is, by identifying $1=2$, then $3=4$, then $5=6$, and so on, up to $15=16$. Thus, we are led to the conclusion in the statement.
\end{proof}

As a consequence of the above result, we have a new look on the Catalan numbers, which is more adapted to our present $SO_3$ considerations, as follows:

\begin{proposition}
The Catalan numbers $C_k=|NC_2(2k)|$ appear as well as
$$C_k=|NC(k)|$$
where $NC(k)$ is the set of all noncrossing partitions of $\{1,\ldots,k\}$.
\end{proposition}

\begin{proof}
This follows indeed from Proposition 12.38.
\end{proof}

Let us formulate now the following definition:

\begin{definition}
The standard Marchenko-Pastur law $\pi_1$ is given by:
$$f\sim\gamma_1\implies f^2\sim\pi_1$$
That is, $\pi_1$ is the law of the square of a variable following the semicircle law $\gamma_1$.
\end{definition}

Here the fact that $\pi_1$ is indeed well-defined comes from the fact that a measure is uniquely determined by its moments. More explicitly now, we have:

\index{Marchenko-Pastur law}

\begin{proposition}
The density of the Marchenko-Pastur law is
$$\pi_1=\frac{1}{2\pi}\sqrt{4x^{-1}-1}\,dx$$
and the moments of this measure are the Catalan numbers.
\end{proposition}

\begin{proof}
The moments of the law in the statement can be computed with the change of variable $x=4\cos^2t$, and we are led to the following formula:
\begin{eqnarray*}
M_k
&=&\frac{1}{2\pi}\int_0^4\sqrt{4x^{-1}-1}\,x^kdx\\
&=&\frac{1}{2\pi}\int_0^{\pi/2}\frac{\sin t}{\cos t}\cdot(4\cos^2t)^k\cdot 2\cos t\sin t\,dt\\
&=&\frac{4^{k+1}}{\pi}\int_0^{\pi/2}\cos^{2k}t\sin^2t\,dt\\
&=&\frac{4^{k+1}}{\pi}\cdot\frac{\pi}{2}\cdot\frac{(2k)!!2!!}{(2k+3)!!}\\
&=&2\cdot 4^k\cdot\frac{(2k)!/2^kk!}{2^{k+1}(k+1)!}\\
&=&C_k
\end{eqnarray*}

Thus, we are led to the conclusion in the statement.
\end{proof}

We can do now the character computation for $SO_3$, as follows:

\begin{theorem}
The main character of $SO_3$, modified by adding $1$ to it, given in standard Euler-Rodrigues coordinates by
$$\chi=3x^2-y^2-z^2-t^2$$
follows a squared semicircle law, or Marchenko-Pastur law $\pi_1$.
\end{theorem}

\begin{proof}
The idea is that this follows by using the canonical quotient map $SU_2\to SO_3$, and the result for $SU_2$ from Theorem 12.37. To be more precise, let us recall from chapter 10 that the elements of $SU_2$ can be parametrized as follows:
$$U=\begin{pmatrix}x+iy&z+it\\ -z+it&x-iy\end{pmatrix}$$

As for the elements of $SO_3$, these can be parametrized as follows:
$$V=\begin{pmatrix}
x^2+y^2-z^2-t^2&2(yz-xt)&2(xz+yt)\\
2(xt+yz)&x^2+z^2-y^2-t^2&2(zt-xy)\\
2(yt-xz)&2(xy+zt)&x^2+t^2-y^2-z^2
\end{pmatrix}$$

The point now is that, by using the above two formulae, in the context of the computation from Theorem 12.37, the main character of $SO_3$ is given by:
\begin{eqnarray*}
\chi
&=&Tr(V)+1\\
&=&3x^2-y^2-z^2-t^2+1\\
&=&4x^2
\end{eqnarray*}

Now recall from the proof of Theorem 12.37 that we have:
$$2x\sim\gamma_1$$

On the other hand, a quick comparison between the moment formulae for the Wigner and Marchenko-Pastur laws, which are very similar, shows that we have:
$$f\sim\gamma_1\implies f^2\sim\pi_1$$

Thus, with $f=2x$, we obtain the result in the statement.
\end{proof}

As an interesting question now, appearing from the above, and which is quite philosophical, we have the problem of understanding how the Wigner and Marchenko-Pastur laws $\gamma_1,\pi_1$ fit in regards with the main limiting laws from classical probability. 

\bigskip

The answer here is quite tricky, the idea being that, with a suitable formalism for freeness, $\gamma_1,\pi_1$ can be thought of as being ``free analogues'' of the Gaussian and Poisson laws $g_1,p_1$. This is something quite subtle, requiring some further knowledge, and we will be back to this in Part IV below, when doing representation theory.

\section*{12e. Exercises}

There has been a lot of technical material in this chapter, with substantial combinatorics, and technical as well will be most of our exercises. First, we have:

\begin{exercise}
Work out the moment formula for Bessel laws, $M_k=|P^s(k)|$, where $P^s(k)$ are the partitions satisfying $\#\circ=\#\bullet(s)$, as a weighted sum, in each block.
\end{exercise}

This is something that we briefly discussed in the above, and the problem is now that of working out all the details, first as $s=1,2,\infty$, and then in general.

\begin{exercise}
Work out all details for the truncated character formula for $H_N^s$,
$$\chi_t\sim b^s_t$$
where $b_t^s=p_{t\varepsilon_s}$, with $\varepsilon_s$ being the uniform measure on the $s$-th roots of unity.
\end{exercise}

As before, this is something that we briefly discussed in the above, and the problem is now that of working out all the details, first as $s=1,2,\infty$, and then in general.

\begin{exercise}
Show that the passage from $H_N^s$ to $H_N^{sd}$ does not change the asymptotic laws of the truncated characters.
\end{exercise}

This is something that we discussed in the previous chapter, in a particular case, namely for the passage from the symmetric group $S_N$ to the alternating group $A_N$.

\begin{exercise}
Compute the asymptotic laws of characters and coordinates for the bistochastic groups $B_N$ and $C_N$, as well as for the symplectic group $Sp_N\subset U_N$.
\end{exercise}

These computations are all quite standard, with the computation for $B_N$ being quite similar to that for $O_{N-1}$, the computation for $C_N$ being quite similar to that for $U_{N-1}$, and the computation for $Sp_N$ being quite similar to that for $O_{N-1}$.

\begin{exercise}
Compute the character laws for the groups $O_1$, $SO_1$, then for the groups $U_1$, $SU_1$, and then for the groups $O_2$, $SO_2$.
\end{exercise}

As before with the previous exercise, the computations here are quite standard. In fact, the more difficult questions of this type concern the next groups in the above series, namely $SU_2$ and $SO_3$, which were discussed in the above.

\begin{exercise}
Work out all the combinatorics and calculus details in relation with the Wigner and Marchenko-Pastur laws, and their moments, the Catalan numbers.
\end{exercise}

This is a very instructive exercise, with lots of nice combinatorics involved. Most of this combinatorics was actually already discussed in the above.

\part{Haar integration}

\ \vskip50mm

\begin{center}
{\em And the band plays Waltzing Matilda

And the old men still answer the call

But year after year, their numbers get fewer

Someday, no one will march there at all}
\end{center}

\chapter{Representations}

\section*{13a. Basic theory}

We have seen so far that some algebraic and probabilistic theory for the finite subgroups $G\subset U_N$, ranging from elementary to quite advanced, can be developed. We have seen as well a few computations for the continuous compact subgroups $G\subset U_N$. In what follows we develop some systematic theory for the arbitrary closed subgroups $G\subset U_N$, covering both the finite and the infinite case. The main examples that we have in mind, and the questions that we would like to solve for them, are as follows:

\medskip

\begin{enumerate}
\item The orthogonal and unitary groups $O_N,U_N$. Here we would like to have an integration formula, and results about character laws, in the $N\to\infty$ limit.

\medskip

\item Various versions of $O_N,U_N$, such as the bistochastic groups $B_N,C_N$, or the symplectic groups $Sp_N$, with similar questions to be solved.

\medskip

\item The reflection groups $H_N^{sd}\subset U_N$, with results about characters extending, or at least putting in a more conceptual framework, what we already have. 
\end{enumerate}

\medskip

There is a lot of theory to be developed, and we will do this gradually. To be more precise, in this chapter and in the next one we will work out algebraic aspects, and then in the chapter afterwards and in the last one we will use these algebraic techniques, in order to work out probabilistic results, and in particular to answer the above questions. As before, the main notion that we will be interested in is that of a representation:

\index{representation}
\index{character}

\begin{definition}
A representation of a compact group $G$ is a continuous group morphism, which can be faithful or not, into a unitary group:
$$u:G\to U_N$$
The character of such a representation is the function $\chi:G\to\mathbb C$ given by
$$g\to Tr(u_g)$$
where $Tr$ is the usual trace of the $N\times N$ matrices, $Tr(M)=\sum_iM_{ii}$.
\end{definition}

As a basic example here, for any compact group we always have available the trivial 1-dimensional representation, or character, which is by definition as follows:
$$u:G\to U_1\quad,\quad 
g\to(1)$$

In fact, talking 1-dimensional representations, we already know about these, from chapter 9, with the summary of our results there being as follows:

\begin{theorem}
The $1$-dimensional representations of $G$ are the morphisms
$$u:G\to\mathbb T$$
and we have $u=\chi$ in this case. These morphisms, or characters, must come via
$$u:G\to G_{ab}\to\mathbb T$$
from the characters $G_{ab}\to\mathbb T$, which themselves form a group, which is the dual $\widehat{G}_{ab}$.
\end{theorem}

\begin{proof}
This is indeed self-explanatory, coming in the finite group case from our discussion from chapter 9, and in general, via a straightforward extension of this.
\end{proof}

Moving now to higher dimensions, as another class of basic examples, we have:

\begin{theorem}
Available for any finite group $G$ is its regular representation
$$u:G\subset S_N\subset O_N\subset U_N$$
with $N=|G|$, obtained via Cayley and permutation matrices, the formula being
$$u_g(e_h)=e_{gh}$$
with $\{e_h|h\in G\}$ being the standard basis of $\mathbb C^N$. Its character is $\chi(g)=N\delta_{g1}$.
\end{theorem}

\begin{proof}
This is again something self-explanatory, coming from our discussion from chapter 9, on the Cayley theorem, permutation matrices and related topics, and with the character computation being something elementary too, as follows:
\begin{eqnarray*}
\chi(g)
&=&Tr(u_g)\\
&=&\sum_{h\in G}<u_g(e_h),e_h>\\
&=&\sum_{h\in G}<e_{gh},e_h>\\
&=&N\delta_{g1}
\end{eqnarray*}

Thus, we are led to the conclusions in the statement.
\end{proof}

Summarizing, we definitely have interesting illustrations for Definition 13.1, and even some beginning of theory on the way, based on our material from chapter 9.

\bigskip

What is next? You guessed it right, more examples. Inspired by the above, let us formulate the following question, which looks like something quite interesting:

\begin{question}
Given a subgroup $G\subset U_N$, besides its fundamental representation
$$u:G\subset U_N\quad,\quad g\to g$$
we can equally talk about its conjugate fundamental representation
$$\bar{u}:G\subset U_N\quad,\quad g\to\bar{g}$$
and probably about many more, coming via other operations. What exactly are these?
\end{question} 

To be more precise here, consider the usual conjugation of the unitary matrices, $(\bar{U})_{ij}=\bar{U}_{ij}$. This can be viewed as a group isomorphism, as follows:
$$U_N\simeq U_N\quad,\quad U\to\bar{U}$$

Now given an embedding $u:G\subset U_N$, we can compose it with this isomorphism $U_N\simeq U_N$, and we obtain another embedding $\bar{u}:G\subset U_N$. And with $\bar{u}$ being in general different from $u$ itself, as the 1D examples, in the context of Theorem 13.2, show.

\bigskip

In order to answer Question 13.4, and see which representations are available, let us first discuss the various operations on the representations. We have here:

\index{sum of representations}
\index{product of representations}
\index{conjugate representation}
\index{spinned representation}

\begin{proposition}
The representations of a given compact group $G$ are subject to the following operations:
\begin{enumerate}
\item Making sums. Given representations $u,v$, having dimensions $N,M$, their sum is the $N+M$-dimensional representation $u+v=diag(u,v)$.

\item Making products. Given representations $u,v$, having dimensions $N,M$, their tensor product is the $NM$-dimensional representation $(u\otimes v)_{ia,jb}=u_{ij}v_{ab}$.

\item Taking conjugates. Given a representation $u$, having dimension $N$, its complex conjugate is the $N$-dimensional representation $(\bar{u})_{ij}=\bar{u}_{ij}$.

\item Spinning by unitaries. Given a representation $u$, having dimension $N$, and a unitary $V\in U_N$, we can spin $u$ by this unitary, $u\to VuV^*$.
\end{enumerate}
\end{proposition}

\begin{proof}
The fact that the operations in the statement are indeed well-defined, among maps from $G$ to unitary groups, can be checked as follows:

\medskip

(1) This follows from the trivial fact that if $g\in U_N$ and $h\in U_M$ are two unitaries, then their diagonal sum is a unitary too, as follows:
$$\begin{pmatrix}g&0\\ 0&h\end{pmatrix}\in U_{N+M}$$

(2) This follows from the fact that if $g\in U_N$ and $h\in U_M$ are two unitaries, then $g\otimes h\in U_{NM}$ is a unitary too. Given unitaries $g,h$, let us set indeed:
$$(g\otimes h)_{ia,jb}=g_{ij}h_{ab}$$

This matrix is then a unitary too, as shown by the following computation:
\begin{eqnarray*}
[(g\otimes h)(g\otimes h)^*]_{ia,jb}
&=&\sum_{kc}(g\otimes h)_{ia,kc}((g\otimes h)^*)_{kc,jb}\\
&=&\sum_{kc}(g\otimes h)_{ia,kc}\overline{(g\otimes h)_{jb,kc}}\\
&=&\sum_{kc}g_{ik}h_{ac}\bar{g}_{jk}\bar{h}_{bc}\\
&=&\sum_kg_{ik}\bar{g}_{jk}\sum_ch_{ac}\bar{h}_{bc}\\
&=&\delta_{ij}\delta_{ab}
\end{eqnarray*}

(3) This simply follows from the fact that if $g\in U_N$ is unitary, then so is its complex conjugate, $\bar{g}\in U_N$, and this due to the following formula, obtained by conjugating:
$$g^*=g^{-1}\implies g^t=\bar{g}^{-1}$$

(4) This is clear as well, because if $g\in U_N$ is unitary, and $V\in U_N$ is another unitary, then we can spin $g$ by this unitary, and we obtain a unitary as follows:
$$VgV^*\in U_N$$

Thus, our operations are well-defined, and this leads to the above conclusions.
\end{proof}

In relation now with characters, we have the following result:

\index{character}

\begin{proposition}
We have the following formulae, regarding characters
$$\chi_{u+v}=\chi_u+\chi_v\quad,\quad 
\chi_{u\otimes v}=\chi_u\chi_v\quad,\quad 
\chi_{\bar{u}}=\bar{\chi}_u\quad,\quad 
\chi_{VuV^*}=\chi_u$$
in relation with the basic operations for the representations.
\end{proposition}

\begin{proof}
All these assertions are elementary, by using the following well-known trace formulae, valid for any two square matrices $g,h$, and any unitary $V$:
$$Tr(diag(g,h))=Tr(g)+Tr(h)\quad,\quad 
Tr(g\otimes h)=Tr(g)Tr(h)$$
$$Tr(\bar{g})=\overline{Tr(g)}\quad,\quad 
Tr(VgV^*)=Tr(g)$$

To be more precise, the first formula is clear from definitions. Regarding now the second formula, the computation here is immediate too, as follows:
\begin{eqnarray*}
Tr(g\otimes h)
&=&\sum_{ia}(g\otimes h)_{ia,ia}\\
&=&\sum_{ia}g_{ii}h_{aa}\\
&=&Tr(g)Tr(h)
\end{eqnarray*}

Regarding now the third formula, this is clear from definitions, by conjugating. Finally, regarding the fourth formula, this can be established as follows:
$$Tr(VgV^*)=Tr(gV^*V)=Tr(g)$$

Thus, we are led to the conclusions in the statement.
\end{proof}

Assume now that we are given a closed subgroup $G\subset U_N$. By using the above operations, we can construct a whole family of representations of $G$, as follows:

\index{Peter-Weyl representations}
\index{colored integer}

\begin{definition}
Given a closed subgroup $G\subset U_N$, its Peter-Weyl representations are the tensor products between the fundamental representation and its conjugate:
$$u:G\subset U_N\quad,\quad 
\bar{u}:G\subset U_N$$ 
We denote these tensor products $u^{\otimes k}$, with $k=\circ\bullet\bullet\circ\ldots$ being a colored integer, with the colored tensor powers being defined according to the rules 
$$u^{\otimes\circ}=u\quad,\quad 
u^{\otimes\bullet}=\bar{u}\quad,\quad
u^{\otimes kl}=u^{\otimes k}\otimes u^{\otimes l}$$
and with the convention that $u^{\otimes\emptyset}$ is the trivial representation $1:G\to U_1$.
\end{definition}

Here are a few examples of such Peter-Weyl representations, namely those coming from the colored integers of length 2, to be often used in what follows:
$$u^{\otimes\circ\circ}=u\otimes u\quad,\quad 
u^{\otimes\circ\bullet}=u\otimes\bar{u}$$
$$u^{\otimes\bullet\circ}=\bar{u}\otimes u\quad,\quad 
u^{\otimes\bullet\bullet}=\bar{u}\otimes\bar{u}$$

In relation now with characters, we have the following result:

\index{colored powers}

\begin{proposition}
The characters of Peter-Weyl representations are given by
$$\chi_{u^{\otimes k}}=(\chi_u)^k$$
with the colored powers of a variable $\chi$ being by definition given by
$$\chi^\circ=\chi\quad,\quad
\chi^\bullet=\bar{\chi}\quad,\quad
\chi^{kl}=\chi^k\chi^l$$
and with the convention that $\chi^\emptyset$ equals by definition $1$.
\end{proposition}

\begin{proof}
This follows indeed from the additivity, multiplicativity and conjugation formulae established in Proposition 13.6, via the conventions in Definition 13.7.
\end{proof}

Getting back now to our motivations, we can see the interest in the above constructions. Indeed, the joint moments of the main character $\chi=\chi_u$ and its adjoint $\bar{\chi}=\chi_{\bar{u}}$ are simply the expectations of the characters of various Peter-Weyl representations: 
$$\int_G\chi^k=\int_G \chi_{u^{\otimes k}}$$

Summarizing, given a closed subgroup $G\subset U_N$, we would like to understand its Peter-Weyl representations, and compute the expectations of the characters of these representations. In order to do so, let us formulate the following key definition:

\index{Hom space}
\index{End space}
\index{Fix space}

\begin{definition}
Given a compact group $G$, and two of its representations,
$$u:G\to U_N\quad,\quad 
v:G\to U_M$$
we define the linear space of intertwiners between these representations as being 
$$Hom(u,v)=\left\{T\in M_{M\times N}(\mathbb C)\Big|Tu_g=v_gT,\forall g\in G\right\}$$
and we use the following conventions:
\begin{enumerate}
\item We use the notations $Fix(u)=Hom(1,u)$, and $End(u)=Hom(u,u)$.

\item We write $u\sim v$ when $Hom(u,v)$ contains an invertible element.

\item We say that $u$ is irreducible, and write $u\in Irr(G)$, when $End(u)=\mathbb C1$.
\end{enumerate}
\end{definition}

The terminology here is very standard, with Hom and End standing for ``homomorphisms'' and ``endomorphisms'', and with Fix standing for ``fixed points''.

\bigskip

In practice, it is useful to think of the representations of $G$ as being the objects of some kind of abstract combinatorial structure associated to $G$, and of the intertwiners between these representations as being the ``arrows'' between these objects. We have in fact the following result, making the link with this viewpoint, called categorical:

\index{tensor category}

\begin{theorem}
The following happen:
\begin{enumerate}
\item The intertwiners are stable under composition:
$$T\in Hom(u,v)\ ,\ 
S\in Hom(v,w)
\implies ST\in Hom(u,w)$$

\item The intertwiners are stable under taking tensor products:
$$S\in Hom(u,v)\ ,\ 
T\in Hom(w,t)\\
\implies S\otimes T\in Hom(u\otimes w,v\otimes t)$$

\item The intertwiners are stable under taking adjoints:
$$T\in Hom(u,v)
\implies T^*\in Hom(v,u)$$

\item Thus, the Hom spaces form a tensor $*$-category.
\end{enumerate}
\end{theorem}

\begin{proof}
All this is clear from definitions, the verifications being as follows:

\medskip

(1) This follows indeed from the following computation, valid for any $g\in G$:
$$STu_g=Sv_gT=w_gST$$

(2) Again, this is clear, because we have the following computation:
\begin{eqnarray*}
(S\otimes T)(u_g\otimes w_g)
&=&Su_g\otimes Tw_g\\
&=&v_gS\otimes t_gT\\
&=&(v_g\otimes t_g)(S\otimes T)
\end{eqnarray*}

(3) This follows from the following computation, valid for any $g\in G$:
\begin{eqnarray*}
Tu_g=v_gT
&\implies&u_g^*T^*=T^*v_g^*\\
&\implies&T^*v_g=u_gT^*
\end{eqnarray*}

(4) This is just a conclusion of (1,2,3), with a tensor $*$-category being by definition an abstract beast satisfying these conditions (1,2,3). We will be back to tensor categories later on, in chapter 14 below, with more details on all this.
\end{proof}

The above result is quite interesting, because it shows that the combinatorics of a compact group $G$ is described by a certain collection of linear spaces, which can be in principle investigated by using tools from linear algebra. Thus, what we have here is a useful ``linearization'' idea. We will heavily use this idea, in what follows.

\section*{13b. Peter-Weyl theory}

In what follows we develop a systematic theory of the representations of the compact groups $G$, with emphasis on the Peter-Weyl representations, in the closed subgroup case $G\subset U_N$, that we are mostly interested in. Let us start with the following fact:

\begin{theorem}
Given a representation of a compact group $u:G\to U_N$, the corresponding linear space of self-intertwiners 
$$End(u)\subset M_N(\mathbb C)$$
is a $*$-algebra, with respect to the usual involution of the matrices.
\end{theorem}

\begin{proof}
By definition, the space $End(u)$ is a linear subspace of $M_N(\mathbb C)$. We know from Theorem 13.10 (1) that this subspace $End(u)$ is a subalgebra of $M_N(\mathbb C)$, and then we know as well from Theorem 13.10 (3) that this subalgebra is stable under the involution $*$. Thus, what we have here is a $*$-subalgebra of $M_N(\mathbb C)$, as claimed.
\end{proof}

The above result is quite interesting, because it gets us into linear algebra. Indeed, associated to any group representation $u:G\to U_N$ is now a quite familiar object, namely the algebra $End(u)\subset M_N(\mathbb C)$. In order to exploit this fact, we will need a well-known result, complementing the basic operator algebra theory from chapter 8, namely:

\index{operator algebra}
\index{finite dimensional algebra}

\begin{theorem}
Let $A\subset M_N(\mathbb C)$ be a $*$-algebra.
\begin{enumerate}
\item We can write $1=p_1+\ldots+p_k$, with $p_i\in A$ being central minimal projections.

\item The linear spaces $A_i=p_iAp_i$ are non-unital $*$-subalgebras of $A$.

\item We have a non-unital $*$-algebra sum decomposition $A=A_1\oplus\ldots\oplus A_k$.

\item We have unital $*$-algebra isomorphisms $A_i\simeq M_{n_i}(\mathbb C)$, with $n_i=rank(p_i)$.

\item Thus, we have a $*$-algebra isomorphism $A\simeq M_{n_1}(\mathbb C)\oplus\ldots\oplus M_{n_k}(\mathbb C)$.
\end{enumerate}
\end{theorem}

\begin{proof}
This is something very standard. Consider indeed an arbitrary $*$-algebra of the $N\times N$ matrices, $A\subset M_N(\mathbb C)$. Let us first look at the center of this algebra, $Z(A)=A\cap A'$. This center, viewed as an algebra, is then of the following form:
$$Z(A)\simeq\mathbb C^k$$

Consider now the standard basis $e_1,\ldots,e_k\in\mathbb C^k$, and let  $p_1,\ldots,p_k\in Z(A)$ be the images of these vectors via the above identification. In other words, these elements $p_1,\ldots,p_k\in A$ are central minimal projections, summing up to 1:
$$p_1+\ldots+p_k=1$$

The idea is then that this partition of the unity will eventually lead to the block decomposition of $A$, as in the statement. We prove this in 4 steps, as follows:

\medskip

\underline{Step 1}. We first construct the matrix blocks, our claim here being that each of the following linear subspaces of $A$ are non-unital $*$-subalgebras of $A$:
$$A_i=p_iAp_i$$

But this is clear, with the fact that each $A_i$ is closed under the various non-unital $*$-subalgebra operations coming from the projection equations $p_i^2=p_i^*=p_i$.

\medskip

\underline{Step 2}. We prove now that the above algebras $A_i\subset A$ are in a direct sum position, in the sense that we have a non-unital $*$-algebra sum decomposition, as follows:
$$A=A_1\oplus\ldots\oplus A_k$$

As with any direct sum question, we have two things to be proved here. First, by using the formula $p_1+\ldots+p_k=1$ and the projection equations $p_i^2=p_i^*=p_i$, we conclude that we have the needed generation property, namely:
$$A_1+\ldots+ A_k=A$$

As for the fact that the sum is indeed direct, this follows as well from the formula $p_1+\ldots+p_k=1$, and from the projection equations $p_i^2=p_i^*=p_i$.

\medskip

\underline{Step 3}. Our claim now, which will finish the proof, is that each of the $*$-subalgebras $A_i=p_iAp_i$ constructed above is in fact a full matrix algebra. To be more precise, with $n_i=rank(p_i)$, our claim is that we have isomorphisms, as follows:
$$A_i\simeq M_{n_i}(\mathbb C)$$

In order to prove this claim, recall that the projections $p_i\in A$ were chosen central and minimal. Thus, the center of each of the algebras $A_i$ reduces to the scalars:
$$Z(A_i)=\mathbb C$$

But this shows, either via a direct computation, or via the bicommutant theorem, that the each of the algebras $A_i$ is a full matrix algebra, as claimed.

\medskip

\underline{Step 4}. We can now obtain the result, by putting together what we have. Indeed, by using the results from Step 2 and Step 3, we obtain an isomorphism as follows:
$$A\simeq M_{n_1}(\mathbb C)\oplus\ldots\oplus M_{n_k}(\mathbb C)$$

In addition to this, a careful look at the isomorphisms established in Step 3 shows that at the global level, of the algebra $A$ itself, the above isomorphism simply comes by twisting the following standard multimatrix embedding, discussed in the beginning of the proof, (1) above, by a certain unitary matrix $U\in U_N$:
$$M_{n_1}(\mathbb C)\oplus\ldots\oplus M_{n_k}(\mathbb C)\subset M_N(\mathbb C)$$

Now by putting everything together, we obtain the result.
\end{proof}

We can now formulate our first Peter-Weyl theorem, as follows:

\index{Peter-Weyl}

\begin{theorem}[PW1]
Let $u:G\to U_N$ be a group representation, consider the algebra $A=End(u)$, and write its unit as above, as follows:
$$1=p_1+\ldots+p_k$$
The representation $u$ decomposes then as a direct sum, as follows,
$$u=u_1+\ldots+u_k$$
with each $u_i$ being an irreducible representation, obtained by restricting $u$ to $Im(p_i)$.
\end{theorem}

\begin{proof}
This basically follows from Theorem 13.11 and Theorem 13.12, as follows:

\medskip

(1) As a first observation, by replacing $G$ with its image $u(G)\subset U_N$, we can assume if we want that our representation $u$ is faithful, $G\subset_uU_N$. However, this replacement will not be really needed, and we will keep using $u:G\to U_N$, as above.

\medskip

(2) In order to prove the result, we will need some preliminaries. We first associate to our representation $u:G\to U_N$ the corresponding action map on $\mathbb C^N$. If a linear subspace $V\subset\mathbb C^N$ is invariant, the restriction of the action map to $V$ is an action map too, which must come from a subrepresentation $v\subset u$. This is clear indeed from definitions, and with the remark that the unitaries, being isometries, restrict indeed into unitaries.

\medskip

(3) Consider now a projection $p\in End(u)$. From $pu=up$ we obtain that the linear space $V=Im(p)$ is invariant under $u$, and so this space must come from a subrepresentation $v\subset u$. It is routine to check that the operation $p\to v$ maps subprojections to subrepresentations, and minimal projections to irreducible representations.

\medskip

(4) To be more precise here, the condition $p\in End(u)$ reformulates as follows:
$$pu_g=u_gp\quad,\quad\forall g\in G$$

As for the condition that $V=Im(p)$ is invariant, this reformulates as follows:
$$pu_gp=u_gp\quad,\quad\forall g\in G$$

Thus, we are in need of a technical linear algebra result, stating that for a projection $P\in M_N(\mathbb C)$ and a unitary $U\in U_N$, the following happens:
$$PUP=UP\implies PU=UP$$

(5) But this can be established with some $C^*$-algebra know-how, as follows:
\begin{eqnarray*}
tr[(PU-UP)(PU-UP)^*]
&=&tr[(PU-UP)(U^*P-PU^*)]\\
&=&tr[P-PUPU^*-UPU^*P+UPU^*]\\
&=&tr[P-UPU^*-UPU^*+UPU^*]\\
&=&tr[P-UPU^*]\\
&=&0
\end{eqnarray*}

Indeed, by positivity this gives $PU-UP=0$, as desired.

\medskip

(6) With these preliminaries in hand, let us decompose the algebra $End(u)$ as in Theorem 13.12, by using the decomposition $1=p_1+\ldots+p_k$ into minimal projections. If we denote by $u_i\subset u$ the subrepresentation coming from the vector space $V_i=Im(p_i)$, then we obtain in this way a decomposition $u=u_1+\ldots+u_k$, as in the statement.
\end{proof}

In order to formulate our second Peter-Weyl theorem, we need to talk about coefficients, and smoothness. Things here are quite tricky, and we can proceed as follows:

\index{smooth representation}
\index{space of coefficients}

\begin{definition}
Given a closed subgroup $G\subset U_N$, and a unitary representation $v:G\to U_M$, the space of coefficients of this representation is:
$$C_v=\left\{f\circ v\Big|f\in M_M(\mathbb C)^*\right\}$$
In other words, by delinearizing, $C_\nu\subset C(G)$ is the following linear space:
$$C_v=span\Big[g\to(v_g)_{ij}\Big]$$
We say that $v$ is smooth if its matrix coefficients $g\to(v_g)_{ij}$ appear as polynomials in the standard matrix coordinates $g\to g_{ij}$, and their conjugates $g\to\overline{g}_{ij}$.
\end{definition}

As a basic example of coefficient we have, besides the matrix coefficients $g\to(v_g)_{ij}$, the character, which appears as the diagonal sum of these coefficients:
$$\chi_v(g)=\sum_i(v_g)_{ii}$$

Regarding the notion of smoothness, things are quite tricky here, the idea being that any closed subgroup $G\subset U_N$ can be shown to be a Lie group, and that, with this result in hand, a representation $v:G\to U_M$ is smooth precisely when the condition on coefficients from the above definition is satisfied. All this is quite technical, and we will not get into it. We will simply use Definition 13.14 as such, and further comment on this later on. 

\bigskip

Here is now our second Peter-Weyl theorem, complementing Theorem 13.13:

\index{Peter-Weyl}

\begin{theorem}[PW2]
Given a closed subgroup $G\subset_uU_N$, any of its irreducible smooth representations 
$$v:G\to U_M$$
appears inside a tensor product of the fundamental representation $u$ and its adjoint $\bar{u}$.
\end{theorem}

\begin{proof}
In order to prove the result, we will use the following three elementary facts, regarding the spaces of coefficients introduced above:

\medskip

(1) The construction $v\to C_v$ is functorial, in the sense that it maps subrepresentations into linear subspaces. This is indeed something which is routine to check.

\medskip

(2) Our smoothness assumption on $v:G\to U_M$, as formulated in Definition 13.14, means that we have an inclusion of linear spaces as follows:
$$C_v\subset<g_{ij}>$$

(3) By definition of the Peter-Weyl representations, as arbitrary tensor products between the fundamental representation $u$ and its conjugate $\bar{u}$, we have:
$$<g_{ij}>=\sum_kC_{u^{\otimes k}}$$

(4) Now by putting together the observations (2,3) we conclude that we must have an inclusion as follows, for certain exponents $k_1,\ldots,k_p$:
$$C_v\subset C_{u^{\otimes k_1}\oplus\ldots\oplus u^{\otimes k_p}}$$

By using now the functoriality result from (1), we deduce from this that we have an inclusion of representations, as follows:
$$v\subset u^{\otimes k_1}\oplus\ldots\oplus u^{\otimes k_p}$$

Together with Theorem 13.13, this leads to the conclusion in the statement.
\end{proof}

As a conclusion to what we have so far, the problem to be solved is that of splitting the Peter-Weyl representations into sums of irreducible representations.

\section*{13c. Haar integration}

In order to further advance, and complete the Peter-Weyl theory, we need to talk about integration over $G$. In the finite group case the situation is trivial, as follows:

\begin{proposition}
Any finite group $G$ has a unique probability measure which is invariant under left and right translations,
$$\mu(E)=\mu(gE)=\mu(Eg)$$
and this is the normalized counting measure on $G$, given by $\mu(E)=|E|/|G|$.
\end{proposition}

\begin{proof}
The uniformity condition in the statement gives, with $E=\{h\}$:
$$\mu\{h\}=\mu\{gh\}=\mu\{hg\}$$

Thus $\mu$ must be the usual counting measure, normalized as to have mass 1.
\end{proof}

In the continuous group case now, the simplest examples, to be studied first, are the compact abelian groups. Here things are standard again, as follows:

\begin{theorem}
Given a compact abelian group $G$, with dual group denoted $\Gamma=\widehat{G}$, we have an isomorphism of commutative algebras
$$C(G)\simeq C^*(\Gamma)$$
and via this isomorphism, the functional defined by linearity and the following formula, 
$$\int_Gg=\delta_{g1}$$
for any $g\in\Gamma$, is the integration with respect to the unique uniform measure on $G$.
\end{theorem}

\begin{proof}
We can indeed apply the Gelfand theorem, from chapter 8, to the group algebra $C^*(\Gamma)$, which is commutative, and this gives all the results.
\end{proof}

Summarizing, we have results in the finite case, and in the compact abelian case. With the remark that the proof in the compact abelian case was quite brief, but this result, coming as an illustration for more general things to follow, is not crucial for us. 

\bigskip

Let us discuss now the construction of the uniform probability measure in general. This is something quite technical, the idea being that the uniform measure $\mu$ over $G$ can be constructed by starting with an arbitrary probability measure $\nu$, and setting:
$$\mu=\lim_{n\to\infty}\frac{1}{n}\sum_{k=1}^n\nu^{*k}$$

Thus, our next task will be that of proving this result. It is convenient, for this purpose, to work with the integration functionals with respect to the various measures on $G$, instead of the measures themselves. Let us begin with the following key result:

\begin{proposition}
Given a unital positive linear form $\varphi:C(G)\to\mathbb C$, the limit
$$\int_\varphi f=\lim_{n\to\infty}\frac{1}{n}\sum_{k=1}^n\varphi^{*k}(f)$$
exists, and for a coefficient of a representation $f=(\tau\otimes id)v$ we have
$$\int_\varphi f=\tau(P)$$
where $P$ is the orthogonal projection onto the $1$-eigenspace of $(id\otimes\varphi)v$.
\end{proposition}

\begin{proof}
By linearity it is enough to prove the first assertion for functions of the following type, where $v$ is a Peter-Weyl representation, and $\tau$ is a linear form:
$$f=(\tau\otimes id)v$$

Thus we are led into the second assertion, and more precisely we can have the whole result proved if we can establish the following formula, with $f=(\tau\otimes id)v$:
$$\lim_{n\to\infty}\frac{1}{n}\sum_{k=1}^n\varphi^{*k}(f)=\tau(P)$$

In order to prove this latter formula, observe that we have:
$$\varphi^{*k}(f)
=(\tau\otimes\varphi^{*k})v
=\tau((id\otimes\varphi^{*k})v)$$

Let us set $M=(id\otimes\varphi)v$. In terms of this matrix, we have:
$$((id\otimes\varphi^{*k})v)_{i_0i_{k+1}}
=\sum_{i_1\ldots i_k}M_{i_0i_1}\ldots M_{i_ki_{k+1}}
=(M^k)_{i_0i_{k+1}}$$

Thus we have the following formula, for any $k\in\mathbb N$:
$$(id\otimes\varphi^{*k})v=M^k$$

It follows that our Ces\`aro limit is given by the following formula:
$$\lim_{n\to\infty}\frac{1}{n}\sum_{k=1}^n\varphi^{*k}(f)
=\lim_{n\to\infty}\frac{1}{n}\sum_{k=1}^n\tau(M^k)
=\tau\left(\lim_{n\to\infty}\frac{1}{n}\sum_{k=1}^nM^k\right)$$

Now since $v$ is unitary we have $||v||=1$, and so $||M||\leq1$. Thus the last Ces\`aro limit converges, and equals the orthogonal projection onto the $1$-eigenspace of $M$:
$$\lim_{n\to\infty}\frac{1}{n}\sum_{k=1}^nM^k=P$$

Thus our initial Ces\`aro limit converges as well, to $\tau(P)$, as desired. 
\end{proof}

The point now is that when the linear form $\varphi\in C(G)^*$ from the above result is chosen to be faithful, we obtain the following finer result:

\begin{proposition}
Given a faithful unital linear form $\varphi\in C(G)^*$, the limit
$$\int_\varphi f=\lim_{n\to\infty}\frac{1}{n}\sum_{k=1}^n\varphi^{*k}(f)$$
exists, and is independent of $\varphi$, given on coefficients of representations by
$$\left(id\otimes\int_\varphi\right)v=P$$
where $P$ is the orthogonal projection onto the space $Fix(v)=\left\{\xi\in\mathbb C^n\big|v\xi=\xi\right\}$.
\end{proposition}

\begin{proof}
In view of Proposition 13.18, it remains to prove that when $\varphi$ is faithful, the $1$-eigenspace of the matrix $M=(id\otimes\varphi)v$ equals the space $Fix(v)$.

\medskip

``$\supset$'' This is clear, and for any $\varphi$, because we have the following implication:
$$v\xi=\xi\implies M\xi=\xi$$

``$\subset$'' Here we must prove that, when $\varphi$ is faithful, we have:
$$M\xi=\xi\implies v\xi=\xi$$

For this purpose, assume that we have $M\xi=\xi$, and consider the following function:
$$f=\sum_i\left(\sum_jv_{ij}\xi_j-\xi_i\right)\left(\sum_kv_{ik}\xi_k-\xi_i\right)^*$$

We must prove that we have $f=0$. Since $v$ is unitary, we have:
\begin{eqnarray*}
f
&=&\sum_{ijk}v_{ij}v_{ik}^*\xi_j\bar{\xi}_k-\frac{1}{N}v_{ij}\xi_j\bar{\xi}_i-\frac{1}{N}v_{ik}^*\xi_i\bar{\xi}_k+\frac{1}{N^2}\xi_i\bar{\xi}_i\\
&=&\sum_j|\xi_j|^2-\sum_{ij}v_{ij}\xi_j\bar{\xi}_i-\sum_{ik}v_{ik}^*\xi_i\bar{\xi}_k+\sum_i|\xi_i|^2\\
&=&||\xi||^2-<v\xi,\xi>-\overline{<v\xi,\xi>}+||\xi||^2\\
&=&2(||\xi||^2-Re(<v\xi,\xi>))
\end{eqnarray*}

By using now our assumption $M\xi=\xi$, we obtain from this:
\begin{eqnarray*}
\varphi(f)
&=&2\varphi(||\xi||^2-Re(<v\xi,\xi>))\\
&=&2(||\xi||^2-Re(<M\xi,\xi>))\\
&=&2(||\xi||^2-||\xi||^2)\\
&=&0
\end{eqnarray*}

Now since $\varphi$ is faithful, this gives $f=0$, and so $v\xi=\xi$, as claimed.
\end{proof}

We can now formulate a main result about Haar integration, as follows:

\index{Haar measure}
\index{Haar integration}

\begin{theorem}
Any compact group $G$ has a unique Haar integration, which can be constructed by starting with any faithful positive unital state $\varphi\in C(G)^*$, and setting:
$$\int_G=\lim_{n\to\infty}\frac{1}{n}\sum_{k=1}^n\varphi^{*k}$$
Moreover, for any representation $v$ we have the formula
$$\left(id\otimes\int_G\right)v=P$$
where $P$ is the orthogonal projection onto $Fix(v)=\left\{\xi\in\mathbb C^n\big|v\xi=\xi\right\}$.
\end{theorem}

\begin{proof}
We can prove this from what we have, in several steps, as follows:

\medskip

(1) Let us first go back to the general context of Proposition 13.18. Since convolving one more time with $\varphi$ will not change the Ces\`aro limit appearing there, the functional $\int_\varphi\in C(G)^*$ constructed there has the following invariance property:
$$\int_\varphi*\,\varphi=\varphi*\int_\varphi=\int_\varphi$$

In the case where $\varphi$ is assumed to be faithful, as in Proposition 13.19, our claim is that we have the following formula, valid this time for any $\psi\in C(G)^*$:
$$\int_\varphi*\,\psi=\psi*\int_\varphi=\psi(1)\int_\varphi$$

Moreover, it is enough to prove this formula on a coefficient of a representation:
$$f=(\tau\otimes id)v$$

(2) In order to do so, consider the following two matrices:
$$P=\left(id\otimes\int_\varphi\right)v\quad,\quad 
Q=(id\otimes\psi)v$$

We have then the following two computations, involving these matrices:
$$\left(\int_\varphi*\,\psi\right)f
=\left(\tau\otimes\int_\varphi\otimes\,\psi\right)(v_{12}v_{13})
=\tau(PQ)$$
$$\left(\psi*\int_\varphi\right)f
=\left(\tau\otimes\psi\otimes\int_\varphi\right)(v_{12}v_{13})
=\tau(QP)$$

Also, regarding the term on the right in our formula in (1), this is given by:
$$\psi(1)\int_\varphi f=\psi(1)\tau(P)$$

We conclude from all this that our claim is equivalent to the following equality:
$$PQ=QP=\psi(1)P$$

(3) But this latter equality holds indeed, coming from the fact, that we know from Proposition 13.19, that $P=(id\otimes\int_\varphi)v$ equals the orthogonal projection onto $Fix(v)$. Thus, we have proved our claim in (1), namely that the following formula holds: 
$$\int_\varphi*\,\psi=\psi*\int_\varphi=\psi(1)\int_\varphi$$

(4) In order to finish now, it is convenient to introduce the following abstract operation, on the continuous functions $f,f':C(G)\to\mathbb C$ on our group:
$$\Delta(f\otimes f')(g\otimes h)=f(g)f'(h)$$

With this convention, the formula that we established above can be written as:
$$\psi\left(\int_\varphi\otimes\, id\right)\Delta
=\psi\left(id\otimes\int_\varphi\right)\Delta
=\psi\int_\varphi(.)1$$

This formula being true for any $\psi\in C(G)^*$, we can simply delete $\psi$. We conclude that the following invariance formula holds indeed, with $\int_G=\int_\varphi$:
$$\left(\int_G\otimes\, id\right)\Delta=\left(id\otimes\int_G\right)\Delta=\int_G(.)1$$

But this is exactly the left and right invariance formula we were looking for.

\medskip

(5) Finally, in order to prove the uniqueness assertion, assuming that we have two invariant integrals $\int_G,\int_G'$, we have, according to the above invariance formula:
$$\left(\int_G\otimes\int_G'\right)\Delta
=\left(\int_G'\otimes\int_G\right)\Delta
=\int_G(.)1
=\int_G'(.)1$$

Thus we have $\int_G=\int_G'$, and this finishes the proof.
\end{proof}

Summarizing, we can now integrate over $G$. As a first application, we have:

\index{main character}
\index{moments}

\begin{theorem}
Given a compact group $G$, we have the following formula, valid for any unitary group representation $v:G\to U_M$:
$$\int_G\chi_v=\dim(Fix(v))$$
In particular, in the unitary matrix group case, $G\subset_uU_N$, the moments of the main character $\chi=\chi_u$ are given by the following formula:
$$\int_G\chi^k=\dim(Fix(u^{\otimes k}))$$
Thus, knowing the law of $\chi$ is the same as knowing the dimensions on the right.
\end{theorem}

\begin{proof}
We have three assertions here, the idea being as follows:

\medskip

(1) Given a unitary representation $v:G\to U_M$ as in the statement, its character $\chi_v$ is a coefficient, so we can use the integration formula for coefficients in Theorem 13.20. If we denote by $P$ the projection onto $Fix(v)$, that formula gives, as desired:
\begin{eqnarray*}
\int_G\chi_v
&=&Tr(P)\\
&=&\dim(Im(P))\\
&=&dim(Fix(v))
\end{eqnarray*}

(2) This comes from (1) applied to the Peter-Weyl representations, as follows:
\begin{eqnarray*}
\int_G\chi^k
&=&\int_G\chi_u^k\\
&=&\int_G\chi_{u^{\otimes k}}\\
&=&\dim(Fix(u^{\otimes k}))
\end{eqnarray*}

(3) This follows from (2), and from the standard fact, which follows from definitions, that a probability measure is uniquely determined by its moments.
\end{proof}

As a key remark now, the integration formula in Theorem 13.20 allows the computation for the truncated characters too, because these truncated characters are coefficients as well. To be more precise, all the probabilistic questions about $G$, regarding characters, or truncated characters, or more complicated variables, require a good knowledge of the integration over $G$, and more precisely, of the various polynomial integrals over $G$:

\index{polynomial integrals}

\begin{definition}
Given a closed subgroup $G\subset U_N$, the quantities
$$I_k=\int_Gg_{i_1j_1}^{e_1}\ldots g_{i_kj_k}^{e_k}\,dg$$
depending on a colored integer $k=e_1\ldots e_k$, are called polynomial integrals over $G$.
\end{definition}

As a first observation, the knowledge of these integrals is the same as the knowledge of the integration functional over $G$. Indeed, since the coordinate functions $g\to g_{ij}$ separate the points of $G$, we can apply the Stone-Weierstrass theorem, and we obtain:
$$C(G)=<g_{ij}>$$

Thus, by linearity, the computation of any functional $f:C(G)\to\mathbb C$, and in particular of the integration functional, reduces to the computation of this functional on the polynomials of the coordinate functions $g\to g_{ij}$ and their conjugates $g\to\bar{g}_{ij}$.

\bigskip

By using now Peter-Weyl theory, everything reduces to algebra, as follows:

\index{Weingarten formula}
\index{Gram matrix}
\index{Weingarten matrix}

\begin{theorem}
The Haar integration over a closed subgroup $G\subset_uU_N$ is given on the dense subalgebra of smooth functions by the Weingarten formula
$$\int_Gg_{i_1j_1}^{e_1}\ldots g_{i_kj_k}^{e_k}\,dg=\sum_{\pi,\sigma\in D_k}\delta_\pi(i)\delta_\sigma(j)W_k(\pi,\sigma)$$
valid for any colored integer $k=e_1\ldots e_k$ and any multi-indices $i,j$, where $D_k$ is a linear basis of $Fix(u^{\otimes k})$, the associated generalized Kronecker symbols are given by
$$\delta_\pi(i)=<\pi,e_{i_1}\otimes\ldots\otimes e_{i_k}>$$
and $W_k=G_k^{-1}$ is the inverse of the Gram matrix, $G_k(\pi,\sigma)=<\pi,\sigma>$.
\end{theorem}

\begin{proof}
We know from Peter-Weyl theory that the integrals in the statement form altogether the orthogonal projection $P^k$ onto the following space:
$$Fix(u^{\otimes k})=span(D_k)$$

Consider now the following linear map, with $D_k=\{\xi_k\}$ being as in the statement:
$$E(x)=\sum_{\pi\in D_k}<x,\xi_\pi>\xi_\pi$$

By a standard linear algebra computation, it follows that we have $P=WE$, where $W$ is the inverse of the restriction of $E$ to the following space:
$$K=span\left(T_\pi\Big|\pi\in D_k\right)$$

But this restriction is precisely the linear map given by the matrix $G_k$, and so $W$ itself is the linear map given by the matrix $W_k$, and this gives the result.
\end{proof}

We will be back to this in chapter 16 below, with some concrete applications.

\section*{13d. More Peter-Weyl}

In order to further develop now the Peter-Weyl theory, which is something very useful, we will need the following result, which is of independent interest:

\index{Frobenius isomorphism}

\begin{proposition}
We have a Frobenius type isomorphism
$$Hom(v,w)\simeq Fix(v\otimes\bar{w})$$
valid for any two representations $v,w$.
\end{proposition}

\begin{proof}
According to the definitions, we have the following equivalences:
\begin{eqnarray*}
T\in Hom(v,w)
&\iff&Tv=wT\\
&\iff&\sum_jT_{aj}v_{ji}=\sum_bw_{ab}T_{bi},\forall a,i
\end{eqnarray*}

On the other hand, we have as well the following equivalences:
\begin{eqnarray*}
T\in Fix(v\otimes\bar{w})
&\iff&(v\otimes\bar{w})T=\xi\\
&\iff&\sum_{jb}v_{ij}w_{ab}^*T_{bj}=T_{ai}\forall a,i
\end{eqnarray*}

With these formulae in hand, both inclusions follow from the unitarity of $v,w$.
\end{proof}

We can now formulate our third Peter-Weyl theorem, as follows:

\index{Peter-Weyl}

\begin{theorem}[PW3]
The norm dense $*$-subalgebra 
$$\mathcal C(G)\subset C(G)$$
generated by the coefficients of the fundamental representation decomposes as 
$$\mathcal C(G)=\bigoplus_{v\in Irr(G)}M_{\dim(v)}(\mathbb C)$$
with the summands being pairwise orthogonal with respect to the scalar product
$$<a,b>=\int_Gab^*$$
where $\int_G$ is the Haar integration over $G$.
\end{theorem}

\begin{proof}
By combining the previous two Peter-Weyl results, we deduce that we have a linear space decomposition as follows:
$$\mathcal C(G)
=\sum_{v\in Irr(G)}C_v
=\sum_{v\in Irr(G)}M_{\dim(v)}(\mathbb C)$$

Thus, in order to conclude, it is enough to prove that for any two irreducible corepresentations $v,w\in Irr(A)$, the corresponding spaces of coefficients are orthogonal:
$$v\not\sim w\implies C_v\perp C_w$$ 

But this follows from Theorem 13.20, via Proposition 13.24. Let us set indeed:
$$P_{ia,jb}=\int_Gv_{ij}w_{ab}^*$$

Then $P$ is the orthogonal projection onto the following vector space:
$$Fix(v\otimes\bar{w})
\simeq Hom(v,w)
=\{0\}$$

Thus we have $P=0$, and this gives the result.
\end{proof}

Finally, we have the following result, completing the Peter-Weyl theory:

\index{Peter-Weyl}
\index{central function}
\index{algebra of characters}

\begin{theorem}[PW4]
The characters of irreducible representations belong to
$$\mathcal C(G)_{central}=\left\{f\in\mathcal C(G)\Big|f(gh)=f(hg),\forall g,h\in G\right\}$$
called algebra of smooth central functions on $G$, and form an orthonormal basis of it.
\end{theorem}

\begin{proof}
We have several things to be proved, the idea being as follows:

\medskip

(1) Observe first that $\mathcal C(G)_{central}$ is indeed an algebra, which contains all the characters. Conversely, consider a function $f\in\mathcal C(G)$, written as follows:
$$f=\sum_{v\in Irr(G)}f_v$$

The condition $f\in\mathcal C(G)_{central}$ states then that for any $v\in Irr(G)$, we must have:
$$f_v\in\mathcal C(G)_{central}$$

But this means precisely that the coefficient $f_v$ must be a scalar multiple of $\chi_v$, and so the characters form a basis of $\mathcal C(G)_{central}$, as stated. 

\medskip

(2) The fact that we have an orthogonal basis follows from Theorem 13.25. 

\medskip

(3) As for the fact that the characters have norm 1, this follows from:
\begin{eqnarray*}
\int_G\chi_v\chi_v^*
&=&\sum_{ij}\int_Gv_{ii}v_{jj}^*\\
&=&\sum_i\frac{1}{N}\\
&=&1
\end{eqnarray*}

Here we have used the fact, coming from Theorem 13.25, that the integrals $\int_Gv_{ij}v_{kl}^*$ form the orthogonal projection onto the following vector space: 
$$Fix(v\otimes\bar{v})
\simeq End(v)
=\mathbb C1$$

Thus, the proof of our theorem is now complete.
\end{proof}

As a key observation now, complementing Theorem 13.26, observe that a function $f:G\to\mathbb C$ is central, in the sense that it satisfies $f(gh)=f(hg)$, precisely when it satisfies the following condition, saying that it must be constant on conjugacy classes:
$$f(ghg^{-1})=f(h),\forall g,h\in G$$

Now the point is that this makes the algebra of central functions something quite easy to compute, via standard algebra, and this puts us on the right track for computing $Irr(G)$. Or at least, this is how the theory goes, because there are many tricks too. 

\bigskip

As a basic illustration for this method, which clarifies some previous considerations from chapter 9, in relation with our study there of the finite abelian groups, we have:

\begin{theorem}
For a finite abelian group $G$ the irreducible representations are all $1$-dimensional, equal to their own characters,
$$\chi:G\to\mathbb T$$
and these characters form the dual discrete abelian group $\widehat{G}$.
\end{theorem}

\begin{proof}
This comes indeed from the Peter-Weyl theory, as follows:

\medskip

(1) Since our group $G$ was assumed to be abelian, any function $f:G\to\mathbb C$ is obviously central, so the algebra of central functions is $C(G)$ itself:
$$C(G)_{central}=C(G)$$

(2) Thus the decomposition of $C(G)$ from Theorem 13.25 reduces in this case to the decomposition of $C(G)_{central}$ from Theorem 13.26, and in particular, the irreducible representations $u\in Irr(G)$ must be all 1-dimensional, equal to their own characters $\chi_u$.

\medskip

(3) Finally, the last assertion is something that we know well from chapter 9, and with the extra comment that we have in fact an isomorphism $\widehat{G}\simeq G$, coming from the structure theorem for the finite abelian groups, as explained there. 

\medskip

(4) As a final comment on this, observe that $\widehat{G}\simeq G$, or the structure theorem for the finite abelian groups, do not come from Peter-Weyl for the abelian groups, whose conclusions reduce to what is said in the statement. Thus, although Peter-Weyl for the finite abelian groups does part of the job that we did in chapter 9, this is not everything, and our arithmetic work there remains something needed, going beyond Peter-Weyl.
\end{proof}

Getting now to the non-abelian case, things here can be quite complicated. For the simplest non-abelian group that we know, namely $S_3=D_3$, the result is as follows:

\begin{theorem}
The group $S_3=D_3$ has $3$ irreducible representations, namely:
\begin{enumerate}
\item The trivial representation, $g\to 1$.

\item The signature representation, $g\to\varepsilon(g)$.

\item The $2D$ representation $u-1$, with $u$ being the standard $3D$ representation.
\end{enumerate}
\end{theorem}

\begin{proof}
We certainly have the representations in (1) and (2), which are obviously irreducible, and non-equivalent. Now let us look at the 3D representation:
$$u:[S_3=D_3]\subset O_3\subset U_3$$

Since this representation appears via the permutation matrices, which sum up to 1 on each row, we conclude that the all-one vector is fixed by this representation:
$$\begin{pmatrix}1\\ 1\\1\end{pmatrix}\in Fix(u)$$

Thus, we can consider the following representation, which is 2-dimensional:
$$v=u-1$$

And we can stop here, because our group being non-abelian, and of order 6, a quick look at Theorem 13.25 shows that the decomposition there must come from:
$$6=1+1+4$$

Thus, $1,\varepsilon,v$ are indeed the irreducible representations, as stated.
\end{proof}

Regarding now more complicated groups, with a bit more work the ideas in the above proof extend to all dihedral groups $D_N$. As for the symmetric groups $S_N$, the situation here is more complicated. We will leave some study and learning here as an exercise

\section*{13e. Exercises}

There has been a lot of theory on this chapter, and as exercises, we will have some more theory, namely an introduction to quantum groups. Let us start with:

\begin{exercise}
Given a finite group $G$, setting $A=C(G)$, prove that the maps
$$\Delta:A\to A\otimes A\quad,\quad 
\varepsilon:A\to\mathbb C\quad,\quad 
S:A\to A$$
which are transpose to the multiplication $m:G\times G\to G$, unit $u:\{.\}\to G$ and inverse map $i:G\to G$, are subject to the following conditions,
$$(\varepsilon\otimes id)\Delta=(id\otimes\varepsilon)\Delta=id$$
$$m(S\otimes id)\Delta=m(id\otimes S)\Delta=\varepsilon(.)1$$
in usual tensor product notation, along with the extra condition $S^2=id$. 
\end{exercise}

This does not look difficult, with the conditions in the statement reminding the usual group axioms, satisfied by $m,u,i$. Up to you to prove this now, with full details.

\begin{exercise}
Given a finite group $H$, setting $A=C^*(H)$, prove that the maps
$$\Delta:A\to A\otimes A\quad,\quad 
\quad,\quad\varepsilon:A\to\mathbb C
\quad,\quad S:A\to A^{opp}$$
given by the formulae $\Delta(g)=g\otimes g$, $\varepsilon(g)=1$, $S(g)=g^{-1}$ and linearity, are subject to the same conditions as above, including the extra condition $S^2=id$. 
\end{exercise}

As before with the previous exercise, this does not look very difficult, with most likely only some elementary algebraic computations being involved.

\begin{exercise}
Let us call finite Hopf algebra a finite dimensional $C^*$-algebra, with maps as follows, called comultiplication, counit and antipode, 
$$\Delta:A\to A\otimes A\quad,\quad 
\varepsilon:A\to\mathbb C\quad,\quad 
S:A\to A^{opp}$$
satisfying the conditions found above. Prove that if $G,H$ are finite abelian groups, dual to each other, we have an isomorphism of finite Hopf algebras as follows:
$$C(G)=C^*(H)$$
Afterwards, based on this, formally write any finite Hopf algebra as 
$$A=C(G)=C^*(H)$$
and call $G,H$ finite quantum groups, dual to each other.
\end{exercise}

Here the thing to be done, namely to establish the identification in the statement, looks like something quite routine, related to many things that we already know. As for the last part, there is nothing to be done here, just enjoying that definition.

\chapter{Tannakian duality}

\section*{14a. Tensor categories}

We have seen that the representations of a closed subgroup $G\subset U_N$ are subject to a number of non-trivial results, collectively known as Peter-Weyl theory. To be more precise, the main ideas of Peter-Weyl theory were as follows:

\medskip

\begin{enumerate}

\item The representations of $G$ split as sums of irreducibles, and the irreducibles can be found inside the tensor products $u^{\otimes k}$ between the fundamental representation $u:G\subset U_N$ and its adjoint $\bar{u}:G\subset U_N$, called Peter-Weyl representations.

\medskip

\item The main problem is therefore that of splitting the various Peter-Weyl representations $u^{\otimes k}$ into irreducibles. Technically speaking, this leads to the question of explicitly computing the corresponding fixed point spaces $Fix(u^{\otimes k})$.

\medskip

\item From a probabilistic perspective, in connection with characters and truncated characters, which require the explicit knowledge of $\int_G$, we are led into the same fundamental question, namely the computation of the spaces $Fix(u^{\otimes k})$.
\end{enumerate}

\medskip

Summarizing, no matter what we want to do with $G$, we must compute the spaces $Fix(u^{\otimes k})$. As a first idea now, it is technically convenient to slightly enlarge the class of spaces to be computed, by talking about Tannakian categories, as follows:

\index{Tannakian category}
\index{Peter-Weyl representations}

\begin{definition}
The Tannakian category associated to a closed subgroup $G\subset_uU_N$ is the collection $C=(C(k,l))$ of vector spaces
$$C(k,l)=Hom(u^{\otimes k},u^{\otimes l})$$
where the representations $u^{\otimes k}$ with $k=\circ\bullet\bullet\circ\ldots$ colored integer, defined by
$$u^{\otimes\emptyset}=1\quad,\quad 
u^{\otimes\circ}=u\quad,\quad
u^{\otimes\bullet}=\bar{u}$$
and multiplicativity, $u^{\otimes kl}=u^{\otimes k}\otimes u^{\otimes l}$, are the Peter-Weyl representations.
\end{definition}

Here are a few examples of such representations, namely those coming from the colored integers of length 2, to be often used in what follows:
$$u^{\otimes\circ\circ}=u\otimes u\quad,\quad 
u^{\otimes\circ\bullet}=u\otimes\bar{u}$$
$$u^{\otimes\bullet\circ}=\bar{u}\otimes u\quad,\quad 
u^{\otimes\bullet\bullet}=\bar{u}\otimes\bar{u}$$

As a first observation, the knowledge of the Tannakian category is more or less the same thing as the knowledge of the fixed point spaces, which appear as:
$$Fix(u^{\otimes k})=C(0,k)$$

Indeed, these latter spaces fully determine all the spaces $C(k,l)$, because of the Frobenius isomorphisms, which for the Peter-Weyl representations read:
\begin{eqnarray*}
C(k,l)
&=&Hom(u^{\otimes k},u^{\otimes l})\\
&\simeq&Hom(1,\bar{u}^{\otimes k}\otimes u^{\otimes l})\\
&=&Hom(1,u^{\otimes\bar{k}l})\\
&=&Fix(u^{\otimes\bar{k}l})
\end{eqnarray*}

In order to get started now, let us make a summary of what we have so far, regarding these spaces $C(k,l)$, coming from the general theory developed in chapter 13. In order to formulate our result, let us start with an abstract definition, as follows:

\index{tensor category}

\begin{definition}
Let $H$ be a finite dimensional Hilbert space. A tensor category over $H$ is a collection $C=(C(k,l))$ of linear spaces 
$$C(k,l)\subset\mathcal L(H^{\otimes k},H^{\otimes l})$$
satisfying the following conditions:
\begin{enumerate}
\item $S,T\in C$ implies $S\otimes T\in C$.

\item If $S,T\in C$ are composable, then $ST\in C$.

\item $T\in C$ implies $T^*\in C$.

\item Each $C(k,k)$ contains the identity operator.

\item $C(\emptyset,k)$ with $k=\circ\bullet,\bullet\circ$ contain the operator $R:1\to\sum_ie_i\otimes e_i$.

\item $C(kl,lk)$ with $k,l=\circ,\bullet$ contain the flip operator $\Sigma:a\otimes b\to b\otimes a$.
\end{enumerate}
\end{definition}

Here the tensor powers $H^{\otimes k}$, which are Hilbert spaces depending on a colored integer $k=\circ\bullet\bullet\circ\ldots\,$, are defined by the following formulae, and multiplicativity:
$$H^{\otimes\emptyset}=\mathbb C\quad,\quad 
H^{\otimes\circ}=H\quad,\quad
H^{\otimes\bullet}=\bar{H}\simeq H$$

With these conventions, we have the following result, summarizing our knowledge on the subject, coming from the results from the previous chapter:

\begin{theorem}
For a closed subgroup $G\subset_uU_N$, the associated Tannakian category
$$C(k,l)=Hom(u^{\otimes k},u^{\otimes l})$$
is a tensor category over the Hilbert space $H=\mathbb C^N$.
\end{theorem}

\begin{proof}
We know that the fundamental representation $u$ acts on the Hilbert space $H=\mathbb C^N$, and that its conjugate $\bar{u}$ acts on the Hilbert space $\bar{H}=\mathbb C^N$. Now by multiplicativity we conclude that any Peter-Weyl representation $u^{\otimes k}$ acts on the Hilbert space $H^{\otimes k}$, so that we have embeddings as in Definition 14.2, as follows:
$$C(k,l)\subset\mathcal L(H^{\otimes k},H^{\otimes l})$$

Regarding now the fact that the axioms (1-6) in Definition 14.2 are indeed satisfied, this is something that we basically already know, as follows:

\medskip

(1,2,3) These results follow from definitions, and were explained in chapter 13.

\medskip

(4) This is something trivial, coming from definitions.

\medskip

(5) This follows from the fact that each element $g\in G$ is a unitary, which can be reformulated as follows, with $R:1\to\sum_ie_i\otimes e_i$ being the map in Definition 14.2:
$$R\in Hom(1,g\otimes\bar{g})\quad,\quad 
R\in Hom(1,\bar{g}\otimes g)$$

Indeed, given an arbitrary matrix $g\in M_N(\mathbb C)$, we have the following computation:
\begin{eqnarray*}
(g\otimes\bar{g})(R(1)\otimes1)
&=&\left(\sum_{ijkl}e_{ij}\otimes e_{kl}\otimes g_{ij}\bar{g}_{kl}\right)\left(\sum_ae_a\otimes e_a\otimes 1\right)\\
&=&\sum_{ika}e_i\otimes e_k\otimes g_{ia}\bar{g}_{ka}^*\\
&=&\sum_{ik}e_i\otimes e_k\otimes(gg^*)_{ik}
\end{eqnarray*}

We conclude from this that we have the following equivalence:
\begin{eqnarray*}
R\in Hom(1,g\otimes\bar{g})
&\iff&gg^*=1
\end{eqnarray*}

By replacing $g$ with its conjugate matrix $\bar{g}$, we have as well:
$$R\in Hom(1,\bar{g}\otimes g)\iff\bar{g}g^t=1$$

Thus, the two intertwining conditions in Definition 14.2 (5) are both equivalent to the fact that $g$ is unitary, and so these conditions are indeed satisfied, as desired.

\medskip

(6) This is again something elementary, coming from the fact that the various matrix coefficients $g\to g_{ij}$ and their complex conjugates $g\to\bar{g}_{ij}$ commute with each other. To be more precise, with $\Sigma:a\otimes b\to b\otimes a$ being the flip operator, we have:
\begin{eqnarray*}
(g\otimes h)(\Sigma\otimes id)(e_a\otimes e_b\otimes 1)
&=&\left(\sum_{ijkl}e_{ij}\otimes e_{kl}\otimes g_{ij}h_{kl}\right)(e_b\otimes e_a\otimes1)\\
&=&\sum_{ik}e_i\otimes e_k\otimes g_{ib}h_{ka}
\end{eqnarray*}

On the other hand, we have as well the following computation:
\begin{eqnarray*}
(\Sigma\otimes id)(h\otimes g)(e_a\otimes e_b\otimes 1)
&=&(\Sigma\otimes id)\left(\sum_{ijkl}e_{ij}\otimes e_{kl}\otimes h_{ij}g_{kl}\right)(e_a\otimes e_b\otimes1)\\
&=&(\Sigma\otimes id)\left(\sum_{ik}e_i\otimes e_k\otimes h_{ia}g_{kb}\right)\\
&=&\sum_{ik}e_k\otimes e_i\otimes h_{ia}g_{kb}\\
&=&\sum_{ik}e_i\otimes e_k\otimes h_{ka}g_{ib}
\end{eqnarray*}

Now since functions commute, $g_{ib}h_{ka}=h_{ka}g_{ib}$, this gives the result.
\end{proof}

Quite remarkably, we have the following result, coming from Peter-Weyl:

\begin{theorem}
Given a compact subgroup $G\subset U_N$, we have
$$G=\left\{g\in U_N\Big|Tg^{\otimes k}=g^{\otimes l}T,\forall k,l,\forall T\in C(k,l)\right\}$$
where $C=(C(k,l))$ is the associated Tannakian category.
\end{theorem}

\begin{proof}
This is something quite standard, the idea being as follows:

\medskip

(1) Consider the set of matrices constructed in the statement, namely:
$$\widetilde{G}=\left\{g\in U_N\Big|Tg^{\otimes k}=g^{\otimes l}T,\forall k,l,\forall T\in C(k,l)\right\}$$

Our first claim is that $\widetilde{G}$ is a group. Indeed, assuming $g,h\in\widetilde{G}$, we have $gh\in\widetilde{G}$, due to the following computation, valid for any $k,l$ and any $T\in C(k,l)$:
\begin{eqnarray*}
T(gh)^{\otimes k}
&=&Tg^{\otimes k}h^{\otimes k}\\
&=&g^{\otimes l}Th^{\otimes k}\\
&=&g^{\otimes l}h^{\otimes l}T\\
&=&(gh)^{\otimes l}T
\end{eqnarray*}

Also, we have $1\in\widetilde{G}$, trivially. Finally, assuming $g\in\widetilde{G}$, we have:
\begin{eqnarray*}
T(g^{-1})^{\otimes k}
&=&(g^{-1})^{\otimes l}[g^{\otimes l}T](g^{-1})^{\otimes k}\\
&=&(g^{-1})^{\otimes l}[Tg^{\otimes k}](g^{-1})^{\otimes k}\\
&=&(g^{-1})^{\otimes l}T
\end{eqnarray*}

Thus we have $g^{-1}\in\widetilde{G}$, and we conclude that $\widetilde{G}$ is a group, as claimed. 

\medskip

(2) Next, observe that this group $\widetilde{G}$ appears as a closed subgroup $\widetilde{G}\subset U_N$, and also that we have an inclusion $G\subset\widetilde{G}$, coming from definitions. Thus, what we have is an intermediate compact group, as follows, that we want to prove to be equal to $G$:
$$G\subset\widetilde{G}\subset U_N$$

(3) In order to prove this, consider the Tannakian category of $\widetilde{G}$, namely:
$$\widetilde{C}_{kl}=\left\{T\in\mathcal L(H^{\otimes k},H^{\otimes l})\Big|Tg^{\otimes k}=g^{\otimes l}T,\forall g\in\widetilde{G}\right\}$$

By functoriality, from $G\subset\widetilde{G}$ we obtain $\widetilde{C}\subset C$. On the other hand, according to the definition of $\widetilde{G}$, we have $C\subset\widetilde{C}$. Thus, we have the following equality:
$$C=\widetilde{C}$$

(4) Assume now by contradiction that $G\subset\widetilde{G}$ is not an equality. Then, at the level of algebras of functions, the following quotient map is not an isomorphism either:
$$C(\widetilde{G})\to C(G)$$ 

On the other hand, we know from Peter-Weyl that we have decompositions as follows, with the sums being over all irreducible unitary representations:
$$C(\widetilde{G})=\overline{\bigoplus}_{v\in Irr(\widetilde{G})}M_{\dim v}(\mathbb C)
\quad,\quad 
C(G)=\overline{\bigoplus}_{w\in Irr(G)}M_{\dim w}(\mathbb C)$$

Now observe that each unitary representation $v:\widetilde{G}\to U_K$ restricts into a certain representation $v':G\to U_K$. Since the quotient map $C(\widetilde{G})\to C(G)$ is not an isomorphism, we conclude that there is at least one representation $v$ satisfying:
$$v\in Irr(\widetilde{G})\quad,\quad v'\notin Irr(G)$$

(5) We are now in position to conclude. By using Peter-Weyl theory again, the above representation $v\in Irr(\widetilde{G})$ appears in a certain tensor power of the fundamental representation $u:\widetilde{G}\subset U_N$. Thus, we have inclusions of representations, as follows:
$$v\in u^{\otimes k}\quad,\quad v'\in u'^{\otimes k}$$

Now since we know that $v$ is irreducible, and that $v'$ is not, by using one more time Peter-Weyl theory, we conclude that we have a strict inequality, as follows:
\begin{eqnarray*}
\dim(\widetilde{C}(k,k))
&=&dim(End(u^{\otimes k}))\\
&<&dim(End(u'^{\otimes k})\\
&=&\dim(C(k,k))
\end{eqnarray*}

But this contradicts the equality $C=\widetilde{C}$ found in (3), which finishes the proof.
\end{proof}

Our purpose now will be that of showing that we have a correspondence as follows, between closed subgroups $G\subset U_N$, and Tannakian categories $C=(C(k,l))$:
$$G\leftrightarrow C$$

This correspondence, known as Tannakian duality, is something quite deep, and very useful. Indeed, the idea is that what we have here is a useful ``linearization'' of $G$, allowing us to do combinatorics, and ultimately reach to very concrete and powerful results, regarding $G$ itself. And as a consequence, solve our probability questions left.

\bigskip

Speaking linearization of the closed subgroups $G\subset U_N$, we should mention that another way of doing this is by considering the tangent space at the origin $\mathfrak g=T_1(G)$, called Lie algebra of $G$. In what follows, we will use instead our Tannakian approach.

\bigskip

Getting started now, we want to construct a correspondence $G\leftrightarrow C$, and we already know from Theorem 14.4 how the correspondence $G\to C$ appears, namely via:
$$C(k,l)=Hom(u^{\otimes k},u^{\otimes l})$$

Regarding now the construction in the other sense, $C\to G$, this is something very simple as well, coming from the following elementary result: 

\begin{theorem}
Given a tensor category $C=(C(k,l))$ over the space $H\simeq\mathbb C^N$,
$$G=\left\{g\in U_N\Big|Tg^{\otimes k}=g^{\otimes l}T\ ,\ \forall k,l,\forall T\in C(k,l)\right\}$$
is a closed subgroup $G\subset U_N$.
\end{theorem}

\begin{proof}
Consider indeed the closed subset $G\subset U_N$ constructed in the statement. We want to prove that $G$ is indeed a group, and the verifications here go as follows:

\medskip

(1) Given two matrices $g,h\in G$, their product satisfies $gh\in G$, due to the following computation, valid for any $k,l$ and any $T\in C(k,l)$:
\begin{eqnarray*}
T(gh)^{\otimes k}
&=&Tg^{\otimes k}h^{\otimes k}\\
&=&g^{\otimes l}Th^{\otimes k}\\
&=&g^{\otimes l}h^{\otimes l}T\\
&=&(gh)^{\otimes l}T
\end{eqnarray*}

(2) Also, we have $1\in G$, trivially. Finally, for $g\in G$ and $T\in C(k,l)$, we have:
\begin{eqnarray*}
T(g^{-1})^{\otimes k}
&=&(g^{-1})^{\otimes l}[g^{\otimes l}T](g^{-1})^{\otimes k}\\
&=&(g^{-1})^{\otimes l}[Tg^{\otimes k}](g^{-1})^{\otimes k}\\
&=&(g^{-1})^{\otimes l}T
\end{eqnarray*}

Thus we have $g^{-1}\in G$, and so $G$ is a group, as claimed.
\end{proof}

Summarizing, we have so far precise axioms for the tensor categories $C=(C(k,l))$, given in Definition 14.2, as well as correspondences as follows:
$$G\to C\quad,\quad 
C\to G$$

We will show in what follows that these correspondences are inverse to each other. In order to get started, we first have the following technical result:

\begin{theorem}
If we denote the correspondences in Theorem 14.4 and 14.5, between closed subgroups $G\subset U_N$ and tensor categories $C=(C(k,l))$ over $H=\mathbb C^N$, as
$$G\to C_G\quad,\quad 
C\to G_C$$
then we have embeddings as follows, for any $G$ and $C$ respectively,
$$G\subset G_{C_G}\quad,\quad 
C\subset C_{G_C}$$
and proving that these correspondences are inverse to each other amounts in proving
$$C_{G_C}\subset C$$
for any tensor category $C=(C(k,l))$ over the space $H=\mathbb C^N$.
\end{theorem}

\begin{proof}
This is something trivial, with the embeddings $G\subset G_{C_G}$ and $C\subset C_{G_C}$ being both clear from definitions, and with the last assertion coming from this.
\end{proof}

In order to establish Tannakian duality, and more specifically in order to prove the embedding $C_{G_C}\subset C$ appearing above, we will need some abstract constructions. 

\bigskip

Following Malacarne \cite{mal}, let us start with the following elementary fact:

\begin{proposition}
Given a tensor category $C=C((k,l))$ over a Hilbert space $H$,
$$E_C
=\bigoplus_{k,l}C(k,l)
\subset\bigoplus_{k,l}B(H^{\otimes k},H^{\otimes l})
\subset B\left(\bigoplus_kH^{\otimes k}\right)$$
is a closed $*$-subalgebra. Also, inside this algebra,
$$E_C^{(s)}
=\bigoplus_{|k|,|l|\leq s}C(k,l)
\subset\bigoplus_{|k|,|l|\leq s}B(H^{\otimes k},H^{\otimes l})
=B\left(\bigoplus_{|k|\leq s}H^{\otimes k}\right)$$
is a finite dimensional $*$-subalgebra.
\end{proposition}

\begin{proof}
This is clear indeed from the categorical axioms from Definition 14.2, which, since satisfied, prove that the various linear spaces in the statement are stable under both the multiplication operation, and under taking the adjoints.
\end{proof}

Now back to our reconstruction question, we want to prove $C=C_{G_C}$, which is the same as proving $E_C=E_{C_{G_C}}$. We will use a standard commutant trick, as follows:

\index{bicommutant theorem}

\begin{theorem}
For any $*$-algebra $A\subset M_N(\mathbb C)$ we have the equality
$$A=A''$$
where prime denotes the commutant, $X'=\left\{T\in M_N(\mathbb C)\big|Tx=xT,\forall x\in X\right\}$.
\end{theorem}

\begin{proof}
This is a particular case of von Neumann's bicommutant theorem, which follows from the explicit description of $A$ worked out in chapter 13, namely:
$$A=M_{n_1}(\mathbb C)\oplus\ldots\oplus M_{n_k}(\mathbb C)$$

Indeed, the center of each matrix algebra being reduced to the scalars, the commutant of this algebra is as follows, with each copy of $\mathbb C$ corresponding to a matrix block:
$$A'=\mathbb C\oplus\ldots\oplus\mathbb C$$

Now when taking once again the commutant, the computation is trivial, and we obtain in this way $A$ itself, and this leads to the conclusion in the statement.
\end{proof}

By using now the bicommutant theorem, we have:

\begin{theorem}
Given a Tannakian category $C$, the following are equivalent:
\begin{enumerate}
\item $C=C_{G_C}$.

\item $E_C=E_{C_{G_C}}$.

\item $E_C^{(s)}=E_{C_{G_C}}^{(s)}$, for any $s\in\mathbb N$.

\item $E_C^{(s)'}=E_{C_{G_C}}^{(s)'}$, for any $s\in\mathbb N$.
\end{enumerate}
In addition, the inclusions $\subset$, $\subset$, $\subset$, $\supset$ are automatically satisfied.
\end{theorem}

\begin{proof}
This follows from the above results, as follows:

\medskip

$(1)\iff(2)$ This is clear from definitions.

\medskip

$(2)\iff(3)$ This is clear from definitions as well.

\medskip

$(3)\iff(4)$ This comes from the bicommutant theorem. As for the last assertion, we have indeed $C\subset C_{G_C}$ from Theorem 14.6, and this shows that we have as well: 
$$E_C\subset E_{C_{G_C}}$$

We therefore obtain by truncating $E_C^{(s)}\subset E_{C_{G_C}}^{(s)}$, and by taking the commutants, this gives $E_C^{(s)}\supset E_{C_{G_C}}^{(s)}$. Thus, we are led to the conclusion in the statement.
\end{proof}

\section*{14b. The correspondence}

Getting to work now, we would like to prove that we have $E_C^{(s)'}\subset E_{C_{G_C}}^{(s)'}$. Let us first study the commutant on the right.  As a first observation, we have:

\begin{proposition}
We have the following equality,
$$E_{C_G}^{(s)}=End\left(\bigoplus_{|k|\leq s}u^{\otimes k}\right)$$
between subalgebras of $B\left(\bigoplus_{|k|\leq s}H^{\otimes k}\right)$.
\end{proposition}

\begin{proof}
We know that the category $C_G$ is by definition given by:
$$C_G(k,l)=Hom(u^{\otimes k},u^{\otimes l})$$

Thus, the corresponding algebra $E_{C_G}^{(s)}$ appears as follows:
$$E_{C_G}^{(s)}
=\bigoplus_{|k|,|l|\leq s}Hom(u^{\otimes k},u^{\otimes l})
\subset\bigoplus_{|k|,|l|\leq s}B(H^{\otimes k},H^{\otimes l})
=B\left(\bigoplus_{|k|\leq s}H^{\otimes k}\right)$$

On the other hand, the algebra of intertwiners of $\bigoplus_{|k|\leq s}u^{\otimes k}$ is given by:
$$End\left(\bigoplus_{|k|\leq s}u^{\otimes k}\right)
=\bigoplus_{|k|,|l|\leq s}Hom(u^{\otimes k},u^{\otimes l})
\subset\bigoplus_{|k|,|l|\leq s}B(H^{\otimes k},H^{\otimes l})
=B\left(\bigoplus_{|k|\leq s}H^{\otimes k}\right)$$

Thus we have indeed the same algebra, and we are done.
\end{proof}

We have to compute the commutant of the above algebra. For this purpose, we can use the following general result, valid for any representation of a compact group:

\begin{proposition}
Given a unitary group representation $v:G\to U_n$ we have an algebra representation as follows,
$$\pi_v:C(G)^*\to M_n(\mathbb C)\quad,\quad 
\varphi\to (\varphi(v_{ij}))_{ij}$$
whose image is given by $Im(\pi_v)=End(v)'$.
\end{proposition}

\begin{proof}
The first assertion is clear, with the multiplicativity claim for $\pi_v$ coming from the following computation, where $\Delta:C(G)\to C(G)\otimes C(G)$ is the comultiplication:
\begin{eqnarray*}
(\pi_v(\varphi*\psi))_{ij}
&=&(\varphi\otimes\psi)\Delta(v_{ij})\\
&=&\sum_k\varphi(v_{ik})\psi(v_{kj})\\
&=&\sum_k(\pi_v(\varphi))_{ik}(\pi_v(\psi))_{kj}\\
&=&(\pi_v(\varphi)\pi_v(\psi))_{ij}
\end{eqnarray*} 

Let us establish now the equality in the statement, namely:
$$Im(\pi_v)=End(v)'$$

Let us first prove the inclusion $\subset$. Given $\varphi\in C(G)^*$ and $T\in End(v)$, we have:
\begin{eqnarray*}
[\pi_v(\varphi),T]=0
&\iff&\sum_k\varphi(v_{ik})T_{kj}=\sum_kT_{ik}\varphi(v_{kj}),\forall i,j\\
&\iff&\varphi\left(\sum_kv_{ik}T_{kj}\right)=\varphi\left(\sum_kT_{ik}v_{kj}\right),\forall i,j\\
&\iff&\varphi((vT)_{ij})=\varphi((Tv)_{ij}),\forall i,j
\end{eqnarray*}

But this latter formula is true, because $T\in End(v)$ means that we have:
$$vT=Tv$$

As for the converse inclusion $\supset$, the proof is quite similar. Indeed, by using the bicommutant theorem, this is the same as proving that we have:
$$Im(\pi_v)'\subset End(v)$$

But, by using the above equivalences, we have the following computation:
\begin{eqnarray*}
T\in Im(\pi_v)'
&\iff&[\pi_v(\varphi),T]=0,\forall\varphi\\
&\iff&\varphi((vT)_{ij})=\varphi((Tv)_{ij}),\forall\varphi,i,j\\
&\iff&vT=Tv
\end{eqnarray*}

Thus, we have obtained the desired inclusion, and we are done.
\end{proof}

By combining the above results, we obtain the following technical statement:

\begin{theorem}
We have $E_{C_G}^{(s)'}=Im(\pi_v)$, where $v$ is the following direct sum,
$$v=\bigoplus_{|k|\leq s}u^{\otimes k}$$
and where the algebra representation $\pi_v:C(G)^*\to M_n(\mathbb C)$ is given by $\varphi\to(\varphi(v_{ij}))_{ij}$.
\end{theorem}

\begin{proof}
This follows indeed by combining the above results, and more precisely by combining Proposition 14.10 and Proposition 14.11.
\end{proof}

We recall that we want to prove that we have $E_C^{(s)'}\subset E_{C_{G_C}}^{(s)'}$, for any $s\in\mathbb N$. And for this purpose, we must first refine Theorem 14.12, in the case $G=G_C$. 

\bigskip

Generally speaking, in order to prove anything about $G_C$, we are in need of an explicit model for this group. In order to construct such a model, let $<u_{ij}>$ be the free $*$-algebra over $\dim(H)^2$ variables, with comultiplication and counit as follows:
$$\Delta(u_{ij})=\sum_ku_{ik}\otimes u_{kj}\quad,\quad 
\varepsilon(u_{ij})=\delta_{ij}$$

Following \cite{mal}, we can model this $*$-bialgebra, in the following way:

\begin{proposition}
Consider the following pair of dual vector spaces,
$$F=\bigoplus_kB\left(H^{\otimes k}\right)\quad,\quad 
F^*=\bigoplus_kB\left(H^{\otimes k}\right)^*$$
and let $f_{ij},f_{ij}^*\in F^*$ be the standard generators of $B(H)^*,B(\bar{H})^*$.
\begin{enumerate}
\item $F^*$ is a $*$-algebra, with multiplication $\otimes$ and involution as follows:
$$f_{ij}\leftrightarrow f_{ij}^*$$

\item $F^*$ is a $*$-bialgebra, with $*$-bialgebra operations as follows:
$$\Delta(f_{ij})=\sum_kf_{ik}\otimes f_{kj}\quad,\quad 
\varepsilon(f_{ij})=\delta_{ij}$$

\item We have a $*$-bialgebra isomorphism $<u_{ij}>\simeq F^*$, given by $u_{ij}\to f_{ij}$.
\end{enumerate}
\end{proposition}

\begin{proof}
Since $F^*$ is spanned by the various tensor products between the variables $f_{ij},f_{ij}^*$, we have a vector space isomorphism as follows:
$$<u_{ij}>\simeq F^*\quad,\quad 
u_{ij}\to f_{ij}\quad,\quad
u_{ij}^*\to f_{ij}^*$$

The corresponding $*$-bialgebra structure induced on the vector space $F^*$ is then the one in the statement, and this gives the result.
\end{proof}

Now back to our group $G_C$, we have the following modeling result for it:

\begin{proposition}
The smooth part of the algebra $A_C=C(G_C)$ is given by 
$$\mathcal A_C\simeq F^*/J$$
where $J\subset F^*$ is the ideal coming from the following relations, for any $i,j$,
$$\sum_{p_1,\ldots,p_k}T_{i_1\ldots i_l,p_1\ldots p_k}f_{p_1j_1}\otimes\ldots\otimes  f_{p_kj_k}
=\sum_{q_1,\ldots,q_l}T_{q_1\ldots q_l,j_1\ldots j_k}f_{i_1q_1}\otimes\ldots\otimes f_{i_lq_l}$$
one for each pair of colored integers $k,l$, and each $T\in C(k,l)$.
\end{proposition}

\begin{proof}
As a first observation, $A_C$ appears as enveloping $C^*$-algebra of the following universal $*$-algebra, where $u=(u_{ij})$ is regarded as a formal corepresentation:
$$\mathcal A_C=\left<(u_{ij})_{i,j=1,\ldots,N}\Big|T\in Hom(u^{\otimes k},u^{\otimes l}),\forall k,l,\forall T\in C(k,l)\right>$$

With this observation in hand, the conclusion is that we have a formula as follows, where $I$ is the ideal coming from the relations $T\in Hom(u^{\otimes k},u^{\otimes l})$, with $T\in C(k,l)$:
$$\mathcal A_C=<u_{ij}>/I$$

Now if we denote by $J\subset F^*$ the image of the ideal $I$ via the $*$-algebra isomorphism $<u_{ij}>\simeq F^*$ from Proposition 14.16, we obtain an identification as follows:
$$\mathcal A_C\simeq F^*/J$$

With standard multi-index notations, and by assuming now that $k,l\in\mathbb N$ are usual integers, for simplifying the presentation, the general case being similar, a relation of type $T\in Hom(u^{\otimes k},u^{\otimes l})$ inside $<u_{ij}>$ is equivalent to the following conditions:
$$\sum_{p_1,\ldots,p_k}T_{i_1\ldots i_l,p_1\ldots p_k}u_{p_1j_1}\ldots u_{p_kj_k}
=\sum_{q_1,\ldots,q_l}T_{q_1\ldots q_l,j_1\ldots j_k}u_{i_1q_1}\ldots u_{i_lq_l}$$

Now by recalling that the isomorphism of $*$-algebras $<u_{ij}>\to F^*$ is given by $u_{ij}\to f_{ij}$, and that the multiplication operation of $F^*$ corresponds to the tensor product operation $\otimes$, we conclude that $J\subset F^*$ is the ideal from the statement.
\end{proof}

With the above result in hand, let us go back to Theorem 14.12. We have:

\begin{proposition}
The linear space $\mathcal A_C^*$ is given by the formula
$$\mathcal A_C^*=\left\{a\in F\Big|Ta_k=a_lT,\forall T\in C(k,l)\right\}$$
and the representation
$$\pi_v:\mathcal A_C^*\to B\left(\bigoplus_{|k|\leq s}H^{\otimes k}\right)$$
appears diagonally, by truncating, $\pi_v:a\to (a_k)_{kk}$.
\end{proposition}

\begin{proof}
We know from Proposition 14.14 that we have an identification of $*$-bialgebras $\mathcal A_C\simeq F^*/J$. But this gives a quotient map, as follows:
$$F^*\to\mathcal A_C$$

At the dual level, this gives $\mathcal A_C^*\subset F$. To be more precise, we have:
$$\mathcal A_C^*=\left\{a\in\ F\Big|f(a)=0,\forall f\in J\right\}$$

Now since $J=<f_T>$, where $f_T$ are the relations in Proposition 14.14, we obtain:
$$\mathcal A_C^*=\left\{a\in F\Big|f_T(a)=0,\forall T\in C\right\}$$

Given $T\in C(k,l)$, for an arbitrary element $a=(a_k)$, we have:
\begin{eqnarray*}
&&f_T(a)=0\\
&\iff&\sum_{p_1,\ldots,p_k}T_{i_1\ldots i_l,p_1\ldots p_k}(a_k)_{p_1\ldots p_k,j_1\ldots j_k}=\sum_{q_1,\ldots,q_l}T_{q_1\ldots q_l,j_1\ldots j_k}(a_l)_{i_1\ldots i_l,q_1\ldots q_l},\forall i,j\\
&\iff&(Ta_k)_{i_1\ldots i_l,j_1\ldots j_k}=(a_lT)_{i_1\ldots i_l,j_1\ldots j_k},\forall i,j\\
&\iff&Ta_k=a_lT
\end{eqnarray*}

Thus, $\mathcal A_C^*$ is given by the formula in the statement. It remains to compute $\pi_v$:
$$\pi_v:\mathcal A_C^*\to B\left(\bigoplus_{|k|\leq s}H^{\otimes k}\right)$$

With $a=(a_k)$, we have the following computation:
\begin{eqnarray*}
\pi_v(a)_{i_1\ldots i_k,j_1\ldots j_k}
&=&a(v_{i_1\ldots i_k,j_1\ldots j_k})\\
&=&(f_{i_1j_1}\otimes\ldots\otimes f_{i_kj_k})(a)\\
&=&(a_k)_{i_1\ldots i_k,j_1\ldots j_k}
\end{eqnarray*}

Thus, our representation $\pi_v$ appears diagonally, by truncating, as claimed.
\end{proof}

In order to further advance, consider the following vector spaces:
$$F_s=\bigoplus_{|k|\leq s}B\left(H^{\otimes k}\right)\quad,\quad 
F^*_s=\bigoplus_{|k|\leq s}B\left(H^{\otimes k}\right)^*$$

We denote by $a\to a_s$ the truncation operation $F\to F_s$. We have:

\begin{proposition}
The following hold:
\begin{enumerate}
\item $E_C^{(s)'}\subset F_s$.

\item $E_C'\subset F$.

\item $\mathcal A_C^*=E_C'$.

\item $Im(\pi_v)=(E_C')_s$.
\end{enumerate}
\end{proposition}

\begin{proof}
These results basically follow from what we have, as follows:

\medskip

(1) We have an inclusion as follows, as a diagonal subalgebra:
$$F_s\subset B\left(\bigoplus_{|k|\leq s}H^{\otimes k}\right)$$

The commutant of this algebra is then given by:
$$F_s'=\left\{b\in F_s\Big|b=(b_k),b_k\in\mathbb C,\forall k\right\}$$

On the other hand, we know from the identity axiom for the category $C$ that we have $F_s'\subset E_C^{(s)}$. Thus, our result follows from the bicommutant theorem, as follows:
$$F_s'\subset E_C^{(s)}\implies F_s\supset E_C^{(s)'}$$

(2) This follows from (1), by taking inductive limits.

\medskip

(3) With the present notations, the formula of $\mathcal A_C^*$ from Proposition 14.15 reads $\mathcal A_C^*=F\cap E_C'$. Now since by (2) we have $E_C'\subset F$, we obtain from this $\mathcal A_C^*=E_C'$.

\medskip

(4) This follows from (3), and from the formula of $\pi_\nu$ in Proposition 14.15.
\end{proof}

Following \cite{mal}, we can now state and prove our main result, as follows:

\index{Tannakian duality}

\begin{theorem}
The Tannakian duality constructions 
$$C\to G_C\quad,\quad 
G\to C_G$$
are inverse to each other.
\end{theorem}

\begin{proof}
According to our various results above, we have to prove that, for any Tannakian category $C$, and any $s\in\mathbb N$, we have an inclusion as follows:
$$E_C^{(s)'}\subset(E_C')_s$$

By taking duals, this is the same as proving that we have:
$$\left\{f\in F_s^*\Big|f_{|(E_C')_s}=0\right\}\subset\left\{f\in F_s^*\Big|f_{|E_C^{(s)'}}=0\right\}$$

In order to do so, we use the following formula, from Proposition 14.16: 
$$\mathcal A_C^*=E_C'$$

We know from the above that we have an identification as follows:
$$\mathcal A_C=F^*/J$$

We conclude that the ideal $J$ is given by the following formula:
$$J=\left\{f\in F^*\Big|f_{|E_C'}=0\right\}$$

Our claim is that we have the following formula, for any $s\in\mathbb N$:
$$J\cap F_s^*=\left\{f\in F_s^*\Big|f_{|E_C^{(s)'}}=0\right\}$$

Indeed, let us denote by $X_s$ the spaces on the right. The axioms for $C$ show that these spaces are increasing, that their union $X=\cup_sX_s$ is an ideal, and that:
$$X_s=X\cap F_s^*$$

We must prove that we have $J=X$, and this can be done as follows:

\medskip

``$\subset$'' This follows from the following fact, for any $T\in C(k,l)$ with $|k|,|l|\leq s$:
\begin{eqnarray*}
(f_T)_{|\{T\}'}=0
&\implies&(f_T)_{|E_C^{(s)'}}=0\\
&\implies&f_T\in X_s
\end{eqnarray*}

``$\supset$'' This follows from our description of $J$, because from $E_C^{(s)}\subset E_C$ we obtain:
$$f_{|E_C^{(s)'}}=0\implies f_{|E_C'}=0$$

Summarizing, we have proved our claim. On the other hand, we have:
\begin{eqnarray*}
J\cap F_s^*
&=&\left\{f\in F^*\Big|f_{|E_C'}=0\right\}\cap F_s^*\\
&=&\left\{f\in F_s^*\Big|f_{|E_C'}=0\right\}\\
&=&\left\{f\in F_s^*\Big|f_{|(E_C')_s}=0\right\}
\end{eqnarray*}

Thus, our claim is exactly the inclusion that we wanted to prove, and we are done.
\end{proof}

Summarizing, we have proved Tannakian duality. We should mention that there are many other versions of this duality, and for more on this, we refer to the quantum algebra literature, where Tannakian duality, in all its forms, is something highly valued.

\section*{14c. Brauer theorems}

As a basic illustration for the Tannakian correspondence, we will work out now Brauer theorems for $O_N,U_N$. These are very classical results, and there are many possible proofs for them. We will follow here the modern approach from \cite{bsp}. Let us start with:

\index{Kronecker symbols}
\index{maps associated to partitions}

\begin{definition}
Given a pairing $\pi\in P_2(k,l)$ and an integer $N\in\mathbb N$, we can construct a linear map between tensor powers of $\mathbb C^N$,
$$T_\pi:(\mathbb C^N)^{\otimes k}\to(\mathbb C^N)^{\otimes l}$$
by the following formula, with $e_1,\ldots,e_N$ being the standard basis of $\mathbb C^N$,
$$T_\pi(e_{i_1}\otimes\ldots\otimes e_{i_k})=\sum_{j_1\ldots j_l}\delta_\pi\begin{pmatrix}i_1&\ldots&i_k\\ j_1&\ldots&j_l\end{pmatrix}e_{j_1}\otimes\ldots\otimes e_{j_l}$$
and with the coefficients on the right being Kronecker type symbols,
$$\delta_\pi\begin{pmatrix}i_1&\ldots&i_k\\ j_1&\ldots&j_l\end{pmatrix}\in\{0,1\}$$
whose values depend on whether the indices fit or not.
\end{definition}

To be more precise here, we put the multi-indices $i=(i_1,\ldots,i_k)$ and $j=(j_1,\ldots,j_l)$ on the legs of our pairing $\pi$, in the obvious way. In the case where all strings of $\pi$ join pairs of equal indices of $i,j$, we set $\delta_\pi(^i_j)=1$. Otherwise, we set $\delta_\pi(^i_j)=0$. 

\bigskip

The point with the above definition comes from the fact that most of the ``familiar'' maps, in the Tannakian context, are of the above form. Here are some examples:

\begin{proposition}
The correspondence $\pi\to T_\pi$ has the following properties:
\begin{enumerate}
\item $T_\cap=(1\to\sum_ie_i\otimes e_i)$.

\item $T_\cup=(e_i\otimes e_j\to\delta_{ij})$.

\item $T_{||\ldots||}=id$.

\item $T_{\slash\hskip-1.5mm\backslash}=(e_a\otimes e_b\to e_b\otimes e_a)$.
\end{enumerate}
\end{proposition}

\begin{proof}
We can assume that all legs of $\pi$ are colored $\circ$, and then:

\medskip

(1) We have $\cap\in P_2(\emptyset,\circ\circ)$, so the corresponding linear map is as follows:
$$T_\cap:\mathbb C\to\mathbb C^N\otimes\mathbb C^N$$

The formula of this linear map is then, as claimed:
\begin{eqnarray*}
T_\cap(1)
&=&\sum_{ij}\delta_\cap(i\ j)e_i\otimes e_j\\
&=&\sum_{ij}\delta_{ij}e_i\otimes e_j\\
&=&\sum_ie_i\otimes e_i
\end{eqnarray*}

(2) Here we have $\cup\in P_2(\circ\circ,\emptyset)$, so the corresponding linear map is as follows:
$$T_\cap:\mathbb C^N\otimes\mathbb C^N\to\mathbb C$$

The formula of this linear form is then as follows:
$$T_\cap(e_i\otimes e_j)
=\delta_\cap(i\ j)
=\delta_{ij}$$

(3) Consider indeed the ``identity'' pairing $||\ldots||\in P_2(k,k)$, with $k=\circ\circ\ldots\circ\circ$. The corresponding linear map is then the identity, because we have:
\begin{eqnarray*}
T_{||\ldots||}(e_{i_1}\otimes\ldots\otimes e_{i_k})
&=&\sum_{j_1\ldots j_k}\delta_{||\ldots||}\begin{pmatrix}i_1&\ldots&i_k\\ j_1&\ldots&j_k\end{pmatrix}e_{j_1}\otimes\ldots\otimes e_{j_k}\\
&=&\sum_{j_1\ldots j_k}\delta_{i_1j_1}\ldots\delta_{i_kj_k}e_{j_1}\otimes\ldots\otimes e_{j_k}\\
&=&e_{i_1}\otimes\ldots\otimes e_{i_k}
\end{eqnarray*}

(4) For the basic crossing $\slash\hskip-2.0mm\backslash\in P_2(\circ\circ,\circ\circ)$, the corresponding linear map is as follows:
$$T_{\slash\hskip-1.5mm\backslash}:\mathbb C^N\otimes\mathbb C^N\to\mathbb C^N\otimes\mathbb C^N$$

This linear map can be computed as follows:
\begin{eqnarray*}
T_{\slash\hskip-1.5mm\backslash}(e_i\otimes e_j)
&=&\sum_{kl}\delta_{\slash\hskip-1.5mm\backslash}\begin{pmatrix}i&j\\ k&l\end{pmatrix}e_k\otimes e_l\\
&=&\sum_{kl}\delta_{il}\delta_{jk}e_k\otimes e_l\\
&=&e_j\otimes e_i
\end{eqnarray*}

Thus we obtain the flip operator $\Sigma(a\otimes b)=b\otimes a$, as claimed.
\end{proof}

The relation with the Tannakian categories comes from the following key result:

\begin{proposition}
The assignement $\pi\to T_\pi$ is categorical, in the sense that
$$T_\pi\otimes T_\sigma=T_{[\pi\sigma]}\quad,\quad 
T_\pi T_\sigma=N^{c(\pi,\sigma)}T_{[^\sigma_\pi]}\quad,\quad 
T_\pi^*=T_{\pi^*}$$
where $c(\pi,\sigma)$ is the number of circles appearing in the middle, when concatenating.
\end{proposition}

\begin{proof}
The concatenation axiom follows from the following computation:
\begin{eqnarray*}
&&(T_\pi\otimes T_\sigma)(e_{i_1}\otimes\ldots\otimes e_{i_p}\otimes e_{k_1}\otimes\ldots\otimes e_{k_r})\\
&=&\sum_{j_1\ldots j_q}\sum_{l_1\ldots l_s}\delta_\pi\begin{pmatrix}i_1&\ldots&i_p\\j_1&\ldots&j_q\end{pmatrix}\delta_\sigma\begin{pmatrix}k_1&\ldots&k_r\\l_1&\ldots&l_s\end{pmatrix}e_{j_1}\otimes\ldots\otimes e_{j_q}\otimes e_{l_1}\otimes\ldots\otimes e_{l_s}\\
&=&\sum_{j_1\ldots j_q}\sum_{l_1\ldots l_s}\delta_{[\pi\sigma]}\begin{pmatrix}i_1&\ldots&i_p&k_1&\ldots&k_r\\j_1&\ldots&j_q&l_1&\ldots&l_s\end{pmatrix}e_{j_1}\otimes\ldots\otimes e_{j_q}\otimes e_{l_1}\otimes\ldots\otimes e_{l_s}\\
&=&T_{[\pi\sigma]}(e_{i_1}\otimes\ldots\otimes e_{i_p}\otimes e_{k_1}\otimes\ldots\otimes e_{k_r})
\end{eqnarray*}

The composition axiom follows from the following computation:
\begin{eqnarray*}
&&T_\pi T_\sigma(e_{i_1}\otimes\ldots\otimes e_{i_p})\\
&=&\sum_{j_1\ldots j_q}\delta_\sigma\begin{pmatrix}i_1&\ldots&i_p\\j_1&\ldots&j_q\end{pmatrix}
\sum_{k_1\ldots k_r}\delta_\pi\begin{pmatrix}j_1&\ldots&j_q\\k_1&\ldots&k_r\end{pmatrix}e_{k_1}\otimes\ldots\otimes e_{k_r}\\
&=&\sum_{k_1\ldots k_r}N^{c(\pi,\sigma)}\delta_{[^\sigma_\pi]}\begin{pmatrix}i_1&\ldots&i_p\\k_1&\ldots&k_r\end{pmatrix}e_{k_1}\otimes\ldots\otimes e_{k_r}\\
&=&N^{c(\pi,\sigma)}T_{[^\sigma_\pi]}(e_{i_1}\otimes\ldots\otimes e_{i_p})
\end{eqnarray*}

Finally, the involution axiom follows from the following computation:
\begin{eqnarray*}
&&T_\pi^*(e_{j_1}\otimes\ldots\otimes e_{j_q})\\
&=&\sum_{i_1\ldots i_p}<T_\pi^*(e_{j_1}\otimes\ldots\otimes e_{j_q}),e_{i_1}\otimes\ldots\otimes e_{i_p}>e_{i_1}\otimes\ldots\otimes e_{i_p}\\
&=&\sum_{i_1\ldots i_p}\delta_\pi\begin{pmatrix}i_1&\ldots&i_p\\ j_1&\ldots& j_q\end{pmatrix}e_{i_1}\otimes\ldots\otimes e_{i_p}\\
&=&T_{\pi^*}(e_{j_1}\otimes\ldots\otimes e_{j_q})
\end{eqnarray*}

Summarizing, our correspondence is indeed categorical.
\end{proof}

The above result suggests the following general definition, from \cite{bsp}:

\index{category of partitions}
\index{matching pairings}

\begin{definition}
Let $P_2(k,l)$ be the set of pairings between an upper colored integer $k$, and a lower colored integer $l$. A collection of subsets 
$$D=\bigsqcup_{k,l}D(k,l)$$
with $D(k,l)\subset P_2(k,l)$ is called a category of pairings when it has the following properties:
\begin{enumerate}
\item Stability under the horizontal concatenation, $(\pi,\sigma)\to[\pi\sigma]$.

\item Stability under vertical concatenation $(\pi,\sigma)\to[^\sigma_\pi]$, with matching middle symbols.

\item Stability under the upside-down turning $*$, with switching of colors, $\circ\leftrightarrow\bullet$.

\item Each set $P(k,k)$ contains the identity partition $||\ldots||$.

\item The sets $P(\emptyset,\circ\bullet)$ and $P(\emptyset,\bullet\circ)$ both contain the semicircle $\cap$.

\item The sets $P(k,\bar{k})$ with $|k|=2$ contain the crossing partition $\slash\hskip-2.0mm\backslash$.
\end{enumerate}
\end{definition} 

Observe the similarity with the axioms for Tannakian categories, from the beginning of this chapter. We will see in a moment that this similarity can be turned into something very precise, with the categories of pairings producing Tannakian categories.

\bigskip

As basic examples of such categories, that we have already met in the above, we have the categories $P_2,\mathcal P_2$ of pairings, and of matching pairings, with the convention that a matching pairing must pair $\circ-\bullet$ on the horizontal, and $\circ-\circ$ or $\bullet-\bullet$ on the vertical. There are many other examples, and we will discuss this gradually, in what follows.

\bigskip

In relation with the compact groups, we have the following result:

\begin{theorem}
Each category of pairings, in the above sense,
$$D=(D(k,l))$$
produces a family of compact groups $G=(G_N)$, one for each $N\in\mathbb N$, via the formula
$$Hom(u^{\otimes k},u^{\otimes l})=span\left(T_\pi\Big|\pi\in D(k,l)\right)$$
and the Tannakian duality correspondence.
\end{theorem}

\begin{proof}
Given an integer $N\in\mathbb N$, consider the correspondence $\pi\to T_\pi$ constructed in Definition 14.18, and then the collection of linear spaces in the statement, namely:
$$C_{kl}=span\left(T_\pi\Big|\pi\in D(k,l)\right)$$

According to Proposition 14.20, and to our axioms for the categories of partitions, from Definition 14.21, this collection of spaces $C=(C_{kl})$ satisfies the axioms for the Tannakian categories, from the beginning of this chapter. Thus the Tannakian duality result applies, and provides us with a closed subgroup $G_N\subset U_N$ such that:
$$C_{kl}=Hom(u^{\otimes k},u^{\otimes l})$$

Thus, we are led to the conclusion in the statement.
\end{proof}

The above result is something fundamental, and suggests formulating:

\begin{definition}
Assuming that a closed subgroup $G\subset_uU_N$ has the property
$$Hom(u^{\otimes k},u^{\otimes l})=span\left(T_\pi\Big|\pi\in D(k,l)\right)$$
for a certain category of pairings $D=(D(k,l))$, we say that $G$ is easy.
\end{definition}

This definition, from \cite{bsp}, is motivated by the fact that, from the point of view of Tannakian duality, the above groups are indeed the ``easiest'' possible ones. Of course, this might sound a bit strange, after all the quite complicated things that we did in this chapter. But hey, there is a beginning for everything. We will get to know better Tannakian duality and easiness, and their applications, in what follows, and please believe me, you will reach too to the conclusion that Definition 14.23 is justified.

\bigskip

As another comment, it is possible to talk about more general easy groups, by using general categories of partitions, instead of just categories of pairings. We will be back to all this, with a systematic study of easiness, in chapter 15 below. 

\bigskip

As a technical remark now, to be always kept in mind, when dealing with easiness, the category of pairings producing an easy group is not unique, for instance because at $N=1$ all the possible categories of pairings produce the same easy group, namely the trivial group $G=\{1\}$. Thus, some subtleties are going on here. More on this later.

\bigskip

Getting back now to concrete things, the point now is that with the above ingredients in hand, and as a first application of Tannakian duality, we can establish a useful result, namely the Brauer theorem for the unitary group $U_N$. The statement is a follows:

\index{Brauer theorem}

\begin{theorem}
For the unitary group $U_N$ we have
$$Hom(u^{\otimes k},u^{\otimes l})=span\left(T_\pi\Big|\pi\in\mathcal P_2(k,l)\right)$$
where $\mathcal P_2$ denotes as usual the category of all matching pairings.
\end{theorem}

\begin{proof}
This is something very old and classical, due to Brauer \cite{bra}, and in what follows we will present a simplified proof for it, based on the easiness technology developed above. Consider the spaces on the right in the statement, namely:
$$C_{kl}=span\left(T_\pi\Big|\pi\in\mathcal P_2(k,l)\right)$$

According to Proposition 14.20 these spaces form a tensor category. Thus, by Tannakian duality, these spaces must come from a certain closed subgroup $G\subset U_N$. To be more precise, if we denote by $v$ the fundamental representation of $G$, then:
$$C_{kl}=Hom(v^{\otimes k},v^{\otimes l})$$

We must prove that we have $G=U_N$. For this purpose, let us recall that the unitary group $U_N$ is defined via the following relations:
$$u^*=u^{-1}\quad,\quad 
u^t=\bar{u}^{-1}$$ 

But these relations tell us precisely that the following two operators must be in the associated Tannakian category $C$:
$$T_\pi\quad:\quad \pi={\ }^{\,\cap}_{\circ\bullet}\ ,\ {\ }^{\,\cap}_{\bullet\circ}$$

Thus the associated Tannakian category is $C=span(T_\pi|\pi\in D)$, with:
$$D
=<{\ }^{\,\cap}_{\circ\bullet}\,\,,{\ }^{\,\cap}_{\bullet\circ}>
=\mathcal P_2$$
  
Thus, we are led to the conclusion in the statement.
\end{proof}

Regarding the orthogonal group $O_N$, we have here a similar result, as follows:

\index{Brauer theorem}

\begin{theorem}
For the orthogonal group $O_N$ we have
$$Hom(u^{\otimes k},u^{\otimes l})=span\left(T_\pi\Big|\pi\in P_2(k,l)\right)$$
where $P_2$ denotes as usual the category of all pairings.
\end{theorem}

\begin{proof}
As before with Theorem 14.24, regarding $U_N$, this is something very old and classical, due to Brauer \cite{bra}, that we can now prove by using the easiness technology developed above. Consider the spaces on the right in the statement, namely:
$$C_{kl}=span\left(T_\pi\Big|\pi\in P_2(k,l)\right)$$

According to Proposition 14.20 these spaces form a tensor category. Thus, by Tannakian duality, these spaces must come from a certain closed subgroup $G\subset U_N$. To be more precise, if we denote by $v$ the fundamental representation of $G$, then:
$$C_{kl}=Hom(v^{\otimes k},v^{\otimes l})$$

We must prove that we have $G=O_N$. For this purpose, let us recall that the orthogonal group $O_N\subset U_N$ is defined by imposing the following relations:
$$u_{ij}=\bar{u}_{ij}$$

But these relations tell us precisely that the following two operators must be in the associated Tannakian category $C$:
$$T_\pi\quad:\quad\pi=|^{\hskip-1.32mm\circ}_{\hskip-1.32mm\bullet}\ ,\ |_{\hskip-1.32mm\circ}^{\hskip-1.32mm\bullet}$$

Thus the associated Tannakian category is $C=span(T_\pi|\pi\in D)$, with:
$$D
=<\mathcal P_2,|^{\hskip-1.32mm\circ}_{\hskip-1.32mm\bullet},|_{\hskip-1.32mm\circ}^{\hskip-1.32mm\bullet}>
=P_2$$
 
Thus, we are led to the conclusion in the statement.
\end{proof}

We will see later, in chapter 16 below, applications of the above results, to integration problems over $O_N,U_N$, by using the Peter-Weyl methods from chapter 13. 

\section*{14d. Clebsch-Gordan rules}

As a last piece of representation theory, we are now in position of dealing, in a quite conceptual way, with $SU_2$ and $SO_3$. Regarding $SU_2$, the result here is as follows:

\index{Wigner law}
\index{semicircle law}
\index{Clebsch-Gordan rules}
\index{fusion rules}

\begin{theorem}
The irreducible representations of $SU_2$ are all self-adjoint, and can be labeled by positive integers, with their fusion rules being as follows,
$$r_k\otimes r_l=r_{|k-l|}+r_{|k-l|+2}+\ldots+r_{k+l}$$
called Clebsch-Gordan rules. The corresponding dimensions are $\dim r_k=k+1$.
\end{theorem}

\begin{proof}
There are several proofs for this fact, the simplest one, with the knowledge that we have, being via purely algebraic methods, as follows:

\medskip

(1) Our first claim is that we have the following estimate, telling us that the even moments of the main character are smaller than the Catalan numbers:
$$\int_{SU_2}\chi^{2k}\leq C_k$$

But this is something that we know from chapter 12, obtained by using $SU_2\simeq S^3_\mathbb R$ and spherical integrals, and with the stronger statement that we have in fact equality $=$. However, for the purposes of what follows, the above $\leq$ estimate will do.

\medskip

(2) Alternatively, the above estimate can be deduced with purely algebraic methods, by using an easiness type argument for $SU_2$, as follows:
\begin{eqnarray*}
\int_{SU_2}\chi^{2k}
&=&\dim(Fix(u^{\otimes 2k}))\\
&=&\dim\left(span\left(T_\pi'\Big|\pi\in NC_2(2k)\right)\right)\\
&\leq&|NC_2(2k)|\\
&=&C_k
\end{eqnarray*}

To be more precise, $SU_2$ is not exactly easy, but rather ``super-easy'', coming from a different implementation $\pi\to T_\pi'$ of the pairings, involving some signs. And with this being proved exactly as the Brauer theorem for $O_N$, with modifications where needed.

\medskip

(3) Long story short, we have our estimate in (1), and this is all that we need. Our claim is that we can construct, by recurrence on $k\in\mathbb N$, a sequence $r_k$ of irreducible, self-adjoint and distinct representations of $SU_2$, satisfying:
$$r_0=1\quad,\quad
r_1=u\quad,\quad 
r_k+r_{k-2}=r_{k-1}\otimes r_1$$

Indeed, assume that $r_0,\ldots,r_{k-1}$ are constructed, and let us construct $r_k$. We have:
$$r_{k-1}+r_{k-3}=r_{k-2}\otimes r_1$$

Thus $r_{k-1}\subset r_{k-2}\otimes r_1$, and since $r_{k-2}$ is irreducible, by Frobenius we have:
$$r_{k-2}\subset r_{k-1}\otimes r_1$$

We conclude there exists a certain representation $r_k$ such that:
$$r_k+r_{k-2}=r_{k-1}\otimes r_1$$

(4) By recurrence, $r_k$ is self-adjoint. Now observe that according to our recurrence formula, we can split $u^{\otimes k}$ as a sum of the following type, with positive coefficients:  
$$u^{\otimes k}=c_kr_k+c_{k-2}r_{k-2}+\ldots$$

We conclude by Peter-Weyl that we have an inequality as follows, with equality precisely when $r_k$ is irreducible, and non-equivalent to the other summands $r_i$:
$$\sum_ic_i^2\leq\dim(End(u^{\otimes k}))$$

(5) But by (1) the number on the right is $\leq C_k$, and some straightforward combinatorics, based on the fusion rules, shows that the number on the left is $C_k$ as well:
$$C_k
=\sum_ic_i^2
\leq\dim(End(u^{\otimes k}))
=\int_{SU_2}\chi^{2k}
\leq C_k$$

Thus we have equality in our estimate, so our representation $r_k$ is irreducible, and non-equivalent to $r_{k-2},r_{k-4},\ldots$ Moreover, this representation $r_k$ is not equivalent to $r_{k-1},r_{k-3},\ldots$ either, with this coming from $r_p\subset u^{\otimes p}$ for any $p$, and from:
$$\dim(Fix(u^{\otimes 2s+1}))=\int_{SU_2}\chi^{2s+1}=0$$

(6) Thus, we proved our claim. Now since each irreducible representation of $SU_2$ appears into some $u^{\otimes k}$, and we know how to decompose each $u^{\otimes k}$ into sums of representations $r_k$, these representations $r_k$ are all the irreducible representations of $SU_2$, and we are done with the main assertion. As for the dimension formula, this is clear.
\end{proof}

Regarding now $SO_3$, we have here a similar result, as follows:

\index{Catalan numbers}
\index{Clebsch-Gordan rules}
\index{fusion rules}
\index{Marchenko-Pastur law}

\begin{theorem}
The irreducible representations of $SO_3$ are all self-adjoint, and can be labeled by positive integers, with their fusion rules being as follows,
$$r_k\otimes r_l=r_{|k-l|}+r_{|k-l|+1}+\ldots+r_{k+l}$$ 
also called Clebsch-Gordan rules. The corresponding dimensions are $\dim r_k=2k+1$.
\end{theorem}

\begin{proof}
As before with $SU_2$, there are many possible proofs here, which are all instructive. Here is our take on the subject, in the spirit of our proof for $SU_2$:

\medskip

(1) Our first claim is that we have the following formula, telling us that the moments of the main character equal the Catalan numbers:
$$\int_{SO_3}\chi^k=C_k$$

But this is something that we know from chapter 12, coming from Euler-Rodrigues. Alternatively, this can be deduced as well from Tannakian duality, a bit as for $SU_2$.

\medskip

(2) Our claim now is that we can construct, by recurrence on $k\in\mathbb N$, a sequence $r_k$ of irreducible, self-adjoint and distinct representations of $SO_3$, satisfying:
$$r_0=1\quad,\quad
r_1=u-1\quad,\quad 
r_k+r_{k-1}+r_{k-2}=r_{k-1}\otimes r_1$$

Indeed, assume that $r_0,\ldots,r_{k-1}$ are constructed, and let us construct $r_k$. The Frobenius trick from the proof for $SU_2$ will no longer work, due to some technical reasons, so we have to invoke (1). To be more precise, by integrating characters we obtain:
$$r_{k-1},r_{k-2}\subset r_{k-1}\otimes r_1$$

Thus there exists a representation $r_k$ such that:
$$r_{k-1}\otimes r_1=r_k+r_{k-1}+r_{k-2}$$

(3) Once again by integrating characters, we conclude that $r_k$ is irreducible, and non-equivalent to $r_1,\ldots,r_{k-1}$, and this proves our claim. Also, since any irreducible representation of $SO_3$ must appear in some tensor power of $u$, and we can decompose each $u^{\otimes k}$ into sums of representations $r_p$, we conclude that these representations $r_p$ are all the irreducible representations of $SO_3$. Finally, the dimension formula is clear.
\end{proof} 

There are of course many other things that can be said about $SU_2$ and $SO_3$. For instance, with the proof of Theorem 14.26 and Theorem 14.27 done in a purely algebraic fashion, by using the super-easiness property of $SU_2$ and $SO_3$, the Euler-Rodrigues formula can be deduced afterwards from this, without any single computation, the argument being that by Peter-Weyl the embedding $PU_2\subset SO_3$ must be indeed an equality.

\section*{14e. Exercises}

With the technology presented above, we can work out a few interesting particular cases of the Tannakian duality. Let us start with something quite elementary:

\begin{exercise}
Check the Brauer theorems for $O_N,U_N$, which are both of type
$$Hom(u^{\otimes k},u^{\otimes l})=span\left(T_\pi\Big|\pi\in D(k,l)\right)$$
for small values of the global length parameter, $k+l\in\{1,2,3\}$.
\end{exercise}

The idea here is to prove these results that we already know directly, by double inclusion, with the inclusion in one sense being something quite elementary.

\begin{exercise}
Write down Brauer theorems for the bistochastic groups 
$$B_N\subset O_N\quad,\quad C_N\subset U_N$$
by identifying first the partition which produces them, as subgroups of $O_N,U_N$.
\end{exercise}

This is actually something that will be discussed later on in this book, but without too much details, so the answer ``done in the book'' will not do.

\begin{exercise}
Look up the original version of Tannakian duality, stating that $G$ can be recovered from the knowledge of its full category of representations $\mathcal R_G$, viewed as subcategory of the category $\mathcal H$ of the finite dimensional Hilbert spaces, with each $\pi\in\mathcal R_G$ corresponding to its Hilbert space $H_\pi\in\mathcal H$, and write down a brief account of this.
\end{exercise}

As already mentioned in the above, the idea is that the group $G$ appears as the group of endomorphisms of the embedding functor $\mathcal R_G\subset\mathcal H$. Time to understand this.

\begin{exercise}
Look up the Doplicher-Roberts and Deligne theorems, stating that the compact group $G$ can be in fact recovered from the sole knowledge of the category $\mathcal R_G$, with no need for the embedding into $\mathcal H$, and write down a brief account of this. 
\end{exercise}

This is obviously something more advanced, and the proof is quite tricky.

\begin{exercise}
Given a closed subgroup $G\subset_uU_N$, understand and then briefly explain, in a short piece of writing, why the $*$-algebras
$$C(k,k)=End(u^{\otimes k})$$
form a planar algebra in the sense of Jones, and then comment as well on the various formulations of Tannakian duality, in the planar algebra setting.
\end{exercise}

This is actually quite difficult. And as a final, bonus exercise, try learning as well some Lie algebras, and their relation with the above, and report on what you learned.

\chapter{Diagrams, easiness}

\section*{15a. Easy groups}

We have seen in the previous chapter that the Tannakian duals of the groups $O_N,U_N$ are very simple objects. To be more precise, the Brauer theorem for these two groups states that we have equalities as follows, with $D=P_2,\mathcal P_2$ respectively:
$$Hom(u^{\otimes k},u^{\otimes l})=span\left(T_\pi\Big|\pi\in D(k,l)\right)$$

Our goal here will be that of axiomatizing and studying the closed subgroups $G\subset U_N$ which are of this type, but with $D$ being allowed to be, more generally, a category of partitions. Let us start our discussion with the following key definition:

\index{Kronecker symbols}
\index{maps associated to partitions}

\begin{definition}
Given a partition $\pi\in P(k,l)$ and an integer $N\in\mathbb N$, we define
$$T_\pi:(\mathbb C^N)^{\otimes k}\to(\mathbb C^N)^{\otimes l}$$
by the following formula, with $e_1,\ldots,e_N$ being the standard basis of $\mathbb C^N$,
$$T_\pi(e_{i_1}\otimes\ldots\otimes e_{i_k})=\sum_{j_1\ldots j_l}\delta_\pi\begin{pmatrix}i_1&\ldots&i_k\\ j_1&\ldots&j_l\end{pmatrix}e_{j_1}\otimes\ldots\otimes e_{j_l}$$
and with the coefficients on the right being Kronecker type symbols.
\end{definition}

To be more precise here, in order to compute the Kronecker type symbols $\delta_\pi(^i_j)\in\{0,1\}$, we proceed exactly as in the pairing case, namely by putting the multi-indices $i=(i_1,\ldots,i_k)$ and $j=(j_1,\ldots,j_l)$ on the legs of $\pi$, in the obvious way. In case all the blocks of $\pi$ contain equal indices of $i,j$, we set $\delta_\pi(^i_j)=1$. Otherwise, we set $\delta_\pi(^i_j)=0$.

\bigskip

With the above notion in hand, we can now formulate the following key definition, from \cite{bsp}, motivated by the Brauer theorems for $O_N,U_N$, as indicated before:

\index{easy group}

\begin{definition}
A closed subgroup $G\subset U_N$ is called easy when
$$Hom(u^{\otimes k},u^{\otimes l})=span\left(T_\pi\Big|\pi\in D(k,l)\right)$$
for any two colored integers $k,l=\circ\bullet\circ\bullet\ldots\,$, for certain sets of partitions
$$D(k,l)\subset P(k,l)$$
where $\pi\to T_\pi$ is the standard implementation of the partitions, as linear maps.
\end{definition}

In other words, we call a group $G$ easy when its Tannakian category appears in the simplest possible way: from the linear maps associated to partitions. The terminology is quite natural, because Tannakian duality is basically our only serious tool.

\bigskip

As basic examples, the orthogonal and unitary groups $O_N,U_N$ are both easy, coming respectively from the following collections of sets of partitions:
$$P_2=\bigsqcup_{k,l}P_2(k,l)\quad,\quad 
\mathcal P_2=\bigsqcup_{k,l}\mathcal P_2(k,l)$$

In the general case now, as an important theoretical remark, in the context of Definition 15.2, consider the following collection of sets of partitions:
$$D=\bigsqcup_{k,l}D(k,l)$$

This collection of sets $D$ determines $G$, but the converse is not true. Indeed, at $N=1$ for instance, both $D=P_2,\mathcal P_2$ produce the same easy group, namely $G=\{1\}$.

\bigskip

Coming next, again inspired from what we did in chapter 14, let us formulate:

\index{category of partitions}

\begin{definition}
Let $P(k,l)$ be the set of partitions between an upper colored integer $k$, and a lower colored integer $l$. A collection of subsets 
$$D=\bigsqcup_{k,l}D(k,l)$$
with $D(k,l)\subset P(k,l)$ is called a category of partitions when it has the following properties:
\begin{enumerate}
\item Stability under the horizontal concatenation, $(\pi,\sigma)\to[\pi\sigma]$.

\item Stability under vertical concatenation $(\pi,\sigma)\to[^\sigma_\pi]$, with matching middle symbols.

\item Stability under the upside-down turning $*$, with switching of colors, $\circ\leftrightarrow\bullet$.

\item Each set $P(k,k)$ contains the identity partition $||\ldots||$.

\item The sets $P(\emptyset,\circ\bullet)$ and $P(\emptyset,\bullet\circ)$ both contain the semicircle $\cap$.

\item The sets $P(k,\bar{k})$ with $|k|=2$ contain the crossing partition $\slash\hskip-2.0mm\backslash$.
\end{enumerate}
\end{definition} 

As before, this is something that we already met in chapter 14, but for the pairings only. Observe the similarity with the axioms for Tannakian categories, also from chapter 14. We will see in a moment that this similarity can be turned into something very precise, the idea being that such a category produces a family of easy quantum groups $(G_N)_{N\in\mathbb N}$, one for each $N\in\mathbb N$, via the formula in Definition 15.1, and Tannakian duality.

\bigskip

As basic examples, that we have already met in chapter 14, in connection with the representation theory of $O_N,U_N$, we have the categories $P_2,\mathcal P_2$ of pairings, and of matching pairings. Further basic examples include the categories $P,P_{even}$ of all partitions, and of all partitions whose blocks have even size. We will see in a moment that these latter categories are related to the symmetric and hyperoctahedral groups $S_N,H_N$.
 
\bigskip

The relation with the Tannakian categories comes from the following result:

\begin{proposition}
The assignement $\pi\to T_\pi$ is categorical, in the sense that
$$T_\pi\otimes T_\sigma=T_{[\pi\sigma]}\quad,\quad 
T_\pi T_\sigma=N^{c(\pi,\sigma)}T_{[^\sigma_\pi]}\quad,\quad 
T_\pi^*=T_{\pi^*}$$
where $c(\pi,\sigma)$ are certain integers, coming from the erased components in the middle.
\end{proposition}

\begin{proof}
This is something that we already know for the pairings, from chapter 14, and the proof in general is similar, with the only axiom where some slight changes appear being the composition one. Here the computation is as follows, as before for pairings, with $c(\pi,\sigma)\in\mathbb N$ counting the middle components, which are not necessarily circles:
\begin{eqnarray*}
&&T_\pi T_\sigma(e_{i_1}\otimes\ldots\otimes e_{i_p})\\
&=&\sum_{j_1\ldots j_q}\delta_\sigma\begin{pmatrix}i_1&\ldots&i_p\\j_1&\ldots&j_q\end{pmatrix}
\sum_{k_1\ldots k_r}\delta_\pi\begin{pmatrix}j_1&\ldots&j_q\\k_1&\ldots&k_r\end{pmatrix}e_{k_1}\otimes\ldots\otimes e_{k_r}\\
&=&\sum_{k_1\ldots k_r}N^{c(\pi,\sigma)}\delta_{[^\sigma_\pi]}\begin{pmatrix}i_1&\ldots&i_p\\k_1&\ldots&k_r\end{pmatrix}e_{k_1}\otimes\ldots\otimes e_{k_r}\\
&=&N^{c(\pi,\sigma)}T_{[^\sigma_\pi]}(e_{i_1}\otimes\ldots\otimes e_{i_p})
\end{eqnarray*}

Thus, our correspondence is indeed categorical, as claimed.
\end{proof}

Time now to put everyting together. All the above was pure combinatorics, and in relation with the compact groups, we have the following result:

\begin{theorem}
Each category of partitions $D=(D(k,l))$ produces a family of compact groups $G=(G_N)$, one for each $N\in\mathbb N$, via the formula
$$Hom(u^{\otimes k},u^{\otimes l})=span\left(T_\pi\Big|\pi\in D(k,l)\right)$$
and the Tannakian duality correspondence.
\end{theorem}

\begin{proof}
Given an integer $N\in\mathbb N$, consider the correspondence $\pi\to T_\pi$ constructed in Definition 15.1, and then the collection of linear spaces in the statement, namely:
$$C_{kl}=span\left(T_\pi\Big|\pi\in D(k,l)\right)$$

According to the formulae in Proposition 15.4, and to our axioms for the categories of partitions, from Definition 15.3, this collection of spaces $C=(C_{kl})$ satisfies the axioms for the Tannakian categories, from chapter 14. Thus the Tannakian duality result there applies, and provides us with a closed subgroup $G_N\subset U_N$ such that:
$$C_{kl}=Hom(u^{\otimes k},u^{\otimes l})$$

Thus, we are led to the conclusion in the statement.
\end{proof}

In relation with the easiness property, we can now formulate a key result, which can serve as an alternative definition for the easy groups, as follows:

\index{easy group}

\begin{theorem}
A closed subgroup $G\subset U_N$ is easy precisely when
$$Hom(u^{\otimes k},u^{\otimes l})=span\left(T_\pi\Big|\pi\in D(k,l)\right)$$
for any colored integers $k,l$, for a certain category of partitions $D\subset P$.
\end{theorem}

\begin{proof}
This basically follows from Theorem 15.5, as follows:

\medskip

(1) In one sense, we know from Theorem 15.5 that any category of partitions $D\subset P$ produces a family of closed groups $G\subset U_N$, one for each $N\in\mathbb N$, according to Tannakian duality and to the Hom space formula there, namely:
$$Hom(u^{\otimes k},u^{\otimes l})=span\left(T_\pi\Big|\pi\in D(k,l)\right)$$

But these groups $G\subset U_N$ are indeed easy, in the sense of Definition 15.2.

\medskip

(2) In the other sense now, assume that $G\subset U_N$ is easy, in the sense of Definition 15.2, coming via the above Hom space formula, from a collection of sets as follows:
$$D=\bigsqcup_{k,l}D(k,l)$$

Consider now the category of partitions $\widetilde{D}=<D>$ generated by this family. This is by definition the smallest category of partitions containing $D$, whose existence follows by starting with $D$, and performing the various categorical operations, namely horizontal and vertical concatenation, and upside-down turning. It follows then, via another application of Tannakian duality, that we have the following formula, for any $k,l$:
$$Hom(u^{\otimes k},u^{\otimes l})=span\left(T_\pi\Big|\pi\in\widetilde{D}(k,l)\right)$$

Thus, our group $G\subset U_N$ can be viewed as well as coming from $\widetilde{D}$, and so appearing as particular case of the construction in Theorem 15.5, and this gives the result. 
\end{proof}

As already mentioned above, Theorem 15.6 can be regarded as an alternative definition for easiness, with the assumption that $D\subset P$ must be a category of partitions being added. In what follows we will rather use this new definition, which is more precise.

\bigskip

The notion of easiness goes back to the results of Brauer in \cite{bra} regarding the orthogonal group $O_N$, and the unitary group $U_N$, which reformulate as follows:

\index{Brauer theorem}
\index{easy group}
\index{orthogonal group}
\index{unitary group}

\begin{theorem}
We have the following results:
\begin{enumerate}
\item The unitary group $U_N$ is easy, coming from the category $\mathcal P_2$.

\item The orthogonal group $O_N$ is easy as well, coming from the category $P_2$.
\end{enumerate}
\end{theorem}

\begin{proof}
This is something that we already know, from chapter 14, based on Tannakian duality, the idea of the proof being as follows:

\medskip

(1) The group $U_N$ being defined via the relations $u^*=u^{-1}$, $u^t=\bar{u}^{-1}$, the associated Tannakian category is $C=span(T_\pi|\pi\in D)$, with:
$$D
=<{\ }^{\,\cap}_{\circ\bullet}\,\,,{\ }^{\,\cap}_{\bullet\circ}>
=\mathcal P_2$$

(2) The group $O_N\subset U_N$ being defined by imposing the relations $u_{ij}=\bar{u}_{ij}$, the associated Tannakian category is $C=span(T_\pi|\pi\in D)$, with:
$$D
=<\mathcal P_2,|^{\hskip-1.32mm\circ}_{\hskip-1.32mm\bullet},|_{\hskip-1.32mm\circ}^{\hskip-1.32mm\bullet}>
=P_2$$
  
Thus, we are led to the conclusion in the statement.
\end{proof}

There are many other examples of easy groups, and we will gradually explore this. To start with, we have the following result, dealing with the groups $B_N,C_N$:

\index{Brauer theorem}
\index{easy group}
\index{bistochastic group}

\begin{theorem}
We have the following results:
\begin{enumerate}
\item The unitary bistochastic group $C_N$ is easy, coming from the category $\mathcal P_{12}$ of matching singletons and pairings.

\item The orthogonal bistochastic group $B_N$ is easy, coming from the category $P_{12}$ of singletons and pairings.
\end{enumerate}
\end{theorem}

\begin{proof}
The proof here is similar to the proof of Theorem 15.7. To be more precise, we can use the results there, and the proof goes as follows:

\medskip

(1) The group $C_N\subset U_N$ is defined by imposing the following relations, with $\xi$ being the all-one vector, which correspond to the bistochasticity condition:
$$u\xi=\xi\quad,\quad 
\bar{u}\xi=\xi$$

But these relations tell us precisely that the following two operators, with the partitions on the right being singletons, must be in the associated Tannakian category $C$:
$$T_\pi\quad:\quad\pi=|_{\hskip-1.32mm\circ}\ ,\ |_{\hskip-1.32mm\bullet}$$

Thus the associated Tannakian category is $C=span(T_\pi|\pi\in D)$, with:
$$D
=<\mathcal P_2,|_{\hskip-1.32mm\circ},|_{\hskip-1.32mm\bullet}>
=\mathcal P_{12}$$
  
Thus, we are led to the conclusion in the statement.

\medskip

(2) In order to deal now with the real bistochastic group $B_N$, we can either use a similar argument, or simply use the following intersection formula: 
$$B_N=C_N\cap O_N$$

Indeed, at the categorical level, this intersection formula tells us that the associated Tannakian category is given by $C=span(T_\pi|\pi\in D)$, with:
$$D
=<\mathcal P_{12},P_2>
=P_{12}$$

Thus, we are led to the conclusion in the statement.
\end{proof}

As a comment here, we have used in the above the fact, which is something quite trivial, that the category of partitions associated to an intersection of easy quantum groups is generated by the corresponding categories of partitions. We will be back to this, and to some other product operations as well, with similar results, later on.

\bigskip

We can put now the results that we have together, as follows:

\begin{theorem}
The basic unitary and bistochastic groups,
$$\xymatrix@R=50pt@C=50pt{
C_N\ar[r]&U_N\\
B_N\ar[u]\ar[r]&O_N\ar[u]}$$
are all easy, coming from the various categories of singletons and pairings.
\end{theorem}

\begin{proof}
We know from the above that the groups in the statement are indeed easy, the corresponding diagram of categories of partitions being as follows:
$$\xymatrix@R=16mm@C=18mm{
\mathcal P_{12}\ar[d]&\mathcal P_2\ar[l]\ar[d]\\
P_{12}&P_2\ar[l]}$$

Thus, we are led to the conclusion in the statement.
\end{proof}

Summarizing, what we have so far is a general notion of easiness, coming from the Brauer theorems for $O_N,U_N$, and their straightforward extensions to $B_N,C_N$.

\section*{15b. Reflection groups}

In view of the above, the notion of easiness is a quite interesting one, deserving a full, systematic investigation. As a first natural question that we would like to solve, we would like to compute the easy group associated to the category of all partitions $P$ itself. And here, no surprise, we are led to the most basic, but non-trivial, classical group that we know, namely the symmetric group $S_N$. To be more precise, we have the following Brauer type theorem for $S_N$, which answers our question formulated above:

\index{Brauer theorem}
\index{symmetric group}

\begin{theorem}
The symmetric group $S_N$, regarded as group of unitary matrices,
$$S_N\subset O_N\subset U_N$$
via the permutation matrices, is easy, coming from the category of all partitions $P$.
\end{theorem}

\begin{proof}
Consider indeed the group $S_N$, regarded as a group of unitary matrices, with each permutation $\sigma\in S_N$ corresponding to the associated permutation matrix:
$$\sigma(e_i)=e_{\sigma(i)}$$

In order to prove the result, consider the one-block ``fork'' partition, namely:
$$\xymatrix@R=1mm@C=2mm{\\ \\ \mu\ \ =\\ \\ }\ \ \ 
\xymatrix@R=2mm@C=3mm{
\circ\ar@/_/@{-}[dr]&&\circ\\
&\ar@/_/@{-}[ur]\ar@{-}[dd]\\
&&&\\
&\circ}$$

The linear map associated to $\mu$ is then given by the following formula:
$$T_\mu(e_i\otimes e_j)=\delta_{ij}e_i$$

In order to do the computations, we use the following formulae:
$$u=(u_{ij})_{ij}\quad,\quad 
u^{\otimes 2}=(u_{ij}u_{kl})_{ik,jl}\quad,\quad 
T_\mu=(\delta_{ijk})_{i,jk}$$

By using these formulae, we obtain the following equality:
$$(T_\mu u^{\otimes 2})_{i,jk}
=\sum_{lm}(T_\mu)_{i,lm}(u^{\otimes 2})_{lm,jk}
=u_{ij}u_{ik}$$

On the other hand, we have as well the following equality:
$$(uT_\mu)_{i,jk}
=\sum_lu_{il}(T_\mu)_{l,jk}
=\delta_{jk}u_{ij}$$

We therefore conclude that we have an equivalence, as follows:
$$T_\mu\in Hom(u^{\otimes 2},u)\iff u_{ij}u_{ik}=\delta_{jk}u_{ij},\forall i,j,k$$

In other words, the elements $u_{ij}$ must be projections, which must be pairwise orthogonal on the rows of $u=(u_{ij})$. But this reformulates into the following equality:
$$C(S_N)=C(O_N)\Big/\Big<T_\mu\in Hom(u^{\otimes 2},u)\Big>$$

According now to our general conventions for easiness, this means that the symmetric group $S_N$ is easy, coming from the following category of partitions:
$$D=<\mu>=P$$

Thus, we are led to the conclusion in the statement.
\end{proof}

Next, regarding the hyperoctahedral group $H_N$, we have the following result:

\index{Brauer theorem}
\index{hyperoctahedral group}

\begin{theorem}
The hyperoctahedral group $H_N$, regarded as group of matrices,
$$H_N\subset O_N\subset U_N$$
is easy, coming from the category of partitions with even blocks $P_{even}$.
\end{theorem}

\begin{proof}
This follows as usual from Tannakian duality. To be more precise, consider the following one-block partition $\chi\in P(2,2)$, which looks like a $\chi$ letter:
$$\xymatrix@R=0.5mm@C=2mm{\\ \\ \\ \chi\ \ =\\ \\ }\ \ \ 
\xymatrix@R=2mm@C=3mm{
\circ\ar@/_/@{-}[dr]&&\circ\\
&\ar@/_/@{-}[ur]\ar@{-}[dd]\\
&&&\\
&\ar@/^/@{-}[dr]\ar@/_/@{-}[dl]\\
\circ&&\circ}$$

The linear map associated to this partition is then given by:
$$T_\chi(e_i\otimes e_j)=\delta_{ij}e_i\otimes e_i$$

By using this formula, we have the following computation:
\begin{eqnarray*}
(T_\chi\otimes id)u^{\otimes 2}(e_a\otimes e_b)
&=&(T_\chi\otimes id)\left(\sum_{ijkl}e_{ij}\otimes e_{kl}\otimes u_{ij}u_{kl}\right)(e_a\otimes e_b)\\
&=&(T_\chi\otimes id)\left(\sum_{ik}e_i\otimes e_k\otimes u_{ia}u_{kb}\right)\\
&=&\sum_ie_i\otimes e_i\otimes u_{ia}u_{ib}
\end{eqnarray*}

On the other hand, we have as well the following computation:
\begin{eqnarray*}
u^{\otimes 2}(T_\chi\otimes id)(e_a\otimes e_b)
&=&\delta_{ab}\left(\sum_{ijkl}e_{ij}\otimes e_{kl}\otimes u_{ij}u_{kl}\right)(e_a\otimes e_a)\\
&=&\delta_{ab}\sum_{ij}e_i\otimes e_k\otimes u_{ia}u_{ka}
\end{eqnarray*}

We conclude from this that we have the following equivalence:
$$T_\chi\in End(u^{\otimes 2})\iff \delta_{ik}u_{ia}u_{ib}=\delta_{ab}u_{ia}u_{ka},\forall i,k,a,b$$

But the relations on the right tell us that the entries of $u=(u_{ij})$ must satisfy $\alpha\beta=0$ on each row and column of $u$, and so that the corresponding closed subgroup $G\subset O_N$ consists of the matrices $g\in O_N$ which are permutation-like, with $\pm1$ nonzero entries. Thus, the corresponding group is $G=H_N$, and as a conclusion to this, we have:
$$C(H_N)=C(O_N)\Big/\Big<T_\chi\in End(u^{\otimes 2})\Big>$$

According now to our conventions for easiness, this means that the hyperoctahedral group $H_N$ is easy, coming from the following category of partitions:
$$D=<\chi>=P_{even}$$

Thus, we are led to the conclusion in the statement.
\end{proof}

Next, regarding the full reflection group $K_N$, we have the following result:

\begin{theorem}
The full reflection group $K_N=\mathbb T\wr S_N$, regarded as subgroup
$$K_N\subset U_N$$
comes from $\mathcal P_{even}$, the partitions satisfying $\#\circ=\#\bullet$, weighted equality, in each block. 
\end{theorem}

\begin{proof}
We are now dealing with unitary matrices, so we must use colored partitions. Consider the following partition $\chi\in P(\circ\bullet\,,\bullet\circ)$, that we already met above, uncolored:
$$\xymatrix@R=0.5mm@C=2mm{\\ \\ \\ \chi\ \ =\\ \\ }\ \ \ 
\xymatrix@R=2mm@C=3mm{
\circ\ar@/_/@{-}[dr]&&\bullet\\
&\ar@/_/@{-}[ur]\ar@{-}[dd]\\
&&&\\
&\ar@/^/@{-}[dr]\ar@/_/@{-}[dl]\\
\bullet&&\circ}$$

Our computations from the previous proof, for the group $H_N$, modify into:
$$(T_\chi\otimes id)(u\otimes\bar{u})(e_a\otimes e_b)=\sum_ie_i\otimes e_i\otimes u_{ia}\bar{u}_{ib}$$
$$(\bar{u}\otimes u)(T_\chi\otimes id)(e_a\otimes e_b)=\delta_{ab}\sum_{ij}e_i\otimes e_k\otimes\bar{u}_{ia}u_{ka}$$

We conclude from this that we have the following equivalence:
$$T_\chi\in Hom(u\otimes\bar{u},\bar{u}\otimes u)\iff \delta_{ik}u_{ia}\bar{u}_{ib}=\delta_{ab}\bar{u}_{ia}u_{ka},\forall i,k,a,b$$

But the relations on the right tell us that the entries of $u=(u_{ij})$ must satisfy $\alpha\beta=0$ on each row and column of $u$, and as a conclusion to this, we have:
$$C(K_N)=C(U_N)\Big/\Big<T_\chi\in Hom(u\otimes\bar{u},\bar{u}\otimes u)\Big>$$

Thus the group $K_N$ is easy, coming from the following category of partitions:
$$D=<\chi>=\mathcal P_{even}$$

We are therefore led to the conclusion in the statement.
\end{proof}

More generally now, we have in fact the following grand result:

\index{Brauer theorem}
\index{complex reflection group}

\begin{theorem}
The complex reflection group $H_N^s=\mathbb Z_s\wr S_N$ is easy, the corresponding category $P^s$ consisting of the partitions satisfying the condition
$$\#\circ=\#\bullet(s)$$
as a weighted sum, in each block. In particular, we have the following results:
\begin{enumerate}
\item $S_N$ is easy, coming from the category $P$.

\item $H_N=\mathbb Z_2\wr S_N$ is easy, coming from the category $P_{even}$.

\item $K_N=\mathbb T\wr S_N$ is easy, coming from the category $\mathcal P_{even}$.
\end{enumerate}
\end{theorem}

\begin{proof}
This is something coming at $s=1,2,\infty$ from Theorems 15.10, 15.11 and 15.12, as indicated in (1,2,3), with this to be discussed in a moment, and in general, the proof is similar. Consider indeed the following partition, with $s+2$ legs:
$$\xymatrix@R=6pt@C=12pt{
&&\ar@{-}[dd]\ar@{-}[rrrrr]&\ar@{-}[dd]&&\ar@{-}[dd]&\ar@{-}[dd]&\ar@{-}[dd]\\
\xi&=&&&\ldots\\
&&\circ&\circ&&\circ&\circ&\bullet}$$

Observe that, up to rotation and some discussion regarding the colors, this coincides with the partitions $\mu,\chi$ that we used before at $s=1,2$. In general now, we have:
$$T_\xi=\sum_je_j^{\otimes s+2}$$

Our claim, which will prove the result, is that we have the following formula:
$$C(H_N^s)=C(K_N)\Big/\Big<T_\xi\in Fix(u^{\otimes s+1}\otimes\bar{u})\Big>$$

Indeed, by using the above formula of $T_\xi$, we have the following computation:
\begin{eqnarray*}
(u^{\otimes s+1}\otimes\bar{u})(T_\xi\otimes 1)
&=&\sum_{ij}e_{i_1}\otimes\ldots\otimes e_{i_{s+2}}\otimes u_{i_1j}\ldots u_{i_{s+1}j}\bar{u}_{i_{s+2}j}\\
&=&\sum_{ij}e_i\otimes\ldots\otimes e_i\otimes u_{ij}^{s+1}\bar{u}_{ij}\\
&=&\sum_ie_i^{\otimes s+2}\otimes\left(\sum_ju_{ij}^{s+1}\bar{u}_{ij}\right)
\end{eqnarray*}

We conclude that, for a subgroup of $K_N$, we have the following equivalence:
$$T_\xi\in Fix(u^{\otimes s+1}\otimes\bar{u})\iff\sum_ju_{ij}^{s+1}\bar{u}_{ij}=1$$

Now the conditions on the right being those defining the subgroup $H_N^s\subset K_N$, we conclude that we have the equality announced above, namely:
$$C(H_N^s)=C(K_N)\Big/\Big<T_\xi\in Fix(u^{\otimes s+1}\otimes\bar{u})\Big>$$

But with this, we can finish the proof of the main assertion. Indeed, it follows that the group $H_N^s$ is easy, coming from the following category of partitions:
$$D=<\mathcal P_{even},\xi>=P^s$$

Summarizing, theorem proved, and in what regards the particular cases, which generalize what we knew from Theorems 15.10, 15.11 and 15.12, these are as follows:

\medskip

(1) At $s=1$ we know that we have $H_N^1=S_N$. Regarding now the corresponding category, here the condition $\#\circ=\#\bullet(1)$ is automatic, and so $P^1=P$.

\medskip

(2) At $s=2$ we know that we have $H_N^2=H_N$. Regarding now the corresponding category, here the condition $\#\circ=\#\bullet(2)$ reformulates as follows:
$$\#\circ+\,\#\bullet=0(2)$$

Thus each block must have even size, and we obtain, as claimed, $P^2=P_{even}$.

\medskip

(3) At $s=\infty$ we know that we have $H_N^\infty=K_N$. Regarding now the corresponding category, here the condition $\#\circ=\#\bullet(\infty)$ reads:
$$\#\circ=\#\bullet$$

But this is the condition defining $\mathcal P_{even}$, and so $P^\infty=\mathcal P_{even}$, as claimed.
\end{proof}

Summarizing, we have many examples. In fact, our list of easy groups has currently become quite big, and here is a selection of the main results that we have so far: 

\begin{theorem}
We have a diagram of compact groups as follows,
$$\xymatrix@R=50pt@C=50pt{
K_N\ar[r]&U_N\\
H_N\ar[u]\ar[r]&O_N\ar[u]}$$
where $H_N=\mathbb Z_2\wr S_N$ and $K_N=\mathbb T\wr S_N$, and all these groups are easy.
\end{theorem}

\begin{proof}
This follows from the above results. To be more precise, we know that the above groups are all easy, the corresponding categories of partitions being as follows:
$$\xymatrix@R=16mm@C=18mm{
\mathcal P_{even}\ar[d]&\mathcal P_2\ar[l]\ar[d]\\
P_{even}&P_2\ar[l]}$$

Thus, we are led to the conclusion in the statement.
\end{proof}

Summarizing, most of the groups that we investigated in this book are covered by the easy group formalism. One exception is the symplectic group $Sp_N$, but this group is covered as well, by a suitable extension of the easy group formalism. See \cite{csn}.

\section*{15c. Basic operations} 

All the above is quite encouraging, so time now to take easiness very seriously, and develop some general abstract theory for the easy groups. Let us first discuss some basic composition operations. We will be mainly interested in the following operations:

\begin{definition}
The closed subgroups of $U_N$ are subject to intersection and generation operations, constructed as follows:
\begin{enumerate}
\item Intersection: $H\cap K$ is the usual intersection of $H,K$.

\item Generation: $<H,K>$ is the closed subgroup generated by $H,K$.
\end{enumerate}
\end{definition}

Alternatively, we can define these operations at the function algebra level, by performing certain operations on the associated ideals, as follows:

\begin{proposition}
Assuming that we have presentation results as follows,
$$C(H)=C(U_N)/I\quad,\quad 
C(K)=C(U_N)/J$$
the groups $H\cap K$ and $<H,K>$ are given by the following formulae,
$$C(H\cap K)=C(U_N)/<I,J>$$
$$C(<H,K>)=C(U_N)/(I\cap J)$$
at the level of the associated algebras of functions.
\end{proposition}

\begin{proof}
This is indeed clear from the definition of the operations $\cap$ and $<\,,\,>$, as formulated above, and from the Stone-Weierstrass theorem.
\end{proof}

In what follows we will need Tannakian formulations of the above two operations. The result here, coming from the general Tannakian duality result established in chapter 14, and that we have in fact already used a couple of times in the above, is as follows:

\begin{theorem}
The intersection and generation operations $\cap$ and $<\,,>$ can be constructed via the Tannakian correspondence $G\to C_G$, as follows:
\begin{enumerate}
\item Intersection: defined via $C_{G\cap H}=<C_G,C_H>$.

\item Generation: defined via $C_{<G,H>}=C_G\cap C_H$.
\end{enumerate}
\end{theorem}

\begin{proof}
This follows from Proposition 15.16, and from Tannakian duality. Indeed, it follows from Tannakian duality that given a closed subgroup $G\subset U_N$, with fundamental representation $v$, the algebra of functions $C(G)$ has the following presentation:
$$C(G)=C(U_N)\Big/\left<T\in Hom(u^{\otimes k},u^{\otimes l})\Big|\forall k,\forall l,\forall T\in Hom(v^{\otimes k},v^{\otimes l})\right>$$

In other words, given a closed subgroup $G\subset U_N$, we have a presentation of the following type, with $I_G$ being the ideal coming from the Tannakian category of $G$:
$$C(G)=C(U_N)/I_G$$

But this leads to the conclusion in the statement.
\end{proof}

In relation now with our easiness questions, we first have the following result:

\begin{proposition}
Assuming that $H,K$ are easy, then so is $H\cap K$, and we have
$$D_{H\cap K}=<D_H,D_K>$$
at the level of the corresponding categories of partitions.
\end{proposition}

\begin{proof}
We have indeed the following computation:
\begin{eqnarray*}
C_{H\cap K}
&=&<C_H,C_K>\\
&=&<span(D_H),span(D_K)>\\
&=&span(<D_H,D_K>)
\end{eqnarray*}

Thus, by Tannakian duality we obtain the result.
\end{proof}

Regarding now the generation operation, the situation here is more complicated, due to a number of technical reasons, and we only have the following statement:

\begin{proposition}
Assuming that $H,K$ are easy, we have an inclusion 
$$<H,K>\subset\{H,K\}$$
coming from an inclusion of Tannakian categories as follows,
$$C_H\cap C_K\supset span(D_H\cap D_K)$$
where $\{H,K\}$ is the easy group having as category of partitions $D_H\cap D_K$.
\end{proposition}

\begin{proof}
This follows from the definition and properties of the generation operation, explained above, and from the following computation:
\begin{eqnarray*}
C_{<H,K>}
&=&C_H\cap C_K\\
&=&span(D_H)\cap span(D_K)\\
&\supset&span(D_H\cap D_K)
\end{eqnarray*}

Indeed, by Tannakian duality we obtain from this all the assertions.
\end{proof}

It is not clear when the inclusions in Proposition 15.19 are isomorphisms or not, and this even under a supplementary $N>>0$ assumption. Technically speaking, the problem comes from the fact that the operation $\pi\to T_\pi$ does not produce linearly independent maps, and so all that we are doing is sensitive to the value of $N\in\mathbb N$. The subject here is quite technical, to be further developed in chapter 16 below, with probabilistic motivations in mind, without however solving the present algebraic questions.

\bigskip

Summarizing, we have some problems here, and we must proceed as follows:

\index{easy generation}

\begin{theorem}
The intersection and easy generation operations $\cap$ and $\{\,,\}$ can be constructed via the Tannakian correspondence $G\to D_G$, as follows:
\begin{enumerate}
\item Intersection: defined via $D_{G\cap H}=<D_G,D_H>$.

\item Easy generation: defined via $D_{\{G,H\}}=D_G\cap D_H$.
\end{enumerate}
\end{theorem}

\begin{proof}
Here the situation is as follows:

\medskip

(1) This is a true and honest result, coming from Proposition 15.18.

\medskip

(2) This is more of an empty statement, coming from Proposition 15.19.
\end{proof}

As already mentioned, there is some interesting mathematics still to be worked out, in relation with all this, and we will be back to this later, with further details. With the above notions in hand, however, even if not fully satisfactory, we can formulate a nice result, which improves our main result so far, namely Theorem 15.14, as follows:

\begin{theorem}
The basic unitary and reflection groups, namely
$$
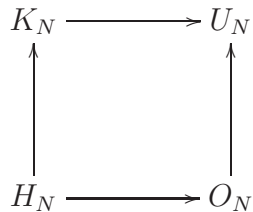
\xymatrix@R=50pt@C=50pt{
K_N\ar[r]&U_N\\
H_N\ar[u]\ar[r]&O_N\ar[u]}$$
are all easy, and they form an intersection and easy generation diagram, in the sense that the above square diagram satisfies $U_N=\{K_N,O_N\}$, and $H_N=K_N\cap O_N$.
\end{theorem}

\begin{proof}
We know from Theorem 15.14 that the groups in the statement are easy, the corresponding categories of partitions being as follows:
$$
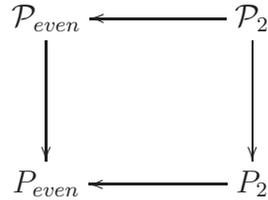
\xymatrix@R=16mm@C=18mm{
\mathcal P_{even}\ar[d]&\mathcal P_2\ar[l]\ar[d]\\
P_{even}&P_2\ar[l]}$$

Now observe that this latter diagram is an intersection and generation diagram. By using Theorem 15.20, this reformulates into the fact that the corresponding diagram of groups is an intersection and easy generation diagram, as claimed.
\end{proof}

It is possible to further improve the above result, by proving that the diagram there is actually a plain generation diagram. However, this is something more technical, and for a discussion here, you can check for instance my group theory book \cite{ba3}.

\bigskip

Moving forward, as a continuation of the above, it is possible to develop some more general theory, along the above lines. Given a closed subgroup $G\subset U_N$, we can talk about its ``easy envelope'', which is the smallest easy group $\widetilde{G}$ containing $G$. This easy envelope appears by definition as an intermediate closed subgroup, as follows:
$$G\subset\widetilde{G}\subset U_N$$

With this notion in hand, Proposition 15.19 can be refined into a result stating that given two easy groups $H,K$, we have inclusions as follows:
$$<H,K>\subset\widetilde{<H,K>}\subset\{H,K\}$$

In order to discuss all this, let us start with the following definition:

\index{homogeneous group}

\begin{definition}
A closed subgroup $G\subset U_N$ is called homogeneous when
$$S_N\subset G\subset U_N$$
with $S_N\subset U_N$ being the standard embedding, via permutation matrices.
\end{definition}

We will be interested in such groups, which cover for instance all the easy groups, and many more. At the Tannakian level, we have the following result:

\begin{theorem}
The homogeneous groups $S_N\subset G\subset U_N$ are in one-to-one correspondence with the intermediate tensor categories
$$span\left(T_\pi\Big|\pi\in\mathcal P_2\right)\subset C\subset span\left(T_\pi\Big|\pi\in P\right)$$
where $P$ is the category of all partitions, $\mathcal P_2$ is the category of the matching pairings, and $\pi\to T_\pi$ is the standard implementation of partitions, as linear maps.
\end{theorem}

\begin{proof}
This follows from Tannakian duality, and from the Brauer type results for $S_N,U_N$. To be more precise, we know from Tannakian duality that each closed subgroup $G\subset U_N$ can be reconstructed from its Tannakian category $C=(C(k,l))$, as follows:
$$C(G)=C(U_N)\Big/\left<T\in Hom(u^{\otimes k},u^{\otimes l})\Big|\forall k,l,\forall T\in C(k,l)\right>$$

Thus we have a one-to-one correspondence $G\leftrightarrow C$, given by Tannakian duality, and since the endpoints $G=S_N,U_N$ are both easy, corresponding to the categories $C=span(T_\pi|\pi\in D)$ with $D=P,\mathcal P_2$, this gives the result.
\end{proof}

Our purpose now will be that of using the Tannakian result in Theorem 15.23, in order to introduce and study a combinatorial notion of ``easiness level'', for the arbitrary intermediate groups $S_N\subset G\subset U_N$.  Let us begin with the following simple fact:

\begin{proposition}
Given a homogeneous group $S_N\subset G\subset U_N$, with associated Tannakian category $C=(C(k,l))$, the sets
$$D^1(k,l)=\left\{\pi\in P(k,l)\Big|T_\pi\in C(k,l)\right\}$$ 
form a category of partitions, in the sense of Definition 15.3.
\end{proposition}

\begin{proof}
We use the basic categorical properties of the correspondence $\pi\to T_\pi$ between partitions and linear maps, that we established in the above, namely:
$$T_{[\pi\sigma]}=T_\pi\otimes T_\sigma\quad,\quad 
T_{[^\sigma_\pi]}\sim T_\pi T_\sigma\quad,\quad 
T_{\pi^*}=T_\pi^*$$

Together with the fact that $C$ is a tensor category, we deduce from these formulae that we have the following implication:
\begin{eqnarray*}
\pi,\sigma\in D^1
&\implies&T_\pi,T_\sigma\in C\\
&\implies&T_\pi\otimes T_\sigma\in C\\
&\implies&T_{[\pi\sigma]}\in C\\
&\implies&[\pi\sigma]\in D^1
\end{eqnarray*}

On the other hand, we have as well the following implication:
\begin{eqnarray*}
\pi,\sigma\in D^1
&\implies&T_\pi,T_\sigma\in C\\
&\implies&T_\pi T_\sigma\in C\\
&\implies&T_{[^\sigma_\pi]}\in C\\
&\implies&[^\sigma_\pi]\in D^1
\end{eqnarray*}

Finally, we have as well the following implication:
\begin{eqnarray*}
\pi\in D^1
&\implies&T_\pi\in C\\
&\implies&T_\pi^*\in C\\
&\implies&T_{\pi^*}\in C\\
&\implies&\pi^*\in D^1
\end{eqnarray*}

Thus $D^1$ is indeed a category of partitions, as claimed.
\end{proof}

We can further refine the above observation, in the following way:

\begin{proposition}
Given a compact group $S_N\subset G\subset U_N$, construct $D^1\subset P$ as above, and let $S_N\subset G^1\subset U_N$ be the easy group associated to $D^1$. Then:
\begin{enumerate}

\item We have $G\subset G^1$, as subgroups of $U_N$.

\item $G^1$ is the smallest easy group containing $G$.

\item $G$ is easy precisely when $G\subset G^1$ is an isomorphism.
\end{enumerate}
\end{proposition}

\begin{proof}
All this is elementary, the proofs being as follows:

\medskip

(1) We know that the Tannakian category of $G^1$ is given by:
$$C_{kl}^1=span\left(T_\pi\Big|\pi\in D^1(k,l)\right)$$

Thus we have $C^1\subset C$, and so $G\subset G^1$, as subgroups of $U_N$.

\medskip

(2) Assuming that we have $G\subset G'$, with $G'$ easy, coming from a Tannakian category $C'=span(D')$, we must have $C'\subset C$, and so $D'\subset D^1$. Thus, $G^1\subset G'$, as desired.

\medskip

(3) This is a trivial consequence of (2).
\end{proof}

Summarizing, we have now a notion of ``easy envelope'', as follows:

\index{easy envelope}

\begin{definition}
The easy envelope of a homogeneous group $S_N\subset G\subset U_N$ is the easy group $S_N\subset G^1\subset U_N$ associated to the category of partitions
$$D^1(k,l)=\left\{\pi\in P(k,l)\Big|T_\pi\in C(k,l)\right\}$$ 
where $C=(C(k,l))$ is the Tannakian category of $G$.
\end{definition}

At the level of examples, most of the known homogeneous groups $S_N\subset G\subset U_N$ are in fact easy. However, there are non-easy interesting examples as well, such as the generic reflection groups $H_N^{sd}$ from chapter 12, and we will certainly have an exercise at the end of this chapter, regarding the computation of the corresponding easy envelopes. 

\bigskip

As a technical observation now, we can in fact generalize the above construction to any closed subgroup $G\subset U_N$, and we have the following result:

\begin{proposition}
Given a closed subgroup $G\subset U_N$, construct $D^1\subset P$ as above, and let $S_N\subset G^1\subset U_N$ be the easy group associated to $D^1$. We have then
$$G^1=(<G,S_N>)^1$$
where $<G,S_N>\subset U_N$ is the smallest closed subgroup containing $G,S_N$.
\end{proposition}

\begin{proof}
According to our Tannakian results, the subgroup $<G,S_N>\subset U_N$ in the statement exists indeed, and can be obtained by intersecting categories, as follows:
$$C_{<G,S_N>}=C_G\cap C_{S_N}$$

We conclude from this that for any $\pi\in P(k,l)$ we have:
$$T_\pi\in C_{<G,S_N>}(k,l)\iff T_\pi\in C_G(k,l)$$

It follows that the $D^1$ categories for the groups $<G,S_N>$ and $G$ coincide, and so the easy envelopes $(<G,S_N>)^1$ and $G^1$ coincide as well, as stated.
\end{proof}

In order now to fine-tune all this, by using an arbitrary parameter $p\in\mathbb N$, which can be thought of as being an ``easiness level'', we can proceed as follows:

\index{easiness level}

\begin{definition}
Given a compact group $S_N\subset G\subset U_N$, and an integer $p\in\mathbb N$, we construct the family of linear spaces
$$E^p(k,l)=\left\{\alpha_1T_{\pi_1}+\ldots+\alpha_pT_{\pi_p}\in C(k,l)\Big|\alpha_i\in\mathbb C,\pi_i\in P(k,l)\right\}$$
and we denote by $C^p$ the smallest tensor category containing $E^p=(E^p(k,l))$, and by $S_N\subset G^p\subset U_N$ the compact group corresponding to this category $C^p$.
\end{definition}

As a first observation, at $p=1$ we have $C^1=E^1=span(D^1)$, where $D^1$ is the category of partitions constructed in Proposition 15.25. Thus the group $G^1$ constructed above coincides with the ``easy envelope'' of $G$, from Definition 15.26.

\bigskip

In the general case, $p\in\mathbb N$, the family $E^p=(E^p(k,l))$ constructed above is not necessarily a tensor category, but we can of course consider the tensor category $C^p$ generated by it, as indicated. Finally, in the above definition we have used of course the Tannakian duality results, in order to perform the operation $C^p\to G^p$.

\bigskip

In practice, the construction in Definition 15.28 is often something quite complicated, and it is convenient to use the following observation:

\begin{proposition}
The category $C^p$ constructed above is generated by the spaces
$$E^p(l)=\left\{\alpha_1T_{\pi_1}+\ldots+\alpha_pT_{\pi_p}\in C(l)\Big|\alpha_i\in\mathbb C,\pi_i\in P(l)\right\}$$
where $C(l)=C(0,l),P(l)=P(0,l)$, with $l$ ranging over the colored integers.
\end{proposition}

\begin{proof}
We use the well-known fact, that we know from chapter 13, that given a closed subgroup $G\subset U_N$, we have a Frobenius type isomorphism, as follows:
$$Hom(u^{\otimes k},u^{\otimes l})\simeq Fix(u^{\otimes\bar{k}l})$$

If we apply this to the group $G^p$, we obtain an isomorphism as follows:
$$C(k,l)\simeq C(\bar{k}l)$$

On the other hand, we have as well an isomorphism $P(k,l)\simeq P(\bar{k}l)$, obtained by performing a counterclockwise rotation to the partitions $\pi\in P(k,l)$. According to the above definition of the spaces $E^p(k,l)$, this induces an isomorphism as follows:
$$E^p(k,l)\simeq E^p(\bar{k}l)$$

We deduce from this that for any partitions $\pi_1,\ldots,\pi_p\in C(k,l)$, having rotated versions $\rho_1,\ldots,\rho_p\in C(\bar{k}l)$, and for any scalars $\alpha_1,\ldots,\alpha_p\in\mathbb C$, we have:
$$\alpha_1T_{\pi_1}+\ldots+\alpha_pT_{\pi_p}\in C(k,l)\iff\alpha_1T_{\rho_1}+\ldots+\alpha_pT_{\rho_p}\in C(\bar{k}l)$$

But this gives the conclusion in the statement, and we are done.
\end{proof}

The main properties of the construction $G\to G^p$ can be summarized as follows:

\begin{theorem}
Given a compact group $S_N\subset G\subset U_N$, the compact groups $G^p$ constructed above form a decreasing family, whose intersection is $G$:
$$G=\bigcap_{p\in\mathbb N}G^p$$
Moreover, $G$ is easy when this decreasing limit is stationary, $G=G^1$.
\end{theorem}

\begin{proof}
By definition of $E^p(k,l)$, and by using Proposition 15.29, these linear spaces form an increasing filtration of $C(k,l)$. The same remains true when completing into tensor categories, and so we have an increasing filtration, as follows:
$$C=\bigcup_{p\in\mathbb N}C^p$$

At the compact group level now, we obtain the decreasing intersection in the statement. Finally, the last assertion is clear from Proposition 15.29.
\end{proof}

As a main consequence of the above results, we can now formulate: 

\begin{definition}
We say that a homogeneous compact group 
$$S_N\subset G\subset U_N$$
is easy at order $p$ when $G=G^p$, with $p$ being chosen minimal with this property.
\end{definition}

Observe that the order 1 notion corresponds to the usual easiness. In general, all this is quite abstract, but there are several explicit examples, that can be worked out. For more on all this, you can check my group theory book \cite{ba3}.

\section*{15d. Classification results}

Let us go back now to plain easiness, and discuss some classification results, following the old paper \cite{bsp}, and then the more recent paper of Tarrago-Weber \cite{twe}. In order to cut from the complexity, we must impose an extra axiom, and we will use here:

\index{uniform group}

\begin{theorem}
For an easy group $G=(G_N)$, coming from a category of partitions $D\subset P$, the following conditions are equivalent:
\begin{enumerate}
\item $G_{N-1}=G_N\cap U_{N-1}$, via the embedding $U_{N-1}\subset U_N$ given by $u\to diag(u,1)$.

\item $G_{N-1}=G_N\cap U_{N-1}$, via the $N$ possible diagonal embeddings $U_{N-1}\subset U_N$.

\item $D$ is stable under the operation which consists in removing blocks.
\end{enumerate}
If these conditions are satisfied, we say that $G=(G_N)$ is uniform.
\end{theorem}

\begin{proof}
We use the general easiness theory explained above, as follows:

\medskip

$(1)\iff(2)$ This is something standard, coming from the inclusion $S_N\subset G_N$, which makes everything $S_N$-invariant. The result follows as well from the proof of $(1)\iff(3)$ below, which can be converted into a proof of $(2)\iff(3)$, in the obvious way.

\medskip

$(1)\iff(3)$ Given a subgroup $K\subset U_{N-1}$, with fundamental representation $u$, consider the $N\times N$ matrix $v=diag(u,1)$. Our claim is that for any $\pi\in P(k)$ we have:
$$\xi_\pi\in Fix(v^{\otimes k})\iff\xi_{\pi'}\in Fix(v^{\otimes k'}),\,\forall\pi'\in P(k'),\pi'\subset\pi$$

In order to prove this, we must study the condition on the left. We have:
\begin{eqnarray*}
\xi_\pi\in Fix(v^{\otimes k})
&\iff&(v^{\otimes k}\xi_\pi)_{i_1\ldots i_k}=(\xi_\pi)_{i_1\ldots i_k},\forall i\\
&\iff&\sum_j(v^{\otimes k})_{i_1\ldots i_k,j_1\ldots j_k}(\xi_\pi)_{j_1\ldots j_k}=(\xi_\pi)_{i_1\ldots i_k},\forall i\\
&\iff&\sum_j\delta_\pi(j_1,\ldots,j_k)v_{i_1j_1}\ldots v_{i_kj_k}=\delta_\pi(i_1,\ldots,i_k),\forall i
\end{eqnarray*}

Now let us recall that our representation has the special form $v=diag(u,1)$. We conclude from this that for any index $a\in\{1,\ldots,k\}$, we must have:
$$i_a=N\implies j_a=N$$

With this observation in hand, if we denote by $i',j'$ the multi-indices obtained from $i,j$ obtained by erasing all the above $i_a=j_a=N$ values, and by $k'\leq k$ the common length of these new multi-indices, our condition becomes:
$$\sum_{j'}\delta_\pi(j_1,\ldots,j_k)(v^{\otimes k'})_{i'j'}=\delta_\pi(i_1,\ldots,i_k),\forall i$$

Here the index $j$ is by definition obtained from $j'$ by filling with $N$ values. In order to finish now, we have two cases, depending on $i$, as follows:

\medskip

\underline{Case 1}. Assume that the index set $\{a|i_a=N\}$ corresponds to a certain subpartition $\pi'\subset\pi$. In this case, the $N$ values will not matter, and our formula becomes:
$$\sum_{j'}\delta_\pi(j'_1,\ldots,j'_{k'})(v^{\otimes k'})_{i'j'}=\delta_\pi(i'_1,\ldots,i'_{k'})$$

\underline{Case 2}. Assume now the opposite, namely that the set $\{a|i_a=N\}$ does not correspond to a subpartition $\pi'\subset\pi$. In this case the indices mix, and our formula reads:
$$0=0$$

Thus, we are led to $\xi_{\pi'}\in Fix(v^{\otimes k'})$, for any subpartition $\pi'\subset\pi$, as claimed.

\medskip

Now with this claim in hand, the result follows from Tannakian duality.
\end{proof}

We can now formulate a first classification result, as follows:

\begin{theorem}
The uniform orthogonal easy groups are as follows,
$$
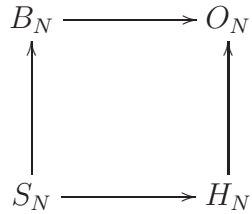
\xymatrix@R=50pt@C=50pt{
B_N\ar[r]&O_N\\
S_N\ar[u]\ar[r]&H_N\ar[u]}$$
and this diagram is an intersection and easy generation diagram.
\end{theorem}

\begin{proof}
We know that the various orthogonal groups in the statement are indeed easy and uniform, the corresponding categories of partitions being as follows:
$$
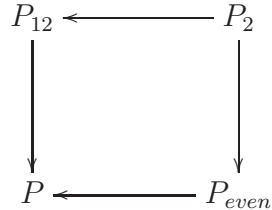
\xymatrix@R=50pt@C50pt{
P_{12}\ar[d]&P_2\ar[d]\ar[l]\\
P&P_{even}\ar[l]}$$

Since this latter diagram is an intersection and generation diagram, we conclude that we have an intersection and easy generation diagram of groups, as stated. Regarding now the classification, consider an arbitrary easy group, as follows:
$$S_N\subset G_N\subset O_N$$

This group must then come from a category of partitions, as follows:
$$P_2\subset D\subset P$$

Now if we assume $G=(G_N)$ to be uniform, this category of partitions $D$ is uniquely determined by the subset $L\subset\mathbb N$ consisting of the sizes of the blocks of the partitions in $D$. Following \cite{bsp}, our claim is that the admissible sets are as follows:

\medskip

\begin{enumerate}
\item $L=\{2\}$, producing $O_N$.

\medskip

\item $L=\{1,2\}$, producing $B_N$.

\medskip

\item $L=\{2,4,6,\ldots\}$, producing $H_N$.

\medskip

\item $L=\{1,2,3,\ldots\}$, producing $S_N$.
\end{enumerate}

\medskip

Indeed, in one sense, this follows from our easiness results for $O_N,B_N,H_N,S_N$. In the other sense now, assume that $L\subset\mathbb N$ is such that the set $P_L$ consisting of partitions whose sizes of the blocks belong to $L$ is a category of partitions. We know from the axioms of the categories of partitions that the semicircle $\cap$ must be in the category, so we have $2\in L$. Our claim is that the following conditions must be satisfied as well:
$$k,l\in L,\,k>l\implies k-l\in L$$
$$k\in L,\,k\geq 2\implies 2k-2\in L$$

Indeed, we will prove that both conditions follow from the axioms of the categories of
partitions. Let us denote by $b_k\in P(0,k)$ the one-block partition, as follows:
$$b_k=\left\{\begin{matrix}\sqcap\hskip-0.7mm \sqcap&\ldots&\sqcap\\
1\hskip2mm 2&\ldots&k\end{matrix} \right\}$$

For $k>l$, we can write $b_{k-l}$ in the following way:
$$b_{k-l}=\left\{\begin{matrix}\sqcap\hskip-0.7mm
\sqcap&\ldots&\ldots&\ldots&\ldots&\sqcap\\ 1\hskip2mm 2&\ldots&l&l+1&\ldots&k\\
\sqcup\hskip-0.7mm \sqcup&\ldots&\sqcup&|&\ldots&|\\ &&&1&\ldots&k-l\end{matrix}\right\}$$

In other words, we have the following formula:
$$b_{k-l}=(b_l^*\otimes |^{\otimes k-l})b_k$$

Since all the terms of this composition are in $P_L$, we have $b_{k-l}\in P_L$, and this proves our first formula. As for the second formula, this can be proved in a similar way, by capping two adjacent $k$-blocks with a $2$-block, in the middle.

\medskip

With the above two formulae in hand, we can conclude in the following way:

\medskip

\underline{Case 1}. Assume $1\in L$. By using the first formula with $l=1$ we get:
$$k\in L\implies k-1\in L$$

This condition shows that we must have $L=\{1,2,\ldots,m\}$, for a certain number $m\in\{1,2,\ldots,\infty\}$. On the other hand, by using the second formula we get:
\begin{eqnarray*}
m\in L
&\implies&2m-2\in L\\
&\implies&2m-2\leq m\\
&\implies&m\in\{1,2,\infty\}
\end{eqnarray*}

The case $m=1$ being excluded by the condition $2\in L$, we reach to one of the two sets producing the groups $S_N,B_N$.

\medskip

\underline{Case 2}. Assume $1\notin L$. By using the first formula with $l=2$ we get:
$$k\in L\implies k-2\in L$$

This condition shows that we must have $L=\{2,4,\ldots,2p\}$, for a certain number $p\in\{1,2,\ldots,\infty\}$. On the other hand, by using the second formula we get:
\begin{eqnarray*}
2p\in L
&\implies&4p-2\in L\\
&\implies&4p-2\leq 2p\\
&\implies&p\in\{1,\infty\}
\end{eqnarray*}

Thus $L$ must be one of the two sets producing $O_N,H_N$, and we are done.
\end{proof}

All the above is very nice, but the continuation of the story is more complicated. When lifting the uniformity assumption, the final classification results become more technical, due to the presence of various copies of $\mathbb Z_2$, that can be added, while keeping the easiness property still true. To be more precise, in the real case, as explained in \cite{bsp}, we have exactly 6 solutions, which are as follows, with the convention $G_N'=G_N\times\mathbb Z_2$:
$$
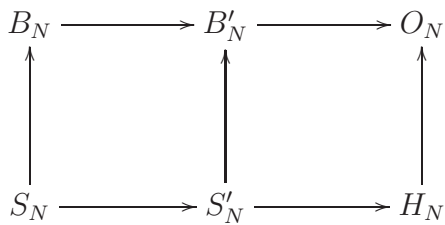
\xymatrix@R=50pt@C=50pt{
B_N\ar[r]&B_N'\ar[r]&O_N\\
S_N\ar[u]\ar[r]&S_N'\ar[u]\ar[r]&H_N\ar[u]}$$

In the unitary case now, the classification is quite similar, but more complicated, as explained in the paper of Tarrago-Weber \cite{twe}. In particular we have:

\begin{theorem}
The uniform easy groups which are purely unitary, in the sense that they appear as complexifications of real easy groups, are as follows,
$$
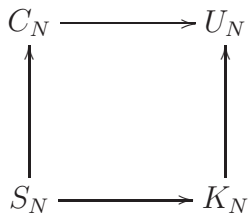
\xymatrix@R=50pt@C=50pt{
C_N\ar[r]&U_N\\
S_N\ar[u]\ar[r]&K_N\ar[u]}$$
and this diagram is an intersection and easy generation diagram.
\end{theorem}

\begin{proof}
We know from the above that the groups in the statement are indeed easy and uniform, the corresponding categories of partitions being as follows:
$$
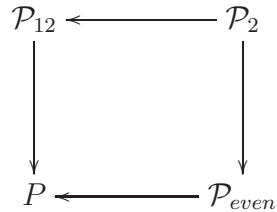
\xymatrix@R=50pt@C50pt{
\mathcal P_{12}\ar[d]&\mathcal P_2\ar[d]\ar[l]\\
P&\mathcal P_{even}\ar[l]}$$

Since this latter diagram is an intersection and generation diagram, we conclude that we have an intersection and easy generation diagram of groups, as stated. As for the uniqueness result, the proof here is similar to the proof from the real case, from Theorem 15.33, by examining the possible sizes of the blocks of the partitions in the category, and doing some direct combinatorics. For details here, we refer to Tarrago-Weber \cite{twe}.
\end{proof}

Finally, let us mention that the easy quantum group formalism can be extended into a ``super-easy'' group formalism, covering as well the symplectic group $Sp_N$. This is something a bit technical, and we refer here to the paper of Collins-\'Sniady \cite{csn}.

\section*{15e. Exercises}

In relation with the notion of easy envelope, we have the following exercise:

\begin{exercise}
Compute the easy envelope of general complex reflection groups
$$H_N^{sd}=\left\{U\in H_N^s\Big|(\square\, U)^d=1\right\}$$
with the symbol $\square$ denoting, as usual, the product of nonzero entries.
\end{exercise}

This is something which does not look very difficult, and you have the choice here, either by using combinatorics, or the universality property of the easy envelope.

\begin{exercise}
Work out the super-easiness property of the symplectic group 
$$Sp_N\subset U_N$$
defined for $N\in\mathbb N$ even, then try as well the groups $SU_2$ and $SO_3$.
\end{exercise}

This is actually a quite difficult exercise. Many things to be done here.

\begin{exercise}
Prove that when lifting the uniformity assumption, the groups
$$\xymatrix@R=50pt@C=50pt{
B_N\ar[r]&B_N'\ar[r]&O_N\\
S_N\ar[u]\ar[r]&S_N'\ar[u]\ar[r]&H_N\ar[u]}$$
with the convention $G_N'=G_N\times\mathbb Z_2$, are the only easy real groups.
\end{exercise}

This is something quite standard, briefly discussed in the above.

\begin{exercise}
Prove that the uniform, purely unitary easy groups are
$$\xymatrix@R=50pt@C=50pt{
C_N\ar[r]&U_N\\
S_N\ar[u]\ar[r]&K_N\ar[u]}$$
with a suitable definition for the notion of pure unitarity.
\end{exercise}

As before, this is something quite standard, briefly discussed in the above, the idea being that of adapting the proof of the classification from the real uniform case.

\chapter{Weingarten calculus}

\section*{16a. Weingarten formula}

Time now to put everything together. We will discuss here applications of the theory developed above, to the computation of the laws of characters, and truncated characters, as to solve the various questions left open in Part III, for the continuous groups. Generally speaking, all these questions require a good knowledge of the integration over $G$, and more precisely, of the various polynomial integrals over $G$, defined as follows:

\index{polynomial integrals}

\begin{definition}
Given a closed subgroup $G\subset U_N$, the quantities
$$I_k=\int_Gg_{i_1j_1}^{e_1}\ldots g_{i_kj_k}^{e_k}\,dg$$
depending on a colored integer $k=e_1\ldots e_k$, are called polynomial integrals over $G$.
\end{definition}

As a first observation, the knowledge of these integrals is the same as the full knowledge of the integration functional over $G$. Indeed, since the coordinate functions $g\to g_{ij}$ separate the points of $G$, we can apply the Stone-Weierstrass theorem, and we obtain:
$$C(G)=<g_{ij}>$$

Thus, by linearity, the computation of any functional $f:C(G)\to\mathbb C$, and in particular of the integration functional, reduces to the computation of this functional on the polynomials of the coordinate functions $g\to g_{ij}$ and their conjugates $g\to\bar{g}_{ij}$.

\bigskip

The point now is that, by using Peter-Weyl, everything reduces to linear algebra, and more specifically to a matrix inversion question, due to the following result:

\index{Weingarten formula}

\begin{theorem}
The Haar integration over a closed subgroup $G\subset U_N$ is given on the dense subalgebra of smooth functions by the Weingarten type formula
$$\int_Gg_{i_1j_1}^{e_1}\ldots g_{i_kj_k}^{e_k}\,dg=\sum_{\pi,\sigma\in D_k}\delta_\pi(i)\delta_\sigma(j)W_k(\pi,\sigma)$$
valid for any colored integer $k=e_1\ldots e_k$ and any multi-indices $i,j$, where $D_k$ is a linear basis of $Fix(u^{\otimes k})$, the associated generalized Kronecker symbols are given by
$$\delta_\pi(i)=<\pi,e_{i_1}\otimes\ldots\otimes e_{i_k}>$$
and $W_k=G_k^{-1}$ is the inverse of the Gram matrix, $G_k(\pi,\sigma)=<\pi,\sigma>$.
\end{theorem}

\begin{proof}
This is something that we know from chapter 13, the idea being that the above integrals form altogether the orthogonal projection $P^k$ onto the following space:
$$Fix(u^{\otimes k})=span(D_k)$$

Consider now the following linear map, with $D_k=\{\xi_k\}$ being as in the statement:
$$E(x)=\sum_{\pi\in D_k}<x,\xi_\pi>\xi_\pi$$

By a standard linear algebra computation, it follows that we have $P=WE$, where $W$ is the inverse of the restriction of $E$ to the following space:
$$K=span\left(T_\pi\Big|\pi\in D_k\right)$$

But this restriction is the linear map given by the matrix $G_k$, and so $W$ is the linear map given by the inverse matrix $W_k=G_k^{-1}$, and this gives the result.
\end{proof}

In the easy case now, we have the following more precise result:

\index{Weingarten formula}
\index{Gram matrix}
\index{Weingarten matrix}

\begin{theorem}
For an easy group $G\subset U_N$, coming from a category of partitions $D=(D(k,l))$, we have the Weingarten integration formula
$$\int_Gu_{i_1j_1}^{e_1}\ldots u_{i_kj_k}^{e_k}=\sum_{\pi,\sigma\in D(k)}\delta_\pi(i)\delta_\sigma(j)W_{kN}(\pi,\sigma)$$
for any multi-indices $i,j$ and any exponent $k=e_1\ldots e_k$, where $D(k)=D(\emptyset,k)$, the $\delta$ numbers are the usual Kronecker type symbols, and $W_{kN}=G_{kN}^{-1}$, with 
$$G_{kN}(\pi,\sigma)=N^{|\pi\vee\sigma|}$$
where $|.|$ is the number of blocks. 
\end{theorem}

\begin{proof}
We use the abstract Weingarten formula, from Theorem 16.2. According to our easiness conventions, the Kronecker symbols are given by:
\begin{eqnarray*}
\delta_{\xi_\pi}(i)
&=&<\xi_\pi,e_{i_1}\otimes\ldots\otimes e_{i_k}>\\
&=&\left<\sum_j\delta_\pi(j_1,\ldots,j_k)e_{j_1}\otimes\ldots\otimes e_{j_k},e_{i_1}\otimes\ldots\otimes e_{i_k}\right>\\
&=&\delta_\pi(i_1,\ldots,i_k)
\end{eqnarray*}

The Gram matrix being as well the correct one, we obtain the result.
\end{proof}

Generally speaking, the above result is something quite powerful, because the main computation there, that of the inverse matrix $W_{kN}=G_{kN}^{-1}$, can be run on an ordinary laptop, after implementing the formula of the Gram matrix, namely $G_{kN}(\pi,\sigma)=N^{|\pi\vee\sigma|}$, which is something quite easy to do. Thus, you can prove theorems about integrals over easy groups just by smoking cigars, and letting your computer do the work.

\bigskip

Let us also mention that there is a long story behind the above results. Generally speaking, such things have been known since ever, and more precisely, since the old work of Weyl \cite{wy2} and Brauer \cite{bra}. However, in what regards the applications of the Weingarten formula, to various questions in mathematics or physics, and the interest in this formula in general, things here have evolved over the time with several ups and lows:

\bigskip

(1) In modern times, this formula has been quite popular among physicists since the 1978 paper of Weingarten \cite{wei}, who was motivated by physics, and among mathematicians, since the 2003 paper of Collins \cite{col}, who was motivated by physics too.

\bigskip

(2) A key step was the 2006 paper of Collins-\'Sniady \cite{csn}, with this formula clearly explained, for the unitary, orthogonal, and symplectic groups as well, and made ready to use, for everyone willing to do so, be them mathematicians or physicists.

\bigskip

(3) This technology has always been something rival to the Lie algebra theory, and a further increase in popularity came from the series of papers \cite{bb+}, \cite{bbc}, \cite{bco}, \cite{bsp}, extending this formula to the quantum group setting, where no Lie theory is available.

\bigskip

(4) Finally, at the level of the applications, there are many of them, but probably the most popular ones, in recent times, came from quantum information theory work of Collins-Nechita, \cite{cne} and subsequent papers, heavily relying on this formula.

\bigskip 

Back to work now, as a first illustration for Theorem 16.3, let us discuss the computation of the Weingarten function for $S_N$. For this purpose, we can use the following result, which actually shows that the Weingarten formula is not really needed for $S_N$:

\begin{theorem}
Consider the symmetric group $S_N\subset O_N$, with coordinates given by: 
$$g_{ij}=\chi\left(\sigma\in S_N\Big|\sigma(j)=i\right)$$
The products of these coordinates span then the algebra of functions $C(S_N)$, and the arbitrary integrals over $S_N$ are given, modulo linearity, by the formula
$$\int_{S_N}g_{i_1j_1}\ldots g_{i_kj_k}=\begin{cases}
\frac{(N-|\ker i|)!}{N!}&{\rm if}\ \ker i=\ker j\\
0&{\rm otherwise}
\end{cases}$$
where $\ker i$ denotes as usual the partition of $\{1,\ldots,k\}$ whose blocks collect the equal indices of $i$, and where $|.|$ denotes the number of blocks.
\end{theorem}

\begin{proof}
This is something that we know from chapter 11, the idea being that, according to the formula of the coordinates $g_{ij}$, the polynomial integrals are given by:
$$\int_{S_N}g_{i_1j_1}\ldots g_{i_kj_k}=\frac{1}{N!}\#\left\{\sigma\in S_N\Big|\sigma(j_1)=i_1,\ldots,\sigma(j_k)=i_k\right\}$$

Now observe that the existence of $\sigma\in S_N$ as above requires:
$$i_m=i_n\iff j_m=j_n$$

Thus, the above integral vanishes when the following happens:
$$\ker i\neq\ker j$$

Regarding now the case $\ker i=\ker j$, if we denote by $b\in\{1,\ldots,k\}$ the number of blocks of this partition $\ker i=\ker j$, we have $N-b$ points to be sent bijectively to $N-b$ points, and so $(N-b)!$ solutions, and the integral is $\frac{(N-b)!}{N!}$, as claimed.
\end{proof}

The above result shows that the integration over $S_N$ is something quite trivial, and no surprise here, and so that the computation of the Weingarten function should be something quite trivial too. In practice now, in order to compute the Weingarten function for $S_N$, by using the above result, we will need some combinatorics, and more specifically the M\"obius inversion formula. Let us begin with some standard definitions, as follows:

\begin{definition}
Let $P(k)$ be the set of partitions of $\{1,\ldots,k\}$, and let $\pi,\sigma\in P(k)$.
\begin{enumerate}
\item We write $\pi\leq\sigma$ if each block of $\pi$ is contained in a block of $\sigma$.

\item We let $\pi\vee\sigma\in P(k)$ be the partition obtained by superposing $\pi,\sigma$.
\end{enumerate}
\end{definition}

As an illustration here, at $k=2$ we have $P(2)=\{||,\sqcap\}$, and we have:
$$||\leq\sqcap$$

Also, at $k=3$ we have $P(3)=\{|||,\sqcap|,\sqcap\hskip-3.2mm{\ }_|\,,|\sqcap,\sqcap\hskip-0.7mm\sqcap\}$, and the order relation is as follows:
$$|||\leq\sqcap|,\sqcap\hskip-3.2mm{\ }_|\,,|\sqcap\leq\sqcap\hskip-0.7mm\sqcap$$

Observe also that we have the following inequalities:
$$\pi,\sigma\leq\pi\vee\sigma$$

In fact, the partition $\pi\vee\sigma$ is by construction the smallest possible one with this property. Due to this fact, this partition $\pi\vee\sigma$ is called supremum of $\pi,\sigma$. 

\bigskip

We can now introduce the M\"obius function, as follows:

\index{M\"obius function}

\begin{definition}
The M\"obius function of any lattice, and so of $P$, is given by
$$\mu(\pi,\sigma)=\begin{cases}
1&{\rm if}\ \pi=\sigma\\
-\sum_{\pi\leq\tau<\sigma}\mu(\pi,\tau)&{\rm if}\ \pi<\sigma\\
0&{\rm if}\ \pi\not\leq\sigma
\end{cases}$$
with the construction being performed by recurrence.
\end{definition}

As an illustration here, let us go back to the set of 2-point partitions, $P(2)=\{||,\sqcap\}$. We have here, by definition of the M\"obius function:
$$\mu(||,||)=\mu(\sqcap,\sqcap)=1$$

Also, we know that we have $||<\sqcap$, with no intermediate partition in between, and so the above recurrence procedure gives the following formulae:
$$\mu(||,\sqcap)=-\mu(||,||)=-1$$

Finally, we have $\sqcap\not\leq||$, and so $\mu(\sqcap,||)=0$. Thus, as a conclusion, the M\"obius matrix $M_{\pi\sigma}=\mu(\pi,\sigma)$ of the lattice $P(2)=\{||,\sqcap\}$ is as follows:
$$M=\begin{pmatrix}1&-1\\ 0&1\end{pmatrix}$$

The interest in the M\"obius function comes from the M\"obius inversion formula:
$$f(\sigma)=\sum_{\pi\leq\sigma}g(\pi)\implies g(\sigma)=\sum_{\pi\leq\sigma}\mu(\pi,\sigma)f(\pi)$$

In linear algebra terms, the statement and proof of this formula are as follows:

\begin{theorem}
The inverse of the adjacency matrix of $P$, given by
$$A_{\pi\sigma}=\begin{cases}
1&{\rm if}\ \pi\leq\sigma\\
0&{\rm if}\ \pi\not\leq\sigma
\end{cases}$$
is the M\"obius matrix of $P$, given by $M_{\pi\sigma}=\mu(\pi,\sigma)$.
\end{theorem}

\begin{proof}
This is well-known, coming for instance from the fact that $A$ is upper triangular. Indeed, when inverting, we are led into the recurrence from Definition 16.6.
\end{proof}

As a first illustration, for $P(2)$ the formula $M=A^{-1}$ appears as follows:
$$\begin{pmatrix}1&-1\\ 0&1\end{pmatrix}=
\begin{pmatrix}1&1\\ 0&1\end{pmatrix}^{-1}$$

Also, for $P(3)=\{|||,\sqcap|,\sqcap\hskip-3.2mm{\ }_|\,,|\sqcap,\sqcap\hskip-0.7mm\sqcap\}$ the formula $M=A^{-1}$ reads:
$$\begin{pmatrix}
1&-1&-1&-1&2\\
0&1&0&0&-1\\
0&0&1&0&-1\\
0&0&0&1&-1\\
0&0&0&0&1
\end{pmatrix}=
\begin{pmatrix}
1&1&1&1&1\\
0&1&0&0&1\\
0&0&1&0&1\\
0&0&0&1&1\\
0&0&0&0&1
\end{pmatrix}^{-1}$$

With the above results in hand, we can now compute the Weingarten function of $S_N$, and also find a precise estimate for it, as follows:

\begin{theorem}
For $S_N$ the Weingarten function is given by
$$W_{kN}(\pi,\sigma)=\sum_{\tau\leq\pi\wedge\sigma}\mu(\tau,\pi)\mu(\tau,\sigma)\frac{(N-|\tau|)!}{N!}$$
and satisfies the folowing estimate,
$$W_{kN}(\pi,\sigma)=N^{-|\pi\wedge\sigma|}(
\mu(\pi\wedge\sigma,\pi)\mu(\pi\wedge\sigma,\sigma)+O(N^{-1}))$$
with $\mu$ being the M\"obius function of $P(k)$.
\end{theorem}

\begin{proof}
The first assertion follows from the Weingarten formula, namely:
$$\int_{S_N}u_{i_1j_1}\ldots u_{i_kj_k}=\sum_{\pi,\sigma\in P(k)}\delta_\pi(i)\delta_\sigma(j)W_{kN}(\pi,\sigma)$$

Indeed, in this formula the integrals on the left are known, from the explicit integration formula over $S_N$ that we established above, namely:
$$\int_{S_N}g_{i_1j_1}\ldots g_{i_kj_k}=\begin{cases}
\frac{(N-|\ker i|)!}{N!}&{\rm if}\ \ker i=\ker j\\
0&{\rm otherwise}
\end{cases}$$

But this allows the computation of the right term, via the M\"obius inversion formula, explained above. As for the second assertion, this follows from the first one. See \cite{bcu}.
\end{proof}

As an illustration, let us record the formulae at $k=2,3$. At $k=2$, with indices $||,\sqcap$, and with the convention that $\approx$ means componentwise dominant term, we have:
$$W_{2N}\approx\begin{pmatrix}N^{-2}&-N^{-2}\\-N^{-2}&N^{-1}\end{pmatrix}$$

At $k=3$ now, with indices $|||,|\sqcap,\sqcap|,\sqcap\hskip-3.2mm{\ }_|,\sqcap\hskip-0.8mm\sqcap$, and same meaning for $\approx$, we have:
$$W_{3N}\approx\begin{pmatrix}
N^{-3}&-N^{-3}&-N^{-3}&-N^{-3}&2N^{-3}\\
-N^{-3}&N^{-2}&N^{-3}&N^{-3}&-N^{-2}\\
-N^{-3}&N^{-3}&N^{-2}&N^{-3}&-N^{-2}\\
-N^{-3}&N^{-3}&N^{-3}&N^{-2}&-N^{-2}\\
2N^{-3}&-N^{-2}&-N^{-2}&-N^{-2}&N^{-1}
\end{pmatrix}$$

We will be back to all this later, with results about the orthogonal group $O_N$ and about some other easy groups as well, where the Weingarten function is in general not explicitly computable, but where some useful estimates are still possible.

\section*{16b. Laws of characters}

As a first concrete application of the above, let us discuss now the computation of the asymptotic laws of truncated characters. We have the following result, to start with:

\index{main character}
\index{moments}

\begin{theorem}
Assuming that $G\subset U_N$ is easy, coming from a category of partitions 
$$D=(D(k,l))$$
the moments of the main character are given by the formula
$$\int_G\chi^k=\dim\Big(span\left(\xi_\pi\big|\pi\in D(k)\right)\Big)$$
where $D(k)=D(\emptyset,k)$, and where for $\pi\in D(k)$ we use the notation $\xi_\pi=T_\pi$.
\end{theorem}

\begin{proof}
We recall that for an easy group $G\subset U_N$, coming from a category of partitions $D=(D(k,l))$, we have by definition equalities as follows:
$$Hom(u^{\otimes k},u^{\otimes l})=span\left(T_\pi\Big|\pi\in D(k,l)\right)$$

By interchanging $k\leftrightarrow l$ in this formula, and then setting $l=\emptyset$, we obtain:
$$Fix(u^{\otimes k})=span\left(\xi_\pi\Big|\pi\in D(k)\right)$$

Now since by the Peter-Weyl theory integrating a character amounts in counting the fixed points, we are led to the conclusion in the statement.
\end{proof}

In order to investigate the linear independence questions for the vectors $\xi_\pi$, we will use the Gram matrix of these vectors. We have the following result, to start with:

\begin{proposition}
The Gram matrix $G_{kN}(\pi,\sigma)=<\xi_\pi,\xi_\sigma>$ is given by
$$G_{kN}(\pi,\sigma)=N^{|\pi\vee\sigma|}$$
where $|.|$ is the number of blocks.
\end{proposition}

\begin{proof}
According to the formula of the vectors $\xi_\pi$, we have:
\begin{eqnarray*}
<\xi_\pi,\xi_\sigma>
&=&\sum_{i_1\ldots i_k}\delta_\pi(i_1,\ldots,i_k)\delta_\sigma(i_1,\ldots,i_k)\\
&=&\sum_{i_1\ldots i_k}\delta_{\pi\vee\sigma}(i_1,\ldots,i_k)\\
&=&N^{|\pi\vee\sigma|}
\end{eqnarray*}

Thus, we have obtained the formula in the statement.
\end{proof}

Next in line, we have the following key result:

\begin{proposition}
The Gram matrix is given by $G_{kN}=AL$, where
$$L(\pi,\sigma)=
\begin{cases}
N(N-1)\ldots(N-|\pi|+1)&{\rm if}\ \sigma\leq\pi\\
0&{\rm otherwise}
\end{cases}$$
and where $A=M^{-1}$ is the adjacency matrix of $P(k)$.
\end{proposition}

\begin{proof}
We have indeed the following computation:
\begin{eqnarray*}
N^{|\pi\vee\sigma|}
&=&\#\left\{i_1,\ldots,i_k\in\{1,\ldots,N\}\Big|\ker i\geq\pi\vee\sigma\right\}\\
&=&\sum_{\tau\geq\pi\vee\sigma}\#\left\{i_1,\ldots,i_k\in\{1,\ldots,N\}\Big|\ker i=\tau\right\}\\
&=&\sum_{\tau\geq\pi\vee\sigma}N(N-1)\ldots(N-|\tau|+1)
\end{eqnarray*}

According to Proposition 16.10 and to the definition of $A,L$, this formula reads:
\begin{eqnarray*}
(G_{kN})_{\pi\sigma}
&=&\sum_{\tau\geq\pi}L_{\tau\sigma}\\
&=&\sum_\tau A_{\pi\tau}L_{\tau\sigma}\\
&=&(AL)_{\pi\sigma}
\end{eqnarray*}

Thus, we obtain in this way the formula in the statement.
\end{proof}

As an illustration for the above result, at $k=2$ we have $P(2)=\{||,\sqcap\}$, and the above formula $G_{kN}=AL$ appears as follows:
$$\begin{pmatrix}N^2&N\\ N&N\end{pmatrix}
=\begin{pmatrix}1&1\\ 0&1\end{pmatrix}
\begin{pmatrix}N^2-N&0\\N&N\end{pmatrix}$$

At $k=3$ now, we have $P(3)=\{|||,\sqcap|,\sqcap\hskip-3.2mm{\ }_|\,,|\sqcap,\sqcap\hskip-0.7mm\sqcap\}$, and the Gram matrix is:
$$G_3=\begin{pmatrix}
N^3&N^2&N^2&N^2&N\\
N^2&N^2&N&N&N\\
N^2&N&N^2&N&N\\
N^2&N&N&N^2&N\\
N&N&N&N&N
\end{pmatrix}$$

Regarding $L_3$, this can be computed by writing down the matrix $E_3(\pi,\sigma)=\delta_{\sigma\leq\pi}|\pi|$, and then replacing each entry by the corresponding polynomial in $N$. We reach to the conclusion that the product $A_3L_3$ is as follows, producing the above matrix $G_3$:
$$A_3L_3=\begin{pmatrix}
1&1&1&1&1\\
0&1&0&0&1\\
0&0&1&0&1\\
0&0&0&1&1\\
0&0&0&0&1
\end{pmatrix}
\begin{pmatrix}
N^3-3N^2+2N&0&0&0&0\\
N^2-N&N^2-N&0&0&0\\
N^2-N&0&N^2-N&0&0\\
N^2-N&0&0&N^2-N&0\\
N&N&N&N&N
\end{pmatrix}$$

In general, the formula $G_k=A_kL_k$ appears a bit in the same way, with $A_k$ being binary and upper triangular, and with $L_k$ depending on $N$, and being lower triangular.

\bigskip

With the above result in hand, we can now investigate the linear independence properties of the vectors $\xi_\pi$. We have here the following result of Lindst\"om \cite{lin}:

\index{linear independence}

\begin{theorem}
The determinant of the Gram matrix $G_{kN}$ is given by
$$\det(G_{kN})=\prod_{\pi\in P(k)}\frac{N!}{(N-|\pi|)!}$$
and in particular, for $N\geq k$, the vectors $\{\xi_\pi|\pi\in P(k)\}$ are linearly independent.
\end{theorem}

\begin{proof}
According to the formula in Proposition 16.11, we have:
$$\det(G_{kN})=\det(A)\det(L)$$

Now if we order $P(k)$ as above, with respect to the number of blocks, and then lexicographically, we see that the matrix $A$ is upper triangular, and that $L$ is lower triangular. Thus $\det(A)$ can be computed simply by making the product on the diagonal, and we obtain $1$. As for $\det(L)$, this can computed as well by making the product on the diagonal, and we obtain the number in the statement, with the technical remark that in the case $N<k$ the convention is that we obtain a vanishing determinant. 
\end{proof}

Now back to the laws of characters, we can formulate:

\begin{theorem}
For an easy group $G=(G_N)$, coming from a category of partitions $D=(D(k,l))$, the asymptotic moments of the main character are given by
$$\lim_{N\to\infty}\int_{G_N}\chi^k=\# D(k)$$
where $D(k)=D(\emptyset,k)$, with the limiting sequence on the left consisting of certain integers, and being stationary at least starting from the $k$-th term.
\end{theorem}

\begin{proof}
This follows indeed from the general formula from Theorem 16.9, by using the linear independence result from Theorem 16.12.
\end{proof}

Our next purpose will be that of understanding what happens for the basic classes of easy groups. We have here the following result, to start with:

\index{asymptotic character}
\index{normal law}

\begin{theorem}
In the $N\to\infty$ limit, the law of the main character 
$$\chi_u=\sum_{i=1}^Nu_{ii}$$ 
for the orthogonal and unitary groups is as follows:
\begin{enumerate}
\item For $O_N$ we obtain a real Gaussian law $g_1$. 

\item For $U_N$ we obtain a complex Gaussian law $G_1$.
\end{enumerate}
\end{theorem}

\begin{proof}
These results follow indeed from the general formula in Theorem 16.13, by using the knowledge of the associated categories of partitions, as follows:

\medskip

(1) For $O_N$ the associated category of partitions is $P_2$, so the asymptotic moments of the main character are as follows, with the convention $k!!=0$ when $k$ is odd:
$$M_k
=\# P_2(k)
=k!!$$

Thus, we obtain the real Gaussian law, as stated.

\medskip

(2) For $U_N$ the associated category of partitions is $\mathcal P_2$, so the asymptotic moments of the main character, with respect to the colored integers, are as follows:
$$M_k=\#\mathcal P_2(k)$$

Thus, we obtain the complex Gaussian law, as stated.
\end{proof}

More generally now, we have the following result:

\begin{theorem}
With $N\to\infty$, the laws of main character is as follows:
\begin{enumerate}
\item For $O_N$ we obtain the Gaussian law $g_1$.

\item For $U_N$ we obtain the complex Gaussian law $G_1$.

\item For $S_N$ we obtain the Poisson law $p_1$.

\item For $H_N$ we obtain the Bessel law $b_1$.

\item For $H_N^s$ we obtain the generalized Bessel law $b_1^s$.

\item For $K_N$ we obtain the complex Bessel law $B_1$.
\end{enumerate}
Also, for $B_N,C_N$ and for $Sp_N$ we obtain modified Gaussian laws.
\end{theorem}

\begin{proof}
We already know the results for $O_N$ and for $U_N$, from Theorem 16.14. In general, the proof is similar, by counting the partitions in the associated category of partitions, and then doing some calculus, based on the various moment results for the laws in the statement, coming from the general theory developed in the above. All this is of course a bit technical, and for details we refer to \cite{bb+}, \cite{csn} and related papers.
\end{proof}

\section*{16c. Truncated characters}

In order to fully solve the various questions left open in Part III, we still have to discuss now the more advanced question of computing the laws of truncated characters. First, we have the following formula, in the general easy group setting:

\begin{proposition}
The moments of truncated characters are given by the formula
$$\int_G(g_{11}+\ldots +g_{ss})^k=Tr(W_{kN}G_{ks})$$
where $G_{kN}$ and $W_{kN}=G_{kN}^{-1}$ are the associated Gram and Weingarten matrices.
\end{proposition}

\begin{proof}
We have indeed the following computation:
\begin{eqnarray*}
\int_G(g_{11}+\ldots +g_{ss})^k
&=&\sum_{i_1=1}^{s}\ldots\sum_{i_k=1}^s\int_Gg_{i_1i_1}\ldots g_{i_ki_k}\\
&=&\sum_{\pi,\sigma\in D(k)}W_{kN}(\pi,\sigma)\sum_{i_1=1}^{s}\ldots\sum_{i_k=1}^s\delta_\pi(i)\delta_\sigma(i)\\
&=&\sum_{\pi,\sigma\in D(k)}W_{kN}(\pi,\sigma)G_{ks}(\sigma,\pi)\\
&=&Tr(W_{kN}G_{ks})
\end{eqnarray*}

Thus, we have obtained the formula in the statement.
\end{proof}

In order to process now the above formula, and reach to concrete results, we can impose the uniformity condition from chapter 15, originally used there for some technical classification purposes. Let us recall indeed from there that we have:

\index{uniform group}

\begin{definition}
An easy group $G=(G_N)$, coming from a category of partitions $D\subset P$, is called uniform if it satisfies the following equivalent conditions:
\begin{enumerate}
\item $G_{N-1}=G_N\cap U_{N-1}$, via the embedding $U_{N-1}\subset U_N$ given by $u\to diag(u,1)$.

\item $G_{N-1}=G_N\cap U_{N-1}$, via the $N$ possible diagonal embeddings $U_{N-1}\subset U_N$.

\item $D$ is stable under the operation which consists in removing blocks.
\end{enumerate}
\end{definition}

Here the equivalence between the above three conditions is something standard, obtained by doing some combinatorics, and this was discussed in chapter 15. We refer as well to chapter 15 for examples and counterexamples of such groups, the idea here being that the most familiar easy groups $G=(G_N)$ that we know are indeed uniform.

\bigskip

In what follows we will be mostly interested in the condition (3) above, which makes the link with our computations for truncated characters, and simplifies them. To be more precise, by imposing the uniformity condition we obtain:

\begin{theorem}
For a uniform easy group $G=(G_N)$, we have the formula
$$\lim_{N\to\infty}\int_{G_N}\chi_t^k=\sum_{\pi\in D(k)}t^{|\pi|}$$
with $D\subset P$ being the associated category of partitions.
\end{theorem}

\begin{proof}
We use the general moment formula from Proposition 16.16, namely:
$$\int_G(g_{11}+\ldots +g_{ss})^k=Tr(W_{kN}G_{ks})$$

By setting $s=[tN]$, with $t>0$ being a given parameter, this formula becomes:
$$\int_{G_N}\chi_t^k=Tr(W_{kN}G_{k[tN]})$$

The point now is that in the uniform case the Gram and Weingarten matrices are asymptotically diagonal, and this leads to the formula in the statement. See \cite{bco}.
\end{proof}

We can now improve our character results, as follows:

\begin{theorem}
With $N\to\infty$, the laws of truncated characters are as follows:
\begin{enumerate}
\item For $O_N$ we obtain the Gaussian law $g_t$.

\item For $U_N$ we obtain the complex Gaussian law $G_t$.

\item For $S_N$ we obtain the Poisson law $p_t$.

\item For $H_N$ we obtain the Bessel law $b_t$.

\item For $H_N^s$ we obtain the generalized Bessel law $b_t^s$.

\item For $K_N$ we obtain the complex Bessel law $B_t$.
\end{enumerate}
Also, for $B_N,C_N$ and for $Sp_N$ we obtain modified normal laws.
\end{theorem}

\begin{proof}
We use the formula that we found in Theorem 16.18, namely:
$$\lim_{N\to\infty}\int_{G_N}\chi_t^k=\sum_{\pi\in D(k)}t^{|\pi|}$$

By doing now some combinatorics, for instance in relation with the cumulants, this gives the results. We refer here to \cite{bco} and various related papers.
\end{proof}

All the above is quite interesting in relation with questions from theoretical probability. Let us recall indeed that we have 4 main limiting results in probability, namely real and complex, and discrete and continuous, which are as follows:
$$
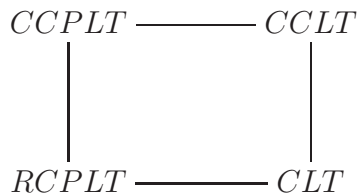
\xymatrix@R=45pt@C=35pt{
CPLT_\mathbb C\ar@{-}[r]\ar@{-}[d]&CCLT\ar@{-}[d]\\
CPLT_\mathbb R\ar@{-}[r]&CLT
}$$

We also know from chapter 12 that the limiting laws in these main limiting theorems are the real and complex Gaussian and Bessel laws, which are as follows:
$$\xymatrix@R=45pt@C=45pt{
B_t\ar@{-}[r]\ar@{-}[d]&G_t\ar@{-}[d]\\
b_t\ar@{-}[r]&g_t
}$$

Moreover, we have also seen in the above that at the level of the moments, these come from certain collections of partitions, as follows:
$$\xymatrix@R=45pt@C=45pt{
\mathcal P_{even}\ar@{-}[r]\ar@{-}[d]&\mathcal P_2\ar@{-}[d]\\
P_{even}\ar@{-}[r]&P_2
}$$

The point now is that, according to our general easiness philosophy, and also to Theorem 16.19, there are some Lie groups behind all this probability theory, namely the basic real and complex rotation and reflection groups, which as follows:
$$\xymatrix@R=45pt@C=45pt{
K_N\ar@{-}[r]\ar@{-}[d]&U_N\ar@{-}[d]\\
H_N\ar@{-}[r]&O_N
}$$

To be more precise, these Lie groups correspond via easiness to the categories of partitions given above, and the corresponding measures can be recaptured as well, as being the asymptotic laws of the corresponding truncated characters, as explained in Theorem 16.19. As for the main probabilistic limiting results themselves, these are of course related too to these Lie groups, but this is something a bit more technical.

\bigskip

All this is very nice. With all this in hand, we are now at a rather advanced level in theoretical probability, and with this knowledge, you can virtually read any article or book in theoretical probability, that you might want to. With our recommendations here being the article of Diaconis-Shahshahani \cite{dsh}, and other texts by Diaconis, which are all quite magic, and no wonder here, because Diaconis used to be a professional magician before doing mathematics, then the classical and lovely random matrix book by Mehta \cite{meh}, and then some fancy theoretical physics from Collins-Nechita \cite{cne}.

\section*{16d. Standard estimates}

We have seen in the above that the Weingarten calculus is something very efficient in dealing with various probability questions over the easy groups $G\subset U_N$. We discuss now, as a continuation of this, a number of more advanced aspects of the Weingarten function combinatorics. We will be mostly interested in the case $G=O_N$. To be more precise, we will be interested in the computation of the polynomial integrals over $O_N$. These polynomial integrals are best introduced in a ``rectangular way'', as follows:

\begin{definition}
Associated to any matrix $a\in M_{p\times q}(\mathbb N)$ is the integral
$$I(a)=\int_{O_N}\prod_{i=1}^p\prod_{j=1}^qu_{ij}^{a_{ij}}\,du$$
with respect to the Haar measure of $O_N$, where $N\geq p,q$.
\end{definition}

As a first observation, we can of course complete our matrix with 0 values, as to always deal with square matrices, $a\in M_N(\mathbb N)$. However, the parameters $p,q$ are very useful, because they measure the ``complexity'' of the problem, so we will keep them.

\bigskip

In order to get familiar with the above integrals, let us do some computations. With the convention $x!!=(x-1)(x-3)(x-5)\ldots\,$, with product ending at $1$ or $2$, we have:

\begin{theorem}
At $p=1$ we have the formula
$$I\begin{pmatrix}a_1&\ldots&a_q\end{pmatrix}=\varepsilon\cdot\frac{(N-1)!!a_1!!\ldots a_q!!}{(N+\Sigma a_i-1)!!}$$
where $\varepsilon=1$ if all $a_i$ are even, and $\varepsilon=0$ otherwise.
\end{theorem}

\begin{proof}
This follows from the fact that the first slice of $O_N$ is isomorphic to the real sphere $S^{N-1}_\mathbb R$. Indeed, this gives the following formula:
$$I\begin{pmatrix}a_1&\ldots&a_q\end{pmatrix}=\int_{S^{N-1}_\mathbb R}x_1^{a_1}\ldots x_q^{a_q}\,dx$$ 

But this latter integral can be computed by using polar coordinates, via the various formulae from chapters 5-6, and we obtain the formula in the statement.
\end{proof}

Another instructive computation, as well of trigonometric nature, is the one at $N=2$. We have here the following result, which completely solves the problem in this case:

\begin{theorem}
At $N=2$ we have the formula
$$I\begin{pmatrix}a&b\\ c&d\end{pmatrix}=\varepsilon\cdot\frac{(a+d)!!(b+c)!!}{(a+b+c+d+1)!!}$$
where $\varepsilon=1$ if $a,b,c,d$ are even, $\varepsilon=-1$ is $a,b,c,d$ are odd, and $\varepsilon=0$ otherwise.
\end{theorem}

\begin{proof}
When computing the integral over $O_2$, we can restrict the integration to $SO_2=\mathbb T$, then further restrict the integration to the first quadrant. We get: 
$$I\begin{pmatrix}a&b\\ c&d\end{pmatrix}
=\varepsilon\cdot\frac{2}{\pi}\int_0^{\pi/2}(\cos t)^{a+d}(\sin t)^{b+c}\,dt$$

By using now the formulae for trigonometric integrals from chapters 5-6, this gives the formula in the statement, with our previous convention for the double factorials.
\end{proof}

The above computations might tend to suggest that $I(a)$ always decomposes as a product of factorials. However, this is far from being true, but in the $2\times 2$ case it is known that $I(a)$ decomposes as a quite reasonable sum of products of factorials. This is something quite technical, from \cite{bcs}, and we will be back to this, later on.

\bigskip

Let us discuss now the representation theory approach to the computation of $I(a)$. The Weingarten formula reformulates, in ``rectangular form'', as follows:

\index{Weingarten formula}

\begin{theorem}
We have the Weingarten formula
$$I(a)=\sum_{\pi,\sigma}\delta_\pi(a_l)\delta_\sigma(a_r)W_{kN}(\pi,\sigma)$$
where $k=\Sigma a_{ij}/2$, and where the multi-indices $a_l/a_r$ are defined as follows:
\begin{enumerate}
\item Start with $a\in M_{p\times q}(\mathbb N)$, and replace each $ij$-entry by $a_{ij}$ copies of $i/j$.

\item Read this matrix in the usual way, as to get the multi-indices $a_l/a_r$.
\end{enumerate}
\end{theorem}

\begin{proof}
This is simply a reformulation of the Weingarten formula. Indeed, according to our definitions, the integral in the statement is given by:
$$I(a)=\int_{O_n}\underbrace{u_{11}\ldots u_{11}}_{a_{11}}\,\underbrace{u_{12}\ldots u_{12}}_{a_{12}}\,\ldots\,\underbrace{u_{pq}\ldots u_{pq}}_{a_{pq}}\,du$$

Thus, what we have here is an integral exactly as in the usual Weingarten formula, the multi-indices which are involved being as follows:
\begin{eqnarray*}
a_l&=&(\underbrace{1\ldots 1}_{a_{11}}\,\underbrace{1\ldots 1}_{a_{12}}\,\ldots\,\underbrace{p\ldots p}_{a_{pq}})\\
a_r&=&(\underbrace{1\ldots 1}_{a_{11}}\,\underbrace{2\ldots 2}_{a_{12}}\,\ldots\,\underbrace{q\ldots q}_{a_{pq}})
\end{eqnarray*}

With this in hand, the result follows now from the Weingarten formula. 
\end{proof}

We are now in position of deriving a first general result from our study. This extends the various vanishing results appearing before, as follows:

\begin{proposition}
We have $I(a)=0$, unless the matrix $a$ is ``admissible'', in the sense that all $p+q$ sums on its rows and columns are even numbers.
\end{proposition}

\begin{proof}
Observe first that the left multi-index associated to $a$ consists of $k_1=\Sigma a_{1j}$ copies of $1$, $k_2=\Sigma a_{2j}$ copies of $2$, and so on, up to $k_p=\Sigma a_{pj}$ copies of $p$. In the case where one of these numbers is odd we have $\delta_\pi(a)=0$ for any $\pi$, and this gives:
$$I(a)=0$$

A similar argument with the right multi-index associated to $a$ shows that the sums on the columns of $a$ must be even as well, and we are done.
\end{proof}

A natural question now is whether the converse of Proposition 16.24 holds, and if so, the question of computing the sign of $I(a)$ appears as well. These are both quite subtle questions, and we begin our investigations with a $N\to\infty$ study. We have here:

\begin{theorem}
The Weingarten matrix is asymptotically diagonal, in the sense that:
$$W_{kN}(\pi,\sigma)=N^{-k}(\delta_{\pi\sigma}+O(N^{-1}))$$

\noindent Moreover, the $O(N^{-1})$ remainder is asymptotically smaller that $(2k/e)^kN^{-1}$.
\end{theorem}

\begin{proof}
It is convenient, for the purposes of this proof, to drop the indices $k,N$. We know that the Gram matrix is given by $G(\pi,\sigma)=N^{|\pi\vee\sigma|}$, so we have:
$$G(\pi,\sigma)=
\begin{cases}
N^k&{\rm for\ }\pi=\sigma\\
N,N^2,\ldots,N^{k-1}&{\rm for\ }\pi\neq\sigma
\end{cases}$$

Thus the Gram matrix is of the following form, with $||H||_\infty\leq N^{-1}$:
$$G=N^k(1+H)$$
 
Now recall that for any $K\times K$ complex matrix $X$, we have the following lineup of standard inequalities, which are all elementary:
$$||X||_\infty
\leq||X||
\leq||X||_2
\leq K||X||_\infty$$

In the case of our matrix $H$, the size is $K=(2k)!!$, and we obtain in this way:
$$||H||\leq KN^{-1}$$

In order to advance, we can use now the following basic inversion formula:
$$(1+H)^{-1}=1-H+H^2-H^3+\ldots$$

We conclude from this that we have the following estimate:
$$||1-(1+H)^{-1}||\leq\frac{||H||}{1-||H||}$$

By putting now everything together, we obtain the following estimate:
\begin{eqnarray*}
||1-N^kW||_\infty
&=&||1-(1+H)^{-1}||_\infty\\
&\leq&||1-(1+H)^{-1}||\\
&\leq&||H||/(1-||H||)\\
&\leq&KN^{-1}/(1-KN^{-1})\\
&=&K/(N-K)
\end{eqnarray*}

Together with the Stirling estimate $K=(2k)!!\approx(2k/e)^k$, this gives the result.
\end{proof}

As a continuation of this, regarding this time integrals over $O_N$, we have:

\begin{theorem}
We have the estimate
$$I(a)=N^{-k}\left(\prod_{i=1}^p\prod_{j=1}^qa_{ij}!!+O(N^{-1})\right)$$
when all $a_{ij}$ are even, and $I(a)=O(N^{-k-1})$ otherwise.
\end{theorem}

\begin{proof}
By using the above results, we obtain the following estimate:
\begin{eqnarray*}
I(a)
&=&\sum_{\pi,\sigma}\delta_\pi(a_l)\delta_\sigma(a_r)W_{kN}(\pi,\sigma)\\
&=&n^{-k}\sum_{\pi,\sigma}\delta_\pi(a_l)\delta_\sigma(a_r)(\delta_{\pi\sigma}+O(N^{-1}))\\
&=&N^{-k}\left(\#\left\{\pi\Big|\delta_\pi(a_l)=\delta_\pi(a_r)=1\right\}+O(N^{-1})\right)
\end{eqnarray*}

In order to count the partitions appearing in the set on the right, it is convenient to view the multi-indices $a_l,a_r$ in a rectangular way, as follows: 
$$a_l=\begin{pmatrix}
\underbrace{1\ldots 1}_{a_{11}}&\ldots&\underbrace{1\ldots 1}_{a_{1q}}\\
\dots&\dots&\dots\\
\underbrace{p\ldots p}_{a_{p1}}&\ldots&\underbrace{p\ldots p}_{a_{pq}}$$
\end{pmatrix}\quad,\quad 
a_r=\begin{pmatrix}
\underbrace{1\ldots 1}_{a_{11}}&\ldots&\underbrace{q\ldots q}_{a_{1q}}\\
\dots&\dots&\dots\\
\underbrace{1\ldots 1}_{a_{p1}}&\ldots&\underbrace{p\ldots p}_{a_{pq}}$$
\end{pmatrix}$$

In other words, the multi-indices $a_l/a_r$ are now simply obtained from the matrix $a$ by ``dropping'' from each entry $a_{ij}$ a sequence of $a_{ij}$ numbers, all equal to $i/j$. These two multi-indices, now in matrix form, have total length $2k=\Sigma a_{ij}$. We agree to view as well any pairing of $\{1,\ldots,2k\}$ in matrix form, by following the same convention. With this picture, the pairings $\pi$ which contribute are simply those interconnecting sequences of indices ``dropped'' from the same $a_{ij}$, and this gives the following results:

\medskip

(1) In the case where one of the entries $a_{ij}$ is odd, there is no pairing that can contribute to the leading term under consideration, so we have $I(a)=O(N^{-k-1})$, and we are done.

\medskip

(2) In the case where all the entries $a_{ij}$ are even, the pairings that contribute to the leading term are those connecting points inside the $pq$ ``dropped'' sets, i.e. are made out of a pairing of $a_{11}$ points, a pairing of $a_{12}$ points, and so on, up to a pairing of $a_{pq}$ points. Now since an $x$-point set has $x!!$ pairings, this gives the formula in the statement. 
\end{proof}

In order to further advance, let us formulate a key definition, as follows:

\index{Brauer space}

\begin{definition}
The Brauer space $D_k$ is defined as follows:
\begin{enumerate}
\item The points are the Brauer diagrams, i.e. the pairings of $\{1,2,\ldots,2k\}$. 

\item The distance function is given by $d(\pi,\sigma)=k-|\pi\vee\sigma|$.
\end{enumerate}
\end{definition}

It is indeed well-known, and elementary to check, that $d$ satisfies the usual axioms for a distance function. This is something standard, and heavily used in probability theory, and for some comments and examples here, we refer to \cite{bco}, \cite{csn} and related papers. Now the point is that we have a series expansion of the Weingarten function in terms of paths on the Brauer space, originally found by Collins in \cite{col} in the unitary case, then by Collins and \'Sniady \cite{csn} in the orthogonal case. We present here a slightly modified statement, along with a complete proof, by using a somewhat lighter formalism: 

\begin{theorem}
The Weingarten function $W_{kN}$ has a series expansion in $N^{-1}$,
$$W_{kN}(\pi,\sigma)=N^{-k-d(\pi,\sigma)}\sum_{g=0}^\infty K_g(\pi,\sigma)N^{-g}$$
where the objects on the right are defined as follows:
\begin{enumerate}
\item A path from $\pi$ to $\sigma$ is a sequence $p=[\pi=\tau_0\neq\tau_1\neq\ldots\neq\tau_r=\sigma]$.

\item The signature of such a path is $+$ when $r$ is even, and $-$ when $r$ is odd.

\item The geodesicity defect of such a path is $g(p)=\Sigma_{i=1}^rd(\tau_{i-1},\tau_i)-d(\pi,\sigma)$.

\item $K_g$ counts the signed paths from $\pi$ to $\sigma$, with geodesicity defect $g$.
\end{enumerate} 
\end{theorem}

\begin{proof}
Let us go back to the proof of our main estimate so far, established in the above. We can write the Gram matrix in the following way:
$$G_{kn}=N^{-k}(1+H)$$

In terms of the Brauer space distance, the formula of $H$ is simply:
$$H(\pi,\sigma)=
\begin{cases}
0&{\rm for\ }\pi=\sigma\\
N^{-d(\pi,\sigma)}&{\rm for\ }\pi\neq\sigma
\end{cases}$$

Consider now the set $P_r(\pi,\sigma)$ of $r$-paths between $\pi$ and $\sigma$. According to the usual rule of matrix multiplication, the powers of $H$ are given by:
$$H^r(\pi,\sigma)
=\sum_{p\in P_r(\pi,\sigma)}H(\tau_0,\tau_1)\ldots H(\tau_{r-1},\tau_r)
=\sum_{p\in P_r(\pi,\sigma)}N^{-d(\pi,\sigma)-g(p)}$$

We can use now the following standard inversion formula:
$$(1+H)^{-1}=1-H+H^2-H^3+\ldots$$

By using this formula, we obtain the following equality:
\begin{eqnarray*}
W_{kN}(\pi,\sigma)
&=&N^{-k}\sum_{r=0}^\infty(-1)^rH^r(\pi,\sigma)\\
&=&N^{-k-d(\pi,\sigma)}\sum_{r=0}^\infty\sum_{p\in P_r(\pi,\sigma)}(-1)^rN^{-g(p)}
\end{eqnarray*}

Now by rearranging the various terms of the double sum according to their geodesicity defect $g=g(p)$, this gives the following formula:
$$W_{kN}(\pi,\sigma)=N^{-k-d(\pi,\sigma)}\sum_{g=0}^\infty K_g(\pi,\sigma)N^{-g}$$

Thus, we have obtained the formula in the statement.
\end{proof}

In order to discuss now the $I(a)$ reformulation of the above result, it is convenient to use the total length of a path, defined as follows:
$$d(p)=\sum_{i=1}^rd(\tau_{i-1},\tau_i)$$

Observe that, in terms of this quantity, we have the following formula:
$$d(p)=d(\pi,\sigma)+g(p)$$

With these conventions, we have the following result:

\begin{theorem}
The integral $I(a)$ has a series expansion in $N^{-1}$ of the form
$$I(a)=N^{-k}\sum_{d=0}^\infty H_d(a)N^{-d}$$
where the coefficient on the right can be interpreted as follows:
\begin{enumerate}
\item Starting from $a\in M_{p\times q}(\mathbb N)$, construct the multi-indices $a_l,a_r$ as usual.

\item Call a path ``$a$-admissible'' if its endpoints satisfy $\delta_\pi(a_l)=1$ and $\delta_\sigma(a_r)=1$.

\item Then $H_d(a)$ counts all $a$-admissible signed paths in $D_k$, of total length $d$.
\end{enumerate}
\end{theorem}

\begin{proof}
We can combine first the above results, in the following way:
\begin{eqnarray*}
I(a)
&=&\sum_{\pi,\sigma}\delta_\pi(a_l)\delta_\sigma(a_r)W_{kN}(\pi,\sigma)\\
&=&N^{-k}\sum_{\pi,\sigma}\delta_\pi(a_l)\delta_\sigma(a_r)\sum_{g=0}^\infty K_g(\pi,\sigma)N^{-d(\pi,\sigma)-g}
\end{eqnarray*}

Let us denote by $H_d(\pi,\sigma)$ the number of signed paths between $\pi$ and $\sigma$, of total length $d$. In terms of the new variable $d=d(\pi,\sigma)+g$, the above expression becomes:
\begin{eqnarray*}
I(a)
&=&N^{-k}\sum_{\pi,\sigma}\delta_\pi(a_l)\delta_\sigma(a_r)\sum_{d=0}^\infty H_d(\pi,\sigma)N^{-d}\\
&=&N^{-k}\sum_{d=0}^\infty\left(\sum_{\pi,\sigma}\delta_\pi(a_l)\delta_\sigma(a_r)H_d(\pi,\sigma)\right)N^{-d}
\end{eqnarray*}

We recognize in the middle the quantity $H_d(a)$, and this gives the result.
\end{proof}

We derive now some concrete consequences from the abstract results in the previous section. First, let us recall the following result, due to Collins and \'Sniady \cite{csn}:

\begin{theorem}
We have the estimate
$$W_{kN}(\pi,\sigma)=N^{-k-d(\pi,\sigma)}(\mu(\pi,\sigma)+O(N^{-1}))$$
where $\mu$ is the M\"obius function.
\end{theorem}

\begin{proof}
We know from the above that we have the following estimate:
$$W_{kN}(\pi,\sigma)=N^{-k-d(\pi,\sigma)}(K_0(\pi,\sigma)+O(N^{-1}))$$

Now since one of the possible definitions of the M\"obius function is that this counts the signed geodesic paths, we have $K_0=\mu$, and we are done.
\end{proof}

Let us go back now to our integrals $I(a)$. We have the following result:

\begin{theorem}
We have the estimate
$$I(a)=N^{-k-e(a)}(\mu(a)+O(N^{-1}))$$
where the objects on the right are as follows:
\begin{enumerate}
\item $e(a)=\min\left\{d(\pi,\sigma)\big|\pi,\sigma\in D_k,\delta_\pi(a_l)=\delta_\sigma(a_r)=1\right\}$. 

\item $\mu(a)$ counts all $a$-admissible signed paths in $D_k$, of total length $e(a)$.
\end{enumerate} 
\end{theorem}

\begin{proof}
We know that we have an estimate of the following type:
$$I(a)=N^{-k-e}(H_e(a)+O(N^{-1}))$$

Here, according to the various notations above, $e\in\mathbb N$ is the smallest total length of an $a$-admissible path, and $H_e(a)$ counts all signed $a$-admissible paths of total length $e$. Now since the smallest total length of such a path is of course attained when the path is just a segment, we have $e=e(a)$ and $H_e(a)=\mu(a)$, and we are done.
\end{proof}

At a more advanced level now, and still on the same topic, integration over $O_N$, we have the following result, due to Collins-Matsumoto \cite{cma} and Zinn-Justin \cite{zin}:

\begin{theorem}
We have the formula
$$W_{kn}(\pi,\sigma)=\frac{\sum_{\lambda\vdash k,\,l(\lambda)\leq k}\chi^{2\lambda}(1_k)w^\lambda(\pi^{-1}\sigma)}{(2k)!!\prod_{(i,j)\in\lambda}(n+2j-i-1)}$$
where the various objects on the right are as follows:
\begin{enumerate}
\item The sum is over all partitions of $\{1,\ldots,2k\}$ of length $l(\lambda)\leq k$.

\item $w^\lambda$ is the corresponding zonal spherical function of $(S_{2k},H_k)$.

\item $\chi^{2\lambda}$ is the character of $S_{2k}$ associated to $2\lambda=(2\lambda_1,2\lambda_2,\ldots)$.

\item The product is over all squares of the Young diagram of $\lambda$.
\end{enumerate}
\end{theorem}

\begin{proof}
This is something quite technical, that we will not attempt to explain here, and for details on all this, we refer to the papers \cite{cma}, \cite{zin}.
\end{proof}

It is of course possible to deduce from this a new a formula for the integrals $I(a)$, just by putting together the various formulae that we have. Let us just record here:

\begin{theorem}
The possible poles of $I(a)$ can be at the numbers
$$-(k-1),-(k-2),\ldots,2k-1,2k$$
where $k\in\mathbb N$, associated to the admissible matrix $a\in M_{p\times q}(\mathbb N)$ is given by $k=\sum a_{ij}/2$.
\end{theorem}

\begin{proof}
We know from the above that the possible poles of $I(a)$ can only come from those of the Weingarten function. On the other hand, Theorem 16.32 tells us that these latter poles are located at the numbers of the form $-2j+i+1$, with $(i,j)$ ranging over all possible squares of all possible Young diagrams, and this gives the result.
\end{proof}

As a last topic, let us discuss Gram determinants. In what regards the symmetric group $S_N$, we have the following result, that we already know, from the above:

\begin{theorem}
The determinant of the Gram matrix of $S_N$ is given by
$$\det(G_{kN})=\prod_{\pi\in P(k)}\frac{N!}{(N-|\pi|)!}$$
with the convention that in the case $N<k$ we obtain $0$.
\end{theorem}

\begin{proof}
This is something that we know, the idea being that $G_{kN}$ naturally decomposes as a product of an upper triangular and lower triangular matrix.
\end{proof}

Let us discuss now the case of the orthogonal group $O_N$. Here the combinatorics is that of the Young diagrams. We denote by $|.|$ the number of boxes, and we use quantity $f^\lambda$, which gives the number of standard Young tableaux of shape $\lambda$. With these conventions, the result, which is something quite technical, is then as follows:

\begin{theorem}
The determinant of the Gram matrix of $O_N$ is given by
$$\det(G_{kN})=\prod_{|\lambda|=k/2}f_N(\lambda)^{f^{2\lambda}}$$
where the quantities on the right are $f_N(\lambda)=\prod_{(i,j)\in\lambda}(N+2j-i-1)$.
\end{theorem}

\begin{proof}
This follows from the results of Zinn-Justin in \cite{zin}. Indeed, it is known from there that the Gram matrix is diagonalizable, as follows:
$$G_{kN}=\sum_{|\lambda|=k/2}f_N(\lambda)P_{2\lambda}$$

Here $1=\sum P_{2\lambda}$ is the standard partition of unity associated to the Young diagrams having $k/2$ boxes, and the coefficients $f_N(\lambda)$ are those in the statement. Now since we have $Tr(P_{2\lambda})=f^{2\lambda}$, this gives the result. See \cite{bcu}, \cite{zin}.
\end{proof}

Finally, since it is late, and time to sleep, and no algebra book would be complete without some quantum groups at the end, let us discuss this. Unfortunately, we are here, with our Gram determinants, into quite advanced things, so we will have to trick a bit, and take some dirty shortcuts. Let us start with a definition, informal as they come:

\begin{definition}
In analogy with the fact that $S_N,O_N$ are easy, coming from $P,P_2$, let us denote by $S_N^+,O_N^+$ the formal objects associated to $NC,NC_2$.
\end{definition}

Observe that $S_N^+,O_N^+$ cannot be groups, because $NC,NC_2$ do not contain the basic crossing $\slash\hskip-2mm\backslash$, and so are not categories of partitions in the sense of chapter 15. This being said, the axiom stating that $\slash\hskip-2mm\backslash$ must be in the category was coming from the fact that the coordinates $u_{ij}:G\to\mathbb C$ of a compact Lie group $G\subset_uU_N$ commute, so in the lack of this axiom, we can only have some kind of ``quantum groups'', which are beasts a bit like groups, save for the fact that the coordinates $u_{ij}:G\to\mathbb C$ do not longer commute.

\bigskip

Anyway. Getting now to business, we would like to compute the  Gram determinants for $S_N^+,O_N^+$. Following Di Francesco \cite{dif}, let us begin with some examples:

\begin{proposition}
At $k=2$ the set of partitions for $S_N^+$ is $NC(2)=\{||,\sqcap\}$, and the corresponding Gram matrix and its determinant are:
$$\det\begin{pmatrix}N^2&N\\N&N\end{pmatrix}=N^2(N-1)$$
Also, at $k=4$ the set of partitions for $O_N^+$ is $NC_2(4)=\{\sqcap\sqcap,\bigcap\hskip-4.9mm{\ }_\cap\,\}$, and the corresponding Gram matrix and its determinant are:
$$\det\begin{pmatrix}N^2&N\\N&N^2\end{pmatrix}=N^2(N^2-1)$$
\end{proposition}

\begin{proof}
This is something which is indeed clear from definitions.
\end{proof}

With a few tricks, we can work out as well the next computation, as follows:

\begin{proposition}
At $k=3$ the partition set for $S_N^+$ is $NC(3)=\{|||,\sqcap|,\sqcap\hskip-3.2mm{\ }_|\,,|\sqcap,\sqcap\hskip-0.7mm\sqcap\}$, and the corresponding Gram matrix and its determinant are:
$$\det\begin{pmatrix}
N^3&N^2&N^2&N^2&N\\
N^2&N^2&N&N&N\\
N^2&N&N^2&N&N\\
N^2&N&N&N^2&N\\
N&N&N&N&N
\end{pmatrix}=N^5(N-1)^4(N-2)$$
Also, at $k=6$ the set of partitions for $O_N^+$ is $NC_2(6)\simeq NC(3)$, and the corresponding Gram matrix and its determinant are:
$$\det\begin{pmatrix}
N^3&N^2&N^2&N^2&N\\
N^2&N^3&N&N&N^2\\
N^2&N&N^3&N&N^2\\
N^2&N&N&N^3&N^2\\
N&N^2&N^2&N^2&N^3
\end{pmatrix}=N^5(N^2-1)^4(N^2-2)$$
\end{proposition}

\begin{proof}
We have two formulae to be proved, the idea being as follows:

\medskip

(1) In what regards $S_N^+$, the set of partitions here is $NC(3)=P(3)$, and so the corresponding Gram matrix is the one in the statement, exactly as for $S_N$. By using the Lindst\"om formula, from Theorem 16.12, the determinant of this matrix is, as claimed:
\begin{eqnarray*}
\det
&=&\prod_{\pi\in P(3)}\frac{N!}{(N-|\pi|)!}\\
&=&\frac{N!}{(N-3)!}\left(\frac{N!}{(N-2)!}\right)^3\frac{N!}{(N-1)!}\\
&=&N(N-1)(N-2)N^3(N-1)^3N\\
&=&N^5(N-1)^4(N-2)
\end{eqnarray*}

(2) Regarding now $O_N^+$, the set of partitions here is $NC_2(6)$, and by using the  fattening/shrinking identification $NC_2(6)\simeq NC(3)$, we obtain, by using (1):
\begin{eqnarray*}
\det
&=&\frac{1}{N^2\sqrt{N}}\times N^{10}(N^2-1)^4(N^2-2)\times\frac{1}{N^2\sqrt{N}}\\
&=&N^5(N^2-1)^4(N^2-2)
\end{eqnarray*}

Thus, we have obtained the formula in the statement.
\end{proof}

In general now, following \cite{dif}, we have the following result:

\index{meander determinant}
\index{Gram determinant}

\begin{theorem}
The determinant of the Gram matrix for $O_N^+$ is given by
$$\det(G_{kN})=\prod_{r=1}^{[k/2]}P_r(N)^{d_{k/2,r}}$$
where $P_r$ are the Chebycheff polynomials, given by
$$P_0=1\quad,\quad 
P_1=X\quad,\quad 
P_{r+1}=XP_r-P_{r-1}$$
and $d_{kr}=f_{kr}-f_{k,r+1}$, with $f_{kr}$ being the following numbers, depending on $k,r\in\mathbb Z$,
$$f_{kr}=\binom{2k}{k-r}-\binom{2k}{k-r-1}$$
with the convention $f_{kr}=0$ for $k\notin\mathbb Z$. 
\end{theorem}

\begin{proof}
This is something quite heavy, and we refer here to Di Francesco \cite{dif}.
\end{proof}

Also following \cite{dif}, we have as well the following result:

\begin{theorem}
The determinant of the Gram matrix for $S_N^+$ is given by
$$\det(G_{kN})=(\sqrt{N})^{a_k}\prod_{r=1}^kP_r(\sqrt{N})^{d_{kr}}$$
where $d_{kr}=f_{kr}-f_{k,r+1}$, with $f_{kr}$ being the following numbers, depending on $k,r\in\mathbb Z$,
$$f_{kr}=\binom{2k}{k-r}-\binom{2k}{k-r-1}$$
with the convention $f_{kr}=0$ for $k\notin\mathbb Z$, and where $a_k=\sum_{\pi\in \mathcal P(k)}(2|\pi|-k)$.
\end{theorem}

\begin{proof}
Again, heavy mathematics, and we refer here to Di Francesco \cite{dif}.
\end{proof}

We refer to \cite{bcu}, \cite{dif}, for a further discussion on these topics. 

\section*{16e. Exercises}

Congratulations for having read this book, and no exercises for this final chapter. But you can try instead to read some of the books and articles referenced below.

\baselineskip=14pt

\printindex

\end{document}